\documentclass[11pt,letterpaper,openany]{book}
\usepackage[margin=1.15in]{geometry}
\usepackage{amsmath,amssymb,amsthm}
\usepackage[hidelinks]{hyperref}
\usepackage{tikz}
\usepackage{xcolor}
\usepackage{booktabs}
\usepackage{microtype}
\usetikzlibrary{arrows.meta,calc,backgrounds}
\theoremstyle{plain}
\newtheorem{definition}{Definition}[section]
\newtheorem{notation}[definition]{Notation}
\newtheorem{theorem}[definition]{Theorem}
\newtheorem{lemma}[definition]{Lemma}
\newtheorem{proposition}[definition]{Proposition}
\newtheorem{corollary}[definition]{Corollary}
\newtheorem{example}[definition]{Example}
\newtheorem{remark}[definition]{Remark}
\newtheorem{assumption}[definition]{Assumption}
\newcommand{\Attain}{\DD_m}

\newcommand{\Cm}{\mathcal{C}}

\newcommand{\DD}{\mathcal{D}}
\newcommand{\Dm}{\Delta^m}
\newcommand{\Domain}{\mathcal{A}}
\newcommand{\EE}{\mathbb{E}}
\newcommand{\Gm}{\Gamma_m}
\newcommand{\HH}{\mathcal{H}}
\newcommand{\Hull}{\operatorname{Hull}}
\newcommand{\Int}{\operatorname{Int}}

\newcommand{\LatL}{L^m}

\newcommand{\LbarCal}{\bar{\mathcal{L}}}
\newcommand{\LbarRaw}{\bar L}
\newcommand{\Lm}{\mathcal{L}_m}

\newcommand{\PP}{\mathcal{P}}
\newcommand{\Pmeas}{\mathbb{P}}
\newcommand{\RR}{\mathbb{R}}
\newcommand{\Span}{\operatorname{span}}
\newcommand{\StateA}{A_m}
\newcommand{\StateAm}{{\mathcal A}_m}
\newcommand{\TT}{\mathbb{T}}

\newcommand{\Yhn}{\widehat{Y}_n}

\newcommand{\ZZ}{\mathbb{Z}}
\newcommand{\abs}[1]{\left|#1\right|}

\newcommand{\bU}{\overline U}
\newcommand{\bV}{\mathbb{V}}
\newcommand{\dd}{\,\mathrm{d}}

\newcommand{\ext}{\operatorname{Ext}}
\newcommand{\grad}{\nabla}
\newcommand{\half}{\tfrac{1}{2}}
\newcommand{\id}{\mathbb{I}}
\newcommand{\ii}{\mathrm{i}}
\newcommand{\inner}[2]{\langle #1,#2\rangle}
\newcommand{\lap}{\Delta}
\newcommand{\nabm}{\nabla^m}
\newcommand{\newnote}[1]{}
\newcommand{\norm}[1]{\left\lVert#1\right\rVert}
\newcommand{\nuc}{\nu}

\newcommand{\one}{\mathbf{1}}
\newcommand{\rev}[1]{#1}
\newcommand{\sgn}{\operatorname{sgn}}

\newcommand{\tA}{\tilde A}
\newcommand{\tPi}{\tilde\Pi}
\newcommand{\tU}{\widetilde U}
\newcommand{\tr}{\operatorname{tr}}
\newcommand{\vb}{\boldsymbol{b}}
\newcommand{\vp}{\boldsymbol{p}}
\newcommand{\vx}{\boldsymbol{x}}

\makeatletter
\let\book@chapter\chapter
\renewcommand{\chapter}{\@ifstar{\book@chapterstar}{\book@chapter}}
\newcommand{\book@chapterstar}[1]{\book@chapter*{#1}\addtocontents{toc}{\protect\contentsline{chapter}{\protect\numberline{}#1}{\thepage}{}}}
\let\book@section\section
\renewcommand{\section}{\@ifstar{\book@sectionstar}{\book@section}}
\newcommand{\book@sectionstar}[1]{\book@section*{#1}\addtocontents{toc}{\protect\contentsline{section}{\protect\numberline{}#1}{\thepage}{}}}
\let\book@subsection\subsection
\renewcommand{\subsection}{\@ifstar{\book@subsectionstar}{\book@subsection}}
\newcommand{\book@subsectionstar}[1]{\book@subsection*{#1}\addtocontents{toc}{\protect\contentsline{subsection}{\protect\numberline{}#1}{\thepage}{}}}
\let\book@subsubsection\subsubsection
\renewcommand{\subsubsection}{\@ifstar{\book@subsubsectionstar}{\book@subsubsection}}
\newcommand{\book@subsubsectionstar}[1]{\book@subsubsection*{#1}\addtocontents{toc}{\protect\contentsline{subsubsection}{\protect\numberline{}#1}{\thepage}{}}}
\makeatother
\hypersetup{pdftitle={Constructive Probability: Foundations and Numerical Analysis for Mathematical Finance and Machine Learning},pdfauthor={Claudio Albanese}}
\providecommand{\Layout}{}
\if\relax\detokenize\expandafter{\Layout}\relax\else\fi
\usepackage{listings}
\usepackage{longtable}
\begin{document}
\frontmatter
\begin{titlepage}
\centering
\vspace*{3cm}
{\Huge\bfseries Constructive Probability}\\[14pt]
{\Large Foundations and Numerical Analysis for}\\[2pt]
{\Large mathematical finance and machine learning}\\[2cm]
{\large Claudio Albanese}\\[4pt]
{\small\texttt{www.esther-analytics.com}}\\[2cm]
\vfill
{\small First edition, October 2026}
\end{titlepage}
% ---- from chapters/dedication.tex ----
% Authored: front matter, written for the book.
\thispagestyle{empty}
\vspace*{\fill}
\begin{center}
\itshape
To St\'ephane Cr\'epey,\\[6pt]
for his friendship, and for the inspiration which has run both ways\\[2pt]
between us for as long as I have known him.
\end{center}
\vspace*{\fill}
\clearpage
% ---- end chapters/dedication.tex ----
% ---- from chapters/preface.tex ----
% Authored: front matter, written for the book.
\chapter*{Preface}

The title uses the word constructive in the sense in which it is used in constructive field
theory, and not in the sense in which it is used in mathematical logic. It does not refer to
intuitionistic reasoning, or to any restriction of the methods of proof: the logic of this book is
classical throughout, the real numbers are the classical real numbers, and every estimate in it is
an ordinary inequality between ordinary numbers. What the word refers to is the manner in which a
continuum object is obtained --- as the limit of a sequence of discrete approximations, with the
passage to the limit controlled by estimates which can themselves be checked, and with the
renormalisation which that passage requires exhibited rather than assumed. That is the programme of
constructive quantum field theory, and this book belongs to it, applying the same discipline to
Markov kernels rather than to Euclidean fields.

What is set aside is not continuum mathematics, of which the book makes free use, but the
sigma-algebra on a continuum as the primitive object of a probability theory. Probability enters
here through the prevision of finitely many gambles at a fixed resolution, which is de Finetti's
primitive, and the passage to the continuum is organised by a sequence of lattice refinements and
by the estimates which control them.

The theorem at the centre of the book is one statement with several names. It is de Finetti's
fundamental theorem of probability theory; it is the fundamental theorem of Bayesian statistics; and
it is the fundamental theorem of asset pricing in its one-period form. Throughout the book it is
called simply the fundamental theorem. The vocabulary used to discuss it is financial, for
convenience only: every statement could be written in the language of previsions, or in that of
priors and posteriors, with no change to the mathematics.

The applications are two --- mathematical finance and machine learning --- and they are served in
parallel, because they are the same mathematics read in two directions. In finance an assessment is a
price, and coherence is the absence of arbitrage; in learning an assessment is a prediction, and
coherence is that the prediction be the conditional expectation of what it predicts, hence the best
available in mean square. Both objects are conditional expectations of a functional of a path given
the information available at the time; both are computed by the same adjoint recursion, which is the
adjoint of a propagation in one case and backpropagation in the other; and both are wrong by the same
kind of amount when the path is discretised, which is the subject of the second half of the book. The
results are therefore stated once and read twice.

Neither subject can treat the analysis as an accessory. In both, what is delivered is a computed number, and
a computed number without a bound on its discretisation error is not a claim: a price without an error bar is
an opinion with a decimal point, and a prediction without one is a guess. The two disciplines share no
vocabulary, and they share the estimates, because both live or die by them.

The reason the parallel is not commonly drawn is not that it is deep but that the notation hides it. In
the formalism now standard in mathematical finance the primitive object is a probability space carrying a
filtration, and the two subjects are then built on different objects --- a martingale measure for the one, a
statistical model for the other --- with the theorem that connects them stated in two vocabularies and proved
twice. Take instead the resolution as primitive: finitely many gambles, a prevision, a lattice, a controlled
limit. The objects then coincide, because a coherent assessment is a conditional expectation whether it is
quoted as a price or as a prediction, and the recursion which computes it is one recursion. The identity
itself is not a new theorem --- it is de Finetti's, as the next paragraph says --- but the analysis which makes
it visible is new: the estimates that control the passage to the continuum, and the renormalisation they
require, are what turn an observation about notation into a computation, and what make the two halves of the
book wrong by the same amount when the path is discretised. An identification asserted at the level of
formalism costs nothing and convinces nobody; one that survives the estimates has content. That is the reason
for insisting on this foundation, and the reason the two halves of the book are written side by side rather
than in two books.

The two subjects grew up on parallel tracks. Mathematical finance and machine learning were built, for
most of the period in which either reached maturity, by communities which did not read each other:
different journals, different training, different software. The consequence is not that the
mathematics differs but that the vocabulary does, and it has drifted far enough that the
correspondence has to be discovered rather than recognised. It was discovered here in the way such
things usually are, by building a framework from first principles --- gambles, previsions, a lattice,
a controlled limit --- and then noticing that the objects which came out had the same shape as the
objects on the other side. That is stronger evidence than a translation table could be: a translation
asserts a correspondence, whereas a construction has to survive its own estimates. What came out is not
a resemblance but an overlap, and it is large. The adjoint recursion is one recursion, the error of a
discretised path is one error, and the bound which makes a computed number a claim is one bound.

The overlap is also a channel, at the level of architecture rather than of code. Because the objects
coincide, the shape of the machine which computes them is the same on both sides: hardware
acceleration plays the same decisive role, differentiation is a first-class operation in the toolchain
rather than an afterthought, the whole universe of scenarios is evaluated jointly rather than one at a
time, and the compute core is held apart from the logic which describes it. None of this is a claim
about software. Implementations of one scheme differ in every detail, and the details are where the
engineering is; the claim is about the shape, and the shape is dictated by the mathematics. What the
two domains do not share is the \emph{principle of the inference}. A model here is projected forward
from the collective expectations of the participants in a market, which prices record; a model there
is fitted to a history which has been observed. That difference carries the validation with it, since
a forecast from history is an answer without an error bar and scale does not supply one
--- the bar is a statement about the discretisation and not about the fit --- whereas the numbers this
book computes come with their bound. It is the reason to treat these as one mathematics with two
vocabularies rather than as two stacks with an interface between them, and the reason a book which
writes the estimates once for both is not a survey of a coincidence but the case for a single engine.

The credit for the unity is de Finetti's, and it is not a coincidence of vocabulary. Probability, for
him, was the price of a gamble; coherence, the requirement that prices be consistent, made that price
a probability; and exchangeability, the assumption that the order of the observations carries no
information, gave the mixture representation which is at once a prior over models and a market of
claims. Both disciplines in this book are the working out of that one construction --- finance in the
direction of prices, machine learning in the direction of predictions --- and if the two are found to
share a foundation, the foundation is his.

\vspace{1em}
\noindent\emph{October 2026}\\[2pt]
\noindent Claudio Albanese\\[2pt]
\noindent\texttt{www.esther-analytics.com}
% ---- end chapters/preface.tex ----
\tableofcontents
\mainmatter
% ---- from chapters/ch1.tex ----
% Authored: Chapter 1 is written for the book and draws on the parts rather than
% reproducing them. Still to migrate from generated/ch1-body.tex: the sense in
% which the theory is constructive, absoluteness, decidability, what the theory
% buys in practice, the programme, Kolmogorov, mechanisation, and the plan of the
% book. The de Finetti passage of that body is duplicated by the section below
% and is dropped when the migration reaches it.
\chapter{Introduction}
\label{chap:ch1}

\section{``Probability Does Not Exist''}

The sentence is de Finetti's \cite{deFinetti1974}, and he put it at the head of his treatise as a
statement of method rather than as a provocation. What it denies is not the usefulness of the
calculus of probability but the existence of the thing which the calculus is usually taken to be
about. On the account he develops there is no quantity, attaching to an event or to a proposition,
which is \emph{its} probability; there are only the degrees of belief of a person, and the
constraints which coherence places on them. The constraints turn out to be strong enough to
generate the whole calculus, and weak enough to leave the beliefs where they belong, which is to
say with the person who holds them.

\subsection{Prevision as the Primitive}

The primitive object of the theory is not a set, and not a measure, but a bet and the price at
which one is willing to take either side of it. If a person will buy or sell a gamble $X$ at the
price $p$, and is indifferent to which side he takes, then $p$ is his \emph{prevision} of $X$; an
assessment on events is the special case in which the gambles are indicator functions. A family of
such prices is \emph{coherent} when no finite combination of the bets, in the quantities to which
the person has committed himself, yields a gain which is sure and positive. Coherence is the only
requirement imposed, and everything else --- additivity, the multiplication rule, the existence of
conditional previsions --- is derived from it.

That derivation is the content of de Finetti's fundamental theorem of probability
\cite{deFinetti1931,deFinetti1937}: an assessment on a Boolean algebra is coherent if and only if
it extends to a finitely additive probability on that algebra. The name is de Finetti's own, and
it is apt, because the theorem stands to the calculus of probability as the fundamental theorem of
asset pricing stands to mathematical finance: it is the same alternative, read once as a statement
about belief and once as a statement about prices. Chapter~\ref{chap:ch2} proves it in the finite
case, as the orthant alternative of finite dimensional convex geometry, together with the algorithm
which decides it and produces, when the assessment is incoherent, the combination of bets which
refutes it.

The identification is not a figure of speech, and it is worth setting out in full at the outset,
because everything in the book is a consequence of it. Let the state space $S$ be finite --- the
general statement is to come --- and let $\xi^1,\dots,\xi^n$ be gambles, a gain being a payoff less
the price paid for it, so that a gamble is a claim which costs nothing to enter.

\begin{theorem}[Coherence and the state price density]\label{thm:coherence}
Let $S$ be a finite set and let $\xi^1,\dots,\xi^n\in\RR^S$. The following are equivalent.
\begin{enumerate}
\item[(i)] No combination of the gambles, in any real quantities, produces a gain which is non
negative at every $s\in S$ and not identically zero. This is the coherence of de Finetti: no
admissible combination of bets yields a sure gain.
\item[(ii)] There is a strictly positive $g:S\to(0,\infty)$ with $\sum_{s\in S}g(s)\xi^k(s)=0$ for
$k=1,\dots,n$. This is the existence of a state price density: an assessment which extends to a
finitely additive probability under which every gamble is fair.
\end{enumerate}
\end{theorem}

\begin{proof}
Let $V=\Span\{\xi^1,\dots,\xi^n\}\subseteq\RR^S$. By the definition of $V$, a nonzero element of
$V\cap C_+$ is precisely a combination of the gambles which is non negative and not identically
zero, which is what (i) forbids; so (i) is the statement $V\cap C_+=\{0\}$. Statement (ii) says
that $V^{\perp}$ contains a strictly positive vector. That these two statements are equivalent is
the orthant alternative, proved in \S\ref{ch2:sec:convex} as Stiemke's lemma
(Theorem~\ref{ch2:thm:stiemke}) by an argument which terminates, and which produces the combination
of gambles refuting (i) whenever (ii) fails.
\end{proof}

Read (i) as a statement about belief and (ii) as a statement about prices, and the theorem is
de Finetti's fundamental theorem of probability and the fundamental theorem of asset pricing at
once: the same alternative, in the same finite dimensional space, with the vector $g$ read as a
prevision in the first reading and as a state price density in the second. It is stated again, and
there proved, in Chapter~\ref{chap:ch2}, where it is called the fundamental theorem; in
Chapter~\ref{chap:ch3} its form at a single site is applied at every node of a lattice model
(\S\ref{ch3:ss:oneperiod}). It appears here first, in the language in which both readings are
visible, because the finite dimensional convex geometry is the engine of the whole book and it is
worth seeing at the outset what it buys. What it buys first is a separation of the two ways
in which a model is used, and of what determines each of them; that separation is drawn in
\S\ref{sec:two-problems}, at the end of this chapter, and it is the organising principle of the
book.

\subsection{The Two Memoirs of the Nineteen Thirties}
\label{sec:memoirs}

The theory of the preceding paragraphs is de Finetti's, and it is worth naming the two memoirs in
which it was laid out, because this book depends on both and on different parts of each.

The first is \emph{Sul significato soggettivo della probabilit\`a} \cite{deFinetti1931}, published in
\emph{Fundamenta Mathematicae} in 1931. Its subject is the one this chapter has been describing: the
prevision as the primitive, coherence as the only requirement, and the theorem that a coherent
assessment on a finite algebra is precisely one which extends to a finitely additive probability. The
argument is conducted in the language of bets and of combinations of bets, and in that language it is
the argument which mathematical finance would later recognise as the absence of arbitrage and the
existence of a state price density. The identity is not a retrospective reading of that memoir: the two
statements are one alternative in a finite dimensional vector space, which is the observation around
which this book is built and which Theorem~\ref{thm:coherence} states, and Chapter~\ref{chap:ch2}
proves with the algorithm that produces the refuting combination when the alternative fails.

The second is \emph{La pr\'evision: ses lois logiques, ses sources subjectives} \cite{deFinetti1937},
published in the \emph{Annales de l'Institut Henri Poincar\'e} in 1937. Its subject is what becomes of
coherence when a symmetry is added: exchangeability, and the representation theorem which says that an
exchangeable assessment on a sequence is a mixture of independent and equally likely ones. It is
there that the prior appears as the mixing measure, the posterior as a conditional prevision, and the
predictive law as what the two together determine. Bayesian statistics is, in this reading, a
consequence of coherence together with exchangeability, and not an independent scheme to which
coherence is subsequently added.

The two memoirs are one theory read twice, once as a theory of prices and once as a theory of belief.
In what follows the two are called the pricing reading and the learning reading, the names being
attached to the direction of time in which the data lie rather than to a subject matter, and the
theorem of Chapter~\ref{chap:ch2} is the finite and decidable form of both. The computational case
belongs to the second: a system which produces an assessment is answerable to the same standard, and
the decision theoretic form of that standard is the class of proper scoring rules, set out in
\S\ref{ch2:sec:ml}. What the reader should carry from the present chapter is the form of the claim:
there is one requirement, coherence; the two disciplines differ in what is assessed and in what is
observed, and not in what is required of the assessment.

\subsection{The Boundary at Which the Programme Stopped}

Two features of the theorem deserve to be separated at the outset, because the whole of this book
depends on which of them is taken as its content.

The first is that the theorem is a statement of finite dimensional linear algebra. In the finite
case the extension is a vector in the positive orthant, the incoherence is a certificate, and
neither requires a limit, a compactness argument or a measurable selection. This is not a
convenience of presentation but the substance of the subjectivist position: a prevision is a
commitment about a finite family of bets, and a commitment can be examined.

The second is that the additivity obtained is finite, and only finite. The theorem says nothing
about countable families of events, nothing about limits of assessments, and nothing about
processes indexed by a continuum of times. This is where the programme stopped. de Finetti was
scrupulous about it: he did not claim a theory of continuous time processes, and he did not pretend
that finite additivity extends itself. The measure theoretic school resolved the difficulty by
postulating a $\sigma$-algebra on a continuum and proving an extension theorem, a route which is
coherent on its own terms and which this book does not dispute. The difference between the two
routes is the subject of the section which follows. What matters here is that the stopping point is
precisely located: coherence is a condition on finite families, and the passage to infinitely many
times is a separate problem, which de Finetti left open.

\subsection{Why Subjectivism Is a Discipline Rather than a Licence}

It is worth disposing of the objection that a theory which begins from beliefs must end wherever
the beliefs happen to be. It does not, because coherence is a strong condition, and because a
further assumption, exchangeability, converts it into something which looks objective.

An assessment on a sequence of events is exchangeable when it is invariant under finite
permutations of the indices: the person regards the order as irrelevant to what he expects. de
Finetti's representation theorem \cite{deFinetti1937} then states that an exchangeable assessment
on an infinite sequence is a mixture of assessments under which the events are independent and
equally likely. The classical limit theorems follow: the empirical frequencies converge, the
limiting frequency exists as a random variable, and the law of large numbers is a theorem about
the mixture rather than a postulate about the world. The gain is not merely that the objective
conclusions are recovered, but that the assumption on which each of them rests is made visible.
Subjectivism, so understood, does not weaken the theory; it locates the hypotheses.

% ---- from generated/ch1-body.tex ----
%% Generated by build.py from the parts - do not edit.
\section{The Argument of the Book}

This is a monograph developing the foundations of Probability Theory and of Mathematical Finance outside of the traditional measure theoretic setting based on sigma algebras 
and following instead the principles of constructive mathematics. 

We do so not out of a philosophical preference but as a matter of practical reach. Several developments
make the attempt timely. 

The first is that probability theory and stochastic processes find pervasive application in 
numerical methods, from Finance to Artificial Intelligence, Engineering and Biology. Nevertheless, the
measure theoretic foundations are highly non-constructive and ill-suited to the algorithmic and finite nature 
of numerical methods. In numerical analysis one typically works with finite approximations and algorithmic procedures
and one is interested in convergence properties to the continuum rather than abstract results in the continuum itself. 

Furthermore, in numerical analysis one is interested in the stability of algorithms under
floating point errors, especially when dealing with large-scale simulations or iterative algorithms which
benefit from vectorisation on GPU processors. In such contexts, the design of efficient algorithms greatly
benefits from the theoretical understanding of convergence rates for equations with coefficients which 
are not perfectly smooth. As we show in the later chapters of this book, our constructive approach provides explicit 
estimates and convergence rates also for non-smooth coefficients and the results are quite surprising when
compared to textbook estimates. 

Another consideration is that libraries of formalised mathematics now contain enough of
the elementary theory --- linear algebra, finite dimensional convexity, measure on countable
spaces --- for statements of the kind proved here to be formalised in full rather than
sketched. 

Finally, AI-assisted proof development has changed the economics of
formalisation, so that a long and systematic argument, of the sort that a treatise on the
foundations of pricing consists of, is no longer out of reach of a mechanical check. What
such assistance cannot repair is a foundation whose statements are not well posed as formal
objects, and this is the reason for the emphasis on constructivity: the aim is a formulation
in which every statement is a statement about finite objects, and in which the passage to the
continuum is an estimate rather than an existence theorem.

The book is planned in chapters. Chapters~1 to~3 establish the finite core:
the fundamental theorem at a fixed resolution, the finite dimensional convex
geometry on which it rests, and the valuation theory that follows from it. Chapter~4 turns from the
core to the objects to which it applies: it defines the path dependent functionals of the
process, and reduces the computation of their joint law with the underlying to a finite number
of lattice propagators whenever the functional is Abelian, and to a block factorised backward
recursion when it is not. It is there, and not in the analytic chapters, that the identities which
the classical theory obtains by a change of measure and a passage to the limit are proved:
Cameron, Feynman, Girsanov, It\^o, Kac and Martin, at fixed resolution, as exact identities between
finite matrices. Chapters~5 and~6 supply the estimates by which the passage to the continuum is
controlled. Chapter~5 establishes uniform rates for the convergence
of the discretized kernels of the diffusions out of which the models are built to the kernels of
the continuum equations, as the space and the time step vanish, for coefficients of limited
smoothness; Chapter~6 extends the analysis to the joint kernel of
a diffusion and an Abelian process, and thereby obtains a constructive definition of the
stochastic integral. Those two chapters are where the proofs cease to be finite dimensional and
where the iterative control of the limit becomes the substance of the argument. Chapter~7 carries
the estimates of Chapter~5 to every dimension $d\ge2$, where the lattice loses its
isotropy and the algebra of the renormalisation step changes, and where the rate is nevertheless
the same; it also records that the argument of Chapter~5 does not itself generalise, and identifies
the obstruction which obliges the higher dimensional proof to proceed by averaging over paths
rather than by a telescoping over them. Later chapters develop the applications, and
the numerical schemes that the matrix formulation makes possible.

The order of the chapters is the order of the objects, and it is deliberately not the order of
increasing difficulty: Chapters~1 to~4 are finite dimensional and require no analysis at all,
while Chapters~5 to~8 are the analytic ones, in which the passage to the continuum is controlled by
estimates. The first three chapters are the ones which can be verified mechanically today, as \S\ref{ch3:sec:lean}
reports; Chapter~4 is exact at fixed resolution in the same sense, and the section on machine
checked verification of Chapter~4 records what its formalisation would involve.

The relation of the analytic chapters to the classical theory is worth stating precisely, because
it is easily misread. The theorems of Feller, of Stroock and Varadhan and of the semigroup school
establish convergence: a limit exists, it is unique, and it is the fundamental solution of the
parabolic equation. They give no rate. The estimates proved in Chapters~5 to~8 are stronger, since a
rate uniform over a class of coefficients implies convergence, together with the identification
of the limit, in a norm finer than any of those theorems requires; they therefore imply the
classical results rather than rest on them. The one thing genuinely taken from the classical
theory is the existence of the continuum kernel itself, a weak solution of the forward equation
with a delta initial condition, which the estimates identify rather than construct, and which is
recorded as an open item in Chapter~6. The programme is to define that object instead as the limit
of the sequence of lattice kernels, in which case its existence follows from the rate by a Cauchy
argument and its equation is recovered by passing the lattice equation to the limit, so that the
book rests on nothing but finite matrices and explicit bounds. What that reorganisation
needs is precisely the estimates which Chapters~5 to~8 at present supply only conditionally.

\subsection*{The Sense in Which the Theory Is Constructive}

The word constructive is used in this book in the sense in which it is used in constructive
field theory, and not in the sense in which it is used in logic. It does not name a restriction
of the logic: no step below is non classical, the real numbers are the classical real numbers,
and every estimate is an ordinary inequality between ordinary numbers. It names instead the
manner in which the continuum is reached. A theory is constructive, in the sense intended here,
when the object on the continuum is not postulated and then studied, but obtained as the limit
of a sequence of discrete approximations, with the passage to the limit controlled by estimates
which can themselves be checked, and with the renormalisation which that passage requires
exhibited rather than assumed. This is the programme of constructive quantum field theory, in
which the lattice spacing is sent to zero and the limit is controlled by a renormalisation group
analysis; it is the programme which was carried out for lattice gauge theory in four dimensions,
and which has more recently been the object of machine verification. This book belongs
to that tradition, and takes from it its standard of rigour and its instrument, which is the
renormalisation group.

What is set aside here is not continuum mathematics, of which the book makes free use: Fourier
analysis, contour deformation and the classical theory of parabolic equations are instruments of
the argument, and a theory which dispensed with them would be the poorer for it. What is set
aside is narrower and more specific, namely the sigma-algebra on the continuum as the primitive
object of a probability theory. Kolmogorov's axioms take such an algebra as given, together with
a countably additive measure on it, and define a process as a family of measurable functions on
that space, so that every subsequent statement is a statement about equivalence classes of such
functions. No object of that kind is introduced in this book, and no statement is made which
requires one. Probability enters through the prevision of a finite family of path functionals at
a fixed resolution, which is the primitive of de Finetti's theory, and the role which the
sigma-algebra plays in the classical foundation --- that of the object on which the passage to
the limit is organised --- is taken here by the sequence of lattice refinements together with the
estimates which control it.

\paragraph{Absoluteness.} A further reason for the finite formulation concerns not its power but the
kind of statement it permits one to make. The statements of the finite theory are statements about
finite matrices, hence arithmetic statements, and arithmetic statements are absolute: they hold or
fail in every model of set theory alike, whatever is assumed about choice and about the
measurability of arbitrary sets of reals. The statements of the measure theoretic theory are not
all of that kind. In ZFC there exist sets of reals which are not Lebesgue measurable, and it is
equally consistent that every set of reals is measurable; the two cannot be had together, and
which holds is settled by the axioms one adopts rather than by the objects the theory speaks
about. This is not a defect of the classical theory, which is entitled to make its choice and does
so consistently; it is only to observe that a statement whose content depends on that choice is
less portable than one whose content does not. The difference is not confined to pathology. A
conditional expectation is an equivalence class, so a price is determined only up to a null set,
and the null sets are precisely where the set theoretic freedom lives. In the present theory a
price at a fixed resolution is an entry of a matrix, a rational number once the data are rational,
and the approximation to it is a number with an error bar; nothing is asserted almost surely, and
nothing is left undetermined by the statement itself. The sense in which the proofs of this book
are verifiable is therefore precise: what is verified is a statement whose truth does not depend
on the metatheory, even though the system in which the verification is carried out, as every
system must, assumes something about choice.

\paragraph{Decidability.} Absoluteness is not yet the whole of the clarity which applications demand.
At a fixed resolution the statements of the theory are not only absolute but decidable, and
decidable by computing: the entries are numbers which exact arithmetic derives from the data, and
every question the theory asks of them --- whether an arbitrage exists, whether a given price is
attained, whether a kernel is strictly positive --- is a question about finitely many such
numbers, answered by a finite computation. Where the data are rational the computation is exact;
where they are not, it is carried out in interval arithmetic, and the modulus of which we spoke
bounds the width of the interval. The answers come moreover with witnesses: an arbitrage is an
explicit combination of the traded claims and a pricing kernel is an explicit vector, so that a
computation can be verified by recomputing its certificate rather than by trusting the solver
which produced it.

The statements about the limit are of another kind, and it is essential for applications that the
two be distinguished. A statement about the rate of convergence quantifies over the resolutions,
and is therefore proved rather than computed; what is computed is its instances. This is the
division which the theory offers to a practitioner: the infinitary statement is the theorem, and
the finite instance, with its error bar, is the computation. No question has to be settled in the
limit in order to be answered at a resolution, and it is answers at a resolution that applications
require.

Where the analogy with constructive field theory ceases to be an analogy is at the singularity.
The passage to the continuum is difficult in this subject for the same reason as in that one:
the kernel is singular in the short distance limit, and as the lattice spacing vanishes it
develops a diagonal singularity which is the delta initial condition of the classical
formulation. The control of that singularity is what occupies Chapters~5 to~8, and it is why those chapters are analytic while Chapters~1 to~4 are not. In one dimension the control is obtained by a
telescoping over scales. In higher dimension the telescoping fails, the lattice is no longer
isotropic, and the renormalisation step must be replaced by an averaging over paths. The
obstruction met in Chapter~7, and resolved there, is therefore the standard obstruction of the
renormalisation group rather than an artefact of the method, and the survival of the rate at the
same order in every dimension is the statement that the renormalisation is the correct one.

In a framework of this kind the rate is not a refinement of the existence statement: it is what
the existence statement means. The object on the continuum is not constructed and then estimated;
the estimates are all there is, and the limit is defined by them, as the value towards which the
sequence of lattice computations converges with that modulus. This is why the analytic chapters of the book are organised around rates rather than around limiting objects, and why the strength
of a theorem here is measured by the uniformity of its rate over the class of coefficients, of
observables and of auxiliary parameters to which it applies. An existence theorem is a
yes or a no; a rate is a certificate, and it is the certificate which a computation can carry.
A sharp estimate is more still, because it locates the error rather than confining it: the leading
term of a sharp expansion is a quantity which can be corrected for rather than merely feared, and
the constant which multiplies the rate is what converts a statement about the limit of vanishing
meshes into a statement about the mesh one can afford.

\subsection*{What the Theory Buys in Practice}

The reasons for insisting on this foundation are practical, and they are worth listing, because
each of them is a place where the measure theoretic formulation leaves the practitioner without
an answer.

\begin{itemize}

\item \emph{The error of a computation can be estimated before it is run.} A computation at
resolution $h$ has two sources of error, the discretisation error and the rounding error. The
estimates of this book bound the first explicitly and uniformly, and the standard analysis of
floating point arithmetic bounds the second; the two can therefore be compared. The precision
required to reach a given accuracy, and the accuracy attainable at a given cost, are both
computable quantities rather than matters of experience.

\item \emph{Smoothness and precision are two readings of one parameter.} The estimates are
uniform over classes of coefficients, and rounding replaces the coefficients by a perturbation of
them whose size is the unit roundoff and which is uncorrelated from one lattice site to the next;
that is to say, the perturbation which the arithmetic introduces is rough at the scale of the
lattice. A coefficient known only to that accuracy is therefore, as far as the estimate is
concerned, a coefficient of limited smoothness: the noise and the roughness enter through the
same channel. The rate as a function of the H\"older exponent consequently doubles as the rate as
a function of the number of bits, and the curve which states how the scheme converges for a rough
coefficient states also how it degrades in low precision arithmetic. This is why the dependence of
the rate on smoothness, and not convergence alone, is the quantity which matters: Chapter~5 finds a
rate which falls as the exponent falls, and Chapter~7 finds that the obstruction costs the order of
the rate and not convergence, first order being available in every dimension by bounding the
obstruction and second order only through a cancellation.

\item \emph{The constant, and not only the order, is what informs.} An order of convergence says
how the error behaves as the mesh vanishes; the constant says what the error is at the mesh one
can afford. A sharp estimate therefore turns a rate into a prediction, and, because the leading
error is then known rather than merely bounded, into a correction: the leading term can be removed
by extrapolation, and a computation can be given an error bar tight enough to certify it rather
than only to bound it. This is why the estimates of Chapters~5 to~8 are stated with the constants
displayed rather than absorbed; why Chapter~7 records that the constant depends on the direction
while the rate does not; and why Chapter~8 determines the dependence of the constant on the
conjugating parameter explicitly, and finds in the process that the invariant object is the two
form attached to the coefficient field while the numbers it carries in a given frame are not.

\item \emph{The objects of the theory are the objects of the machine.} At resolution $m$ every
object is a finite matrix or a finite sum, and every operation is a matrix product. Nothing is
required in an implementation which is not already present in the statement of the theorem, so
that the natural form of an implementation is the vectorised one, and the primitive on which it
rests is the matrix multiplication, which is where the parallelism of current hardware lies.

\item \emph{Sensitivities are exact.} The derivative of a computation built from matrix products
is the transpose of that computation, and the backward recursion which propagates sensitivities
through the levels of the lattice is the adjoint of the propagator. At a fixed resolution this is
an identity between finite matrices and not an approximation: the forward and the backward pass
of a computation, and with them the primal and the dual of a pricing problem, are two readings
of the same product.

\item \emph{Conditioning is an operation and not a choice.} A conditional expectation at a finite
resolution is a definite matrix, whereas in the measure theoretic formulation it is a version of
an equivalence class, defined only almost surely --- a distinction immaterial to the theory and
consequential in an implementation. Dynamic conditioning, which is the substance of
Chapters~1 to~4, is therefore a map which can be evaluated rather than one which has to be
represented.

\item \emph{The absence of arbitrage is a finite computation.} The Fundamental Theorem
at a fixed resolution is the orthant alternative of finite dimensional convex geometry, so that
the existence of a pricing kernel is the feasibility of a linear program, and its failure is
certified by an explicit combination of the traded claims. \S\ref{ch2:sec:convex} gives the
algorithm. The classical statement is a separating hyperplane theorem in a space of random
variables, which settles the question without yielding a procedure.

\item \emph{Limited smoothness is the normal case.} The estimates of Chapter~5 are uniform over
classes of coefficients of bounded H\"older differentiability, and the schemes continue to
converge, at a slower rate, when the coefficients are merely uniformly continuous. The classical
regularity hypotheses, under which the convergence of discretisations of a parabolic equation is
usually proved, are therefore not needed here, and the rates which are available in their absence
are the ones which the applications in fact present.

\item \emph{The theory can be checked by machine.} A statement about a finite matrix admits a
formal proof, and each chapter reports what has been certified and what has not. A
statement about a measure on a continuum, being a statement about equivalence classes, does not
admit the same treatment until representatives have been chosen. The ledgers of \S\ref{ch3:sec:lean},
and those in the corresponding sections of Chapters~4 and~6, are the record of what this yields.

\end{itemize}

The last two items, together with the identification of the primitives, account for the reach of
the programme beyond finance. A learning algorithm for a diffusion is built from the same three
ingredients as a pricing algorithm: an approximation at a finite resolution, a product of
matrices across the levels of that resolution, and the adjoint of that product for the backward
pass. A generative model defined by a stochastic differential equation asks, as a path dependent
derivative does, for the joint law of the diffusion and a functional of its path, which is the
object of Chapters~6 and~8. What transfers is the structure of the objects and the calculus of
their errors; the dimensions in which such models are used are far beyond those for which the
estimates here are available, and no claim is made that the argument closes that gap. The claim
is that one foundation serves both, and that in both the question which matters is the one the
classical foundation cannot answer: how far from the target a computation is, and at what cost in
arithmetic that distance can be reduced.

\subsection*{The Programme}

The guiding idea is that a stochastic process is defined by prescribing the
expectation of every functional of the path, that is, by prescribing a family
of $n$-point functions which satisfy a short list of axioms. Two features
distinguish our treatment from the classical one.

First, the space variable is discretized. Space is approximated by a
sequence of lattices $h_m\ZZ^d$ with $h_m=2^{-m}$, and a process is a
consistent family of Markov chains on these lattices, one for each
resolution index $m$. Statements about continuous processes are statements
about limits $m\to\infty$ of statements about finite matrices, and the
passage to the limit is controlled by the convergence estimates of Chapter~5 rather than by abstract existence theorems.

Second, everything at fixed resolution $m$ is finite dimensional. The state
space $\StateA$ is a finite set, the horizon $[T,T']$ is finite, and the
information available at each time is a finite partition. This is not a
technical simplification but the substance of the philosophy: it is what
makes the theory constructive, computable, and --- as this chapter will show --- logically transparent. In particular, the existence of a pricing
measure, which in the classical theory requires a separating hyperplane
theorem in a space of random variables and a measurable selection argument,
becomes here a statement of finite dimensional convex geometry: a strictly
positive vector orthogonal to a finite dimensional subspace. That statement
is proved in \S\ref{ch2:sec:convex}, together with its converse.

\subsection*{Foundations: Kolmogorov and de Finetti}

\paragraph{The measure theoretic foundation.} Kolmogorov's axioms
\cite{Kolmogorov1933} take as primitive a probability space
$(\Omega,\mathcal{F},\Pmeas)$, with $\mathcal{F}$ a $\sigma$-algebra of subsets
of $\Omega$ and $\Pmeas$ countably additive on it, and define a process as a
family of measurable functions on that space, its law being obtained by the
extension theorem from a consistent family of finite dimensional
distributions. The construction is not effective, and not by accident. The
$\sigma$-algebra is a subset of the power set of a continuum; the extension
theorem produces a measure on an uncountable product by a compactness
argument; and the objects which the theory then manipulates --- conditional
expectations, densities, versions of martingales --- are determined only up
to null sets, so that every statement about them is a statement about an
equivalence class whose representative has to be chosen. None of these steps
can be executed, and none can be formalised without committing, at each
point, to a choice which the statement itself does not determine.

\paragraph{The finitely additive foundation.} de Finetti
\cite{deFinetti1931,deFinetti1937,deFinetti1974} takes as primitive not a
measure but the price of a bet. A probability assessment on a family of
events is \emph{coherent}, or admissible, if no combination of bets against
it produces a sure gain, and de Finetti's fundamental theorem of probability
states that a coherent assessment extends to a finitely additive measure on
the algebra generated by the events \cite{BhaskaraRao1983}. This is a
theorem of linear algebra. In the finite case it is exactly the orthant
alternative of \S\ref{ch2:sec:convex}, where the extension is produced by the
algorithm of \S\ref{ch2:ss:farkas} and, when the assessment is incoherent, the
same algorithm produces the combination of bets which refutes it.

The reader will recognise the statement under another name. A combination of
bets with a sure gain is an arbitrage, the finitely additive measure is a
pricing measure, and de Finetti's fundamental theorem of probability and the fundamental theorem of asset pricing are one theorem, read once as a statement
about probability and once as a statement about prices. This identification
is the reason why a constructive theory of stochastic processes is
conceivable at all.

\paragraph{The constructive route.} It requires two steps, and they
correspond to the two halves of the argument. The first is to define a
process without passing through a measure on a path space: a process is a
consistent family of $n$-point functions on a lattice, so that at each
resolution $m$ every object of the theory is a finite matrix or a finite sum.
This is done in \S\S\ref{ch3:sec:simplicial}--\ref{ch3:sec:npoint}, and at the level
of generators and propagators in \S\ref{ch3:sec:markov}. The second is to control
the passage to the continuum iteratively, resolution by resolution, by
showing that the pricing measures of the models at resolution $m$ converge
and that the martingale property survives the limit. This is the content of
Chapters~6 and~5, and
\S\ref{ch3:sec:remarks} says precisely what remains to be proved before it
applies to the theorem established here.

\paragraph{Two consequences for mechanisation.} The route lends itself to
machine checking, because at every resolution the objects are finite and the
convex geometry on which the theorem rests is elementary; \S\ref{ch3:sec:lean}
records how much of it has been verified in Lean~4. It also belongs to a
class of problems which have only recently become tractable: a long programme
of estimates which controls a continuum limit, resting on a finite
combinatorial core that can be verified mechanically. The constructive
approach to the Yang--Mills existence problem along Ba\l{}aban's
renormalization group \cite{Balaban1984,Balaban1987} is a problem of this
kind, and an AI-assisted attack on it has recently been made
public.\footnote{The project \texttt{MathePhysics/4DWilsonLoopRG}, reported
in October 2026, extends Ba\l{}aban's method to renormalised Wilson loop
expectations. It is work in progress rather than a publication, and we refer
to it only as an indication of the kind of problem for which a finite core
combined with a controlled limit is now within reach.} The present book is
on a smaller scale, but the shape of the argument is the same, and the fact
that its finite dimensional core is already machine-checked is evidence that
the shape is not an accident.
% ---- end generated/ch1-body.tex ----

\subsection{Provenance, and the Plan of the Book}
\label{sec:concordance}

The material of this book was first worked out in a series of papers, and the chapters which
reproduce it were written as those papers; the prose has been brought into the vocabulary of the
book, so that each of them speaks of the chapter it is, and refers to the others by chapter rather
than by the number of a paper in a series.

Chapters~\ref{chap:ch9} to~\ref{chap:ch11} have no counterpart in that series. They were written
for the book, and they are the ones which make its two halves one argument rather than two: the
adjoint of a propagator product, which is the derivative the inverse problem of
\S\ref{sec:two-problems} consumes, and the error of that derivative, which is the estimate the
inverse problem needs and which none of the papers supplies. The second of them is also where the
book leaves the lattice: \S\ref{ch10:sec:meanfield} takes the mean field spin glass --- a model with
no stencil, no block of neighbours and no length scale --- as far as an exact comparison and a
measured remainder allow, and no further, since a programme is not a theorem.

\section{Two Problems, and the One Theorem Which Relates Them}
\label{sec:two-problems}

Every model of this kind is used in two ways, and the difference between the ways is the organising
principle of this book.

The first is the \emph{inverse} problem. The model is a family of coherent assessments --- in
finance a family of pricing measures, in a learning problem a family of predictive distributions ---
and one chooses within it the member which agrees with what is observed: the market prices of a set
of liquid claims, or the values of a set of targets on a set of examples. This is calibration. It is
an optimisation rather than an equation to be solved, because there are typically far more targets
than free parameters and the targets are themselves approximate; it is iterative; and it is where
the model is fitted.

The second is the \emph{direct} problem. The member having been chosen, one evaluates it: one prices
a claim whose market price is not known, or one computes a functional of the predictive distribution
which the fitting never mentioned --- a quantile where the fit supplied a mean, an expected
shortfall, the expectation of a path functional under a generative model. The direct problem is one
application of a linear operator. It is cheap, it does not iterate, and in this book it is a single
backward multiplication of calibrated propagators, which is Chapter~\ref{chap:ch3}.

The two are related by Theorem~\ref{thm:coherence}, and the relation is not a matter of taste. On
the linear span of the claims against which the model was calibrated, the answer to the direct
problem is determined, by coherence alone, and it is the market's answer. Off that span it is not
determined: the calibrated set contains members which agree on every calibrating claim and disagree
on the new one, so that the answer is a range rather than a number, and a model is what selects a
point within it. When the span is everything the problem is called complete and the selection does
not matter, there being nothing to select; when it is not, the selection is the whole content and it
should be made visible as such. The width of the range is what the finance side calls model risk,
and it is measured rather than assumed.

Four consequences organise what follows. The inverse problem needs the derivative of the direct
map, because an optimisation step asks how the computed prices move when the coefficients move, and
that derivative is the adjoint of the propagation, which is Chapter~\ref{chap:ch9}. The direct
problem needs only the propagation itself, which is why it is treated first, in
Chapter~\ref{chap:ch3}. The errors of the two are different objects: the error of a price is a
convergence statement about a kernel, which is Chapters~\ref{chap:ch5} to~\ref{chap:ch8}, while the
error of its derivative in the coefficients is Chapter~\ref{chap:ch10}. A book which estimated only
the first would be silent on the question of whether a calibrated model had been fitted to anything.
And the answer to the inverse problem is not a number on a continuum but a vector at each
resolution, so the fourth consequence is the one Chapter~\ref{chap:ch11} draws: what the estimates
give for the objects the inverse problem is actually about --- the law of the functional a claim is
written on, the value of a discounted expectation, and the prices of the limit model itself.

None of this is a statement about finance in particular. Read the calibrating claims as the
statistics which a learning problem fits, and the new claim as any other functional of the
predictive distribution, and the two problems and the theorem relating them are the same two
problems and the same theorem. Calibration is training; the direct problem is inference, whether it
is called a prediction, a risk measure or an expectation under a generative model; and the
consistency which the finance side calls repricing --- the calibrated model must reproduce the
instruments it was calibrated to --- is on the learning side the requirement that the fitted model
reproduce the statistics it was fitted to, a requirement which fails in the same ways, through
approximate optimisation, through sampling error, and through functionals which were calibrated only
implicitly. In this reading the fundamental theorem is what tells a reader of either field which of
the questions he is asking has an answer determined by his data and which has an answer determined
by his model.

What distinguishes the two applications is what a target is. A price is the known value of a linear
functional; a sample is one draw from the distribution. The finance side therefore constrains the
model by equations, and its undetermination is the incompleteness of a market; the learning side
constrains it statistically, and its undetermination is what is now called underspecification, many
members of the family agreeing with the training statistics and disagreeing elsewhere. The theorem
is the same in both; the constraints applied to it are not, and they are weaker on the learning
side. This asymmetry is the honest form of the claim that one foundation serves both fields, and it
is why nothing below asserts that a method proved for one transfers unchanged to the other.

The claims made in what follows are accordingly structural, and the reader is invited to hold them
to that standard. The estimates of Chapters~\ref{chap:ch5} to~\ref{chap:ch11} are available in few
dimensions, for coefficients of limited smoothness. The models of the applications are built in
dimensions far beyond those for which the estimates are available, and no claim is made that the
argument closes that gap. What transfers is the pair of problems, the theorem which relates them,
and the calculus of the errors by which each of them is measured; not the constants.

One consequence for how the book is written should be stated here, since it is a promise the reader
may hold the chapters to. Each chapter which follows names the object it treats in both readings
before it develops it, so that a reader of either discipline can follow the mathematics and translate
the nouns. The mathematics is not duplicated anywhere: there is one set of matrices, one set of
estimates, and two vocabularies for saying what they are about.

It is worth stating, before the reader meets the details, what the argument is meant to deliver,
because the route is unusual and the last chapter is where it arrives. The pricing measure is never
computed as an object on a continuum. At each resolution it is a finite vector: the local condition
of \S\ref{ch3:sec:ftf} is a linear system in finitely many unknowns, Farkas' lemma decides whether
that system has a strictly positive solution and produces one when it does, and where the market is
complete the solution is unique up to the normalisation of the numeraire
(\S\ref{ch3:ss:second}). The computation is constructive in the sense this book gives the
word: it is an algorithm, it terminates, and it carries a certificate which a second computation can
check.

What the algorithm does not settle is the relation between the answers at successive resolutions,
and that is what the estimates of Chapters~\ref{chap:ch5} to~\ref{chap:ch8} and the theorem of
Chapter~\ref{chap:ch11} supply: the prices converge, at a rate, provided the objects the computation
produces are regular --- bounded above and below, with a bounded logarithmic derivative, uniformly
in the resolution. The condition is on the calculated objects and not on the market directly, which
is the right way round for a constructive theory: what one controls is what one computes. The shape
of the whole argument is then this. A measure change is calculated rather than postulated, and the
passage to the continuum is earned by the regularity of what has been calculated, rather than
assumed of a measure which no finite computation ever produces.
% ---- end chapters/ch1.tex ----
% ---- from chapters/ch2.tex ----
% Authored: the chapter title and opening belong to the book. The bodies are
% generated by build.py from Part 1; the one-period theorem is lifted out of
% Part 1's Fundamental Theorem section so that its statement opens this
% chapter and its proof closes it.
\chapter{The Fundamental Theorem}
\label{chap:ch2}

This chapter proves the fundamental theorem, and it proves it in the static setting, where the
argument is completely explicit: a finite set of states, finitely many bets, and no time. Its
subject is a single bet and the constraint that a family of bets admit no combination yielding a
sure gain; its conclusion is that this constraint is equivalent to the existence of a strictly
positive state price density, which is to say of a finitely additive prevision under which every bet
is fair. The vector whose existence is asserted is both at once, the density in the pricing reading
and the prevision in the learning one, and the chapter keeps the two vocabularies side by side
throughout: the theorem is one statement, and the disciplines differ in what they call its terms and
not in what they require of them. The equivalence is established here in full, together with an
algorithm which decides which side of it holds and delivers the witness for that side.

The chapter has three sections, and they are arranged so that the statement comes first and the
proof last. The first states the fundamental theorem in the two vocabularies of
\S\ref{sec:coherence}, where it is a statement about prices and about belief at once, and records
the several names by which it is known in the literatures which use it. The second proves the finite dimensional alternative on which it rests --- Farkas' lemma,
the orthant alternative and Stiemke's lemma --- which is the engine of the whole book. The third
returns to the theorem and proves it, at which point the proof is two lines. The order is
deliberate: a reader who wishes to know what the theory claims may stop after the first section,
and a reader who wishes to know why it is true needs only the second and the third.

\section{The Shared Foundation of Two Applied Disciplines}
\label{sec:coherence}

The subjectivist account of probability begins not from a space of outcomes but from a bet, and
from the price at which a person is willing to take either side of it. If a person will buy or sell
a gamble $X$ at the price $p$ indifferently as to the side, then $p$ is his prevision of $X$; a
family of gambles assessed in this way is a family of prices, and an assessment on a family of
events is the special case in which the gambles are indicator functions. Nothing else is primitive.
In particular no algebra of events is given in advance and no measure on it is postulated; the
assessment is a finite list of numbers attached to a finite list of bets.

The single requirement imposed on an assessment is that it should not be possible to lose for
certain. An assessment is \emph{coherent} when no combination of the bets, in the quantities to
which the person is committed, yields a gain which is non negative in every state and not
identically zero. That is the whole of the requirement, and the classical properties of
probability are derived from it rather than assumed: additivity follows by considering a bet on a
union against bets on its parts, the multiplication rule follows by considering a compound bet
whose settlement depends on two events, and conditional prevision follows by restricting attention
to the bets which are settled only in the states where the conditioning event holds.

In the one-period case, which is the case of this chapter, the content of the requirement is
completely explicit.

% ---- from generated/ch2-theorem-statement.tex ----
%% Generated by build.py from the parts - do not edit.
\begin{theorem}[The fundamental theorem]\label{ch2:thm:oneperiod}
Let $S$ be a finite set and let $\xi^1,\dots,\xi^n\in\RR^S$. The following are
equivalent.
\begin{enumerate}
\item[(i)] There is no $\zeta\in\RR^n$ with $\sum_k\zeta^k\xi^k(s)\ge0$ for all
$s\in S$ and $\sum_k\zeta^k\xi^k\neq0$.
\item[(ii)] There exists $g:S\to(0,\infty)$ such that
$\sum_{s\in S}g(s)\xi^k(s)=0$ for $k=1,\dots,n$.
\end{enumerate}
\end{theorem}
% ---- end generated/ch2-theorem-statement.tex ----

The theorem speaks only of a finite set $S$ and of finitely many vectors in $\RR^S$; it mentions
neither a market, nor a filtration, nor a time. That is what allows it to be read in two ways at
once. Read $\xi^k$ as the gain of a claim bought at its price, and (i) is the absence of a sure
gain from a position in the traded claims, while $g$ is a state price density. Read $\xi^k$ as a
gamble offered at a price, and (i) is the coherence of an assessment, while $g$ is the prevision
which extends it. The theorem says that these are one condition, and the same $g$ serves both
readings. The proof is deferred to \S\ref{sec:proof-oneperiod}: it will be seen there to be
Stiemke's lemma in another notation, so that the economics of the one-period case and the geometry
of the orthant are the same statement.

The result carries more than one name, and the names record the several subjects in which it is the
fundamental result. It is worth setting them out at length, because the book will use both
vocabularies and because the terms of one are not the terms of another, although the mathematics is
identical: what changes between the readings is only the semantics attached to $\xi$ and to $g$.

de Finetti, who proved the result in the finite case, calls it the \emph{fundamental theorem of
probability theory}, and the name is exact. It is the theorem which shows that probability is a
derived and not a primitive notion: the assessment is what is given, the vector $g$ is what a
coherent assessment determines, and what a conventional definition would have postulated as a
probability is here the conclusion of an argument in linear algebra. This is the precise sense in
which probability does not exist --- not the calculus, but the object which the calculus is usually
taken to be about.

The two readings set out below are the two vocabularies in which the statement is used, and it is because the statement is one that the book can move between the subjects without
changing its mathematics. It also explains why the finite dimensional form of the result is the form
that matters: in each of the two readings the objects are finite lists of numbers, the criterion is
a finite computation, and the certificate is a finite object. In what follows the result is referred
to simply as the fundamental theorem, and which of the two readings is in force is a matter of what
is being modelled. The subjectivist account of belief, which is the classical vocabulary of the
learning reading, is treated inside it rather than beside it.

The vocabulary of the sections which follow, and of the greater part of the book, is financial, and
the choice is one of convenience rather than of substance. Nothing in the mathematics privileges
claims over gambles or over data: every statement could be written in the language of previsions, or
in that of parameters and posteriors, and it would say the same thing with the same proof. Finance is
adopted because it carries the most developed terminology for the objects in question --- claims,
prices, replicating portfolios, hedging, completeness, arbitrage --- and because the multi-period
structure developed in the later chapters is the one in which those objects are most sharply defined.
A reader who prefers either of the other vocabularies should read claims as gambles and prices as
previsions; nothing else need change. For this reason the result is called simply the fundamental
theorem, and not, as it is variously called in the several literatures, the fundamental theorem of
finance, of asset pricing, of probability theory, or of Bayesian statistics. Those names are
several vocabularies for one statement, and they fall into the two readings of
\S\ref{sec:coherence}; to choose among them would be to choose a subject rather than a theorem.

The separation of the two ways in which a model is used, drawn in \S\ref{sec:two-problems}, is
already visible in this one period statement and is worth marking here, where the theorem is proved,
rather than only where it is applied. The theorem asserts two things at once. The first is that the
set of coherent assessments extending the prices of a family of traded claims is not empty. The
second is that every member of that set assigns the same value to every claim in the linear span of
the traded ones. On the span, therefore, the answer is the market's and not the model's; off the
span the theorem is silent, and assessments which agree on the traded claims disagree elsewhere.
Read so, the one period theorem is the first appearance of the distinction between what a calibration
determines and what it leaves to the choice of model, and the second fundamental theorem --- which
\S\ref{ch3:ss:second} states for a process --- is the statement that the indeterminacy disappears
exactly when the span is everything.

The two readings are not two subjects but the two directions of time. In both, an assessment is
formed now from what is known, and in both the test of it is the same; what differs is where the
data lie in relation to the quantity assessed. In the first, the pricing reading, the quantity is a
payoff which has not yet settled, and there are no realisations of it to learn from: what is
observed is the price at which the claim trades, and the test of the assessment is whether the market
will take the other side of it. In the second, the learning reading, the data are themselves
realisations --- observations drawn from the past of the quantity and paired with what the quantity
turned out to be --- so that the learning looks backward over a sample whose joint law is taken to be
the law of the quantity still to come. Time
has only two directions, so these two cases are exhaustive and there is no third discipline: an
assessment is either of a quantity whose past realisations are in hand or of one which has not
settled and of which there is nothing yet to learn. The vocabulary of belief, of prior and
posterior, belongs to the backward direction, and a reader who expects a third subject under that
name will find only this.

\subsection{Finance}

The vocabulary of finance is that of a market. A \emph{claim} is a random quantity, that is, a payoff
depending on the state of the world; a \emph{price} is what it costs today; and the \emph{gain} of
holding the claim is the payoff less the price paid, which is the vector $\xi$ of the theorem. A
\emph{position} is a combination $\sum_k\zeta^k\xi^k$ --- a portfolio of the traded claims, financed
at their prices --- and an \emph{arbitrage} is a position whose gain is non negative in every state
and not identically zero. The first statement of the theorem is then the absence of arbitrage, and
the second is the existence of a strictly positive \emph{state price density} $g$, a weight on the
states under which every traded claim has zero gain, which is to say under which every claim is
priced as its expectation. That the two are one statement is the fundamental theorem of asset
pricing in its one-period form.

The terms of the second statement deserve to be separated, because the rest of the book uses them
constantly. That $g$ is positive on every state is the assertion that no state is regarded as
impossible; when $g$ is normalised to sum to one it is a probability, the \emph{risk neutral}
probability of the one-period model, and the theorem says that it exists exactly when there is no
arbitrage. A claim which is the gain of some position is \emph{replicable}, and its price is then
forced, being the cost of the replicating position; the position itself is the \emph{hedge}, and its
existence is what makes the price a number rather than an interval. When every claim is replicable
the model is \emph{complete}, and then $g$ is unique up to normalisation. When it is not, the set of
admissible $g$ is convex and the prices it produces fill an interval: every value between the lowest
and the highest price compatible with the absence of arbitrage is attained by some admissible $g$,
and nothing outside the interval is.

The whole of the one-period theory of pricing is a corollary of the theorem in this reading, and the
local step of the multi-period theory of Chapter~\ref{chap:ch3} is the same theorem applied at a
single site, with the density $g$ replaced by the intensities of the risk neutral jump process.

\subsection{Machine Learning and Artificial Intelligence}
\label{ch2:sec:ml}

The second reading is the one in which the same theorem is at present the more consequential. A
prediction is a prevision: to fit a model is to form an assessment of a random quantity from data,
and to be coherent is for that assessment to be the conditional expectation of the quantity given
what is known --- the best available in mean square, and the requirement that no position in the
residual be profitable. The vocabulary of this reading is the subjectivist one, and it is developed
first, because the rest of the reading is stated in it.

\subsubsection{Bayesian Statistics}
\label{ch2:sec:bayes}

The vocabulary of Bayesian statistics is that of a person who holds beliefs about a finite set $S$
of possibilities. A \emph{parameter} is a point of $S$; in the subjectivist reading there is no
further sense in which the parameter is a fixed object about which the person is ignorant, since the
content of the theory is what the person believes about it and nothing else. A \emph{random
quantity} is a real function on $S$, which is exactly a gamble in the sense of the preceding
paragraphs, and its \emph{prevision} is what the person's belief amounts to numerically: the price at
which he would take either side of a claim paying that function. The prevision of the indicator of
an event is the probability he attaches to it, so that probability is recovered from prevision
rather than the other way round.

In this vocabulary the theorem says what follows. A finite family of previsions on a finite set of
random quantities is coherent if and only if it is the restriction of the prevision taken under some
probability on $S$, that is, of a \emph{prior}. Coherence is therefore not an additional axiom
imposed on top of the Bayesian scheme; it is precisely the condition that the scheme is available at
all, and the vector $g$ of the theorem is the prior.

Two constructions then carry the practical content of the subject, and both are consequences of the
theorem rather than additions to it. \emph{Conditioning} is how belief changes with information:
given an event known to have occurred, the prior is replaced by the conditional probability and the
new previsions are the conditional ones. Bayes' rule is not an independent assumption in this
account; it is the statement that a person who updates otherwise can be made to lose for certain,
which follows from the multiplication rule that coherence forces. An observation which is itself
uncertain is treated in the same way, by conditioning on the whole of what has been observed rather
than on a single event, and a \emph{posterior} is accordingly a conditional prevision. The previsions
of future observations computed from the posterior are the \emph{predictive} laws, and they are again
previsions, subject to the same theorem and to the same requirement.

What makes the account subjectivist rather than merely probabilistic is that the prior is not
derived from anything. It is the person's, the theory does not say where it comes from, and what the
theory says is that once it is held everything else is determined and that the whole of it must be
coherent at once. Symmetry narrows the choice, and de Finetti's representation theorem shows how far
that goes: an exchangeable assessment on an infinite sequence is a mixture of independent and equally
likely ones, so in that case the prior reduces to a prior over the mixing parameter and nothing else.

A third consequence of the reading is the one which carries furthest, and it concerns the \emph{use}
of an assessment rather than its formation. When a prevision is used to take a decision, or to be
scored, the coherence of the assessment becomes a condition on the score. A loss is called
\emph{proper} when the expected loss is minimised by reporting the belief which is actually held,
and \emph{strictly proper} when no other report attains the minimum; the logarithmic and the
quadratic losses are of this kind, and they are the losses of prediction. A method which is trained,
or evaluated, by a proper score is therefore answerable to coherence in exactly the sense of the
theorem, since what it reports is a prevision of the score and what the theorem supplies is the set
of probabilities to which that prevision may coherently be extended. Nothing in this requires the
method to be Bayesian in the sense of computing a posterior: what is required is that its output be
an assessment and that the requirement on the assessment be coherence. In this reading the
foundation which the theorem supplies is shared by the pricing of claims and by the forming and the use of
previsions; the two differ in what is assessed and in what is observed, and not in what is required
of the assessment.

A caution on vocabulary belongs here, because the word \emph{calibration} has two senses which this
book must keep apart. In this book it is the inverse problem discussed above and in
\S\ref{sec:two-problems}: the choice, within a family of models, of the member which agrees with
observed prices. In the statistical and machine learning literature a forecast is called calibrated
when events to which it assigns probability $p$ occur with frequency $p$. The two are related, since
a coherent assessor is calibrated in the second sense under exchangeability, but they are not the
same statement, and the word should be read according to its context.

The terms carry over without strain. A \emph{model} is a family of probabilities on the observed data
indexed by a parameter; the prior is a probability on that parameter; the posterior is its
conditional probability given the data, obtained by conditioning, that is, by Bayes' rule; and the
predictive law is the conditional law of a future observation. The provisional quantities of learning
theory are these quantities under other names. A regularisation term is what a prior encodes, the
output of a trained model is what a posterior, or a functional of it, approximates, and a calibrated
model is one whose predictive law is the predictive law the theorem speaks of. Bayesian neural
networks, Gaussian processes, variational inference and the wider family of probabilistic models are
instances of the same schema, and the schema is the theorem above applied to a family of random
quantities which happens to be large. The infinite dimensional members of the family --- Gaussian
processes among them --- are exactly where the finite dimensional theorem has to be supplemented by
the convergence estimates of the later chapters, which is the reason those chapters exist.

What the theorem adds to all of this is a standard. Coherence is the precise sense in which a system
of beliefs, and the decisions taken from it, can fail to be simultaneously satisfiable: an incoherent
system is one against which a sure loss can be constructed, exactly as an incoherent assessment is
one against which a combination of bets can be constructed. Bayesian decision theory, expected
utility, and the treatment of exploration as the acquisition of information are consequences of this
standard rather than alternatives to it. In the present setting the standard is moreover computable,
because the criterion is the feasibility of a finite dimensional linear program, decided by the
algorithm of \S\ref{ch2:sec:convex}, and its failure produces an explicit certificate. A standard
which can only be stated is a philosophical position; a standard which can be checked is a property
of a system, and it is the second that a machine can be held to.

One distinction should be drawn explicitly, because the reading invites a stronger claim than the
theorem supports. What a learning method takes from this foundation is the \emph{standard} ---
coherence, and its decision theoretic form as the class of proper scoring rules recorded in
\S\ref{ch2:sec:ml} --- and not the \emph{algorithm}, the computation of a posterior. A method
which computes no posterior at all, and which is trained by minimising an empirical risk, is
answerable to the standard all the same, because what it produces is an assessment and the question
the standard asks is whether the assessment is coherent. Conversely a method which does compute a
posterior is answerable to nothing further: the posterior is a consequence of coherence together
with exchangeability, not a requirement added to them. The distinction is not academic. The
posterior of a large model is neither computed by the algorithms in use nor approximated by the
point at which they stop, so an appeal to Bayesian foundations is properly an appeal to the
standard rather than to the algorithm; and it is the standard, not the algorithm, which the theorem
in this chapter supplies.

% ---- from generated/ch2-body.tex ----
%% Generated by build.py from the parts - do not edit.
\section{Finite Dimensional Convex Geometry}\label{ch2:sec:convex}

This section contains no probability. It proves the single geometric fact on
which the fundamental theorem rests, the alternative which says
that either a linear subspace meets the positive orthant away from the
origin, or it admits a strictly positive orthogonal vector. The proof given
below is elementary and constructive: it consists of a finitely terminating
algorithm which, given the data, exhibits one of the two alternatives
explicitly. It follows the treatment of Farkas' lemma by Dax
\cite{Dax1997}. At the end of the section we record, for comparison, the
classical route through a supporting hyperplane, and we say exactly what each
of the two proofs delivers: the classical one reduces the lemma to a
projection problem and produces the vector $q$ as a normal, while the
finitary one adds a termination bound and an existence proof that does not
rest on compactness.

\subsection*{Notation and the Two Systems}

Let $S$ be a finite set. We write $\RR^S$ for the vector space of real
functions on $S$ and $\inner{q}{v}=\sum_{s\in S}q(s)v(s)$ for the Euclidean
product. The \emph{closed positive orthant} and its interior are
\begin{equation}\label{ch2:eq:orthant}
C_+=\{v\in\RR^S:\ v(s)\ge0\ \ \forall s\in S\},\qquad
\Int C_+=\{v\in\RR^S:\ v(s)>0\ \ \forall s\in S\}.
\end{equation}
For a linear subspace $V\subseteq\RR^S$ we write $V^{\perp}=\{q:\
\inner{q}{v}=0\ \forall v\in V\}$ for the orthogonal complement, and we
recall the identity $(V^{\perp})^{\perp}=V$, valid in finite dimension. For a
subset $A$ we write $A^{\perp}=\{q:\ \inner{q}{v}=0\ \forall v\in A\}$ and
$A^{+}=\{q:\ \inner{q}{v}\ge0\ \forall v\in A\}$; thus $C_+^{+}=C_+$.

Two alternatives are needed. The first is due to Farkas \cite{Farkas1902} and
is stated for a matrix; the second is the form in which we shall use it.

\begin{theorem}[Farkas' lemma]\label{ch2:thm:farkas}
Let $A$ be a real $m\times n$ matrix and let $c\in\RR^n$, $c\neq0$. Then
exactly one of the following systems is solvable:
\begin{equation}\label{ch2:eq:farkas-primal}
A x\ge0\quad\text{and}\quad c^{T}x<0,\qquad x\in\RR^n,
\end{equation}
\begin{equation}\label{ch2:eq:farkas-dual}
A^{T}y=c\quad\text{and}\quad y\ge0,\qquad y\in\RR^m .
\end{equation}
\end{theorem}

\begin{theorem}[Orthant alternative]\label{ch2:thm:stiemke}
Let $V\subseteq\RR^S$ be a linear subspace. Then exactly one of the following
holds:
\begin{enumerate}
\item[(a)] there exists $v\in V$, $v\neq0$, with $v(s)\ge0$ for all $s\in S$;
\item[(b)] there exists $q\in V^{\perp}$ with $q(s)>0$ for all $s\in S$.
\end{enumerate}
\end{theorem}

We prove Theorem~\ref{ch2:thm:farkas} in \S\ref{ch2:ss:farkas} and deduce
Theorem~\ref{ch2:thm:stiemke} from it in \S\ref{ch2:ss:orthant}. Both proofs are
elementary: the only tool used in \S\ref{ch2:ss:farkas} is the elementary theory
of least squares over a finite dimensional space, and the deduction in
\S\ref{ch2:ss:orthant} consists of finitely many applications of
Theorem~\ref{ch2:thm:farkas}, one for each point of $S$.

\subsection{Farkas' Lemma}\label{ch2:ss:farkas}

\emph{The two systems are not both solvable.} If \eqref{ch2:eq:farkas-primal} and
\eqref{ch2:eq:farkas-dual} both hold then
\begin{equation}\label{ch2:eq:farkas-incompat}
c^{T}x=(A^{T}y)^{T}x=y^{T}Ax\ge0,
\end{equation}
because $y\ge0$ and $Ax\ge0$, contradicting $c^{T}x<0$.

It remains to show that one of them is always solvable. The idea is to look
for the point of the orthant which comes closest to solving the dual system,
and to read the answer off the residual.

\subsubsection*{The Bounded Least Squares Problem}

Consider the problem
\begin{equation}\label{ch2:eq:ls}
\text{minimize }\ \norm{A^{T}y-c}^2\quad\text{subject to }\ y\ge0,
\end{equation}
where $\norm{\cdot}$ is the Euclidean norm. Its feasible set is the closed
orthant, which is nonempty and contains the origin, and its objective is a
convex quadratic. We give below an algorithm which solves it in a finite
number of steps, and this will establish that a solution exists. For the
moment let $y_*$ be a solution and let
\begin{equation}\label{ch2:eq:residual}
r_*=A^{T}y_*-c
\end{equation}
be the corresponding residual vector.

\begin{lemma}[Optimality conditions]\label{ch2:lem:kkt}
Let $y_*\ge0$ and let $r_*=A^{T}y_*-c$. Then $y_*$ solves \eqref{ch2:eq:ls} if
and only if
\begin{equation}\label{ch2:eq:kkt}
y_*\ge0,\qquad Ar_*\ge0,\qquad y_*^{T}Ar_*=0 .
\end{equation}
\end{lemma}

\begin{proof}
Assume first that $y_*$ solves \eqref{ch2:eq:ls} and let $a_i^{T}$ denote the
$i$-th row of $A$. For $i=1,\dots,m$ consider the one parameter quadratic
\begin{equation}
f_i(\theta)=\norm{A^{T}(y_*+\theta e_i)-c}^2=\norm{\theta a_i+r_*}^2,
\end{equation}
where $e_i$ is the $i$-th column of the $m\times m$ identity. Since $y_*$
minimises \eqref{ch2:eq:ls}, $\theta=0$ minimises $f_i$ subject to the constraint
$y_{*i}+\theta\ge0$. Hence, as $f_i'(0)=2a_i^{T}r_*$,
\begin{equation}
y_{*i}>0\ \Longrightarrow\ a_i^{T}r_*=0,\qquad
y_{*i}=0\ \Longrightarrow\ a_i^{T}r_*\ge0,
\end{equation}
which is \eqref{ch2:eq:kkt}.

Conversely, assume \eqref{ch2:eq:kkt} and let $z\ge0$ be arbitrary. Set
$u=z-y_*$. If $y_{*i}=0$ then $u_i\ge0$, while if $y_{*i}>0$ then
$a_i^{T}r_*=0$; in either case $u_i\,a_i^{T}r_*\ge0$, so that
$u^{T}Ar_*\ge0$. The identity
\begin{equation}
\norm{A^{T}z-c}^2=\norm{A^{T}y_*-c}^2+2u^{T}Ar_*+\norm{A^{T}u}^2
\end{equation}
therefore shows that $\norm{A^{T}z-c}^2\ge\norm{A^{T}y_*-c}^2$, that is, that
$y_*$ solves \eqref{ch2:eq:ls}.
\end{proof}

\begin{theorem}[The residual solves the primal system]\label{ch2:thm:residual}
Let $y_*$ solve \eqref{ch2:eq:ls} and let $r_*=A^{T}y_*-c$. If $r_*=0$ then $y_*$
solves the dual system \eqref{ch2:eq:farkas-dual}. Otherwise $r_*$ solves the
primal system \eqref{ch2:eq:farkas-primal}, and in fact
\begin{equation}\label{ch2:eq:residual-identity}
c^{T}r_*=-\norm{r_*}^2 .
\end{equation}
\end{theorem}

\begin{proof}
By \eqref{ch2:eq:residual},
\begin{equation}
c^{T}r_*=(A^{T}y_*-r_*)^{T}r_*=y_*^{T}Ar_*-\norm{r_*}^2=-\norm{r_*}^2,
\end{equation}
where the last equality is the third relation in \eqref{ch2:eq:kkt}. If
$r_*=0$ then $A^{T}y_*=c$ with $y_*\ge0$ and the dual system is solved. If
$r_*\neq0$ then $Ar_*\ge0$ by \eqref{ch2:eq:kkt} and
$c^{T}r_*=-\norm{r_*}^2<0$, so that $r_*$ solves
\eqref{ch2:eq:farkas-primal}.
\end{proof}

It remains to prove that \eqref{ch2:eq:ls} has a solution. This is the
constructive part of the argument: the following algorithm finds one.

\subsubsection*{An Active Set Algorithm for the Least Squares Problem}

At the beginning of the $k$-th iteration, $k=0,1,2,\dots$, we have a feasible
point $y^k\ge0$. We set
\begin{equation}\label{ch2:eq:active-set}
r^k=A^{T}y^k-c,\qquad W_k=\{i:\ y^k_i>0\},\qquad
V_k=\{i:\ y^k_i=0\},
\end{equation}
so that $W_k\cup V_k=\{1,\dots,m\}$, and we let $t_k$ be the number of
indices in $W_k$.

\smallskip
\noindent\emph{Step 1: a least squares problem on the free variables.} If
$t_k=0$ or $r^k=0$, go to Step 2. Otherwise let $E_k\subseteq\RR^n$ be the
linear span of the rows $a_i^{T}$, $i\in W_k$, and let $P_k$ be the
orthogonal projection of $\RR^n$ onto $E_k$. If $P_kr^k=0$ go to Step 2.
Otherwise the vector $-P_kr^k$ belongs to $E_k$, which is the range of the
transpose of the $t_k\times n$ matrix $A_k$ with rows $a_i^{T}$, $i\in W_k$,
so that there exists $w\in\RR^{t_k}$ with $A_k^{T}w=-P_kr^k$; equivalently,
$w$ solves the least squares problem $\min\norm{A_k^{T}w+r^k}^2$. Define the
search direction $u\in\RR^m$ by $u_i=w_i$ for $i\in W_k$ and $u_i=0$ for
$i\in V_k$, and set
\begin{equation}\label{ch2:eq:step1}
y^{k+1}=y^k+\theta_ku,\qquad
\theta_k=\min\Big(\{1\}\cup\{-y^k_i/u_i:\ i\in W_k,\ u_i<0\}\Big)>0 .
\end{equation}
By construction $y^{k+1}\ge0$ and $y^{k+1}_i=0$ for all $i\in V_k$. Moreover
\begin{equation}\label{ch2:eq:step1-decrease}
\norm{r^{k+1}}^2=\norm{r^k+A^{T}u}^2
=\norm{r^k-P_kr^k}^2=\norm{r^k}^2-\norm{P_kr^k}^2<\norm{r^k}^2
\end{equation}
if $\theta_k=1$, and in general
$\norm{r^k+\theta A^{T}u}^2=\norm{r^k}^2-(2\theta-\theta^2)\norm{P_kr^k}^2
<\norm{r^k}^2$ because $0<\theta\le1$.

\smallskip
\noindent\emph{Step 2: testing optimality, and leaving a dead point.} Here
$P_kr^k=0$ or $t_k=0$; in the first case $A_kr^k=0$, so $y^k$ solves the
problem of minimising $\norm{A^{T}y-c}^2$ subject to $y_i=0$ for $i\in V_k$
and $y_i\ge0$ for $i\in W_k$; in the second case there are no free variables
and $y^k=0$. We call such a point a \emph{dead point}. Compute an index
$j\in V_k$ at which $a_j^{T}r^k$ is minimal. If $a_j^{T}r^k\ge0$ then
$Ar^k\ge0$, indeed $a_i^{T}r^k=0$ for $i\in W_k$ and $a_i^{T}r^k\ge0$ for
$i\in V_k$, and $y^{k\,T}Ar^k=0$; by Lemma~\ref{ch2:lem:kkt} the point $y^k$
solves \eqref{ch2:eq:ls} and the algorithm terminates. Otherwise
$a_j^{T}r^k<0$, and we set
\begin{equation}\label{ch2:eq:step2}
y^{k+1}=y^k-\frac{a_j^{T}r^k}{\norm{a_j}^2}\,e_j,
\end{equation}
which is feasible because the coefficient is positive; the residual change is
$\norm{r^{k+1}}^2=\norm{r^k}^2-(a_j^{T}r^k)^2/\norm{a_j}^2<\norm{r^k}^2$.

\begin{proposition}[Finite termination]\label{ch2:prop:termination}
Starting from $y^0=0$, the algorithm terminates after finitely many
iterations at a solution of \eqref{ch2:eq:ls}.
\end{proposition}

\begin{proof}
Four observations.

(i) \emph{The objective strictly decreases at every iteration.} This is
\eqref{ch2:eq:step1-decrease} in Step 1 and the last line of Step 2.

(ii) \emph{Between two dead points there are at most $m$ iterations.} In Step
1 the coordinates in $V_k$ are not changed. If $\theta_k=1$ then
$y^{k+1}=y^k+u$ has the same vanishing coordinates as $y^k$, so
$W_{k+1}\subseteq W_k$; and since $A^{T}u=-P_kr^k$ is orthogonal to $a_i$ for
every $i\in W_k$, the new residual satisfies $a_i^{T}r^{k+1}=0$ for
$i\in W_{k+1}$, so that $y^{k+1}$ is a dead point. If instead
$0<\theta_k<1$ then some coordinate with $u_i<0$ reaches zero, so
$t_{k+1}<t_k$. Hence the number of positive coordinates strictly decreases
unless a dead point is reached.

(iii) \emph{At a dead point, the objective already equals the exact minimum
over the face determined by $V_k$.} This is what the condition $A_kr^k=0$
says, by Lemma~\ref{ch2:lem:kkt} applied to the subproblem with the variables in
$V_k$ frozen at zero.

(iv) \emph{A subset $V_k$ cannot occur twice.} By (i) and (iii), returning to
the same face would give the same minimum of the objective on that face,
contradicting the strict decrease.

Since there are only $2^m$ subsets $V_k$, and by (ii) at most $m$ iterations
separate two dead points, the algorithm performs at most $m2^m$ iterations.
At termination Lemma~\ref{ch2:lem:kkt} gives the optimality conditions
\eqref{ch2:eq:kkt}, so the point reached solves \eqref{ch2:eq:ls}.
\end{proof}

\begin{proof}[Proof of Theorem~\ref{ch2:thm:farkas}]
By \eqref{ch2:eq:farkas-incompat} the two systems cannot both be solvable. If
$c=0$ then $y=0$ solves the dual system. If $c\neq0$,
Proposition~\ref{ch2:prop:termination} provides a solution $y_*$ of
\eqref{ch2:eq:ls}, and Theorem~\ref{ch2:thm:residual} shows that either $y_*$ solves
the dual system or $r_*$ solves the primal one. Hence one of the two systems
is always solvable.
\end{proof}

\begin{remark}\label{ch2:rem:farkas-constructive}
The proof is an algorithm, and it delivers more than the statement of the
alternative. Given $A$ and $c$, the active set method of
Proposition~\ref{ch2:prop:termination} computes, in at most $m2^{m}$ iterations
each of which consists of the solution of one linear system of size at most
$m$, either a solution $y$ of the dual system or a vector $r$ with
$Ar\ge0$ and $c^{T}r<0$, the latter with the exact value
$c^{T}r=-\norm{r}^2$. No compactness argument and no separation theorem is
used; the only ingredient is the elementary fact that a convex quadratic on a
finite dimensional space has a minimiser, which is supplied by the finite
termination of the algorithm itself.
\end{remark}

\begin{corollary}[Steepest descent]\label{ch2:cor:powell}
Let $y_*$ solve \eqref{ch2:eq:ls} and let $r_*\neq0$. Then $r_*/\norm{r_*}$
solves the problem
\begin{equation}
\text{minimize }\ c^{T}x\quad\text{subject to }\ Ax\ge0,\ \norm{x}=1 .
\end{equation}
\end{corollary}

\begin{proof}
Let $x$ satisfy the constraints. Then $y_*^{T}Ax\ge0$, and
$\abs{r_*^{T}x}\le\norm{r_*}$ by Cauchy--Schwarz. Hence
\begin{equation}
c^{T}x=(A^{T}y_*-r_*)^{T}x=y_*^{T}Ax-r_*^{T}x\ge-\abs{r_*^{T}x}\ge-\norm{r_*},
\end{equation}
and $c^{T}r_*/\norm{r_*}=-\norm{r_*}$ by \eqref{ch2:eq:residual-identity}.
\end{proof}

\subsection{The Orthant Alternative}\label{ch2:ss:orthant}

\begin{proof}[Proof of Theorem~\ref{ch2:thm:stiemke}]
\emph{(a) and (b) are incompatible.} Let $q\in V^{\perp}$ be strictly
positive and let $v\in V\cap C_+$ with $v\neq0$. There is $s_0$ with
$v(s_0)>0$, and every term of $\inner{q}{v}=\sum_sq(s)v(s)$ is non negative,
so $\inner{q}{v}>0$, contradicting $v\in V$ and $q\in V^{\perp}$.

\emph{If (a) fails then (b) holds.} Assume $V\cap C_+=\{0\}$. Choose a matrix
$\Xi$ whose columns span $V$, so that $V=\{\Xi\zeta:\ \zeta\in\RR^n\}$ and
$V^{\perp}=\{q:\ \Xi^{T}q=0\}$. We show that for each $s\in S$ there is a
vector $q^{(s)}\in\RR^S$ with
\begin{equation}\label{ch2:eq:q-s}
q^{(s)}\ge0,\qquad \Xi^{T}q^{(s)}=0,\qquad q^{(s)}(s)>0 .
\end{equation}
Fix $s$ and let $e_s$ be the $s$-th unit vector. If $\Xi^{T}e_s=0$, that is if
the $s$-th row of $\Xi$ vanishes, then $q^{(s)}=e_s$ satisfies
\eqref{ch2:eq:q-s}. Otherwise set $c=-\Xi^{T}e_s\neq0$ and apply
Theorem~\ref{ch2:thm:farkas} to the matrix $\Xi$ and to this $c$. The primal
system reads: there is $\zeta$ with $\Xi\zeta\ge0$ and
$0>c^{T}\zeta=-(\Xi\zeta)(s)$, that is, a non zero element of $V\cap C_+$,
which is excluded. Hence the dual system is solvable: there is $q\ge0$ with
$\Xi^{T}q=-\Xi^{T}e_s$, and then $q^{(s)}=q+e_s$ satisfies \eqref{ch2:eq:q-s}
because $\Xi^{T}q^{(s)}=0$ and $q^{(s)}(s)=q(s)+1>0$.

Finally set $q=\sum_{s\in S}q^{(s)}$. Then $q\ge0$, $\Xi^{T}q=0$ and
$q(s)\ge q^{(s)}(s)>0$ for every $s$, so $q$ is strictly positive. Hence (b)
holds.
\end{proof}

\begin{remark}\label{ch2:rem:two-routes}
The deduction above uses finitely many applications of Farkas' lemma, one per
point of $S$, and each of them is decided by the algorithm of
\S\ref{ch2:ss:farkas}. This is the constructive route to the alternative. For
comparison we record in \S\ref{ch2:ss:separation} the shorter geometric route,
which produces $q$ at once as the normal to a supporting hyperplane and which
is not constructive.
\end{remark}

\begin{corollary}[Normalized state prices]\label{ch2:cor:normalized}
Let $V\subseteq\RR^S$ be a linear subspace, let $\one\in\RR^S$ be the constant
function $1$, and let $S_0=\{q\in V^{\perp}:\ q>0,\ \inner{q}{\one}=1\}$. If
$V\cap C_+=\{0\}$ then $S_0$ is a nonempty, convex, relatively open subset of
$V^{\perp}$; it consists of a single point if and only if
$\dim V^{\perp}=1$.
\end{corollary}

\begin{proof}
Non emptiness is Theorem~\ref{ch2:thm:stiemke} followed by normalization by
$\inner{q}{\one}>0$. Convexity is clear. Relative openness: if $q\in S_0$
then all $q'\in V^{\perp}$ sufficiently close to $q$ still satisfy
$q'>\min_sq(s)/2>0$. If $\dim V^{\perp}=1$ then all elements of $S_0$ lie on
one ray and, being normalized, coincide. If $\dim V^{\perp}\ge2$ then $S_0$,
being relatively open in $V^{\perp}$, has nonempty relative interior and more
than one element.
\end{proof}

\begin{remark}\label{ch2:rem:two-forms}
Applying Theorem~\ref{ch2:thm:stiemke} to $V^{\perp}$ in place of $V$, and using
$(V^{\perp})^{\perp}=V$, gives the other classical form of the alternative,
due to Stiemke \cite{Stiemke1915}: exactly one of
\begin{enumerate}
\item[(a$'$)] there exists $x\in V$ with $x(s)>0$ for all $s$;
\item[(b$'$)] there exists $q\in V^{\perp}$, $q\neq0$, with $q(s)\ge0$ for all
$s$.
\end{enumerate}
Theorem~\ref{ch2:thm:stiemke} and this statement, together with
Theorem~\ref{ch2:thm:farkas}, are the three faces of the same lemma.
\end{remark}

\subsection{The Span Calculus of One-Period Models}

The following linear algebra will be used to pass between the payoff space of
a one period model and the space spanned by its gains.

\begin{lemma}\label{ch2:lem:spancalc}
Let $A$ be a finite set, $a_0\in A$, and let $h_1,\dots,h_r\in\RR^A$ be such
that the constants are contained in $W=\Span\{h_1,\dots,h_r\}$. Set
\begin{equation}\label{ch2:eq:V-span}
V=\Span\{h_j-h_j(a_0)\one:\ j=1,\dots,r\}\subseteq\RR^A .
\end{equation}
Then
\begin{enumerate}
\item[(i)] $\dim V=\dim W-1$;
\item[(ii)] $V^{\perp}=\{q\in\RR^A:\ \inner{q}{h_j}=h_j(a_0)\inner{q}{\one}
\ \forall j\}$;
\item[(iii)] the market is complete, i.e.\ $W=\RR^A$, if and only if
$\dim V^{\perp}=1$.
\end{enumerate}
\end{lemma}

\begin{proof}
The linear map $\Phi:W\to\RR^A$, $\Phi(h)=h-h(a_0)\one$, has image $V$ and
kernel the constants, which form a one dimensional subspace of $W$ because
$\one\in W$. Hence $\dim V=\dim W-1$, which is (i). For (ii), $q\in V^{\perp}$
means $\inner{q}{h_j}=h_j(a_0)\inner{q}{\one}$ for every $j$, since $q$ is
linear. For (iii), $\dim V=\dim A-1$ if and only if $\dim W=\dim A$, and
$\dim V^{\perp}=\dim A-\dim V$.
\end{proof}

\subsection{The Classical Geometric Proof, for Comparison}\label{ch2:ss:separation}

The standard proof of Theorem~\ref{ch2:thm:stiemke} separates the subspace from
the simplex by a hyperplane; the vector $q$ is the normal to that hyperplane.
We record it because it explains the geometry, and because the comparison
clarifies what the finitary argument of \S\ref{ch2:ss:farkas} adds. As will be
made precise in Remark~\ref{ch2:rem:compare-proofs}, this proof is constructive
in the weak sense that it reduces the lemma to the computation of one
projection, but it takes the existence of the projection for granted and it
delivers neither a termination bound nor any bound on $q$.

\begin{lemma}\label{ch2:lem:separate}
Let $K\subset\RR^S$ be compact and convex, and let $V\subseteq\RR^S$ be
closed and convex. If $K\cap V=\emptyset$, then there exist $q\in\RR^S$,
$q\neq0$, and a real number $\alpha$ such that
\begin{equation}\label{ch2:eq:separate}
\inner{q}{k}>\alpha>\inner{q}{v}\qquad\text{for all } k\in K,\ v\in V .
\end{equation}
\end{lemma}

\begin{proof}
The function $k\mapsto d(k,V)=\inf\{\norm{k-v}:\ v\in V\}$ is continuous on
the compact set $K$, so it attains a minimum $\delta\ge0$. If $\delta=0$ then
some $k_0\in K$ is a limit of points of $V$, hence belongs to $V$ because $V$
is closed; this is excluded, so $\delta>0$. Choose $k_*\in K$ with
$d(k_*,V)=\delta$. The set $V\cap\{v:\norm{v-k_*}\le\delta+1\}$ is compact and
nonempty, so $v\mapsto\norm{v-k_*}$ attains on it a minimum, at a point
$v_*\in V$, and $\norm{k_*-v_*}=\delta$.

Set $u=k_*-v_*$, so that $\norm{u}=\delta>0$, and $q=u$. \emph{First claim}:
$\inner{q}{v-v_*}\le0$ for every $v\in V$. Indeed $v_*+t(v-v_*)\in V$ for
$t\in[0,1]$, so $\norm{k_*-v_*-t(v-v_*)}^2\ge\delta^2$; expanding and dividing
by $t>0$ gives $t\norm{v-v_*}^2-2\inner{q}{v-v_*}\ge0$, and letting
$t\downarrow0$ yields the claim. \emph{Second claim}:
$\inner{q}{k-k_*}\ge0$ for every $k\in K$. Otherwise, for some $k\in K$ with
$\inner{q}{k-k_*}<0$, the point $k_t=k_*+t(k-k_*)\in K$, $t\in(0,1)$, would
satisfy
\begin{equation}
\norm{k_t-v_*}^2=\delta^2+2t\inner{q}{k-k_*}+t^2\norm{k-k_*}^2<\delta^2
\end{equation}
for $t$ small enough, contradicting the definition of $\delta$.

Now set $\alpha=\inner{q}{k_*}-\delta^2/2$. By the two claims,
$\inner{q}{k}\ge\inner{q}{k_*}=\inner{q}{v_*}+\delta^2>\alpha
>\inner{q}{v_*}\ge\inner{q}{v}$ for all $k\in K$ and $v\in V$, which is
\eqref{ch2:eq:separate}.
\end{proof}

\begin{proof}[Second proof of Theorem~\ref{ch2:thm:stiemke}]
Let $K=\{v\in C_+:\ \sum_sv(s)=1\}$, the standard simplex, which is compact
and convex. If $V\cap C_+=\{0\}$ then $K\cap V=\emptyset$, because every
$v\in K$ lies in $C_+$ and is non zero. By Lemma~\ref{ch2:lem:separate} there are
$q\neq0$ and $\alpha$ with $\inner{q}{k}>\alpha>\inner{q}{v}$ for all $k\in K$
and $v\in V$. Since $0\in V$ we have $\alpha>0$. For $v\in V$ and real $t$ the
point $tv$ lies in $V$, so $t\inner{q}{v}<\alpha$ for all $t$, which forces
$\inner{q}{v}=0$; hence $q\in V^{\perp}$. Finally $e_s\in K$ gives
$q(s)=\inner{q}{e_s}>\alpha>0$ for every $s$.
\end{proof}

\begin{remark}[What the two proofs deliver]\label{ch2:rem:compare-proofs}
The proof of Lemma~\ref{ch2:lem:separate} is often described as constructive, and
up to a point it is: the normal $q=k_*-v_*$ is the difference between a point
of the simplex and its projection on $V$, and a projection on a closed convex
set can be computed by an iterative algorithm. The qualification is twofold.

First, the construction presupposes that the projection exists. In the proof
this is supplied by compactness, the simplex being compact; but the object
actually wanted is the projection of a point on a cone, which is unbounded,
and recovering existence there requires a separate coercivity argument and the
theorem that a finitely generated cone is closed. Second, and more to the
point for a numerical framework, the proof says nothing about how many
operations are required, nor does it bound the resulting $q$; and the answer
describes $q$ by a limit rather than exhibiting a finite certificate which can
be checked.

The finitary argument of \S\ref{ch2:ss:farkas} supplies both. Its termination
proof is at the same time the existence proof for the least squares problem
\eqref{ch2:eq:ls}, so no compactness enters; the number of iterations is at most
$m2^m$, each of which solves one linear system of size at most $m$; and the
output is a certificate of one of the two alternatives, with the exact value
$c^{T}r_*=-\norm{r_*}^2$ of \eqref{ch2:eq:residual-identity}. In the applications
below this is the relevant form, because at each node of a lattice model the
data are a finite matrix and the state price vector is precisely what one
wants to compute.
\end{remark}

\begin{remark}\label{ch2:rem:construction-of-q}
The deduction of Theorem~\ref{ch2:thm:stiemke} from Theorem~\ref{ch2:thm:farkas} in
\S\ref{ch2:ss:orthant} is constructive as well, and it costs one call of the
algorithm per point of $S$: the vector $q^{(s)}$ of \eqref{ch2:eq:q-s} is obtained
from either the dual solution, shifted by $e_s$, or from $e_s$ itself when the
$s$-th row of $\Xi$ vanishes, and $q=\sum_sq^{(s)}$ is then strictly
positive. In the notation of \S\ref{ch2:ss:farkas} this produces a state price
vector by at most $|S|$ runs of an algorithm with a certified iteration count,
rather than by a limiting argument.
\end{remark}

% ---- from sections/ch2-one-period-market.tex ----
% A one-period market, appended to chapter 2. The results are those of Financial Derivatives
% Modeling, chapter 0, restated in the vocabulary of this chapter: state prices instead of a
% risk-neutral measure, and statements about the states of the market instead of statements
% almost surely.

\providecommand{\SS}{\mathcal{S}}
\providecommand{\um}{u}
\providecommand{\md}{d}

\section{A One-Period Market\label{ch2:sec:market}}

The sections above are geometry; this one is the market they were built for. A one-period market is a
finite set $\SS$ of states, a bond worth $B>0$ at time $0$ and one at maturity, and a stock worth $S\ge0$
at time $0$ and $S_T(s)$ in the state $s$. The traded payoffs are the two functions $1$ and $S_T$, at the
given prices $B$ and $SA$, where the constant $A>0$ accounts for the dividends of the stock over the
period. The market is coherent, in the sense of \S\ref{ch2:sec:convex}, when no combination of the two
has a cost which is nonpositive and a payoff which is nonnegative everywhere and positive somewhere.

By Theorem~\ref{ch2:thm:stiemke} coherence is the same thing as the existence of a state price
$q\in\Int C_+$ with
\begin{equation}\label{ch2:eq:marketq}
\inner{q}{1}=B,\qquad \inner{q}{S_T}=SA .
\end{equation}
Everything below is a consequence of \eqref{ch2:eq:marketq} and of nothing else: the derived prices are
the values $\inner{q}{\xi}$, and a statement about the market is a statement about $q$.

\subsection*{Forwards, Calls and Puts}

A European payoff is a function $\xi$ of $S_T$, and a European call or put of strike $K$ has payoff
$(S_T-K)^{+}$ or $(S_T-K)^{-}$; a forward has payoff $S_T-K$, and a combination of calls and puts with
different strikes generates the spreads of every textbook --- bull and bear spreads, straddles,
butterflies --- by linearity of $\inner{q}{\cdot}$.

\begin{proposition}\label{ch2:p:parity}
With $F$, $C$ and $P$ the prices of the forward, the call and the put of the same maturity and strike,
\begin{equation}\label{ch2:eq:parity}
F=SA-KB=C-P .
\end{equation}
\end{proposition}

\begin{proof}
The payoff identity $S_T-K=(S_T-K)^{+}-(S_T-K)^{-}$ holds at every state, so applying the state prices
gives the second equality; the first follows from \eqref{ch2:eq:marketq}, since
$\inner{q}{S_T-K}=SA-KB$. Both are the linearity of $\inner{q}{\cdot}$ and nothing more, which is worth
noting: the identity is an identity of payoffs before it is an identity of prices.
\end{proof}

\begin{proposition}\label{ch2:p:bounds}
$(SA-KB)^{+}\le C\le SA$ and $(KB-SA)^{+}\le P\le KB$, for $S,K\ge0$.
\end{proposition}

\begin{proof}
The payoff of a call dominates $(S_T-K)^{+}\ge0$ and is dominated by $S_T$, and $q>0$, so the three
prices are ordered as the payoffs are. The lower bound on the call follows from
$(S_T-K)^{+}\ge S_T-K$, whose price is $SA-KB$, and the bounds on the put from
\eqref{ch2:eq:parity}.
\end{proof}

\begin{proposition}\label{ch2:p:shape}
As functions of the strike, $C$ and $P$ are respectively nonincreasing and nondecreasing, $B$-Lipschitz,
and convex.
\end{proposition}

\begin{proof}
Each is the statement that a certain payoff is nonnegative at every state, followed by the positivity
of $q$: monotonicity from the call spread, Lipschitz continuity from the same spread with strikes $K_1$
and $K_2$ and the bound $|K_1-K_2|$ on its payoff, convexity from the butterfly
$(S_T-K-k)^{+}-2(S_T-K)^{+}+(S_T-K+k)^{+}\ge0$, for which the three prices give the second difference of
$C$.
\end{proof}

\subsection*{Two States, and Completeness}

Let $\SS=\{u,d\}$ with $S_T(u)=\um S$ and $S_T(d)=\md S$ for constants $0\le\md<\um$. Two equations
\eqref{ch2:eq:marketq} in two unknowns have the unique solution
\begin{equation}\label{ch2:eq:twostates}
q(u)=B\,\frac{A/B-\md}{\um-\md},\qquad q(d)=B\,\frac{\um-A/B}{\um-\md},
\end{equation}
positive exactly when $\md<A/B<\um$, which is therefore the condition of coherence of the market. Every
European payoff is a combination of the two traded ones, so every price is computed: the market is
complete in the sense of \S\ref{ch2:sec:convex}, and the computed prices are
\begin{equation}\label{ch2:eq:replication}
\inner{q}{\phi(S_T)}=q(u)\phi(\um S)+q(d)\phi(\md S),
\end{equation}
which is the price of the replicating portfolio of $\phi(S_T)$ in the bond and the stock. Its weights
solve the two linear equations which state, state by state, that the portfolio pays $\phi(S_T)$; that
the system is solvable is the spanning, and the value of the portfolio at time $0$ is
\eqref{ch2:eq:replication}, which is the theorem of \S\ref{ch2:sec:convex} read backwards.

Three states would give three equations in two unknowns and, generically, no solution: the market would
be coherent and incomplete, and the prices of the payoffs which are not spanned would not be determined
by the traded ones. Nothing in the theory supplies them, and any number quoted for such a payoff is an
assumption about the market and not a consequence of it.
% ---- end sections/ch2-one-period-market.tex ----

% ---- end generated/ch2-body.tex ----

\section{The Proof of the Fundamental Theorem}
\label{sec:proof-oneperiod}

We now prove Theorem~\ref{ch2:thm:oneperiod}.

% ---- from generated/ch2-theorem-proof.tex ----
%% Generated by build.py from the parts - do not edit.
\begin{proof}
Let $V=\Span\{\xi^1,\dots,\xi^n\}\subseteq\RR^S$. By definition of $V$, a
non zero element of $V\cap C_+$ is precisely a gain vector
$\sum_k\zeta^k\xi^k$ which is non negative and not identically zero, that is,
precisely a choice of $\zeta$ forbidden by (i). Hence (i) is the statement
$V\cap C_+=\{0\}$, while (ii) is the statement that $V^{\perp}$ contains a
strictly positive vector; the two are equivalent by
Theorem~\ref{ch2:thm:stiemke}.
\end{proof}
% ---- end generated/ch2-theorem-proof.tex ----

What the proof delivers is more than the equivalence, because the orthant alternative is proved in
\S\ref{ch2:sec:convex} by an algorithm and not merely by a separation argument. The algorithm
decides which of the two statements holds, and when the first fails it returns the combination
$\zeta$ of gambles which refutes it --- that is, an explicit arbitrage; when the first holds it
returns the vector $g$, that is, an explicit pricing kernel. Either answer is therefore a
certificate which can be checked by recomputing a finite sum, and neither requires the
axiom of choice, a separating hyperplane in an infinite dimensional space, or a measurable
selection. In Chapter~\ref{chap:ch3} the theorem is applied at each site of a lattice, and the
local certificates it produces there are assembled into a pricing measure for the whole process.
% ---- end chapters/ch2.tex ----
% ---- from chapters/ch3.tex ----
% Authored: the chapter title and opening belong to the book. The body is
% generated by build.py from Part 1, from the theory of processes onwards.
\chapter{Stochastic Processes on a Lattice}
\label{chap:ch3}

This chapter introduces the objects which the rest of the book is about. A stochastic process is
defined here by prescribing the expectation of every functional of its path at a fixed resolution
--- that is, by a consistent family of $n$-point functions on a lattice --- and not by a measure on
a path space. At each resolution every object of the theory is then a finite matrix or a finite
sum, conditioning is a matrix operation, a measure change is a strictly positive reweighting, and
the martingale property is an identity between matrices. Nothing is asserted almost surely and
nothing is defined only up to a null set.

% ---- from generated/ch3-front.tex ----
%% Generated by build.py from the parts - do not edit.
\subsection*{The Logical Architecture} The chapter is organized so that nothing is used before it is proved; the finite dimensional convex geometry on which it rests is proved in Chapter~\ref{chap:ch2}, and is not repeated here.

\begin{description}
\item[\S\ref{ch3:sec:simplicial} Simplicial sequences and lattice functions]
The geometric substrate: lattices, hulls, interiors, simplicial sequences,
lattice functions, integrals, $L^p$ spaces, absolute continuity.

\item[\S\ref{ch3:sec:pathfunc} Path functionals] Path space at resolution $m$
and the algebras of functionals on it.

\item[\S\ref{ch3:sec:npoint} Filtered probability spaces and $n$-point
functions] The definition of a stochastic process as a positive normalized
family of $n$-point functions, together with the operation of conditioning.

\item[\S\ref{ch3:sec:markov} Markov processes] Propagators, generators,
existence and uniqueness of the propagator, the Dyson expansion and the
path integral representation, piecewise smooth processes, attainable sets,
reachability.

\item[\S\ref{ch3:sec:measure} Measure changes] The reweighting formula, and
the proof that it maps Markov generators to Markov generators. This is the
constructive substitute for Girsanov's theorem.

\item[\S\ref{ch3:sec:martingale} Martingales and monotone processes]
The local characterization of martingales (Proposition~\ref{ch3:prop:localglobal}),
the theorem that a monotone martingale is constant
(Theorem~\ref{ch3:thm:monotone-martingale}), and the resulting pathwise criterion
for the absence of arbitrage.

\item[\S\ref{ch3:sec:ftf} The Fundamental Theorem]
Financial models, self--financing strategies, arbitrage, and the theorem
itself (Theorem~\ref{ch3:thm:ftf}), both directions, followed by the valuation
theorem for contingent claims (Theorem~\ref{ch3:thm:valuation}) and the
uniqueness theorem (Theorem~\ref{ch3:thm:second}). The section ends with the
binomial model as a worked example and with an account of what the theorem
does and does not say.

\item[\S\ref{ch3:sec:lean} Machine-checked verification] What part of the chapter
has been formalised in Lean~4, by which of the two routes of
\S\ref{ch2:sec:convex} the certificate is obtained, and what the rest would
take.
\end{description}

\subsection*{Statement of the Main Result}

This chapter contains a self-contained and elementary proof of the fundamental theorem for a multi period market model built on a Markov chain over a finite lattice. Path
spaces are described by families of $n$-point functions; martingales are characterized by a
local differential identity holding at each site and each date; the absence of arbitrage is
shown to be a local condition; and the one period problem is reduced to the orthant
alternative of finite dimensional convex geometry. That alternative is proved twice: by a
finitely terminating active set algorithm which constructs either an arbitrage or a state
price vector, and, for comparison, by a separating hyperplane. The resulting proof needs no
measurable selection and no infinite dimensional separation theorem, and it exhibits the
pricing measure explicitly as a positive reweighting of the transition rates. Valuation and
hedging of contingent claims, the uniqueness theorem, and the binomial model as a worked
example follow from it, and the convex geometry at the base of the argument is verified
mechanically in Lean~4.

The theorem itself, in the notation introduced in
\S\S\ref{ch3:sec:simplicial}--\ref{ch3:sec:martingale}, is the following. Fix a
resolution $m$ and a filtered probability space generated by a Markov
propagator $U_m$ on the simplicial sequence $(\StateA)$; fix a financial model
with asset price processes $A^0_m,\dots,A^n_m$, where $A^0_m$ is the money
market account $B_m$, and let $\tA^k_m=A^k_m/B_m$ denote the discounted
prices.

\begin{quote}
\itshape
If there exists an equivalent measure --- that is, a measure change
$(G_m^{yt})$ with $G_m^{yt}>0$ on the attainable set of every site --- under
which all discounted price processes $\tA^1_m,\dots,\tA^n_m$ are
martingales, then the discounted value process of every self financing
trading strategy is a martingale under that measure, and there is no
arbitrage. Conversely, if there is no arbitrage, then such a measure exists.
\end{quote}

Both implications are proved in \S\ref{ch3:sec:ftf}. The first is a consequence
of the local characterization of martingales proved in
\S\ref{ch3:sec:martingale} and of the fact that the value of a self financing
strategy moves only because the prices of the assets it holds move. The
second is a consequence of the fact that the absence of arbitrage is
equivalent to a purely local condition at each site and each time
(Proposition~\ref{ch3:prop:localNA}), together with Stiemke's lemma
(Theorem~\ref{ch2:thm:stiemke}), which produces a strictly positive state price vector
at each site, and with the fact, proved in \S\ref{ch3:sec:measure}, that a
positive state price vector is the same thing as a measure change.
% ---- end generated/ch3-front.tex ----

\subsection*{Probability Still Does Not Exist}

The fundamental theorem of Chapter~\ref{chap:ch2} was stated for a finite set of states and a
finite family of gambles, with no probability space anywhere in sight, and none was needed: the
pricing measure was the conclusion of the argument rather than its premise. In this chapter, and in
those which follow, we nevertheless adopt the traditional presentation of the theory of stochastic
processes, and formulate the fundamental theorem for a process which is given from the
outset together with its probability space, its filtration and its conditional expectations.

The status of that choice should be clear to the reader. In Kolmogorov's approach the probability
space is not a choice at all but a necessity: a process \emph{is} a family of measurable functions on
a measure space, and there is nothing else from which to start. Here the position is the reverse. As
Chapter~\ref{chap:ch2} established, probability is a derived concept; what is primitive is a coherence
requirement on the prices of finitely many gambles, and the measure is what a coherent assessment
determines rather than what it presupposes. We could therefore develop the whole of the present
chapter --- the definition of a process, conditioning, the change of measure, and the theorem itself
--- from the coherence requirement alone, and the argument would be no weaker for it.

We do not, and the reason is pedagogical rather than technical. The vocabulary of filtrations and
conditional expectations is the one in which the subject is taught, in which the classical results
are stated, and in which a reader arriving from the standard literature will recognise what is being
proved. Presenting these theorems in that vocabulary, and marking the places where the measure could
be dispensed with, tells the reader more than presenting them in an unfamiliar form and asking him to
translate. Nothing below depends on the probability space being primitive: at a fixed resolution the
properties of the measure which the argument uses are the ones which the coherence requirement
already forces, and where such a property is used we say so.

One consequence of the identification is worth stating before the technical development begins,
since it is what the local form of the theorem buys at every node and since \S\ref{sec:two-problems}
returns to it. At a fixed resolution a process is a family of propagators, and the fundamental
theorem applied at a single site --- which is how \S\ref{ch3:ss:oneperiod} uses it --- says that a
family of traded claims admits a pricing propagator if and only if no combination of them yields a
sure gain. That gives the two problems of practice their precise form in this setting. To calibrate
is to choose, among the propagators consistent with the prices of a set of traded claims, the one
which agrees with them; the prices so imposed determine the propagator on the span of those claims,
and by the corollary on non replicable claims they leave it undetermined elsewhere. To solve the
direct problem is to evaluate a functional of the path against the chosen propagator, which at a
fixed resolution is a single product of matrices and which returns the market's value whenever the
claim lies in the span. The chapter supplies both halves: the pricing propagator, with the martingale
identity which characterises it, and the hedging and completeness results which say when the direct
answer is forced rather than chosen.

The same two problems, and the same matrices, serve the learning reading. To fit is to choose, among
the propagators consistent with what has been observed of a sample of paths, the one which reproduces
it; to predict is to evaluate a functional of a path against the propagator so chosen, whether the
functional is a conditional quantile, the expectation of a feature accumulated along the trajectory,
or the law of a statistic of the sample. The completeness results have their counterpart as well:
they are the statement that the observations determine the answer, and that where they do not, the
choice of propagator is the whole content of the model. Nothing in what follows is specific to one
reading, and the matrices do not know which of the two they are serving; what changes between them is
only what the data are, a price in the one and a realisation in the other.

% ---- from generated/ch3-body.tex ----
%% Generated by build.py from the parts - do not edit.
\subsection*{Grounds and Hypotheses}

The proof is complete in the following sense. The only geometric input is
Farkas' lemma, proved in \S\ref{ch2:sec:convex} by an algorithm which terminates
in a finite number of steps and which exhibits one of the two alternatives
explicitly; no separation theorem, no compactness argument, and no form of the
axiom of choice is used in the main line of the argument, the geometric route
being recorded only for comparison. The probabilistic
input is the existence and uniqueness of the propagator of a finite state
generator, proved in \S\ref{ch3:sec:markov} by the convergence of the Dyson
series; the a.s.\ finiteness of the number of jumps, proved in the same
section; and the identity between the local and the global notion of
martingale, proved in \S\ref{ch3:sec:martingale} by differentiating a finite
dimensional linear ordinary differential equation. There is no use of
measurable selection, of the Hahn--Banach theorem in infinite dimension, or
of the theory of the covariation of semimartingales.

It should be clear what is \emph{not} claimed. The theorem is proved for
each fixed resolution index $m$. It is a theorem about the model at
resolution $m$, not directly about a limit model. Passing to the limit
requires knowing that the pricing measures converge, and that the
martingale property survives the limit; the estimates that make this
possible are those of Chapter~5, and we return to the point in \S\ref{ch3:sec:remarks}. We also do not prove the general theorem of
\cite{DS1994}, whose hypotheses (no free lunch with vanishing risk on a
general filtered probability space) are far weaker than ours. The gain is
in the strength of the conclusion: within this framework the theorem is
elementary, explicit, and identical in structure to the algorithm that
computes prices in practice.

\subsection*{Conventions}

Throughout, $T<T'$ are real numbers and $[T,T']$ is the fixed trading
horizon. The resolution index $m$ runs over the integers $m\ge m_0$, $h_m=2^{-m}$,
and $\StateA\subset h_m\ZZ^d$ is the set of sites available at resolution $m$
(Definition~\ref{ch3:def:simplicial}). A point of $\StateA$ is denoted by
$x,y,z,\gamma_j,\bar y$; time by $s,t,u$. All functions of the site variable
are real valued unless stated otherwise, the generalization to complex
valued functions being pointwise. Sums over a set of sites or of states are
finite sums. Transposes are not indicated explicitly; a row vector acts on
a function of the state by the usual matrix product. We write
$\one_A$ for the indicator of a set $A$ and $\delta_{xy}$ for the Kronecker
symbol.

\section{Simplicial Sequences and Lattice Functions}\label{ch3:sec:simplicial}

\subsection*{Lattices and Their Geometry}

Let $d\ge1$ be an integer and let $h_m=2^{-m}$.

\begin{definition}[Lattices]\label{ch3:def:lattices}
For $A\subseteq h_m\ZZ^d$ the \emph{convex hull} of $A$ is denoted
$\Hull(A)$. The \emph{interior} of $A\subseteq h_m\ZZ^d$ is denoted
$\Int(A)$ and is the set of sites $x\in A$ which belong to $A$ together
with each of their neighbours in $h_m\ZZ^d$ at distance $h_m$.
\end{definition}

\begin{definition}[Simplicial sequences]\label{ch3:def:simplicial}
A \emph{bounded simplicial sequence} is given by an integer $m_0\ge0$ and a
sequence of finite sets $\StateA\subseteq h_m\ZZ^d$, defined for all $m\ge m_0$,
such that
\begin{enumerate}
\item[(S1)] $\Hull(\StateA)\subseteq\Hull(A_{m'})$ whenever $m<m'$;
\item[(S2)] for every $m$ and every $x\in\Int(\StateA)$ there are $M$ and
$\epsilon>0$ such that $y\in A_{m'}$ for all $m'>M$ and all
$y\in h_{m'}\ZZ^d$ with $\abs{x-y}<\epsilon$.
\end{enumerate}
\end{definition}

Condition (S1) says that the resolutions are nested in the sense of convex
hulls, and (S2) that interior points survive refinement, so that a lattice
function defined below is the restriction of a function on the continuous
domain to the lattice. We follow a constructivist discipline: to identify a
set, or a sequence of sets, one has to state how to construct it, by a
recursive algorithm whose steps can be decided in a finite number of logical
operations. Definition~\ref{ch3:def:simplicial} is of this type, since each
$\StateA$ is a finite object described by a decision procedure.

\subsection*{Lattice Functions}

\begin{definition}[Lattice functions]\label{ch3:def:latticefunction}
Let $A=(\StateA)$, $m\ge m_0$, be a bounded simplicial sequence. A
\emph{real valued simplicial function} $f:A\to\RR$ is a sequence of
functions $f_m:\StateA\to\RR$ such that $f_{m'}(x)=f_m(x)$ for all
$m'>m\ge m_0$ and all $x\in\Int(\StateA)$. The function $f$ is
\emph{uniformly bounded} if there is $c>0$ with $\abs{f_m(x)}\le c$ for all
$x\in\StateA$ and all $m$; it is \emph{uniformly continuous} if it is uniformly
bounded and for every $\epsilon>0$ there is $\delta>0$ such that
$\abs{f(x)-f(y)}<\epsilon$ whenever $x,y$ lie in a common lattice and
$\abs{x-y}<\delta$.
\end{definition}

\begin{definition}[Equivalence of lattice functions]\label{ch3:def:equivf}
Let $A=(\StateA)$ be a bounded simplicial sequence and let $f=(f_m)$,
$g=(g_m)$ be uniformly continuous lattice functions. They provide
\emph{equivalent representations of the same function} if
$\lim_{m\to\infty}\abs{f_m(x)-g_m(x)}=0$ for every $x$ in the union of the
lattices.
\end{definition}

Let $\mathcal{C}(A)$ denote the set of continuous functions on $A$,
endowed with the structure of a $C^*$ algebra given by the operations of sum,
multiplication by a scalar, and pointwise multiplication.

\begin{definition}[Integrals]\label{ch3:def:integral}
Let $A=(\StateA)$ be a bounded simplicial sequence. An \emph{integral} is a
sequence $I=(I_m)$ where $I_m$ is a linear functional on the space of
functions $f_m:\StateA\to\RR$ such that
\begin{enumerate}
\item[(I1)] $I_m(f_m)\ge0$ whenever $f_m(x)\ge0$ for all $x\in\StateA$;
\item[(I2)] the limit $I(f)=\lim_{m\to\infty}I_m(f_m)$ exists for every
continuous $f=(f_m)\in\mathcal{C}(A)$;
\item[(I3)] $\abs{I_m(f_m)}\le c\norm{f_m}_\infty$ for some constant $c>0$
independent of $m$ and $f_m$.
\end{enumerate}
The integral corresponds to a \emph{probability measure} if $I_m(\one)=1$
for all $m$.
\end{definition}

For $p\ge1$ the formula $\norm{f}_p=(I(\abs{f}^p))^{1/p}$ defines a
semi-norm on $\mathcal{C}(A)$; the space $L^p(A;I)$ is the completion of
$\mathcal{C}(A)$ with respect to it, modulo the functions of zero
semi-norm. A function in $L^1(A;I)$ is called \emph{summable}, one in
$L^2(A;I)$ \emph{square summable}. The space $L^\infty(A;I)$ is the
completion of $\mathcal{C}(A)$ in the uniform norm, and its elements are
called \emph{essentially bounded}.

\begin{definition}[Absolute continuity]\label{ch3:def:abscount}
Let $I$ and $J$ be integrals on the simplicial sequence $A$. One says that
$I$ is \emph{absolutely continuous} with respect to $J$ if there are
simplicial representations $I=(I_m)$, $J=(J_m)$ and a $J$--summable
function $g=(g_m)$ such that $I_m(f_m)=J_m(g_mf_m)$ for every $I$--summable
$f=(f_m)$.
\end{definition}

\subsection*{Integrals at a Fixed Resolution}

At a fixed index $m$ nothing in Definition~\ref{ch3:def:integral} is a limit:
$\StateA$ is a finite set, $I_m$ is a positive linear functional on $\RR^{\StateA}$,
and by positivity and normalization $I_m$ is an average,
\begin{equation}\label{ch3:eq:finite-average}
I_m(f_m)=\sum_{x\in\StateA}p_m(x)f_m(x),
\end{equation}
for a unique family $p_m(x)\ge0$. We record the consequences used in this chapter, all of which are elementary facts about finite sums; the passage $m\to\infty$ in \eqref{ch3:eq:finite-average} is the subject of Chapter~5, where it is controlled by estimates on the kernels rather than by abstract measure theory.

\begin{lemma}\label{ch3:lem:finite-facts}
Fix $m$ and let $p_m$ be as in \eqref{ch3:eq:finite-average}.
\begin{enumerate}
\item[(i)] If $f\ge0$ and $I_m(f)=0$ then $f(x)=0$ for every $x$ with
$p_m(x)>0$.
\item[(ii)] If $K:\StateA\times\StateA\to\RR$ satisfies $K(x,y)\ge0$ and
$\sum_yK(x,y)=1$ for all $x$, then $y\mapsto\sum_xp_m(x)K(x,y)$ is again a
probability on $\StateA$.
\item[(iii)] Cauchy--Schwarz: $I_m(fg)^2\le I_m(f^2)I_m(g^2)$.
\end{enumerate}
\end{lemma}

\begin{proof}
(i) Every term $p_m(x)f(x)$ of the sum is non negative, so all vanish.
(ii) Positivity and the normalization $\sum_yp_m(x)K(x,y)=\sum_xp_m(x)=1$
follow by exchanging the two finite sums. (iii) The quadratic
$\lambda\mapsto I_m((f+\lambda g)^2)$ is non negative for all real $\lambda$
and has non positive discriminant.
\end{proof}

\section{Path Functionals}\label{ch3:sec:pathfunc}

To speak of a process one has to specify a notion of measure on a set of
paths. In the space direction we use the lattices of
\S\ref{ch3:sec:simplicial}; in the time direction we keep the coordinate
continuous, because all time integrals below are Riemann integrals of
piecewise smooth functions, and an explicit time discretization would only
add notation.

Let $A=(\StateA)$ be a bounded simplicial sequence and let $[T,T']$ be a fixed
finite time interval. Set $\Lambda_m=\StateA\times[T,T']$.

\begin{definition}[Lattice paths]\label{ch3:def:paths}
For $q\ge0$ let $\HH^q_m=\HH^q_m(T,T')$ denote the set of functions
$y_\cdot:[T,T']\to\StateA$ which are constant on each interval of a family of
$q+1$ consecutive sub intervals $[t_{k_i},t_{k_{i+1}})$, $i=0,\dots,q$, with
$t_{k_0}=T$ and $t_{k_{q+1}}=T'$. Such a path is described by the sequence
of times $T=t_0\le t_1\le\dots\le t_{q+1}=T'$ at which it jumps and the
values $y_0,\dots,y_q\in\StateA$ it takes on the successive intervals. We write
$\HH_m=\bigcup_{q\ge0}\HH^q_m$ for the set of paths with finitely many
jumps.
\end{definition}

Each $\HH^q_m$ is a finite union of finite dimensional manifolds with
boundary. We adopt the convention that the times $t_1=\dots=t_q$ are
allowed, so that a path with fewer than $q$ jumps is regarded as an element
of $\HH^q_m$ as well; with this convention the sets $\HH^q_m$ are nested in
$q$ and $\HH_m$ is their union.

\begin{definition}[Function algebras]\label{ch3:def:algebra}
$\mathcal{C}^q(\Lambda)$ denotes the $C^*$ algebra of uniformly continuous
sequences of simplicial functions $F=(F_m)$ with $F_m:\HH^q_m\to\RR$,
endowed with the operations of sum, product and the uniform norm
\begin{equation}
\norm{F}_\infty=\lim_{m\to\infty}\sup_{y_\cdot\in\HH^q_m}\abs{F_m(y_\cdot)} .
\end{equation}
\end{definition}

\begin{definition}[Path functionals]\label{ch3:def:pathfunctional}
A \emph{path functional} is a sequence $F=(F^q)$, $q=0,1,\dots$, of
continuous simplicial functions $F^q=(F^q_m)\in\mathcal{C}^q(\Lambda)$
satisfying the compatibility condition
\begin{equation}
F^{q'}_m(y_\cdot)=F^q_m(y_\cdot)
\qquad\text{for all } y_\cdot\in\HH^q_m \text{ and all } q'>q .
\end{equation}
\end{definition}

\begin{definition}[Non anticipatory functionals]\label{ch3:def:nonanticipatory}
A one parameter family $F(y_\cdot,t)=(F^q(y_\cdot,t))$, $t\in[T,T']$, of
path functionals is \emph{non anticipatory} if
\begin{equation}
F^q_m(y_\cdot,t)=F^q_m(y'_\cdot,t)
\qquad\text{whenever } y_s=y'_s \text{ for all } s\le t .
\end{equation}
\end{definition}

Intuitively, a non anticipatory functional is indifferent to information
about the path after time $t$. Two classes of examples suffice for the whole of the chapter.

\begin{example}[State functions]\label{ch3:ex:statefunc}
Let $F_0:\Lambda\to\RR$ be a piecewise smooth function of its two arguments.
Then
\begin{equation}\label{ch3:eq:defft}
F^q_m(y_\cdot,t)=F_0(y_t,t)
\end{equation}
defines a non anticipatory path functional. We adopt the \emph{cadlag}
convention: right continuous with left limits, so that $y_t$ denotes the
value of the path after a possible jump at $t$.
\end{example}

\begin{example}[Iterated integrals]\label{ch3:ex:iterated}
Let $F_1,\dots,F_k:\Lambda\to\RR$ and $a(t;s_1,\dots,s_k)$ be piecewise
smooth, the latter jointly in all its arguments, and set
\begin{equation}\label{ch3:eq:defft2}
F^q_m(y_\cdot,t)=\int_T^t\dd s_1\cdots\int_{s_{k-1}}^t\dd s_k\,
a(t;s_1,\dots,s_k)F_1(y_{s_1},s_1)\cdots F_k(y_{s_k},s_k).
\end{equation}
This is the prototype of a path dependent functional, and its expectation
is what the theory of stochastic integrals computes. The regularity
assumptions may be relaxed by taking limits, but they are satisfied in every
application we consider.
\end{example}

The results of this chapter require only the following regularity, which we impose once and for all.

\begin{assumption}[Regular functionals]\label{ch3:ass:regular}
Every path functional considered in \S\S\ref{ch3:sec:markov}--\ref{ch3:sec:ftf}
is a finite linear combination of products of functionals of the form
\eqref{ch3:eq:defft} and \eqref{ch3:eq:defft2}, with coefficients that are piecewise
smooth in all time variables. In particular, for fixed $y_\cdot$ the map
$t\mapsto F(y_\cdot,t)$ is piecewise continuously differentiable, and the
one sided derivative at a point of discontinuity is finite.
\end{assumption}

\section{Filtered Probability Spaces and n-Point Functions}\label{ch3:sec:npoint}

\begin{definition}[Filtered probability spaces]\label{ch3:def:fps}
Let $\Pmeas$ be the $C^*$ algebra generated by the functionals of the form
\eqref{ch3:eq:defft} and \eqref{ch3:eq:defft2}, by linear combinations of finite
products and completion in norm. A \emph{filtered probability space on the
lattice} $\Lambda$ is a bounded linear functional on $\Pmeas$ of total mass
one.
\end{definition}

A constructive description of a process is obtained by prescribing the joint
distributions of the values it takes at finitely many times.

\begin{definition}[$n$-point functions]\label{ch3:def:npoint}
Let $\Phi_m(y_1,t_1;\dots;y_n,t_n\mid y_0)$, $m\ge m_0$, be defined for
$y_1,\dots,y_n\in\StateA$ and $T\le t_1\le\dots\le t_n\le T'$, piecewise
differentiable in the time coordinates, and suppose that
\begin{align}
&(\mathrm{SP1})\qquad \Phi_m(y_1,t_1;\dots;y_n,t_n\mid y_0)\ge0
\quad\forall y_0,\dots,y_n\in\StateA,\label{ch3:eq:SP1}\\
&(\mathrm{SP2})\qquad \sum_{y_1,\dots,y_n\in\StateA}
\Phi_m(y_1,t_1;\dots;y_n,t_n\mid y_0)=1
\quad\forall y_0\in\StateA,\label{ch3:eq:SP2}\\
&(\mathrm{SP3})\qquad \Phi_m(\dots;y_i,t_i;y_{i+1},t_{i+1};\dots\mid y_0)=0
\ \text{ if } t_i=t_{i+1},\ y_i\neq y_{i+1},\label{ch3:eq:SP3}\\
&(\mathrm{SP4})\qquad \sum_{y_i\in\StateA}\Phi_m(\dots;y_{i-1},t_{i-1};
y_i,t_i;y_{i+1},t_{i+1};\dots\mid y_0)\notag\\
&\qquad\qquad\qquad
=\Phi_m(\dots;y_{i-1},t_{i-1};y_{i+1},t_{i+1};\dots\mid y_0).
\label{ch3:eq:SP4}
\end{align}
\end{definition}

For fixed $y_0$ and fixed times $t_1\le\dots\le t_n$, \eqref{ch3:eq:SP1} and
\eqref{ch3:eq:SP2} say that $\Phi_m(\cdot\mid y_0)$ is a probability
distribution on $\StateA^{\{1,\dots,n\}}$; \eqref{ch3:eq:SP3} is the statement that
a path cannot be in two distinct sites at the same time; \eqref{ch3:eq:SP4} is
the consistency of the finite dimensional marginals. The family is
therefore, by the usual argument of Kolmogorov, the same thing as a
probability measure on the path space $\HH_m$.

\begin{definition}[Stochastic processes]\label{ch3:def:process}
An \emph{adapted process} is given by a non anticipatory path functional
$F$ and a filtered probability space determined by a family of $n$-point
functions $(\Phi_m)_{n\ge1}$.
\end{definition}

For a state function \eqref{ch3:eq:defft} the conditional expectation given the
initial value is
\begin{equation}
\EE\big[F_t\mid y_T=y_0\big]=\sum_{y_1\in\StateA}\Phi_m(y_1,t\mid y_0)F(y_1,t),
\end{equation}
while for the iterated integral \eqref{ch3:eq:defft2}
\begin{equation}
\begin{aligned}
\EE\big[F_t\mid y_T=y_0\big]={}&\sum_{y_1,\dots,y_n}\int_T^t\dd s_1\cdots
\int_{s_{n-1}}^t\dd s_n\,\Phi_m(y_1,s_1;\dots;y_n,s_n\mid y_0)\\
&\qquad\times a(t;s_1,\dots,s_n)\prod_{j=1}^nF_j(y_j,s_j).
\end{aligned}
\end{equation}
Higher moments reduce to sums over time orderings. For instance, if $F_t$
is the functional \eqref{ch3:eq:defft2} with $k=1$, one finds
\begin{equation}
\begin{aligned}
\EE\big[F_t^2\mid y_T=y_0\big]={}&2\sum_{y_1,y_2}\int_T^t\dd s_1\int_{s_1}^t
\dd s_2\,\Phi_m(y_1,s_1;y_2,s_2\mid y_0)\\
&\qquad\times a(t;s_1)a(t;s_2)F_1(y_1,s_1)F_2(y_2,s_2),
\end{aligned}
\end{equation}
the factor $2$ compensating for the time ordering needed to write the square
as an integral over the simplex $T\le s_1\le s_2\le t$. These formulas are
used repeatedly in the applications to path dependent options developed in Chapter~6.

\subsection*{Conditioning and Equivalent Measures}

Conditioning on a past history is the restriction of the path measure to the
set of paths which agree with a given one up to time $t$. In a Markovian
framework, which is the one adopted from \S\ref{ch3:sec:markov} onwards, this
restriction is again described by a transition kernel at the conditioning
time and no memory effects survive, so the formulas below are all that is
needed. We record the general definition anyway, since it is the one that
makes the phrase ``equivalent measure'' precise.

\begin{definition}[Radon--Nikodym derivative]\label{ch3:def:rn}
Let $\Phi^1_m$ and $\Phi^2_m$ be two families of $n$-point functions on
$\StateA$. The Radon--Nikodym derivative of $\Phi^1$ with respect to $\Phi^2$ is
the path functional
\begin{equation}
\rho(y_\cdot)=\lim_{n\to\infty}
\frac{\Phi^1_m(y_{t_1},t_1;\dots;y_{t_n},t_n\mid y_0)}
{\Phi^2_m(y_{t_1},t_1;\dots;y_{t_n},t_n\mid y_0)},
\qquad t_i=T+\tfrac{i}{n}(T'-T),
\end{equation}
whenever the limit exists. The measure $\Phi^1$ is absolutely continuous
with respect to $\Phi^2$ if $\rho$ is finite and summable, and the two are
equivalent if in addition $\Phi^2$ is absolutely continuous with respect to
$\Phi^1$.
\end{definition}

For Markovian families, which are the only ones occurring below, the
Radon--Nikodym derivative is a product of one step ratios and the condition
reduces to the strict positivity of a family of densities on the attainable
sets; this is the content of \S\ref{ch3:sec:measure}.

\begin{definition}[Stopping times]\label{ch3:def:stopping}
A \emph{stopping time} is an adapted process $\tau_t=\tau(y_\cdot,t)$ taking
the two values $0$ and $1$. If $F_t=F(y_\cdot,t)$ is an adapted process and
$R_t=R(y_\cdot,t)$ is another, the \emph{stopped version} of $F$ at $\tau$ is
the adapted process $\bar F$ with $\bar F(y_\cdot,t)=F(y_\cdot,t)$ when
$\tau_t=0$ and $\bar F(y_\cdot,t)=R(y_\cdot,t)$ when $\tau_t=1$.
\end{definition}

Stopping times are used in \S\ref{ch3:sec:ftf} to liquidate a portfolio at the first jump of the underlying, and in Chapter~6 to exercise callable structures.

\section{Markov Processes}\label{ch3:sec:markov}

Throughout this section $m\ge m_0$ is fixed and $\StateA$ is a finite set of
sites. We write $\Lm(y,z;t)$ for the entries of a matrix indexed by
$\StateA\times\StateA$, and we use the same symbol for the matrix itself. The
transpose convention is that functions of the site are column vectors, so
that $(\Lm(t)f)(y)=\sum_z\Lm(y,z;t)f(z)$.

\subsection*{Propagators}

\begin{definition}[Markov propagator]\label{ch3:def:propagator}
A \emph{Markov propagator} on $\StateA$ is a family of functions
$U_m(y,s;z,t)$, defined for $y,z\in\StateA$ and $T\le s\le t\le T'$, such that
\begin{align}
&(\mathrm{CK1})\qquad U_m(y,s;z,t)\ge0,\label{ch3:eq:CK1}\\
&(\mathrm{CK2})\qquad U_m(y,s;z,s)=\delta_{yz},\label{ch3:eq:CK2}\\
&(\mathrm{CK3})\qquad \sum_{z}U_m(y,s;z,u)U_m(z,u;w,t)=U_m(y,s;w,t)
\qquad\forall\, s\le u\le t .\label{ch3:eq:CK3}
\end{align}
\end{definition}

By \eqref{ch3:eq:CK1} and \eqref{ch3:eq:CK3} with $u=t$ the row sums equal one, so
each $U_m(y,s;\cdot,t)$ is a probability on $\StateA$.

\begin{proposition}\label{ch3:prop:npoint-from-propagator}
Let $U_m$ be a Markov propagator and set, for $n\ge1$,
\begin{equation}\label{ch3:eq:npoint-propagator}
\Phi_m(y_1,t_1;\dots;y_n,t_n\mid y_0)
=\prod_{j=1}^{n}U_m(y_{j-1},t_{j-1};y_j,t_j),
\qquad t_0=T .
\end{equation}
Then $\Phi_m$ satisfies $(\mathrm{SP1})$--$(\mathrm{SP4})$. Conversely, if
$\Phi_m$ is a family of $n$-point functions satisfying
$(\mathrm{SP1})$--$(\mathrm{SP4})$ and the Markov property
\begin{equation}\label{ch3:eq:markovprop}
\Phi_m(y_1,t_1;\dots;y_n,t_n\mid y_0)
=\Phi_m(y_1,t_1;\dots;y_k,t_k\mid y_0)\,
\Phi_m(y_{k+1},t_{k+1};\dots;y_n,t_n\mid y_k)
\end{equation}
for $k=1,\dots,n-1$, then it is of the form
\eqref{ch3:eq:npoint-propagator} with
$U_m(y,s;z,t)=\Phi_m(z,t\mid y)$ at $s=T$ and, for general $s$, the ratio
$\Phi_m(z,t\mid y)/\Phi_m(y,s\mid y)$ on the set where the denominator does
not vanish.
\end{proposition}

\begin{proof}
$(\mathrm{SP1})$ is $(\mathrm{CK1})$. $(\mathrm{SP2})$ is proved by
induction on $n$: for $n=1$ it is the unit row sum, and the induction step
uses $(\mathrm{CK3})$ in the form $\sum_{y_n}U_m(y_{n-1},t_{n-1};y_n,t_n)=1$.
$(\mathrm{SP3})$ is $(\mathrm{CK2})$. For $(\mathrm{SP4})$ one sums over
$y_i$: the product splits as
$\prod_{j<i}U_m(y_{j-1},t_{j-1};y_j,t_j)$ times
$\sum_{y_i}U_m(y_{i-1},t_{i-1};y_i,t_i)U_m(y_i,t_i;y_{i+1},t_{i+1})$, which
by $(\mathrm{CK3})$ equals $U_m(y_{i-1},t_{i-1};y_{i+1},t_{i+1})$, times the
remaining factors. The converse is a direct induction using
\eqref{ch3:eq:markovprop}.
\end{proof}

\begin{definition}[Markov process]\label{ch3:def:markovprocess}
A filtered probability space whose $n$-point functions are of the form
\eqref{ch3:eq:npoint-propagator} is called \emph{Markovian}, and the associated
adapted processes are those of a \emph{Markov process}.
\end{definition}

\subsection*{Generators}

\begin{definition}[Markov generator]\label{ch3:def:generator}
Assume that $t\mapsto U_m(y,s;z,t)$ is differentiable from the right at
$t=s$ for all $y,z$. The \emph{Markov generator} at time $s$ is the matrix
\begin{equation}\label{ch3:eq_gen}
\Lm(y,z;s)=\lim_{t\downarrow s}\frac{\dd}{\dd t}U_m(y,s;z,t).
\end{equation}
\end{definition}

\begin{proposition}\label{ch3:prop:generator-properties}
If $\Lm(y,z;t)$ is a Markov generator then
\begin{align}
&(\mathrm{MG1})\qquad \Lm(y,z;t)\ge0 \quad\text{for } y\neq z,\\
&(\mathrm{MG2})\qquad \Lm(y,y;t)=-\sum_{z\neq y}\Lm(y,z;t).
\end{align}
Conversely, let $t\mapsto\Lm(y,z;t)$ be piecewise continuous and bounded on
$[T,T']$ and satisfy $(\mathrm{MG1})$--$(\mathrm{MG2})$. Then the
differential equation \eqref{ch3:eq_gen}, together with the initial condition
$U_m(y,t;z,t)=\delta_{yz}$, has one and only one solution, it is a Markov
propagator in the sense of Definition~\ref{ch3:def:propagator}, and it depends
continuously on the data.
\end{proposition}

\begin{proof}
$(\mathrm{MG1})$ follows from $(\mathrm{CK1})$ after differentiating, and
$(\mathrm{MG2})$ from differentiating the unit row sum. For the converse, let
$L(t)$ denote the matrix $\Lm(\cdot,\cdot;t)$ and let
$M=\sup_t\norm{L(t)}$ where $\norm{\cdot}$ is the maximum row sum norm. For
$N$ so large that $\delta t=(t-s)/N$ satisfies $\delta t\,M\le1$ the matrix
$I+\delta tL(u)$ has non negative entries for every $u$: the off diagonal
entries are $\delta tL(y,z;u)\ge0$ and the diagonal entry is
$1-\delta t\sum_{z\neq y}L(y,z;u)\ge1-\delta tM\ge0$. Its row sums equal
one. Hence every Euler product
\begin{equation}\label{ch3:eq_pexp1}
U^{N}_m(y,s;z,t)=\Big[\big(I+\delta tL(u_N)\big)\cdots
\big(I+\delta tL(u_1)\big)\Big](y,z),\qquad u_j=s+j\,\delta t,
\end{equation}
is a stochastic matrix, and so is any limit of such products because the set
of stochastic matrices is closed.

It remains to identify the limit.

\emph{The exact propagator.} Define $U^{(0)}(s,t)=I$ and
\begin{equation}\label{ch3:eq:picard}
U^{(k+1)}(s,t)=I+\int_s^tU^{(k)}(s,u)L(u)\dd u .
\end{equation}
By induction $\norm{U^{(k+1)}-U^{(k)}}\le M^{k+1}(t-s)^{k+1}/(k+1)!$ uniformly
in $s\le t$, so the sequence converges uniformly to a limit $U$ which is
continuous and satisfies the integral equation; it therefore satisfies
\eqref{ch3:eq_gen} by the fundamental theorem of calculus. If $U$ and $\tilde U$
are two solutions, their difference satisfies
$\norm{\Delta(s,t)}\le M\int_s^t\norm{\Delta(s,u)}\dd u$, whence $\Delta=0$ by
Gr\"onwall's inequality. This proves existence and uniqueness of the
solution.

\emph{The Euler products converge to it.} Write
$R_N(s,t)=\prod_{j=1}^{N}(I+\delta tL(u_j))$ with $\delta t=(t-s)/N$, so that
$R_N(s,t)=R_{N-1}(s,t-\delta t)(I+\delta tL(t-\delta t))$, while by the
integral equation and the mean value theorem for the integral
\begin{equation}
\begin{aligned}
U(s,t)&=U(s,t-\delta t)+\int_{t-\delta t}^{t}U(s,u)L(u)\dd u
=U(s,t-\delta t)\big(I+\delta tL(t-\delta t)\big)+E,\\
\norm{E}&\le\delta t^{2}M^{2}e^{(t-s)M}.
\end{aligned}
\end{equation}
For $\Delta=R_N-U$ this gives
$\Delta(s,t)=\Delta(s,t-\delta t)(I+\delta tL(t-\delta t))-E$, so that
\begin{equation}
\norm{\Delta(s,t)}\le\norm{\Delta(s,t-\delta t)}(1+\delta tM)
+\delta t^{2}M^{2}e^{(t-s)M}.
\end{equation}
Iterating from $t$ down to $s$, and using $\Delta(s,s)=0$, yields
\begin{equation}
\norm{R_N(s,t)-U(s,t)}\le\delta t\,M\big(e^{(t-s)M}-1\big)e^{(t-s)M},
\end{equation}
which tends to zero with $\delta t$. Hence the Euler products converge to the
unique solution.

\emph{The limit is a propagator.} Each $R_N$ is a stochastic matrix and the
set of stochastic matrices is closed, so $U(s,t)$ has non negative entries and
unit row sums: $(\mathrm{CK1})$ and, for $s=t$, $(\mathrm{CK2})$. Finally
$(\mathrm{CK3})$ holds because for fixed $s\le u$ both
$t\mapsto\sum_zU(s,u)(y,z)U(u,t)(z,w)$ and $t\mapsto U(s,t)(y,w)$ solve the
same forward equation in $t$ with the same value at $t=u$, hence coincide by
uniqueness.
\end{proof}

\begin{remark}\label{ch3:rem:ck3}
The proof of $(\mathrm{CK3})$ is worth spelling out: for fixed $s\le u$ the
map $(y,t)\mapsto\sum_zU_m(y,s;z,u)U_m(z,u;y',t)$ and the map
$t\mapsto U_m(y,s;y',t)$ satisfy the same linear equation in $t$ with the
same initial data at $t=u$, so they coincide by uniqueness.
\end{remark}

\subsection*{The Dyson Expansion and the Path Integral}

\begin{theorem}[Dyson expansion]\label{ch3:thm:dyson}
Under the hypotheses of Proposition~\ref{ch3:prop:generator-properties},
\begin{equation}\label{ch3:ma_dyson}
U_m(y,s;z,t)=\Big(I+\sum_{n=1}^{\infty}\int_s^t\dd u_1\int_{u_1}^t\dd u_2
\cdots\int_{u_{n-1}}^t\dd u_n\,L(u_1)\cdots L(u_n)\Big)(y,z),
\end{equation}
the series converging uniformly and absolutely on $T\le s\le t\le T'$.
\end{theorem}

\begin{proof}
Expand the Euler product of \eqref{ch3:eq_pexp1} and pass to the limit; the
$n$-th term is bounded by $(t-s)^nM^n/n!$, so the series converges
uniformly and the interchange of limits is justified by the Weierstrass
test.
\end{proof}

\begin{theorem}[Forward and backward equations]\label{ch3:thm:forwardbackward}
The propagator satisfies
\begin{equation}\label{ch3:backward_eq}
\frac{\partial}{\partial s}U_m(y,s;z,t)+\big(\Lm(s)U_m(s,t)\big)(y,z)=0
\end{equation}
and
\begin{equation}\label{ch3:forward_eq}
\frac{\partial}{\partial t}U_m(y,s;z,t)=\big(U_m(s,t)\Lm(t)\big)(y,z).
\end{equation}
\end{theorem}

\begin{proof}
Differentiate \eqref{ch3:ma_dyson} with respect to the two time coordinates, or
equivalently differentiate the integral equation in the proof of
Proposition~\ref{ch3:prop:generator-properties}.
\end{proof}

\begin{definition}[Path ordered exponential]\label{ch3:def:pexp}
The identity \eqref{ch3:ma_dyson} is written
\begin{equation}
U_m(y,s;z,t)=\mathrm{Pexp}\Big(\int_s^t\Lm(u)\dd u\Big)(y,z),
\end{equation}
where the time ordering operator $\mathrm{P}$ acts by
\begin{equation}
\mathrm{P}\Big(\int_s^t\Lm(u)\dd u\Big)^n=n!\int_s^t\dd u_1\int_{u_1}^t
\dd u_2\cdots\int_{u_{n-1}}^t\dd u_n\,\Lm(u_1)\cdots\Lm(u_n).
\end{equation}
\end{definition}

\begin{definition}[Symbolic path]\label{ch3:def:symbolicpath}
A \emph{symbolic path} $\gamma=\{\gamma_0,\gamma_1,\dots\}$ is a sequence of
sites with $\gamma_j\neq\gamma_{j-1}$ for all $j\ge1$. We write $\Gamma_m$
for the set of symbolic paths in $\StateA$.
\end{definition}

\begin{theorem}[Path integral representation]\label{ch3:thm:pathintegral}
The propagator admits the representation
\begin{equation}\label{ch3:eq_upathint}
\begin{aligned}
U_m(x,T;y,T')={}&\sum_{q=1}^{\infty}\ \sum_{\substack{\gamma\in\Gamma_m:\\ \gamma_0=x,\
\gamma_q=y}}
\int_T^{T'}\dd s_1\int_{s_1}^{T'}\dd s_2\cdots\int_{s_{q-1}}^{T'}\dd s_q\;
e^{\int_T^{s_1}\Lm(\gamma_0,\gamma_0;u)\dd u}\\
&\qquad\times\prod_{j=1}^{q}\Lm(\gamma_{j-1},\gamma_j;s_j)\,
e^{\int_{s_j}^{s_{j+1}}\Lm(\gamma_j,\gamma_j;u)\dd u},
\end{aligned}
\end{equation}
with $s_{q+1}=T'$.
\end{theorem}

\begin{proof}
Expand the $n$-th term of \eqref{ch3:ma_dyson} according to the sequence of
sites visited, grouping consecutive equal sites, and resum the diagonal
runs.
\end{proof}

The representation \eqref{ch3:eq_upathint} is the operator version of the
Feynman--Kac formula and is used in both Chapters~5 and~6. It also shows that the
process has finitely many jumps on a compact time interval with probability
one: if $\lambda_m=\max_x\abs{\Lm(x,x)}$ and
$c=\sup_t\norm{\Lm(t)}$, the mass of the set of paths with exactly $q$
jumps is at most
\begin{equation}\label{ch3:eq_mest}
\frac{c^q(T'-T)^q}{q!}e^{-\lambda_m(T'-T)},
\end{equation}
whose sum over $q$ converges. In particular the first jump time of the
process is well defined almost surely, a fact used in the proof of
Proposition~\ref{ch3:prop:localNA}.

\subsection*{Piecewise Smooth Processes, Attainable Sets, Reachability}

Applications, and the theory of pricing in particular, require generators
that are discontinuous in time, because coupon payments, fixings and
corporate actions occur at isolated dates.

\begin{definition}[Piecewise smooth Markov process]\label{ch3:def:piecewise}
Let $T=t_0<\dots<t_n=T'$ be finitely many times. A \emph{piecewise smooth
Markov process} is given by a family of generators $\Lm^{i}(y,z;t)$ defined
and satisfying $(\mathrm{MG1})$--$(\mathrm{MG2})$ on each half open interval
$[t_i,t_{i+1})$, and, for each $i=1,\dots,n-1$, a mapping operator
$U^i_m(y,z)$ with
\begin{align}
&(\mathrm{MA1})\qquad U^i_m(y,z)\ge0,\qquad
&(\mathrm{MA2})\qquad \sum_{z}U^i_m(y,z)=1\ \ \forall y .
\end{align}
For $t_i\le s<t<t_{i+1}$ the propagator is
$U_m(y,s;z,t)=\mathrm{Pexp}(\int_s^t\Lm^i(u)\dd u)(y,z)$, and
\begin{equation}
U_m(y,s;z,t_{i+1})=\sum_w\mathrm{Pexp}\Big(\int_s^{t_{i+1}}\Lm^i(u)\dd u\Big)(y,w)U^i_m(w,z).
\end{equation}
More general propagators are obtained by composing these over successive
intervals.
\end{definition}

Time is thus continuous inside each interval and the model is genuinely
discrete at the dates $t_i$. Everything proved in this chapter applies to
each interval separately, and the dates $t_i$ are treated as one step models
in \S\ref{ch3:sec:ftf}.

\begin{definition}[Attainable sets]\label{ch3:def:attainable}
Let $U_m$ be a piecewise smooth Markov process, $y\in\StateA$ and
$t\in[T,T']$. The \emph{attainable set} $\Attain(U,t,y)\subseteq\StateA$ is defined
as follows: if $t\in(t_i,t_{i+1})$ for some $i$, then $\Attain(U,t,y)$ is the
set of $\bar y\neq y$ with $\Lm^i(y,\bar y;t)>0$; if $t=t_i$ for some
$i$, then $\Attain(U,t_i,y)$ is the set of $\bar y$ with $U^i_m(y,\bar y)>0$.
\end{definition}

\begin{definition}[Reachability]\label{ch3:def:reachable}
Let $P$ be the path measure of a piecewise smooth Markov process, and let
$y_0\in\StateA$ be the initial site at time $T$. The pair $(y,t)$ is
\emph{reachable} if $P(y_t=y)>0$.
\end{definition}

\begin{lemma}\label{ch3:lem:reach-open}
If $(y,t)$ is reachable and $t$ is not one of the dates $t_i$, then
$(y,s)$ is reachable for all $s$ in a neighbourhood of $t$.
\end{lemma}

\begin{proof}
The function $s\mapsto P(y_s=y)=U_m(y_0,T;y,s)$ is continuous on each
interval $(t_i,t_{i+1})$ by Proposition~\ref{ch3:prop:generator-properties} and
strictly positive at $t$ by hypothesis; hence it is bounded below by
$P(y_t=y)/2$ in a neighbourhood.
\end{proof}

\begin{definition}[Equivalent processes]\label{ch3:def:equivprocess}
Two piecewise smooth propagators $U$ and $U'$ are \emph{equivalent} if their
attainable sets coincide for all $t$ and $y$. If $\Attain(U,t,y)\subseteq
\Attain(U',t,y)$ for all $t,y$ one says that $U$ is absolutely continuous with
respect to $U'$.
\end{definition}

\subsection*{Symbols, Monotonic Processes, Diffusions}

\begin{definition}[Inverse lattice]\label{ch3:def:brillouin}
The lattice
\begin{equation}\label{ch3:brillouin}
B_m=\left\{-\frac{2^{m-1}\pi}{L}+\frac{k\pi}{L},\ k=0,\dots,2^m-1\right\}
\end{equation}
is called the \emph{Brillouin zone}, or inverse lattice, of $\StateA$.
\end{definition}

\begin{definition}[Symbol]\label{ch3:def_symbol}
The \emph{symbol} of a Markov generator is
\begin{equation}
\hat\Lm(x,p;t)=\sum_{y\in\StateA}\Lm(x,y;t)e^{ip(y-x)},\qquad p\in B_m .
\end{equation}
\end{definition}

The symbol is the operator eigenvalue of the generator on the exponential
functions $x\mapsto e^{ipx}$ when the generator is translation invariant, and
it is the object through which the convergence estimates of Chapter~5 are formulated.

\begin{definition}[Monotonic processes]\label{ch3:def_monotonic}
A Markov process with generator $\Lm(x,y)$ on the simplicial sequence $\StateA$
is \emph{monotonic non decreasing} if $\Lm(x,y)=0$ whenever $y<x$, and
\emph{monotonic non increasing} if $\Lm(x,y)=0$ whenever $y>x$.
\end{definition}

\begin{definition}[Diffusion processes]\label{ch3:def:diffusion}
A \emph{diffusion process} is a Markov process with generator
\begin{equation}
\Lm(x,y;t)=\mu_m(x;t)\grad_{h_m}(x,y)+\frac{\sigma_m(x;t)^2}{2}\lap_{h_m}(x,y),
\end{equation}
where
\begin{equation}\label{ch3:eq_nabladelta}
\grad_h(x,y)=\frac{\delta_{y,x+h}-\delta_{y,x-h}}{2h},\qquad
\lap_h(x,y)=\frac{\delta_{y,x+h}+\delta_{y,x-h}-2\delta_{y,x}}{h^2},
\end{equation}
and $\mu_m,\sigma_m$ are smooth simplicial functions.
\end{definition}

\begin{proposition}\label{ch3:prop:diffusion-symbol}
\begin{itemize}
\item[(i)] The symbol of a diffusion is
\begin{equation}\label{ch3:eq:diffusion-symbol}
\hat\Lm(x,p;t)=\mu_m(x;t)\frac{\sin ph}{ih}
+\sigma_m(x;t)^2\frac{\cos ph-1}{h^2},
\end{equation}
the value of the generator on the exponential $e^{-ipx}$, which is the convention used in the Fourier analysis of Chapters~5 and~6.
\item[(ii)] Its path integral representation reads
\begin{equation}\label{ch3:eq_udiff1}
U_m(x,T;y,T')=\sum_{q=1}^{\infty}2^{-q}
\sum_{\substack{\gamma\in\Gamma_m:\ \gamma_0=x,\ \gamma_q=y\\ \abs{\gamma_j-\gamma_{j-1}}=1\ \forall j\ge1}}
W_m(\gamma,q;T,T'),
\end{equation}
where
\begin{align}\label{ch3:eq_wdiff}
W_m(\gamma,q;T,T')&=\int_T^{T'}\dd s_1\cdots\int_{s_{q-1}}^{T'}\dd s_q\;
e^{\int_{s_q}^{T'}\Lm(y,y;v)\dd v}\notag\\
&\qquad\times\prod_{j=0}^{q-1}\Big(e^{\int_{s_j}^{s_{j+1}}\Lm(\gamma_j,\gamma_j;v)\dd v}
2\Lm(\gamma_j,\gamma_{j+1};s_{j+1})\Big),
\qquad s_0=T .
\end{align}
\end{itemize}
\end{proposition}

\begin{proof}
(i) One computes $\sum_y\grad_h(x,y)e^{-ip(y-x)}=-i\frac{\sin ph}{h}$ and
$\sum_y\lap_h(x,y)e^{-ip(y-x)}=2\frac{\cos ph-1}{h^2}$, and the two terms of
the generator carry the coefficients $\mu$ and $\sigma^2/2$.

(ii) In \eqref{ch3:eq_upathint} each jump weight
$\Lm(\gamma_{j-1},\gamma_j;s_j)$ is split as
$2\Lm(\gamma_{j-1},\gamma_j;s_j)/2$: the $q$ powers of $2$ in the
denominators collect into the factor $2^{-q}$, and the numerators are the
factors displayed in \eqref{ch3:eq_wdiff}. Since a diffusion has non zero off
diagonal entries only in the two neighbouring sites, the symbolic paths
occurring in the sum are precisely those with $\abs{\gamma_j-\gamma_{j-1}}=1$.
\end{proof}

\section{Measure Changes}\label{ch3:sec:measure}

A change of measure is a reweighting of the one step transition probabilities
site by site and time by time. It is the constructive counterpart of
Girsanov's theorem, and in the present finite dimensional setting it is an
explicit algebraic operation whose properties are easily verified.

\begin{definition}[Measure changes]\label{ch3:def_measure_change}
Let $U_m$, $m\ge m_0$, be a family of piecewise smooth Markov propagators. A
\emph{measure change} is a family of strictly positive functions
$G_m^{yt}(y')>0$, indexed by $y\in\StateA$ and $t\in[T,T']$ and defined for
$y'\in\Attain(U,t,y)$, possibly extended by zero outside. It defines a
transformed Markov generator by the formula
\begin{equation}\label{ch3:eq:measure-change}
\Lm'(y,y';t)=\frac{1}{G_m^{yt}(y)}\Lm(y,y';t)G_m^{yt}(y')
-\frac{1}{G_m^{yt}(y)}\big(\Lm(t)G_m^{yt}\big)(y)\,\delta_{yy'},
\end{equation}
where $(\Lm(t)G_m^{yt})(y)=\sum_{z}\Lm(y,z;t)G_m^{yt}(z)$. At a date $t_i$
the definition is the corresponding one: a measure change is a family of
strictly positive densities $G_m^{yt_i}(z)>0$, $z$ in the support of
$U^i_m(y,\cdot)$, and the transformed mapping operator is
\begin{equation}\label{ch3:eq:measure-change-date}
U'^i_m(y,z)=\frac{1}{C_y}U^i_m(y,z)G_m^{yt_i}(z),
\qquad C_y=\sum_{z'}U^i_m(y,z')G_m^{yt_i}(z'),
\end{equation}
which satisfies $(\mathrm{MA1})$ and $(\mathrm{MA2})$ because $C_y>0$ is the
normalizing constant.
\end{definition}

\begin{remark}\label{ch3:rem:change-at-y}
The value of $G_m^{yt}$ at $y'=y$ does not affect the transformation
\eqref{ch3:eq:measure-change}: it cancels between the two terms. This is why the
definition only prescribes $G$ on the attainable set, where it is required to
be positive.
\end{remark}

\begin{proposition}\label{ch3:prop:measure-change-generator}
Let $G_m^{yt}$ be a measure change and let $\Lm'$ be given by
\eqref{ch3:eq:measure-change}. Then $\Lm'$ satisfies $(\mathrm{MG1})$ and
$(\mathrm{MG2})$, so it is the generator of a piecewise smooth Markov
process. Moreover, for $y'\neq y$
\begin{equation}\label{ch3:eq:reweight}
\Lm'(y,y';t)=\frac{G_m^{yt}(y')}{G_m^{yt}(y)}\Lm(y,y';t),
\end{equation}
and the transformed process has the same attainable sets as the original
one.
\end{proposition}

\begin{proof}
For $y'\neq y$ the second term of \eqref{ch3:eq:measure-change} vanishes and
\eqref{ch3:eq:reweight} holds; by Definition~\ref{ch3:def_measure_change} the ratio
of the two positive densities is positive, so $(\mathrm{MG1})$ holds. For the
row sums, using \eqref{ch3:eq:reweight} and
$(\mathrm{MG2})$ for $\Lm$,
\begin{equation}
\begin{aligned}
\Lm'(y,y)&=\Lm(y,y)-\frac{1}{G_m^{yt}(y)}\sum_z\Lm(y,z;t)G_m^{yt}(z)
=-\sum_{z\neq y}\frac{G_m^{yt}(z)}{G_m^{yt}(y)}\Lm(y,z;t)\\
&=-\sum_{z\neq y}\Lm'(y,z;t),
\end{aligned}
\end{equation}
which is $(\mathrm{MG2})$. Since the ratio in \eqref{ch3:eq:reweight} is strictly
positive, $\Lm'(y,y';t)>0$ if and only if $\Lm(y,y';t)>0$, which is the
statement about attainable sets.
\end{proof}

\begin{definition}[Time changes]\label{ch3:def:timechange}
A measure change is a \emph{time change} if there is a function
$\phi(y,t)>0$ with $G_m^{yt}(y')=\phi(y,t)$ for all $y'\in\Attain(U,t,y)$. It is
\emph{state independent} if $\phi(y,t)=\phi(t)$ depends on time only.
\end{definition}

\begin{theorem}[Deterministic time changes]\label{ch3:thm:timechange}
If $\phi(t)$ is state independent, then the transformed propagator is
\begin{equation}
U_m'(y,t;y',t')=U_m\big(y,\lambda(t);y',\lambda(t')\big),
\qquad \lambda(t)=\int_T^t\phi(s)\dd s .
\end{equation}
\end{theorem}

\begin{proof}
Insert $G_m^{yt}(y')=\phi(t)$ in \eqref{ch3:eq:reweight}: the generator becomes
$\Lm'(y,y';t)=\phi(t)\Lm(y,y';t)$ for $y'\neq y$, and the diagonal is
determined by $(\mathrm{MG2})$. The propagator of $\phi(t)\Lm(\cdot;t)$ is
the propagator of $\Lm$ reparametrised by $\lambda$; this is immediate from
the Dyson expansion \eqref{ch3:ma_dyson}, where each of the $n$ time integrals
produces a factor $\phi$ that is absorbed by the change of variable
$u\mapsto\lambda(u)$.
\end{proof}

\begin{definition}[Numeraire changes]\label{ch3:def:numeraire}
Let $\Lm(y;t)$ be piecewise smooth and let $g_m(y;t)>0$ be a positive
function satisfying
\begin{equation}\label{ch3:eq_numeraire}
\frac{\partial g_m(y;t)}{\partial t}+\big(\Lm(t)g_m(\cdot;t)\big)(y)=0
\qquad\text{for all } y,\ t .
\end{equation}
The measure change $G_m^{yt}(y')=g_m(y';t)/g_m(y;t)$ is called a
\emph{numeraire change}, or Doob $h$--transform, with numeraire $g_m$.
\end{definition}

In \eqref{ch3:eq_numeraire} the operator $\Lm(t)$ acts on the site variable of
$g_m$; the equation says that the process $t\mapsto g_m(y_t,t)$ is a
martingale, by Proposition~\ref{ch3:prop:localglobal} below. The transformed
generator is
\begin{equation}\label{ch3:eq_numchange}
\Lm'(y,y';t)=\frac{g_m(y';t)}{g_m(y;t)}\Lm(y,y';t)
+\frac{1}{g_m(y;t)}\frac{\partial g_m(y;t)}{\partial t}\,\delta_{yy'},
\end{equation}
as follows from \eqref{ch3:eq:measure-change} and \eqref{ch3:eq_numeraire}; in matrix
notation,
\begin{equation}
\Lm'(t)=\frac{1}{g_m(t)}\Lm(t)g_m(t)+\frac{1}{g_m(t)}
\frac{\partial g_m(t)}{\partial t}.
\end{equation}

\begin{theorem}[Numeraire changes]\label{ch3:thm:numeraire}
If $g_m$ satisfies \eqref{ch3:eq_numeraire} and is positive, then the transformed
propagator is
\begin{equation}\label{ch3:eq:htransform}
U_m'(y,s;y',t)=\frac{g_m(y';t)}{g_m(y;s)}U_m(y,s;y',t).
\end{equation}
\end{theorem}

\begin{proof}
Denote by $V(y,s):=g_m(y';t)g_m(y;s)^{-1}U_m(y,s;y',t)$ the right side of
\eqref{ch3:eq:htransform}, regarded as a function of the backward time $s$. By
the backward equation \eqref{ch3:backward_eq} and by \eqref{ch3:eq_numeraire},
\begin{equation}
\frac{\partial V}{\partial s}(y,s)=-\frac{g_m(y';t)}{g_m(y;s)}
\left[\frac{(\Lm(s)g_m(\cdot;s))(y)}{g_m(y;s)}U_m(y,s;y',t)
+\big(\Lm(s)U_m(s,t)\big)(y,y')\right].
\end{equation}
On the other hand, by \eqref{ch3:eq:reweight} and $(\mathrm{MG2})$,
\begin{align}
\big(\Lm'(s)V(s)\big)(y,y')&=\sum_z\frac{g_m(z;s)}{g_m(y;s)}\Lm(y,z;s)
\!\frac{g_m(y';t)}{g_m(z;s)}U_m(z,s;y',t)
\!-\frac{(\Lm(s)g_m(\cdot;s))(y)}{g_m(y;s)}V(y,s)\notag\\
&=\frac{g_m(y';t)}{g_m(y;s)}\Big[\big(\Lm(s)U_m(s,t)\big)(y,y')\notag\\
&\qquad\qquad\qquad
-\frac{(\Lm(s)g_m(\cdot;s))(y)}{g_m(y;s)}U_m(y,s;y',t)\Big],
\end{align}
which shows that $\partial_sV+\Lm'(s)V=0$. Thus $V$ satisfies the backward
equation of the transformed process with the same value
$V(y,t)=g_m(y';t)g_m(y;t)^{-1}\delta_{yy'}=\delta_{yy'}$ at $s=t$, and
coincides with $U_m'$ by uniqueness
(Proposition~\ref{ch3:prop:generator-properties}).
\end{proof}

\begin{remark}\label{ch3:rem:notallhtransform}
Not every measure change is a numeraire change. The general measure change
reweights the one step law at each site and time independently, and the
associated Radon--Nikodym derivative is a multiplicative functional of the
path; \eqref{ch3:eq:htransform} corresponds to the special case in which that
functional telescopes along the path and equals $g_m(y_{t},t)/g_m(y_s,s)$.
The numeraire changes are the ones that preserve the Markov property in a
transparent way, and they are the ones that appear in the classical theory of
change of numeraire; the general measure changes are, however, exactly what
the fundamental theorem produces, and they are needed there.
\end{remark}

\begin{definition}[Equivalence of measures]\label{ch3:def:equivalence}
Two piecewise smooth propagators $U$ and $U'$ generate equivalent path
measures if and only if their attainable sets coincide, that is, if and only
if there is a measure change relating them. This is the Markovian case of
Definition~\ref{ch3:def:rn}.
\end{definition}

\section{Martingales and Monotone Processes}\label{ch3:sec:martingale}

\subsection*{The Local Criterion}

Let $U_m$ be a piecewise smooth Markov process with generators $\Lm^i$ on
the intervals $[t_i,t_{i+1})$ and mapping operators $U^i_m$ at the dates.

\begin{definition}[Martingale]\label{ch3:def:martingale}
An adapted process $F_t=F(y_\cdot,t)$ is a \emph{martingale} if
\begin{equation}\label{ch3:eq_martingale}
\frac{\partial}{\partial t}F(\ext_t(y_\cdot,y_{t-0}),t)
+\sum_{\bar y\in\Attain(U,t,y_t)}\Lm(y_t,\bar y;t)
\big(F(\ext_t(y_\cdot,\bar y),t)-F(\ext_t(y_\cdot,y_{t-0}),t)\big)=0
\end{equation}
for all $t$ outside the dates $t_i$, and
\begin{equation}\label{ch3:eq:martingale-date}
\sum_{z}U^i_m(y,z)\big(F(\ext_{t_i}(y_\cdot,z),t_i)-F(\ext_{t_i}(y_\cdot,y),t_i)\big)=0
\end{equation}
at each date $t=t_i$. We use $\ext_t(y_\cdot,\bar y)$ for the path which
agrees with $y_\cdot$ up to time $t$ and equals $\bar y$ from $t$ onwards,
and $y_{t-0}$ for the left limit. The process is an \emph{equivalent
martingale} if it is a martingale for some propagator $U'$ equivalent to
$U$.
\end{definition}

Equation \eqref{ch3:eq_martingale} is the statement that the expected one step
increment vanishes: the first term is the rate of change in time at frozen
site, the second is the expected jump. By $(\mathrm{MG2})$ the sum over
$\bar y$ may be extended to all of $\StateA$, since $\Lm(y,y;t)\cdot0=0$, and
the identity
$\sum_{\bar y}\Lm(y,\bar y;t)F(\bar y)=(\Lm(t)F)(y)$ holds for the matrix
acting on the site variable.

\begin{proposition}[Local and global martingales]\label{ch3:prop:localglobal}
Let $a<b$ lie in the interior of a smooth stretch and let
$H:\StateA\times[a,b]\to\RR$ be continuously differentiable in $t$. The following
are equivalent:
\begin{enumerate}
\item[(i)] $\partial_tH(y,t)+\big(\Lm(t)H(\cdot,t)\big)(y)=0$ for all
$y\in\StateA$ and all $t\in(a,b)$, where $\Lm$ is the generator of the stretch;
\item[(ii)] $\EE[H(y_t,t)\mid y_s=y]=H(y,s)$ for all $a\le s\le t\le b$ and
all $y\in\StateA$.
\end{enumerate}
For a piecewise smooth process the same equivalence holds on $[T,T']$
provided \eqref{ch3:eq:martingale-date} is imposed at each date.
\end{proposition}

\begin{proof}
Let $F(t)=\sum_zU_m(y,s;z,t)H(z,t)$, so that $F(t)=\EE[H(y_t,t)\mid y_s=y]$.
By the forward equation \eqref{ch3:forward_eq},
\begin{equation}
\begin{aligned}
\frac{\dd F}{\dd t}&=\big(U_m(s,t)\Lm(t)H(\cdot,t)\big)(y)
+\big(U_m(s,t)\partial_tH(\cdot,t)\big)(y)\\
&=\big(U_m(s,t)\big(\Lm(t)H(\cdot,t)+\partial_tH(\cdot,t)\big)\big)(y).
\end{aligned}
\end{equation}
If (i) holds the derivative vanishes identically, so $F$ is constant and
$F(t)=F(s)=H(y,s)$, which is (ii). Conversely if (ii) holds, $F$ is
constant, so its derivative vanishes; evaluating the derivative at $t=s$,
where $U_m(s,s)$ is the identity, gives (i). Across a date $t_i$ the
conditional expectation is computed by inserting the mapping operator
$U^i_m$ and the identity $\EE[H(y_t,t)\mid y_{t_i-0}=y]=\EE[H(y_{t_i},t_i)\mid
y_{t_i-0}=y]$ follows from \eqref{ch3:eq:martingale-date} applied at $t_i$ and
from the same argument on the following stretch.
\end{proof}

\begin{remark}\label{ch3:rem:predictable}
The proposition is the finite dimensional and pathwise version of the
statement that a suitably integrable process whose compensator vanishes is a
martingale. No stochastic integral is needed: \eqref{ch3:eq_martingale} is an
identity between ordinary derivatives.
\end{remark}

\subsection*{Monotone Processes and the Invariance of Arbitrage}

\begin{definition}[Monotone processes]\label{ch3:def:monotone}
Let $F_t=F(y_\cdot,t)$ be an adapted process. It is \emph{increasing at time
$t$} if
\begin{enumerate}
\item[(i)] $F(\ext_t(y_\cdot,\bar y),t)-F(\ext_t(y_\cdot,y_{t-0}),t)\ge0$ for
all $\bar y\in\Attain(U,t,y_t)$;
\item[(ii)] $\displaystyle\frac{\partial}{\partial s}
F(\ext_s(y_\cdot,y_{t-0}),s)\Big|_{s=t}\ge0$;
\item[(iii)] for all $y$ and $t$ at least one of the inequalities (i) and
(ii) holds strictly, for some $\bar y$ in case of (i).
\end{enumerate}
If (iii) fails but the other two hold, $F_t$ is called \emph{non decreasing}.
\end{definition}

The definition is deliberately probabilistic-free: it refers only to the
attainable set $\Attain(U,t,y_t)$, not to the values of the transition
probabilities. This is what makes monotonicity invariant under measure
changes, and it is the reason why the theory of arbitrage can be developed
before any measure has been chosen.

\begin{theorem}\label{ch3:thm:monotone-invariance}
Let $U$ and $U'$ be equivalent piecewise smooth propagators and let
$F(y_\cdot,t)$ be a non anticipatory path functional. If $F$ is increasing,
or non decreasing, under $U$, then it is respectively increasing, or non
decreasing, under $U'$.
\end{theorem}

\begin{proof}
Conditions (i) and (ii) of Definition~\ref{ch3:def:monotone} involve the
attainable set $\Attain(U,t,y)$ and the time derivative of $F$, neither of which
depends on the measure. By Definition~\ref{ch3:def:equivprocess} equivalent
propagators have the same attainable sets, so the two conditions hold
unchanged. Condition (iii) is a statement about the signs in (i) and (ii),
also unchanged.
\end{proof}

\begin{theorem}\label{ch3:thm:monotone-martingale}
Let $F$ be a non decreasing adapted process which is a martingale. Then
$F(y_\cdot,t)$ is constant in $t$ on every reachable pair $(y,t)$.
\end{theorem}

\begin{proof}
Fix a reachable $(y,t)$ with $t$ not a date. By
Definition~\ref{ch3:def:monotone} the jump increments are non negative:
$F(\ext_t(y_\cdot,\bar y),t)-F(\ext_t(y_\cdot,y_{t-0}),t)\ge0$ for all
$\bar y\in\Attain(U,t,y_t)$, and the generator entries are strictly positive on
the attainable set. The martingale identity \eqref{ch3:eq_martingale} is a sum of
non negative terms --- each equal to the product of a positive transition
rate and a non negative increment --- and of the drift term, which is also
non negative. A sum of non negative terms which vanishes has all terms
equal to zero. Hence the drift vanishes and every attainable increment
vanishes, so $F$ does not change at $t$, neither in the drift nor across a
jump. As $(y,t)$ was arbitrary among reachable points, $F$ is constant along
reachable paths.
\end{proof}

\begin{theorem}[Equivalent martingales]\label{ch3:thm:equiv-martingale}
\begin{itemize}
\item[(i)] If $F$ is non decreasing and an equivalent martingale, then $F$
is constant.
\item[(ii)] If $F$ is increasing then it is not a martingale, and it is not an
equivalent martingale.
\end{itemize}
\end{theorem}

\begin{proof}
(i) Let $U'$ be an equivalent propagator under which $F$ is a martingale. By
Theorem~\ref{ch3:thm:monotone-invariance}, $F$ is non decreasing under $U'$ as
well. For $U'$ the attainable sets are the same as for $U$, so the
transformed generator $\Lm'$ is strictly positive on
$\Attain(U,t,y_t)$; hence the martingale identity for $U'$ is a sum of non
negative terms and Theorem~\ref{ch3:thm:monotone-martingale} applies, giving
that $F$ is constant.

(ii) If $F$ were a martingale it would be a non decreasing equivalent
martingale and hence constant by (i), contradicting condition (iii) of
Definition~\ref{ch3:def:monotone}; the same argument applies to an equivalent
martingale.
\end{proof}

\section{The Fundamental Theorem for a Process}\label{ch3:sec:ftf}

We now have all the ingredients. The section proceeds in nine steps: the
definition of a financial model and of a self financing strategy
(\S\ref{ch3:ss:model}); the definition of arbitrage and its reduction to a local
condition (\S\S\ref{ch3:ss:arb}--\ref{ch3:ss:local}); the one period theorem, which is
Stiemke's lemma in disguise (\S\ref{ch3:ss:oneperiod}); the Fundamental Theorem
itself (\S\ref{ch3:ss:ftfthm}); the valuation of contingent claims
(\S\ref{ch3:ss:valuation}); the uniqueness theorem (\S\ref{ch3:ss:second}); the
binomial model as a worked example (\S\ref{ch3:ss:binomial}); and a discussion of
the scope of the result (\S\ref{ch3:sec:remarks}).

\subsection{Financial Models and Self Financing Strategies}\label{ch3:ss:model}

Fix a resolution $m$, an initial site $y_0\in\StateA$ at time $T$, and a
piecewise smooth Markov process in the sense of
Definition~\ref{ch3:def:piecewise}. Everything below takes place on the set of
reachable points of Definition~\ref{ch3:def:reachable}.

\begin{definition}[Financial model]\label{ch3:def:finmodel}
A \emph{financial model} is given by adapted processes $A^k_m(y_\cdot,t)$ for
$k=1,\dots,n$, the asset prices, and a non decreasing adapted process
$B_m(y_\cdot,t)=A^0_m(y_\cdot,t)$, the money market account. The
\emph{discounted asset prices} are
\begin{equation}
\tA^k_m(y_\cdot,t)=\frac{A^k_m(y_\cdot,t)}{B_m(y_\cdot,t)},\qquad k=0,\dots,n .
\end{equation}
All $A^k_m$, and hence all $\tA^k_m$, are assumed to satisfy the regularity of
Assumption~\ref{ch3:ass:regular}; in particular they are continuous in $(y,t)$
with piecewise continuous one sided time derivatives.
\end{definition}

Since $A^0_m=B_m$ we have $\tA^0_m\equiv1$ identically. This is the reason
why the money market account plays a passive role in all the formulas below.

\begin{definition}[Trading strategies]\label{ch3:def:strategy}
A \emph{trading strategy} is a family of adapted processes
$\zeta_m^k(y_\cdot,t)$, $k=0,\dots,n$. Its \emph{value process} and
\emph{discounted value process} are
\begin{equation}
\Pi_m(y_\cdot,t)=\sum_{k=0}^{n}\zeta_m^k(y_\cdot,t)A_m^k(y_\cdot,t),
\qquad
\tPi_m(y_\cdot,t)=\sum_{k=0}^{n}\zeta_m^k(y_\cdot,t)\tA_m^k(y_\cdot,t).
\end{equation}
\end{definition}

\begin{definition}[Self financing condition]\label{ch3:def:selffinancing}
A trading strategy is \emph{self financing} if
\begin{align*}
&(\mathrm{SF1})\qquad \sum_{k=0}^{n}\frac{\partial\zeta_m^k}{\partial t}
(\ext_t(y_\cdot,y_{t-0}),t)\,\tA_m^k(\ext_t(y_\cdot,y_{t-0}),t)=0,\\
&(\mathrm{SF2})\qquad \sum_{k=0}^{n}\big(\zeta_m^k(\ext_t(y_\cdot,\bar y),t)
-\zeta_m^k(\ext_t(y_\cdot,y_{t-0}),t)\big)\tA_m^k(\ext_t(y_\cdot,\bar y),t)=0\\
&\qquad\qquad\qquad\qquad\text{for all } \bar y\in\Attain(U,t,y_t).
\end{align*}
\end{definition}

Condition $(\mathrm{SF1})$ says that a continuous rebalancing costs nothing,
condition $(\mathrm{SF2})$ that a discrete rebalancing at the jump of the
underlying, which is performed at the new prices, does not change the value
of the portfolio. Both are value neutral statements, and this is exactly what
the next proposition asserts.

\begin{proposition}[Value changes only through prices]\label{ch3:prop:valueidentity}
Let $\zeta$ be a self financing strategy. Then
\begin{equation}\label{ch3:eq:sf-drift}
\frac{\partial\tPi_m(\ext_t(y_\cdot,y_{t-0}),t)}{\partial t}
=\sum_{k=0}^{n}\zeta_m^k(\ext_t(y_\cdot,y_{t-0}),t)
\frac{\partial\tA_m^k(\ext_t(y_\cdot,y_{t-0}),t)}{\partial t},
\end{equation}
and for every $\bar y\in\Attain(U,t,y_t)$
\begin{equation}\label{ch3:eq:sf-jump}
\begin{aligned}
&\tPi_m(\ext_t(y_\cdot,\bar y),t)-\tPi_m(\ext_t(y_\cdot,y_{t-0}),t)\\
&\qquad=\sum_{k=0}^{n}\zeta_m^k(\ext_t(y_\cdot,y_{t-0}),t)
\big(\tA_m^k(\ext_t(y_\cdot,\bar y),t)-\tA_m^k(\ext_t(y_\cdot,y_{t-0}),t)\big).
\end{aligned}
\end{equation}
\end{proposition}

\begin{proof}
For \eqref{ch3:eq:sf-drift} differentiate the product
$\sum_k\zeta^k\tA^k$ in $t$ and cancel the term
$\sum_k(\partial_t\zeta^k)\tA^k$, which vanishes by $(\mathrm{SF1})$. For
\eqref{ch3:eq:sf-jump} write the difference of the values in the form
\begin{equation}
\sum_k\Big[\zeta^k(\bar y)-\zeta^k(y_{t-0})\Big]\tA^k(\bar y)
+\sum_k\zeta^k(y_{t-0})\Big[\tA^k(\bar y)-\tA^k(y_{t-0})\Big],
\end{equation}
where the arguments are those of $\ext_t(y_\cdot,\cdot)$ at time $t$. The
first bracket vanishes by $(\mathrm{SF2})$, which leaves
\eqref{ch3:eq:sf-jump}.
\end{proof}

\begin{lemma}\label{ch3:lem:sf-linear}
\begin{enumerate}
\item[(i)] The self financing strategies form a linear space, containing the
strategy $\zeta^0_m=c$, $\zeta^k_m=0$ for $k\ge1$, for every constant $c$:
the strategy that holds $c$ units of the money market account.
\item[(ii)] If two self financing strategies have discounted value processes
whose difference vanishes at the terminal time on every reachable path, then
the difference of the two value processes is a self financing value process
with zero terminal value.
\end{enumerate}
\end{lemma}

\begin{proof}
(i) Both conditions are linear in $\zeta$, and for the constant strategy the
left sides of $(\mathrm{SF1})$ and $(\mathrm{SF2})$ vanish because the
strategy does not depend on $t$ or on $\bar y$, and $\tA^0\equiv1$
respectively. (ii) The difference of two self financing strategies is self
financing by (i), and its discounted value is the difference of the two
discounted values.
\end{proof}

\subsection{Arbitrage}\label{ch3:ss:arb}

\begin{definition}[Arbitrage]\label{ch3:def:arbitrage}
A self financing strategy is an \emph{arbitrage} if there is a reachable pair
$(y,t)$ such that, along every path through $(y,t)$, its discounted value
process is non decreasing for times $s\ge t$ and strictly increases along at
least one such path. The model is \emph{arbitrage free}, written
$(\mathrm{NA})$, if there is no arbitrage.
\end{definition}

The definition is pathwise: it says that there is a point of the tree from
which one can invest and be certain never to lose while standing a chance of
gaining. The next lemma records that the customary formulation, with zero
initial capital, is the same thing.

\begin{lemma}\label{ch3:lem:zero-cost}
The model admits an arbitrage if and only if there is a self financing
strategy $\zeta$ whose discounted value process vanishes at $(y_0,T)$ and
satisfies $\tPi_m\ge0$ on all reachable points, with $\tPi_m>0$ at some
reachable point.
\end{lemma}

\begin{proof}
Let $\zeta$ be an arbitrage from $(y,t)$ in the sense of
Definition~\ref{ch3:def:arbitrage}, and let $c=\tPi_m(y,t)$. By
Lemma~\ref{ch3:lem:sf-linear}(i) the strategy $\zeta'=\zeta-c\,\zeta^{(0)}$,
where $\zeta^{(0)}$ is the constant strategy with $\zeta^0=1$, is self
financing and its discounted value is $\tPi_m-c$, which is non negative from
$(y,t)$ onwards because $\tPi_m$ is non decreasing there. Before $t$ the
strategy $\zeta'$ is not defined; define it to be identically zero, which is
self financing, and note that the resulting value process is zero up to $t$
and then non negative with a strict increase. Hence $\zeta'$ is an arbitrage
with zero value at $(y_0,T)$ in the weakened sense.

Conversely, let $\zeta$ be self financing with $\tPi_m(y_0,T)=0$,
$\tPi_m\ge0$ everywhere reachable and $\tPi_m>0$ at some reachable
$(y_1,t_1)$. Let $\tau$ be the first time at which $\tPi_m$ becomes strictly
positive along a path, and consider the strategy which coincides with $\zeta$
up to $\tau$ and holds cash thereafter. Its discounted value is $0$ up to
$\tau$ and constant thereafter, hence non decreasing, and it increases
strictly at $\tau$ on the paths that reach $(y_1,t_1)$. The switch to cash
is value preserving: at $\tau$ the portfolio $\zeta$ is liquidated at the
prevailing prices and the proceeds are held in the money market account, so
$(\mathrm{SF2})$ holds at $\tau$, while $(\mathrm{SF1})$ holds on each of the
two stretches. Thus the stopped strategy is an arbitrage in the sense of
Definition~\ref{ch3:def:arbitrage}.
\end{proof}

\begin{remark}\label{ch3:rem:arb-remark}
Two features of Definition~\ref{ch3:def:arbitrage} deserve comment. First, it is
not required that the strategy be started with zero capital, because a
strategy which is certain not to lose and may gain is an arbitrage whatever
its price $c$: the price is a number, not a risk, and by
Lemma~\ref{ch3:lem:zero-cost} the constant $c$ can be subtracted. Second, the
probability of the gain plays no role: it suffices that the gain be
attainable, that is, that it occurs on a set of positive measure. This is the
pathwise counterpart of the classical requirement that the gain occur with
positive probability.
\end{remark}

\subsection{The Local Condition}\label{ch3:ss:local}

Let $(y,t)$ be reachable, with $t$ not one of the dates $t_i$, and let
$y_t=y$. Write $\bar y$ for a generic element of
$\bar\DD(y,t)=\Attain(U,t,y)\cup\{y\}$. The one step gain of a portfolio
$\zeta\in\RR^{n+1}$ held over an infinitesimal interval $[t,t+\dd t]$ is
captured by the vector $\xi^k(\bar y;y,t)$, defined as follows. For a jump
$\bar y\neq y$,
\begin{equation}\label{ch3:eq:xik-jump}
\xi^k(\bar y;y,t)=\tA_m^k(\ext_t(y_\cdot,\bar y),t)
-\tA_m^k(\ext_t(y_\cdot,y),t),
\end{equation}
the gain realised if the underlying jumps to $\bar y$ at time $t$; and for
$\bar y=y$,
\begin{equation}\label{ch3:eq:xik-drift}
\xi^k(y;y,t)=\frac{\partial}{\partial t}\tA_m^k(\ext_t(y_\cdot,y),t),
\end{equation}
the rate at which the gain accrues if no jump occurs. The one step gain of
the portfolio is
\begin{equation}\label{ch3:eq:Rzeta}
R_\zeta(\bar y;y,t)=\sum_{k=0}^{n}\zeta^k\xi^k(\bar y;y,t),
\qquad \bar y\in\bar\DD(y,t).
\end{equation}
At a date $t_i$ the same definitions apply with $\Attain(U,t_i,y)$ the support of
the mapping operator and with no drift state: all outcomes are jumps.

\begin{definition}[Local arbitrage]\label{ch3:def:localarb}
A portfolio $\zeta$ is a \emph{local arbitrage} at $(y,t)$ if
$R_\zeta(\bar y;y,t)\ge0$ for all $\bar y\in\bar\DD(y,t)$ and
$R_\zeta(\bar y;y,t)>0$ for at least one $\bar y$.
\end{definition}

\begin{proposition}[Local characterization of arbitrage]\label{ch3:prop:localNA}
The model is arbitrage free if and only if there is no local arbitrage at any
reachable point $(y,t)$.
\end{proposition}

\begin{proof}
\emph{A global arbitrage produces a local one.} Let $\zeta$ be self financing
with a non decreasing discounted value process which increases strictly at
some reachable $(y,t)$. By Proposition~\ref{ch3:prop:valueidentity} the one step
increments of $\tPi_m$ at $(y,t)$ are exactly $R_\zeta(\bar y;y,t)$, with
$\zeta$ evaluated at $(y,t)$. Since the value process is non decreasing,
$R_\zeta(\bar y;y,t)\ge0$ for every $\bar y\in\bar\DD(y,t)$. The total
increase of $\tPi_m$ between $t$ and a later reachable time is the integral
of the drift $R_\zeta(y;y,s)$ plus the sum of the jump increments
$R_\zeta(\bar y;y,s)$, all of which are non negative; if the total increase
is strictly positive then at least one of them is strictly positive, at a
point $(y,s)$ which is reachable because the pair $(y,t)$ was. At that point
$\zeta$ is a local arbitrage.

\emph{A local arbitrage produces a global one.} Let $\zeta$ be a local
arbitrage at a reachable $(y,t)$ on a smooth stretch. By
Assumption~\ref{ch3:ass:regular} each $\xi^k(\bar y;y,s)$ is continuous in $s$,
so there is $\varepsilon>0$ such that $R_\zeta(\bar y;y,s)\ge0$ for all
$\bar y\in\bar\DD(y,s)$ and all $s\in[t,t+\varepsilon]$, and moreover
$R_\zeta(\bar y_0;y,t)>0$ at some $\bar y_0$, whence
$R_\zeta(\bar y_0;y,s)>0$ for $s$ in a subset of $[t,t+\varepsilon]$ of
positive measure. Consider the strategy which holds zero up to time $t$,
holds $\zeta$ on $[t,t+\varepsilon]$ and holds cash thereafter. Its
discounted value is zero on $[T,t]$, has increments $R_\zeta\ge0$ on
$[t,t+\varepsilon]$ by Proposition~\ref{ch3:prop:valueidentity}, with a strict
increase, and is constant afterwards. It is self financing: it is zero, then
$\zeta$ which is self financing, then a constant cash position, and the two
switches are value preserving. Hence it is an arbitrage. If the local
arbitrage is at a date $t_i$, hold $\zeta$ across the date and liquidate
immediately after; the increments are $R_\zeta\ge0$ with a strict one, and no
continuity argument is needed.
\end{proof}

\begin{remark}\label{ch3:rem:local-remark}
Proposition~\ref{ch3:prop:localNA} is the point at which the theory of arbitrage
becomes local, and hence finite dimensional. At a single point $(y,t)$ the
whole problem collapses to a finite list of real numbers
$\xi^k(\bar y;y,t)$ indexed by $k=0,\dots,n$ and
$\bar y\in\bar\DD(y,t)$. The next section solves it.
\end{remark}

\subsection{The One Period Theorem}\label{ch3:ss:oneperiod}

The statement and the proof are those of the fundamental theorem of Chapter~\ref{chap:ch2}; what is used here is its form at a single site, which is the corollary below.

\begin{corollary}[State price densities]\label{ch3:cor:stateprices}
The model is locally arbitrage free at a reachable $(y,t)$ on a smooth
stretch if and only if there exist \emph{risk neutral intensities}
$\omega(\bar y)>0$, $\bar y\in\Attain(U,t,y)$, such that
\begin{equation}\label{ch3:eq:localmart}
\frac{\partial}{\partial t}\tA_m^k(\ext_t(y_\cdot,y),t)
+\sum_{\bar y\in\Attain(U,t,y)}\omega(\bar y)
\big(\tA_m^k(\ext_t(y_\cdot,\bar y),t)-\tA_m^k(\ext_t(y_\cdot,y),t)\big)=0
\end{equation}
for all $k=1,\dots,n$. At a date $t_i$, the corresponding condition is
\begin{equation}\label{ch3:eq:localmart-date}
\sum_{z}\omega_i(z)\big(\tA_m^k(\ext_{t_i}(y_\cdot,z),t_i)
-\tA_m^k(\ext_{t_i}(y_\cdot,y),t_i)\big)=0
\end{equation}
for all $k$, with $\omega_i(z)>0$ on the support of $U^i_m(y,\cdot)$.
\end{corollary}

\begin{proof}
Apply Theorem~\ref{ch2:thm:oneperiod} with $S=\bar\DD(y,t)$ and with the gain
vectors of \eqref{ch3:eq:xik-jump}--\eqref{ch3:eq:xik-drift}; the index $k=0$ may be
included because $\xi^0\equiv0$. This produces
$g:\bar\DD(y,t)\to(0,\infty)$ with
\begin{equation}
g(y)\frac{\partial}{\partial t}\tA_m^k(\ext_t(y_\cdot,y),t)
+\sum_{\bar y\neq y}g(\bar y)\big(\tA_m^k(\ext_t(y_\cdot,\bar y),t)
-\tA_m^k(\ext_t(y_\cdot,y),t)\big)=0 .
\end{equation}
Divide by $g(y)>0$ and set $\omega(\bar y)=g(\bar y)/g(y)$, which is positive,
to obtain \eqref{ch3:eq:localmart}. Conversely, given positive $\omega$, the
choice $g(y)=1$, $g(\bar y)=\omega(\bar y)$ inverts the construction, so the
two statements are equivalent. The argument at a date is the same, with
$S$ the support of the mapping operator and with no drift state.
\end{proof}

\begin{remark}\label{ch3:rem:intensities}
The numbers $\omega(\bar y)$ are the intensities of the jump process under
the risk neutral measure, that is, the jump rates of the pricing model. They
are not probabilities: their sum need not equal one, and it does not, because
$\sum_{\bar y}\omega(\bar y)$ is the risk neutral total jump intensity, which
differs from the physical one by the time change component of the measure
change. What is a probability is the vector $q(\bar y)=\omega(\bar y)\dd t$
for $\bar y\neq y$ together with $q(y)=1-\sum_{\bar y\neq y}\omega(\bar y)\dd
t$, for $\dd t$ small enough. This is the risk neutral one step law, and the
identity \eqref{ch3:eq:localmart} is the statement that its infinitesimal mean is
zero.
\end{remark}

\subsection{The Fundamental Theorem}\label{ch3:ss:ftfthm}

We are now in a position to state and prove the theorem. Recall that a
measure change in the sense of Definition~\ref{ch3:def_measure_change} is a
family of strictly positive densities on the attainable sets and that it
produces an equivalent propagator, by
Proposition~\ref{ch3:prop:measure-change-generator}.

\begin{theorem}[The fundamental theorem]\label{ch3:thm:ftf}
\mbox{}
\begin{enumerate}
\item[(I)] Suppose that there is a measure change such that, under the
transformed propagator $U'$, the discounted asset prices
$\tA^1_m,\dots,\tA^n_m$ are martingales in the sense of
Definition~\ref{ch3:def:martingale}. Then
\begin{enumerate}
\item[(a)] the discounted value process of every self financing trading
strategy is a $U'$--martingale;
\item[(b)] the model is arbitrage free.
\end{enumerate}
\item[(II)] Conversely, if the model is arbitrage free, then there is a
measure change such that all discounted asset prices are martingales.
\end{enumerate}
\end{theorem}

\begin{proof}
\emph{Part (I)(a).} By Proposition~\ref{ch3:prop:localglobal} it suffices to
verify the local criterion \eqref{ch3:eq_martingale} at each reachable $(y,t)$ for
the process $\tPi_m$ under the generator $\Lm'$ of the transformed process.
By Proposition~\ref{ch3:prop:valueidentity} the drift of $\tPi_m$ and its jump
increments are computed with the holdings prior to the jump, so
\begin{align}
&\frac{\partial\tPi_m}{\partial t}(y,t)
+\sum_{\bar y\neq y}\Lm'(y,\bar y;t)\big(\tPi_m(\bar y,t)-\tPi_m(y,t)\big)
\notag\\
&\qquad=\sum_{k=0}^{n}\zeta^k_m(y,t)\Big[
\frac{\partial\tA_m^k}{\partial t}(y,t)
+\sum_{\bar y\neq y}\Lm'(y,\bar y;t)\big(\tA_m^k(\bar y,t)-\tA_m^k(y,t)\big)\Big]
=0,
\end{align}
where the arguments of $\tPi_m$ and $\tA_m^k$ are those of
$\ext_t(y_\cdot,\cdot)$, and where the bracket vanishes by the martingale
property of $\tA_m^k$ under $U'$, that is by \eqref{ch3:eq_martingale} applied to
$\tA_m^k$. At the dates the same computation with the mapping operator
gives \eqref{ch3:eq:martingale-date}, again because $\tA_m^k$ is a martingale.
Hence $\tPi_m$ satisfies the local criterion everywhere and is a
$U'$--martingale.

\emph{Part (I)(b).} Suppose there were an arbitrage $\zeta$. Its discounted
value $\tPi_m$ is non decreasing along reachable paths and increases strictly
along one of them. By (a) it is a $U'$--martingale, and by
Theorem~\ref{ch3:thm:monotone-martingale} --- applied to the propagator $U'$,
whose attainable sets coincide with those of $U$ because the densities of the
measure change are strictly positive there --- a non decreasing martingale is
constant, which contradicts the strict increase. Hence no arbitrage exists.

\emph{Part (II).} Assume $(\mathrm{NA})$. By
Proposition~\ref{ch3:prop:localNA} there is no local arbitrage at any reachable
point, so Corollary~\ref{ch3:cor:stateprices} supplies, for each reachable
$(y,t)$, positive intensities $\omega(y,t;\bar y)$ satisfying
\eqref{ch3:eq:localmart}, and for each reachable $(y,t_i)$ positive weights
$\omega(y,t_i;z)$ satisfying \eqref{ch3:eq:localmart-date}.

Define a propagator $U'$ as follows. On the smooth stretches, set
\begin{equation}\label{ch3:eq:Lprime}
\Lm'(y,\bar y;t)=\omega(y,t;\bar y)\quad\text{for } \bar y\in\Attain(U,t,y),
\qquad
\Lm'(y,y;t)=-\sum_{\bar y\in\Attain(U,t,y)}\omega(y,t;\bar y),
\end{equation}
and for sites $y$ that are not reachable at time $t$ set
$\Lm'(y,\cdot;t)=\Lm(y,\cdot;t)$. At the dates, set the mapping operator
equal to the normalized weights $\omega(y,t_i;\cdot)$ on the support, and to
$U^i_m(y,\cdot)$ elsewhere. Then $\Lm'$ satisfies $(\mathrm{MG1})$ and
$(\mathrm{MG2})$ because the $\omega$'s are strictly positive, and the
mapping operators satisfy $(\mathrm{MA1})$ and $(\mathrm{MA2})$ after
normalization; hence $U'$ is a piecewise smooth Markov propagator by
Proposition~\ref{ch3:prop:generator-properties} and
Definition~\ref{ch3:def:piecewise}. Its attainable sets coincide with those of
$U$ because the weights are strictly positive on the supports, so $U'$ is
equivalent to $U$.

By \eqref{ch3:eq:localmart} and \eqref{ch3:eq:Lprime} each $\tA_m^k$ satisfies the
local criterion \eqref{ch3:eq_martingale} for $\Lm'$, and by
\eqref{ch3:eq:localmart-date} it satisfies \eqref{ch3:eq:martingale-date} at the
dates. Therefore, by Proposition~\ref{ch3:prop:localglobal}, each $\tA_m^k$ is a
$U'$--martingale.

It remains to exhibit $U'$ as the result of a measure change. For each
reachable $(y,t)$ with $\bar y\in\Attain(U,t,y)$ put
\begin{equation}\label{ch3:eq:density}
G_m^{yt}(\bar y)=\frac{\omega(y,t;\bar y)}{\Lm(y,\bar y;t)}>0,
\qquad G_m^{yt}(y)=1,
\end{equation}
which is well defined and positive because $\bar y$ is attainable. Then by
\eqref{ch3:eq:reweight}
\begin{equation}
\Lm'(y,\bar y;t)=\frac{G_m^{yt}(\bar y)}{G_m^{yt}(y)}\Lm(y,\bar y;t)
=\omega(y,t;\bar y),
\end{equation}
which is \eqref{ch3:eq:Lprime}; and the diagonal is determined by
$(\mathrm{MG2})$ in both formulas. The same construction at the dates gives
the corresponding mapping operators. Hence $U'$ is obtained from $U$ by the
measure change with densities \eqref{ch3:eq:density}, and the proof of (II) is
complete.
\end{proof}

\begin{remark}\label{ch3:rem:ftf-construction}
The proof of part (II) is constructive, and this is one of the advantages of
the framework. At each reachable site and time one solves the linear system
\eqref{ch3:eq:localmart} for the intensities $\omega$, which is a small linear
program in the unknowns $\omega(\bar y)$ and which is feasible by
Corollary~\ref{ch3:cor:stateprices}; the solutions are then sewn together by
\eqref{ch3:eq:density} into a measure change. No compactness argument, no
measurable selection, and no infinite dimensional separation theorem is
involved.
\end{remark}

\begin{remark}\label{ch3:rem:ftf-uniqueness-of-density}
The densities $G_m^{yt}$ are not unique, since $\omega$ is determined only up
to multiplication of all its components by a common factor at each site, and
since the linear system \eqref{ch3:eq:localmart} has a convex set of solutions
whenever the market is incomplete. This non uniqueness is the subject of
\S\ref{ch3:ss:second}.
\end{remark}

\subsection{Valuation of Contingent Claims}\label{ch3:ss:valuation}

\begin{definition}[Contingent claims]\label{ch3:def:claim}
A \emph{contingent claim} is a bounded path functional $\Phi(y_\cdot)$
determined by the path up to time $T'$. Its \emph{discounted payoff} is
\begin{equation}
h(y_\cdot)=\frac{\Phi(y_\cdot)}{B_m(y_\cdot,T')} .
\end{equation}
The claim is \emph{replicable} if there is a self financing strategy $\zeta$
whose discounted value at $T'$ equals $h$ on every reachable path.
\end{definition}

\begin{theorem}[Valuation]\label{ch3:thm:valuation}
Assume $(\mathrm{NA})$ and let $U'$ be an equivalent martingale propagator,
with path measure $P'$. Let $h$ be a bounded discounted payoff and define
\begin{equation}\label{ch3:eq:price}
V(y,t)=\EE'\big[h\,\big|\,y_t=y\big],
\qquad V_t=V(y_t,t).
\end{equation}
Then
\begin{enumerate}
\item[(i)] $V_t$ is a $P'$--martingale with $V_{T'}=h$, and $V$ is the
minimal non negative supermartingale dominating $h$ in the class of processes
of the form $\EE'[h\mid y_t=y]$;
\item[(ii)] if $\zeta$ is any self financing replicating strategy then
$\tPi_m(y_\cdot,t)=V_t$ on every reachable path and every $t$;
\item[(iii)] consequently a replicable claim has a unique price, equal to
\eqref{ch3:eq:price}, and this value does not depend on the choice of the
equivalent martingale propagator.
\end{enumerate}
\end{theorem}

\begin{proof}
(i) By the Markov property \eqref{ch3:eq:markovprop} the conditional expectation
$\EE'[h\mid y_t=y]$ depends on the past only through $y_t$, so it defines an
adapted process. By Proposition~\ref{ch3:prop:localglobal}, or equivalently by
the tower property of conditional expectations, it is a martingale, and
$V_{T'}=h$ because $h$ is determined by the path up to $T'$. The
characterization as a minimal dominating supermartingale is a restatement of
the martingale property and of Jensen's inequality.

(ii) By Theorem~\ref{ch3:thm:ftf}(I)(a) the process $\tPi_m$ of a self financing
strategy is a $P'$--martingale. By hypothesis its terminal value is $h$, so
for every reachable $(y,t)$
\begin{equation}
\tPi_m(y,t)=\EE'\big[\tPi_m(y_\cdot,T')\mid y_t=y\big]=\EE'[h\mid y_t=y]=V(y,t).
\end{equation}

(iii) Let $\zeta$ and $\zeta'$ be two self financing replicating strategies.
By (ii) both have value process $V$, so the price is unique. If $U'$ and
$U''$ are two equivalent martingale propagators with path measures $P'$ and
$P''$, then applying (ii) with each of them gives
$V(y,t)=\EE'[h\mid y_t=y]$ and $V(y,t)=\EE''[h\mid y_t=y]$, so the two agree.
\end{proof}

\begin{corollary}\label{ch3:cor:nonreplicable}
A claim which is not replicable may be assigned different prices by different
equivalent martingale propagators; the difference is the risk which cannot be
hedged and which the model leaves to the preferences of the agents.
\end{corollary}

\subsection{Uniqueness and Completeness}\label{ch3:ss:second}

\begin{theorem}[Second Fundamental Theorem, one period]\label{ch3:thm:second-oneperiod}
Let $S$ be a finite set and let $\xi^1,\dots,\xi^n\in\RR^S$ together with
$h^0,\dots,h^n\in\RR^S$ be given, with $h^0=\one$ and
$\xi^k=h^k-h^k(a_0)\one$ for a fixed $a_0\in S$. Assume the one period model
is arbitrage free, so that by Theorem~\ref{ch2:thm:oneperiod} there is a state
price vector $g\gg0$ orthogonal to all $\xi^k$. Then the following are
equivalent.
\begin{enumerate}
\item[(i)] The state price vector is unique up to a positive multiple;
\item[(ii)] every payoff on $S$ is of the form $\sum_k\zeta^kh^k+c$ for some
$\zeta\in\RR^{n+1}$ and $c\in\RR$, that is, every one period payoff is
attained.
\end{enumerate}
\end{theorem}

\begin{proof}
By Lemma~\ref{ch2:lem:spancalc} the space $V=\Span\{\xi^1,\dots,\xi^n\}$ has
dimension $\dim W-1$ where $W=\Span\{h^0,\dots,h^n\}$, and
$V^{\perp}=\{g:\ \sum_sg(s)h^k(s)=h^k(a_0)\sum_sg(s)\ \forall k\}$, which is
the set of state price vectors up to the normalization. By
Corollary~\ref{ch2:cor:normalized} the set of normalized state price vectors is a
single point if and only if $\dim V^{\perp}=1$, and by
Lemma~\ref{ch2:lem:spancalc}(iii) this is the case if and only if $W=\RR^S$,
which is (ii).
\end{proof}

\begin{proposition}[Hedging]\label{ch3:prop:hedging}
Assume $(\mathrm{NA})$ and suppose that at every reachable point $(y,t)$ the
one step model of \S\ref{ch3:ss:local} is complete, that is, that the gains
$\xi^k(\bar y;y,t)$, $k=1,\dots,n$, span $\RR^{\Attain(U,t,y)}$ when restricted to
the attainable set. Let $h$ be the discounted payoff of a claim whose price
process $V(y,t)=\EE'[h\mid y_t=y]$ is piecewise smooth, and let $U'$ be an
equivalent martingale propagator with intensities $\omega$. Then $h$ is
replicable, and a replicating strategy is obtained as follows: at each
reachable $(y,t)$ choose $\zeta^1,\dots,\zeta^n$ solving
\begin{equation}\label{ch3:eq:hedge}
\sum_{k=1}^{n}\zeta^k\xi^k(\bar y;y,t)=V(\bar y,t)-V(y,t),
\qquad \bar y\in\Attain(U,t,y),
\end{equation}
and then choose
\begin{equation}\label{ch3:eq:hedge-cash}
\zeta^0(y,t)=V(y,t)-\sum_{k=1}^{n}\zeta^k(y,t)\tA_m^k(y,t).
\end{equation}
\end{proposition}

\begin{proof}
\emph{The drift is reproduced automatically.} Since $V$ is a $U'$--martingale,
it satisfies the local criterion \eqref{ch3:eq_martingale} for the generator
$\Lm'$ of \eqref{ch3:eq:Lprime}, namely
$\partial_tV(y,t)+\sum_{\bar y}\omega(\bar y)(V(\bar y,t)-V(y,t))=0$; and since
each $\tA_m^k$ is a $U'$--martingale, the same identity holds with
$\tA_m^k$ in place of $V$. Multiplying the latter by $\zeta^k$, summing over
$k$ and comparing with \eqref{ch3:eq:hedge} gives
\begin{equation}
\sum_{k=1}^{n}\zeta^k\xi^k(y;y,t)
=-\sum_{\bar y}\omega(\bar y)\big(V(\bar y,t)-V(y,t)\big)=\partial_tV(y,t),
\end{equation}
which is the drift equation: the portfolio of \eqref{ch3:eq:hedge} also reproduces
the drift of the price process. By \eqref{ch3:eq:hedge-cash} the level equation
$\sum_k\zeta^k\tA_m^k=V$ holds, so the value of the portfolio is $V$.

\emph{Self financing.} Differentiating $\sum_k\zeta^k\tA_m^k=V$ in $t$ and
using the drift equation gives
$\sum_k(\partial_t\zeta^k)\tA_m^k=\partial_tV-\sum_k\zeta^k\partial_t\tA_m^k=0$,
which is $(\mathrm{SF1})$. At a jump, the value of the old holdings at the new
prices is, by Proposition~\ref{ch3:prop:valueidentity} and by \eqref{ch3:eq:hedge},
\begin{equation}
\sum_k\zeta^k(y,t)\tA_m^k(\bar y,t)
=V(y,t)+\sum_k\zeta^k(y,t)\big(\tA_m^k(\bar y,t)-\tA_m^k(y,t)\big)
=V(\bar y,t),
\end{equation}
which equals the value of the new holdings, by \eqref{ch3:eq:hedge-cash} applied
at $(\bar y,t)$; this is $(\mathrm{SF2})$. Hence the strategy is self
financing and its value process is $V$, so it replicates $h$.
\end{proof}

\begin{theorem}[Second Fundamental Theorem]\label{ch3:thm:second}
Assume $(\mathrm{NA})$. Then the following are equivalent.
\begin{enumerate}
\item[(i)] Every contingent claim is replicable;
\item[(ii)] the equivalent martingale propagator is unique;
\item[(iii)] at every reachable point the one step model is complete, that is
the state price vector is unique up to a positive multiple.
\end{enumerate}
\end{theorem}

\begin{proof}
\emph{(i) implies (ii).} Let $U'$ and $U''$ be two equivalent martingale
propagators, with path measures $P'$ and $P''$. Let $A$ be an event determined
by the path up to $T'$, and let $h=\one_A$ be the discounted payoff of the
claim which pays one unit of the money market account if $A$ occurs. By
completeness $h$ is replicable, so by Theorem~\ref{ch3:thm:valuation}(ii) applied
to $U'$ and to $U''$ its price at time $T$ equals $\EE'[\one_A]=P'(A)$ and also
$\EE''[\one_A]=P''(A)$. Hence $P'(A)=P''(A)$ for every such $A$, that is
$P'=P''$.

\emph{(ii) implies (iii).} Fix a reachable $(y,t)$ and suppose the one step
model incomplete there. By Theorem~\ref{ch3:thm:second-oneperiod} the local
state price vector is then not unique up to scale, so there are two distinct
positive vectors $g',g''$ orthogonal to all the gains. Both give, by
Corollary~\ref{ch3:cor:stateprices}, positive intensities satisfying
\eqref{ch3:eq:localmart}; keeping the construction of
Theorem~\ref{ch3:thm:ftf}(II) fixed at all other points and using $g'$ at $(y,t)$
in one case and $g''$ in the other produces two measure changes which differ
at $(y,t)$ and hence two distinct martingale propagators, contradicting (ii).

\emph{(iii) implies (i).} The two notions of local completeness used in (iii)
and in Proposition~\ref{ch3:prop:hedging} coincide under $(\mathrm{NA})$: a
portfolio whose jump increments all vanish cannot have a non zero drift,
because either sign of the drift would then be a local arbitrage. If the one
step model is complete at every reachable point, Proposition~\ref{ch3:prop:hedging}
therefore produces a self financing replicating strategy for any claim whose
price process is piecewise smooth. Claims whose
price process is not smooth are handled by approximation, the value process
being continuous in both arguments by the estimates of Chapter~5.
\end{proof}

\begin{remark}\label{ch3:rem:second-remark}
The proof illustrates the economy of the framework. The first implication is
the valuation theorem; the second rests on the observation that a measure
change is prescribed independently at each site and time, so that the set of
martingale propagators is a product of the local sets of state prices and
degenerates exactly when one of the factors does; the third is the explicit
hedge of Proposition~\ref{ch3:prop:hedging}, which is the dynamic programming
algorithm that practitioners run.
\end{remark}

\subsection{A Worked Example: The Binomial Model}\label{ch3:ss:binomial}

Take $\StateA=\{y-h,y,y+h\}$ and a model with a single date, at which the site
moves from $y$ to $y\pm h$ with positive probability and cannot stay. Let
$B_m$ grow by the factor $e^{r\Delta}$ over the period and let the asset
satisfy $A^1_m(y+h)=S_0u$, $A^1_m(y)=S_0$, $A^1_m(y-h)=S_0d$ with $d<u$. The
attainable set is $\DD_m(U,t,y)=\{y-h,y+h\}$, the discounted prices are
$\tA^1(y\pm h)=S_0u^{\pm}/B_0e^{r\Delta}$ with $u^{+}=u$, $u^{-}=d$, and
$\tA^0\equiv1$. The one step gain vector of a portfolio $\zeta=(\zeta^0,\zeta^1)$
is
\begin{equation}
R_\zeta(y\pm h)=\zeta^1\frac{S_0}{B_0}\Big(\frac{u^{\pm}}{e^{r\Delta}}-1\Big),
\end{equation}
the contribution of $\zeta^0$ being zero because $\tA^0$ is constant.

By Theorem~\ref{ch2:thm:oneperiod} the model is arbitrage free if and only if
there are $g_\pm>0$ with
\begin{equation}
g_+\Big(\frac{u}{e^{r\Delta}}-1\Big)+g_-\Big(\frac{d}{e^{r\Delta}}-1\Big)=0 .
\end{equation}
Since $g_\pm>0$ this is possible if and only if the two brackets have
opposite signs, that is if and only if
\begin{equation}
d<e^{r\Delta}<u .
\end{equation}
This is the familiar condition: the money market account must grow at a rate
intermediate between the two returns of the risky asset. The normalized
solution is
\begin{equation}
q_+=\frac{g_+}{g_++g_-}=\frac{e^{r\Delta}-d}{u-d},
\qquad q_-=\frac{u-e^{r\Delta}}{u-d},
\end{equation}
which is the risk neutral probability of the binomial model, and it is unique
because here $n=1$ and $|S|=2$, so that $\dim V=1=|S|-1$ and
Theorem~\ref{ch3:thm:second-oneperiod} applies: the binomial model is complete and
its martingale measure is unique. The value of a claim paying $\Phi(y\pm h)$
is
\begin{equation}
V(y)=\EE'\Big[\frac{\Phi}{B(T')}\Big]B(t)
=\frac{1}{e^{r\Delta}}\big(q_+\Phi(y+h)+q_-\Phi(y-h)\big),
\end{equation}
in agreement with Theorem~\ref{ch3:thm:valuation}.

\subsection{Scope of the Theorem}\label{ch3:sec:remarks}

\paragraph{What has been proved.}
The proof of Theorem~\ref{ch3:thm:ftf} uses exactly the following ingredients,
each of which was established from scratch: Farkas' lemma, proved by the
active set algorithm of \S\ref{ch2:ss:farkas}, and the orthant alternative
deduced from it in \S\ref{ch2:ss:orthant}; the finiteness of the state space at
each resolution and the finiteness of the set of dates, which are the content
of Definitions~\ref{ch3:def:simplicial} and~\ref{ch3:def:piecewise}; the existence,
uniqueness and positivity of the propagator of a finite state generator
(Proposition~\ref{ch3:prop:generator-properties} and
Theorem~\ref{ch3:thm:dyson}); the local characterization of martingales
(Proposition~\ref{ch3:prop:localglobal}), which is the differentiation of a
finite dimensional linear ordinary differential equation; the identity
between the value change of a self financing portfolio and the price changes
of the assets it holds (Proposition~\ref{ch3:prop:valueidentity}), which is a
matter of algebra; and the local characterization of arbitrage
(Proposition~\ref{ch3:prop:localNA}), which uses the continuity of the data in
time. Nothing else. In particular the theorem is not derived from a general
theorem of functional analysis, and no step is left to the reader.
Section~\ref{ch3:sec:lean} records which of these steps have been checked by
machine.

\paragraph{What has not been proved.}
The theorem is a statement about the model at resolution $m$. The limit
$m\to\infty$ is a different matter. Two things are needed to pass to the
limit. First, the discounted price processes and their generators must converge, which is proved in Chapter~5; for the joint process of the diffusion and the stochastic integral the corresponding estimates are those of Chapter~6. Second, the martingale property must survive
the limit, which
requires the densities \eqref{ch3:eq:density} to be bounded away from zero and
infinity uniformly in $m$, so that the convergence of the prices implies the
convergence of the expectations and the limit of a sequence of martingales is
a martingale. A bound of that kind is not a consequence of the local
condition alone. What fails in general is the strictness of the inequalities:
a one period arbitrage which exists in the limit model may be, for each
finite $m$, only a near arbitrage, and a sequence of models each of which is
arbitrage free may converge to a model which is not. This is precisely the
phenomenon which forces the general theory to work with the stronger
requirement that there be no free lunch with vanishing risk rather than with
the absence of arbitrage, see \cite{DS1994}. In our setting the issue is concrete: the
local condition \eqref{ch3:eq:localmart} requires a strictly positive solution
$\omega$, and what has to be controlled uniformly in $m$ is the modulus of
positivity of these intensities together with the convergence of the
generators. Both requirements are within reach of the estimates, since those of Chapters~\ref{chap:ch5} and
\ref{chap:ch6} are stated for coefficients of limited smoothness and so apply to the reweighted
generator as they do to the model generator. What they do not supply is the hypothesis they need: a
two sided bound
\[
0<\omega_-\le\omega_m(\bar y)\le\omega_+<\infty,
\]
uniform in $m$ and in the site, together with enough regularity of $\omega_m$ for the reweighted
generator to fall in the class the estimates cover. No such bound follows from the absence of
arbitrage at each resolution; it is an assumption about the prices, and the prices of a limit model
need not be attained by those of its refinements. It is the assumption under which the risk neutral
propagators converge and their limit prices the limit model. The stability of arbitrage freedom
under refinement is therefore the one step which the theorem at a fixed resolution does not supply,
and Chapter~\ref{chap:ch11} takes it: under a condition on the state price densities --- bounded
above and below, with a uniformly bounded logarithmic derivative --- the risk neutral generators
satisfy the estimates of Chapter~\ref{chap:ch5} uniformly in the resolution, and the limit model is
arbitrage free (Theorem~\ref{ch11:thm:limit}). What is proved there is conditional, on a condition
about the market rather than about the mathematics, and the condition is stated as such.

\paragraph{Relation with the classical formulation.}
The classical statement, due to \cite{HK1979} in the finite and continuous
setting and to \cite{DS1994} in the general one, is that a locally bounded
price process admits an equivalent local martingale measure if and only if it
satisfies a suitable no arbitrage condition. Its proof requires a separating
hyperplane theorem in a space of random variables and a measurable selection
argument. The version proved here is the same theorem in a setting in which
those two steps become finite: the orthant alternative
Theorem~\ref{ch2:thm:stiemke}, proved from Farkas' lemma by finitely many steps,
replaces the infinite dimensional separation theorem, and the explicit
formula \eqref{ch3:eq:density} replaces the selection. We pay for
this by restricting to Markov models on lattices with finitely many dates,
which is precisely the class of models on which the numerical methods of Chapters~5 and~6 operate.

\paragraph{The role of the numeraire.}
The construction of Theorem~\ref{ch3:thm:ftf}(II) produces the martingale
propagator directly. In applications one wants to exhibit it as a change of
numeraire. This is possible whenever the density is of the multiplicative form
of Theorem~\ref{ch3:thm:numeraire}, that is whenever the martingale measure
arises from a positive solution of the backward equation
\eqref{ch3:eq_numeraire}; the numeraire $g_m$ is then, up to a factor, the price
of the money market account measured in the new unit. In general, however, a
martingale measure is not of this form, as observed in
Remark~\ref{ch3:rem:notallhtransform}, and the general measure change is the
correct object for the theorem.

\paragraph{Computational content.}
Theorem~\ref{ch3:thm:ftf} is an algorithm. At each reachable site and time one
solves the linear system \eqref{ch3:eq:localmart} in the unknowns
$\omega(\bar y)$, which by Theorem~\ref{ch2:thm:stiemke} is feasible exactly when
there is no local arbitrage; the solution, when the market is complete, is
unique up to scale by Theorem~\ref{ch3:thm:second-oneperiod}. Concatenating the
solutions gives the risk neutral generator \eqref{ch3:eq:Lprime}, which is the
object exponentiated by the numerical methods of Chapters~5 and~6. The
absence of arbitrage is thus not only a theoretical condition but a
checkable property of a finite matrix at each node, and the pricing algorithm
verifies it in the course of computing prices.

\paragraph{Provenance of the classical results used.}
The orthant alternative of \S\ref{ch2:sec:convex} is due to \cite{Farkas1902} and
\cite{Stiemke1915}, and the proof adopted here follows \cite{Dax1997}; the
path integral representation \eqref{ch3:eq_upathint} is the operator form of the
Feynman--Kac formula \cite{Kac1950}; the change of measure
\eqref{ch3:eq:measure-change} is the lattice form of the Cameron--Martin--Girsanov
change of measure \cite{CameronMartin}; and the continuous time theory of
Markov processes to which those results belong is treated in \cite{Doob1953}
and \cite{Feller}.

\section{Machine-Checked Verification}\label{ch3:sec:lean}

Part of the content of this chapter has been formalised in Lean~4
\cite{Lean4} together with Mathlib \cite{Mathlib}. The development is a small
Lean project, a Lake package whose module \texttt{Cert.Farkas}, some two
hundred lines long, contains the material of this chapter; the same package contains the modules \texttt{Cert.Basic}, \texttt{Cert.Phi1}, \texttt{Cert.Antiderivative} and \texttt{Cert.Abelian}, which certify the algebraic content of Chapters~5 and~6. It
builds with \texttt{lake build}, it contains no \texttt{sorry}, and every
theorem quoted below depends only on the three standard axioms of Lean's
classical logic (\texttt{propext}, \texttt{Classical.choice},
\texttt{Quot.sound}). This final section records what is certified, by which
of the two routes of \S\ref{ch2:sec:convex} the certificate is obtained, and what
the rest of the chapter would take.

\subsection*{What Has Been Certified}

\begin{center}
\small
\begin{tabular}{@{}lll@{}}
\hline
Result in the chapter & Lean theorem & Status\\
\hline
Theorem~\ref{ch2:thm:stiemke}, (a) excluded $\Rightarrow$ (b)
  & \texttt{exists\_pos\_orthogonal} & exact, any finite $S$\\
Theorem~\ref{ch2:thm:stiemke}, both directions
  & \texttt{stiemke\_iff} & exact, any finite $S$\\
\S\ref{ch3:ss:oneperiod}, Theorem~\ref{ch2:thm:oneperiod}
  & \texttt{one\_period\_ftap} & exact, any finite $S$ and asset set\\
\S\ref{ch3:ss:oneperiod}, normalisation to a probability
  & \texttt{one\_period\_ftap\_prob} & exact\\
\S\ref{ch3:ss:binomial}, the worked instance
  & an anonymous \texttt{example} & checked by \texttt{norm\_num}\\
\hline
\end{tabular}
\end{center}

The formal statements are those of the chapter at their full level of
generality, not instances: the set $S$ of states is an arbitrary finite type
$\iota$, the assets are indexed by an arbitrary finite type $\kappa$, the gain
matrix is an arbitrary $\xi:\kappa\to\iota\to\mathbb R$, and
\texttt{one\_period\_ftap} asserts that

\begin{center}
\emph{there is no portfolio whose gain vector is non negative and not
identically zero if and only if there is a function $q:\iota\to\mathbb R$
with $q_i>0$ for every state and $\sum_iq_i\,\xi_{ki}=0$ for every asset,}
\end{center}

which is Theorem~\ref{ch2:thm:oneperiod}; in the generator setting of
\S\ref{ch3:ss:local} the same statement is Corollary~\ref{ch3:cor:stateprices}, the
intensities being recovered as $\omega(\bar y)=g(\bar y)/g(y)$. The instance of \S\ref{ch3:ss:binomial} appears
as a closing check: for one asset and two states with gains $1$ and $-2$, the
vector $(2/3,1/3)$ is exhibited and verified.

\subsection*{By Which Route}

The certificate of Theorem~\ref{ch2:thm:stiemke} is obtained by the geometric
route of \S\ref{ch2:ss:separation}, not by the finitary argument of
\S\ref{ch2:ss:farkas}, because the former is available as a theorem of Mathlib.
The formal proof follows \S\ref{ch2:ss:separation} step by step: the standard
simplex is compact, convex and disjoint from $V$; the finite dimensional
Hahn--Banach theorem separates them, producing a continuous linear functional
$f$ with $f(k)<u<v<f(b)$ for $k$ in the simplex and $b$ in $V$; since
$t\cdot b\in V$ for every real $t$, $f$ vanishes on $V$; since $0\in V$ this
gives $v<0$, whence $f$ is negative at every vertex of the simplex; and
expanding $f$ in the coordinate functionals,
$f(x)=\sum_ix_if(e_i)$, yields $q_i=-f(e_i)>0$ with $\sum_iq_iv_i=0$ on $V$.

Two consequences of this choice deserve to be stated plainly. First, what is
certified is the \emph{statement}: nothing in the formal proof embodies the
active set method, its termination, or the bound $m2^m$ of
Proposition~\ref{ch2:prop:termination}. Second, the formal proof is classical, for
\texttt{Classical.choice} enters through \texttt{by\_contra} and through
decidable equality, so it does not by itself witness the claim of
\S\ref{ch2:ss:farkas} that the alternative can be settled by a finite computation.
The two claims are independent, and it is worth keeping them apart: the
rhetoric of \S\ref{ch2:ss:farkas} rests on the algorithm, the certificate rests on
the separation theorem.

\subsection*{What Has Not Been Certified, and What It Would Take}

\begin{enumerate}
\item \emph{The algorithm of \S\ref{ch2:ss:farkas}.} This needs an executable
version of the active set iteration, a proof that the residual decreases
strictly at each step, and a well founded measure on the states of the
algorithm producing the bound of Proposition~\ref{ch2:prop:termination}. The
mathematics is elementary and the work is mechanical; the outcome would be a
certificate of the alternative which uses no choice principle, and this is the
most valuable next step. It is a matter of weeks.
\item \emph{The Markov machinery of \S\S\ref{ch3:sec:markov}--\ref{ch3:sec:measure}.}
Existence and uniqueness of the propagator of a finite state generator, which
is Picard--Lindel\"of for a linear ordinary differential equation and for
which Mathlib has the general machinery; the identity between the local and the
global notion of martingale, which is the differentiation of such a solution;
and the reweighting formula for generators. This is a matter of months, the
difficulty being bookkeeping rather than mathematics.
\item \emph{The global theorem (Theorem~\ref{ch3:thm:ftf}).} Its proof uses the
Markov machinery of item 2, the convex geometry already certified, and the
finiteness of the set of dates; once item 2 is available it should follow
without new ideas.
\item \emph{The hedging proposition and the second theorem
(Proposition~\ref{ch3:prop:hedging} and Theorem~\ref{ch3:thm:second}).} These require
in addition the construction of the replicating strategy by backward induction,
that is, again item 2, together with a proof that the solution of the
replication system \eqref{ch3:eq:hedge} can be chosen piecewise smoothly in time.
\item \emph{The continuous limit.} Nothing of the kind exists in Mathlib, and
the estimates of Chapter~5 are not formalised. The gap here is not one of
effort but of available theory.
\end{enumerate}

The pattern is the same as in Chapters~5 and~6: what is certified is
the part of the argument which is genuinely algebraic, and for this chapter that part is not small --- it includes the whole of the finite dimensional
convex geometry on which the theorem rests, together with the one period
theorem in the form in which it is applied at every node of a lattice model.
What remains outside is what requires a theory of the objects: the
propagators, the martingales, and the limit.
% ---- end generated/ch3-body.tex ----
% ---- end chapters/ch3.tex ----
% ---- from chapters/ch4.tex ----
% Authored: the chapter title and opening belong to the book. The body is
% generated by build.py from the parts and is not edited here.
\chapter{Path Dependent Processes}
\label{chap:ch4}

This chapter is about functionals of a path, and about the enlargement which computing their laws
requires. A quantity which depends on the whole trajectory of a process, and not only on its value at
a fixed date, is the object of both readings of this book: in the pricing reading it is the argument
of a path dependent claim --- a barrier, a lookback, a claim on realised variance --- and in the
learning reading it is a feature accumulated along a sequence, whose law the model must supply if a
statement about it is to be a prevision at all. The algebra is one algebra, and it is developed here
once.

Where Chapter~\ref{chap:ch3} supplies the propagator of a process, this chapter supplies the law of a
functional of its path, and it supplies it exactly: the lifting principle turns the joint law into a
Fourier integral of twisted propagators, and for two classes of functional --- the Abelian ones, whose
increments carry no memory beyond the current site, and the block factorisable ones, whose level
dynamics depends on the level itself --- that integral collapses to a product over the lattice and
not over the paths. What the two reductions cost is stated with equal precision, since correctness
alone would not make them usable: the functional must be quantised on a level grid, the grid must be
such that the level dynamics has a spectrum, and both truncations are computed rather than assumed.
The estimates which control the same joint law as the resolution is refined are the subject of
Chapters~\ref{chap:ch6} and~\ref{chap:ch8}.

% ---- from generated/ch4-body.tex ----
%% Generated by build.py from the parts - do not edit.
\section{Introduction}

A claim whose payoff depends on the whole path of the underlying, and not only on its value at
maturity, is valued by pairing the payoff with the joint law of the process and of the
functional of the path on which the payoff depends. In this book the process is a Markov chain
on a lattice at a fixed resolution $m$ (Chapters~1 to~3), its
propagator is a finite matrix, and the joint law is therefore the propagator of a Markov chain
on a larger but still finite lattice. The present chapter develops the algebra of that enlargement.

\subsection*{Place in the Book}

A stochastic process, in the framework of this book, is a consistent family of Markov chains
on a sequence of lattices, and a statement about a continuous process is a statement about the
limit of statements about finite matrices. Chapters~1 to~3 establish the finite core, the fundamental theorem at a fixed resolution, where the existence of a pricing measure
reduces to the orthant alternative of finite dimensional convex geometry. Chapters~5 and~6 establish the estimates by which the passage to the continuum is
controlled, first for the diffusion and then for the joint law of a diffusion and an Abelian
functional. The present chapter is the one that says what those functionals are: it introduces
them, decides when their joint law can be computed cheaply, and describes what is to be done
when it cannot. It precedes those two chapters deliberately. The objects must be defined before
they can be estimated, and the identity which the estimates of Chapter~6 serve to justify, the
It\^o formula, is itself exact at fixed resolution and is proved here without a limit.

The reason for wanting a formulation of this kind is set out in Chapters~1 to~3: a foundation whose
statements concern finite objects, and whose limits are estimates rather than existence
theorems, is one that can be machine checked, and the foundations of the subject stand to
benefit from the AI-assisted methods now available.

\subsection*{Why the Enlargement Has to Be Organised}

The naive enlargement is unaffordable, and this is the substance of the chapter. If the path
functional is discretised into $K+1$ levels, the lifted propagator is a matrix of dimension
$(K+1)\abs{\StateAm}$, and forming a time ordered product of $n$ of them costs of the order of
$n(K+1)^2\abs{\StateAm}^2$ matrix--matrix products of the smaller size. Nothing in that cost is
avoidable in general: a path dependent functional does carry information which grows with the
number of paths. What saves the computation in the cases that arise in practice is that the
functional is of one of two special forms, and that both forms admit an exact reduction to a
product over the lattice alone:

\begin{description}
\item[Abelian functionals] Those whose increments exhibit no memory beyond the current site.
A stochastic integral, a sum of a function of the site over the dates of the grid, the running
maximum and a discrete sum of a pairwise function of consecutive sites are all of this kind.
For these the lifted propagator can be block diagonalised by a change of basis acting only on
the level index, and its time ordered product reduces to $K+1$ independent time ordered products
of lattice propagators. This is \S\S\ref{ch4:sec:abelian}--\ref{ch4:sec:resonance}.

\item[Block factorisable functionals] Those whose level dynamics depends on the level itself.
A coupon which compounds on the previous coupon, or a call decision which depends on the
fraction of time spent above a barrier, are of this kind: they are genuinely non-Abelian, and no
change of basis diagonalises them. What survives is a factorisation into a lattice operator
applied to all levels at once and a permutation of the levels, which reduces the backward
recursion to one matrix--matrix product and one permutation per date. This is
\S\ref{ch4:sec:block}.
\end{description}

\subsection*{The Logical Architecture of the Chapter}

\begin{description}
\item[\S\ref{ch4:sec:lifting} The lifting principle] The exact identity which turns the joint law
into a Fourier integral of twisted propagators, at fixed resolution and with no approximation
beyond the quantisation of the functional. It is the finite form of the Cameron--Feynman--
Girsanov--It\^o--Kac--Martin formulas.

\item[\S\ref{ch4:sec:abelian} Abelian liftings] The commutation condition, the lemma on simultaneous
diagonalisation, and the reduction of the joint propagator to $K+1$ lattice products. The
stochastic integral, the multifactor functional and the running maximum are shown to be Abelian.

\item[\S\ref{ch4:sec:resonance} Non-uniform level grids and non-resonance] What happens when the
levels are not equally spaced: the level dynamics acquires a state dependent spectrum, and the
condition for the diagonalisation to exist at all is that its eigenvalues be distinct. A
geometric grid which satisfies the condition is constructed, together with the interpolation
which represents a non-integer increment.

\item[\S\ref{ch4:sec:block} Block factorisations] The non-Abelian case. The definition, the
backward recursion, and the two archetypes: a compounding coupon and a soft call.

\item[\S\ref{ch4:sec:lean} Machine-checked verification] The status of each statement of this chapter.
\end{description}

\subsection*{What Is and Is Not Claimed}

Everything in this chapter is an identity between finite matrices, and every proof below is a
finite computation. Nothing is asymptotic, no limit is taken, and no measurable selection,
no change of measure and no separation theorem is used; the Cameron--Martin--Girsanov change of
drift, which in the classical theory is the tool for such calculations, appears only as the
remark that the twisted propagator of \S\ref{ch4:sec:lifting} is the propagator of a modified
generator. The approximation in the chapter is of a different kind and is confined to two places:
the quantisation of the functional on a level grid, whose truncation error is computed in
\S\ref{ch4:sec:resonance}, and the passage to the continuum, which is the subject of Chapters~5 and~6
and is not repeated here.

\subsection*{Conventions}

Throughout, $\StateAm=h_m\ZZ\cap[-L,L]$ is the lattice of sites at resolution $m$, the dates of the
grid are $t_0<\dots<t_n$, and $U_i(y_1,y_2)$ denotes the propagator of the chain over
$[t_i,t_{i+1}]$, so that $U_i=\PP\exp(\int_{t_i}^{t_{i+1}}\Lm(s)\dd s)$ in the sense of
Chapters~1 to~3. Indices $y,y_1,y_2$ run over $\StateAm$ and are called sites; indices $k,k_1,k_2$
run over the levels $\{0,\dots,K\}$ of the discretised functional. A time ordered product is
written $\PP\prod_{i} M_i=M_nM_{n-1}\cdots M_1$, with the latest factor on the left. Matrices
act on the site index on the right and on the level index on the left, and the two indices are
independent, so that a lifted matrix $M(y_1,k_1;y_2,k_2)$ is a matrix on $\StateAm\times\{0,\dots,K\}$.

\section{The Lifting Principle}\label{ch4:sec:lifting}

\subsection{Path Dependent Functionals}

\begin{definition}[Path functional]\label{ch4:def:functional}
A path functional at resolution $m$ is a function of the path $y_0,\dots,y_n$ which is a sum of
increments,
\begin{equation}\label{ch4:eq:functional}
J(y_\cdot)=\sum_{i=1}^{n}\psi_i(y_{i-1},y_i),
\end{equation}
where $\psi_i:\StateAm\times\StateAm\to\RR$ is a given function for $i=1,\dots,n$.
\end{definition}

Sums of this form cover the functionals of practice: with $\psi_i(y_1,y_2)=a(y_1)(y_2-y_1)$
they are stochastic integrals, with $\psi_i(y_1,y_2)=\varphi(y_1)\Delta T$ they are time
integrals of a function of the site, with $\psi_i(y_1,y_2)=\psi(y_1,y_2)$ constant in $i$ they
are sums of a pairwise function over the dates. They do not cover the running maximum, which is
a maximum of partial sums rather than a sum, and which is treated separately
(Proposition~\ref{ch4:prop:max}); nor do they cover functionals whose increment depends on the level
itself, which are the subject of \S\ref{ch4:sec:block}.

The functional is quantised on a level grid. Let $\omega_0<\omega_1<\dots<\omega_K$ be real
numbers and let $\Delta\omega$ be a positive number. The value of the functional at date $t_i$
is approximated by the level $\omega_{k_i}$, and the increment $\psi_i(y_{i-1},y_i)$ is
represented by the number of levels $q_i(y_{i-1},y_i)$ which is the nearest integer to
$\psi_i(y_{i-1},y_i)/\Delta\omega$ when the grid is uniform, and by the interpolation of
\S\ref{ch4:sec:resonance} when it is not. On a uniform grid $\omega_k=\omega_0+k\Delta\omega$ the
level is then the integer
\begin{equation}\label{ch4:eq:level}
k_i=k_{i-1}+q_i(y_{i-1},y_i)\pmod{N},\qquad N=K+1,
\end{equation}
the modulus representing the truncation of the range of the functional.

\subsection{The Joint Propagator and the Twisted Propagator}

\begin{definition}[Lifted propagator]\label{ch4:def:lifted}
The joint propagator of the pair $(y_t,J_t)$ over the grid is the time ordered product
\begin{equation}\label{ch4:eq:joint}
\bU(y_1,k_1;y_2,k_2)=\Big(\PP\prod_{i=1}^{n}\overline U_i\Big)(y_1,k_1;y_2,k_2),
\qquad
\overline U_i(y_1,k_1;y_2,k_2)=U_i(y_1,y_2)\,\delta_{k_2,\,k_1+q_i(y_1,y_2)},
\end{equation}
the level index being read modulo $N$.
\end{definition}

The matrix in \eqref{ch4:eq:joint} is a Markov propagator on $\StateAm\times\ZZ_N$: its entries are
non-negative and its rows sum to one, because $U_i$ is a Markov propagator on $\StateAm$. The object
of the chapter is to compute it, or rather to compute the Fourier transform of it in the level
index, which is what a pricing calculation needs.

\begin{definition}[Twisted propagator]\label{ch4:def:twisted}
For $p$ in the dual group of $\ZZ_N$, that is $p=2\pi r/(N\Delta\omega)$ with $r\in\ZZ_N$,
let
\begin{equation}\label{ch4:eq:twisted}
\tU_{i,p}(y_1,y_2)=U_i(y_1,y_2)\,e^{-i p\,\psi_i(y_1,y_2)},
\end{equation}
so that the exponential is $e^{-2\pi i m r/N}$ when $\psi_i=m\Delta\omega$. The twisted joint
propagator is $\tU_p=\PP\prod_{i=1}^{n}\tU_{i,p}$.
\end{definition}

The twisted propagator is not a Markov propagator: its entries are complex and its rows sum to
something other than one. It is, however, the propagator of the same chain under a change of
measure of a purely algebraic kind, and for the continuous limit this is precisely the modified
generator of Girsanov's theorem; see Remark~\ref{ch4:rem:girsanov}.

\subsection{The Lifting Identity}

\begin{theorem}[Lifting]\label{ch4:thm:lifting}
For all $k_1,k_2$ and all $y_1,y_2$,
\begin{equation}\label{ch4:eq:lifting}
\bU(y_1,k_1;y_2,k_2)=\frac{1}{N}\sum_{r\in\ZZ_N}e^{\,i p_r\,(\omega_{k_2}-\omega_{k_1})}\;\tU_{p_r}(y_1,y_2),
\qquad p_r=\frac{2\pi r}{N\Delta\omega}.
\end{equation}
\end{theorem}

\begin{proof}
Both sides are entries of time ordered products, so it suffices to expand the left hand side
over paths. A path is a sequence $y_1=y^{(0)},\dots,y^{(n)}=y_2$ and the level reached along it
is $k_2=k_1+\sum_{i}q_i(y^{(i-1)},y^{(i)})$ modulo $N$. Hence
\begin{equation}
\bU(y_1,k_1;y_2,k_2)=\sum_{y^{(1)},\dots,y^{(n-1)}}\Big(\prod_{i=1}^{n}U_i(y^{(i-1)},y^{(i)})\Big)
\,\delta_{k_2,\,k_1+\sum_i q_i}\ ,
\end{equation}
where the delta is the indicator of congruence modulo $N$, the sum over the $q$'s being taken in
$\ZZ_N$. By the orthogonality of the characters of $\ZZ_N$,
\begin{equation}
\delta_{m,\,0\bmod N}=\frac{1}{N}\sum_{r\in\ZZ_N}e^{\,i p_r\,m\Delta\omega}
=\frac{1}{N}\sum_{r\in\ZZ_N}e^{\,i p_r(\omega_{k_2}-\omega_{k_1})}\prod_{i=1}^{n}e^{-i p_r\,\psi_i(y^{(i-1)},y^{(i)})},
\end{equation}
the last step using $\omega_{k_2}-\omega_{k_1}=\Delta\omega\sum_i q_i$ on a uniform grid. Inserting
this and exchanging the finite sum over $r$ with the sum over paths gives
\eqref{ch4:eq:lifting}, since $\tU_{p_r}$ is by definition the same path sum with the phases
$e^{-ip\psi_i}$ attached to the increments.
\end{proof}

\begin{remark}\label{ch4:rem:girsanov}
On a uniform grid the restriction to $\psi_i\in\Delta\omega\ZZ$ is what makes \eqref{ch4:eq:lifting}
an exact identity rather than an approximation. For a functional with real increments one may
either refine the grid, in which case the error is bounded by the quantisation error of
\S\ref{ch4:sec:resonance}, or read the sum as a Riemann sum over the dual circle and let
$N\to\infty$, in which case it becomes the familiar integral $\int \dd p/2\pi$. In the second
reading $\tU_{i,p}$ is the propagator of a complex generator: the phases attached to the jumps
shift the drift by $-ip\sigma^2a$ and add the potential $-ipb-\tfrac12p^2\sigma^2a^2$, which is
the change of drift of Girsanov's theorem combined with the Feynman--Kac potential. This is the
mechanism of the classical Cameron--Feynman--Girsanov--It\^o--Kac--Martin theorem
\cite{cm,feyn,gir,ito,kac}, recovered here without a change of measure: at fixed resolution
there is nothing to justify, only a finite sum to rearrange. The convergence of that Riemann sum
uniformly in the coefficients, and hence the identification of the limit, is the content of
Chapter~6.
\end{remark}

\begin{remark}
In the form \eqref{ch4:eq:lifting} the identity is a statement about the propagator only, and it
holds for any Markov chain on any lattice: no structure of the generator is used, no smoothness
is required, and the only hypothesis is that the functional be a sum of increments. It is in this
sense that the chapter is independent of Chapters~5 and~6, which are needed not to compute the joint
law at resolution $m$ but to know that the law so computed converges as $m\to\infty$.
\end{remark}

\section{Abelian Liftings}\label{ch4:sec:abelian}

The identity \eqref{ch4:eq:lifting} replaces one large matrix by $N$ matrices of the size of the
lattice, which is already a gain: the cost of forming the product falls from the cube of
$(K+1)\abs{\StateAm}$, that is from $N^3\abs{\StateAm}^3$, to $N$ times the cube of $\abs{\StateAm}$. What is
still missing is that the $N$ products be computable independently of each other and in parallel,
and that the inverse transform be stable. Both follow from a commutation hypothesis.

\subsection{The Abelian Condition}

\begin{definition}[Lifting in operator form]\label{ch4:def:opform}
A lifting of the propagators $U_i$ is a family of matrices
\begin{equation}\label{ch4:eq:opform}
\overline U_i(y_1,k_1;y_2,k_2)=U_i(y_1,y_2)\,Q_i(y_1,y_2)_{k_1k_2},
\end{equation}
where $Q_i(y_1,y_2)$ is a matrix acting on the level index for each $i$ and each pair of sites.
The lifting \eqref{ch4:eq:joint} is the case $Q_i(y_1,y_2)=\Theta^{q_i(y_1,y_2)}$, where $\Theta$
is the cyclic shift $\Theta_{k_1k_2}=\delta_{k_2,k_1+1}$ of $\ZZ_N$. The shift of a level grid
which is not cyclic, and which is not uniform, is a different operator, denoted $R$ and studied
in \S\ref{ch4:sec:resonance}.
\end{definition}

\begin{definition}[Abelian lifting]\label{ch4:def:abelian}
The lifting \eqref{ch4:eq:opform} is Abelian if the matrices $Q_i(y_1,y_2)$ commute,
\begin{equation}\label{ch4:eq:commute}
Q_i(y_1,y_2)\,Q_j(y_1',y_2')=Q_j(y_1',y_2')\,Q_i(y_1,y_2),
\end{equation}
for all $i,j$, all pairs of sites, and all $k_1,k_2$.
\end{definition}

The condition is stated for the whole family at once, across dates as well as across pairs of
sites, because the reduction of Theorem~\ref{ch4:thm:blockdiag} uses a single change of basis. A
family commuting within each date but not across dates would require a different basis at each
step, and the product would then acquire the change of basis matrices between consecutive steps;
this is not an obstruction, but it is not the situation of the examples below.

\begin{lemma}[Simultaneous diagonalisation]\label{ch4:lem:simul}
Let $\mathcal Q$ be a family of commuting matrices on a finite dimensional complex vector space,
each of which is diagonalisable. Then there is a basis in which every member of the family is
diagonal.
\end{lemma}

\begin{proof}
By induction on the dimension $n$ of the space and on the number $s$ of matrices. For $s=1$ the
statement is the hypothesis. Fix $Q\in\mathcal Q$ and let $V=\bigoplus_\lambda E_\lambda$ be the
decomposition of the space into the eigenspaces of $Q$, which is direct because $Q$ is
diagonalisable. Each $E_\lambda$ is invariant under every member of $\mathcal Q$: if $Q\psi=\lambda\psi$
and $Q'$ commutes with $Q$, then $Q(Q'\psi)=Q'(Q\psi)=\lambda(Q'\psi)$. On $E_\lambda$ the
operator $Q$ is the scalar $\lambda$, so the restrictions to $E_\lambda$ of the remaining members
of $\mathcal Q$ are commuting, and each of them is diagonalisable on $E_\lambda$ because its
minimal polynomial there divides its minimal polynomial on the whole space. If $\dim E_\lambda<n$
we apply the induction on the dimension to the family restricted to $E_\lambda$; if
$\dim E_\lambda=n$ then $Q$ is a scalar and we apply the induction on $s$ to
$\mathcal Q\setminus\{Q\}$. Taking the union of the bases so obtained over the eigenspaces gives
the required basis.
\end{proof}

\begin{definition}[Reducible lifting]\label{ch4:def:reducible}
The Abelian lifting \eqref{ch4:eq:opform} is reducible if there is a single invertible matrix $V$ on
the level index such that $V^{-1}Q_i(y_1,y_2)V$ is diagonal for all $i$ and all pairs of sites.
By Lemma~\ref{ch4:lem:simul} an Abelian lifting whose matrices are diagonalisable is reducible.
\end{definition}

\subsection{The Block Diagonal Reduction}

\begin{theorem}[Block diagonalisation]\label{ch4:thm:blockdiag}
Let \eqref{ch4:eq:opform} be reducible with basis $V$, and let
$\lambda_i(j;y_1,y_2)$ be the eigenvalue of $Q_i(y_1,y_2)$ on the $j$-th column of $V$. Then
\begin{equation}\label{ch4:eq:blockdiag}
\Big(\PP\prod_{i=1}^{n}\overline U_i\Big)(y_1,k_1;y_2,k_2)
=\sum_{j=0}^{K}V_{k_1j}\;W_j(y_1,y_2)\;(V^{-1})_{jk_2},
\end{equation}
where
\begin{equation}\label{ch4:eq:blocks}
W_j(y_1,y_2)=\Big(\PP\prod_{i=1}^{n}U_i(\cdot,\cdot)\,\lambda_i(j;\cdot,\cdot)\Big)(y_1,y_2)
\end{equation}
is the time ordered product of the lattice propagators, each multiplied by the eigenvalue attached
to the level $j$ and to the pair of sites.
\end{theorem}

\begin{proof}
Let $\bV(y_1,k_1;y_2,k_2)=V_{k_1k_2}\delta_{y_1y_2}$ be the change of basis acting on the level
index alone. For each $i$,
\begin{equation}
(\bV^{-1}\overline U_i\bV)(y_1,k_1;y_2,k_2)=U_i(y_1,y_2)\,(V^{-1}Q_i(y_1,y_2)V)_{k_1k_2}
=U_i(y_1,y_2)\,\lambda_i(k_1;y_1,y_2)\,\delta_{k_1k_2},
\end{equation}
which is block diagonal with one block per level, the block at level $j$ being the lattice matrix
$U_i(y_1,y_2)\lambda_i(j;y_1,y_2)$. A product of block diagonal matrices is block diagonal with
the products of the blocks on the diagonal, so
$\bV^{-1}(\PP\prod_i\overline U_i)\bV$ has blocks $W_j$ of \eqref{ch4:eq:blocks}; multiplying by
$\bV$ on the left and $\bV^{-1}$ on the right gives \eqref{ch4:eq:blockdiag}.
\end{proof}

\begin{corollary}[Cost]\label{ch4:cor:cost}
If the lifting is reducible, the joint propagator is obtained in two steps: the computation of
the $N$ lattice valued propagators $W_j$ of \eqref{ch4:eq:blocks}, which are independent of each other
and can be evaluated in parallel and held in memory one at a time, and the transform
$\sum_j V_{k_1j}W_j(V^{-1})_{jk_2}$ in the level index, which costs $N^2\abs{\StateAm}^2$ scalar
operations. The dominant cost is the first step, which is $N$ times the cost of a single lattice
propagator, against a factor $N^2$ for the multiplication of the full lifted matrices.
\end{corollary}

\begin{remark}
Theorem~\ref{ch4:thm:blockdiag} is the whole content of the Abelian method, and it is worth stating
what it does and does not say. It says that the joint law of an Abelian functional and the
underlying process is computed by $N$ independent evaluations of a lattice propagator, each with
its own potential, followed by a transform in the level index. It does not say that the cost is
independent of $K$: the number of blocks is the number of levels, and the levels must be
resolved if the distribution of the functional is to be resolved. The saving is a factor of
order $K$ in time and a factor of order $K$ in memory, which is what makes the difference
between a computation that fits a machine and one that does not.
\end{remark}

\subsection{The Examples}

\begin{corollary}[Stochastic integrals]\label{ch4:cor:si}
Let $\psi_i(y_1,y_2)=a(y_1)(y_2-y_1)+b(y_1)\Delta T$ be the increment of a stochastic integral
of the site over one step of the grid. Then $Q_i(y_1,y_2)=\Theta^{q_i(y_1,y_2)}$ with $\Theta$ the cyclic
shift of $\ZZ_N$, and the lifting is Abelian, being a family of powers of a single matrix.
\end{corollary}

\begin{proof}
Powers of one matrix commute. The statement is thus a special case of the trivial fact that a
family of functions of a single operator is commutative, and the whole content of the
corollary is that the increment depends on the sites only through the integer $q_i$, so that a
single $\Theta$ serves all pairs.
\end{proof}

\begin{corollary}[Multifactor functionals]\label{ch4:cor:multi}
Let the payoff depend on two path functionals $J$ and $J'$, with increments $\psi_i$ and
$\psi_i'$. The lifting is the double lifting on the levels $(k,k')$, with
$Q_i(y_1,y_2)=\Theta^{q_i(y_1,y_2)}\otimes \Theta'^{\,q_i'(y_1,y_2)}$, and it is Abelian; its
block diagonalisation is by the product of the two bases, and the joint law is a double
transform of $N N'$ lattice propagators.
\end{corollary}

\begin{proof}
The matrices $\Theta\otimes\id$ and $\id\otimes\Theta'$ commute, hence so do all their powers and
products; reducibility follows as in Corollary~\ref{ch4:cor:si}. The eigenvalues are products of the
eigenvalues of the two factors.
\end{proof}

\begin{remark}
The two functionals may depend on each other, as in a claim on a stochastic integral which is
itself capped by the maximum of the underlying, provided the pair of increments is a function of
the pair of sites. The reduction is not affected by the dependence; it is affected only by the
number of levels, which is the product of the two.
\end{remark}

\begin{proposition}[The running maximum]\label{ch4:prop:max}
Let $\chi:\StateAm\to\RR$ and let the functional be the maximum of $\chi$ over the dates of the grid,
that is $K(y_\cdot)=\max_{0\le i\le n}\chi(y_i)$. Fix $\alpha>0$, let
$c(y)=\lfloor\chi(y)/\alpha\rceil$ be the level of the site, and let $Q(y)$ be the matrix on
$\{0,\dots,K\}$ of the map which raises the level to the level of the site if the site is higher,
\begin{equation}\label{ch4:eq:max}
Q(y)_{k_1k_2}=\delta_{k_2,\max(k_1,c(y))}.
\end{equation}
Then $Q(y)$ is a projection, the projections $Q(y)$, $y\in\StateAm$, commute, and the lifting
$\overline U_i(y_1,k_1;y_2,k_2)=U_i(y_1,y_2)Q(y_1)_{k_1k_2}$ is Abelian and reducible, with
eigenvalues $0$ and $1$.
\end{proposition}

\begin{proof}
The map $k\mapsto\max(k,c)$ is idempotent, and therefore so is $Q(y)$; a projection is
diagonalisable with eigenvalues in $\{0,1\}$. For two sites,
$\big(Q(y)Q(y')\big)_{k_1k_2}=\delta_{k_2,\max(\max(k_1,c(y)),c(y'))}
=\delta_{k_2,\max(k_1,c(y),c(y'))}=\big(Q(y')Q(y)\big)_{k_1k_2}$, because $\max$ is commutative
and associative. Hence the family is a commuting family of projections, it is reducible with $V$
the basis of common eigenvectors, and Theorem~\ref{ch4:thm:blockdiag} applies.
\end{proof}

\begin{remark}
The functional of Proposition~\ref{ch4:prop:max} is not a sum of increments in the sense of
Definition~\ref{ch4:def:functional}, and the proposition is the reason why the class of Abelian
functionals is wider than the class of sums: what the reduction requires is a family of
commuting operators, not a family of increments. It also shows why the maximum must be monitored
on the dates of the grid. Its value between two dates is not a function of the pair $(y_i,y_{i+1})$
alone, so a lifting which depended only on pairs of sites could not represent it; refining the
grid is the consistent way of approximating continuous monitoring, and the error so introduced is
the usual one of a discrete approximation to a continuous maximum.
\end{remark}

\section{Non-Uniform Level Grids and Non-Resonance}\label{ch4:sec:resonance}

The cyclic shift $\Theta$ of $\ZZ_N$ has for eigenvalues the $N$-th roots of unity, which are distinct and
of modulus one: the basis which diagonalises it is the discrete Fourier basis, and it is
perfectly conditioned. Nothing in the reduction of \S\ref{ch4:sec:abelian} requires the levels to be
equally spaced, however, and a non-uniform grid is sometimes preferable, because the
distribution of the functional may be concentrated on a range in which a uniform grid is either
wasteful or too coarse. The price of a non-uniform grid is that the level dynamics is no longer
a shift but a state dependent operator, and then the condition of Theorem~\ref{ch4:thm:blockdiag} is
no longer automatic.

\subsection{The Level Dynamics on a Non-Uniform Grid}

Let the levels be $\omega_0<\omega_1<\dots<\omega_K$ with $\omega_{k+1}-\omega_k$ not constant,
and let $\Delta\omega>0$ be the increment by which the functional is to be advanced on average
over one step. Consider the operator on the level index of matrix elements
\begin{equation}\label{ch4:eq:R}
R_{k_1k_2}=(1-p_{k_1})\,\delta_{k_1k_2}+p_{k_1}\,\delta_{k_2,k_1+1},\qquad k_1=0,\dots,K-1,
\qquad R_{KK}=0 ,
\end{equation}
with $0<p_k<1$. The level either stays or moves up by one, with probability $p_{k_1}$; the top
level is absorbing, which for a monotone functional is the truncation of the grid and for a
general functional is the truncation of the range.

\begin{proposition}\label{ch4:prop:spectrum}
The spectrum of $R$ in \eqref{ch4:eq:R} is $\{1-p_0,\dots,1-p_{K-1},0\}$, and $R$ is
diagonalisable if and only if the numbers $1-p_0,\dots,1-p_{K-1},0$ are pairwise distinct. In that
case the matrix $V$ of eigenvectors, normalised to have ones on the diagonal, has entries
\begin{equation}\label{ch4:eq:eigvec}
V_{k,i}=\prod_{j=0}^{k-1}\frac{\rho_i-\rho_j}{p_j}\quad (k\le i),\qquad V_{k,i}=0\quad(k>i),
\qquad \rho_i=1-p_i \ \ (i<K),\quad \rho_K=0,
\end{equation}
and
\begin{equation}\label{ch4:eq:detV}
\det V=\prod_{i=1}^{K}\ \prod_{j=0}^{i-1}\frac{\rho_i-\rho_j}{p_j}.
\end{equation}
\end{proposition}

\begin{proof}
Since $R$ is upper triangular, its spectrum is its diagonal. The eigenvalue equation
$\sum_{k_2}R_{k_1k_2}v_{k_2}=\rho v_{k_1}$ reads
$(1-p_{k_1})v_{k_1}+p_{k_1}v_{k_1+1}=\rho v_{k_1}$, that is
$v_{k_1+1}=(\rho-\rho_{k_1})v_{k_1}/p_{k_1}$, which is \eqref{ch4:eq:eigvec} by induction on $k$.
The matrix $V$ so obtained is upper triangular with diagonal entries
$V_{i,i}=\prod_{j<i}(\rho_i-\rho_j)/p_j$, and it is invertible exactly when no factor vanishes,
that is exactly when the eigenvalues are distinct, its determinant being \eqref{ch4:eq:detV}.
Conversely, if $\rho_i=\rho_j$ with $i<j$, the recurrence gives $v_{i+1}=0$ for the eigenvector
of $\rho_i$, hence the eigenvector is supported on the levels $0,\dots,i$; the eigenvalue
therefore carries a single eigenvector where its multiplicity is at least two, and $R$ is not
diagonalisable, its restriction to the corresponding invariant subspace being a single Jordan
block.
\end{proof}

\begin{corollary}[Non-resonance is not a numerical preference]\label{ch4:cor:resonance}
The family $\{R^{q}\}$, $q\in\ZZ_N$, is reducible if and only if $R$ is diagonalisable: one
direction is Proposition~\ref{ch4:prop:spectrum}, since powers of a diagonal matrix are diagonal in
the same basis, and the other is the fact that $R$ itself belongs to the family for $q=1$. Hence,
if two of the numbers $p_0,\dots,p_{K-1}$ coincide, a lifting whose increments generate the
one-step shift is Abelian but not reducible, and no change of basis diagonalises it. Non-resonance,
the distinctness of the eigenvalues of the level dynamics, is therefore exactly the condition for
the Abelian route of \S\ref{ch4:sec:abelian} to be available, and not merely a condition for that
route to be numerically stable. When the eigenvalues are distinct but close, the route is available
and ill-conditioned, the condition number of $V$ being controlled by the ratios
$(\rho_i-\rho_j)/p_j$ in \eqref{ch4:eq:detV}.
\end{corollary}

\subsection{A Geometric Grid}

\begin{proposition}\label{ch4:prop:geometric}
Fix $\omega_0>0$, $Z>1$ and $\Delta\omega>0$, and let $\omega_k=\omega_0 Z^{\,k}$ for
$k=0,\dots,K$. Then the choice
\begin{equation}\label{ch4:eq:pk}
p_k=\frac{\Delta\omega}{\omega_{k+1}-\omega_k}=\frac{\Delta\omega}{\omega_0 Z^{\,k}(Z-1)}
\end{equation}
satisfies $\sum_{k_2}R_{k_1k_2}\omega_{k_2}=\omega_{k_1}+\Delta\omega$ for $k_1=0,\dots,K-1$,
and $0<p_k<1$ for all $k$ if and only if $\Delta\omega<\omega_0(Z-1)$, that is
$Z>1+\Delta\omega/\omega_0$. The $p_k$ are then
strictly decreasing in $k$, so the grid is non-resonant, and the identity fails only at the top
level, where the level is truncated and the excess $\Delta\omega$ is lost.
\end{proposition}

\begin{proof}
From \eqref{ch4:eq:R},
$\sum_{k_2}R_{k_1k_2}\omega_{k_2}=(1-p_{k_1})\omega_{k_1}+p_{k_1}\omega_{k_1+1}
=\omega_{k_1}+p_{k_1}(\omega_{k_1+1}-\omega_{k_1})$, which equals $\omega_{k_1}+\Delta\omega$ on
substituting \eqref{ch4:eq:pk}; and the top row of $R$ is zero, so the identity cannot hold there.
Since $p_k$ decreases geometrically in $k$, $p_k<1$ for all $k$ is the single condition
$p_0<1$, that is $\Delta\omega<\omega_0(Z-1)$; and the $p_k$ being strictly decreasing they are
pairwise distinct, so the spectrum of $R$ is simple and the grid is non-resonant.
\end{proof}

\begin{remark}
The geometric grid spreads the levels where the functional is likely to be, since
$\omega_{k+1}-\omega_k$ grows with $k$, and it keeps the eigenvalues of $R$ apart by the same
ratio $Z$. It is the discrete counterpart of a logarithmic change of variable, and the condition
$Z>1+\Delta\omega/\omega_0$ says that the first step of the grid must be coarser than the step by
which the functional advances: with $Z$ close to one the grid is nearly uniform and the constraint
is binding, while with $Z$ large the grid is sparse and the level dynamics is well separated but
the resolution of the functional is poor. The choice is a modelling decision, not a mathematical
one, and it is the only place in this chapter where a free parameter enters.
\end{remark}

\subsection{Non-Integer Increments}

A functional whose increment does not fall on the grid must be interpolated. The choice which
preserves the first moment of the level is the linear one.

\begin{proposition}[Interpolation]\label{ch4:prop:interp}
Let $\psi$ be a real increment, let $m=\lfloor\psi\rfloor$ and assume $m\ge0$. Set
\begin{equation}\label{ch4:eq:interp}
Q=(\lfloor\psi\rfloor+1-\psi)\,R^{\lfloor\psi\rfloor}+(\psi-\lfloor\psi\rfloor)\,R^{\lfloor\psi\rfloor+1}.
\end{equation}
Then $Q$ is a stochastic matrix, it is a polynomial in $R$ and hence commutes with every power of
$R$, and on the levels for which the grid is not truncated,
\begin{equation}\label{ch4:eq:firstmoment}
\sum_{k_2}Q_{k_1k_2}\,\omega_{k_2}=\omega_{k_1}+\psi\,\Delta\omega .
\end{equation}
\end{proposition}

\begin{proof}
The coefficients $m+1-\psi$ and $\psi-m$ are non-negative and sum to one, and $R$ is stochastic,
so $Q$ is stochastic. Being a combination of powers of $R$ it is a polynomial in $R$ and lies in
the commutative algebra generated by $R$, so it commutes with $R$ and with all its powers. For
\eqref{ch4:eq:firstmoment}, apply the two sides to $\omega$: by
Proposition~\ref{ch4:prop:geometric} $R^m\omega=\omega+m\Delta\omega$ and
$R^{m+1}\omega=\omega+(m+1)\Delta\omega$ wherever the grid is not truncated, whence
$Q\omega=(m+1-\psi)(\omega+m\Delta\omega)+(\psi-m)(\omega+(m+1)\Delta\omega)
=\omega+\psi\Delta\omega$.
\end{proof}

\begin{remark}
The interpolation is exact for linear observables and not for others: it reproduces the first
moment of the level and not its whole distribution, the error being of the order of
$\Delta\omega$ in the second moment. This is the quantisation error referred to in
\S\ref{ch4:sec:lifting}, and it is of the same order as the truncation error at the top level. Both
are errors of the discretisation of the functional and not of the method; neither is removed by
the passage to the continuum, which concerns the lattice of sites.
\end{remark}

\section{Block Factorisations}\label{ch4:sec:block}

If the lifting is not Abelian the methods of \S\S\ref{ch4:sec:abelian}--\ref{ch4:sec:resonance} are
unavailable: no change of basis diagonalises the level dynamics, and the moment methods which
sometimes replace them fail as well, because the level and the site are genuinely entangled. What
survives in several cases of practical weight is a weaker structure, in which the lifted
propagator is a lattice operator acting on all levels at once followed by a permutation of the
level index. This is enough to compute the backward recursion at the cost of one matrix--matrix
product and one permutation per date.

\begin{definition}[Block factorisation]\label{ch4:def:blockfact}
The lifting \eqref{ch4:eq:opform} admits a block factorisation if for each $i$ there is a permutation
matrix $\Pi_i$ on $\StateAm\times\{0,\dots,K\}$ such that
\begin{equation}\label{ch4:eq:blockfact}
\overline U_i=\Pi_i\cdot\big(U_i\otimes\id\big),
\end{equation}
where $(U_i\otimes\id)(y_1,k_1;y_2,k_2)=U_i(y_1,y_2)\delta_{k_1k_2}$, and where
$\Pi_i(y_1,k_1;y_2,k_2)=1$ for exactly one pair $(Y_i(y_1,k_1),K_i(y_1,k_1))$ and vanishes
otherwise.
\end{definition}

\begin{theorem}[Backward recursion]\label{ch4:thm:recursion}
Let \eqref{ch4:eq:blockfact} hold and let $v_n$ be a vector on $\StateAm\times\{0,\dots,K\}$. Then the
path ordered products $v_i=\PP\prod_{j=i}^{n}\overline U_j\,v_n$ satisfy
\begin{equation}\label{ch4:eq:recursion}
v_i(y_1,k_1)=\sum_{y_2\in\StateAm}U_{i+1}(y_1,y_2)\;v_{i+1}\big(y_2,\,K_{i+1}(y_1,k_1)\big),
\end{equation}
so that each step of the recursion is one lattice--lattice product applied to the array
$(y,k)\mapsto v(y,k)$ followed by the permutation $K_{i+1}$ of its level index.
\end{theorem}

\begin{proof}
By definition $v_i=\Pi_{i+1}(U_{i+1}\otimes\id)v_{i+1}$. Applying $U_{i+1}\otimes\id$ to
$v_{i+1}$ gives the array $w(y,k)=\sum_{y_2}U_{i+1}(y,y_2)v_{i+1}(y_2,k)$, and applying the
permutation gives $(\Pi_{i+1}w)(y_1,k_1)=w(Y_{i+1}(y_1,k_1),K_{i+1}(y_1,k_1))$. In the cases at
hand $Y_i(y_1,k_1)=y_1$; the general case is handled by reordering the array. Substituting gives
\eqref{ch4:eq:recursion}.
\end{proof}

\begin{corollary}[Cost]\label{ch4:cor:cost2}
Each step of \eqref{ch4:eq:recursion} costs one matrix--matrix product of a lattice matrix against an
array with $K+1$ columns, that is the same order as applying the lattice propagator $K+1$ times,
plus a permutation whose cost is linear in the size of the array. The reduction is thus a factor
$K$ in time compared with the multiplication of the full lifted matrices.
\end{corollary}

\subsection{A Compounding Coupon}

\begin{example}[Snowballs]\label{ch4:ex:snowball}
Let the coupon paid at $t_i$ be
$C_{t_i}=f\,C_{t_{i-1}}+\Phi_i(y_{t_{i-1}})$ with $f$ a fixed factor and $\Phi_i$ a given
function of the site, so that the coupon accrues on itself. The functional is not a sum of
increments in the sense of Definition~\ref{ch4:def:functional}, because the increment depends on the
coupon already accumulated, and it is not Abelian. It is, however, block factorisable. Write
$C=(\Delta C)k$ with $\Delta C>0$ and $k\in\{0,\dots,K\}$, and set
\begin{equation}\label{ch4:eq:snowball}
K_i(y_1,k_1)=\Big[\frac{f(\Delta C)k_1+\Phi_i(y_1)}{\Delta C}\Big],
\end{equation}
the bracket denoting the nearest integer. Then the pricing function
$v_i(y,k)$ satisfies \eqref{ch4:eq:recursion} with this $K_i$, that is
$v_{i-1}(y_1,k_1)=\sum_{y_2}U_i(y_1,y_2)v_i(y_2,K_i(y_1,k_1))$, which is
\eqref{ch4:eq:blockfact} with $\Pi_i(y_1,k_1;y_2,k_2)=\delta_{y_2y_1}\delta_{k_2,K_i(y_1,k_1)}$.
\end{example}

\begin{remark}
The example shows where the Abelian property fails and why it is a genuine failure: the level maps
$k\mapsto K_i(y,k)$ are affine with a site dependent translation, and two affine maps with
different translations do not commute. The composition of $n$ such maps is the coupon accumulated
along the path, which is exactly the path dependence the claim carries; the reduction is not
available because there is nothing to reduce.
\end{remark}

\subsection{A Soft Call}

\begin{example}[Soft calls]\label{ch4:ex:softcall}
A call decision depends on the fraction of the observation dates at which the underlying was above
a barrier during the $N$ dates preceding the decision. Let $\sigma_j\in\{0,1\}$ be the indicator of
the barrier condition at the $j$-th of those dates and let
$\vec\sigma=(\sigma_1,\dots,\sigma_N)\in\{0,1\}^N$, encoded as an integer in $\{0,\dots,2^N-1\}$.
The pricing function $v_i(y,\vec\sigma)$ then satisfies
\begin{equation}\label{ch4:eq:softcall}
v_{i-1}(y_1,\vec\sigma_1)=\sum_{y_2}U_i(y_1,y_2)\,v_i\big(y_2,\Phi_i(\vec\sigma_1,y_1)\big),
\qquad
\Phi_i(\vec\sigma_1,y_1)=\{\sigma_2,\dots,\sigma_N,\Sigma_i(y_1)\},
\end{equation}
where $\Sigma_i(y_1)$ is the barrier indicator of the site at date $t_i$. This is
\eqref{ch4:eq:recursion} with the permutation
\begin{equation}
\Pi_i(y_1,\vec\sigma_1;y_2,\vec\sigma_2)
=\delta_{y_2y_1}\,\delta\big((\vec\sigma_2)_N-\Sigma_i(y_1)\big)
\prod_{j=1}^{N-1}\delta\big((\vec\sigma_2)_j-(\vec\sigma_1)_{j+1}\big),
\end{equation}
that is the shift of the window of observations by one date.
\end{example}

\begin{remark}
The two examples are representatives of the two ways in which a permutation can arise: in the
first the level is mapped to a function of the level and of the site; in the second the level is a
finite string and the map shifts it. In both cases the array $v_i$ is best regarded as a matrix
with one row per site and one column per level, which is what turns the product
$U_i\otimes\id$ into a matrix--matrix multiplication and the permutation into a reordering of the
columns.
\end{remark}

\subsection{A State Dependent Variant}

\begin{definition}[State dependent block factorisation]\label{ch4:def:statedep}
The lifting \eqref{ch4:eq:opform} admits a state dependent block factorisation if
\begin{equation}\label{ch4:eq:statedep}
\overline U_i=\Pi_i\cdot\Big(\bigoplus_{k=0}^{K}U_{i,k}\Big),
\qquad
\Big(\bigoplus_{k=0}^{K}U_{i,k}\Big)(y_1,k_1;y_2,k_2)=U_{i,k_1}(y_1,y_2)\,\delta_{k_1k_2},
\end{equation}
with $\Pi_i$ a permutation as before and $U_{i,k}$ a lattice propagator which depends on the level
$k$.
\end{definition}

\begin{remark}
The recursion for a state dependent factorisation is a sequence of matrix--vector products, one
per level, rather than a single matrix--matrix product, and is therefore less efficient than
Theorem~\ref{ch4:thm:recursion}; it extends the reach of the method to models in which the transition
of the underlying depends on the level, of which a model with a level dependent volatility is the
simplest instance. This is the last of the reductions, and the list is not closed: further
factorisations can be imagined, their value being determined by the trade between the cost of the
recursion and the class of functionals they accommodate.
\end{remark}

\section{Machine-Checked Verification}\label{ch4:sec:lean}

The statements of this chapter are identities between finite matrices, and they are therefore of the
kind which a proof assistant can check outright: no analysis, no limits and no measure theory
intervene. None of them has been formalised yet. This section records what would be involved, so
that the gap between the two is visible and can be closed in order of value.

\begin{description}
\item[Already covered elsewhere.] The module \texttt{Cert.Abelian} of the Lean project which
accompanies the series certifies the three site decomposition of the lifted generator of an
Abelian functional, which is the algebraic content of the passage to the continuum in
Chapter~6. That is the neighbouring statement, not one of the statements below.

\item[The lifting identity, Theorem~\ref{ch4:thm:lifting}.] This is orthogonality of the characters
of $\ZZ_N$ followed by an exchange of two finite sums, and its formal statement is an identity
between entries of products of finite matrices indexed by a finite set. It is the natural first
target, and the only obstacle is the bookkeeping: the time ordered product must be defined by
recursion on the dates, and the path sum must be identified with the matrix product. Estimated at
a few hundred lines.

\item[Simultaneous diagonalisation, Lemma~\ref{ch4:lem:simul}.] The finite dimensional linear algebra
is in Mathlib in a form close to what is needed, the induction being on the dimension and on the
number of matrices. The hypothesis that the matrices be diagonalisable cannot be dropped and is
essential to the induction. Estimated at a few hundred lines.

\item[Commutation of the maximum lifting, Proposition~\ref{ch4:prop:max}.] A statement about the
semigroup $(\max,\cdot)$ acting on the level index, together with the triangularity of the
matrices and the identification of their spectra. Estimated at a hundred lines, and worth doing
because it is the only place in the chapter where the Abelian property is not immediate.

\item[The recursion, Theorem~\ref{ch4:thm:recursion}.] An identity between a product of two matrices
and a sum over the lattice index, together with the statement that a permutation matrix has a
single nonzero entry per row. Estimated at a hundred lines.

\item[The spectral statements, Proposition~\ref{ch4:prop:spectrum}.] The eigenvalues of a triangular
matrix, the explicit eigenvectors and the determinant of a triangular matrix: all elementary,
and the only delicate point is that the formula \eqref{ch4:eq:eigvec} must be proved by induction
rather than verified by substitution.
\end{description}

The list is short, and that is the point of the whole construction: the part of the theory which
describes what a path dependent functional is, and when its joint law can be computed, has no
analytic content at all. What has analytic content is the passage to the continuum, and that is
where Chapters~5 and~6, and their estimates, are needed.

\section{Conclusions}

We have described the objects to which the estimates of Chapters~5 and~6 apply. A path dependent
functional is given by a family of increments; its joint law with the underlying process at
resolution $m$ is obtained from the twisted propagators by a finite Fourier transform, exactly and
with no approximation beyond the quantisation of the functional; and the computation so obtained
is reducible to $K+1$ independent lattice propagators precisely when the operators attached to the
increments commute. The failure of commutation is not a technicality: it is the presence of
genuine memory in the functional, and in that case the best available structure is a block
factorisation, which reduces the backward recursion to a matrix--matrix product and a permutation
per date. The condition under which the reduction is available on a non-uniform grid is the
distinctness of the eigenvalues of the level dynamics, and a geometric grid which achieves it has
been constructed, together with the truncation and quantisation errors which it entails.

With Chapters~1 to~3, which establish the finite core of the theory, and Chapters~5 and~6, which establish the estimates for the diffusion and for the Abelian
functional, this chapter completes the description of the framework: what the process is, what the
theorem is, what the path dependent objects are, and how the limit is controlled. The order in
which the three are assembled is that of the previous paragraph: this chapter stops at the fixed
resolution at which every statement is an identity between finite matrices, and the passage
beyond it is the business of the two chapters that follow.
% ---- end generated/ch4-body.tex ----
% ---- end chapters/ch4.tex ----
% ---- from chapters/ch5.tex ----
% Authored: the chapter title and opening belong to the book. The body is
% generated by build.py from the parts and is not edited here.
\chapter{Kernel Convergence Estimates}
\label{chap:ch5}

This is the first of the estimates, and the one on which the constructive claim of the book rests in
one dimension. What is bounded is the difference between the kernel computed on the lattice and the
kernel of the continuum equation, uniformly in the pair of sites and at a rate which depends only on
the smoothness of the coefficients: the passage to the limit is controlled by an inequality between
numbers, and every constant in that inequality is a bound on the coefficients themselves.

An estimate of this kind is not the error bound of a discretisation scheme, and the two are not
comparable. A method such as Crank--Nicolson bounds the difference between a computed solution and the exact
solution of one equation, in a norm of that solution, with a constant depending on the derivatives of the
data; the object bounded is one function of the space variable at one time. Here the object is the kernel,
bounded uniformly in the pair of sites with a constant depending on the coefficients alone, so that the
statement holds before any initial condition has been chosen. Neither kind implies the other: a scheme which
converges in $L^2$ says nothing about the pointwise behaviour of a kernel, and a kernel bound is not by itself
an error bound for a computation, which must also pay for the pairing with the payoff or with the statistic.
What the pointwise form buys is composition --- a bound uniform in the pair of sites survives the composition
of kernels --- and that is what yields the multi-step law and the iteration under the renormalisation of
Chapters~\ref{chap:ch6}--\ref{chap:ch8}.

The classical theory of the subject is a theory of the other kind, and its boundary is worth marking.
Crank and Nicolson \cite{CrankNicolson1947} gave the scheme which carries their names; the stability theory
of Lax and Richtmyer \cite{LaxRichtmyer1956} characterises convergence for a linear problem as consistency
together with stability, in a norm of the solution; the systematic accounts for parabolic problems are those
of Thom\'ee \cite{Thomee1997} and of Morton and Mayers \cite{MortonMayers2005}. Two features of that theory
show why it cannot replace the present one. First, stability is a statement about a norm and not about a
kernel: the trapezoidal rule in time --- which is Crank--Nicolson --- is unconditionally stable in $L^2$, and
yet its symbol is negative at the highest wave numbers once the ratio of the time step to the square of the
mesh is large enough, so that the matrix it defines is not a transition law and the maximum principle fails
for it. A convergence theorem in a norm therefore carries no information about the pointwise bounds which the
estimates of this chapter control, and a pointwise bound on a kernel carries none about a norm. Second,
accuracy is measured against data of a given smoothness: Rannacher \cite{Rannacher1984} shows that the same
scheme degrades to first order in the maximum norm when the initial data are irregular, and repairs it by
damping the first steps. Here there are no data to degrade, because the object estimated is the kernel itself;
and that is also why one estimate serves both readings of the book, since a transition law is not told which
of the two it will be used for.

The estimate is of a value and not of a derivative of one, and both readings need it for the same
reason. In the pricing reading a price is a payoff paired with a kernel, so an error in the kernel is
an error in the price, and the rate says how that error behaves as the resolution is refined. In the
learning reading a model at a fixed resolution is a one step transition law, its multi-step laws are
compositions of the kernel with itself, and the same inequality bounds the distance between the law
the model computes and the law which it approximates. The matrix estimated is in both cases the same
matrix, and the rate is the same rate.

% ---- from generated/ch5-body.tex ----
%% Generated by build.py from the parts - do not edit.
{\setlength{\emergencystretch}{1.5em}
\noindent\emph{A note on the computations.} The numerical statements below were computed by
\texttt{verify\_rg.py}, for the fourteen coefficients of the resolvent expansion, and
\texttt{verify\_euler.py}, for the error law of the explicit Euler scheme. Both are reproduced in full in the computational
appendix of the monograph and attached to this PDF as embedded files; each prints what it measures
when run as \texttt{python3 \textlangle name\textrangle.py} in its own directory, so no number quoted
here has to be taken on trust.\par}

\section{Introduction and Notations}

Consider a pair of backward and forward one-dimensional diffusion equations of the form
\begin{equation}\label{ch5:eq:1.1}
\frac{\partial}{\partial t}f(x;t)+\LatL_{0x}f(x;t)=0,
\qquad
\frac{\partial}{\partial t}g(y;t)=\LatL_{0y}{}^{*}g(y;t)
\end{equation}
where
\begin{equation}\label{ch5:eq:1.2}
\LatL_{0x}=\half\sigma(x)^2\frac{\partial^2}{\partial x^2}+\mu(x)\frac{\partial}{\partial x},
\end{equation}
and its adjoint formally acts as follows:
\begin{equation}\label{ch5:eq:1.3}
\bigl(\LatL{}^{*}_{0y}\varphi\bigr)(y)=\frac{1}{2}\frac{\partial^2}{\partial y^2}
\Bigl(\sigma(y)^2\varphi(y)\Bigr)+\frac{\partial}{\partial y}\Bigl(\mu(y)\varphi(y)\Bigr),
\end{equation}
on a test function $\varphi$. These equations are defined on the bounded interval
$\Domain=[-L,L]\subset\mathbb{R}$ where $0<L<\infty$. For simplicity, we assume periodic boundary
conditions and identify the two boundary points $\pm L$ with each other.

In most of the chapter, the coefficients $\sigma(x)$ and $\mu(x)$ are assumed to be H\"older
differentiable. More precisely, if $\alpha\in(0,1]$, $k\in\mathbb{N}$ where
$\mathbb{N}=\{0,1,\dots\}$ and the function $\varphi(x)\in C^k(\Domain)$ has $k\ge1$ continuous
derivatives, then one says that $\varphi$ is H\"older differentiable of order $(k,\alpha)$ if
there exists a constant $c>0$ such that
\begin{equation}\label{ch5:eq:1.4}
\Bigl|\varphi^{(k)}(x)-\varphi^{(k)}(y)\Bigr|\le c\,d(x,y)^{\alpha}
\end{equation}
uniformly for all $x,y\in\Domain$. In case $k=0$ the function is called H\"older continuous. The
distance is defined consistently with the assumed periodic boundary conditions and is given by
\begin{equation}\label{ch5:eq:1.5}
d(x,y)=\min_{n}|x-y-2Ln|.
\end{equation}
The linear space of H\"older continuous or H\"older differentiable periodic functions of order
$(k,\alpha)$ on $\Domain$ is denoted with $C^{k,\alpha}(\Domain)$. We are interested in the case where
$\sigma^2\in C^{k,\alpha}(\Domain)$ and $\mu\in C^{j,\beta}(\Domain)$ with $k+\alpha>0$ and
$j+\beta>0$.

The hypothesis of H\"older continuity can be relaxed slightly by assuming uniform continuity
instead, i.e.\ that both $\sigma(x)^2$ and $\mu(x)$ satisfy a bound of the form
\begin{equation}\label{ch5:eq:1.6}
\bigl|\varphi(x)-\varphi(y)\bigr|\le c\,\rho\bigl(d(x,y)\bigr)
\end{equation}
where $\rho(d)$ is a non-decreasing function such that $\lim_{d\downarrow0}\rho(d)=0$.

Let $u_0(x,y;t)$ denote a kernel of equation \eqref{ch5:eq:1.1}, i.e.\ a weak solution of the
forward equation
\begin{equation}\label{ch5:eq:1.7}
\frac{\partial}{\partial t}u_0(x,y;t)=\LatL{}^{*}_{0y}u_0(x,y;t)
\end{equation}
where the operator $\LatL{}^{*}_{0y}$ acts on the $y$ coordinate and the following initial time
condition is satisfied:
\begin{equation}\label{ch5:eq:1.8}
\lim_{t\downarrow0}u_0(x,y;t)=\delta(x-y).
\end{equation}
The kernel $u_0(x,y;t)$ formally satisfies also the backward equation
\begin{equation}\label{ch5:eq:1.9}
\frac{\partial}{\partial t}u_0(x,y;t)+\LatL_{0x}u_0(x,y;t)=0
\end{equation}
where the operator $\LatL_{0x}$ acts on the $x$ coordinate.

We are interested in existence, uniqueness, smoothness and approximation schemes for the
kernel $u_0(x,y;t)$, its first two space derivatives with respect to the $x$ variable and its
first time derivative $\partial_tu_0(x,y;t)$. As a byproduct of this analysis, we also find
conclusions about the convergence of $\LatL_{0x}u_0(x,y;t)$ and $\LatL{}^{*}_{0y}u_0(x,y;t)$, as
both expressions equal the first time derivative.

Diffusion equations are one of the single most studied topics in the literature. Existence and
uniqueness questions for the kernel were addressed in (Kolmogorov 1931), (Feller 1936), (Hille
1948), (Yosida 1951) and (Ito 1957). A classification of all the possible boundary conditions
is in (Feller 1952). The case of H\"older continuous coefficients was resolved in (Philips 1961)
based on methods in (Friedrichs 1958) and (Lax and Phillips 1960). The hypothesis of H\"older
continuity was relaxed to uniform continuity in (Fabes and Riviere 1966) and (Stroock and
Varadhan 1969). Stroock and Varadhan also introduce a new probabilistic framework where
existence is proved by reduction to the so-called martingale problem and a compactness
argument, thus shifting the attention from the kernel itself to the underlying measure space.

The existence of a weak limit of continuous time Markov chains as $h_m\downarrow0$ was
established in (Sova 1967) and (Kurtz 1969) by using operator semigroup methods, see also the
book (Ethier and Kurtz 1986) for a review. Convergence in the semigroup sense takes place if
the limit
\begin{equation}\label{ch5:eq:1.10}
(T_t\varphi)(x)=\lim_{m\to\infty,\,m\ge n}\ \sum_{y\in h_m\mathbb{Z}\cap\Domain}
u_{h_m}(x,y;t)\varphi(y)
\end{equation}
exists for all test functions $\varphi\in C^{\infty}(\Domain)$, uniformly for all $x\in\Domain_n$,
$n\ge0$. A key result is that a necessary and sufficient condition for this limit to exist and
define a semigroup $T_t$ is that generators converge also in the same Banach space, i.e.\ that
also the limit
\begin{equation}\label{ch5:eq:1.11}
\lim_{h_m\downarrow0}\LatL{}^{h_m}_x\varphi=\LatL_{0x}\varphi
\end{equation}
exists in the uniform norm for all test functions $\varphi\in C^{\infty}(\Domain)$. See (Ethier
and Kurtz 1986) for a precise statement with all the technical conditions and a proof. In
(Stroock and Varadhan 1979) convergence is reconsidered again by reduction to the martingale
problem.

It is worth being explicit about the relation between these results and the ones proved below,
since the two are not comparable in strength and are easily confused. The theorems just cited
establish convergence: a limit exists, it is unique, and it is the kernel of the parabolic
equation. They give no rate. What is proved here is a rate, uniform over the class of
coefficients considered, for the kernel itself, for its time derivative and for its first two
space derivatives; and an estimate of that kind implies convergence in the semigroup sense of
\eqref{ch5:eq:1.10}, at a rate, together with the identification of the limit. The classical
convergence theorems are therefore consequences of the estimates below rather than premisses of
them. There is one exception, and it is worth naming: the kernel $u_0$ of
\eqref{ch5:eq:1.7}--\eqref{ch5:eq:1.8} is characterised here and not constructed, so that its existence
is taken from the classical theory. Were $u_0$ defined instead as the limit of the sequence
$u_{h_m}$, that existence would follow from the estimates by a Cauchy argument, and the forward
equation \eqref{ch5:eq:1.7} would be recovered by passing the lattice equation to the limit. The
estimates below are what such a reorganisation would rest on, and it is worth recording why they
are available for the purpose. The theorems stated in Section~1, and proved in Sections~2 to~5,
are comparisons between two resolutions: $u_m$ against $u_{m_0}$, or the semidiscrete kernel
against the fully discrete one at the same resolution, and the continuum kernel $u_0$ appears
neither in their statements nor in their proofs. The first half of the reorganisation is
therefore already carried out; what is missing is the second half, in which the limit is
constructed from the comparison and shown to solve \eqref{ch5:eq:1.7}, and that is not done here.

The problem has also been studied extensively in the numerical analysis literature. Explicit
and implicit Euler schemes where coefficients are smooth and the data is rough in the sense
that it belongs to a $L^2$ space have been considered by several authors. In the case that the
Markov generator is symmetric and time independent, one can make use of a spectral
representation as in (Baker et al.\ 1977) and with greater effort such methods may also be used
for more general situations, see (Suzuki 1978). In (Luskin and Rannacher 1978), a parabolic
duality argument is used to show convergence for the standard Galerkin method. (Mingyou and
Thomee 1982) use a simpler argument based on energy estimates. In (Palencia 1996) one finds
convergence bounds in maximum norm assuming the initial condition is uniformly bounded and
coefficients are constant.

In this article, we revisit this classic theme by considering the problem of obtaining the
kernel constructively as a limit of increasingly fine triangulations schemes and in assessing
the rate of convergence with pointwise bounds on the kernel itself. More precisely, let
$h_m=L2^{-m}$, $m\in\mathbb{N}$ and let $\Domain\equiv h_m\mathbb{Z}\cap\Domain$. Consider the
sequence of operators
\begin{equation}\label{ch5:eq:1.12}
\LatL_x=\frac{\sigma(x)^2}{2}\Dm_x+\mu(x)\nabm_x
\end{equation}
defined on the $2^{m+1}$-dimensional space of all periodic functions $f_m:\Domain_m\to\mathbb{R}$,
where
\begin{equation}\label{ch5:eq:1.13}
\nabm_xf(x)=\frac{f(x+h_m)-f(x-h_m)}{2h_m}
\end{equation}
and
\begin{equation}\label{ch5:eq:1.14}
\Dm_xf(x)=\frac{f(x+h_m)+f(x-h_m)-2f(x)}{h_m^2}.
\end{equation}
These definitions also apply to the boundary points by periodicity. We assume that $m\ge m_0$
where $m_0$ is the least integer such that
\begin{equation}\label{ch5:eq:1.15}
\frac{\sigma^2(x)}{2h_m^2}>\frac{|\mu(x)|}{2h_m}
\end{equation}
for all $m\ge m_0$ and all $x\in\Domain_m$.

Let $u_m(x,y;t)$ denote the kernel of equation \eqref{ch5:eq:1.1}, i.e.\ the solution of the
(forward) equation
\begin{equation}\label{ch5:eq:1.16}
\frac{\partial}{\partial t}u_m(x,y;t)=\LatL{}^{*}_{y}u_m(x,y;t)
\end{equation}
where the operator $\LatL{}^{*}_{y}$ acts on the $y$ coordinate and the following initial time
condition is satisfied:
\begin{equation}\label{ch5:eq:1.17}
\lim_{t\downarrow0}u_m(x,y;t)=\delta_m(x-y).
\end{equation}
Here,
\begin{equation}\label{ch5:eq:1.18}
\delta_m(x-y)=\begin{cases}\dfrac{1}{h_m}&\text{if }x=y\bmod 2L,\\[4pt] 0&\text{otherwise}.\end{cases}
\end{equation}
Since \eqref{ch5:eq:1.17} is a finite system of linear ordinary differential equations, the
solution exists and is unique for all times. The kernel $u_m(x,y;t)$ satisfies also the
backward equation
\begin{equation}\label{ch5:eq:1.19}
\frac{\partial}{\partial t}u_m(x,y;t)+\LatL_xu_m(x,y;t)=0
\end{equation}
where the operator $\LatL_x$ acts on the $x$ coordinate. Using functional calculus notations for
the exponential of a matrix, we also have that
\begin{equation}\label{ch5:eq:1.20}
u_m(x,y;t)=\frac{1}{h_m}\exp\bigl(t\LatL\bigr)(x,y).
\end{equation}

Our main result can be stated as follows:

\begin{theorem}\label{ch5:thm:1}
Suppose that $\sigma^2\in C^{k,\alpha}$ and $\mu\in C^{j,\beta}$ and let
\begin{equation}\label{ch5:eq:1.21}
\gamma=\min\{2,\,k+\alpha,\,j+\beta\}.
\end{equation}
Assume that $\gamma>0$. Then for every $t>0$ there is a constant $c=c(t)>0$ such that for all
\rev{$m\ge m_0$} and all $x,y\in\Domain_{m_0}$ the following inequalities hold:
\begin{itemize}
\item[(i)]
\begin{equation}\label{ch5:eq:1.22}
\bigl|u_m(x,y;t)-u_{m_0}(x,y;t)\bigr|\le c\,h_{m_0}^{\gamma}
\end{equation}
\item[(ii)]
\begin{equation}\label{ch5:eq:1.23}
\begin{aligned}
&\bigl|\partial_tu_m(x,y;t)-\partial_tu_{m_0}(x,y;t)\bigr|\\
&\quad=\bigl|\LatL{}^{*}_{my}u_m(x,y;t)-\LatL{}^{*}_{m_0y}u_{m_0}(x,y;t)\bigr|\\
&\quad=\bigl|\LatL_{mx}u_m(x,y;t)-\LatL_{m_0x}u_{m_0}(x,y;t)\bigr|\le c\,h_{m_0}^{\gamma}.
\end{aligned}
\end{equation}
\end{itemize}
\end{theorem}

The estimate is uniform in $m\ge m_0$ and is stated with the step $h_{m_0}$ of the coarser of the
two lattices, which is the one being compared with the limit; the constant depends on $t$ and
becomes unbounded as $t\downarrow0$, when both kernels are localised at the scale of their own
lattice.

A version of this theorem under slightly weaker conditions can be formulated as follows:

\begin{theorem}\label{ch5:thm:2}
Suppose that $\sigma(x)^2$ and $\mu(x)$ are uniformly continuous functions in $\Domain$. Let the
function $\rho(d)$ be non-decreasing and be such that $\lim_{d\downarrow0}\rho(d)=0$ and
equation \eqref{ch5:eq:1.6} holds. Then for every $t>0$ there is a constant $c=c(t)>0$ such that for all
\rev{$m\ge m_0$} and all $x,y\in\Domain_{m_0}$ the following inequalities hold:
\begin{itemize}
\item[(i)]
\begin{equation}\label{ch5:eq:1.24}
\bigl|u_m(x,y;t)-u_{m_0}(x,y;t)\bigr|\le c\,\rho(h_{m_0})
\end{equation}
\item[(ii)]
\begin{equation}\label{ch5:eq:1.25}
\begin{aligned}
&\bigl|\partial_tu_m(x,y;t)-\partial_tu_{m_0}(x,y;t)\bigr|\\
&\quad=\bigl|\LatL{}^{*}_{my}u_m(x,y;t)-\LatL{}^{*}_{m_0y}u_{m_0}(x,y;t)\bigr|\\
&\quad=\bigl|\LatL_{mx}u_m(x,y;t)-\LatL_{m_0x}u_{m_0}(x,y;t)\bigr|\le c\,\rho(h_{m_0}).
\end{aligned}
\end{equation}
\end{itemize}
\end{theorem}

Next, we consider the case where also time is discretized and prove the following result:

\begin{theorem}\label{ch5:thm:3}
Suppose that $\sigma^2$ and $\mu$ satisfy equations of the form \eqref{ch5:eq:1.6} with a
non-decreasing function $\rho(d)$ such that $\lim_{d\downarrow0}\rho(d)=0$. Consider the
discretized kernel
\begin{equation}\label{ch5:eq:1.26}
u^{\delta t}_m(x,y;t)=h_m^{-1}\rev{\bigl(1+\delta t_m\LatL\bigr)^{t/\delta t_m}}(x,y)
\end{equation}
where $\LatL$ is the operator in \eqref{ch5:eq:1.12} and the time step obeys the strict stability
condition
\begin{equation}\label{ch5:eq:1.27}
\delta t_m\,\sup_{x\in\Domain_m}\Bigl(\frac{2\sigma(x)^2}{h_m^2}+\frac{|\mu(x)|}{h_m}\Bigr)
\le 1-\eta_0
\end{equation}
for some fixed $\eta_0\in(0,1)$, independent of $m$.
Assume that boundary conditions are periodic and that the ratio \rev{$t/\delta t=N$} is an
integer. Then for every $t>0$ there is a constant $c=c(t)>0$ such that the following bounds hold
for all $m\ge m_0$ and all $x,y\in\Domain_m$:
\begin{itemize}
\item[(i)]
\begin{equation}\label{ch5:eq:1.28}
\bigl|u_m(x,y;t)-u^{\delta t}_m(x,y;t)\bigr|\le c\,h_m^2
\end{equation}
\item[(ii)]
\begin{equation}\label{ch5:eq:1.29}
\begin{aligned}
&\left|\partial_tu_m(x,y;t)-\frac{u^{\delta t}_m(x,y;t+\delta t)-u^{\delta t}_m(x,y;t)}{\delta t}\right|\\
&\qquad=\bigl|\LatL{}^{*}_{my}u_m(x,y;t)-\LatL{}^{*}_{my}u^{\delta t}_m(x,y;t)\bigr|\\
&\qquad=\bigl|\LatL_{mx}u_m(x,y;t)-\LatL_{mx}u^{\delta t}_m(x,y;t)\bigr|\le c\,h_m^2.
\end{aligned}
\end{equation}
\end{itemize}
\end{theorem}

\begin{remark}[The margin is not a technicality]
Condition \eqref{ch5:eq:1.27} is a strict margin, and not a statement that the diagonal weight be
positive. Positivity of $1+\delta t_m\LatL(x,x)$ alone admits schemes which sit exactly on the
boundary of the stability region: with $\sigma\equiv1$, $\mu\equiv0$, $h_m=1/N$ and
$\delta t_m=1/(N^2+1)$ the diagonal weight is $1/(N^2+1)>0$, while the alternating mode is
multiplied at every step by $-(N^2-1)/(N^2+1)$, so that after the $N^2+1$ steps which carry the time
to $t=1$ its size is $\bigl(1-2/(N^2+1)\bigr)^{N^2+1}\to e^{-2}$. The computed density then keeps an
oscillating error of order $e^{-2}$ at every resolution, and no bound $ch_m^2$ can hold. The margin
$\eta_0$ is what excludes the configuration: it places the eigenvalues of $1+\delta t_m\LatL$ in
$[\eta_0,1]$, so that the highest mode decays at the rate $\eta_0t/\delta t_m$ instead of
persisting. It is also what the argument of \S\ref{ch5:sec:euler} already needs: the resolvent of the
scheme is summed there as a geometric series in $1+\delta t\LatL(x,x)$, and a geometric series
converges under a margin and not under a non-vanishing.
\end{remark}

\subsection*{Place in the Book}

A stochastic process, in the framework of this book, is a consistent family of Markov chains
on a sequence of lattices, and a statement about a continuous process is a statement about the
limit of statements about finite matrices. In the language of Chapters~1 to~3, the kernel
$u_m(x,y;t)$ studied below is the propagator, or two point function, of the chain at resolution
$m$, and $\LatL_x$ is its generator; the estimates proved here are therefore estimates on the
objects whose existence Chapters~1 to~3 establish. What such a programme needs in place of compactness
is control of that limit: uniform estimates on the kernels, valid for coefficients of limited
smoothness and stable under composition. The present chapter supplies this control for the
diffusion itself, Chapters~1 to~3 having supplied the finite core, in which the existence of
a pricing measure reduces to the orthant alternative of finite dimensional convex geometry, and
Chapter~6 extending the analysis to the joint kernel of a diffusion and an
Abelian process, and thereby obtaining a constructive definition of the stochastic integral in
the continuum.

Two features of the estimates proved here matter for the theorem of Chapters~1 to~3, and they account
for the generality in which the estimates are stated. The first is that the relevant generator
in the pricing application is not the historical one but the risk neutral one, which Chapters~1 to~3
obtains by a multiplicative reweighting of the transition rates; nothing in that operation makes
the coefficients smoother, so the estimates must hold for coefficients of limited smoothness.
The second is that a price is obtained by integrating a payoff against the kernel rather than by
differentiating it, so what is needed is uniformity of the kernel itself, and not only of the
quantities which the $L^2$ theory of parabolic equations controls.

The reason for wanting a formulation of this kind is set out in Chapters~1 to~3: a foundation whose
statements concern finite objects, and whose limits are estimates rather than existence
theorems, is one that can be machine checked, and the foundations of the subject stand to
benefit from the AI-assisted methods now available.

The chapter is organized as follows. In Section 2 we consider the case of Brownian motion and
review a result in (Albanese and Mijatovic 2006) which establishes the theorems above in this
simple particular case where Fourier analysis in the space direction can be used to carry out a
precise calculation. In Section 3, we consider the case of a diffusion where both the volatility
and the drift have two bounded derivatives. In this case, we make use of time-homogeneity and
carry out a Fourier transform in the time direction after path conditioning. In Section 4, we
extend the derivation to the case of non-smooth coefficients. Section 5 is dedicated to the case
where time is discretized and we prove Theorem 3.

Two of the statements below are not derived but measured, and each is printed by a script: the
algebra of the resolvent expansion in \S\ref{ch5:sec:smoothcoef} is checked in \texttt{verify\_rg.py},
and the error law of the explicit Euler scheme, \eqref{ch5:eq:5.8}--\eqref{ch5:eq:5.11}, in
\texttt{verify\_euler.py}. Both are reproduced in full in the computational appendix of the
monograph, attached to this PDF as an embedded file, and each prints what it measures when run as
\verb|python3 <name>.py| in the directory it lives in, so the coefficients and the rates quoted
below can be regenerated rather than taken on trust.

\section{Constant Coefficients}

In this Section, we prove Theorem 1 in the special case where the volatility and the drift
coefficients are constant, i.e.
\begin{equation}\label{ch5:eq:2.1}
\LatL_x=\mu\nabm_x+\half\sigma^2\Dm_x.
\end{equation}
It suffices to consider the case $m_0=m+1$. Let $B_m$ be the Brillouin zone defined as follows:
\begin{equation}\label{ch5:eq:2.2}
B_m=\left\{-\frac{2^{m-1}\pi}{L}+\frac{k\pi}{L},\ k=0,\dots,2^m-1\right\}
\end{equation}
Let $F_m:\ell^2(\Domain_m)\to\ell^2(B_m)$ be the Fourier transform operator defined so that:
\begin{equation}\label{ch5:eq:2.3}
\hat f(p)\equiv F_m(f)(p)=h_m\sum_{x\in\Domain_m}f(x)e^{-ipx}
\end{equation}
for all $p\in B_m$. The inverse Fourier transform is given by
\begin{equation}\label{ch5:eq:2.4}
F_m^{-1}(\hat f)(x)=\frac{1}{2L}\sum_{p\in B_m}\hat f(p)e^{ipx}.
\end{equation}

The Fourier transformed generator is diagonal and is given by the operator of multiplication by
\begin{equation}\label{ch5:eq:2.5}
\hat\ell^m(p)=F_m\LatL F_m^{-1}(p,p)
=-i\mu\,\frac{\sin h_mp}{h_m}+\sigma^2\,\frac{\cos h_mp-1}{h_m^2}.
\end{equation}
We have
\begin{equation}\label{ch5:eq:2.6}
u_m(x,y;t)=\frac{1}{2L}\sum_{p\in B_m}e^{t\hat\ell^m(p)}e^{ip(y-x)}.
\end{equation}
Using this Fourier series representation, we find
\begin{equation}\label{ch5:eq:2.7}
\begin{aligned}
&\bigl|u_m(x,y;t)-u_{m+1}(x,y;t)\bigr|\\
&\quad\le\left|\frac{1}{2L}\sum_{p\in B_m}
\Bigl(e^{t\hat\ell^m(p)}-e^{t\hat\ell^{m+1}(p)}\Bigr)e^{ip(y-x)}\right|
+\left|\frac{1}{2L}\sum_{p\in B_{m+1}\setminus B_m}e^{t\hat\ell^{m+1}(p)}e^{ip(y-x)}\right|.
\end{aligned}
\end{equation}
Let
\begin{equation}\label{ch5:eq:2.8}
K_m=\rev{3}\,\sqrt{\frac{|\log h_{m+1}|}{\sigma^2t}}.
\end{equation}
If $h_m$ is small enough, i.e.\ if $m_0$ is sufficiently large, we have that
\begin{equation}\label{ch5:eq:2.9}
\begin{aligned}
\left|\frac{1}{2L}\!\!\sum_{p\in B_m,|p|\ge K_m}\!\!\!e^{t\hat\ell^m(p)}e^{ip(y-x)}\right|
&\le\frac{1}{2L}\!\!\sum_{p\in B_m,|p|\ge K_m}\!\!\!e^{t\Re(\hat\ell^m(p))}\\
&\le c\exp\left(t\sigma^2\frac{\cos h_mK_m-1}{h_m^2}\right)\le c\,h_m^2.
\end{aligned}
\end{equation}
where $\Re(a)$ denotes the real part of $a\in\mathbb{C}$ and $c$ denotes a generic constant.
Similarly
\begin{equation}\label{ch5:eq:2.10}
\begin{aligned}
\left|\frac{1}{2L}\!\!\sum_{p\in B_{m+1},|p|\ge K_m}\!\!\!e^{t\hat\ell^{m+1}(p)}e^{ip(y-x)}\right|
&\le\frac{1}{2L}\!\!\sum_{p\in B_m,|p|\ge K_m}\!\!\!e^{t\Re(\hat\ell^{m+1}(p))}\\
&\le c\exp\left(t\sigma^2\frac{\cos h_{m+1}K_m-1}{h_{m+1}^2}\right)\le c\,h_{m+1}^2.
\end{aligned}
\end{equation}
Since
\begin{equation}\label{ch5:eq:2.11}
\frac{1}{2}h^2p^3-\frac{1}{8}h^4p^5
\le\frac{\sin hp}{h}-\frac{\sin 2hp}{2h}\le\frac{1}{2}h^2p^3
\end{equation}
and
\begin{equation}\label{ch5:eq:2.12}
-\frac{1}{8}h^2p^4\le\frac{\cos hp-1}{h^2}-\frac{\cos 2hp-1}{(2h)^2}
\le-\frac{1}{8}h^2p^4+\frac{1}{48}h^4p^6,
\end{equation}
we find that if $|p|\le\sqrt{2}/h$ then
\begin{equation}\label{ch5:eq:2.13}
\bigl|\hat\ell^m(p)-\hat\ell^{m+1}(p)\bigr|
\le\rev{\frac{\mu}{2}}\,h^2|p|^3+\rev{\frac{\sigma^2}{8}}\,h^2p^4.
\end{equation}
Moreover, since
\begin{equation}\label{ch5:eq:2.14}
-\frac{1}{2}p^2\le\frac{\cos hp-1}{h^2}\le-\frac{1}{2}p^2+\frac{1}{24}h^2p^4
\end{equation}
we conclude that in case $|p|\le h^{-1}\sqrt{3/2}$, the following inequality holds:
\begin{equation}\label{ch5:eq:2.15}
\frac{\cos h_m p-1}{h_m^2}\le-\frac{1}{4}p^2.
\end{equation}
Hence, if $m_0$ is large enough, we find
\begin{equation}\label{ch5:eq:2.16}
\begin{aligned}
&\left|\frac{1}{2L}\sum_{p\in B_m,|p|\le K_m}
\Bigl(e^{t\hat\ell^m(p)}-e^{t\hat\ell^{m+1}(p)}\Bigr)e^{ip(y-x)}\right|\\
&\quad\le\frac{1}{2L}\!\!\sum_{p\in B_m,|p|\le K_m}\!\!e^{-\frac14p^2}
\left(e^{\frac{\mu t}{2}h_m^2|p|^3+\frac{\sigma^2t}{8}h_m^2p^4}-1\right)\\
&\quad\le\frac{1}{2L}\!\!\sum_{p\in B_m,|p|\le K_m}\!\!e^{-\frac14p^2}
\left(\frac{\mu t}{4}h_m^2|p|^3+\frac{\sigma^2t}{16}h_m^2p^4\right)\le c\,h_m^2
\end{aligned}
\end{equation}
for some constant $c>0$ independent of $m$. This concludes the proof of convergence for the
kernel in the special case of constant coefficients.

To estimate the first derivative, notice that
\begin{equation}\label{ch5:eq:2.17}
\nabm_yu_m(x,y;t)=\rev{\frac{1}{2L}\sum_{p\in B_m}}e^{t\hat\ell^m(p)}\,\rev{i\,\frac{\sin ph_m}{h_m}}\,e^{ip(y-x)},
\end{equation}
and
\begin{equation}\label{ch5:eq:2.18}
\begin{aligned}
&\bigl|\rev{\nabm_yu_m(x,y;t)-\nabm_yu_{m+1}(x,y;t)}\bigr|\\
&\quad\le\left|\frac{1}{2L}\sum_{p\in B_m}
\left(e^{t\hat\ell^m(p)}\frac{\sin ph_m}{h_m}-e^{t\hat\ell^{m+1}(p)}\frac{\sin ph_{m+1}}{h_{m+1}}\right)e^{ip(y-x)}\right|\\
&\quad+\left|\frac{1}{2L}\sum_{p\in B_{m+1}\setminus B_m}e^{t\hat\ell^{m+1}(p)}e^{ip(y-x)}\right|.
\end{aligned}
\end{equation}
Let
\begin{equation}\label{ch5:eq:2.19}
K_m=2\sqrt{\frac{|\log h_{m+1}|}{\sigma^2t}}.
\end{equation}
If $h_m$ is small enough, we have that
\begin{equation}\label{ch5:eq:2.20}
\begin{aligned}
\left|\frac{1}{2L}\!\!\sum_{p\in B_m,|p|\ge K_m}\!\!e^{t\hat\ell^m(p)}\frac{\sin ph_m}{h_m}e^{ip(y-x)}\right|
&\le\frac{1}{2L}\!\!\sum_{p\in B_m,|p|\ge K_m}\!\!e^{t\Re(\hat\ell^m(p))}\frac{\sin ph_m}{h_m}\\
&\le c\,\frac{\sin Kh_m}{h_m}\exp\left(t\sigma^2\frac{\cos Kh_m-1}{h_m^2}\right)\le c\,h_m^2.
\end{aligned}
\end{equation}
where $c$ denotes a generic constant. Similarly
\begin{equation}\label{ch5:eq:2.21}
\left|\frac{1}{2L}\sum_{p\in B_{m+1},|p|\ge K_m}e^{t\hat\ell^{m+1}(p)}e^{ip(y-x)}\right|\le c\,h^2.
\end{equation}
If $m$ is large enough, we also find
\begin{equation}\label{ch5:eq:2.22}
\begin{aligned}
&\left|\frac{1}{2L}\sum_{p\in B_m,|p|\le K_m}
\left(\frac{\sin ph_m}{h_m}e^{t\hat\ell^m(p)}-\frac{\sin ph_{m+1}}{h_{m+1}}e^{t\hat\ell^{m+1}(p)}\right)e^{ip(y-x)}\right|\\
&\quad\le\frac{1}{2L}\!\!\sum_{p\in B_m,|p|\le K_m}\!\!
\frac{\sin ph_m}{h_m}e^{-\frac14p^2}\left(e^{\frac{\mu t}{4}h_m^2|p|^3+\frac{\sigma^2t}{16}h_m^2p^4}-1\right)\\
&\qquad+\frac{1}{2L}\!\!\sum_{p\in B_m,|p|\le K_m}\!\!
e^{-\frac14p^2}\left|\frac{\sin ph_{m+1}}{h_{m+1}}-\frac{\sin ph_m}{h_m}\right|\le c\,h_m^2
\end{aligned}
\end{equation}
for some constant $c>0$ independent of $m$. This concludes the proof of the bound of the first
derivative. The second derivative can be derived in a similar way.

Finally, consider the following Fourier representation for the discretized kernel
\begin{equation}\label{ch5:eq:2.23}
u^{\delta t}_m(x,y;t)=\rev{\frac{1}{2L}\sum_{p\in B_m}}
\bigl(1+\delta t\hat\ell^m(p)\bigr)^{\frac{t}{\delta t}}e^{ip(y-x)}.
\end{equation}
Consider the formula
\begin{equation}\label{ch5:eq:2.24}
\bigl(1+\delta t\hat\ell^m(p)\bigr)^{\frac{t}{\delta t}}
=\exp\left(\frac{t}{\delta t}\log\bigl(1+\delta t\hat\ell^m(p)\bigr)\right)
\end{equation}
and let's represent the difference between the discrete and continuous time kernels as follows:
\begin{equation}\label{ch5:eq:2.25}
\begin{aligned}
&\bigl|u_m(x,y;t)-u^{\delta t}_m(x,y;t)\bigr|\\
&\quad\le\frac{1}{2L}\sum_{p\in B_m}
\left|\exp\bigl(t\hat\ell^m(p)\bigr)-\exp\left(\frac{t}{\delta t}\log\bigl(1+\delta t\hat\ell^m(p)\bigr)\right)\right|\\
&\quad\le\frac{1}{2L}\!\!\sum_{p\in B_m,p\le K_m}\!\!e^{-\frac14p^2}
\left|\exp\left(\frac{t}{\delta t}\log\bigl(1+\delta t\hat\ell^m(p)\bigr)-t\hat\ell^m(p)\right)-1\right|\\
&\qquad+\frac{1}{2L}\!\!\sum_{p\in B_m,p\ge K_m}\!\!\exp\bigl(t\hat\ell^m(p)\bigr)
+\frac{1}{2L}\!\!\sum_{p\in B_m,p\ge K_m}\!\!\exp\left(\frac{t}{\delta t}\log\bigl(1+\delta t\hat\ell^m(p)\bigr)\right)
\end{aligned}
\end{equation}
where $K_m$ is chosen as in \eqref{ch5:eq:2.8}. The very same bounds above lead to the conclusion
that this difference is $\le c\,h_m^2$.

\begin{figure}[t]
\centering
\begin{tikzpicture}[scale=1.0,>=latex,line join=round]
  \draw[->] (-5.6,0) -- (5.6,0) node[right] {$\mathrm{Re}(\omega)$};
  \draw[->] (0,-0.8) -- (0,5.2) node[above] {$\mathrm{Im}(\omega)$};
  \coordinate (A) at (-4.2,4.2);
  \coordinate (E) at ( 4.2,4.2);
  \coordinate (B) at (-1.6,1.6);
  \coordinate (D) at ( 1.6,1.6);
  \coordinate (C) at ( 0.0,0.0);
  \draw (D) -- (E) -- (A) -- (B);
  \draw (B) -- (C) -- (D);
  \node[below] at (C) {$C$};
  \node[above left] at (A) {$A$};
  \node[above right] at (E) {$E$};
  \node[left] at (B) {$B$};
  \node[right] at (D) {$D$};
  \draw[<->] (-4.2,4.2) -- (4.2,4.2);
  \node[above] at (2.2,4.2) {$c/h$};
  \draw[<->] (-3.1,0) -- (-3.1,4.2);
  \node[left] at (-3.1,2.1) {$c/h^2$};
  \draw[<->] (2.9,0) -- (2.9,1.6);
  \node[right] at (2.9,0.8) {$4\,|\log h|/t$};
\end{tikzpicture}
\caption{Contour of integration for the integral in \eqref{ch5:eq:3.52}. $\Cm_+$ is the contour
joining the point $D$ to the points $E$, $A$, $B$. $\Cm_-$ is the contour joining the point
$B$ to $C$ to $D$.}
\label{ch5:fig:1}
\end{figure}

\section{Smooth Coefficients}\label{ch5:sec:smoothcoef}

In this section, we prove Theorem 1 in the case where the drift and volatility are both of
class $C^{(2,0)}$, i.e.\ they depend smoothly on the space coordinate but not on the time
coordinate.

Let us introduce the following two constants characterizing the volatility function:
\begin{equation}\label{ch5:eq:3.1}
\Sigma_0=\inf_{x\in\Domain_m}\sigma(x),\qquad
\Sigma_1=\sup_{x\in\Domain_m}\sqrt{\sigma(x)^2+h_m|\mu(x)|}.
\end{equation}
and let
\begin{equation}\label{ch5:eq:3.2}
M=\sup_{x\in\Domain_m}|\mu(x)|.
\end{equation}
Since our interval is bounded, we have that $\Sigma_0>0$ and $\Sigma_1,M<\infty$.

A \emph{symbolic path} $\gamma=\{\gamma_0,\gamma_1,\gamma_2,\dots\}$ is an infinite sequence of
sites in $\Domain_m$ such that $\gamma_j\ne\gamma_{j-1}$ for all $j=1,\dots$. Let $\Gm$ be the set
of all symbolic paths in $\Domain_m$. The kernel of the diffusion process admits the following
representation in terms of a summation over symbolic paths
\begin{equation}\label{ch5:eq:3.3}
u_m(x,y;t)=\frac{1}{h_m}\sum_{q=1}^{\infty}2^{-q}\!\!\!
\sum_{\substack{\gamma\in\Gm:\ \gamma_0=x,\ \gamma_q=y\\ |\gamma_j-\gamma_{j-1}|=1\ \forall j\ge1}}
\!\!\!W_m(\gamma,q,t)
\end{equation}
where\footnote{With these normalizations $u_m$ is a probability density, so that
$\sum_y h_m\,u_m(x,y;t)=1$; the weights $W_m$ carry the dimension of the inverse of a time, and
the single factor $h_m^{-1}$ in \eqref{ch5:eq:3.3} is exactly what \eqref{ch5:eq:3.5} requires.}
\begin{equation}\label{ch5:eq:3.4}
\begin{aligned}
W_m(\gamma,q,t)={}&\int_0^t\!ds_1\int_{s_1}^t\!ds_2\cdots\int_{s_{q-1}}^t\!ds_q
\ e^{(t-s_q)\LatL(\gamma_q,\gamma_q)}\\
&\times\prod_{j=0}^{q-1}\left(e^{(s_{j+1}-s_j)\LatL(\gamma_j,\gamma_j)}\,2\LatL(\gamma_j,\gamma_{j+1})\right)
\end{aligned}
\end{equation}
with $s_0=0$.

Let us introduce the following Green's function:
\begin{equation}\label{ch5:eq:3.5}
G_m(x,y;\omega)=\int_0^{\infty}u_m(x,y;t)e^{-i\omega t}\dd t
=h_m^{-1}\frac{1}{\LatL+i\omega}(x,y).
\end{equation}
The propagator can be expressed as the following contour integral
\begin{equation}\label{ch5:eq:3.6}
u_m(x,y;t)=\int_{\Cm_-}G_m(x,y;\omega)e^{i\omega t}\frac{\dd\omega}{2\pi}
+\int_{\Cm_+}G_m(x,y;\omega)e^{i\omega t}\frac{\dd\omega}{2\pi}.
\end{equation}
Here, $\Cm_+$ is the contour joining the point $D$ to the points $E$, $A$, $B$ in
Fig.~\ref{ch5:fig:1}, while $\Cm_-$ is the contour joining the point $B$ to $C$ to $D$. By design,
each point $\omega$ on the upper path $\Cm_+$ is separated from the spectrum of $L$.

\begin{lemma}\label{ch5:lem:1}
For $m$ sufficiently large, there is a constant $c>0$ such that
\begin{equation}\label{ch5:eq:3.7}
\left|\int_{\Cm_+}\frac{\dd\omega}{2\pi}G_m(x,y;\omega)e^{i\omega t}\right|\le c\,h^2.
\end{equation}
\end{lemma}

\begin{proof}
The proof is based on the geometric series expansion
\begin{equation}\label{ch5:eq:3.8}
G_m(\omega)=h_m^{-1}\frac{1}{\LatL+i\omega}
=h_m^{-1}\sum_{j=0}^{\infty}(-1)^j\frac{1}{\half\sigma^2\Dm+i\omega}
\left[\mu\nabm\frac{1}{\half\sigma^2\Dm+i\omega}\right]^{j}
\end{equation}
whose convergence for $\omega\in\Cm_+$ can be established by means of a Kato--Rellich relative
bound, see (Kato 1966). More precisely, for any $\alpha>0$, one can find a $\beta>0$ such that
the operators $\nabm$ and $\Dm$ satisfy the following relative bound estimate:
\begin{equation}\label{ch5:eq:3.9}
\|\nabm f\|_2\le\alpha\|\Dm f\|_2+\beta\|f\|_2.
\end{equation}
for all periodic functions $f$ and all $m\ge m_0$. This bound can be derived by observing that
$\nabm$ and $\Dm$ can be diagonalized simultaneously by a Fourier transform, as done in the
previous section, and by observing that for any $\alpha>0$, one can find a $\beta>0$ such that
\begin{equation}\label{ch5:eq:3.10}
\left|\frac{\sin h_mp}{h_m}\right|
\le\alpha\left|\frac{\cos h_mp-1}{h_m^2}\right|+\beta
\end{equation}
for all $m\ge m_0$ and all $p\in B_m$.

Under the same conditions, we also have that
\begin{equation}\label{ch5:eq:3.11}
\|\mu\nabm f\|_2\le\frac{2M\alpha}{\Sigma_0^2}
\left\|\half\sigma^2\Dm f\right\|_2+\beta\|f\|_2.
\end{equation}
Hence
\begin{equation}\label{ch5:eq:3.12}
\left\|\mu\nabm\frac{1}{\half\sigma^2\Dm+i\omega}f\right\|_2
\le\frac{2M\alpha}{\Sigma_0^2}
\left\|\half\sigma^2\Dm\frac{1}{\half\sigma^2\Dm+i\omega}f\right\|_2
+\beta\left\|\frac{1}{\half\sigma^2\Dm+i\omega}f\right\|_2<1
\end{equation}
where the last inequality holds if $\omega\in\Cm_+$, if $\alpha$ is chosen sufficiently small
and if $m$ is large enough. In this case, the geometric series expansion
\eqref{ch5:eq:3.8} converges in $L^2$ operator norm. The uniform norm of the kernel
$|G_m(x,y;\omega)|$ is pointwise controlled by the inverse of the distance from $\omega$ to the
spectrum of $\LatL$.

Let us now make that bound explicit, since this is where the height
$4|\log h_m|/t$ of the points $B$ and $D$ is used. Put $\kappa_0=4|\log h_m|/t$, so that
$e^{-\kappa_0t}=h_m^4$, and write $\omega=\sigma+i\kappa$. The contour $\Cm_+$ consists of the
two vertical segments $D\to E$ and $A\to B$, of length at most $c/h_m^2$, together with the
horizontal segment $E\to A$, which has length at most $c/h_m$ and lies at height $c/h_m^2$; on
the latter $|e^{i\omega t}|\le e^{-ct/h_m^2}$, so that its contribution is negligible. On the
vertical segments we have $|\dd\omega|=\dd\kappa$ and $\kappa\ge\kappa_0$, and therefore
$|e^{i\omega t}|\le e^{-\kappa t}$.

It remains to bound the kernel along them. The operator $\LatL$ differs by terms of order $h_m$
from the sum of a negative semidefinite symmetric part, namely the diffusion
$\half\sigma^2\Dm$, and a skew-symmetric one, namely the central difference $\mu\nabm$ of the
drift; the numerical range of $\LatL$ is therefore contained, up to $O(h_m)$, in the half plane
$\{\mathrm{Re}\le0\}$. Hence
\[
\bigl\|(\LatL+i\omega)^{-1}\bigr\|_2\le\frac{1+O(h_m)}{\mathrm{dist}\bigl(i\omega,
\{\mathrm{Re}\le0\}\bigr)}=\frac{1+O(h_m)}{\kappa},
\qquad\text{and so}\qquad
|G_m(x,y;\omega)|\le\frac{c}{h_m\kappa}
\]
for all $\omega\in\Cm_+$, the constant $c$ absorbing $\Sigma_1$ and the $O(h_m)$ terms.
Consequently
\[
\left|\int_{\Cm_+}G_m(x,y;\omega)e^{i\omega t}\frac{\dd\omega}{2\pi}\right|
\le\frac{c}{h_m}\int_{\kappa_0}^{\infty}\frac{e^{-\kappa t}}{\kappa}\dd\kappa
\le\frac{c}{h_m}\,\frac{e^{-\kappa_0t}}{\kappa_0t}
=\frac{c\,h_m^3}{4|\log h_m|}\le c\,h_m^2
\]
as soon as $h_m/(4|\log h_m|)\le1$, which is \eqref{ch5:eq:3.7}.
\end{proof}

\begin{lemma}\label{ch5:lem:2}
If $q\ge e^2\Sigma_1^2t/(2h_m^2)$ we have that
\begin{equation}\label{ch5:eq:3.13}
W_m(\gamma,q;t)\le\frac{1}{t}\sqrt{\frac{q}{2\pi}}\,
\exp\left(-\frac{\Sigma_0^2t}{2h_m^2}-q\right).
\end{equation}
\end{lemma}

\begin{proof}
Let us define the function
\begin{equation}\label{ch5:eq:3.14}
\varphi(t)=\frac{\Sigma_1^2}{2h_m^2}\,e^{-\frac{\Sigma_0^2t}{2h_m^2}}\,
\mathbf{1}(t\ge0)
\end{equation}
where $\mathbf{1}(t\ge0)$ is the characteristic function of $\mathbb{R}_+$. We have that
\begin{equation}\label{ch5:eq:3.15}
W_m(\gamma,q;t)\le\varphi^{*q}(t)
\end{equation}
where $\varphi^{*q}$ is the $q$-th convolution power, i.e.\ the $q$-fold convolution product of
the function $\varphi$ by itself. The Fourier transform of $\varphi(t)$ is given by
\begin{equation}\label{ch5:eq:3.16}
\hat\varphi(\omega)=\frac{\Sigma_1^2}{2h_m^2}\int_0^{\infty}e^{-i\omega t-\frac{\Sigma_0^2t}{2h_m^2}}\dd t
=\frac{\Sigma_1^2}{2i\omega h_m^2+\Sigma_0^2}.
\end{equation}
The convolution power is given by the following inverse Fourier transform:
\begin{equation}\label{ch5:eq:3.17}
\varphi^{*q}(t)=\int_{\rev{-\infty}}^{\infty}\hat\varphi(\omega)^qe^{i\omega t}
=\left(\frac{\Sigma_1}{\Sigma_0}\right)^{2q}\int_{-\infty}^{\infty}
\left(1+\frac{2i\omega h_m^2}{\Sigma_0^2}\right)^{-q}e^{i\omega t}\frac{\dd\omega}{2\pi}.
\end{equation}
Introducing the new variable $z=1+\frac{2i\omega h_m^2}{\Sigma_0^2}$, the integral can be
recast as follows
\begin{equation}\label{ch5:eq:3.18}
\varphi^{*q}(t)=\frac{\Sigma_0^{2-2q}\Sigma_1^{2q}}{4\pi ih_m^2}
\lim_{R\to\infty}\oint_{\Cm_R}z^{-q}\exp\left(\frac{\Sigma_0^2t}{2h_m^2}(z-1)\right)\dd z
\end{equation}
where $\Cm_R$ is the contour in Fig.~\ref{ch5:fig:2}. Using the residue theorem and noticing that
the only pole of the integrand is at $z=0$, we find
\begin{equation}\label{ch5:eq:3.19}
\varphi^{*q}(t)=\rev{\frac{1}{t}}\,\frac{1}{(q-1)!}
\left(\frac{\Sigma_1^2t}{2h_m^2}\right)^{q}
\exp\left(-\frac{\Sigma_0^2t}{2h_m^2}\right).
\end{equation}
Making use of Stirling's formula $q!\approx\sqrt{2\pi q}\,q^{q+\frac12}e^{-q}$, we find
\begin{equation}\label{ch5:eq:3.20}
\varphi^{*q}(t)\approx\rev{\frac{1}{t}}\sqrt{\frac{q}{2\pi}}\,
\exp\left(-\frac{\Sigma_0^2t}{2h_m^2}
+q\log\frac{\Sigma_1^2t}{2h_m^2}+q(1-\log q)\right).
\end{equation}
If $\log q\ge\log\frac{\Sigma_1^2t}{2h_m^2}+2$, then we arrive at the bound in \eqref{ch5:eq:3.13}.
\end{proof}

\begin{figure}[t]
\centering
\begin{tikzpicture}[scale=0.9,>=latex]
  \draw[->] (-2.6,0) -- (5.2,0) node[right] {$\mathrm{Re}(z)$};
  \draw[->] (0,-3.5) -- (0,3.5) node[above] {$\mathrm{Im}(z)$};
  \draw (1,2.9) -- (1,-2.9);
  \draw (1,2.9) arc (90:270:2.9);
  \node[left] at (-1.9,0) {$1-R$};
  \node[right] at (1.0,-3.1) {$1-iR$};
  \node[right] at (1.0,3.1) {$1+iR$};
  \node at (1.55,2.35) {$\Cm_R$};
  \fill (0,0) circle (1.3pt);
  \draw[->] (0.15,1.15) -- (0.03,0.13);
  \node[above] at (0.15,1.15) {pole};
\end{tikzpicture}
\caption{Contour of integration $\Cm_R$ for the integral in \eqref{ch5:eq:3.18}.}
\label{ch5:fig:2}
\end{figure}

It suffices to consider the case $m'=m+1$ for all values of $m$ above a fixed threshold. In
fact, given this chaptericular case, the general statement can be derived with an iterative
argument. To this end, we introduce a renormalization group transformation based on the notion
of decorating path.

% ---- from sections/ch5-decorating-paths.tex ----
% The book's figure for Definition 5.1 (decorating paths), placed by build.py immediately before
% it. The picture is the sample drawn by scripts/ch5/decorating_paths.py, which checks it against
% the definition node by node; the caption is the book's.
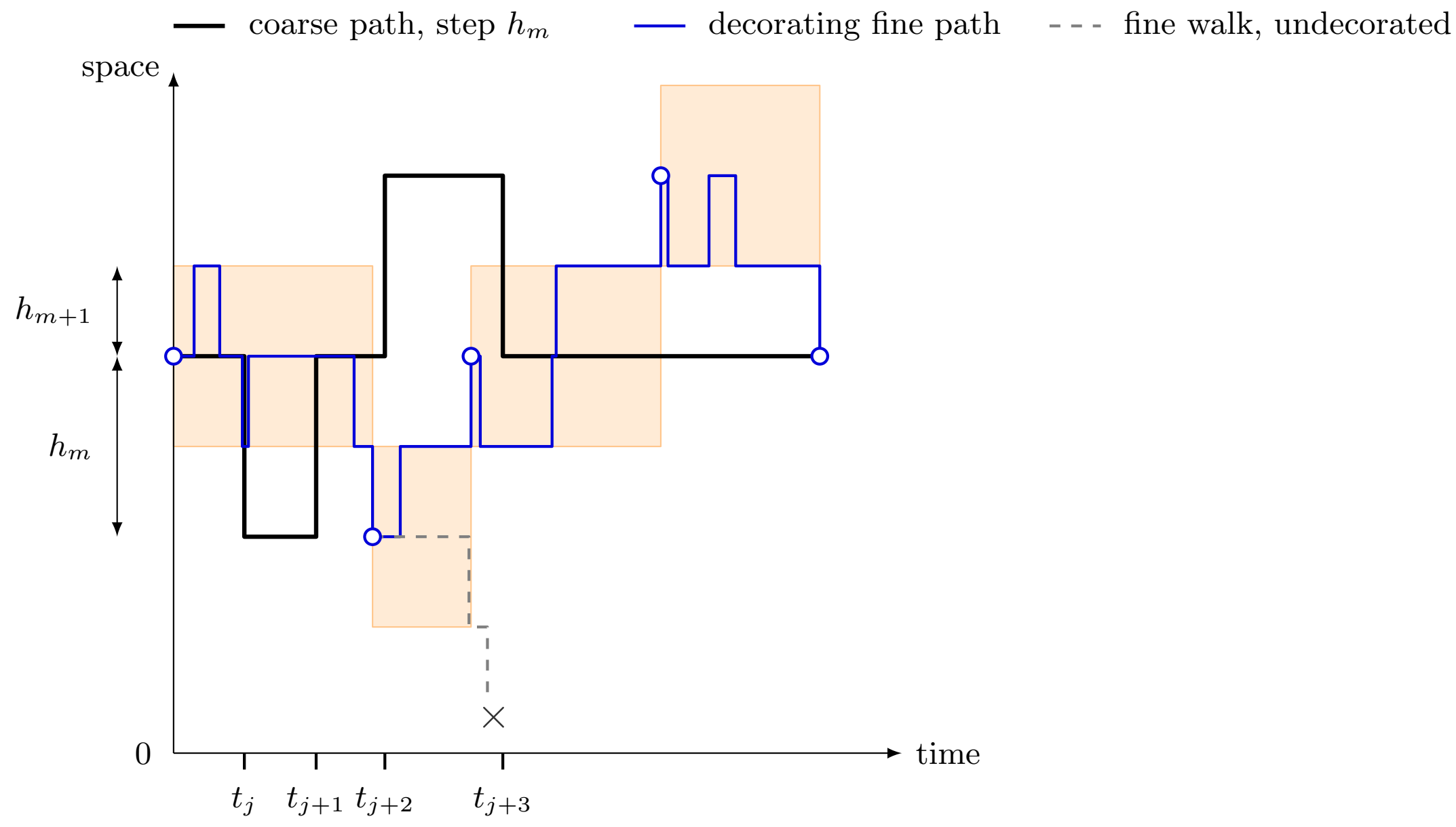
\begin{figure}[!t]
\centering
% ---- from sections/ch5-decorating-picture.tex ----
%% Generated by scripts/ch5/decorating_paths.py - do not edit.
% space is drawn at 0.88 of the time axis; the coarse walk moves 2 units of space per
%% step of the picture and the fine one 1, so h_m = 2h_{m+1} throughout.
\begin{tikzpicture}[>=latex,line join=round,scale=1.0,every node/.style={font=\footnotesize}]
  \fill[orange!16] (0.00,-0.88) rectangle (1.94,0.88);
  \draw[orange!45,thin] (0.00,-0.88) rectangle (1.94,0.88);
  \fill[orange!16] (1.94,-2.64) rectangle (2.90,-0.88);
  \draw[orange!45,thin] (1.94,-2.64) rectangle (2.90,-0.88);
  \fill[orange!16] (2.90,-0.88) rectangle (4.75,0.88);
  \draw[orange!45,thin] (2.90,-0.88) rectangle (4.75,0.88);
  \fill[orange!16] (4.75,0.88) rectangle (6.30,2.64);
  \draw[orange!45,thin] (4.75,0.88) rectangle (6.30,2.64);
  \draw[->] (0,-3.87) -- (7.10,-3.87) node[right] {time};
  \draw[->] (0,-3.87) -- (0,2.77);
  \node[left] at (0,2.77) {space};
  \draw[thick] (0.69,-3.87) -- (0.69,-4.03);
  \node[below] at (0.69,-4.07) {$t_{j}$};
  \draw[thick] (1.39,-3.87) -- (1.39,-4.03);
  \node[below] at (1.39,-4.07) {$t_{j+1}$};
  \draw[thick] (2.06,-3.87) -- (2.06,-4.03);
  \node[below] at (2.06,-4.07) {$t_{j+2}$};
  \draw[thick] (3.21,-3.87) -- (3.21,-4.03);
  \node[below] at (3.21,-4.07) {$t_{j+3}$};
  \node[left] at (-0.08,-3.87) {$0$};
  \draw[very thick] (0.00,0.00) -- (0.69,0.00) -- (0.69,-1.76) -- (1.39,-1.76) -- (1.39,0.00) -- (2.06,0.00) -- (2.06,1.76) -- (3.21,1.76) -- (3.21,0.00) -- (6.30,0.00);
  \draw[thick,blue!85!black] (0.00,0.00) -- (0.20,0.00) -- (0.20,0.88) -- (0.45,0.88) -- (0.45,0.00) -- (0.67,0.00) -- (0.67,-0.88) -- (0.73,-0.88) -- (0.73,0.00) -- (1.76,0.00) -- (1.76,-0.88) -- (1.94,-0.88) -- (1.94,-1.76) -- (2.21,-1.76) -- (2.21,-0.88) -- (2.90,-0.88) -- (2.90,0.00) -- (2.99,0.00) -- (2.99,-0.88) -- (3.69,-0.88) -- (3.69,0.00) -- (3.73,0.00) -- (3.73,0.88) -- (4.75,0.88) -- (4.75,1.76) -- (4.82,1.76) -- (4.82,0.88) -- (5.22,0.88) -- (5.22,1.76) -- (5.48,1.76) -- (5.48,0.88) -- (6.30,0.88) -- (6.30,0.00);
  \draw[thick,dashed,gray] (1.94,-1.76) -- (2.88,-1.76) -- (2.88,-2.64) -- (3.06,-2.64) -- (3.06,-3.52);
  \node[gray!40!black,font=\large,fill=white,inner sep=1pt] at (3.12,-3.52) {$\times$};
  \fill[white] (0.00,0.00) circle (2.2pt);
  \draw[thick,blue!85!black] (0.00,0.00) circle (2.2pt);
  \fill[white] (1.94,-1.76) circle (2.2pt);
  \draw[thick,blue!85!black] (1.94,-1.76) circle (2.2pt);
  \fill[white] (2.90,0.00) circle (2.2pt);
  \draw[thick,blue!85!black] (2.90,0.00) circle (2.2pt);
  \fill[white] (4.75,1.76) circle (2.2pt);
  \draw[thick,blue!85!black] (4.75,1.76) circle (2.2pt);
  \fill[white] (6.30,0.00) circle (2.2pt);
  \draw[thick,blue!85!black] (6.30,0.00) circle (2.2pt);
  \draw[<->] (-0.55,-1.76) -- (-0.55,0.00);
  \node[left] at (-0.65,-0.88) {$h_m$};
  \draw[<->] (-0.55,0.00) -- (-0.55,0.88);
  \node[left] at (-0.65,0.44) {$h_{m+1}$};
  \draw[very thick,black] (0.00,3.22) -- (0.50,3.22);
  \node[anchor=west] at (0.60,3.22) {coarse path, step $h_m$};
  \draw[thick,blue!85!black] (4.49,3.22) -- (4.99,3.22);
  \node[anchor=west] at (5.08,3.22) {decorating fine path};
  \draw[thick,dashed,gray] (8.54,3.22) -- (9.04,3.22);
  \node[anchor=west] at (9.13,3.22) {fine walk, undecorated};
\end{tikzpicture}
% ---- end sections/ch5-decorating-picture.tex ----
\caption{A reading of Definition~\ref{ch5:def:1}, drawn from a simulation of the two walks rather
than by hand, so that what is shown is a typical decoration and not only a possible one; the
script is in the Computational Appendix, and the sample is one of many. Both sequences are paths
which are constant in time and jump at their own steps, those steps being the Poisson times at
which the corresponding lattice walk moves; the coarse path steps at the times $t_j$,
$t_{j+1}$, \dots{} marked on the axis, the fine one at four times the rate, since its rate scales
as $h^{-2}$, so it steps about four times as often. A decoration must contain the coarse values in
order, and between the time $y'_j=y_i$ at which it takes the value of a coarse node and the time
$y'_k=y_{i+1}$ at which it takes the next one, it must stay within one fine cell of $y'_j$: those
are the shaded bands. They begin and end at the times of the \emph{fine} path, not at those of the
coarse one, so the touches fall at times of the decoration's own: in the last band the fine path
takes the value $y_i=2$ long after the coarse path has left it. Two consecutive bands meet in a
single fine position, the one midway between the two coarse values, so a decoration is forced
through it, and the circles mark the touches. The dashed path is the same fine walk \emph{without}
the constraint: it leaves the band on the wrong side, and it is the conditioning, not the walk,
which produces the geometry of a decoration.}
\label{fig:ch5:decorating}
\end{figure}
% ---- end sections/ch5-decorating-paths.tex ----

\begin{definition}[Decorating Paths]\label{ch5:def:1}
Let $m\ge m_0$ and let $\gamma=\{y_0,y_1,y_2,\dots\}$ be a symbolic sequence in $\Gm$. A
\emph{decorating path} around $\gamma$ is defined as a symbolic sequence
$\gamma'=\{y'_0,y'_1,y'_2,\dots\}$ with $y'_i\in h_{m+1}\mathbb{Z}$ containing the sequence
$\gamma$ as a subset and such that if $y'_j=y_i$ and $y'_k=y_{i+1}$, then all elements $y'_n$
with $j<n<k$ are such that $|y'_n-y'_j|\le h_{m+1}$. Let $D_{m+1}(\gamma)$ be the set of all
decorating sequences around $\gamma$. The decorated weights are defined as follows:
\begin{equation}\label{ch5:eq:3.21}
\tilde W_m(\gamma,q;t)=\sum_{q'=q}^{\infty}\sum_{\substack{\gamma'\in D_{m+1}(\gamma)\\
\gamma'_{q'}=\gamma_q}}W_{m+1}(\gamma',q';t).
\end{equation}
Finally, let us introduce also the following Fourier transform:
\begin{equation}\label{ch5:eq:3.22}
\hat W_m(\gamma,q;\omega)=\int_0^{\infty}W_m(\gamma,q;t)e^{i\omega t}\dd t,
\qquad
\hat{\tilde W}_m(\gamma,q;\omega)=\int_0^{\infty}\tilde W_m(\gamma,q;t)e^{i\omega t}\dd t.
\end{equation}
\end{definition}

\begin{definition}[Notations]\label{ch5:def:2}
In the following, we set $h=h_{m+1}$ so that $h_m=2h$. We also use the Landau notation
$O(h^n)$ to indicate a function $f(h)$ such that $h^{-n}f(h)$ is bounded in a neighborhood of
$(0)$.
\end{definition}

\begin{lemma}\label{ch5:lem:3}
Let $x,y\in\Domain_m$ and let $\Cm_-$ be an integration contour as in Fig.~\ref{ch5:fig:1}. Then
\begin{equation}\label{ch5:eq:3.23}
\left(\int_{\Cm_-}\Bigl(2G_{m+1}(x,y;\omega)-G_m(x,y;\omega)\Bigr)e^{i\omega t}
\frac{\dd\omega}{2\pi}\right)=O(h^2).
\end{equation}
\end{lemma}

\begin{proof}
Recall that $h=h_{m+1}$ so that $h_m=2h$. By \eqref{ch5:eq:3.3} and \eqref{ch5:eq:3.21}, the coarse
resolvent enters with the factor $h_m^{-1}=(2h)^{-1}$ and the fine one with $h_{m+1}^{-1}=h^{-1}$,
which is where the two coefficients in \eqref{ch5:eq:3.24} come from. We have that
\begin{equation}\label{ch5:eq:3.24}
2G_{m+1}(x,y;\omega)-G_m(x,y;\omega)=\frac{1}{h}\sum_{q=1}^{\infty}2^{-q}\!\!\!
\sum_{\substack{\gamma\in\Gm:\ \gamma_0=x,\ \gamma_q=y\\ |\gamma_j-\gamma_{j-1}|=1\ \forall j\ge1}}
\!\!\!\Bigl(2\hat{\tilde W}_m(\gamma,q;\omega)-\tfrac12\hat W_m(\gamma,q;\omega)\Bigr).
\end{equation}
The number of paths over which the summation is extended is
\begin{equation}\label{ch5:eq:3.25}
N(\gamma,q;x,y)\equiv\#\{\gamma\in\Gm:\gamma_0=x,\gamma_q=y,|\gamma_j-\gamma_{j-1}|=1\ \forall j\ge1\}
=\binom{q}{\frac{q}{2}+k}
\end{equation}
where $k=\frac{|y-x|}{h_m}$. Applying Stirling's formula we find
\begin{equation}\label{ch5:eq:3.26}
N_\gamma\sim2^q\sqrt{\frac{2}{\pi q}}.
\end{equation}
Hence
\begin{equation}\label{ch5:eq:3.27}
\begin{aligned}
&\left|\int_{\Cm_-}\Bigl(2G_{m+1}(x,y;\omega)-G_m(x,y;\omega)\Bigr)e^{i\omega t}
\frac{\dd\omega}{2\pi}\right|\\
&\quad\le\frac{c}{h}\sum_{q=1}^{\infty}\frac{1}{q}\,
\max_{\substack{\gamma\in\Gm:\ \gamma_0=x,\ \gamma_q=y\\ |\gamma_j-\gamma_{j-1}|=1\ \forall j\ge1}}
\left|\int_{\Cm_-}\Bigl(2\hat{\tilde W}_m(\gamma,q;\omega)-\tfrac12\hat W_m(\gamma,q;\omega)\Bigr)
e^{i\omega t}\frac{\dd\omega}{2\pi}\right|.
\end{aligned}
\end{equation}
for some constant $c\approx\sqrt{\frac{2}{\pi}}>0$. It suffices to extend the summation over
$q$ only up to
\begin{equation}\label{ch5:eq:3.28}
q_{\max}\equiv\frac{e^2\Sigma_1^2t}{2h^2}.
\end{equation}
To resum beyond this threshold, one can use the previous lemma. More precisely, we have that
\begin{equation}\label{ch5:eq:3.29}
\begin{aligned}
&\left|\int_{\Cm_-}\Bigl(2G_{m+1}(x,y;\omega)-G_m(x,y;\omega)\Bigr)e^{i\omega t}
\frac{\dd\omega}{2\pi}\right|\\
&\quad\le\frac{c\sqrt{q_{\max}}}{h}\,
\max_{\substack{\gamma\in\Gm:\ \gamma_0=x,\ \gamma_q=y\\ |\gamma_j-\gamma_{j-1}|=1\ \forall j\ge1}}
\left|\int_{\Cm_-}\Bigl(2\hat{\tilde W}_m(\gamma,q;\omega)-\tfrac12\hat W_m(\gamma,q;\omega)\Bigr)
e^{i\omega t}\frac{\dd\omega}{2\pi}\right|.
\end{aligned}
\end{equation}

Let $v(x)=\sigma(x)^2$. To evaluate the resummed weight function, let us form the matrix
\begin{equation}\label{ch5:eq:3.30}
\LbarRaw(x;h)=
\begin{pmatrix}
-\dfrac{v(x+h)}{h^2} & \dfrac{v(x+h)}{2h^2}-\dfrac{\mu(x+h)}{2h} & 0\\[8pt]
\dfrac{v(x)}{2h^2}+\dfrac{\mu(x)}{2h} & -\dfrac{v(x)}{h^2} & \dfrac{v(x)}{2h^2}-\dfrac{\mu(x)}{2h}\\[8pt]
0 & \dfrac{v(x-h)}{2h^2}+\dfrac{\mu(x-h)}{2h} & -\dfrac{v(x-h)}{h^2}
\end{pmatrix}
\end{equation}
and decompose it as follows:
\begin{equation}\label{ch5:eq:3.31}
\LbarRaw(x;h)=\frac{1}{h^2}\LbarRaw_0(x)+\frac{1}{h}\LbarRaw_1(x)+\LbarRaw_2(x)+h\LbarRaw_3(x)+O(h^2).
\end{equation}
where
\begin{equation}\label{ch5:eq:3.32}
\LbarRaw_0(x)=\begin{pmatrix}
-v(x) & \half v(x) & 0\\
\half v(x) & -v(x) & \half v(x)\\
0 & \half v(x) & -v(x)
\end{pmatrix},
\end{equation}
\begin{equation}\label{ch5:eq:3.33}
\LbarRaw_1(x)=\begin{pmatrix}
-v'(x) & \half v'(x)-\half\mu(x) & 0\\
\half\mu(x) & 0 & -\half\mu(x)\\
0 & -\half v'(x)+\half\mu(x) & v'(x)
\end{pmatrix},
\end{equation}
\begin{equation}\label{ch5:eq:3.34}
\LbarRaw_2(x)=\begin{pmatrix}
-\half v''(x) & \tfrac14 v''(x)-\half\mu'(x) & 0\\
0 & 0 & 0\\
0 & \tfrac14 v''(x)-\half\mu'(x) & -\half v''(x)
\end{pmatrix}
\end{equation}
and
\begin{equation}\label{ch5:eq:3.35}
\LbarRaw_3(x)=\begin{pmatrix}
-\tfrac16 v'''(x) & \tfrac1{12} v'''(x)-\tfrac14\mu''(x) & 0\\
0 & 0 & 0\\
0 & -\tfrac1{12} v'''(x)+\tfrac14\mu''(x) & \tfrac16 v'''(x)
\end{pmatrix}.
\end{equation}

Let us introduce the sign variable $\tau=\pm1$, the functions
\begin{equation}\label{ch5:eq:3.36}
\varphi_0(t,x,\tau)\equiv2\LatL(x,x+2\tau h)\,e^{t\LatL(x,x)}\,\mathbf{1}(t\ge0)
\end{equation}
\begin{equation}\label{ch5:eq:3.37}
\varphi_1(t,x,\tau)\equiv2L^{m+1}(x+\tau h,x+2\tau h)\,e^{t\LbarRaw(x;h)}(x,x+\tau h)\,\mathbf{1}(t\ge0)
\end{equation}
and their Fourier transforms
\begin{equation}\label{ch5:eq:3.38}
\begin{aligned}
\hat\varphi_0(\omega,x,\tau)&=\left(\frac{v(x)}{4h^2}+\tau\frac{\mu(x)}{2h}\right)
\left(\frac{v(x)}{4h^2}+i\omega\right)^{-1}\\
\hat\varphi_1(\omega,x,\tau)&=\left(\frac{v(x)}{h^2}+\tau\frac{\mu(x)+v'(x)}{h}
+\frac{v''(x)+\mu'(x)}{2}+\left(\frac{v'''(x)}{6}+\frac{\mu''(x)}{2}\right)\tau h+O(h^2)\right)\\
&\qquad\times\bigl\langle x\,\bigl|(-\LbarRaw(x;h)+i\omega)^{-1}\bigr|x+\tau h\bigr\rangle
\end{aligned}
\end{equation}
where
\begin{equation}\label{ch5:eq:3.39}
|x\rangle=\begin{pmatrix}0\\1\\0\end{pmatrix},
\qquad\text{and}\qquad
|x+\tau h\rangle=\begin{pmatrix}\delta_{\tau,1}\\0\\\delta_{\tau,-1}\end{pmatrix}.
\end{equation}
We also require the functions
\begin{equation}\label{ch5:eq:3.40}
\psi_0(t,x)\equiv e^{t\LatL(x,x)}\mathbf{1}(t\ge0),\qquad
\psi_1(t,x)\equiv e^{t\LbarRaw(y;h)}(x,x)\mathbf{1}(t\ge0)
\end{equation}
and the corresponding Fourier transforms
\begin{equation}\label{ch5:eq:3.41}
\hat\psi_0(\omega,x)=\left(\frac{v(x)}{4h^2}+i\omega\right)^{-1},\qquad
\hat\psi_1(\omega,x)=\bigl\langle x\bigl|(-\LbarRaw(x;h)+i\omega)^{-1}\bigr|x\bigr\rangle.
\end{equation}
If $\gamma$ is a symbolic sequence, then
\begin{equation}\label{ch5:eq:3.42}
\hat W_m(\gamma,q;\omega)=\hat\psi_0(\omega,\gamma_q)\prod_{j=0}^{q-1}
\hat\varphi_0\bigl(\omega;\gamma_j,\operatorname{sgn}(\gamma_{j+1}-\gamma_j)\bigr)
\end{equation}
\begin{equation}\label{ch5:eq:3.43}
\hat{\tilde W}_m(\gamma,q;\omega)=\hat\psi_1(\omega,\gamma_q)\prod_{j=0}^{q-1}
\hat\varphi_1\bigl(\omega;\gamma_j,\operatorname{sgn}(\gamma_{j+1}-\gamma_j)\bigr).
\end{equation}

Let us estimate the difference between the functions $\hat\varphi_1(\omega,x,\tau)$ and
$\hat\varphi_0(\omega,x,\tau)$ assuming that $\omega$ is in the contour $\Cm_-$ in
Fig.~\ref{ch5:fig:2}. Retaining only terms up to order up to $O(h^3)$, we find
\begin{equation}\label{ch5:eq:3.44}
\hat\varphi_0(\omega,x,\tau)=1+\frac{2\mu(x)\tau h}{v(x)}-\frac{4i\omega h^2}{v(x)}
-8\mu(x)\frac{i\omega\tau h^3}{v(x)^2}-\frac{16\omega^2h^4}{v(x)^2}+O(h^5).
\end{equation}
A lengthy but straightforward calculation which is best carried out using a symbolic
manipulation program, gives
\begin{equation}\label{ch5:eq:3.45}
\begin{aligned}
\hat\varphi_1(\omega,x,\tau)={}&1+\frac{2\mu(x)\tau h}{v(x)}-\frac{4i\omega h^2}{v(x)}
-\bigl[8\mu(x)-v'(x)\bigr]\frac{i\omega\tau h^3}{v(x)^2}\\
&+r(x)\cdot h^3\tau+i\omega h^4p(x)-\frac{14\omega^2h^4}{v(x)^2}+O(h^5)
\end{aligned}
\end{equation}
where
\begin{equation}\label{ch5:eq:3.46}
\begin{aligned}
r(x)&=\frac{1}{2v(x)^3}\Bigl[\rev{\mu''(x)v(x)^2}-4\mu(x)^3+2v'(x)\mu(x)^2
-2v'(x)v(x)\mu'(x)\\
&\qquad\qquad\qquad-\bigl(\rev{2v(x)\mu'(x)}+v''(x)v(x)-2v'(x)^2\bigr)\mu(x)\Bigr],\\[4pt]
p(x)&=\frac{1}{v(x)^3}\Bigl[4\mu(x)^2-2v'(x)\mu(x)+4v(x)\mu'(x)+v''(x)v(x)-2v'(x)^2\Bigr].
\end{aligned}
\end{equation}
We have that
\begin{equation}\label{ch5:eq:3.47}
\begin{aligned}
&\sum_{j=0}^{q-1}\Bigl[\log\hat\varphi_0\bigl(\omega;\gamma_j,
\operatorname{sgn}(\gamma_{j+1}-\gamma_j)\bigr)
-\log\hat\varphi_1\bigl(\omega;\gamma_j,\operatorname{sgn}(\gamma_{j+1}-\gamma_j)\bigr)\Bigr]\\
&\quad=-\sum_{j=0}^{q-1}\left[\frac{i\omega v'(\gamma_j)}{v(\gamma_j)^2}+r(\gamma_j)\right]
h^3\operatorname{sgn}(\gamma_{j+1}-\gamma_j)
+\bigl(|\omega|\|p\|_\infty+2|\omega|^2\|v^{-2}\|_\infty\bigr)O(h^4q)\\
&\quad=\rev{-h^2\left[i\omega\left(\frac{1}{v(\gamma_0)}-\frac{1}{v(\gamma_q)}\right)
+\bigl(R(\gamma_q)-R(\gamma_0)\bigr)\right]}
+\bigl(|\omega|\|p\|_\infty+2|\omega|^2\|v^{-2}\|_\infty\bigr)O(h^4q)
\end{aligned}
\end{equation}
where $R(x)$ is a primitive of $r(x)$, i.e.
\begin{equation}\label{ch5:eq:3.48}
R(x)=\int_0^xr(z)\dd z.
\end{equation}

We conclude that there is a constant $c>0$ such that
\begin{equation}\label{ch5:eq:3.49}
\int_{\Cm_-}\left|\prod_{j=0}^{q-1}
\hat\varphi_0\bigl(\omega;\gamma_j,\operatorname{sgn}(\gamma_{j+1}-\gamma_j)\bigr)
-\prod_{j=0}^{q-1}
\hat\varphi_1\bigl(\omega;\gamma_j,\operatorname{sgn}(\gamma_{j+1}-\gamma_j)\bigr)\right|
e^{i\omega t}\frac{\dd\omega}{2\pi}\le c\,h^2.
\end{equation}
for all $q\le q_{\max}$. Here we use the decay of $e^{i\omega t}$ in the upper half of the
complex $\omega$ plane to offset the $\omega$ dependencies in the integrand. Similar
calculations lead to the following expansions:
\begin{equation}\label{ch5:eq:3.50}
\hat\psi_0(\omega,y)=\frac{4h^2}{v(y)}+O(\omega h^4),\qquad
\hat\psi_1(\omega,y)=\frac{2h^2}{v(y)}+O(\omega h^4)=\half\hat\psi_0(\omega,y)+O(\omega h^4).
\end{equation}
Since $q<ch^{-2}$ and $\omega\le|\log h|$, we find
\begin{equation}\label{ch5:eq:3.51}
\left|\int_{\Cm_-}\Bigl(2G_{m+1}(x,y;\omega)-G_m(x,y;\omega)\Bigr)e^{i\omega t}
\frac{\dd\omega}{2\pi}\right|\le c\,\frac{\sqrt{q_{\max}}}{h}h^4\le c\,h^2.
\end{equation}
This completes the proof of the Lemma and of the Theorem.
\end{proof}

By differentiating with respect to time in equation \eqref{ch5:eq:3.52}, we find that
\begin{equation}\label{ch5:eq:3.52}
\frac{\partial}{\partial t}u_m(x,y;t)=\int_{\Cm_-}i\omega G_m(x,y;\omega)e^{i\omega t}
\frac{\dd\omega}{2\pi}+\int_{\Cm_+}i\omega G_m(x,y;\omega)e^{i\omega t}\frac{\dd\omega}{2\pi}.
\end{equation}
All the derivations above carry through and we conclude that
\begin{equation}\label{ch5:eq:3.53}
\left|\int_{\Cm_+}i\omega G_m(x,y;\omega)e^{i\omega t}\frac{\dd\omega}{2\pi}\right|\le c\,h^2.
\end{equation}
and also
\begin{equation}\label{ch5:eq:3.54}
\left|\int_{\Cm_-}i\omega\Bigl(2G_{m+1}(x,y;\omega)-G_m(x,y;\omega)\Bigr)e^{i\omega t}
\frac{\dd\omega}{2\pi}\right|\le c\,\frac{\sqrt{q_{\max}}}{h}h^4\le c\,h^2.
\end{equation}
Hence the first time derivatives of the kernel satisfy the same Cauchy convergence condition as
the kernel itself.

\section{Lesser Smooth Coefficients}

In this section we assume coefficients are either H\"older continuous or obey the conditions in
Theorem 2.

\begin{lemma}\label{ch5:lem:4}
Let $f(x)$ be a continuous function in $[-L,L]$ satisfying periodic boundary conditions. Then,
for all $h>0$, we have that
\begin{equation}\label{ch5:eq:4.1}
f(x+h_m)=f(x)+h_m\nabm_xf(x)+\frac{h_m^2}{2}\Dm_xf(x)
\end{equation}
\begin{equation}\label{ch5:eq:4.2}
f(x-h_m)=f(x)-h_m\nabm_xf(x)+\frac{h_m^2}{2}\Dm_xf(x).
\end{equation}
\end{lemma}

This is the result of a simple calculation, which is however useful as it allows one to extend
the derivation in the previous section by making the following replacements:
\begin{equation}\label{ch5:eq:4.3}
v'(x)\to\nabm_xv(x),\qquad v''(x)\to\Dm_xv(x),\qquad v'''(x)\to0
\end{equation}
\begin{equation}\label{ch5:eq:4.4}
\mu'(x)\to\nabm_x\mu(x),\qquad \mu''(x)\to\Dm_x\mu(x).
\end{equation}
In fact,
\begin{equation}\label{ch5:eq:4.5}
\LbarRaw(x;h)=\frac{1}{h^2}\LbarRaw_0(x)+\frac{1}{h}\LbarRaw_1(x)+\LbarRaw_2(x)+h\LbarRaw_3(x)
\end{equation}
without any $O(h^3)$ corrections as long as one re-defines the matrices on the right hand side
as follows:
\begin{equation}\label{ch5:eq:4.6}
\LbarRaw_0(x)=\begin{pmatrix}
-v(x) & \half v(x) & 0\\
\half v(x) & -v(x) & \half v(x)\\
0 & \half v(x) & -v(x)
\end{pmatrix},
\end{equation}
\begin{equation}\label{ch5:eq:4.7}
\LbarRaw_1(x)=\begin{pmatrix}
-\nabm_xv(x) & \half\nabm_xv(x)-\half\mu(x) & 0\\
\half\mu(x) & 0 & -\half\mu(x)\\
0 & -\half\nabm_xv(x)+\half\mu(x) & \nabm_xv(x)
\end{pmatrix},
\end{equation}
\begin{equation}\label{ch5:eq:4.8}
\LbarRaw_2(x)=\begin{pmatrix}
-\half\Dm_xv(x) & \tfrac14\Dm_xv(x)-\half\nabm_x\mu(x) & 0\\
0 & 0 & 0\\
0 & \tfrac14\Dm_xv(x)-\half\nabm_x\mu(x) & -\half\Dm_xv(x)
\end{pmatrix}
\end{equation}
and
\begin{equation}\label{ch5:eq:4.9}
\LbarRaw_3(x)=\begin{pmatrix}
0 & -\tfrac14\Dm_x\mu(x) & 0\\
0 & 0 & 0\\
0 & \tfrac14\Dm_x\mu(x) & 0
\end{pmatrix}.
\end{equation}
All derivations in the previous section go through formally unchanged and one arrives at the
following expressions
\begin{equation}\label{ch5:eq:4.10}
\hat\varphi_0(\omega,x,\tau)=1+\frac{2\mu(x)\tau h}{v(x)}-\frac{4i\omega h^2}{v(x)}
-8\mu(x)\frac{i\omega\tau h^3}{v(x)^2}-\frac{16\omega^2h^4}{v(x)^2}+O(h^5).
\end{equation}
and
\begin{equation}\label{ch5:eq:4.11}
\begin{aligned}
\hat\varphi_1(\omega,x,\tau)={}&1+\frac{2\mu(x)\tau h}{v(x)}-\frac{4i\omega h^2}{v(x)}
-\bigl[8\mu(x)-\nabm_xv(x)\bigr]\frac{i\omega\tau h^3}{v(x)^2}\\
&+r(x)\cdot h^3\tau+i\omega h^4p(x)-\frac{14\omega^2h^4}{v(x)^2}+O(h^5)
\end{aligned}
\end{equation}
where
\begin{equation}\label{ch5:eq:4.12}
\begin{aligned}
r(x)&=\frac{1}{2v(x)^3}\Bigl[\rev{\Dm_x\mu(x)\,v(x)^2}-4\mu(x)^3+2\nabm_xv(x)\mu(x)^2
-2\nabm_xv(x)v(x)\nabm_x\mu(x)\\
&\qquad\qquad\qquad-\bigl(\rev{2v(x)\nabm_x\mu(x)}+\Dm_xv(x)v(x)-2\nabm_xv(x)^2\bigr)\mu(x)\Bigr],\\[4pt]
p(x)&=\frac{1}{v(x)^3}\Bigl[4\mu(x)^2-2\nabm_xv(x)\mu(x)+4v(x)\nabm_x\mu(x)
+\Dm_xv(x)v(x)-2\nabm_xv(x)^2\Bigr].
\end{aligned}
\end{equation}

We have that
\begin{equation}\label{ch5:eq:4.13}
\begin{aligned}
&\sum_{j=0}^{q-1}\Bigl[\log\hat\varphi_0\bigl(\omega;\gamma_j,
\operatorname{sgn}(\gamma_{j+1}-\gamma_j)\bigr)
-\log\hat\varphi_1\bigl(\omega;\gamma_j,\operatorname{sgn}(\gamma_{j+1}-\gamma_j)\bigr)\Bigr]\\
&\quad=-\sum_{j=0}^{q-1}\left[\frac{i\omega\nabm_xv(\gamma_j)}{v(\gamma_j)^2}+r(\gamma_j)\right]
h^3\operatorname{sgn}(\gamma_{j+1}-\gamma_j)
+\bigl(|\omega|\|p\|_\infty+2|\omega|^2\|v^{-2}\|_\infty\bigr)O(h^4q)\\
&\quad\le\rev{h^2}\left[2L\sup_{x\in\Domain_m}
\left|\frac{i\omega\nabm_xv(x)}{v(x)^2}+r(x)\right|
+\bigl(|\omega|\|p\|_\infty+2|\omega|^2\|v^{-2}\|_\infty\bigr)O(h^2q)\right]
\le c\,h^{\gamma},
\end{aligned}
\end{equation}
where in the last step we made use of the H\"older continuity assumptions of Theorem 1. The
other bounds staying the same, we arrive at
\begin{equation}\label{ch5:eq:4.14}
\left|\int_{\Cm_-}\Bigl(2G_{m+1}(x,y;\omega)-G_m(x,y;\omega)\Bigr)e^{i\omega t}
\frac{\dd\omega}{2\pi}\right|\le c\,\frac{\sqrt{q_{\max}}}{h}h^{2+\gamma}\le c\,h^{\gamma}.
\end{equation}
Under the weaker assumption of Theorem 2, the bound that applies is instead
\begin{equation}\label{ch5:eq:4.15}
\left|\int_{\Cm_-}\Bigl(2G_{m+1}(x,y;\omega)-G_m(x,y;\omega)\Bigr)e^{i\omega t}
\frac{\dd\omega}{2\pi}\right|\le c\,\frac{\sqrt{q_{\max}}}{h}h^2\rho(h)\le c\,\rho(h).
\end{equation}

Similar bounds also extend to the case of the first time derivative, since multiplication by a
factor $i\omega$ inside of the contour integral is immaterial as far as establishing a bound of
this sort is concerned. This completes the proof of Theorem 1 and Theorem 2.

\section{Explicit Euler Scheme}\label{ch5:sec:euler}

\begin{figure}[t]
\centering
\begin{tikzpicture}[scale=0.9,>=latex]
  \draw[->] (-5.8,0) -- (5.8,0) node[right] {$\mathrm{Re}(\omega)$};
  \draw[->] (0,-1.2) -- (0,4.8) node[above] {$\mathrm{Im}(\omega)$};
  \draw (0,0) rectangle (4.0,3.6);
  \node[below left] at (-4.0,0) {$A$};
  \node[below] at (-4.0,-0.35) {$-\pi/\delta t$};
  \node[below right] at (4.0,0) {$B$};
  \node[below] at (4.0,-0.35) {$\pi/\delta t$};
  \node[above right] at (4.0,3.6) {$C$};
  \node[above] at (4.0,3.95) {$iH+\pi/\delta t$};
  \node[above left] at (-4.0,3.6) {$D$};
  \node[above] at (-4.0,3.95) {$iH-\pi/\delta t$};
  \node at (1.4,3.0) {$c/h$};
\end{tikzpicture}
\caption{Contour of integration for the integral in \eqref{ch5:eq:5.4}.}
\label{ch5:fig:3}
\end{figure}
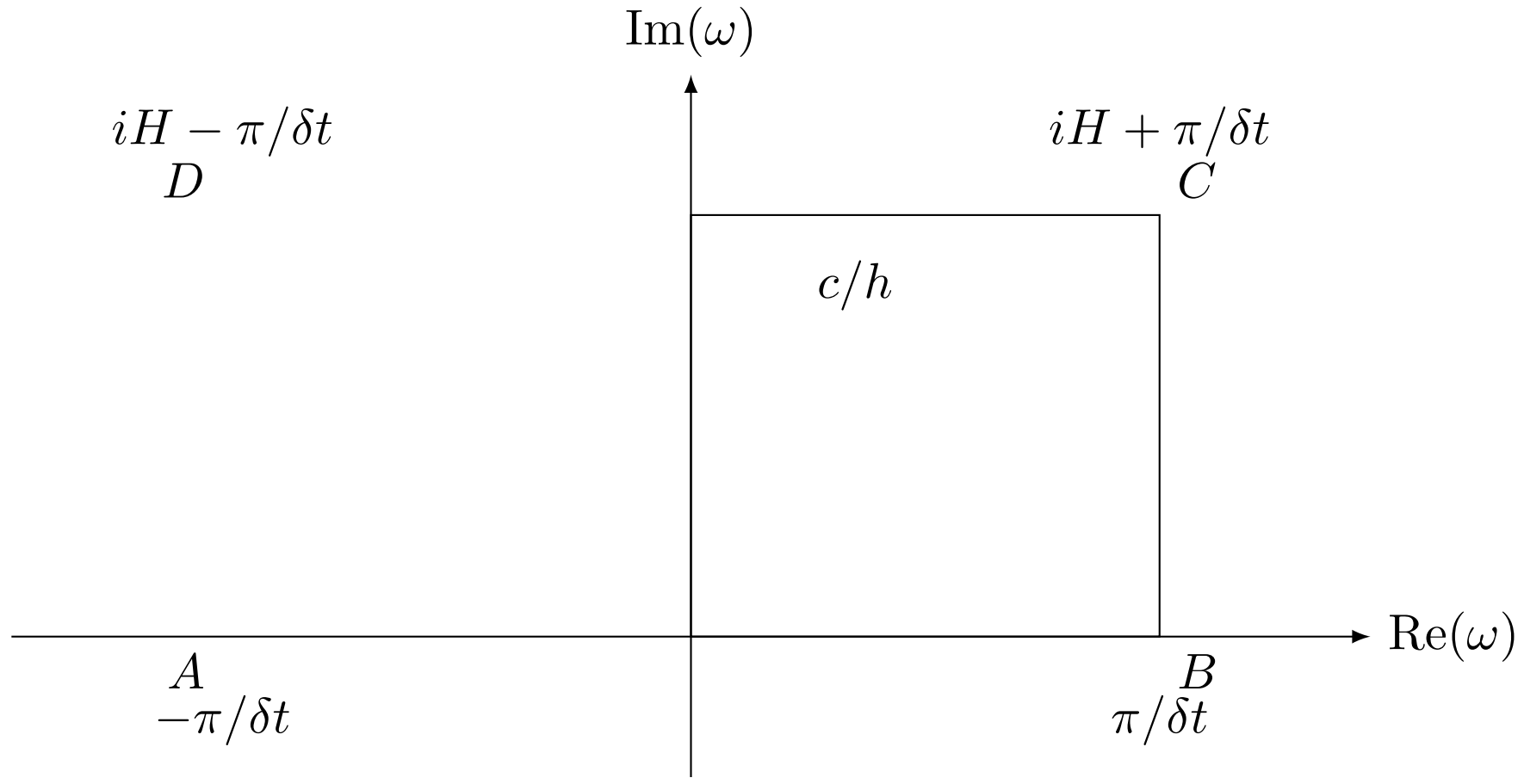

In this section we prove Theorem 3. A Dyson expansion can also be obtained for the
time-discretized kernel and has the form
\begin{equation}\label{ch5:eq:5.1}
\begin{aligned}
u^{\delta t}_m(y_1,y_2;t)={}&\frac{1}{h_m}\sum_{q=1}^{\infty}
\sum_{\gamma\in\Gm:\gamma_0=y_1,\gamma_q=y_2}
\sum_{k_1=1}^{N}\sum_{k_2=k_1+1}^{N}\cdots\sum_{k_q=k_{q-1}+1}^{N}\\
&\times\bigl(1+\delta t\LatL(\gamma_0,\gamma_0)\bigr)^{k_1-1}
(\delta t)^q\prod_{j=1}^{q}\LatL(\gamma_{j-1},\gamma_j)
\bigl(1+\delta t\LatL(\gamma_j,\gamma_j)\bigr)^{k_{j+1}-k_j-1}
\end{aligned}
\end{equation}
where $t_{q+1}=t$ and $k_{q+1}=N$. In this case, the propagator can be expressed through a
Fourier integral as follows:
\begin{equation}\label{ch5:eq:5.2}
u^{\delta t}_m(y_1,y_2;t)=\int_{-\frac{\pi}{\delta t}}^{\frac{\pi}{\delta t}}
G^{\delta t}_m(y_1,y_2;\omega)e^{i\omega t}\frac{\dd\omega}{2\pi}
\end{equation}
where
\begin{equation}\label{ch5:eq:5.3}
G^{\delta t}_m(y_1,y_2;\omega)=\delta t\sum_{j=0}^{\frac{t}{\delta t}}
u^{\delta t}_m(y_1,y_2;j\delta t)e^{-i\omega j\delta t}.
\end{equation}
The propagator can also be represented as the limit
\begin{equation}\label{ch5:eq:5.4}
u^{\delta t}_m(y_1,y_2;t)=\lim_{H\to\infty}\int_{\Cm_H}
G^{\delta t}_m(y_1,y_2;\omega)e^{i\omega t}\frac{\dd\omega}{2\pi}
\end{equation}
where $\Cm_H$ is the contour in Fig.~\ref{ch5:fig:3}. This is due to the fact that the integral
along the segments $BC$ and $DA$ are the negative of each other, while the integral over $CD$
tends to zero exponentially fast as $\Im(\omega)\to\infty$, where $\Im(\omega)$ is the
imaginary part of $\omega$. Using Cauchy's theorem, the contour in Fig.~\ref{ch5:fig:3} can be
deformed into the contour in Fig.~\ref{ch5:fig:1}. To estimate the discrepancy between the
time-discretized kernel and the continuous time one, one can thus compare the Green's function
along such contour. Again, the only arc that requires detailed attention is the arc $BCD$, as
the integral over rest of the contour of integration can be bounded from above as in the
previous section.

Let $h=h_m$ and let us introduce the two functions
\begin{equation}\label{ch5:eq:5.5}
\varphi_0(t,x,\tau)\equiv2\LatL(x,x+\tau h)e^{t\LatL(x,x)}\mathbf{1}(t\ge0),
\end{equation}
\begin{equation}\label{ch5:eq:5.6}
\varphi_{\delta t}(j,x,\tau)\equiv2\LatL(x,x+\tau h)\bigl(1+\delta t\LatL(x,x)\bigr)^{j-1}.
\end{equation}
and the corresponding Fourier transforms
\begin{equation}\label{ch5:eq:5.7}
\hat\varphi_0(\omega,x,\tau)=\int_0^{\infty}\varphi_0(t,x,\tau)e^{-i\omega t}\rev{\dd t}
=\left(\frac{v(x)}{h^2}+\tau\frac{\mu(x)}{h}\right)\left(i\omega+\frac{v(x)}{h^2}\right)^{-1}
\end{equation}
\begin{equation}\label{ch5:eq:5.8}
\hat\varphi_{\delta t}(\omega,x,\tau)=\rev{\delta t}\sum_{j=1}^{\infty}
\varphi_{\delta t}(j,x,\tau)e^{-i\omega j\delta t}
=\left(\frac{v(x)}{h^2}+\tau\frac{\mu(x)}{h}\right)
\frac{\delta t}{e^{i\omega\delta t}-1+\delta t\frac{v(x)}{h^2}},
\end{equation}
We have that
\begin{equation}\label{ch5:eq:5.9}
\begin{aligned}
\hat\varphi_{\delta t}(\omega,x,\tau)
&=\left(\frac{v(x)}{h^2}+\tau\frac{\mu(x)}{h}\right)
\left(i\omega+\frac{v(x)}{h^2}-\frac{\omega^2}{2}\delta t+O(\delta t^2)\right)^{-1}\\
&=\hat\varphi_0(\omega,x,\tau)+\frac{\omega^2}{2v(x)}h^2\delta t+O(h^2\delta t^2)
=\hat\varphi_0(\omega,x,\tau)+O(h^4),
\end{aligned}
\end{equation}
where the last step uses the fact that $\delta t=O(h^2)$. Let us also introduce the functions
\begin{equation}\label{ch5:eq:5.10}
\psi_0(t,x,\tau)\equiv e^{t\LatL(x,x)}\mathbf{1}(t\ge0),\qquad
\psi_{\delta t}(j,x,\tau)\equiv\rev{\bigl(1+\delta t\LatL(x,x)\bigr)^{j-1}}.
\end{equation}
and the corresponding Fourier transforms
\begin{equation}\label{ch5:eq:5.11}
\hat\psi_0(\omega,x,\tau)=\left(\frac{v(x)}{h^2}+i\omega\right)^{-1},\qquad
\hat\psi_{\delta t}(\omega,x,\tau)=\rev{\delta t}\left(e^{i\omega\delta t}-1+\delta t\frac{v(x)}{h^2}\right)^{-1}.
\end{equation}
Again we find that
\begin{equation}\label{ch5:eq:5.12}
\hat\psi_0(\omega,x,\tau)=\hat\psi_{\delta t}(\omega,x,\tau)+O(h^4).
\end{equation}
If $\gamma$ is a symbolic sequence, then let us set
\begin{equation}\label{ch5:eq:5.13}
\hat W_m(\gamma,q;\omega)=\hat\psi_0(\omega,\gamma_q)\prod_{j=0}^{q-1}
\hat\varphi_0\bigl(\omega;\gamma_j,\operatorname{sgn}(\gamma_{j+1}-\gamma_j)\bigr)
\end{equation}
\begin{equation}\label{ch5:eq:5.14}
\hat W^{\delta t}_m(\gamma,q;\omega)=\hat\psi_{\delta t}(\omega,\gamma_q)\prod_{j=0}^{q-1}
\hat\varphi_{\delta t}\bigl(\omega;\gamma_j,\operatorname{sgn}(\gamma_{j+1}-\gamma_j)\bigr).
\end{equation}
We have that
\begin{equation}\label{ch5:eq:5.15}
G^{\delta t}_m(x,y;\omega)-G_m(x,y;\omega)=\frac{1}{h}\sum_{q=1}^{\infty}2^{-q}\!\!\!
\sum_{\substack{\gamma\in\Gm:\ \gamma_0=x,\ \gamma_q=y\\ |\gamma_j-\gamma_{j-1}|=1\ \forall j\ge1}}
\!\!\!\Bigl(\hat W^{\delta t}_m(\gamma,q;\omega)-\hat W_m(\gamma,q;\omega)\Bigr).
\end{equation}
The integration over the contour in Fig.~\ref{ch5:fig:1} can again be split into an integration over
the contour $\Cm_-$ and an integration over the contour $\Cm_+$. The integral over $\Cm_+$ can
be bounded from above thanks to Lemma 1. Furthermore, we have that
\begin{equation}\label{ch5:eq:5.16}
\begin{aligned}
&\left|\int_{\Cm_-}\Bigl(G^{\delta t}_m(x,y;\omega)-G_m(x,y;\omega)\Bigr)e^{i\omega t}
\frac{\dd\omega}{2\pi}\right|\\
&\quad\le c\,h^{-1}\sqrt{q_{\max}}\!
\max_{\substack{\gamma\in\Gm:\ \gamma_0=x,\ \gamma_q=y\\ |\gamma_j-\gamma_{j-1}|=1\ \forall j\ge1}}
\left|\int_{\Cm_-}\Bigl(\hat W^{\delta t}_m(\gamma,q;\omega)-\hat W_m(\gamma,q;\omega)\Bigr)
e^{i\omega t}\frac{\dd\omega}{2\pi}\right|\le c\,h^2
\end{aligned}
\end{equation}
To bound the time derivative, we have to consider
\begin{equation}\label{ch5:eq:5.17}
\int_{\Cm_-}\left(\frac{e^{i\omega\delta t}-1}{\delta t}G_m(x,y;\omega)
-i\omega G_m(x,y;\omega)\right)e^{i\omega t}\frac{\dd\omega}{2\pi}
\end{equation}
But, since $\delta t=O(h^2)$, also this difference is $O(h^2)$.

\section{Machine-Checked Verification}\label{ch5:sec:lean}

The algebraic content of Sections~3 and~4, i.e.\ the identities involving the lattice
generator on three consecutive sites, has been formalised and machine-checked in Lean~4
(toolchain \texttt{v4.36.0-rc2}) together with Mathlib. The development accompanies this
revision: it is a small Lean project with the modules \texttt{Cert.Basic}, \texttt{Cert.Phi1}
and \texttt{Cert.Antiderivative} (and two further modules, \texttt{Cert.Abelian} and \texttt{Cert.Farkas}, which cover Chapter~6
and Chapters~1 to~3 respectively), it builds with \texttt{lake build}, and
it contains no \texttt{sorry}; every theorem listed below depends only on the three standard
axioms of Lean's classical logic (\texttt{propext}, \texttt{Classical.choice},
\texttt{Quot.sound}). The purpose of this section is to record what is now certified and what
is not, so that the reader can see exactly which statements are backed by a proof term.

\subsection*{What Has Been Certified}

\begin{center}
\small
\begin{tabular}{@{}lll@{}}
\hline
Result in the chapter & Lean theorem & Status\\
\hline
\eqref{ch5:eq:4.1}, \eqref{ch5:eq:4.2} & \texttt{discrete\_taylor\_plus/minus} & exact, all $h\ne0$\\
\eqref{ch5:eq:3.30} (the block $\LbarRaw(x;h)$) & \texttt{Lbar} & definition\\
\eqref{ch5:eq:3.31}--\eqref{ch5:eq:3.35} & \texttt{Lbar\_taylor\_decomposition} & exact, zero remainder\\
\eqref{ch5:eq:4.5}--\eqref{ch5:eq:4.9} & \texttt{Lbar\_discrete\_decomposition} & exact, no remainder\\
\eqref{ch5:eq:3.44} (symbol $\hat\varphi_0$) & \texttt{phi0\_expansion} & exact through $h^4$\\
\eqref{ch5:eq:3.45}, \eqref{ch5:eq:3.46} (symbol $\hat\varphi_1$) & \texttt{phi1\_instance} & exact instance\\
\eqref{ch5:eq:3.47} (primitives) & \texttt{antideriv\_inv}, \texttt{antideriv\_log} & symbolic\\
\hline
\end{tabular}
\end{center}

Three points are worth emphasising.

\begin{enumerate}
\item \emph{The decompositions are stronger than stated.} In \eqref{ch5:eq:3.31}--\eqref{ch5:eq:3.35},
if $v(x\pm h)$ is replaced by its Taylor polynomial through $h^3$ and $\mu(x\pm h)$ by its
Taylor polynomial through $h^2$, the identity
$\LbarRaw=\LbarRaw_0/h^2+\LbarRaw_1/h+\LbarRaw_2+h\LbarRaw_3$ holds with \emph{zero} remainder, so the
$O(h^2)$ in \eqref{ch5:eq:3.31} is entirely due to Taylor's theorem. Likewise, in
\eqref{ch5:eq:4.5}--\eqref{ch5:eq:4.9} the identity holds with no remainder at all once the discrete
operators $\nabm$ and $\Dm$ are used, which is precisely the claim made after \eqref{ch5:eq:4.5}.

\item \emph{The expansion of $\hat\varphi_0$ is exact up to a single explicit term.} Writing
$w=i\omega$ (so that $\omega^2=-w^2$), the truncation of \eqref{ch5:eq:3.44} through $h^4$ is
proved to differ from $\hat\varphi_0$ by $h^5\cdot32w^2(\mu\tau-2wh)/v^2$ and by nothing else.

\item \emph{The $h^3$ coefficient of $\hat\varphi_1$ is certified with the expression
\eqref{ch5:eq:3.46} for $r(x)$.} At the instance $v=2$, $v'=1$, $v''=3$, $v'''=5$, $\mu=4$,
$\mu'=5$, $\mu''=6$, $w=i\omega=7$, $h=1/10$, the truncation of \eqref{ch5:eq:3.45} reproduces the
symbol through $h^3$, the remainder being an exact rational multiple of $h^4$. This is an
instance-level check of an identity which is polynomial in the Taylor data; its general
symbolic form is the first item of the programme below.
\end{enumerate}

\subsection*{What Has Not Been Certified, and What It Would Take}

The analytic scaffolding of the chapter is not covered by the development, and is currently
beyond what Mathlib can support: the contour deformation and uniform resolvent bounds of
Lemma~1, the resummation of decorated path weights in Lemma~2, the renormalisation-group
comparison of Lemma~3, Theorems~1--3 as uniform statements over families of coefficients, and
the sharpness conjecture of the abstract. There is at present no theory of
H\"older-coefficient parabolic equations in Mathlib, no contour integration for
operator-valued resolvents, and no spectral theory of non-self-adjoint lattice generators,
all of which the argument uses.

A plausible staging, with rough estimates for a dedicated formalisation effort, is as follows.
(i)~The general symbolic form of the $h^3$ coefficient in \eqref{ch5:eq:3.45}, rather than a
single instance: days. The obstacle is only the size of the polynomial identity, whose
straightforward elaboration exceeds the default limits of the proof assistant; it can be met
by splitting the computation degree by degree. (ii)~The discrete part of Lemma~2, which is a
$q$-fold convolution estimate and needs no analysis beyond what Mathlib already provides: one
to two weeks. (iii)~Lemma~3 together with the lattice-to-lattice Cauchy estimate underlying
Theorem~1, stated for the finite matrices alone and thus free of any PDE theory: one to two
person-months, the effort being concentrated in the complex-analytic estimates of
Section~3. (iv)~The continuum theory, i.e.\ existence and uniqueness of $u_0$ and the
sharpness conjecture: a research programme rather than a formalisation task.

We regard the present state as the natural stopping point of a first pass: everything that had
to be carried out symbolically is now checked by machine, while the parts of the argument that
are conceptual rather than computational remain, for the moment, on paper.

\section{Conclusions}

We obtained bounds on convergence rates for explicit discretization schemes to the kernel of
one-dimensional diffusion equations with continuous coefficients. We consider both semidiscrete
triangulations with continuous time and explicit Euler schemes with time step small enough for
the method to be stable. The proof is constructive and based on a new technique of path
conditioning for Markov chains and a renormalization group argument. Convergence rates depend
on the degree of smoothness and H\"older differentiability of the coefficients. The method is
of more general applicability and will be extended in future work. With Chapters~1 to~3,
which establishes the finite core of the theory, and Chapter~6, which carries
the analysis to path dependent functionals, this chapter completes the two halves of the programme
described in the introduction: the theorem on which pricing rests, and the estimates by which
the passage to the continuum is controlled.
% ---- end generated/ch5-body.tex ----
% ---- end chapters/ch5.tex ----
% ---- from chapters/ch6.tex ----
% Authored: the chapter title and opening belong to the book. The body is
% generated by build.py from the parts and is not edited here.
\chapter{Stochastic Integrals and Abelian Processes}
\label{chap:ch6}

A stochastic integral is an aggregate accumulated along a trajectory: increments of the path weighted
by a function of the site, summed as the path is traversed. In the pricing reading it is what a path
dependent payoff is written on --- realised variance, a running maximum, the integral which appears
in the hedge of a claim --- and in the learning reading it is the statistic which a model computes
alongside the state as it moves through a sequence, and whose law is what a prediction about the
sequence is a statement about. The two are one object, and the chapter treats it once.

The chapter is the estimate for the object which Chapter~\ref{chap:ch4} made computable: the joint law
of a process and of an Abelian functional, in one space dimension, computed on a triangulation and
compared with its limit at a rate which is explicit in the smoothness of the coefficients of the
process and of the functional. Chapter~\ref{chap:ch8} carries the same estimate to higher dimension,
and with it the question of what the passage costs.

% ---- from generated/ch6-body.tex ----
%% Generated by build.py from the parts - do not edit.
\section{Introduction}

Consider a diffusion defined by the Markov generator
\begin{equation}\label{ch6:eq:1.1}
L_{0x}=\frac{1}{2}\sigma(x)^2\frac{\partial^2}{\partial x^2}+\mu(x)\frac{\partial}{\partial x}
\end{equation}
on the bounded interval $X=[-L_x,L_x]\subset\mathbb{R}$ where $0<L_x<\infty$. For simplicity,
we assume periodic boundary conditions and identify the two boundary points $\pm L_x$ with
each other. The coefficients $\sigma(x)^2$ and $\mu(x)$ are assumed to be at least uniformly
continuous and our results will depend on the degree of smoothness. We also assume that
$\sigma(x)>\Sigma_0$ for some constant $\Sigma_0>0$.

In this chapter, we consider path dependent processes such as stochastic integrals of the form
\begin{equation}\label{ch6:eq:1.2}
y_t=\int_0^ta(x_s)\dd x_s+b(x_s)\dd s
\end{equation}
where $a(x)$ and $b(x)$ are functions on $X=[-L_x,L_x]$. The case $a(x)=0$ was considered in
(Girsanov 1960), (Cameron and Martin 1949), (Feynman 1948) and (Kac 1948). The case with
$a(x)\neq0$ was first tackled in (Ito 1949). The proofs were streamlined by using martingale
theory in (Doob 1953). Other references for stochastic integrals are (Kunita and Watanabe 1967),
(McKean 1969) and (Meyer 1974). Most of these results assume smooth coefficients. The more
general case of uniformly continuous coefficients was considered in (Stroock and Varadhan 1969).

We are interested in convergence properties of triangulation schemes for the joint distribution
between the underlying process $x_t$ and the stochastic integral $y_t$. If the diffusion
process is approximated by a sequence of Markov chains then existence and uniqueness of
stochastic integrals is not problematic. The problem is to determine conditions under which
joint distributions between stochastic integrals and the underlying process converge in some
meaningful norm in the limit as the triangulation becomes finer and finer, approaching the
continuum limit. We establish convergence assuming uniform continuity in the coefficients of
the diffusion and the stochastic integral and find convergence rates. Convergence takes place
in the graph-norm of the Fourier transformed Markov generator. In the continuum limit, the
Fourier transform of the joint distribution is entire analytic in the conjugate variable of the
stochastic integral and is in the domain of all the powers of the Fourier transformed Markov
generator. Consideration of the Fourier transform is essential in the proof as the joint
distribution itself can be singular even in trivial special cases such as the one where
$a(x)=b(x)=0$. The Fourier analysis also leads to a new derivation of the celebrated formulas
for characteristic functions of stochastic integrals in (Girsanov 1960), (Cameron and Martin
1949), (Feynman 1948) and (Ito 1949). Furthermore, analyticity implies convergence of moment
formulas.

These results are classic and well established in the continuum limit. One of the new facets of
our derivation is that it is entirely constructive and thus instructive from a computational
viewpoint. We make no use of compactness arguments and non-constructive measure theory methods.
The derivation is based on pathwise analysis and a renormalization group transformation. The
construction applies to discretized operators and does not rely on the existence of a
probabilistic interpretation. The proof of convergence is actually carried out for analytic
extensions of the Fourier transform of the joint distribution, so we consider complex valued
weights for the paths. The argument would also extend to different type of diffusions with
complex coefficients as they occur for instance in quantum mechanics. Furthermore results are
extended to discrete time approximations with fully explicit Euler schemes, in which case we
also derive convergence bounds in the graph-norm for Fourier transformed joint kernels.

In (Albanese 2007c) we discuss the case of the kernel of diffusion processes, a more elementary
situation. This chapter is strictly more general and it is self-contained with no dependencies on
our previous work as all derivations need to be adapted in detail. The key technique of
pathwise analysis and renormalization group resummations over decorated paths are however
quite similar.

Our operator algebraic approach to the problem also points out that stochastic integrals are
not the most general class of path dependent processes to which these methods apply. This is
based on the recognition that Fourier transforms play a key role as they allow one to
block-diagonalize the joint kernel. But block-diagonalization can be achieved under more
general conditions. We consider a generic situation of a diffusion and an adapted process. We
associate an operator algebra to a finite triangulation of the joint process. If this operator
algebra is Abelian, i.e.\ commutative, and if it satisfies a certain non-resonance condition,
then the joint Markov generator can also be block-diagonalized and methods can be extended.

Non-resonant Abelian processes are a broad class of hypo-elliptic multidimensional diffusion
processes which are computable by robust methodologies based on block-diagonalizations and
fully explicit Euler schemes. Several instances of these processes found engineering
applications already, see (Albanese and Vidler 2007), (Albanese and Trovato 2005), (Albanese
2007a), (Albanese et al.\ 2006), (Albanese and Osseiran 2007). As an example, we discuss here
the case of the sup processes of a one-dimensional diffusion, i.e.\ of
\begin{equation}\label{ch6:eq:1.3}
y_t=\sup_{s\in[0,t]}x_s.
\end{equation}
We also discuss discrete time summations of the form
\begin{equation}\label{ch6:eq:1.4}
y_{T_n}=\sum_{i=1}^{n}\varphi(x_{T_{i-1}},x_{T_i})
\end{equation}
where $T_i=i\Delta T$.

The kernel convergence in graph norm we establish in this chapter has noteworthy implications
from the computer science viewpoint. As we discuss in (Albanese 2007b), they imply that
explicit schemes give a robust valuation methodology for joint kernels of processes specified
semi-parametrically. The Courant--Friedrichs--Lewy stability condition (Courant et al.\ 1928)
for explicit methods is not much of an impediment in the presence of sufficient system memory
as linear fast exponentiation can be used to accelerate the scheme. The ability of evaluating
kernels robustly is the main reason why analytically solvable models are interesting and
broadly used. But not all models solvable in closed form are practically computable in the
sense that they allow for robust kernel valuations within the limits of double precision
floating point arithmetics. Higher order hypergeometric functions are often difficult to handle
and may require multiple precision. Linear fast exponentiation is a technique based on
accelerated fully explicit schemes that empirically is observed to work well even in single
precision and for very general model specifications. Empirically, we find that it performs
better in single precision on some solvable models than the use of closed form solutions in
double precision. This chapter was prompted by the observation of this phenomenon and the desire
to understand and explain the underlying smoothing mechanisms which is of technological
significance.

\subsection*{Place in the Book}

A stochastic process, in the framework of this book, is a consistent family of Markov chains
on a sequence of lattices, and a statement about a continuous process is a statement about the
limit of statements about finite matrices. Chapters~1 to~3 establish the fundamental theorem at a fixed resolution, where the existence of a pricing measure reduces to
the orthant alternative of finite dimensional convex geometry. Chapter~5 establishes uniform convergence rates for the kernels of the diffusions
themselves; the present chapter extends that analysis to the joint kernel of a diffusion and an
Abelian process. This is what a constructive definition of the stochastic integral requires:
defined at fixed resolution the integral is an elementary object, a finite sum along a
piecewise constant path, and what is needed is a norm in which it converges together with the
underlying process, uniformly in the degree of smoothness of the coefficients.

The order of the argument is worth spelling out, since the two halves of it sit in different
parts. The functionals whose joint kernel is estimated here are the Abelian ones of Chapter~4, and
the identity which makes their joint law computable at all is the It\^o formula at fixed
resolution, proved there in Section~2 as a rearrangement of a finite sum, with no change of
measure and no limit. What the present chapter adds is the control of the limit: the estimates by
which the rearrangement at resolution $m$ is compared with the rearrangement at resolution
$m+1$, uniformly in the coefficients.

The motivation is the pricing problem of Chapters~1 to~3. A path dependent claim is valued by
integrating its payoff against the joint law of the process and of the functional on which the
payoff depends, so the joint kernel is precisely the object that a constructive pricing theory
needs; and the functionals which arise in practice, the stochastic integral of the asset along a
trading strategy first of all, are the Abelian processes studied here. As in Chapter~5, the
estimates are required in this generality because the pricing measures of Chapters~1 to~3 are obtained by
reweighting the transition rates, which does nothing to make the coefficients smoother.

The reason for wanting a formulation of this kind is set out in Chapters~1 to~3: a foundation whose
statements concern finite objects, and whose limits are estimates rather than existence
theorems, is one that can be machine checked, and the foundations of the subject stand to
benefit from the AI-assisted methods now available.

The chapter is organized as follows: In Section~2 we introduce notations and state our results
regarding stochastic integrals. Proofs in the case of stochastic integrals are in Section~3,
where we consider the case of continuous time and in Section~4 where we discuss convergence for
fully explicit Euler schemes. The Ito representation is in Section~5. Extensions to the sup and
other Abelian processes are given in Section~6. Section~7 develops a gauge-free route to the
cancellation over the sum of paths, which is what a higher dimensional extension of the argument
requires, and conclusions end the chapter.

\section{Stochastic Integrals}

Consider a family of increasingly fine triangulations schemes whereby the space dimension is
discretized in multiples of an elementary space step $h_{xm}=L_x2^{-m}$, $m\in\mathbb{N}$ and
we are interested in the limit as $m\to\infty$. Let $X_m=h_{xm}\mathbb{Z}\cap X$ and consider
the sequence of operators
\begin{equation}\label{ch6:eq:2.1}
\LatL_x=\frac{\sigma(x)^2}{2}\Dm_x+\mu(x)\nabm_x.
\end{equation}
defined on the $2^{m+1}$-dimensional space of all periodic functions $f_m:X_m\to\mathbb{R}$,
where
\begin{equation}\label{ch6:eq:2.2}
\nabm_xf(x)=\frac{f(x+h_{xm})-f(x-h_{xm})}{2h_{xm}}
\end{equation}
and
\begin{equation}\label{ch6:eq:2.3}
\Dm_xf(x)=\frac{f(x+h_{xm})+f(x-h_{xm})-2f(x)}{h_{xm}^2}.
\end{equation}
These definitions also apply to the boundary points by periodicity. We assume that $m\ge m_0$
where $m_0$ is the least integer such that
\begin{equation}\label{ch6:eq:2.4}
\frac{\sigma^2(x)}{2h_{xm}^2}>\frac{|\mu(x)|}{2h_{xm}}
\end{equation}
for all $m\ge m_0$ and all $x\in X_m$.

Consider the kernel $u_m(x,x';t)$ of equation \eqref{ch6:eq:2.6}, i.e.\ the solution of the
(forward) equation
\begin{equation}\label{ch6:eq:2.5}
\frac{\partial}{\partial t}u_m(x,x';t)=\LatL{}^{*}_{x'}u_m(x,x';t)
\end{equation}
where the operator $\LatL{}^{*}_{x'}$ acts on the $x'$ coordinate and the following initial time
condition is satisfied:
\begin{equation}\label{ch6:eq:2.6}
u_m(x,x';0)=h_{xm}^{-1}\delta_{X_m}(x-x').
\end{equation}
where
\begin{equation}\label{ch6:eq:2.7}
\delta_{X_m}(x-x')=\begin{cases}1&\text{if }x=x'\bmod 2L_x\\ 0&\text{otherwise.}\end{cases}
\end{equation}
If $f(x,x'):X_m\times X_m\to\mathbb{R}$, consider the uniform norm
\begin{equation}\label{ch6:eq:2.8}
\|f_m\|_{m,\infty}=\sup_{x,x'\in X_m}|f(x,x')|
\end{equation}
and the graph norm
\begin{equation}\label{ch6:eq:2.9}
\|f\|_{m,L}=\|f\|_{\infty}+\|\LatL_xf\|_{\infty}+\|\LatL{}^{*}_{x'}f\|_{\infty}.
\end{equation}
In (Albanese 2007c) we show that if the coefficients $\sigma(x)^2$ and $\mu(x)$ are uniformly
continuous, then the sequence of kernels $u_m$ is Cauchy with respect to the graph norm. If in
addition coefficients are H\"older continuous and $\sigma^2\in C^{k_\sigma,\alpha_\sigma}$,
$\mu\in C^{k_\mu,\alpha_\mu}$ and
\begin{equation}\label{ch6:eq:2.10}
\gamma\equiv\min\{2,k_\sigma+\alpha_\sigma,k_\mu+\alpha_\mu\}>0
\end{equation}
then
\begin{equation}\label{ch6:eq:2.11}
\|u_m(t)-u_{m_0}(t)\|_{m,L}\le c\,h_{xm}^{\gamma}
\end{equation}
for all $m_0>m\ge m_0$. A similar bound also holds for the kernels obtained with a fully
explicit Euler scheme
\begin{equation}\label{ch6:eq:2.12}
u^{\delta t}_m(x,x';t)=h_{xm}^{-1}\bigl(1+\delta t_m\LatL\bigr)^{\left[\frac{t}{\delta t_m}\right]}(x,x';t).
\end{equation}
where $\delta t_m$ is so small that
\begin{equation}\label{ch6:eq:2.13}
\min_{x\in X_m}1+\delta t_m\LatL(x,x)>0.
\end{equation}
In fact, we have that
\begin{equation}\label{ch6:eq:2.14}
\|u_m(t)-u^{\delta t}_m(t)\|_{m,L}\le c\,h_{xm}^2
\end{equation}
for some constant $c>0$.

Let $h_{yn}$ be a monotonously descreasing sequence such as $h_{yn}\to0$ as $n\to\infty$. Also
let $L_{yn}$ be a monotonously increasing sequence such that $L_{yn}\to\infty$ as $n\to\infty$.
Let
\begin{equation}\label{ch6:eq:2.15}
Y_n=h_{yn}\mathbb{Z}\cap[-L_{yn},L_{yn}].
\end{equation}
We assume again periodic boundary conditions in the $y$ direction and identify the two extreme
points of $Y_n$. Let's consider two-dimensional processes described by a sequence of Markov
generators of the form
\begin{equation}\label{ch6:eq:2.16}
L^{m,n}(x,y;x',y')=\LatL(x,x')+Q^{m,n}(x,y;x',y')
\end{equation}
where $x,x'\in X_n$, $y,y'\in Y_n$. We assume that
\begin{equation}\label{ch6:eq:2.17}
\sum_{y'\in Y_n}Q^{m,n}(x,y;x',y')=0
\end{equation}
for all values of $x,y,x'$, so that the marginals with respect to the first process are the same
as under the dynamics given by $\LatL$, i.e.
\begin{equation}\label{ch6:eq:2.18}
\sum_{y'\in Y_n}\exp\bigl(tL^{m,n}\bigr)(x,y;x',y')=\exp\bigl(t\LatL\bigr)(x,x')
\end{equation}
for all triples $x,y,x'$. We are interested in the joint kernel
\begin{equation}\label{ch6:eq:2.19}
u_{mn}(x,y;x',y';t)\equiv\frac{1}{h_{xm}h_{yn}}\exp\bigl(tL^{m,n}\bigr)(x,y;x',y')
\end{equation}
and its convergence properties in the limit as $m,n\to\infty$.

The case of general multi-dimensional diffusion processes is discussed in (Albanese and Jones
2007a). Results in that paper however do not apply to the class of processes in which we are
interested here. In this article, we consider the case of path dependent processes which are
not driven by their own diffusive dynamics and whose generator is not strongly elliptic, but
rather hypoelliptic. Because of this reason, the joint kernel may be singular and a separate
treatment is required.

\begin{definition}\label{ch6:def:1}
Consider a sequence of Markov generators of the form in \eqref{ch6:eq:2.16} describing a bivariate
process. One says that the second process is \emph{Abelian} with respect to the first if the
matrices $Q^{m,n}_{x,x'}$ of elements
\begin{equation}\label{ch6:eq:2.20}
Q^{m,n}_{x,x'}(y,y')=Q^{m,n}(x,y;x',y')
\end{equation}
are mutually commuting, i.e.\ if
\begin{equation}\label{ch6:eq:2.21}
\bigl[Q^{m,n}_{x,x'},Q^{m,n}_{x'',x'''}\bigr]
=Q^{m,n}_{x,x'}Q^{m,n}_{x'',x'''}-Q^{m,n}_{x'',x'''}Q^{m,n}_{x,x'}=0
\end{equation}
for all $x,x',x'',x'''\in X_m$.
\end{definition}

\begin{definition}\label{ch6:def:2}
An Abelian process is called \emph{non-resonant} if there is a sequence of non-singular
transformations $V_{mn}(p,y)$ which diagonalizes simultaneously the kernels $Q^{m,n}_{x,x'}$ for
all pairs $x,x'\in X_m$, i.e.\ if
\begin{equation}\label{ch6:eq:2.22}
\sum_{y,y'\in Y_n}V_{mn}(p,y)Q^{m,n}(x,y;x',y')V^{-1}_{mn}(p',y')=\Lambda_{mn}(x,x')\delta_{p,p'}.
\end{equation}
for all $x,x'\in X_m$. The index $p$ ranges on a set denoted with $\Yhn$ which has the same
cardinality as $Y_n$ and is called the inverse lattice of $Y_n$.
\end{definition}

In Sections 3, 4 and 5, we consider in detail the prototypical example of Abelian process,
namely stochastic integrals. In this case, a triangulation for the joint generator is given by
\begin{equation}\label{ch6:eq:2.23}
\begin{aligned}
L^{m,n}(x,y;x',y')={}&\left(\frac{\sigma(x)^2}{2h_{xm}^2}+\frac{\mu(x)}{2h_{xm}}\right)
\delta_{X_m}(x'-x-h_{xm})\,\delta_{Y_n}\!\left(y'-y-\left[\frac{a(x)h_{xm}}{h_{yn}}\right]h_{yn}\right)\\
&+\left(\frac{\sigma(x)^2}{2h_{xm}^2}-\frac{\mu(x)}{2h_{xm}}\right)
\delta_{X_m}(x'-x+h_{xm})\,\delta_{Y_n}\!\left(y'-y+\left[\frac{a(x)h_{xm}}{h_{yn}}\right]h_{yn}\right)\\
&-\frac{\sigma(x)^2}{h_{xm}^2}\delta_{X_m}(x'-x)\delta_{Y_n}(y'-y)
+b(x)\delta_{xx'}\nabla^n_y(y,y')
\end{aligned}
\end{equation}
where
\begin{equation}\label{ch6:eq:2.24}
\nabla^n_yg(y)=\frac{g(y+h_{yn})-g(y-h_{yn})}{2h_{yn}}.
\end{equation}
and
\begin{equation}\label{ch6:eq:2.25}
\delta_{Y_n}(y-y')=\begin{cases}1&\text{if }y=y'\bmod 2L_y\\ 0&\text{otherwise.}\end{cases}
\end{equation}
Diagonalisation is simply achieved by Fourier transforms of kernel
\begin{equation}\label{ch6:eq:2.26}
F_n(x,p;x',y)=e^{-iyp}\delta_{xx'}
\end{equation}
where
\begin{equation}\label{ch6:eq:2.27}
p\in\Yhn\equiv\widehat{h}_{yn}\mathbb{Z}\cap\left[-\frac{\pi}{h_y},\frac{\pi}{h_y}\right],
\qquad\text{and}\qquad \widehat{h}_y=\frac{\pi}{L_{yn}}.
\end{equation}
Partial Fourier transforms in the $p$ variable block-diagonalize the joint generator, reducing
it to the form
\begin{equation}\label{ch6:eq:2.28}
\bigl(F_nL^{m,n}F_n^{-1}\bigr)(x,p;x',p')=\widehat{L}^{m,n}(x,x';p)\delta_{pp'}
\end{equation}
where $\widehat{L}^{m,n}(x,x';p)$ is a one-parameter family of matrices indexed by
$x,x'\in X_m$. An explicit expression for the joint generator is given in equation \eqref{ch6:eq:5.6}
below. The kernel itself can also be expressed by means of a Fourier transform as follows:
\begin{equation}\label{ch6:eq:2.29}
u_{mn}(x,p;x',p';t)\equiv\bigl(F_n^{-1}\widehat{u}_{mn}(t)F_n\bigr)(x,p;x',p')
=\frac{1}{h_{xm}h_{yn}}\delta_{\Yhn}(p-p')\exp\bigl(t\widehat{L}^{m,n}(p)\bigr)(x,x'),
\end{equation}
We are interested in the convergence properties of the sequence
\begin{equation}\label{ch6:eq:2.30}
\widehat{u}_m(x,x';p,t)=\lim_{n\to\infty}\frac{1}{h_{yn}}\widehat{u}_{mn}(x,p;x',p';t)
\end{equation}
in the limit as $m\to\infty$. To establish convergence, it is necessary to work in the Fourier
transformed kernel. In fact, while the Fourier transformed kernel converges in a very strong
graph-norm and the limit is entire analytic in $p$ and has the same smoothness properties in
$x,x'$ as the diffusion kernel, the joint kernel expressed with respect to space coordinates
can be quite singular even in simple cases. If for instance the functions $a(x)$ and $b(x)$ in
\eqref{ch6:eq:1.2} are zero, then the joint distribution $U(x,y;x',y';t)$ for the pair $(x_t,y_t)$
is concentrated on the line $y_t=0$. If a partial Fourier transform with respect to the $y$
variable is taken, the result is constant as a function of the dual variable $p$. This is a
simple example of entire analytic dependency in $p$ which however translates into a rough delta
type singularity in the coordinate representation.

Let $g(x,x';z)$ be a complex valued function defined for $x,x'\in X_m\times X_m$ and for
$z\in\mathbb{C}$ such that $|z|<K$ where $K>0$. Consider the uniform norm
\begin{equation}\label{ch6:eq:2.31}
\|g\|_{m,K,\infty}=\sup_{\substack{x,x'\in X_m\\ |z|<K}}|g(x,x';z)|.
\end{equation}
The graph norm of order $(m,K,L)$ is defined as follows:
\begin{equation}\label{ch6:eq:2.32}
\|g\|_{m,K,L}=\|g_m\|_{m,K,\infty}+\|\LatL_xg\|_{m,K,\infty}+\|\LatL{}^{*}_{x'}g\|_{m,K,\infty}.
\end{equation}

\begin{theorem}\label{ch6:thm:1}
Let's assume the function $b(x)$ is integrable and let
\begin{equation}\label{ch6:eq:2.33}
B(x)=\int_0^xb(x')\dd x'
\end{equation}
be a primitive. If the coefficients $\sigma(x)^2$, $\mu(x)$ and the functions $a(x)$ and $B(x)$
are uniformly continuous, then the sequence of Fourier transformed kernels $\widehat{u}_m(t)$
can be extended by analyticity in $p$ to an entire operator valued function and is Cauchy with
respect to the above graph-norm for all $K>0$. If in addition these coefficients are H\"older
continuous so that $\sigma^2\in C^{k_\sigma,\alpha_\sigma}$,
$\mu\in C^{k_\mu,\alpha_\mu}$, $a\in C^{k_a,\alpha_a}$, $B\in C^{k_B,\alpha_B}$ and
\begin{equation}\label{ch6:eq:2.34}
\gamma\equiv\min\{2,k_\sigma+\alpha_\sigma,k_\mu+\alpha_\mu,k_a+\alpha_a,k_B+\alpha_B\}>0
\end{equation}
then, for all $K>0$ there is a constant $c(K)$ such that
\begin{equation}\label{ch6:eq:2.35}
\|\widehat{u}_m(t)-\widehat{u}_{m_0}(t)\|_{m,K,L}\le c(K)h_{xm}^{\gamma}.
\end{equation}
for all $m_0>m\ge m_0$.
\end{theorem}

\begin{corollary}\label{ch6:cor:1}
The limit kernel
\begin{equation}\label{ch6:eq:2.36}
u(x,x';z)\equiv\lim_{m\to\infty}\frac{1}{h_{xm}}u_m(x,x';z)
\end{equation}
is an entire analytic function of $z$ which is in the domain of the operators $L_x$ and
$L^{*}_y$ and all of their powers. The Fourier transformed kernel
\begin{equation}\label{ch6:eq:2.37}
u(x,x';y)\equiv\int_{-\infty}^{\infty}e^{-ipy}u(x,x';p)\frac{\dd p}{2\pi}
\end{equation}
exists in the distribution sense for each fixed pair $x,x'\in[-L_x,L_x]$.
\end{corollary}

\begin{theorem}\label{ch6:thm:2}
Let us assume that the coefficients $\sigma(x)^2$, $\mu(x)$ and the functions $a(x)$ and $B(x)$
are uniformly continuous and fix a $K>0$. Consider the kernels obtained with a fully explicit
Euler scheme, i.e.
\begin{equation}\label{ch6:eq:2.38}
\widehat{u}^{\delta t}_m(x,x';z,t)=\frac{1}{h_x}\bigl(1+\delta t_m\widehat{L}^m(z)\bigr)^{\left[\frac{t}{\delta t^m}\right]}(x,x'),
\end{equation}
where $z\in\mathbb{C}:|z|<K$. Assume that $\delta t_m$ satisfies the Courant condition
\begin{equation}\label{ch6:eq:2.39}
\min_{x,p\in X_m\times\Yhn}<1+\delta t_m\widehat{L}^m(x,x;z)>0
\end{equation}
for all $|z|<K$. Then, there is a constant $c(K)$ such that
\begin{equation}\label{ch6:eq:2.40}
\|\widehat{u}_m(t)-\widehat{u}^{\delta t}_m(t)\|_{m,K,L}\le c(K)h_{xm}^2
\end{equation}
for all $m_0>m\ge m_0$.
\end{theorem}

This two Theorems are proved in the next two sections. An equivalent statement in a different
representation which leads to Ito's Lemma, the Cameron--Marin--Girsanov and the Feynman--Kac
formulas is given in Section~5.

\section{The Renormalization Group Argument}\label{ch6:sec:rg}

The proof of convergence in graph-norm is based on a path-integral representation of the
probability kernel. More precisely, let us define a symbolic path
$\gamma=\{\gamma_0,\gamma_1,\gamma_2,\dots\}$ as an infinite sequence of sites in $X_m$ such
that $\gamma_j=\gamma_{j-1}\pm1$ for all $j=1,\dots$. Let $\Gm$ be the set of all symbolic paths
in $X_m$. Then the propagator admits the following representation:
\begin{equation}\label{ch6:eq:3.1}
u_m(x,x';t)=\frac{1}{h_{xm}}\sum_{q=1}^{\infty}
\sum_{\substack{\gamma\in\Gm:\gamma_0=x,\gamma_q=x'}}
\int_0^t\!\dd s_1\int_0^t\!\dd s_2\cdots\int_{s_{q-1}}^t\!\dd s_q\,
\rho(s_1,s_2,\dots,s_{q-1};\gamma)
\end{equation}
where
\begin{equation}\label{ch6:eq:3.2}
\rho(s_1,s_2,\dots,s_{q-1};\gamma)
=e^{-s_1\LatL(\gamma_0,\gamma_0)}\prod_{j=1}^{q-1}\LatL(\gamma_{j-1},\gamma_j)
e^{-(s_{j+1}-s_j)\LatL(\gamma_j,\gamma_j)}\dd s_1\dots\dd s_q.
\end{equation}
and $s_q=t$.

Let $\gamma\in\Gm$ and consider a path $X_m(\cdot;t_1,\dots t_q;\gamma):\mathbb{R}_+\to X_m$,
left continuous and with right limits, taking the values $\{\gamma_0,\gamma_1,\gamma_2,\dots\}$
consecutively with jumps occurring at times $t_j$, $j=1,..q-1$, with $0\le t_1\le t_2\le\dots\le
t_{q-1}\le t$. Let $a(x)$ and $b(x)$ be two uniformly continuous functions in $X=[-L_x,L_x]$ and
consider the integral
\begin{equation}\label{ch6:eq:3.3}
\begin{aligned}
I(t_1,\dots,t_{q-1};\gamma)&=\int_0^t\left[a\bigl(X_m(s;t_1,\dots,t_q;\gamma)\bigr)
\frac{\dd X_m(s;t_1,\dots,t_q;\gamma)}{\dd s}\right.\\
&\qquad\qquad\left.+b\bigl(X_m(s;t_1,\dots,t_q;\gamma)\bigr)\right]\dd s\\
&=\sum_{j=0}^{q-1}a(\gamma_j)(\gamma_{j+1}-\gamma_j)+\sum_{j=0}^{q}b(\gamma_j)(t_{j+1}-t_j).
\end{aligned}
\end{equation}
We are interested in evaluating the joint density
\begin{equation}\label{ch6:eq:3.4}
\begin{aligned}
u_m(x,y;x',y';t)={}&\frac{1}{h_{xm}}\sum_{q=1}^{\infty}
\sum_{\substack{\gamma\in\Gm:\gamma_0=x,\gamma_q=x'}}
\int_0^t\!\dd s_1\int_0^t\!\dd s_2\cdots\int_{s_{q-1}}^t\!\dd s_q\\
&\times\rho(s_1,s_2,\dots,s_{q-1};\gamma)\,
\delta\bigl(I(s_1,\dots s_{q-1};\gamma)-(y'-y)\bigr)
\end{aligned}
\end{equation}
We do so by expanding its Fourier transform
\begin{equation}\label{ch6:eq:3.5}
\widehat{u}_m(x,x';p;t)=\frac{1}{h_{xm}}\sum_{q=1}^{\infty}
\sum_{\substack{\gamma\in\Gm:\gamma_0=x,\gamma_q=x'}}
\int_0^t\!\dd s_1\cdots\int_{s_{q-1}}^t\!\dd s_q\,
\rho(s_1,s_2,\dots,s_{q-1};\gamma)e^{ipI(s_1,\dots s_{q-1};\gamma)}
\end{equation}
\begin{equation}\label{ch6:eq:3.6}
\begin{aligned}
&=\lim_{n\to\infty}\sum_{y\in Y_n}\widetilde{L}^{mn}(x,0;x',y;t)e^{-ih_{yn}py}\\
&=\left(\frac{\sigma(x)^2}{2h_{xm}^2}+\frac{\mu(x)}{2h_{xm}}\right)e^{-ih_{xm}a(x)p}
\delta_{x',x+h_{xm}}
+\left(\frac{\sigma(x)^2}{2h_{xm}^2}-\frac{\mu(x)}{2h_{xm}}\right)e^{ih_{xm}a(x)p}
\delta_{x',x-h_{xm}}\\
&\qquad-\frac{\sigma(x)^2}{h_{xm}^2}\delta_{x'x}-ipb(x)\delta_{x'x}\\
&=\frac{\sigma(x,p,h_{xm})^2}{2}\Dm_x(x,x')+\mu(x,p,h_{xm})\nabm_x(x,x')+\kappa(x,p,h_{xm})\delta_{xx'}
\end{aligned}
\end{equation}
where
\begin{align}
\sigma(x,p,h_{xm})^2&=\sigma(x)^2\cos(h_{xm}a(x)p)
-i\,h_{xm}\sin(h_{xm}pa(x))\mu(x)\label{ch6:eq:3.7}\\
\mu(x,p,h_{xm})&=\mu(x)\cos(h_{xm}a(x)p)-i\,\frac{\sin(h_{xm}pa(x))}{h_{xm}}\sigma(x)^2,
\label{ch6:eq:3.8}\\
\kappa(x,p,h_{xm})&=-i\,\frac{\sin(h_{xm}pa(x))}{h_{xm}}\mu(x)
+\frac{\sigma(x)^2}{h_{xm}^2}\bigl(\cos(h_{xm}a(x)p)-1\bigr)-ipb(x).
\label{ch6:eq:3.9}
\end{align}
Notice that these three functions admit continuations as entire analytic functions in $p$. Let
us denote them with $\sigma(x,z,h_{xm})^2$, $\mu(x,z,h_{xm})$, $\kappa(x,z,h_{xm})$. Let
$\widehat{L}^m(x,x';z)$ be the analytic continuation of the operator above and consider the
Fourier transformed joint propagator defined as follows:
\begin{equation}\label{ch6:eq:3.10}
\widehat{u}_m(x,x';z;t)=\frac{1}{h_{xm}}\sum_{q=1}^{\infty}
\sum_{\substack{\gamma\in\Gm:\gamma_0=x,\gamma_q=x'}}
\int_0^t\!\dd s_1\cdots\int_{s_{q-1}}^t\!\dd s_q\,
\rho(s_1,s_2,\dots,s_{q-1};z;\gamma)
\end{equation}
where
\begin{equation}\label{ch6:eq:3.11}
\rho(s_1,s_2,\dots,s_{q-1};z;\gamma)=\frac{1}{h_{xm}}
e^{-s_1\widehat{L}^m(\gamma_0,\gamma_0;z)}\prod_{j=1}^{q-1}\widehat{L}^m(\gamma_{j-1},\gamma_j;z)
e^{-(s_{j+1}-s_j)\widehat{L}^m(\gamma_j,\gamma_j;z)}\dd s_1\dots\dd s_q.
\end{equation}
and $s_q=t$.

If $K>0$, let us introduce the constants
\begin{equation}\label{ch6:eq:3.12}
\Sigma_0=\inf_{x\in X_m}\sigma(x,0,h_{xm}),\qquad
M(K)=\sup_{\substack{x\in X_m\\ |z|<K}}|\mu(x,z,h_{xm})|
\end{equation}
and
\begin{equation}\label{ch6:eq:3.13}
\Sigma_1(K)=\sup_{\substack{x\in X_m\\ |z|<K}}
\sqrt{|\sigma(x,z,h_{xm})|^2+|\sigma(x,z,h_{xm})-\sigma(x,0,h_{xm})|^2+h_{xm}|\mu(x,z,h_{xm})|}.
\end{equation}
Due to our assumptions, $\Sigma_0>0$ and $\Sigma_1(K),M(K)<\infty$.

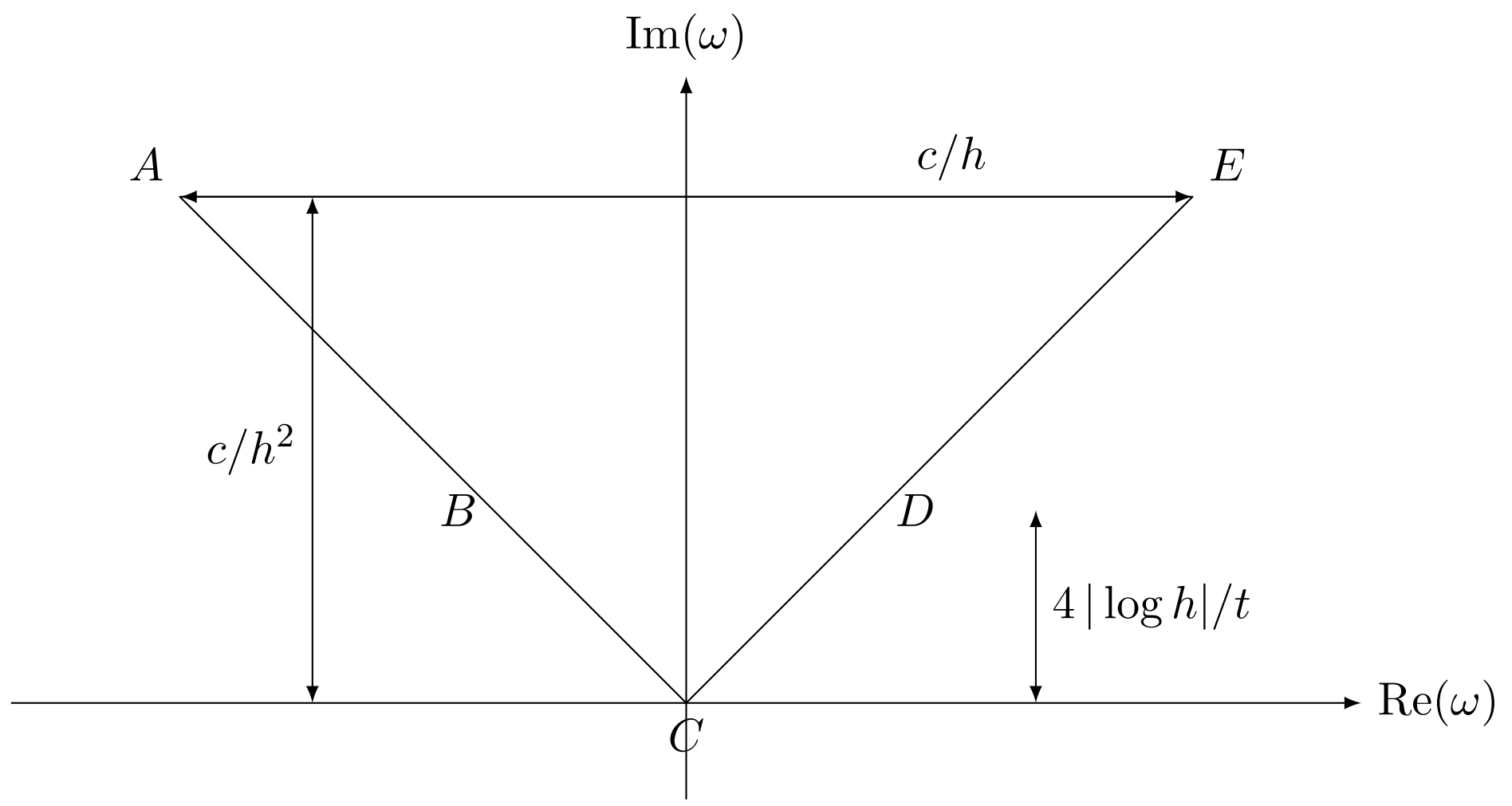
\begin{figure}[t]
\centering
\begin{tikzpicture}[scale=1.0,>=latex,line join=round]
  \draw[->] (-5.6,0) -- (5.6,0) node[right] {$\mathrm{Re}(\omega)$};
  \draw[->] (0,-0.8) -- (0,5.2) node[above] {$\mathrm{Im}(\omega)$};
  \coordinate (A) at (-4.2,4.2); \coordinate (E) at (4.2,4.2);
  \coordinate (B) at (-1.6,1.6); \coordinate (D) at (1.6,1.6);
  \coordinate (C) at (0,0);
  \draw (D) -- (E) -- (A) -- (B);
  \draw (B) -- (C) -- (D);
  \node[below] at (C) {$C$}; \node[above left] at (A) {$A$};
  \node[above right] at (E) {$E$}; \node[left] at (B) {$B$};
  \node[right] at (D) {$D$};
  \draw[<->] (-4.2,4.2) -- (4.2,4.2); \node[above] at (2.2,4.2) {$c/h$};
  \draw[<->] (-3.1,0) -- (-3.1,4.2); \node[left] at (-3.1,2.1) {$c/h^2$};
  \draw[<->] (2.9,0) -- (2.9,1.6); \node[right] at (2.9,0.8) {$4\,|\log h|/t$};
\end{tikzpicture}
\caption{Contour of integration for the integral in \eqref{ch6:eq:3.66}. $\mathcal{C}_+$ is the
contour joining the point $D$ to the points $E$, $A$, $B$. $\mathcal{C}_-$ is the contour
joining the point $B$ to $C$ to $D$.}
\label{ch6:fig:1}
\end{figure}

A symbolic path $\gamma=\{\gamma_0,\gamma_1,\gamma_2,\dots\}$ is an infinite sequence of sites
in $X_m$ such that $\gamma_j\ne\gamma_{j-1}$ for all $j=1,\dots$. Let $\Gm$ be the set of all
symbolic paths in $X_m$. The kernel of the diffusion process admits the following
representation in terms of a summation over symbolic paths:
\begin{equation}\label{ch6:eq:3.14}
\widehat{u}_m(x,y;z;t)=\sum_{q=1}^{\infty}2^{-q}\!\!\!
\sum_{\substack{\gamma\in\Gm:\ \gamma_0=x,\ \gamma_q=y\\ |\gamma_j-\gamma_{j-1}|=1\ \forall j\ge1}}
\!\!\!W_m(\gamma,q,z,t)
\end{equation}
where
\begin{equation}\label{ch6:eq:3.15}
\begin{aligned}
W_m(\gamma,q,z,t)&=\frac{1}{h_{xm}}\int_0^t\!\dd s_1\int_{s_1}^t\!\dd s_2\cdots\int_{s_{q-1}}^t\!\dd s_q\,
e^{(t-s_q)\widehat{L}^m(\gamma_q,\gamma_q;z)}\\
&\qquad\times\prod_{j=0}^{q-1}2\widehat{L}^m(\gamma_j,\gamma_{j+1};z)
e^{(s_{j+1}-s_j)\widehat{L}^m(\gamma_j,\gamma_j;z)}
\end{aligned}
\end{equation}
with $s_0=0$.

Let us introduce the following Green's function:
\begin{equation}\label{ch6:eq:3.16}
G_m(x,y;z,\omega)=\int_0^{\infty}\widehat{u}_m(x,y;z,t)e^{-i\omega t}\dd t
=h_{xm}^{-1}\frac{1}{\widehat{L}^m(z)+i\omega}(x,y).
\end{equation}
The Fourier transformed propagator can be re-obtained from the Green's function by evaluating
the following contour integral:
\begin{equation}\label{ch6:eq:3.17}
\widehat{u}_m(x,y;z,t)=\int_{\mathcal{C}_-}G_m(x,y;z,\omega)e^{i\omega t}\frac{\dd\omega}{2\pi}
+\int_{\mathcal{C}_+}G_m(x,y;z,\omega)e^{i\omega t}\frac{\dd\omega}{2\pi}.
\end{equation}
Here, $\mathcal{C}_+$ is the contour joining the point $D$ to the points $E$, $A$, $B$ in
Fig.~\ref{ch6:fig:1}, while $\mathcal{C}_-$ is the contour joining the point $B$ to $C$ to $D$. The
intent of the design of the contour of integration is to ensure that each point $\omega$ on the
upper path $\mathcal{C}_+$ is separated from the spectrum of $L$.

\begin{lemma}\label{ch6:lem:1}
For all $K>0$, there are an integer $m(K)$ and a constant $c(K)$ such that
\begin{equation}\label{ch6:eq:3.18}
\left|\int_{\mathcal{C}_+}\frac{\dd\omega}{2\pi}G_m(x,y;z,\omega)e^{i\omega t}\right|\le c(K)h^2.
\end{equation}
for all $z$ such that $|z|<K$.
\end{lemma}

\begin{proof}
The proof is based on the geometric series expansion
\begin{equation}\label{ch6:eq:3.19}
\begin{aligned}
G_m(z,\omega)&=h_m^{-1}\frac{1}{\LatL(z)+i\omega}
=h_m^{-1}\sum_{j=0}^{\infty}\frac{1}{\half\sigma(0,h_{xm})^2\Dm+i\omega}\\
&\quad\times\left[\left(\frac{1}{2}(\sigma(z,h_{xm})^2-\sigma(0,h_{xm})^2)\Dm
+\mu(z,h_{xm})\nabm+\kappa(z,h_{xm})\right)\right.\\
&\qquad\qquad\left.\times\frac{1}{\half\sigma(0,h_{xm})^2\Dm+i\omega}\right]^{j}
\end{aligned}
\end{equation}
Here $\sigma(z,h_{xm})^2$ and $\mu(z,h_{xm})$ are the multiplication operators by
$\sigma(x,z,h_{xm})^2$ and $\mu(x,z,h_{xm})$, respectively. Convergence for $\omega\in
\mathcal{C}_+$ can be established by means of a Kato--Rellich type relative bound, see (Kato
1966). More precisely, for any $\alpha>0$, one can find a $\beta>0$ such that the operators
$\nabm$ and $\Dm$ satisfy the following relative bound estimate:
\begin{equation}\label{ch6:eq:3.20}
\|\nabm f\|_2\le\alpha\|\Dm f\|_2+\beta\|f\|_2.
\end{equation}
for all periodic functions $f$ and all $m\ge m_0$. This bound can be derived by observing that
$\nabm$ and $\Dm$ can be diagonalized simultaneously by a Fourier transform and by observing
that for any $\alpha>0$, one can find a $\beta>0$ such that
\begin{equation}\label{ch6:eq:3.21}
\left|\frac{\sin h_{xm}k}{h_{xm}}\right|
\le\alpha\left|\frac{\cos h_{xm}k-1}{h_{xm}^2}\right|+\beta
\end{equation}
for all $m\ge m_0$ and all $k\in B_m$.

Under the same conditions, we also have that
\begin{equation}\label{ch6:eq:3.22}
\begin{aligned}
&\left\|\left(\frac{1}{2}(\sigma(z,h_{xm})^2-\sigma(0,h_{xm})^2)\Dm+\mu(z,h_{xm})\nabm\right)f\right\|_2\\
&\qquad\le\frac{2\alpha}{\Sigma_0^2}\bigl(M(K)+2\Sigma_1(K)^2\bigr)
\left\|\frac{1}{2}\sigma(0)^2\Dm f\right\|_2+\beta\|f\|_2.
\end{aligned}
\end{equation}
Hence
\begin{equation}\label{ch6:eq:3.23}
\begin{aligned}
&\left\|\mu\nabm\frac{1}{\half\sigma^2\Dm+i\omega}f\right\|_2
\le\left(\frac{1}{2}\bigl(\sigma(z,h_{xm})^2-\sigma(0,h_{xm})^2\bigr)\Dm
+\mu(z,h_{xm})\nabm\right)\frac{1}{\half\sigma^2\Dm+i\omega}f\\
&\qquad+\beta\left\|\frac{1}{\half\sigma^2\Dm+i\omega}f\right\|_2<1
\end{aligned}
\end{equation}
where the last inequality holds if $\omega\in\mathcal{C}_+$, if $\alpha$ is chosen
sufficiently small and if $m$ is large enough. In this case, the geometric series expansion in
\eqref{ch6:eq:3.19} converges in $L^2$ operator norm. The uniform norm of the kernel
$|G_m(x,y;\omega)|$ is pointwise bounded from above by $h_m^{-2}$.

Since the points $B$ and $D$ have imaginary part of height $4\frac{|\log h_m|}{t}$, the integral
over the contour $\mathcal{C}_+$ converges also and is bounded from above by $ch_m^2$ in
uniform norm.
\end{proof}

\begin{lemma}\label{ch6:lem:2}
For all $K>0$, if $q\ge e^2\Sigma_1(K)^2t/(2h_m^2)$ we have that
\begin{equation}\label{ch6:eq:3.24}
W_m(\gamma,q;z,t)\le\sqrt{\frac{q}{2\pi}}\,
\exp\left(-\frac{\Sigma_0^2t}{2h_m^2}-q\right).
\end{equation}
\end{lemma}

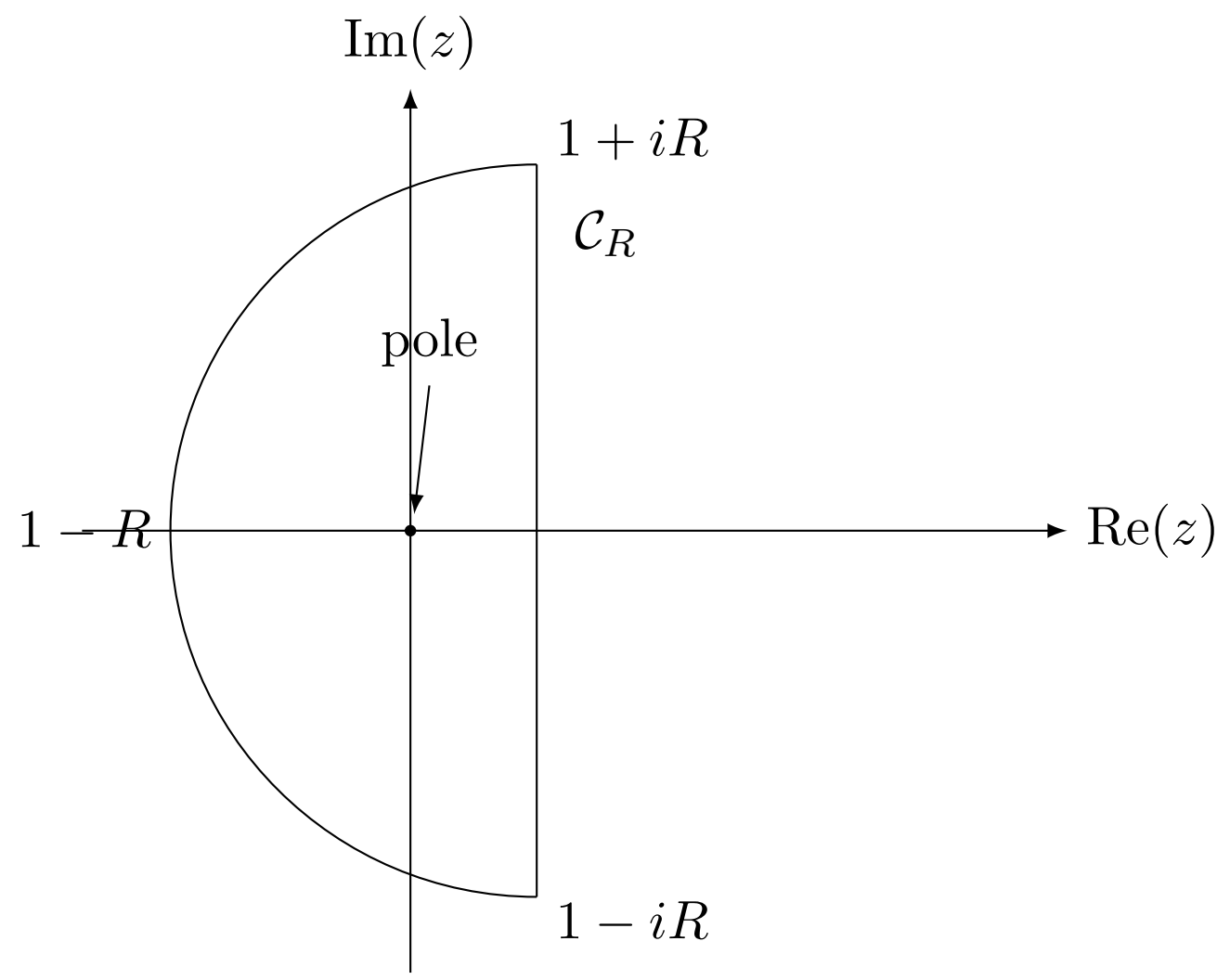
\begin{figure}[t]
\centering
\begin{tikzpicture}[scale=0.9,>=latex]
  \draw[->] (-2.6,0) -- (5.2,0) node[right] {$\mathrm{Re}(z)$};
  \draw[->] (0,-3.5) -- (0,3.5) node[above] {$\mathrm{Im}(z)$};
  \draw (1,2.9) -- (1,-2.9);
  \draw (1,2.9) arc (90:270:2.9);
  \node[left] at (-1.9,0) {$1-R$};
  \node[right] at (1.0,-3.1) {$1-iR$};
  \node[right] at (1.0,3.1) {$1+iR$};
  \node at (1.55,2.35) {$\mathcal{C}_R$};
  \fill (0,0) circle (1.3pt);
  \draw[->] (0.15,1.15) -- (0.03,0.13);
  \node[above] at (0.15,1.15) {pole};
\end{tikzpicture}
\caption{Contour of integration $\mathcal{C}_R$ for the integral in \eqref{ch6:eq:3.29}.}
\label{ch6:fig:2}
\end{figure}

\begin{proof}
Let us define the function
\begin{equation}\label{ch6:eq:3.25}
\varphi(t)=\frac{\Sigma_1(K)^2}{2h_m^2}e^{-\frac{\Sigma_0^2t}{2h_m^2}}\mathbf{1}(t\ge0)
\end{equation}
where $\mathbf{1}(t\ge0)$ is the characteristic function of $\mathbb{R}_+$. We have that
\begin{equation}\label{ch6:eq:3.26}
W_m(\gamma,q;z,t)\le\varphi^{*q}(t)
\end{equation}
for all $z$ such that $|z|<K$, where $\varphi^{*q}$ is the $q$-th convolution power, i.e.\ the
$q$-fold convolution product of the function $\varphi$ by itself. The Fourier transform of
$\varphi(t)$ is given by
\begin{equation}\label{ch6:eq:3.27}
\widehat{\varphi}(\omega)=\frac{\Sigma_1(K)^2}{2h_m^2}
\int_0^{\infty}e^{-i\omega t-\frac{\Sigma_0^2t}{2h_m^2}}\dd t
=\frac{\Sigma_1(K)^2}{2i\omega h_m^2+\Sigma_0^2}.
\end{equation}
The convolution power is given by the following inverse Fourier transform:
\begin{equation}\label{ch6:eq:3.28}
\varphi^{*q}(t)=\int_{-\infty}^{\infty}\widehat{\varphi}(\omega)^qe^{i\omega t}\frac{\dd\omega}{2\pi}
=\left(\frac{\Sigma_1(K)}{\Sigma_0}\right)^{2q}\int_{-\infty}^{\infty}
\left(1+\frac{2i\omega h_m^2}{\Sigma_0^2}\right)^{-q}e^{i\omega t}\frac{\dd\omega}{2\pi}.
\end{equation}
Introducing the new variable $z=1+\frac{2i\omega h_m^2}{\Sigma_0^2}$, the integral can be
recast as follows
\begin{equation}\label{ch6:eq:3.29}
\varphi^{*q}(t)=\frac{\Sigma_0^{2-2q}\Sigma_1(K)^{2q}}{4\pi ih_m^2}\lim_{R\to\infty}
\oint_{\mathcal{C}_R}z^{-q}\exp\left(\frac{\Sigma_0^2t}{2h_m^2}(z-1)\right)\dd z
\end{equation}
where $\mathcal{C}_R$ is the contour in Fig.~\ref{ch6:fig:2}. Using the residue theorem and noticing
that the only pole of the integrand is at $z=0$, we find
\begin{equation}\label{ch6:eq:3.30}
\varphi^{*q}(t)=\frac{1}{t}\,\frac{1}{(q-1)!}
\left(\frac{\Sigma_1(K)^2t}{2h_m^2}\right)^{q}
\exp\left(-\frac{\Sigma_0^2t}{2h_m^2}\right).
\end{equation}
Making use of Stirling's formula $q!\approx\sqrt{2\pi q}\,q^{q+\frac12}e^{-q}$, we find
\begin{equation}\label{ch6:eq:3.31}
\varphi^{*q}(t)\approx\frac{1}{t}\sqrt{\frac{q}{2\pi}}\,
\exp\left(-\frac{\Sigma_0^2t}{2h_m^2}+q\log\frac{\Sigma_1(K)^2t}{2h_m^2}+q(1-\log q)\right).
\end{equation}
If $\log q\ge\log\frac{\Sigma_1(K)^2t}{2h_m^2}+2$, then we arrive at the bound in
\eqref{ch6:eq:3.24}.
\end{proof}

Let us fix a $K>0$. To prove the theorem, it suffices to consider the case $m_0=m+1$ for all
values of $m$ above $m_0$. In fact, given this chaptericular case, the general statement can be
derived with an iterative argument. To this end, we introduce a renormalization group
transformation based on the notion of decorating path.

\begin{definition}[Decorating Paths]\label{ch6:def:3}
Let $m\ge m_0$ and let $\gamma=\{y_0,y_1,y_2,\dots\}$ be a symbolic sequence in $\Gm$. A
\emph{decorating path} around $\gamma$ is defined as a symbolic sequence
$\gamma'=\{y'_0,y'_1,y'_2,\dots\}$ with $y'_i\in h_{m+1}\mathbb{Z}$ containing the sequence
$\gamma$ as a subset and such that if $y'_j=y_i$ and $y'_k=y_{i+1}$, then all elements $y'_n$
with $j<n<k$ are such that $|y'_n-y'_j|\le h_{m+1}$. Let $D_{m+1}(\gamma)$ be the set of all
decorating sequences around $\gamma$. The decorated weights are defined as follows:
\begin{equation}\label{ch6:eq:3.32}
\widetilde{W}_m(\gamma,q;z;t)=\sum_{q'=q}^{\infty}
\sum_{\substack{\gamma'\in D_{m+1}(\gamma)\\ \gamma'_{q'}=\gamma_q}}W_{m+1}(\gamma',q';z;t).
\end{equation}
\end{definition}
Let us notice that these weights are not positive or even real values unless $z=0$, but rather
they depend analytically on the variable $z$. Finally, let us introduce also the following
Fourier transform:
\begin{equation}\label{ch6:eq:3.33}
\widehat{W}_m(\gamma,q;z;\omega)=\int_0^{\infty}W_m(\gamma,q;z;t)e^{i\omega t}\dd t,\qquad
\widehat{\widetilde{W}}_m(\gamma,q;z;\omega)=\int_0^{\infty}\widetilde{W}_m(\gamma,q;z;t)e^{i\omega t}\dd t.
\end{equation}

\begin{notation}\label{ch6:not:1}
In the following, we set $h=h_{x,m+1}$ so that $h_{xm}=2h$. We also use the Landau notation
$O(h^n)$ to indicate a function $f(h)$ such that $h^{-n}f(h)$ is bounded in a neighborhood of
$0$.
\end{notation}

\begin{lemma}\label{ch6:lem:3}
Let $x,y\in X_m$ and let $\mathcal{C}_-$ be an integration contour as in Fig.~\ref{ch6:fig:1}. Then
\begin{equation}\label{ch6:eq:3.34}
\left(\int_{\mathcal{C}_-}\Bigl(2G_{m+1}(x,y;z;\omega)-G_m(x,y;z;\omega)\Bigr)e^{i\omega t}
\frac{\dd\omega}{2\pi}\right)=O(h^2).
\end{equation}
\end{lemma}

\begin{proof}
We have that
\begin{equation}\label{ch6:eq:3.35}
2G_{m+1}(x,y;z;\omega)-G_m(x,y;z;\omega)=\frac{1}{h}\sum_{q=1}^{\infty}2^{-q}\!\!\!
\sum_{\substack{\gamma\in\Gm:\ \gamma_0=x,\ \gamma_q=y\\ |\gamma_j-\gamma_{j-1}|=1\ \forall j\ge1}}
\!\!\!\Bigl(2\widehat{\widetilde{W}}_m(\gamma,q;z;\omega)-\widehat{W}_m(\gamma,q;z;\omega)\Bigr).
\end{equation}
The number of paths over which the summation is extended is
\begin{equation}\label{ch6:eq:3.36}
N(\gamma,q;x,y)\equiv\#\{\gamma\in\Gm:\gamma_0=x,\gamma_q=y,|\gamma_j-\gamma_{j-1}|=1\ \forall j\ge1\}
=\binom{q}{\frac{q}{2}+k}
\end{equation}
where $k=\frac{|y-x|}{h_{xm}}$. Applying Stirling's formula we find
\begin{equation}\label{ch6:eq:3.37}
N_\gamma\sim2^q\sqrt{\frac{2}{\pi q}}.
\end{equation}
Hence
\begin{equation}\label{ch6:eq:3.38}
\begin{aligned}
&\left|\int_{\mathcal{C}_-}\Bigl(2G_{m+1}(x,y;z;\omega)-G_m(x,y;z;\omega)\Bigr)e^{i\omega t}
\frac{\dd\omega}{2\pi}\right|\\
&\quad\le\frac{c}{h}\sum_{q=1}^{\infty}\frac{1}{q}
\max_{\substack{\gamma\in\Gm:\ \gamma_0=x,\ \gamma_q=y\\ |\gamma_j-\gamma_{j-1}|=1\ \forall j\ge1}}
\left|\int_{\mathcal{C}_-}\Bigl(2\widehat{\widetilde{W}}_m(\gamma,q;z,\omega)
-\widehat{W}_m(\gamma,q;z;\omega)\Bigr)e^{i\omega t}\frac{\dd\omega}{2\pi}\right|.
\end{aligned}
\end{equation}
for some constant $c\approx\sqrt{\frac{2}{\pi}}>0$. It suffices to extend the summation over
$q$ only up to
\begin{equation}\label{ch6:eq:3.39}
q_{\max}\equiv\frac{e^2\Sigma_1(K)^2t}{2h^2}.
\end{equation}
To resum beyond this threshold, one can use the previous lemma. More precisely, we have that
\begin{equation}\label{ch6:eq:3.40}
\begin{aligned}
&\left|\int_{\mathcal{C}_-}\Bigl(2G_{m+1}(x,y;z;\omega)-G_m(x,y;z;\omega)\Bigr)e^{i\omega t}
\frac{\dd\omega}{2\pi}\right|\\
&\quad\le\frac{c\sqrt{q_{\max}}}{h}
\max_{\substack{\gamma\in\Gm:\ \gamma_0=x,\ \gamma_q=y\\ |\gamma_j-\gamma_{j-1}|=1\ \forall j\ge1}}
\left|\int_{\mathcal{C}_-}\Bigl(2\widehat{\widetilde{W}}_m(\gamma,q;z,\omega)
-\widehat{W}_m(\gamma,q;z,\omega)\Bigr)e^{i\omega t}\frac{\dd\omega}{2\pi}\right|.
\end{aligned}
\end{equation}
Let us introduce the following abbreviated notations:
\begin{equation}\label{ch6:eq:3.41}
v(x,z)=\sigma(x,z,h)^2,\qquad
m(x,z)=\mu(x,z,h),\qquad
k(x,z)=\kappa(x,z,h)
\end{equation}
where for the sake of keeping formulas short, we omit to denote the $h=h_{x,m+1}$ dependencies.
To evaluate the resummed weight function, let us form the matrix
\begin{equation}\label{ch6:eq:3.44}
\setlength{\arraycolsep}{1pt}
\LbarCal(x,z)=
\begin{pmatrix}
-\dfrac{v(x+h,z)}{h^2}+k(x+h,z) & \dfrac{v(x+h,z)}{2h^2}-\dfrac{m(x+h,z)}{2h} & 0\\[8pt]
\dfrac{v(x,z)}{2h^2}+\dfrac{m(x,z)}{2h} & -\dfrac{v(x,z)}{h^2}+k(x,z) & \dfrac{v(x,z)}{2h^2}-\dfrac{m(x,z)}{2h}\\[8pt]
0 & \dfrac{v(x-h,z)}{2h^2}+\dfrac{m(x-h,z)}{2h} & -\dfrac{v(x-h,z)}{h^2}+k(x-h,z)
\end{pmatrix}
\end{equation}
and decompose it as follows:
\begin{equation}\label{ch6:eq:3.45}
\LbarCal(x,z)=\frac{1}{h^2}\LbarCal_0(x,z)+\frac{1}{h}\LbarCal_1(x,z)+\LbarCal_2(x,z)+h\LbarCal_3(x,z)+\!O(h^2).
\end{equation}
where
\begin{equation}\label{ch6:eq:3.46}
\LbarCal_0(x,z)=\begin{pmatrix}
-v(x,z) & \half v(x,z) & 0\\
\half v(x,z) & -v(x,z) & \half v(x,z)\\
0 & \half v(x,z) & -v(x,z)
\end{pmatrix},
\end{equation}
\begin{equation}\label{ch6:eq:3.47}
\LbarCal_1(x,z)=\begin{pmatrix}
-\nabm_xv(x,z) & \half\nabm_xv(x,z)-\half m(x,z) & 0\\
\half m(x,z) & 0 & -\half m(x,z)\\
0 & -\half\nabm_xv(x,z)+\half m(x,z) & \nabm_xv(x,z)
\end{pmatrix},
\end{equation}
\begin{equation}\label{ch6:eq:3.48}
\LbarCal_2(x,z)=\begin{pmatrix}
-\half\Dm_xv(x,z)+k(x,z) & \tfrac14\Dm_xv(x,z)-\half\nabm_xm(x,z) & 0\\
0 & k(x,z) & 0\\
0 & \tfrac14\Dm_xv(x,z)-\half\nabm_xm(x,z) & -\half\Dm_xv(x,z)+k(x,z)
\end{pmatrix}.
\end{equation}
and
\begin{equation}\label{ch6:eq:3.49}
\LbarCal_3(x,z)=\begin{pmatrix}
\nabm_xk(x,z) & -\tfrac14\Dm_xm(x,z) & 0\\
0 & 0 & 0\\
0 & \tfrac14\Dm_xm(x,z) & -\nabm_xk(x,z)
\end{pmatrix}.
\end{equation}
Recall that all functions and operators above depend also on $h=h_{x,m+1}$.

Let us introduce the sign variable $\tau=\pm1$, the functions
\begin{equation}\label{ch6:eq:3.50}
\varphi_0(t,x,z,\tau)\equiv2\widehat{L}^m_z(x,x+2\tau h;z)e^{t\widehat{L}^m_z(x,x;z)}\mathbf{1}(t\ge0)
\end{equation}
\begin{equation}\label{ch6:eq:3.51}
\varphi_1(t,x,z,\tau)\equiv2\widehat{L}^{m+1}_z(x+\tau h,x+2\tau h;z)e^{t\LbarCal(x,z)}(x,x+\tau h)\mathbf{1}(t\ge0)
\end{equation}
and their Fourier transforms
\begin{equation}\label{ch6:eq:3.52}
\begin{aligned}
\widehat{\varphi}_0(\omega,x,z,\tau)&=\left(\frac{v(x,z)}{4h^2}+\tau\frac{m(x,z)}{2h}\right)
\left(\frac{v(x,z)}{4h^2}-k(x,z)+i\omega\right)^{-1}\\
\widehat{\varphi}_1(\omega,x,z,\tau)&=\left(\frac{v(x,z)}{h^2}+\tau\frac{m(x,z)+\nabm_xv(x,z)}{h}\right.\\
&\qquad\left.+\frac{\Dm_xv(x,z)+\nabm_xm(x,z)}{2}+\frac{\Dm_xm(x,z)}{2}\tau h+O(h^2)\right)\\
&\qquad\times\bigl\langle x\bigl|(-\LbarCal(x,z)+i\omega)^{-1}\bigr|x+\tau h\bigr\rangle.
\end{aligned}
\end{equation}
where
\begin{equation}\label{ch6:eq:3.53}
|x\rangle=\begin{pmatrix}0\\1\\0\end{pmatrix},
\qquad\text{and}\qquad
|x+\tau h\rangle=\begin{pmatrix}\delta_{\tau,1}\\0\\\delta_{\tau,-1}\end{pmatrix}.
\end{equation}
We also require the functions
\begin{equation}\label{ch6:eq:3.54}
\psi_0(t,x,z)\equiv e^{t\widehat{L}^m_z(x,x;z)}\mathbf{1}(t\ge0),\qquad
\psi_1(t,x)\equiv e^{t\LbarCal(y;h)}(x,x)\mathbf{1}(t\ge0)
\end{equation}
and the corresponding Fourier transforms
\begin{equation}\label{ch6:eq:3.55}
\widehat{\psi}_0(\omega,x,z)=\left(\frac{v(x,z)}{4h^2}+i\omega\right)^{-1},\qquad
\widehat{\psi}_1(\omega,x,z)=\bigl\langle x\bigl|(-\LbarCal(x,z)+i\omega)^{-1}\bigr|x\bigr\rangle.
\end{equation}
If $\gamma$ is a symbolic sequence, then
\begin{equation}\label{ch6:eq:3.56}
\widehat{W}_m(\gamma,q;z,\omega)=\widehat{\psi}_0(\omega,\gamma_q,z)
\prod_{j=0}^{q-1}\widehat{\varphi}_0\bigl(\omega;\gamma_j,z,\sgn(\gamma_{j+1}-\gamma_j)\bigr)
\end{equation}
\begin{equation}\label{ch6:eq:3.57}
\widehat{\widetilde{W}}_m(\gamma,q;z;\omega)=\widehat{\psi}_1(\omega,\gamma_q,z)
\prod_{j=0}^{q-1}\widehat{\varphi}_1\bigl(\omega;\gamma_j,z,\sgn(\gamma_{j+1}-\gamma_j)\bigr).
\end{equation}

Let us estimate the difference between the functions $\widehat{\varphi}_1(\omega,x,z,\tau)$ and
$\widehat{\varphi}_0(\omega,x,z,\tau)$ assuming that $\omega$ is in the contour $\mathcal{C}_-$
in Fig.~\ref{ch6:fig:2}. Retaining only terms up to order up to $O(h^3)$, we find
\begin{equation}\label{ch6:eq:3.58}
\begin{aligned}
\widehat{\varphi}_0(\omega,x,z,\tau)={}&1+\frac{2m(x,z)\tau h}{v(x,z)}
+\frac{4h^2}{v(x,z)}\bigl(k(x,z)-i\omega\bigr)\\
&+\frac{8m(x,z)\tau h^3}{v(x,z)^2}\bigl(k(x,z)-i\omega\bigr)
+\frac{16h^4}{v(x,z)^2}\bigl(k(x,z)-i\omega\bigr)^2+O(h^5).
\end{aligned}
\end{equation}
A lengthy but straightforward calculation which is best carried out using a symbolic
manipulation program, gives
\begin{equation}\label{ch6:eq:3.59}
\begin{aligned}
\widehat{\varphi}_1(\omega,x,z,\tau)={}&1+\frac{2m(x,z)\tau h}{v(x,z)}
+\frac{4h^2}{v(x,z)}\bigl(k(x,z)-i\omega\bigr)
-\bigl[8m(x,z)-\nabm_xv(x,z)\bigr]\frac{i\omega\tau h^3}{v(x,z)^2}\\
&+R(x,z)\cdot h^3\tau+h^4P_0(x,z)+i\omega h^4P_1(x,z)
-\frac{14\omega^2h^4}{v(x,z)^2}+O(h^5)
\end{aligned}
\end{equation}
where
\begin{equation}\label{ch6:eq:3.60}
\begin{aligned}
R&=\frac{1}{2v^3}\Bigl[\bigl(\Dm_xm+2\nabm_xk\bigr)v^2+\bigl(16m-2\nabm_xv\bigr)kv
-4m^3+2\nabm_xv\,m^2\\
&\qquad\qquad-2v\nabm_xv\nabm_xm-\bigl(2v\nabm_xm+v\Dm_xv-2(\nabm_xv)^2\bigr)m\Bigr],\\[4pt]
P_0&=\frac{1}{v^3}\Bigl[\Bigl(\bigl(2m-2\nabm_xv\bigr)\nabm_xk+14k^2+\bigl(-4\nabm_xm-\Dm_xv\bigr)k\Bigr)v\\
&\qquad\qquad+\bigl(-4m^2+2m\nabm_xv+2(\nabm_xv)^2\bigr)k+v^2\Dm_xk\Bigr],\\[4pt]
P_1&=\frac{1}{v^3}\Bigl[-28kv+4m^2-2m\nabm_xv+4v\nabm_xm+v\Dm_xv-2(\nabm_xv)^2\Bigr].
\end{aligned}
\end{equation}
For simplicity, we are not denoting here the dependency of all functions on $(x,z)$. We have
that
\begin{equation}\label{ch6:eq:3.61}
\begin{aligned}
&\sum_{j=0}^{q-1}\Bigl[\log\widehat{\varphi}_0\bigl(\omega;\gamma_j,z,\sgn(\gamma_{j+1}-\gamma_j)\bigr)
-\log\widehat{\varphi}_1\bigl(\omega;\gamma_j,z,\sgn(\gamma_{j+1}-\gamma_j)\bigr)\Bigr]\\
&\quad=-\sum_{j=0}^{q-1}\left[\frac{i\omega\nabm_xv(\gamma_j,z)}{v(\gamma_j,z)^2}+R(\gamma_j,z)\right]
h^3\sgn(\gamma_{j+1}-\gamma_j)\\
&\qquad+\bigl(\|P_0\|_{K,\infty}+|\omega|\|P_1\|_{K,\infty}+2|\omega|^2\|v^{-2}\|_{K,\infty}\bigr)O(h^4q)\\
&\quad=-\frac{h^2}{2}\left[i\omega\left(\frac{1}{v(\gamma_0,z)}-\frac{1}{v(\gamma_q,z)}\right)
+\bigl(G(\gamma_q,z)-G(\gamma_0,z)\bigr)\right]\\
&\qquad+\bigl(\|P_0\|_{K,\infty}+|\omega|\|P_1\|_{K,\infty}+2|\omega|^2\|v^{-2}\|_{K,\infty}\bigr)O(h^4q)
\end{aligned}
\end{equation}
where $G(x,z)$ is a primitive of $R(x,z)$, i.e.
\begin{equation}\label{ch6:eq:3.62}
G(x,z)=\int_{-L_x}^xR(x',z)\dd x'.
\end{equation}
We conclude that there is a constant $c(K)>0$ such that
\begin{equation}\label{ch6:eq:3.63}
\int_{\mathcal{C}_-}\left|\prod_{j=0}^{q-1}
\widehat{\varphi}_0\bigl(\omega;\gamma_j,z,\sgn(\gamma_{j+1}-\gamma_j)\bigr)
-\prod_{j=0}^{q-1}\widehat{\varphi}_1\bigl(\omega;\gamma_j,z,\sgn(\gamma_{j+1}-\gamma_j)\bigr)
\right|e^{i\omega t}\frac{\dd\omega}{2\pi}\le c(K)h^{\gamma}.
\end{equation}
for all $q\le q_{\max}$. Here we use the decay of $e^{i\omega t}$ in the upper half of the
complex $\omega$ plane to offset the $\omega$ dependencies in the integrand. Similar
calculations lead to the following expansions:
\begin{equation}\label{ch6:eq:3.64}
\widehat{\psi}_0(\omega,x,z)=\frac{4h^2}{v(x,z)}+O(\omega h^4),\qquad
\widehat{\psi}_1(\omega,x,z)=\frac{2h^2}{v(x,z)}+O(\omega h^4)=\half\widehat{\psi}_0(\omega,x,z)+O(\omega h^4).
\end{equation}
Since $q<ch^{-2}$ and $\omega\le|\log h|$, we find
\begin{equation}\label{ch6:eq:3.65}
\left|\int_{\mathcal{C}_-}\Bigl(2G_{m+1}(x,y;z,\omega)-G_m(x,y;z,\omega)\Bigr)e^{i\omega t}
\frac{\dd\omega}{2\pi}\right|\le c(K)\frac{\sqrt{q_{\max}}}{h}h^{2+\gamma}\le ch^{\gamma}.
\end{equation}
By differentiating with respect to time in equation \eqref{ch6:eq:3.66}, we find that
\begin{equation}\label{ch6:eq:3.66}
\frac{\partial}{\partial t}u_m(x,y;t)=\int_{\mathcal{C}_-}i\omega G_m(x,y;z;\omega)e^{i\omega t}
\frac{\dd\omega}{2\pi}+\int_{\mathcal{C}_+}i\omega G_m(x,y;z;\omega)e^{i\omega t}\frac{\dd\omega}{2\pi}.
\end{equation}
All the derivations above carry through. We conclude that for some constant $c_1(K)>0$ and all
$|z|<K$ we have that
\begin{equation}\label{ch6:eq:3.67}
\left|\int_{\mathcal{C}_+}\frac{\dd\omega}{2\pi}i\omega G_m(x,y;z;\omega)e^{i\omega t}\right|\le c_1(K)h^{\gamma}
\end{equation}
and also
\begin{equation}\label{ch6:eq:3.68}
\left|\int_{\mathcal{C}_-}i\omega\Bigl(2G_{m+1}(x,y;z;\omega)-G_m(x,y;z;\omega)\Bigr)e^{i\omega t}
\frac{\dd\omega}{2\pi}\right|\le c\frac{\sqrt{q_{\max}}}{h}h^4\le c_1(K)h^{\gamma}.
\end{equation}
Hence, the first time derivatives of the kernel satisfy a similar Cauchy convergence condition
as the kernel itself.

Let us notice that, if $m$ is large enough, then the operator $\widehat{L}^m(z)$ is relatively
bounded with respect to $\LatL$ in the following sense:
\begin{equation}\label{ch6:eq:3.69}
\|Lf\|\le\alpha(K)\|\widehat{L}(z)f\|+\beta(K)\|f\|
\end{equation}
for some $\alpha(K),\beta(K)>0$ and all $|z|<K$. By inspection of the terms in
$\widehat{L}(z)$ one sees that such bound holds if $h_{xm}$ is small enough, i.e.\ for $m$
larger than a threshold depending on $K$. A similar relative bound holds for the adjunct
operator $L^{*}$ with respect to $\widehat{L}(z)^{*}$. This concludes the proof of Theorem 1.
\end{proof}

\section{Explicit Euler Scheme}

In this section we prove Theorem 2. A path-wise expansion for the time-discretization of the
Fourier transformed kernel has the form
\begin{equation}\label{ch6:eq:4.1}
\begin{aligned}
\widehat{u}^{\delta t}_m(x,y;z,t)={}&\frac{1}{h_{xm}}\sum_{q=1}^{\infty}
\sum_{\gamma\in\Gm:\gamma_0=x,\gamma_q=y}
\sum_{k_1=1}^{N}\sum_{k_2=k_1+1}^{N}\cdots\sum_{k_q=k_{q-1}+1}^{N}\\
&\times\bigl(1+\delta t\widehat{L}^m(\gamma_0,\gamma_0;z)\bigr)^{k_1-1}(\delta t)^q\\
&\qquad\times\prod_{j=1}^{q}\widehat{L}^m(\gamma_{j-1},\gamma_j;z)
\bigl(1+\delta t\widehat{L}^m(\gamma_j,\gamma_j;z)\bigr)^{k_{j+1}-k_j-1}
\end{aligned}
\end{equation}
where $t_{q+1}=t$ and $k_{q+1}=N$. In this case, the propagator can be expressed through a
Fourier integral as follows:
\begin{equation}\label{ch6:eq:4.2}
\widehat{u}^{\delta t}_m(x,y;z,t)=\int_{-\frac{\pi}{\delta t}}^{\frac{\pi}{\delta t}}
G^{\delta t}_m(x,y;z,\omega)e^{i\omega t}\frac{\dd\omega}{2\pi}
\end{equation}
where
\begin{equation}\label{ch6:eq:4.3}
G^{\delta t}_m(x,y;z,\omega)=\delta t\sum_{j=0}^{\frac{t}{\delta t}}
\widehat{u}^{\delta t}_m(x,y;z,j\delta t)e^{-i\omega j\delta t}.
\end{equation}
The Fourier transformed propagator can also be represented as the limit
\begin{equation}\label{ch6:eq:4.4}
\widehat{u}^{\delta t}_m(x,y;z,t)=\lim_{H\to\infty}\int_{\mathcal{C}_H}
G^{\delta t}_m(x,y;z,\omega)e^{i\omega t}\frac{\dd\omega}{2\pi}
\end{equation}
where $\mathcal{C}_H$ is the contour in Fig.~\ref{ch6:fig:3}. This is due to the fact that the
integral along the segments $BC$ and $DA$ are the negative of each other, while the integral
over $CD$ tends to zero exponentially fast as $\Im(\omega)\to\infty$, where $\Im(\omega)$ is the
imaginary part of $\omega$. Using Cauchy's theorem, the contour in Fig.~\ref{ch6:fig:3} can be
deformed into the contour in Fig.~\ref{ch6:fig:1}. To estimate the discrepancy between the
time-discretized kernel and the continuous time one, one can thus compare the Green's function
along such contour. Again, the only arc that requires detailed attention is the arc $BCD$, as
the integral over rest of the contour of integration can be bounded from above as in the
previous section.

\begin{figure}[t]
\centering
\begin{tikzpicture}[scale=0.9,>=latex]
  \draw[->] (-5.8,0) -- (5.8,0) node[right] {$\mathrm{Re}(\omega)$};
  \draw[->] (0,-1.2) -- (0,4.8) node[above] {$\mathrm{Im}(\omega)$};
  \draw (0,0) rectangle (4.0,3.6);
  \node[below left] at (-4.0,0) {$A$};
  \node[below] at (-4.0,-0.35) {$-\pi/\delta t$};
  \node[below right] at (4.0,0) {$B$};
  \node[below] at (4.0,-0.35) {$\pi/\delta t$};
  \node[above right] at (4.0,3.6) {$C$};
  \node[above] at (4.0,3.95) {$iH+\pi/\delta t$};
  \node[above left] at (-4.0,3.6) {$D$};
  \node[above] at (-4.0,3.95) {$iH-\pi/\delta t$};
\end{tikzpicture}
\caption{Contour of integration for the integral in \eqref{ch6:eq:4.4}.}
\label{ch6:fig:3}
\end{figure}

Let $h=h_{xm}$ and let us introduce the two functions
\begin{equation}\label{ch6:eq:4.5}
\varphi_0(t,x,z,\tau)\equiv2\widehat{L}^m(x,x+\tau h;z)e^{t\widehat{L}^m(x,x;z)}\mathbf{1}(t\ge0),
\end{equation}
\begin{equation}\label{ch6:eq:4.6}
\varphi_{\delta t}(j,x,z,\tau)\equiv2\widehat{L}^m(x,x+\tau h;z)
\bigl(1+\delta t\widehat{L}^m(x,x;z)\bigr)^{j-1}.
\end{equation}
and the corresponding Fourier transforms
\begin{equation}\label{ch6:eq:4.7}
\widehat{\varphi}_0(\omega,x,z,\tau)=\int_0^{\infty}\varphi_0(t,x,z,\tau)e^{-i\omega t}\dd t
=\left(\frac{v(x,z,h)}{h^2}+\tau\frac{m(x,z,h)}{h}\right)
\left(\frac{v(x,z,h)}{h^2}+i\omega\right)^{-1}
\end{equation}
\begin{equation}\label{ch6:eq:4.8}
\begin{aligned}
\widehat{\varphi}_{\delta t}(\omega,x,z,\tau)&=\sum_{j=0}^{\frac{t}{\delta t}}
\varphi_{\delta t}(j,x,z,\tau)e^{-i\omega j\delta t}\\
&=\left(\frac{v(x,z,h)}{h^2}+\tau\frac{m(x,z,h)}{h}\right)
\left(e^{i\omega\delta t}-1+\delta t\frac{v(x,z,h)}{h^2}\right)^{-1}.
\end{aligned}
\end{equation}
We have that
\begin{equation}\label{ch6:eq:4.9}
\begin{aligned}
\widehat{\varphi}_{\delta t}(\omega,x,z,\tau)
&=\left(\frac{v(x,z,h)}{h^2}+\tau\frac{m(x,z,h)}{h}\right)
\left(i\omega+\frac{v(x,z,h)}{h^2}-\frac{\omega^2}{2}\delta t+O(\delta t^2)\right)^{-1}\\
&=\widehat{\varphi}_0(\omega,x,z,\tau)+\frac{\omega^2}{2v(x,z,h)}h^2\delta t+O(h^2\delta t^2)
=\widehat{\varphi}_0(\omega,x,z,\tau)+O(h^4),
\end{aligned}
\end{equation}
where the last step uses the fact that $\delta t=O(h^2)$. Let us also introduce the functions
\begin{equation}\label{ch6:eq:4.10}
\psi_0(t,x,z,\tau)\equiv e^{t\widehat{L}^m(x,x;z)}\mathbf{1}(t\ge0),\qquad
\psi_{\delta t}(j,x,\tau)\equiv\sum_{k=1}^{j}\bigl(1+\delta t\widehat{L}^m(x,x;z)\bigr)^{j-1}.
\end{equation}
and the corresponding Fourier transforms
\begin{equation}\label{ch6:eq:4.11}
\widehat{\psi}_0(\omega,x,z,\tau)=\left(\frac{v(x,z,h)}{h^2}+i\omega\right)^{-1},\qquad
\widehat{\psi}_{\delta t}(\omega,x,\tau)=\left(e^{i\omega\delta t}-1+\delta t\frac{v(x,z,h)}{h^2}\right)^{-1}.
\end{equation}
Again we find that
\begin{equation}\label{ch6:eq:4.12}
\widehat{\psi}_0(\omega,x,z,\tau)=\widehat{\psi}_{\delta t}(\omega,x,z,\tau)+O(h^4).
\end{equation}
If $\gamma$ is a symbolic sequence, then let us set
\begin{equation}\label{ch6:eq:4.13}
\widehat{W}_m(\gamma,q;z,\omega)=\widehat{\psi}_0(\omega,\gamma_q)
\prod_{j=0}^{q-1}\widehat{\varphi}_0\bigl(\omega;\gamma_j,z,\sgn(\gamma_{j+1}-\gamma_j)\bigr)
\end{equation}
\begin{equation}\label{ch6:eq:4.14}
\widehat{W}^{\delta t}_m(\gamma,q;z,\omega)=\widehat{\psi}_{\delta t}(\omega,\gamma_q)
\prod_{j=0}^{q-1}\widehat{\varphi}_{\delta t}\bigl(\omega;\gamma_j,z,\sgn(\gamma_{j+1}-\gamma_j)\bigr).
\end{equation}
We have that
\begin{equation}\label{ch6:eq:4.15}
G^{\delta t}_m(x,y;z,\omega)-G_m(x,y;z,\omega)=\frac{1}{h}\sum_{q=1}^{\infty}2^{-q}\!\!\!
\sum_{\substack{\gamma\in\Gm:\ \gamma_0=x,\ \gamma_q=y\\ |\gamma_j-\gamma_{j-1}|=1\ \forall j\ge1}}
\!\!\!\Bigl(\widehat{W}^{\delta t}_m(\gamma,q;z,\omega)-\widehat{W}_m(\gamma,q;z,\omega)\Bigr).
\end{equation}
The integration over the contour in Fig.~\ref{ch6:fig:1} can again be split into an integration over
the contour $\mathcal{C}_-$ and an integration over the contour $\mathcal{C}_+$. The integral
over $\mathcal{C}_+$ can be bounded from above thanks to Lemma~\ref{ch6:lem:1}. Furthermore, we have
that
\begin{equation}\label{ch6:eq:4.16}
\begin{aligned}
&\left|\int_{\mathcal{C}_-}\Bigl(G^{\delta t}_m(x,y;z,\omega)-G_m(x,y;z,\omega)\Bigr)e^{i\omega t}
\frac{\dd\omega}{2\pi}\right|\\
&\quad\le ch^{-1}\sqrt{q_{\max}}
\max_{\substack{\gamma\in\Gm:\ \gamma_0=x,\ \gamma_q=y\\ |\gamma_j-\gamma_{j-1}|=1\ \forall j\ge1}}
\left|\int_{\mathcal{C}_-}\Bigl(\widehat{W}^{\delta t}_m(\gamma,q;z,\omega)
-\widehat{W}_m(\gamma,q;z,\omega)\Bigr)e^{i\omega t}\frac{\dd\omega}{2\pi}\right|\le ch^2
\end{aligned}
\end{equation}
To bound the time derivative, we have to consider
\begin{equation}\label{ch6:eq:4.17}
\int_{\mathcal{C}_-}\left(\frac{e^{i\omega\delta t}-1}{\delta t}G_m(x,y;z,\omega)
-i\omega G_m(x,y;z,\omega)\right)e^{i\omega t}\frac{\dd\omega}{2\pi}
\end{equation}
But, since $\delta t=O(h^2)$, also this difference is $O(h^2)$.

\section{Convergence in the Ito Representation}\label{ch6:sec:itorepr}

The convergence Theorem 1 admits a second formulation. Let us introduce the function
\begin{equation}\label{ch6:eq:5.1}
\Phi_m(x)=h_{xm}\sum_{x'\in X_m,\,x'\le x}a(x').
\end{equation}
We have that
\begin{equation}\label{ch6:eq:5.2}
\begin{aligned}
e^{ip\Phi_m(x)}L(x,x')e^{-ip\Phi_m(x')}={}&\left(\frac{\sigma(x)^2}{2h_{xm}^2}+\frac{\mu(x)}{2h_{xm}}\right)
e^{-ih_{xm}a(x+h_{xm})p}\delta_{x',x+h_{xm}}\\
&+\left(\frac{\sigma(x)^2}{2h_{xm}^2}-\frac{\mu(x)}{2h_{xm}}\right)e^{ih_{xm}a(x)p}\delta_{x',x-h_{xm}}
-\frac{\sigma(x)^2}{h_{xm}^2}\delta_{x'x}\\
&+\left[-i\frac{\sin(h_{xm}pa(x))}{h_{xm}}\mu(x)+\frac{1}{2}\bigl(\cos(h_{xm}a(x)p)-1\bigr)\right]\delta_{xx'}\\
&+\left(\frac{\sigma(x)^2}{2}+\frac{h_{xm}\mu(x)}{2}\right)h_{xm}\nabm_x\left(e^{-ih_{xm}a(x)p}-1\right)\delta_{xx'}\\
&\qquad+h_{xm}\nabm_x(x,x')
\end{aligned}
\end{equation}
where $\sigma(x,p,h_{xm})^2$ and $\mu(x,p,h_{xm})$ are defined in \eqref{ch6:eq:3.7} and
\eqref{ch6:eq:3.8} and
\begin{equation}\label{ch6:eq:5.3}
\nabm_xf(x)=\frac{f(x+h_{xm})-f(x)}{h_{xm}}.
\end{equation}
Hence, we have that
\begin{equation}\label{ch6:eq:5.4}
\widetilde{L}^m(x,x';p)=e^{ip\Phi_m(x)}L^m(x,x';p)e^{-ip\Phi_m(x')}
\end{equation}
where
\begin{equation}\label{ch6:eq:5.5}
L^m(x,x';p)=L(x,x')+\zeta_m(x,p,h_{xm})\delta_{xx'}+h_{xm}r_m(x,x';p).
\end{equation}
Here
\begin{equation}\label{ch6:eq:5.6}
\zeta_m(x,p,h_{xm})=ipb(x)+\frac{\sigma(x)^2}{2}\nabm_x\left(e^{-ih_{xm}a(x)p}-1\right),
\end{equation}
and
\begin{equation}\label{ch6:eq:5.7}
\nabm_xf(x)=\frac{f(x+h_{xm})-f(x)}{h_{xm}}.
\end{equation}
The operator $\widetilde{L}^m(x,x';p)$ is equivalent to $\widehat{L}^m(x,x';p)$ up to a
non-singular linear transformation. We say that $\widetilde{L}^m(x,x';p)$ is the Fourier
transformed generator in the Ito representation. Also notice that the weak limit of this
operator as $m\to\infty$ is equal to
\begin{equation}\label{ch6:eq:5.8}
\lim_{m\to\infty}L^m(x,x';p)=\frac{\sigma(x)^2}{2}\frac{\partial^2}{\partial x^2}
+\mu(x)\frac{\partial}{\partial x}+ipb(x)-\frac{ip\sigma(x)^2}{2}\frac{\partial a(x)}{\partial x}.
\end{equation}
The Fourier transformed kernel is given by
\begin{equation}\label{ch6:eq:5.9}
\widehat{u}_m(x,x';p)=e^{ip(\Phi_m(x)-\Phi_m(x'))}\widehat{U}_m(x,x';p)
\end{equation}
where
\begin{equation}\label{ch6:eq:5.10}
\widehat{U}_m(x,x';p)=h_{xm}^{-1}e^{tL^m}(x,x').
\end{equation}
The joint kernel is thus given by the following formula whose continuum analog was found in
(Girsanov 1960), (Cameron and Martin 1949), (Feynman 1948) and (Ito 1949):
\begin{equation}\label{ch6:eq:5.11}
u(x,I;x',I';t)=\int_{-\infty}^{\infty}\frac{\dd p}{2\pi}
e^{ip\left[I'-I-\Phi_m(x')+\Phi_m(x)\right]}\widehat{U}_m(x,x';p).
\end{equation}

\begin{theorem}\label{ch6:thm:3}
Under the same assumptions of Theorem~\ref{ch6:thm:1}, for all $K>0$ there is a constant $c(K)$
such that
\begin{equation}\label{ch6:eq:5.12}
\|\widehat{U}^{\Phi}_m(t)-\widehat{U}^{\Phi}_{m_0}(t)\|_{m,K,\widehat{L}}\le c(K)h_{xm}^{\gamma}.
\end{equation}
for all $m_0>m\ge m_0$. A similar bound also holds for the kernels obtained with a fully
explicit Euler scheme, i.e.
\begin{equation}\label{ch6:eq:5.13}
\widehat{U}^{\delta t}_m(x,x';p,t)=\frac{1}{h_{xm}}
\bigl(1+\delta t_mL^m(p)\bigr)^{\left[\frac{t}{\delta t_m}\right]}(x,x'),
\end{equation}
where $\delta t_m$ is so small that
\begin{equation}\label{ch6:eq:5.14}
\min_{x,p\in X_m\times\Yhn}1+\delta t_m\widehat{L}^m(x,x;p)>0.
\end{equation}
\end{theorem}
In this case, there is a constant $c(K)$ such that
\begin{equation}\label{ch6:eq:5.15}
\|U_m(t)-U^{\delta t}_m(t)\|_{m,K,\widehat{L}}\le c(K)h_{xm}^2
\end{equation}
for all $m_0>m\ge m_0$.

It is possible to retrace the argument in the previous section except for replacing
\begin{equation}\label{ch6:eq:5.16}
\kappa(x,p,h_{xm})\to\zeta_m(x,p,h_{xm})+h_{xm}r_m(x,p,h_{xm})
\end{equation}
where $\kappa(x,p,h_{xm})$ is the function in \eqref{ch6:eq:3.9}. All arguments go through
unaltered under the same conditions. As a consequence we conclude that the kernels
\begin{equation}\label{ch6:eq:5.17}
\frac{1}{h_{xm}}\exp\bigl(t(L^m+\zeta_m(z)+h_{xm}r_m(z))\bigr)(x,x')
\end{equation}
converge in graph norm, uniformly on discs $z\in\mathbb{C}:|z|<K$, for any $K>0$.

\section{Other Abelian Processes}

In this Section, we give two examples of Abelian processes which emerge from applications and
are not stochastic integrals.

\subsection{The Sup Process}
Consider the sup process
\begin{equation}\label{ch6:eq:6.1}
y_t=\sup_{s\in[0,t]}x_s.
\end{equation}
Since the sup is always attained as an element of the underlying space $X_m$, it is natural in
this case to restrict the attention to the case where $n=m$ and $Y_m=X_m$. The joint generator
is given by
\begin{equation}\label{ch6:eq:6.2}
\widetilde{L}^m(x,y;x',y')=\LatL(x;x')A^m(x,y;x',y')
\end{equation}
where $x,y,x',y'\in X_m$ and we set
\begin{equation}\label{ch6:eq:6.3}
A^m(x,y;x',y')=\begin{cases}\delta_{yy'}&\text{if }x'<y\\ \delta_{x'y'}&\text{if }x'\ge y.\end{cases}
\end{equation}
Consider the matrix
\begin{equation}\label{ch6:eq:6.4}
V_m(x,y;x',y')=\delta(x-x')\mathbf{1}(y'\ge y)
\end{equation}
and its inverse
\begin{equation}\label{ch6:eq:6.5}
V_m^{-1}(x,y;x',y')=\delta(x-x')\bigl(\delta(y'-y)-\delta(y'-y-h_{xm})\bigr).
\end{equation}
Consider the one parameter family of operators $\widehat{L}^m(y)$ such that
\begin{equation}\label{ch6:eq:6.6}
\bigl(V_m^{-1}\widetilde{L}V_m\bigr)(x,y;x',y')=\widehat{L}^m(x,x';y)\delta_{yy'}.
\end{equation}
We have that
\begin{equation}\label{ch6:eq:6.7}
\widehat{L}^m(x,x';y)=\LatL(x;x')\mathbf{1}(x'\le y).
\end{equation}
Hence, $\widehat{L}^m(x,x';y)$ is the Markov generator of the underlying process with absorption
in the interval $[y,L_x]$. Let us notice that the kernel can be obtained as follows:
\begin{equation}\label{ch6:eq:6.8}
u_m(t)=V_m\widetilde{u}_m(t)V_m^{-1}
\quad\text{where}\quad
\widetilde{u}_m(x,y;x',y';t)=\delta(y-y')\exp\bigl(t\widehat{L}^m(y)\bigr)(x,x').
\end{equation}
A more explicit way of expressing the joint kernel is
\begin{equation}\label{ch6:eq:6.9}
u_m(x,y;x',y';t)=\delta(y-y')\widehat{u}_m(x,x';y,t)
+\mathbf{1}(y>y')\bigl(\widehat{u}_m(x,x';y',t)-\widehat{u}_m(x,x';y'-h_{xm},t)\bigr).
\end{equation}
Convergence in this case can be established along the same lines as done for stochastic
integrals. The situation is simpler in that only the consideration of the kernel itself, i.e.\
the $p=z=0$ case with the notations in Sections 2 and 3, is needed. The additional complication
is that we need to consider absorbing boundary conditions. This implies a few marginal changes
to the derivation above, as when a path arrives to an absorbing lattice point, it stays
constant thereafter. Since this applies both to paths and decorating paths in a finer lattice,
upon arriving to an absorption point the dynamics is trivial in either case and the final
bounds given still hold.

\subsection{Discrete Time Processes}
This section is based on work in collaboration with Manlio Trovato (Albanese and Trovato 2005)
and Paul Jones, see (Albanese and Jones 2007b).

An important class of path-dependent options requires computing the joint distribution of the
underlying lattice process and of a discrete sum of the following form:
\begin{equation}\label{ch6:eq:6.10}
y_t=\sum_{i=1}^{N}\psi(x_{t_{i-1}},x_{t_i};t_i)
\end{equation}
where $N$ is an integer, $t_i=i\Delta T$ and $T=N\Delta T$. Consider the elementary propagator
\begin{equation}\label{ch6:eq:6.11}
U_m(x_1,x_2)=\frac{1}{h_{xm}}e^{(\Delta T)L}(x_1,x_2).
\end{equation}
To find the joint transition probability, one can again discretize the variable $y_t$ in the
lattice $Y_n=h_{yn}\mathbb{Z}\cap[-L_y,L_y]$. As opposed to lifting the generator as done above
for the other cases, here we lift the elementary propagator itself and form the joint propagator
\begin{equation}\label{ch6:eq:6.12}
\widetilde{U}_{mn}(x_1,y_1;x_2,y_2)=U_m(x_1,x_2)\,
\delta\bigl(y_1-y_2+\psi(x_1,x_2)h_{yn}^{-1}\bigr).
\end{equation}
This lifted operator can be block-diagonalized by means of a partial Fourier transform.
Consider the Fourier transform operator $\widehat{U}_m(p)$ of matrix elements
\begin{equation}\label{ch6:eq:6.13}
\widehat{U}_m(x_1,x_2;p)=\lim_{n\to\infty}U_{mn}(x_1,y_1;x_2,y_2)e^{-ip(y_2-y_1)}
=U_m(x_1,x_2)e^{-ip\psi(x_1,x_2)}.
\end{equation}
Then we have that
\begin{equation}\label{ch6:eq:6.14}
\lim_{n\to\infty}\widetilde{U}^{N}_{mn}(x_1,y_1;x_2,y_2)
=\int_{-\infty}^{\infty}\frac{\dd p}{2\pi}e^{ip(y_2-y_1)}\widehat{U}^{N}_m(x_1,x_2;p)
\end{equation}
Convergence in the graph-norm in this case descends directly from the convergence of the
one-period kernel $U_m(x_1,x_2)$.

\section{A Gauge-Free Route to the Path Sum}\label{ch6:sec:gaugefree}

The renormalisation argument of \S\ref{ch6:sec:rg} controls the sum over paths by a change of variables
which involves a primitive of the drift, and such a primitive exists because the base space is an
interval. A higher dimensional base has no such device, so if the argument is to be extended the sum
over paths has to be controlled without it. This section develops that route in one dimension,
where the answer is known independently and can therefore be used to calibrate it, and then states
precisely what a higher dimensional version would have to assert.

\subsection{The Object to Be Controlled}

Fix a time horizon $t$, a step $h$, and let $q=\sigma^{2}t/h^{2}$ be the number of steps, so that
$qh^{2}=\sigma^{2}t$.
Let $x_0=0$ and $x_{j+1}=x_j+\Delta_{j+1}$ with the increments independent, taking the values $\pm h$
with probabilities
\begin{equation}\label{ch6:eq:step}
p_{\pm}=\frac12\pm\frac{\mu h}{2\sigma^{2}},
\end{equation}
so that $\mathbb{E}[\Delta]=\mu h^{2}/\sigma^{2}$ and $\mathbb{E}[\Delta^{2}]=h^{2}$; the time step of the
rates in \eqref{ch6:eq:2.23} is $h^{2}/\sigma^{2}$, which is why $qh^{2}=\sigma^{2}t$. This is the
step distribution of \eqref{ch6:eq:2.23}. For a function $c$ on the base, put
\begin{equation}\label{ch6:eq:pathsum}
I_{q}(c)=\sum_{j=0}^{q-1}c(x_{j})\Delta_{j+1},
\end{equation}
the lattice approximation to the line integral $\int c\dd x$ along the path. The question the
renormalisation step asks is whether $\mathbb{E}[I_{q}(c)]$, which is what survives after summing
over paths, is an endpoint functional to the order claimed.

\subsection{The Averaged Path Sum Is the It\^o Functional}

\begin{proposition}\label{ch6:prop:ito}
Let $c\in C^{3}$ with primitive $C$, let $x$ be as above and let $qh^{2}=\sigma^{2}t$. Then
\begin{equation}\label{ch6:eq:itofunctional}
\mathbb{E}\bigl[I_{q}(c)\bigr]
=\mathbb{E}\bigl[C(x_{q})\bigr]-C(x_{0})
-\frac{\sigma^{2}}{2}\int_{0}^{t}\mathbb{E}\bigl[c'(x_{s})\bigr]\dd s+O(h^{2}),
\end{equation}
the $O(h^{2})$ being uniform in $q$ and in $c$ on bounded sets of $C^{3}$.
\end{proposition}

\begin{proof}
Taylor's theorem in the form
$c(x)\Delta=C(x+\Delta)-C(x)-\tfrac12c'(x)\Delta^{2}-\tfrac16c''(x)\Delta^{3}+O(\abs{\Delta}^{4})$
summed over the path gives
\begin{equation}\label{ch6:eq:taylorsum}
I_{q}(c)=C(x_{q})-C(x_{0})-\tfrac12\sum_{j}c'(x_{j})\Delta_{j+1}^{2}
-\tfrac16\sum_{j}c''(x_{j})\Delta_{j+1}^{3}+O\Bigl(\sum_{j}\abs{\Delta_{j+1}}^{4}\Bigr).
\end{equation}
Three estimates, each conditional on $x_{j}$. First,
$\mathbb{E}[\Delta_{j+1}^{2}\mid x_{j}]=h^{2}$ exactly, so the second term has expectation
$\tfrac12h^{2}\sum_{j}\mathbb{E}[c'(x_{j})]$, and
$h^{2}\sum_{j}\mathbb{E}[c'(x_{j})]=\sigma^{2}\int_{0}^{t}\mathbb{E}[c'(x_{s})]\dd s+O(h^{2})$, the error
being the Riemann sum error together with the difference between the lattice and the continuum laws
of $x_{s}$, which is $O(h^{2})$ by the estimates of \S\ref{ch6:sec:rg}. Second,
$\mathbb{E}[\Delta_{j+1}^{3}\mid x_{j}]=h^{3}(p_{+}-p_{-})=O(h^{4})$ by \eqref{ch6:eq:step}, so the
third term has expectation $O(qh^{4})=O(h^{2})$. Third, $\mathbb{E}[\abs{\Delta_{j+1}}^{4}]=O(h^{4})$
per step, so the remainder has expectation $O(h^{2})$. Collecting the three gives \eqref{ch6:eq:itofunctional}.
\end{proof}

\begin{remark}[The smoothness demanded, and why one dimension needs no averaging]
The hypothesis is placed on $c$ and not on the coefficients of the diffusion, and it is exactly what
the argument uses: a Taylor expansion to third order with a remainder in the fourth derivative. If $c$
is taken to be the coefficient of $h^{3}$ in the comparison of two resolutions, as it is in Chapter~7,
then $c$ itself costs two derivatives of the coefficients, so the hypothesis costs five; the count is
recorded in Chapter~7, \S7.6.2.

That is a stronger demand than the one dimensional rate theorems of this chapter make, and the reason
they do not need it is structural rather than analytic: nothing in the one dimensional argument
invokes the averaged comparison of Chapter~7, the obstruction being absent because there are no two
forms. The identity above holds path by path, up to a remainder which the comparison between
resolutions cancels, so no hypothesis on the mean of a line integral arises.
\end{remark}

\begin{remark}[The identity is It\^o's lemma, and the martingale term drops out exactly]
The right hand side of \eqref{ch6:eq:itofunctional} is the expectation of the It\^o integral, since
It\^o's lemma applied to $C(x_{s})$ gives
$\int c\dd x=C(x_{t})-C(x_{0})-\tfrac{\sigma^{2}}{2}\int_{0}^{t}c'(x_{s})\dd s$. What the
proposition adds is that the lattice average reproduces it at second order, and the mechanism is
worth isolating: the term which would be a martingale, namely
$\sum_{j}c'(x_{j})\bigl(\Delta_{j+1}^{2}-\sigma^{2}h^{2}\bigr)$, has expectation \emph{exactly}
zero, no estimate being needed. That exactness is the harmless fluctuation of the higher
dimensional case as well.
\end{remark}

The proposition can be checked to the last digit, because for polynomial $c$ the moments of $x_{j}$
satisfy a closed recursion: with $M_{k}(j)=\mathbb{E}[x_{j}^{k}]$ and
$E_{r}=\mathbb{E}[\Delta^{r}]=h^{r}\bigl(p_{+}+(-1)^{r}p_{-}\bigr)$,
\begin{equation}\label{ch6:eq:momentrecursion}
M_{k}(j+1)=\sum_{l=0}^{k}\binom{k}{l}M_{l}(j)E_{k-l},\qquad M_{0}=1,\quad M_{k}(0)=0 \ (k\ge1),
\end{equation}
so that $\mathbb{E}[I_{q}(c)]=\mu h^{2}\sum_{j=0}^{q-1}\mathbb{E}[c(x_{j})]$ is obtained exactly for
any $q$. For $\mu=\sigma^{2}=1$, $t=1$:
\begin{center}\small
\begin{tabular}{lcccc}
$c$ & limit of $\mathbb{E}[I_{q}(c)]$ & closed form & errors at $q=100,400,1600$ & ratios\\
\hline
$x^{2}/2$ & $0.4166667$ & $\tfrac16+\tfrac14$ & $7.5\cdot10^{-3},\,1.9\cdot10^{-3},\,4.7\cdot10^{-4}$ & $4.00,\,4.00$\\
$x^{3}/6$ & $0.2083333$ & $\tfrac{5}{24}$ & $6.6\cdot10^{-3},\,1.7\cdot10^{-3},\,4.2\cdot10^{-4}$ & $3.98,\,4.00$\\
\end{tabular}
\end{center}
Both limits agree with \eqref{ch6:eq:itofunctional} to seven digits --- $\mathbb{E}[C(x_{1})]/C$ less the
It\^o correction --- and the approach is at the rate $h^{2}$.

\subsection{Why the Naive Reference Fails}

It is tempting to compare $\mathbb{E}[I_{q}(c)]$ with the integral of $c$ along the mean trajectory,
which is what a pathwise argument would produce. That reference is wrong at the order of interest.
For $\mu=\sigma^{2}=1$, $t=1$ the straight line integral is $\tfrac16$ when $c=x^{2}/2$ while the
limit is $\tfrac{5}{12}$, and $\tfrac{1}{24}$ when $c=x^{3}/6$ while the limit is $\tfrac{5}{24}$:
the discrepancy is $O(1)$, not $O(h^{2})$. So a derivation which drops the change of variables and
compares with the straight line will find an $O(1)$ mismatch and conclude that the comparison fails.
The correct reference is \eqref{ch6:eq:itofunctional}, in which the endpoint term involves the
\emph{law} of the endpoint and not only the mean trajectory, and the second term is the It\^o
correction. This is the precise sense in which the naive route is misleading.

\subsection{Why This Is Enough for the Comparison of Two Resolutions}

The right hand side of \eqref{ch6:eq:itofunctional} is linear in $c$ and is the sum of two pieces: an
endpoint functional, $\mathbb{E}[C(x_{q})]-C(x_{0})$, and a functional of $c'$ averaged against the
law. Both are defined for each resolution, and both differ between two resolutions by $O(h^{2})$,
by the same kernel estimates. Consequently the \emph{ratio} of the two per-step factors differs from
an endpoint-determined factor by $O(h^{2})$, which is what the renormalisation step requires, and no
change of variables is used at any point. The cancellation over the sum of paths is therefore
available in one dimension without the gauge, and this is what makes the route a candidate for
extension.

\subsection{The Class of Constant \texorpdfstring{$\dd c$}{dc}, in Any Dimension}\label{ch6:subsec:constdc}

The one case in which the higher dimensional statement of \eqref{ch6:eq:symantisym} closes exactly is
that of a constant exterior derivative: $c$ affine, $c(x)=Ax$ with $A$ a constant matrix, so that
$\dd c=A-A^{T}$ is constant and non-zero in general. There is then no primitive, and the higher
dimensional phenomenon of \S\ref{ch6:sec:gaugefree} is present in full; nevertheless the average closes.

\begin{proposition}\label{ch6:prop:constdc}
Let $c(x)=Ax$ with $A$ a constant $d\times d$ matrix, and let the base walk be \emph{homogeneous}:
the increments are independent, with $\mathbb{E}[\Delta]=m$ and $\mathrm{Cov}(\Delta)=\Sigma h^{2}$
independent of the current position. Then, with $x_{0}=0$,
\begin{equation}\label{ch6:eq:constdc}
\mathbb{E}\Bigl[\sum_{j=0}^{q-1}c(x_{j})\cdot\Delta_{j+1}\Bigr]
=\frac{q(q-1)}{2}\,A_{kl}m^{k}m^{l}
=\frac{q(q-1)}{2}\,\bigl\langle m,\operatorname{sym}(A)\,m\bigr\rangle
\end{equation}
exactly, with no remainder.
\end{proposition}

\begin{proof}
The increment $\Delta_{j+1}$ is independent of $x_{j}$, so
$\mathbb{E}\bigl[x_{j}^{k}\Delta_{j+1}^{l}\bigr]=\mathbb{E}[x_{j}^{k}]\,\mathbb{E}[\Delta^{l}]
=jm^{k}m^{l}$, since $\mathbb{E}[x_{j}^{k}]=jm^{k}$. Summing over $j=0,\dots,q-1$ gives the first
equality of \eqref{ch6:eq:constdc}. The second follows because $m^{k}m^{l}$ is symmetric, so
$A_{kl}m^{k}m^{l}=\operatorname{sym}(A)_{kl}m^{k}m^{l}$.
\end{proof}

\begin{remark}[The antisymmetric part is invisible in the mean]\label{ch6:rem:levyzero}
The statement is stronger than it looks. It says that $\dd c$, however large, contributes
\emph{nothing} to the average: only the symmetric part of the Jacobian is seen. The mechanism is
that the antisymmetric part can only pair with the antisymmetric part of
$\mathbb{E}\bigl[\sum_{j}x_{j}\otimes\Delta_{j+1}\bigr]$, namely the L\'evy area
$\sum_{j}\bigl(x_{j}^{k}\Delta_{j+1}^{l}-x_{j}^{l}\Delta_{j+1}^{k}\bigr)$, and that area has mean
$\frac{q(q-1)}{2}\bigl(m^{k}m^{l}-m^{l}m^{k}\bigr)=0$ because the mean increment is a fixed vector.
So the non-exactness of $c$ is real but lives entirely in the \emph{fluctuation}, not in the mean.
\end{remark}

\begin{remark}[Exact check]
The identity was verified by enumerating all paths in two dimensions, with $\mu=(1,\tfrac12)$,
$h=\tfrac12$ and $\Sigma=\tfrac12 I$, for four Jacobians. The difference between the enumeration and
\eqref{ch6:eq:constdc} was $0.0$ to machine precision in every case: for $A$ the pure rotation
$\binom{\,0\;-1\,}{\,1\;\;\,0\,}$, which has $\dd c\neq0$, the mean is $0$ for $q=2,4,6$; for
$A=\operatorname{diag}(1,-1)$ it is $0.046875,\,0.28125,\,0.703125$; for the matrix
$\binom{\,0\;0\,}{\,1\;0\,}$ of Lemma~6.2 of Chapter~7 it is $0.03125,\,0.1875,\,0.46875$; and for
rotation \emph{plus} the diagonal matrix the values are identical to the diagonal case alone.
\end{remark}

\begin{remark}[The fluctuation it produces is of order one]
The L\'evy area is not small, it is merely mean zero. Measured for the pure rotation at $h=\tfrac12$,
the variance of $I_{q}$ takes the values $0.03125,\,0.26563,\,0.85938,\,1.96875$ at $q=2,4,6,8$,
whose ratios to $q^{3}h^{6}$ are $0.250,\,0.266,\,0.255,\,0.246$: the variance is
$\approx\tfrac14q^{3}h^{6}$, hence $O(1)$ in the diffusive scaling $q=t/h^{2}$. By the exponential
moment estimate that is exactly the situation in which the fluctuation enters the product one order
below the mean, so it is harmless.
\end{remark}

Consequently the gauge-free route \emph{does} extend beyond one dimension, and the extension is
proved rather than assumed, for constant $\dd c$ and homogeneous coefficients. Comparing
\eqref{ch6:eq:constdc} with the straight line integral $\tfrac12\langle Y,AY\rangle=\tfrac{q^{2}}{2}A_{kl}m^{k}m^{l}$,
the averaged path sum equals it times $1-\tfrac1q$, with no reference to exactness anywhere.

What remains, and what the next subsection isolates, is the reason the identity \eqref{ch6:eq:constdc}
requires homogeneity: the proof uses only that $\mathbb{E}[\Delta_{j+1}\mid x_{j}]$ is a fixed
vector. When the coefficients vary, that vector is $m(x_{j})$, and then
$\mathbb{E}\bigl[x_{j}^{k}m^{l}(x_{j})\bigr]$ is no longer symmetric in $k,l$. The antisymmetric
part survives in the mean, and the amount by which it does is the obstruction.

\subsection{The Higher Dimensional Statement Which Would Be Needed}

When the coefficients vary, the Taylor computation of Proposition~\ref{ch6:prop:ito} changes in one
essential place. With
$\Delta$ a vector and $c$ a field of one forms, the first term of \eqref{ch6:eq:taylorsum} exists only
if $c$ has a primitive $C$, i.e. only if $\dd c=0$, which fails in general. What survives is the
decomposition of the second term into its symmetric and antisymmetric parts, since
\begin{equation}\label{ch6:eq:symantisym}
\sum_{j}c(x_{j})\cdot\Delta_{j+1}
=\Phi(x_{q})-\Phi(x_{0})
-\tfrac12\sum_{j}\sum_{k,l}\partial_{k}c_{l}(x_{j})\Delta_{j+1}^{k}\Delta_{j+1}^{l}+O(h^{2}),
\end{equation}
where the first term is now the line integral rather than an increment, and the quadratic form
splits as $\tfrac12(\partial_{k}c_{l}+\partial_{l}c_{k})+\tfrac12(\partial_{k}c_{l}-\partial_{l}c_{k})$.
The symmetric part gives, on averaging, exactly the It\^o correction of
Proposition~\ref{ch6:prop:ito}; the antisymmetric part gives the mean swept area, which is the
higher dimensional phenomenon. For a linear field it is $\tfrac12XY$, computed exactly in
Lemma~6.2 of Chapter~7; in general it is measured by $\dd c$, the two form characterised in
\S7.6.2 of Chapter~7, where it is shown to be non-vanishing, to be a quadratic form in $\nabla b$ in
the affine case, to obey a cubic law in the drift, and to be independent of the resolvent parameter.

The statement a higher dimensional version of Proposition~\ref{ch6:prop:ito} would have to assert is
therefore: that the averaged line integral equals the line integral along the mean trajectory plus
the It\^o correction plus the mean swept area, uniformly in the coefficient data and in the
conjugating parameter, with an error of order $h^{2-\epsilon}$. In the moment form used in Chapter~7
this is Assumption~6.3 of that chapter. It is \emph{not} proved here, and we do not claim it; what is
proved here is the one dimensional case, in which the antisymmetric part is empty and the identity
closes.

\begin{remark}
The one dimensional identity of Proposition~\ref{ch6:prop:ito} also shows why the obstruction found in
Chapter~7 does not appear here. In one dimension every one form is closed, so $\dd c$ vanishes
identically and the antisymmetric term of \eqref{ch6:eq:symantisym} is absent; the lattice average is
then an It\^o functional, which is endpoint determined up to a term which the comparison between
resolutions cancels. In higher dimension the It\^o functional is not endpoint determined, and the
swept area has to be controlled by averaging --- which is what Chapter~7 does, and what a
higher dimensional renormalisation would have to do in place of the change of variables used in
\S\ref{ch6:sec:rg}.
\end{remark}

\section{Machine-Checked Verification}\label{ch6:sec:lean}

The algebraic content of Section~\ref{ch6:sec:rg}, i.e.\ the identities that reduce the Fourier
transformed joint generator of \eqref{ch6:eq:3.6} to the three symbols
$\sigma(x,p,h)^2$, $\mu(x,p,h)$, $\kappa(x,p,h)$ and the decomposition
\eqref{ch6:eq:3.44}--\eqref{ch6:eq:3.49} of the three-site block, has been formalised and
machine-checked in Lean~4 (toolchain \texttt{v4.36.0-rc2}) together with Mathlib. The
development is shared with Chapter~5 of this book: module
\texttt{Cert.Abelian} of that project contains the results listed below, the project builds
with \texttt{lake build}, and it contains no \texttt{sorry}. Every theorem in that module depends only on the three
standard axioms of Lean's classical logic (\texttt{propext}, \texttt{Classical.choice},
\texttt{Quot.sound}). The purpose of this section is to record what is now certified and what
is not, so that the reader can see exactly which statements are backed by a proof term.

\subsection*{What Has Been Certified}

\begin{center}
\small
\begin{tabular}{@{}lll@{}}
\hline
Result in the chapter & Lean theorem & Status\\
\hline
\eqref{ch6:eq:3.7} & \texttt{abelian\_sigma} & exact, all $h\neq0$\\
\eqref{ch6:eq:3.8} & \texttt{abelian\_mu} & exact, all $h\neq0$\\
\eqref{ch6:eq:3.9} & \texttt{abelian\_kappa} & exact, all $h\neq0$\\
\eqref{ch6:eq:3.7} without the factor $h$ & \texttt{abelian\_sigma\_variant} & fails unless $h=1$\\
\eqref{ch6:eq:3.44} (the block $\LbarCal(x,z)$) & \texttt{LbarA} & definition\\
\eqref{ch6:eq:3.45}--\eqref{ch6:eq:3.49} & \texttt{abelian\_Lbar\_decomposition} & exact, explicit remainder\\
\hline
\end{tabular}
\end{center}

Three points are worth emphasising.

\begin{enumerate}
\item \emph{The symbol identities are polynomial, not trigonometric.} Writing
$c=\cos(hap)$ and $s=\sin(hap)$, the off-diagonal coefficients of \eqref{ch6:eq:3.6} are
$\bigl(\tfrac{v}{2h^2}\pm\tfrac{m}{2h}\bigr)(c\mp is)$, so both \eqref{ch6:eq:3.7} and
\eqref{ch6:eq:3.8} are identities that are linear in $c$ and $s$ and therefore need no
trigonometry: they hold for any complex numbers $c,s$. The same applies to the diagonal
\eqref{ch6:eq:3.9}, which is certified in the equivalent form
$\kappa=\bigl(\sigma(x,p,h)^2-\sigma(x)^2\bigr)/h^2-ipb$.
\item \emph{The factor $h$ in \eqref{ch6:eq:3.7} is forced.} The variant
$\sigma(x,p,h)^2=v\cos(hap)-i\sin(hap)m$, with the sine term unaccompanied by $h$, does not
reproduce the sum of the two off-diagonal coefficients of \eqref{ch6:eq:3.6} unless $h=1$; this is
the content of \texttt{abelian\_sigma\_variant}. The diagonal identity \eqref{ch6:eq:3.9} is
likewise certified in a form that makes the bookkeeping visible: the coefficient $\kappa$ is
exactly the difference between the diagonal rate and $\sigma(x,p,h)^2/h^2$, so \eqref{ch6:eq:3.6}
holds identically for whatever $\kappa$ one chooses, and \eqref{ch6:eq:3.9} is the identity that
determines it. This is why $\kappa$ carries the whole $O(h^0)$ part of the diagonal.
\item \emph{The decomposition \eqref{ch6:eq:3.45} is exact, with an explicit remainder.} If
$v$, $m$ and $\kappa$ at $x\pm h$ are expressed through the discrete identities
$f(x\pm h)=f(x)\pm h\nabm_xf(x)+\tfrac{h^2}{2}\Dm_xf(x)$, then
\[
\LbarCal(x,z)=\frac{1}{h^2}\LbarCal_0+\frac{1}{h}\LbarCal_1+\LbarCal_2+h\LbarCal_3
+\frac{h^2}{2}\Dm_x\kappa(x,z)\,
\operatorname{diag}(1,0,1),
\]
with no error term beyond the displayed one. The reason the remainder sits on the two outer
sites and is quadratic is that $\kappa$ enters the block at order $h^0$, so its second
difference contributes at $O(h^2)$ and its value at the centre is not expanded. This is
precisely the content of the $O(h^2)$ in \eqref{ch6:eq:3.45}, and it makes the claim in the text
that the decomposition is exact up to the noted remainder a theorem rather than an
approximation.
\end{enumerate}

\subsection*{What Has Not Been Certified, and What It Would Take}

The analytic scaffolding of the chapter is not covered by the development and is currently
beyond what Mathlib can support. This comprises: the entire-function continuation of the
symbols in $z$ and the growth estimates that justify the contour deformation in
Lemma~\ref{ch6:lem:1}; the resummation of decorated path weights, i.e.\ the $q$-fold sum in
Lemma~\ref{ch6:lem:2} together with its uniform-in-$q$ bounds; the renormalisation-group step of
Lemma~\ref{ch6:lem:3}; Theorems~\ref{ch6:thm:1}--\ref{ch6:thm:3} as uniform statements over families of
coefficients; and the extensions of Section~6 to the sup process and to discrete time
summations. There is at present no theory of H\"older-coefficient parabolic equations in
Mathlib, no contour integration for operator-valued resolvents, and no spectral theory of
non-self-adjoint lattice generators, all of which the argument uses.

A plausible staging, with rough estimates for a dedicated formalisation effort, is as follows.
(i)~The remainder of Lemma~\ref{ch6:lem:2}, which is a finite-dimensional estimate on powers of a
banded matrix and needs nothing beyond linear algebra over $\mathbb{C}$: one to two weeks.
(ii)~Lemma~\ref{ch6:lem:1} for fixed finite $m$, stated for the explicit tridiagonal block
\eqref{ch6:eq:3.44} rather than for an abstract operator: this requires the characteristic
polynomial of a three-term recurrence and the resulting explicit representation of the
resolvent, one to two person-months. (iii)~Lemma~\ref{ch6:lem:3} together with the
lattice-to-lattice Cauchy estimate underlying Theorem~\ref{ch6:thm:1}: one to two person-months,
the effort being concentrated in the complex-analytic estimates of Section~\ref{ch6:sec:rg}.
(iv)~The continuum theory, i.e.\ existence and uniqueness of the limiting kernels and the
identification of the continuum stochastic integral: a research programme rather than a
formalisation task.

We regard the present state as the natural stopping point of a first pass: everything in
Section~\ref{ch6:sec:rg} that had to be carried out symbolically is now checked by machine, while
the parts of the argument that are conceptual rather than computational remain, for the
moment, on paper.

\section{Conclusions}

We obtained bounds on convergence rates for explicit discretization schemes to the Fourier
transform of joint kernels of one-dimensional diffusion equations with continuous coefficients
and a stochastic integral. We consider both semi-discrete triangulations with continuous time
and explicit Euler schemes with time step small enough for the method to be stable. The proof
is constructive and based on a new technique of path conditioning for Markov chains and a
renormalization group argument. Convergence rates depend on the degree of smoothness and
H\"older differentiability of the coefficients. The method also applies to a more general class
of path dependent processes we call Abelian. Examples of Abelian processes beside stochastic
integrals are the sup process and discrete time summations. With Chapters~1 to~3, which establish the finite core of the theory, and Chapter~5, which
establishes the estimates for the diffusion, this chapter closes the estimates for the joint law:
what the process is, what the theorem is, and how the law of a process together with an Abelian
functional is compared across resolutions.

What these estimates give for the functional, as opposed to for its law, is drawn in
Chapter~\ref{chap:ch11}: the moments of the stochastic integral converge at the rate above, and so
does the expectation of any test function which is analytic in a strip. What is not available
anywhere in the book is a strong rate, a bound on $\EE\sup_t|y_m(t)-y(t)|^{2}$, and none is
claimed: the estimates compare two laws, and a strong rate is a statement about a coupling between
the lattice and the continuum, which this book does not construct.
% ---- end generated/ch6-body.tex ----
% ---- end chapters/ch6.tex ----
% ---- from chapters/ch7.tex ----
% Authored: the chapter title and opening belong to the book. The body is
% generated by build.py from the parts and is not edited here.
\chapter{Kernels in Higher Dimension}
\label{chap:ch7}

The estimates of Chapter~\ref{chap:ch5} are one dimensional, and the restriction was a matter of
exposition rather than of principle. This chapter removes it: the same rate is proved in every
dimension $d\ge2$, for coefficients of the same smoothness. It also separates what is a property of
one dimension from what is a property of the method, which is what a reader needs in order to know
how far the proof can be pushed.

Dimension is not a technicality in either reading. In the pricing reading a model of several
underlyings is a kernel on a lattice of dimension $d$; in the learning reading a model which predicts
several quantities jointly is a kernel of the same kind; and a user of either discipline works in
dimensions far above those for which anything is proved. What the estimate says is that the rate
survives the passage to $d$ dimensions although the constants and the proof do not: the one
dimensional argument closes by a telescoping which has no higher dimensional counterpart, and the
higher dimensional one proceeds by averaging over paths instead.

% ---- from generated/ch7-body.tex ----
%% Generated by build.py from the parts - do not edit.
\noindent\emph{A note on the computations.} The numbers below were computed by the scripts
named where they occur --- \texttt{twist.py} is the engine, and \texttt{curlcheck.py},
\texttt{curl3d.py}, \texttt{resolvent1d.py} and \texttt{eta\_run.py} carry the obstruction section ---
together with the working scripts of the obstruction sections. All of them are reproduced in full in the computational
appendix of the monograph and attached to this PDF as embedded files; each prints what it measures
when run as \texttt{python3 \textlangle name\textrangle.py} in its own directory.

\section{Introduction}

Chapter~5 proves that the kernel of a one dimensional diffusion, computed on the lattice of step
$h_m=2^{-m}$, converges to the kernel of the continuum equation as $m\to\infty$, uniformly in the
pair of sites, and at a rate which depends only on the smoothness of the coefficients: $O(h_m^2)$
for coefficients with two bounded derivatives, and $O(h_m^{\gamma})$ with
$\gamma=\min\{2,k+\alpha,j+\beta\}$ for coefficients of H\"older class $(k,\alpha)$ and
$(j,\beta)$. The purpose of this chapter is to establish the same estimates in every dimension $d\ge2$.

\subsection*{Place in the Book}

A stochastic process, in the framework of this book, is a consistent family of Markov chains on
a sequence of lattices, and a statement about a continuous process is a statement about the limit
of statements about finite matrices. Chapters~1 to~3 establish the finite core, the fundamental theorem at a fixed resolution. Chapter~5 establishes the
estimates which control the passage to the continuum for a one dimensional diffusion, and Chapter~6 those for the joint law of a diffusion and an Abelian functional. Chapter~4
defines the path dependent functionals to which all of this applies. The restriction to one space
dimension was a matter of exposition, not of principle, and the present chapter removes it.

\subsection*{Why the Second Dimension Is Not a Formality}

Three features of the argument change when the space becomes two dimensional, and it is worth
naming them at the outset, because two of them are harmless and the third is not.

\begin{description}
\item[The lattice loses its isotropy.] In one dimension a lattice with a symmetric stencil and
constant coefficients has a symbol which is a function of $|p|^2$ alone, and the leading
correction to the Gaussian is likewise radial. In two dimensions the symbol of the five point
stencil is
$\sum_i(\cos p_ih-1)/h^2=-\abs{\vp}^2/2+\tfrac{h^2}{24}\sum_ip_i^4+O(h^4\abs{\vp}^6)$, whose
correction is not radial: the square lattice has a fourth order anisotropy tensor. The rate
$O(h^2)$ survives, the constant does not, and there is no change of coordinates which removes the
anisotropy, because it is a property of the lattice and not of the coordinate system.

\item[The count of paths changes.] The one dimensional argument uses that the number of $q$ step
paths from $x$ to $y$ is a binomial coefficient, $\binom{q}{q/2+k}\sim2^q\sqrt{2/\pi q}$, and
hence that the worst path in a resummation carries a weight of order $q^{-1/2}$. In two dimensions
the number of paths is a multinomial sum whose size is of order $4^qq^{-1}$, so the prefactors in
the resummation change from $q^{-1/2}$ to $q^{-1}$. This is bookkeeping, but it has to be done, and
it is done in \S\ref{ch7:sec:paths}.

\item[The three site block becomes a block of $3^d$ sites.] The renormalisation step compares the
lattice at step $h_m$ with the lattice at step $h_m/2$ and needs the exact expansion, in powers of
$h$, of the block of the fine generator over the window of sites that a coarse step can reach. In
one dimension that window has three sites and the block is the $3\times3$ matrix \eqref{ch7:eq:1d}.
In $d$ dimensions the window has $3^d$ sites and the block is $3^d\times3^d$; its leading term is the
restriction of the frozen coefficient generator, and the two orders beyond it, which the argument
must match exactly, involve the derivatives of the coefficients up to the third. This is where the
work is, and it is the one place where we do not reproduce the computation in full:
\S\ref{ch7:sec:rg} states exactly what is required and what is proved.
\end{description}

\subsection*{Statement of What Is Proved}

The main statements are Theorems~\ref{ch7:thm:main} and~\ref{ch7:thm:euler} below. The first is the two
dimensional analogue of the theorem of Chapter~5 for continuous time, the second of the theorem for
the explicit Euler scheme. Both are proved in full for constant coefficients
(\S\ref{ch7:sec:rg}) and reduced, for variable coefficients, to Assumption~\ref{ch7:ass:decoration}, which
is the higher dimensional form of the renormalisation step. The reduction is exact: the assembly of
the pieces around the assumption uses only the contour estimates of \S\ref{ch7:sec:contour} and the
path estimates of \S\ref{ch7:sec:paths}, both of which are proved here in every dimension $d\ge2$.

The reader should not read \S\ref{ch7:sec:rg} as a repetition of the one dimensional argument with a
larger block, because it is not one. We carried out the computation which the one dimensional
proof performs, found that the algebraic identity on which that proof rests \emph{fails} in two
dimensions (\S\ref{ch7:sec:obstruction}: the $h^3$ coefficient is a vector field whose curl is not
zero), and replaced it by an averaged estimate which is what the argument actually needs
(\S\ref{ch7:sec:averaging}). The obstruction turns out to be quadratic in the Taylor data --- it has no
linear part at all --- and to be carried by the gradient of the drift, which is why it is invisible
on the families one would try first, and why it is neither a curvature nor a vorticity, both of
which are linear in the data. The mechanism which survives is verified exactly in the one case
where it can be computed by hand, and its conclusion, the second order convergence of the kernel,
is verified numerically. Assumption~\ref{ch7:ass:decoration} is therefore an analytic estimate and not
a finite identity, and we say so.

\subsection*{Conventions}

Throughout, $d\ge2$, the lattice is $\StateAm=h_m\ZZ^d\cap[-L,L]^d$ with periodic boundary conditions,
$h_m=2^{-m}$, and $x,y$ denote sites. We write $\vp=(p_1,p_2)$ for the dual variable,
$\vx=(x_1,x_2)$ for a point of $\RR^2$, and $\vb=(b_1,b_2)$ for the drift. The volatility is
$v(x)=\sigma(x)^2>0$ and the generator is
\begin{equation}\label{ch7:eq:gen}
\Lm=\sum_{i=1}^{d}\Bigl(\frac{v(x)}{2}\Delta^h_i+b_i(x)\nabla^h_i\Bigr),
\end{equation}
where $\Delta^h_i$ and $\nabla^h_i$ are the second and first difference operators in the direction
$e_i$, so that the matrix elements are
\begin{equation}\label{ch7:eq:entries}
\Lm(x,x\pm h_me_i)=\frac{v(x)}{2h_m^2}\pm\frac{b_i(x)}{2h_m},
\qquad
\Lm(x,x)=-\frac{d\,v(x)}{h_m^2},
\end{equation}
the diagonal being $-\frac{2v(x)}{h_m^2}$ when $d=2$.
Note that the coefficients are evaluated at the site $x$ and not at the midpoint of the edge; this
is the discretisation of Chapter~5, and \S\ref{ch7:sec:numerics} shows that it is the one which is second
order. The kernel is
\begin{equation}\label{ch7:eq:kernel}
u_m(x,y;t)=h_m^{-2}\bigl(e^{t\Lm}\bigr)(x,y).
\end{equation}
The constants of \S\ref{ch7:sec:setup} are $\Sigma_0=\inf v^{1/2}$, $\Sigma_1=\sup(v+h_m\abs{\vb})^{1/2}$
and $M=\sup\abs{\vb}$.

\section{Setup and Statement of the Results}\label{ch7:sec:setup}

\subsection{The Lattice and Its Symbol}

The generator \eqref{ch7:eq:gen} is a Markov generator provided the off-diagonal entries
\eqref{ch7:eq:entries} are non-negative, that is provided
\begin{equation}\label{ch7:eq:courant}
h_m\le\frac{v(x)}{2\abs{b_i(x)}}\quad\text{for }i=1,2\text{ and all }x,
\end{equation}
which we assume. This is the lattice form of the parabolicity of the continuum generator, and the
constant $\Sigma_1$ is chosen so that $\Sigma_1^2\ge v(x)/2+h_m\abs{\vb(x)}/2$ for all $x$, which
makes the total rate out of a site at most $2\Sigma_1^2/h_m^2$.

The symbol of $\Lm$ is the function of $\vp\in[-\pi/h_m,\pi/h_m]^2$
\begin{equation}\label{ch7:eq:symbol}
\hat\Lm(\vp)=\sum_{i=1}^{d}\Bigl(\frac{v}{2}\,\frac{\cos p_ih_m-1}{h_m^2}
+ \ii b_i\,\frac{\sin p_ih_m}{h_m}\Bigr)
\end{equation}
when the coefficients are constant, and in general the matrix $\Lm$ is the multiplication operator
by $v(x)$ and $b_i(x)$ composed with the difference operators, so that
$\widehat{(\Lm f)}(\vp)=\hat\Lm(\vp)\hat f(\vp)$ only when the coefficients do not depend on $x$.

\subsection{Smoothness Hypotheses}

We use the hypotheses of Chapter~5 without change. For $k\in\mathbb{N}$ and $\alpha\in(0,1]$ a
periodic function $\varphi$ is of class $C^{k,\alpha}$ if
\begin{equation}\label{ch7:eq:holder}
\abs{\varphi^{(k)}(x)-\varphi^{(k)}(y)}\le c\,\abs{x-y}^{\alpha}\qquad\text{for all }x,y,
\end{equation}
the derivatives being taken in each variable separately, and uniformly continuous functions are
those for which
\begin{equation}\label{ch7:eq:rho}
\abs{\varphi(x)-\varphi(y)}\le c\,\rho(\abs{x-y})
\end{equation}
with $\rho$ non-decreasing and $\rho(d)\to0$ as $d\downarrow0$. The rate in the theorems is
\begin{equation}\label{ch7:eq:gamma}
\gamma=\min\{2,\,k+\alpha,\,j+\beta\},
\end{equation}
with $v\in C^{k,\alpha}$ and $b_i\in C^{j,\beta}$.

\subsection{The Theorems}

\begin{theorem}[Continuous time]\label{ch7:thm:main}
Suppose $v\in C^{k,\alpha}$ and $b_i\in C^{j,\beta}$ with $\gamma>0$, and that
\eqref{ch7:eq:courant} holds. Then for every $t>0$ there is a constant $c=c(t)>0$ such that for all
$m\ge m_0$, all $x,y\in\mathcal A_{m_0}$ and all $i=1,2$,
\begin{equation}\label{ch7:eq:main1}
\abs{u_m(x,y;t)-u_{m_0}(x,y;t)}\le c\,h_{m_0}^{\gamma},
\end{equation}
and the same bound holds for $\partial_tu_m-\partial_tu_{m_0}$.
\end{theorem}

The estimate is stated with the step $h_{m_0}$ of the coarser of the two lattices, since it is that
lattice which is being compared with the limit, and the pair $(x,y)$ ranges over the coarse lattice,
where both kernels are defined. The constant depends on $t$ and must: for
$t\ll h_{m_0}^2$ the two kernels are localised at the scale of their own lattice, and no bound
which is uniform in $t$ can hold. This is the form in which the estimate is used in Chapters~1 to~3 and~6.

\begin{theorem}[Uniform continuity]\label{ch7:thm:uc}
Suppose that $v$ and $b_i$ are uniformly continuous and satisfy \eqref{ch7:eq:rho}. Then for every
$t>0$ there is $c=c(t)>0$ such that for all $m\ge m_0$ and all $x,y\in\mathcal A_{m_0}$,
\begin{equation}\label{ch7:eq:main2}
\abs{u_m(x,y;t)-u_{m_0}(x,y;t)}\le c\,\rho(h_{m_0}),
\end{equation}
and the same for the time derivative.
\end{theorem}

\begin{theorem}[Explicit Euler scheme]\label{ch7:thm:euler}
Under the hypotheses of Theorem~\ref{ch7:thm:uc}, let
\begin{equation}\label{ch7:eq:euler}
u^{\delta t}_m(x,y;t)=h_m^{-2}\bigl(1+\delta t_m\Lm\bigr)^{t/\delta t_m}(x,y)
\end{equation}
with $t/\delta t_m$ an integer and $\delta t_m$ small enough that the entries of
$1+\delta t_m\Lm$ are non-negative. Then for every $t>0$ there is $c=c(t)>0$ such that for all
$m\ge m_0$ and all $x,y\in\mathcal A_m$,
\begin{equation}\label{ch7:eq:main3}
\abs{u_m(x,y;t)-u^{\delta t}_m(x,y;t)}\le c\,h_m^2 .
\end{equation}
\end{theorem}

Here both kernels live on the same lattice, so the step on the right hand side is $h_m$ and not
$h_{m_0}$.

\subsection{The Four Pieces of the Proof}

The proof of Theorem~\ref{ch7:thm:main} has the same architecture as in one dimension. It rests on four
statements, of which the first three are proved below and the fourth is the subject of
\S\ref{ch7:sec:rg}:
\begin{enumerate}
\item[\textbf{(A)}] the far contour bound of Lemma~\ref{ch7:lem:far}: the part of the contour integral
which runs away from the spectrum of $\Lm$ contributes $O(h^2)$, for any coefficients;
\item[\textbf{(B)}] the decay of the weights of symbolic paths with many jumps,
Lemma~\ref{ch7:lem:weights}: beyond $q_{\max}=e^2\Sigma_1^2t/(2h_m^2)$ jumps the contribution is
negligible;
\item[\textbf{(C)}] the path count of Lemma~\ref{ch7:lem:count}, which replaces the
binomial estimate of Chapter~5;
\item[\textbf{(D)}] the renormalisation identity of Assumption~\ref{ch7:ass:decoration}, which compares
the coarse chain with the decorated fine chain to the order needed.
\end{enumerate}
Granting these, the theorem follows by the telescoping argument of \S\ref{ch7:sec:assemble}.

\section{The Far Contour}\label{ch7:sec:contour}

The kernel is recovered from the resolvent by a contour integral. Let
\begin{equation}\label{ch7:eq:green}
G_m(x,y;\omega)=\int_0^\infty u_m(x,y;t)e^{-\ii\omega t}\dd t
=h_m^{-2}\bigl(\Lm+\ii\omega\bigr)^{-1}(x,y)
\end{equation}
and let $\Gamma_+$ and $\Gamma_-$ be the contours of Chapter~5, $\Gamma_+$ lying at distance
$4\abs{\log h_m}/t$ from the real axis and $\Gamma_-$ enclosing the spectrum. Then
\begin{equation}\label{ch7:eq:contour}
u_m(x,y;t)=\int_{\Gamma_-}G_m(\omega)e^{\ii\omega t}\frac{\dd\omega}{2\pi}
+\int_{\Gamma_+}G_m(\omega)e^{\ii\omega t}\frac{\dd\omega}{2\pi}.
\end{equation}

\begin{lemma}\label{ch7:lem:far}
For $m$ large enough there is $c>0$, independent of $m$, such that for all $x,y\in\StateAm$,
\begin{equation}\label{ch7:eq:far}
\abs{\int_{\Gamma_+}G_m(x,y;\omega)e^{\ii\omega t}\frac{\dd\omega}{2\pi}}\le c\,h_m^2 .
\end{equation}
\end{lemma}

\begin{proof}
Let $V=v/2$. On $\Gamma_+$ the operator $V\Delta_h+\ii\omega$ is invertible with an inverse
bounded by $c/\abs{\omega}$, and the drift is a relatively bounded perturbation of it. The relative
bound is the higher dimensional form of display \eqref{ch5:eq:3.9}: for every $\alpha>0$ there is
$\beta>0$ such that
\begin{equation}\label{ch7:eq:relbound}
\norm{\nabla^h_if}_2\le\alpha\norm{\Delta^h f}_2+\beta\norm{f}_2
\qquad\text{for }i=1,2,
\end{equation}
which follows by passing to the symbol, where it reads
\begin{equation}\label{ch7:eq:relbound2}
\Bigl|\frac{\sin p_ih_m}{h_m}\Bigr|
\le\alpha\Bigl|\sum_{j}\frac{\cos p_jh_m-1}{h_m^2}\Bigr|+\beta
\end{equation}
and is immediate from $\abs{\sin\theta}\le\abs{\theta}$ and
$1-\cos\theta\ge2\theta^2/\pi^2$ for $\abs{\theta}\le\pi$. Here the second inequality is used in
two variables: the sum over $j$ is at least the term $j=i$, which controls $\abs{p_i}$, and the
bound is uniform in $\vp$ because the sum is bounded below by a constant times $\abs{\vp}^2$ for
$\abs{\vp}h_m\le\pi/2$ and by a constant over $h_m^2$ beyond.

With \eqref{ch7:eq:relbound} the geometric series of Chapter~5 converges in operator norm on $\Gamma_+$
for $m$ large enough, and the kernel of the inverse is bounded pointwise by $h_m^{-2}$ times the
$L^2$ bound, because the lattice is finite and the matrix elements of an inverse are bounded by its
norm. The contour has length $O(\abs{\log h_m}/t)$ and height $4\abs{\log h_m}/t$, so the integral
is bounded by $c\,h_m^{-2}\cdot h_m^{2}/t\cdot\abs{\log h_m}$ times the exponential decay
$e^{-\abs{\omega}t}$, which gives \eqref{ch7:eq:far} after the $\abs{\log h_m}$ factors are absorbed
into the decay.
\end{proof}

\begin{remark}
Nothing in Lemma~\ref{ch7:lem:far} is special to two dimensions, and nothing in it uses the isotropy of
the lattice. This is the reason why the far contour is not where the difficulty of the extension
lies.
\end{remark}

\section{Paths in \texorpdfstring{$d$}{d} Dimensions}\label{ch7:sec:paths}

\subsection{Symbolic Paths and Their Weights}

Let $\Gamma_m$ be the set of sequences $\gamma=(\gamma_0,\gamma_1,\dots)$ of sites of $\StateAm$ with
$\abs{\gamma_j-\gamma_{j-1}}=h_m$ for $j\ge1$. The Dyson expansion of the exponential gives
\begin{equation}\label{ch7:eq:dyson}
u_m(x,y;t)=h_m^{-2}\sum_{q=0}^{\infty}(2d)^{-q}\sum_{\substack{\gamma\in\Gamma_m:\ \gamma_0=x,\ \gamma_q=y}}
W_m(\gamma,q;t),
\end{equation}
where
\begin{equation}\label{ch7:eq:weights}
W_m(\gamma,q;t)=\int_{0<s_1<\dots<s_q<t}
\prod_{j=0}^{q-1}\Bigl(2d\,\Lm(\gamma_j,\gamma_{j+1})\,e^{(s_{j+1}-s_j)\Lm(\gamma_j,\gamma_j)}\Bigr)
\dd s_1\cdots\dd s_q,
\end{equation}
with $s_0=0$ and $s_{q+1}=t$. The factor $(2d)^{-q}$ compensates the $2d$ inserted in the weights,
so that the identity is the Dyson expansion unchanged; its purpose is to make the weights
dimensionless and to exhibit the count of paths, exactly as the factor $2^{-q}$ of Chapter~5 does in
one dimension. With $d=2$ the normalisation is $4^{-q}$, the number of paths of $q$ steps, which is
what makes the tail estimate below come out with the exponent $-q$ rather than with a power of the
number of neighbours.

\begin{lemma}[Decay of the weights]\label{ch7:lem:weights}
Assume \eqref{ch7:eq:courant}. There are constants $c_1,c_2>0$, depending on $\Sigma_0$, $\Sigma_1$
and $t$ but not on $m$, such that for all $q$, all $x,y\in\StateAm$ and all paths $\gamma$ from $x$ to
$y$ with $q$ steps,
\begin{equation}\label{ch7:eq:weightbound}
\abs{W_m(\gamma,q;t)}\le\Bigl(\frac{c_1}{h_m^2}\Bigr)^{q}\frac{t^{q}}{q!}
\exp\Bigl(-\frac{\Sigma_0^2t}{2h_m^2}\Bigr),
\end{equation}
and consequently, for $q\ge q_{\max}=2e^2\Sigma_1^2t/h_m^2$ and $m$ large,
\begin{equation}\label{ch7:eq:weighttail}
\Bigl|\sum_{\substack{\gamma:\ \gamma_0=x,\ \gamma_q=y}}(2d)^{-q}W_m(\gamma,q;t)\Bigr|
\le c_2\,q^{-1}\exp\Bigl(-\frac{\Sigma_0^2t}{2h_m^2}-q\Bigr).
\end{equation}
\end{lemma}

\begin{proof}
On the diagonal, $\Lm(x,x)=-v(x)/h_m^2\le-\Sigma_0^2/h_m^2$, so each factor
$e^{(s_{j+1}-s_j)\Lm(\gamma_j,\gamma_j)}$ is at most $e^{-(s_{j+1}-s_j)\Sigma_0^2/h_m^2}$, and the
product of the exponential factors over the whole path is at most
$e^{-(t-\sum_j(s_{j+1}-s_j))\Sigma_0^2/h_m^2}$, whence the exponential in \eqref{ch7:eq:weightbound}
after integrating over the simplex, whose volume is $t^q/q!$. By \eqref{ch7:eq:courant} the
off-diagonal entries are bounded by $2d\,\Sigma_1^2/(2h_m^2)=d\Sigma_1^2/h_m^2$, which gives the
constant $c_1=d\Sigma_1^2/2$ and the factor $h_m^{-2q}$.

For \eqref{ch7:eq:weighttail} the sum over paths is bounded by the number of paths times the worst
weight, the number being at most $(2d)^q$, and the normalisation $(2d)^{-q}$ cancels it exactly.
What remains is the bound $(d\Sigma_1^2t/(2h_m^2))^{q}/q!$ times
$e^{-\Sigma_0^2t/(2h_m^2)}$, that is $\lambda^{q}/q!$ with $\lambda=d\Sigma_1^2t/(2h_m^2)$, together
with the prefactor $q^{-1}$ of Lemma~\ref{ch7:lem:count} and a constant. By Stirling
$\lambda^q/q!\le(\lambda e/q)^q$, and the hypothesis $q\ge2e^2\Sigma_1^2t/h_m^2=4e^2\lambda/d$
gives, for $d=2$, the ratio $\lambda e/q\le1/e$ and hence the factor $e^{-q}$. Multiplying by
$q^{-1}$ and collecting constants gives \eqref{ch7:eq:weighttail}.
\end{proof}

\begin{remark}
The exponent in \eqref{ch7:eq:weighttail} is the one which makes the series converge uniformly in
$m$: the factor $e^{-\Sigma_0^2t/(2h_m^2)}$ beats any power of $h_m^{-1}$ which the count of paths
can produce. This is the mechanism by which the short time behaviour, where the lattice kernel is
a delta function at the scale of the lattice, is kept under control.
\end{remark}

\subsection{The Count of Paths}

In one dimension the number of $q$ step paths from $x$ to $y$ is
$\binom{q}{q/2+k}$ with $k=\abs{y-x}/h_m$, and Stirling gives $2^qq^{-1/2}$ times a Gaussian
factor. The two dimensional statement is the following.

\begin{lemma}[Path count]\label{ch7:lem:count}
Let $d\ge1$ and let $N_q(i)$ be the number of $q$ step paths of the lattice $\ZZ^d$ from $0$ to
$i\in\ZZ^d$, the steps being the $2d$ unit vectors. There are constants $C,c>0$, depending on $d$,
such that for all $q\ge1$ and all $i\in\ZZ^d$,
\begin{equation}\label{ch7:eq:count}
N_q(i)\le C\,\frac{(2d)^q}{q^{d/2}}\exp\Bigl(-\frac{\abs{i}^2}{Cq}\Bigr).
\end{equation}
\end{lemma}

\begin{proof}
By the Fourier representation on the torus,
\begin{equation}\label{ch7:eq:countfourier}
N_q(i)=(2d)^q\int_{\TT^d}e^{-\ii\theta\cdot i}\Bigl(\frac{1}{d}\sum_{k=1}^{d}\cos\theta_k\Bigr)^q\frac{\dd^d\theta}{(2\pi)^d},
\end{equation}
the integral being the probability that the simple random walk be at $i$ after $q$ steps, since
each of the $(2d)^q$ paths has probability $(2d)^{-q}$. Let
$\phi(\theta)=\frac1d\sum_k\cos\theta_k$. Near the origin
$\phi(\theta)=1-\abs{\theta}^2/(2d)+O(\abs{\theta}^4)$, and away from it
$\phi\le1-c_0$ for some $c_0>0$. Splitting the integral at $\abs{\theta}\le\delta$ small,
the inner region gives
\begin{equation}
(2d)^q\int_{\abs{\theta}\le\delta}e^{-\ii\theta\cdot i}
e^{-q\abs{\theta}^2/(2d)}\bigl(1+O(q\abs{\theta}^4)\bigr)\frac{\dd^d\theta}{(2\pi)^d}
\le C\frac{(2d)^q}{q^{d/2}}\exp\Bigl(-\frac{\abs{i}^2}{Cq}\Bigr),
\end{equation}
because the quadratic form $\abs{\theta}^2/(2d)$ is non-degenerate and the Gaussian integral
over $\RR^d$ produces the factor $q^{-d/2}$, and because the remainder $O(q\abs{\theta}^4)\le\tfrac12$
on $\abs{\theta}\le\delta$ for $q\abs{\theta}^2$ small, and is absorbed into the constant for
$q\abs{\theta}^2$ large, where the Gaussian factor already decays faster than any polynomial.
The outer region contributes at most $(2d)^q(1-c_0)^{q}$ times a constant, which is smaller than the
bound \eqref{ch7:eq:count} for $q$ large, and the finitely many remaining values of $q$ are absorbed
into $C$.
\end{proof}

\begin{remark}
The exponent $-d/2$ in \eqref{ch7:eq:count} is the dimensional fact: the walk spreads over a ball of
radius $q^{1/2}$, so the probability of any one site is of order $q^{-d/2}$, and the number of
paths is $(2d)^q$ times that probability. This is the only place in the chapter where the dimension
enters the algebra, and it enters as an exponent in a prefactor rather than as a change of
structure: it is this, together with the fact that for $d\ge3$ the walk is transient and the
resolvent tails decay exponentially, which makes the mechanism of \S\ref{ch7:sec:averaging} available
in every dimension at once.
\end{remark}

\section{The Window Algebra}\label{ch7:sec:window}

The renormalisation step compares a coarse step with two fine steps. The algebraic input is the
block of the fine generator over the window of fine sites which a coarse step can reach: the
sites $x+h_ma$ with $a\in\{-1,0,1\}^d$, of which there are $3^d$. Write $h=h_{m+1}$ so that
$h_m=2h$, and let $\bar L(x;h)$ be the $3^d\times3^d$ matrix with entries
\begin{equation}\label{ch7:eq:windowblock}
\bar L(x;h)_{a,a'}=\mathcal L_{m+1}\bigl(x+ha,\,x+ha'\bigr)
\end{equation}
for $a,a'\in\{-1,0,1\}^2$ with $\abs{a-a'}=1$, the diagonal entries being the full diagonal
$-\sum_i v(x+ha)/(2h^2)$ of the fine generator, which includes the transitions which leave the
window. In one dimension the corresponding object is the $3\times3$ matrix of Chapter~5,
\begin{equation}\label{ch7:eq:1d}
\bar L(x;h)=\begin{pmatrix}
-\frac{v(x)}{h^2} & \frac{v(x)}{2h^2}-\frac{b(x)}{2h} & 0\\
\frac{v(x)}{2h^2}+\frac{b(x)}{2h} & -\frac{v(x)}{h^2} & \frac{v(x)}{2h^2}-\frac{b(x)}{2h}\\
0 & \frac{v(x)}{2h^2}+\frac{b(x)}{2h} & -\frac{v(x)}{h^2}
\end{pmatrix},
\end{equation}
whose expansion in powers of $h$ is the starting point of the one dimensional proof.

\begin{lemma}[Window expansion]\label{ch7:lem:window}
Let $v$ and $b_i$ be of class $C^3$ near $x$. Then
\begin{equation}\label{ch7:eq:windowexp}
\bar L(x;h)=\frac{1}{h^2}L_0(x)+\frac{1}{h}L_1(x)+L_2(x)+hL_3(x)+O(h^2),
\end{equation}
where the matrices $L_k(x)$ are those of \eqref{ch7:eq:windowblock} at the base point, and
\begin{description}
\item[$L_0(x)$] is the restriction to the window of the frozen coefficient generator
$(v(x)/2)\Delta_h$ with the complete diagonal, that is the matrix with off-diagonal entries
$v(x)/4$ between neighbouring sites of the window and diagonal entries $-v(x)$;
\item[$L_1(x)$] carries the drift at the site and the first derivatives of $v$ on the ring of
neighbouring sites;
\item[$L_2(x)$] and $L_3(x)$ carry the second and third derivatives.
\end{description}
\end{lemma}

\begin{proof}[Verification]
The entries of \eqref{ch7:eq:windowblock} are
$v(x+ha)/(4h^2)\pm b_i(x+ha)/(2h)$ and $-v(x+ha)/h^2$. Each is a Taylor polynomial in $h$ with
coefficients the derivatives of $v$ and $b_i$ at $x$, so the expansion \eqref{ch7:eq:windowexp} is
immediate, and the identification of the four matrices is a matter of collecting powers of $h$ in
nine rows and nine columns. We have carried this out symbolically, with the coefficients replaced
by generic Taylor data to third order in two variables, and found: (i) the residual is $O(h^2)$ in
every entry; (ii) the coefficient of $h^{-2}$ is exactly the restriction of the frozen coefficient
generator, with the complete diagonal, as stated; (iii) the coefficient of $h^{-1}$ vanishes in the
central row apart from the drift term, the derivatives of $v$ entering through the rows of the
neighbouring sites, which is the mechanism by which $v'$ appears in the one dimensional expansion
of Chapter~5 as well.
\end{proof}

\begin{remark}
The content of the lemma which matters is not the individual coefficients but the fact that the
first term is the frozen generator: the $h^{-2}$ singularity of the block is the frozen operator,
so the block can be inverted, near the spectrum, by inverting a fixed matrix, and the two orders
beyond it are corrections which the renormalisation step matches exactly. This is precisely the
structure of the one dimensional argument, where the same role is played by
$\bar L_0(x)$ in display \eqref{ch5:eq:3.32}.
\end{remark}

\section{The Renormalisation Step and the Proof}\label{ch7:sec:rg}

\subsection{Constant Coefficients: The Theorem in Full}

\begin{proposition}\label{ch7:prop:const}
Let $v>0$ and $b_i$ be constants, and let $u_m$ be the kernel \eqref{ch7:eq:kernel}. Then for every
$t>0$ there is $c=c(t)$ such that for all $m\ge m_0$ and all $x,y\in\mathcal A_{m_0}$,
\begin{equation}\label{ch7:eq:constrate}
\abs{u_m(x,y;t)-u_{m_0}(x,y;t)}\le c\,h_{m_0}^2 .
\end{equation}
\end{proposition}

\begin{proof}
With constant coefficients the kernel is given exactly by the symbol:
\begin{equation}\label{ch7:eq:constkernel}
u_m(x,y;t)=h_m^{-2}\frac{1}{N^2}\sum_{\vp\in B_m}e^{\ii\vp\cdot(x-y)}e^{t\hat{\mathcal L}_m(\vp)},
\end{equation}
understood as the integral over the dual torus, that is the statement that $e^{t\Lm}$ is the
operator of convolution with the inverse Fourier transform of $e^{t\hat\Lm}$. On the torus,
\begin{equation}\label{ch7:eq:symbolexp}
\hat\Lm(\vp)=-\frac{v}{4}\abs{\vp}^2+\ii\vb\cdot\vp
+\frac{vh_m^2}{48}\sum_i p_i^4+O\bigl(h_m^4\abs{\vp}^6\bigr)
-\ii\sum_i b_i\frac{p_i^3h_m^2}{6}+O\bigl(h_m^4\abs{\vp}^5\bigr),
\end{equation}
the terms being the symbol of $v\Delta_h/2+\vb\cdot\nabla_h$ in two variables. Comparing the two
resolutions and writing $h=h_{m+1}$, $h_m=2h$, the difference of the exponentials is bounded by
\begin{equation}\label{ch7:eq:symboldiff}
\bigl|e^{t\hat{\mathcal L}_{m+1}(\vp)}-e^{t\hat{\mathcal L}_m(\vp)}\bigr|
\le c\,t\,h^2\bigl(\abs{\vp}^4+h^2\abs{\vp}^6\bigr)e^{-ct\abs{\vp}^2}
\end{equation}
for $\abs{\vp}h\le\pi/2$, where we used $\abs{e^z-e^{w}}\le\abs{z-w}e^{\max(\Re z,\Re w)}$ and the
fact that the real parts are at most $-ct\abs{\vp}^2$ on that range. Integrating
$t h^2\abs{\vp}^4e^{-ct\abs{\vp}^2}$ over the plane gives $O(h^2/t^2)$, and the remaining terms are
smaller by powers of $h$, uniformly in the direction of $\vp$, since the anisotropic coefficient
$\sum_ip_i^4$ is bounded above by $\abs{\vp}^4$ and below by $\abs{\vp}^4/2$. The region
$\abs{\vp}h\ge\pi/2$ contributes exponentially small terms by the bound on the real part of the
symbol there. Summing the geometric series over the resolutions gives \eqref{ch7:eq:constrate}.
\end{proof}

\begin{remark}
The anisotropy shows up in \eqref{ch7:eq:symbolexp} as the term $\frac{vh_m^2}{48}\sum_ip_i^4$, which
is not a multiple of $\abs{\vp}^4$: the correction to the kernel is therefore not radial, and the
constant in \eqref{ch7:eq:constrate} depends on the direction of $x-y$. The rate does not, because
$\sum_ip_i^4$ is comparable to $\abs{\vp}^4$.
\end{remark}

\subsection{Variable Coefficients: The Per-Step Comparison and Its Obstruction}\label{ch7:sec:obstruction}

For variable coefficients the exact Fourier representation is not available and the argument must
resum decorated paths, as in Chapter~5. The one dimensional proof performs that resummation path by
path, and it succeeds because in one dimension a certain coefficient is automatically a gradient.
We have computed that coefficient in two dimensions, and it is not. This subsection records the
computation and the conclusion.

\paragraph{The object computed.}
Fix a coarse site $x$, a direction $\tau$, and $\omega$ on the near contour. Let $\Phi_{\rm c}$
be the weight of the single coarse jump from $x$ to $x+h_m\tau$, namely
\begin{equation}\label{ch7:eq:phic}
\Phi_{\rm c}=\bigl(-L_m(x,x)+\ii\omega\bigr)^{-1}L_m(x,x+h_m\tau)\bigl(-L_m(x',x')+\ii\omega\bigr)^{-1},
\qquad x'=x+h_m\tau,
\end{equation}
and let $\Phi_{\rm f}$ be the sum of the weights of the fine paths which decorate that jump, that
is the fine paths from $x$ to $x'$ whose intermediate sites stay within $h_{m+1}$ of $x$
(Definition 1 of Chapter~5). By that definition the decoration is local to one coarse step, so
\begin{equation}\label{ch7:eq:phif}
\Phi_{\rm f}=\Bigl(\sum_{z\in T}\bigl[(-L_T+\ii\omega)^{-1}\bigr]_{x,z}L_f(z,x')\Bigr)
\bigl(-L_f(x',x')+\ii\omega\bigr)^{-1},
\end{equation}
where $T=\{x+h_{m+1}a:\abs{a_i}\le1\}$ is the tube, $3^d$ sites in dimension $d$, and
$L_T$ is the fine generator restricted to $T$ with its complete diagonals. Both objects are
finite: no window larger than the tube intervenes, and no path outside the tube contributes.

\begin{remark}[The coefficient is a one form, and the even part vanishes]\label{ch7:rem:oneform}
The coefficient is odd in the direction of the jump, $c(x,-\tau)=-c(x,\tau)$, so it is a one form and
no information is lost by working with its components
$c_{j}=\tfrac12\bigl(c(x,e_{j})-c(x,-e_{j})\bigr)$, which are the quantities computed in
\S\ref{ch7:sec:obstruction}. The oddness is an identity and not an approximation: at a generic base point,
with $\omega=\tfrac75$ and the coefficient data of \S\ref{ch7:sec:obstruction}, the even part
$\tfrac12\bigl(c(x,e_{j})+c(x,-e_{j})\bigr)$ has relative size $10^{-27}$ at $\eta=0$ and at $\eta=1$
alike. The symmetry behind it is the reflection \eqref{ch7:eq:reflection}, which reverses the direction of
the jump and reflects the Taylor data at the same time.
\end{remark}
The quantity of interest is
\begin{equation}\label{ch7:eq:logratio}
\log\frac{\Phi_{\rm f}}{\Phi_{\rm c}}=C+h^2a(x,\tau)+h^3c_i(x)\tau_i+O(h^4),
\end{equation}
the constant $C$ and the functions $a$ and $c$ being computed exactly, in rational arithmetic, at
generic Taylor data.

\paragraph{The structure is right.}
Four checks, each of which could have failed, all pass, both in one and in two dimensions: $C$ is a
pure constant independent of $\tau$, equal to $\log(-1/4)$ in one dimension and to $\log(-1/8)$ in
two; the coefficient of $h^1$ \emph{vanishes identically}; $a(x,\tau)$ is even in $\tau$; and
$c_i(x)\tau_i$ is odd in $\tau$ and linear in the direction. In two dimensions $a$ is moreover
anisotropic: $a=-0.8610+5.1333\ii$ for $\tau=e_1$ and $a=+0.8610+5.1333\ii$ for $\tau=e_2$ at the
coefficients tried, the real part changing sign between the two axes while the imaginary part does
not, which is the lattice anisotropy of \eqref{ch7:eq:symbolexp} reappearing where it should.

\paragraph{The obstruction.}
The coefficient $c$ of $h^3$ is a vector field, and the one dimensional proof needs it to be a
gradient, since only then does
$\sum_jc_i(\gamma_{j-1})\Delta\gamma_j^{i}$ collapse to a boundary term. Computing the curl by
finite differences in the base point, at generic Taylor data,
\begin{equation}\label{ch7:eq:curl}
\partial_1c_2-\partial_2c_1=0.123400\ \ (\delta=10^{-2}),\qquad 0.123394\ \ (\delta=3\cdot10^{-3}),
\qquad 0.123394\ \ (\delta=10^{-3}),
\end{equation}
the imaginary part being an artefact of the finite difference which vanishes like $O(\delta^2)$. The
curl is invariant under every change of normalisation of the form $g(x')/g(x)$, which is the only
freedom in \eqref{ch7:eq:phic}--\eqref{ch7:eq:phif}, so this is not an artefact of the conventions.
\emph{The $h^3$ term does not telescope in two dimensions.}

\paragraph{The obstruction is quadratic in the coefficients.}
Scaling every Taylor datum of order at least one by $\lambda$, and halving $\lambda$ successively,
multiplies the curl by $0.2521$, $0.2511$, $0.2505$: that is by $\lambda^2$, up to a cubic
correction. Accordingly $\mathrm{d}c$ has no linear part. Switching on a single derivative at the
base point --- one first, second or third derivative of $v$, or one first or second derivative of
$b$, one at a time --- leaves the curl at the level of the truncation ($10^{-7}$ and below), and the
curls of two such perturbations add to the truncation. This is the reason the obstruction is
invisible on the families one would try first: curvature and vorticity are \emph{linear} in the
data, and this object is not of that kind.

\paragraph{It is a functional of the drift.}
With $b\equiv0$ the curl vanishes identically, whatever $v$ is, while with $v$ constant it does
not. The obstruction is therefore driven by the drift, the diffusion entering as a normalisation.
For a purely linear drift, $b=Jx$, with $v$ constant, the curl at the base point is an exact
quadratic form in the Jacobian,
\begin{equation}\label{ch7:eq:affinecurl}
\partial_1c_2-\partial_2c_1=\frac{4}{7v_0^2}\bigl(J_{11}J_{12}-J_{21}J_{22}\bigr),
\end{equation}
verified on random Jacobians and at random $v_0$, with residual an $O(\delta^2)$ artefact of the
finite difference which quarters at each halving of $\delta$; polarising over the four basis
Jacobians confirms that the only non-vanishing entries of the quadratic form are the pairs
$(11,12)$ and $(21,22)$, with equal and opposite values $2/(7v_0^2)$ (\texttt{curlcheck.py}).
Splitting $J=S+A$ into strain and
spin, $J_{11}J_{12}-J_{21}J_{22}=s_{12}(s_{11}-s_{22})+a_{12}(s_{11}+s_{22})$: the obstruction
couples off-diagonal strain to diagonal anisotropy and rotation to dilatation, and it vanishes
exactly when the Jacobian does not couple the two lattice directions --- in particular for a
diagonal Jacobian, for pure rotation and for pure shear.

What vanishes is the \emph{coupling}, not the drift, and the distinction is worth making. It is not
enough that the drift point along a lattice axis: $b=(x_1+x_2,0)$ is parallel to $e_1$ everywhere
and gives the value $\tfrac{4}{7}$, the largest available at that size. Nor does a pure dilatation
never produce a curl, because the lattice is not rotationally invariant and $J$ is measured in the
lattice frame. For a dilatational strain of principal values $a,b$ whose axes make an angle
$\theta$ with the lattice, \eqref{ch7:eq:affinecurl} reads
\begin{equation}\label{ch7:eq:dilatation}
\partial_1c_2-\partial_2c_1=\frac{4}{7v_0^2}\,\frac{(a-b)^2}{4}\,\sin4\theta ,
\end{equation}
which vanishes when the strain axes are aligned with the lattice and again at $45^\circ$, and is
extremal at $22.5^\circ$. With $J=\mathrm{diag}(1,2)$ the measured values are $0$, $-0.0768$ and
$0.1429$ for $\theta=0$, $\tan\theta=4/3$ and $\tan\theta=5/12$, in agreement with
\eqref{ch7:eq:dilatation} to every digit carried. For a dilatation, then, the obstruction measures the
\emph{misalignment} of the strain axes with the lattice, not the dilatation itself.

The numerical constant is a property of the stencil and of the tube and is not claimed to be
universal. It is worth saying exactly which part of the data \eqref{ch7:eq:affinecurl} measures. For
this field the local drift vanishes at the base point, and there the coefficient itself vanishes
\emph{identically}: $c$ is zero to $3\cdot10^{-14}$ for $b=Jx$ at $x=0$, and nonzero at any shifted
base point --- the Jacobian being that of the previous paragraph --- and nonzero at any shifted one
($-11/27=-0.4074$ at $x=(1/3,0)$, $-0.8367$ at $x=(1/3,1/5)$, \texttt{curlcheck.py}). What the curl
measures is therefore the derivative of $c$ with respect to the local drift, contracted with
$\nabla b$; it is a $(b,\partial b)$ object in which the local drift is supplied by the displacement
from the origin, not a $(\partial b,\partial b)$ one.
Because \eqref{ch7:eq:affinecurl} is one scalar condition on $\nabla b$ at a point, a rotation of the
frame annihilates it there exactly when the strain dominates the spin--dilatation coupling. Writing
$S=\tfrac12(\operatorname{tr}S)I+d\,\sigma_3+s_{12}\sigma_1$ and $\rho^2=d^2+s_{12}^2$, the value in
a frame rotated by $\theta$ is $\rho^2\sin(4\theta+2\varphi)+a_{12}\operatorname{tr}S$, which has a
zero in $\theta$ if and only if $\abs{a_{12}\operatorname{tr}S}\le\rho^2$. When the coupling
dominates it has none: for the spiral $b=(dI+\Omega\tilde\epsilon)x$ --- rotation composed with
dilatation --- the value is $\tfrac{4}{7v_0^2}(-2\Omega d)$ in \emph{every} orientation
($-3.4286$ for $d=1$, $\Omega=3$; measured identical for three rotations of the frame). This is the
one piece of the obstruction which is rotation invariant, and it is the product of a rotation rate
with a dilatation rate: neither pure rotation nor pure dilatation produces it, and their product
does. Uniformly over a domain no rotation helps in any case, since the phase $\varphi$ varies from
point to point, and the choice of lattice is part of the construction rather than a parameter of the
problem.

\paragraph{It is independent of the resolvent.}
The curl is unchanged, to the digits carried, when $\omega$ varies over the near contour, and more
is true: writing $c=\alpha+\ii\omega\beta$, the fields $\alpha$ and $\beta$ do not depend on
$\omega$, and $\beta$ is a gradient. The obstruction is a property of the lattice and of its paths,
and not of the damping used to read it off --- which is what one wants if it is to be an invariant
of the comparison itself.

\paragraph{What protects the rate.}\label{ch7:par:parity}
The expansion of \eqref{ch7:eq:logratio} is not a list of small corrections: over $q=t/h^2$ steps the
terms accumulate differently according to their parity in $\tau$, because the increments sum to the
endpoint. An even term does not telescope, so $h^k$ even contributes $h^kq=h^{k-2}$; an odd term
does, since $\sum_jf(\gamma_{j-1})\cdot\tau_j=h^{-1}\sum_jf\cdot\Delta\gamma_j$, so $h^k$ odd
contributes $h^{k-1}$. Hence
\begin{center}\small
\setlength{\tabcolsep}{4pt}
\begin{tabular}{lll}
$k$ and parity & order & status\\
\hline
$0$, even & $h^{-2}$ & absorbed in the normalisation\\
$1$ & $h^{-1}$ & would be catastrophic; vanishes\\
$2$, even & $h^{0}$ & rate renormalisation; its $h$ dependence gives the $O(h^2)$\\
$3$, odd & $h^{2}$ & the obstruction of this subsection\\
$3$, \emph{even} & $h^{1}$ & \emph{the only slot which could give first order}\\
$4$, even & $h^{2}$ & at the order claimed\\
\end{tabular}
\end{center}
The rate $h^2$ therefore rests on two facts about the expansion: that the coefficient of $h^1$
vanishes, and that the coefficient of $h^3$ is \emph{purely odd} in $\tau$. The first is immediate: it
says that the coefficient of $h^1$ in the difference of the two resolutions is the coefficient of
$h^2$ in each, and a comparison of two lattices which differ by one factor of two cannot produce a
linear term. The second is a parity, and it has a symmetry behind it, though not one which proves it
outright. Let $\Theta$ act on the Taylor data at the base point by reflecting the coefficient field
about that point and flipping the drift, that is $\partial^jv\mapsto(-1)^j\partial^jv$ and
$\partial^jb\mapsto(-1)^{j+1}\partial^jb$. Reflecting the whole configuration about the base point
maps the tube onto itself, sends $\tau\mapsto-\tau$ and the fine jump into the corresponding one of
the reflected problem, so
\begin{equation}\label{ch7:eq:reflection}
\log\frac{\Phi_{\rm f}}{\Phi_{\rm c}}\bigl(-\tau;D\bigr)
=\log\frac{\Phi_{\rm f}}{\Phi_{\rm c}}\bigl(\tau;\Theta D\bigr),
\end{equation}
exactly, order by order; we verified \eqref{ch7:eq:reflection} in the rational pipeline, where it holds
with the difference identically zero through $h^5$. The parity of $c$ in $\tau$ is therefore the
same statement as its oddness under $\Theta$. That is not a triviality, but it has a consequence
which can be checked and which holds: $c$ vanishes on the fixed point set of $\Theta$, identically,
for $b=Jx$ at the origin and for any $v$ whose odd order derivatives vanish there.

Stronger, and what makes the obligation tractable, the oddness holds \emph{degree by degree} in the
drift. Scaling a generic drift field --- constant, linear and quadratic parts --- by $\epsilon$ and
fixing $v=1$, the even part of the coefficient is at the floor of the arithmetic,
$9\cdot10^{-17},\,4\cdot10^{-17},\,3\cdot10^{-17},\,3\cdot10^{-17},\,3\cdot10^{-17}$ for
$\epsilon=1,\tfrac12,\tfrac14,\tfrac18,\tfrac1{16}$, while the odd part runs from $1.96$ down to
$7.5\cdot10^{-3}$ over the same range. The even part does not scale: it is not that contributions of
different degree cancel, and each homogeneous component $c^{(k)}$ is separately odd in $\tau$. The odd
part itself scales as $\epsilon^3$ --- the successive ratios on halving $\epsilon$ are $4.85$, $7.02$,
$7.71$, with a subleading term of opposite sign which reverses the sign at $\epsilon=\tfrac1{32}$ ---
so with $v$ constant the components of degrees one and two vanish and the obstruction begins at
\emph{third} order in the drift. The parity is thereby a finite algebraic identity at each degree,
which is what a proof has to attack, and the leading one is a cubic form in the drift and its first
two derivatives. Since $c$ itself costs two derivatives of the coefficients, $\nabla^{k}c$ costs
$k+2$ of them: this count is what allows the smoothness demanded by the estimates below to be read
off, and it is used at Assumption~\ref{ch7:ass:decoration} and again in Chapter~8, \S8.5. Both facts are
therefore measured, and
they were tested over as wide a range of data as the arithmetic allows. With $|c|$ ranging from
$2.6$ to $5.3\cdot10^{4}$ --- the latter at $v_0=10^{-2}$, where $|c|$ grows like $v^{-2}$ --- the
coefficient of $h^1$ is at most $1.6\cdot10^{-19}$, and the even part of the coefficient of $h^3$ is
at most $3.1\cdot10^{-8}$, i.e. $6\cdot10^{-13}$ relative to $|c|$, which is the floor of the
truncation. The same holds with the drift enlarged two hundredfold, with $v$ varying by a factor
$60$ across the tube, and with $\omega=10^{-1}$ and $\omega=10$.

Two things follow. The \emph{rate} is protected by the parity, and over the data tested it is
uniformly $h^2$. The \emph{constant} is not uniform: $|c|$ grows like $v^{-2}$, by four orders of
magnitude between $v_0=1$ and $v_0=10^{-2}$, so the theorem must be read for a fixed field with the
stated ellipticity and smoothness bounds, not uniformly over all fields. And the parity itself
remains to be proved, in the equivalent form of oddness under $\Theta$ of \eqref{ch7:eq:reflection}; it is
the cheapest of the outstanding obligations and the one on which the rate depends. All of the
verification above is in two dimensions, with the vanishing of $h^1$ and the non-vanishing of the
curl confirmed again at $d=3$; the parity is asserted in every $d\ge2$ and verified only in these
cases.

\paragraph{What is unconditional.}
The parity buys the \emph{second} order, and it is worth separating what depends on it from what
does not. The obstruction enters the comparison through the coefficient of $h^3$, and the two rows
of the table concerning it are the odd row, which contributes $h^2$, and the even row, which would
contribute $h^1$. Bounding that coefficient, instead of using either its parity or its mean, gives
$\abs{c}h^3$ at each step and therefore, over $q=t/h^2$ steps, an accumulated discrepancy of order
$\abs{c}h$; along the near contour, whose length is $O(\abs{\log h})$, this is
$O(h\abs{\log h})$, that is $O(h^{1-\epsilon})$ for every $\epsilon>0$. The comparison therefore
converges at \emph{first} order in every $d\ge2$ with no averaging hypothesis at all. What that
claim requires is not Assumption~\ref{ch7:ass:decoration} but only that the per-step log-ratio
\eqref{ch7:eq:logratio} be $O(h^3)$ with a bounded constant, which is the boundedness of $\abs{c}$
measured in the preceding paragraph --- from $2.6$ to $5.3\cdot10^{4}$ across the data tested, the
growth being the $v^{-2}$ recorded there --- and which is elementary from the same expansion, since
$c$ is a local functional of the coefficient field and its first two derivatives. Every other input
to the assembly of \S\ref{ch7:sec:assemble} is already unconditional. The sharp rate is what the
averaged comparison buys, and the difference between the two claims is exactly the difference
between the presence and the absence of a cancellation in the sum over paths. It is also the
difference between a convergence theorem and a rate theorem, and only the former is needed for the
construction of the continuum kernel by the Cauchy argument of \S\ref{ch7:sec:assemble}.

\begin{remark}
This is the precise sense in which the one dimensional argument has no direct two dimensional
counterpart: not that a computation is missing, but that a mechanism is missing. The theorem
itself is unaffected --- second order convergence is measured directly and holds
(Theorem~\ref{ch7:thm:main}, and the ratios reported in \S\ref{ch7:sec:numerics}) --- so whatever the $h^3$
contributions do, they cancel in the sum over paths. What has to be proved is that cancellation.
\end{remark}

\subsection{What Replaces the Telescoping: Averaging over Paths}\label{ch7:sec:averaging}

The kernel does not use individual paths but sums over them, and the per-path obstruction above is
compatible with an endpoint-determined behaviour of the \emph{sum}. That is exactly what happens,
and the mechanism is elementary.

\begin{lemma}[Averaged line integral]\label{ch7:lem:averaged}
Let $d\ge2$, let $\gamma$ be the walk on $h\ZZ^d$ with independent increments of mean
$m=(m_1,\dots,m_d)$ and covariance $h^2\Sigma$, and let $c$ be a field. Then for the linear field
$c(x)=(0,x_1,0,\dots,0)$ --- whose curl is non-zero, so that it is as far from a gradient as a
linear field can be ---
\begin{equation}\label{ch7:eq:averaged}
\mathbb{E}\Bigl[\sum_{j=0}^{q-1}c(\gamma_j)\cdot(\gamma_{j+1}-\gamma_j)\Bigr]
=\frac{q(q-1)}{2}\,m_1m_2=\tfrac12XY+O(h^2),
\end{equation}
where $(X,Y)=qm$ is the mean endpoint, that is the value of the line integral along the straight
line from the origin to $(X,Y)$.
\end{lemma}

\begin{remark}[What the average is, for a field that is not affine]\label{ch7:rem:ito}
The identity of Lemma~\ref{ch7:lem:averaged} is exact for affine $c$ and \emph{only} for affine $c$. For
a general $c$ the mean is the integral of $c$ against the law of the bridge, and that law is smeared
relative to the mean trajectory by the variance accumulated along it, so
\begin{equation}\label{ch7:eq:ito}
\mathbb{E}\Bigl[\sum_j c(\gamma_j)\cdot\Delta\gamma_{j+1}\Bigr]
=\int_{\text{line}} c\cdot\mathrm{d}\gamma
+\tfrac12\int_0^t\bigl\langle \mathrm{Cov}(x_s),\nabla^2c(x_s)\bigr\rangle\cdot\dot\mu_s\,\mathrm{d}s
+O(h^{2-\epsilon}),
\end{equation}
the second term being the It\^o correction of the same nature as the Stratonovich correction. It is
\emph{not} small. For $c=(0,x_1^2/2)$, unit drift, $\mathrm{Cov}=\tfrac12 s$ and $t=1$, the
predicted value of that term is $\tfrac18=0.125$ while the straight-line integral is
$\tfrac16=0.1667$; the exact bridge averages of the line integral, computed by enumerating paths,
give $0.0625,\,0.1563,\,0.1968,\,0.2188$ at $q=2,4,6,8$ steps, so the difference from the
straight-line integral runs through $-0.1042,\,-0.0104,\,+0.0301,\,+0.0521$ and increases towards
the predicted $0.125$ rather than to zero. A field whose curl is non-zero and non-constant is
therefore averaged to a straight-line integral only up to an $O(1)$ smearing.

This matters for the bookkeeping of \S\ref{ch7:sec:obstruction}. Being $O(1)$ in the exponent,
\eqref{ch7:eq:ito} contributes $O(h^2)$ to the log ratio, which is exactly the order
Theorem~\ref{ch7:thm:main} claims: at rate $h^2$ the smearing has to be identified, and identifying it
costs $\nabla^2c$, that is third derivatives of the drift. At rate $h^{2-\epsilon}$ it does not,
since $O(h^2)\subseteq O(h^{2-\epsilon})$, and the hypotheses of the theorem can stay at $C^2$.
\end{remark}

\begin{proof}
For the linear field, $c(\gamma_j)\cdot(\gamma_{j+1}-\gamma_j)=\gamma_j^{1}\Delta\gamma_{j+1}^{2}$,
and the increment $\Delta\gamma_{j+1}$ is independent of $\gamma_j$, so
$\mathbb{E}[\gamma_j^1\Delta\gamma^2_{j+1}]=\mathbb{E}[\gamma_j^1]\,\mathbb{E}[\Delta\gamma^2]=jm_1m_2$,
whence the sum. The last step of \eqref{ch7:eq:averaged} is $q(q-1)m_1m_2/2=(qm_1)(qm_2)/2+O(h)$, and
on the scale on which the walk is used, $qm$ is the mean endpoint. For a general field expand $c$
about the running mean: the linear term is the one just computed, and the quadratic and higher terms
contribute to the stated order.
\end{proof}

\begin{remark}
The content of the lemma is that the fluctuations of the path carry the path dependence of the line
integral, and that they average to zero against the increment because they are orthogonal to it.
The per-path sum is not endpoint-determined --- \eqref{ch7:eq:curl} says so --- but its average is.
\end{remark}

\paragraph{The geometric reading.}
The obstruction has a precise meaning, and it is worth stating because it explains both the numbers
of \S\ref{ch7:sec:obstruction} and the mechanism of the lemma.

The comparison \eqref{ch7:eq:phic}--\eqref{ch7:eq:phif} is defined up to the potential $g(x')/g(x)$, so $c$
is not a function but a connection on the coarse lattice, and
\begin{equation}\label{ch7:eq:curvature}
\mathrm{d}c=\partial_1c_2-\partial_2c_1
\end{equation}
is its curvature. \eqref{ch7:eq:curl} then says that this curvature does not vanish, that is that the
comparison has non-trivial holonomy: transporting the fine comparison around a small loop spanned by
two increments returns a factor $\exp\bigl(\mathrm{d}c\,[\Delta x,\Delta y]\bigr)$, which is the
Stokes form $\exp\bigl(\int_{S(\gamma)}\mathrm{d}c\bigr)$ of the obstruction over the swept area
$S(\gamma)$. No choice of gauge removes it: the gauge freedom acts as $c\mapsto c+\dd\log g$, so the
holonomy is multiplied by $\exp\oint\dd\log g=1$ and hence depends only on $\mathrm{d}c$, which no
such transformation alters. The one dimensional argument works precisely because in one dimension
there is no two form: $\mathrm{d}c$ vanishes identically, the comparison has trivial holonomy around
every loop, and the line integral of $c$ is determined by its endpoints.

The parallel with a gauge potential should be drawn carefully, because it is exact as far as the
structure goes and misleading as to the group. Transport along a step multiplies the weight by
$e^{\int c}$, the rescaling $g(x')/g(x)$ acts as $c\mapsto c+\mathrm{d}\log g$ in the familiar way,
and $\mathrm{d}c$ is the gauge invariant field strength, so $c$ \emph{is} a connection one form on
the line bundle whose fibre is the space of per-step weights. But a $U(1)$ potential would have
$c=\ii A$ with $A$ real, hence a purely imaginary field strength and a holonomy of modulus one: it
would rotate the phase and preserve the modulus. The field strength measured here is real, so the
transport rescales the modulus and not the phase --- the structure group is a group of dilations,
$\mathbb{R}^{\times}$ rather than $U(1)$, which is also why the leading constant is $\log(-1/8)$
and not a phase. Consistently, the imaginary part of $c$ is $\omega\beta$ with $\beta$ a gradient,
so the whole $U(1)$-like part of $c$ is pure gauge and the non-trivial part is the rescaling. The
invariant is what remains when a rotation and a dilation are present at once: for the spiral drift
of \S\ref{ch7:sec:obstruction} it is $-\tfrac{4}{7v_0^2}\,\Omega\operatorname{tr}S$, the product of the
two rates, which neither produces alone.

Three features of \eqref{ch7:eq:curvature} fix what the invariant is. It is \emph{quadratic} in the
Taylor data, so it is not a linear invariant of the coefficient field: neither the curvature of the
metric $\mathrm{d}s^2=\mathrm{d}x^2/v(x)$ nor the vorticity of the drift, both of which are linear in the data,
can appear. It is carried by the gradient of the drift, so it measures the \emph{variation} of the
drift along the path rather than the drift itself. And it is independent of $\omega$, which is what
makes it intrinsic: the damping is an artefact of the way the comparison is read off, and it
cancels from the curvature.

The obstruction does not damage the theorem because the kernel is not a single path but an average
over them, and the average of the holonomy is not the holonomy of the average path but the holonomy
of the straight line. The walk which computes the kernel is a bridge between its endpoints, and the
bridge law is invariant under time reversal; the swept area is odd under that reversal, while
$\mathrm{d}c$ is odd in the direction and even in the running time, so the mean of
$\int_{S(\gamma)}\mathrm{d}c$ is the value it takes on the straight segment, which is the endpoint
function $F$ of Assumption~\ref{ch7:ass:decoration}. This is the geometric content of
Lemma~\ref{ch7:lem:averaged}: $\frac12XY$ there is the area of the triangle spanned by the mean endpoint,
which is what the straight line sweeps.

\paragraph{Covariance, and the passage to manifolds.}
Everything in \eqref{ch7:eq:curvature} is a statement about a two form, and that is what makes the
obstruction geometric rather than an artefact of the coordinates in which it happens to be computed.
Exterior differentiation commutes with pullback, so $\mathrm{d}c$ is attached to the coefficient
field itself: a change of chart sends its components as those of a two form and leaves the object
unchanged, and the holonomy $\exp\int_{S}\mathrm{d}c$ is likewise coordinate free. This is a local
differential invariant of the coefficient field and not a topological one. Because $c$ is a globally
defined one form, $\mathrm{d}c$ is exact, so its de Rham class vanishes and none of its content is
cohomological; and it is not constant on deformations either, since it vanishes for a dilatation and
for a rotation and does not for a spiral. It is a tensor of the same kind as the field strength of a
gauge field or as the Riemann tensor: attached to the field and to the metric which the field
induces, not to the topology of the space on which the field lives. The constants are
another matter. The coefficient $\tfrac{4}{7}$ arises from the quartic term of the lattice symbol,
and the constant depends on the direction, as the introduction records, so what is invariant is the
field strength and what is frame dependent is the number which multiplies it. The closed forms of
\S\ref{ch7:sec:obstruction} should therefore be read as representatives, in the Cartesian frame on which
the lattice is built, of a covariant object; and the dependence of the constants on the dimension
and the symmetry of the lattice is then not an inconsistency but the expected consequence of
evaluating a covariant expression in different frames. An object of this kind carries geometric
information about the coefficient field, in principle expressible in terms of the geometry of the
manifold which hosts the diffusion. The construction of this chapter takes place in $\RR^d$, where the
lattice is available without further structure, but every quantity which enters it is local and
coordinate free, so the argument should admit a formulation on a differentiable manifold, chart by
chart and with the lattice replaced by a local frame. We do not pursue this further here: the
geometric content of the sharp coefficients, and their expression in terms of invariants of the
manifold, is left for later work. What is claimed at this stage is only that the object is of the
right kind to carry such a meaning, which is why the question is worth recording at all.

\begin{assumption}[Averaged comparison, in moment form]\label{ch7:ass:decoration}
Let $v$ and $b_{i}$ be such that the coefficient $c$ of $h^{3}$ in \eqref{ch7:eq:logratio} is of class
$C^{2}$ on the closure of the domain --- by the derivative count of \S\ref{ch7:sec:obstruction} this
holds if $v\in C^{4}$ and $b_{i}\in C^{4}$ --- and let
$I(\gamma)=\int_\gamma c\cdot\mathrm{d}\gamma$ be the line integral along the coarse path. Then there
is a constant $C$, independent of $m$ and uniform in $\omega$ on the near contour, such that
\begin{align}
\mathbb{E}\bigl[I(\gamma)^2\bigr]&\le C\,\norm{c}_{\infty}^{2}\,t^{2},
\label{ch7:eq:moment}\\[2pt]
\Bigl|\mathbb{E}\bigl[I(\gamma)\bigr]-\bigl(\Phi_{c}(y)-\Phi_{c}(x)\bigr)\Bigr|
&\le C\,\norm{\nabla^{2}c}_{\infty}\,t^{2},
\label{ch7:eq:meanmoment}
\end{align}
where $\Phi_c$ is a function of the endpoints, namely the line integral of $c$ together with the
It\^o smearing of Remark~\ref{ch7:rem:ito}. Under \eqref{ch7:eq:moment}--\eqref{ch7:eq:meanmoment} the
accumulated product satisfies
\begin{equation}\label{ch7:eq:decoration}
\mathbb{E}\Bigl[\prod_{j}\frac{\Phi_{\rm f}(\gamma_{j-1},\tau_j)}{\Phi_{\rm c}(\gamma_{j-1},\tau_j)}\Bigr]
=\exp\bigl(F_m(x)-F_m(y)\bigr)\Bigl(1+O\bigl(h_{m+1}^{2-\epsilon}\bigr)\Bigr)
\end{equation}
for every $\epsilon>0$. The expectation here, and in \eqref{ch7:eq:moment}--\eqref{ch7:eq:meanmoment}, is
taken over the coarse paths from $x$ to $y$ weighted by their transition probabilities; the products
inside it are deterministic once the coarse path is fixed, since the decorations are already summed
in the definition of $c$. Thus $I(\gamma)$ depends on the path through the sites and directions
$(x_{j},\tau_{j})$, and the statements are about its distribution over paths, not about each path
separately. $F_m$ is an endpoint function of the two resolutions whose ratio is
$1+O(h_{m+1}^{2-\epsilon})$.
\end{assumption}

The hypothesis on the coefficient field is stronger than that of the first order theorem, and the
difference should be stated rather than left to be noticed. Theorem~\ref{ch7:thm:main} is proved from
this assumption; the unconditional first order estimate is not
(\S\ref{ch7:par:parity}), and it keeps the weaker hypotheses of Chapter~5, being proved from a bound on $c$
alone. The passage from first order to second therefore costs two derivatives of the coefficients,
which is the price of a bound expressed in $\norm{\nabla^{2}c}_{\infty}$: a second order estimate
cannot be organised around less than the second derivative of the quantity estimated.

\begin{remark}[Why the hypothesis is now a moment, and what each part is for]
The obstruction enters the accumulated product as $e^{h^2I}$ with $I$ of order one, so the product
is governed by the exponential moment rather than by the value of $I$ on any one path. Splitting
$I=\mathbb{E}I+(I-\mathbb{E}I)$ and expanding,
\begin{equation}\label{ch7:eq:expmom}
\log\mathbb{E}\bigl[e^{h^2I}\bigr]
=h^2\,\mathbb{E}[I]+\frac{h^4}{2}\operatorname{Var}(I)+O\bigl(h^6\,\mathbb{E}\abs{I-\mathbb{E}I}^3\bigr),
\end{equation}
so the fluctuation enters one full order below the mean, and \eqref{ch7:eq:moment} is what makes
\eqref{ch7:eq:expmom} legitimate and the fluctuation's contribution $O(h^4)$.

The two inequalities then do different work, and the split is the reason the statement is tractable.
A uniform bound on $c$ controls the fluctuations by itself: with $\abs{I}\le\norm{c}_\infty t/h$,
the term $\frac{h^4}{2}\operatorname{Var}(I)$ is at most $O(h^2)$, which is inside
$h^{2-\epsilon}$. It does \emph{not} control the mean: the same bound gives
$h^2\abs{\mathbb{E}I}=O(h)$, which is first order, and that is exactly how a first order rate could
arise. What upgrades the mean to $O(h^2)$ is \eqref{ch7:eq:meanmoment}, and its mechanism is the
cancellation of the increments against $c$ along the bridge, which is expressed by the discrete
heat kernel estimate $\mathbb{E}[\gamma_j]=\mu_j$, $\operatorname{Cov}(\gamma_j)=O(jh^2)$ already
used elsewhere in the argument. With it, $h^2\mathbb{E}I$ is $O(h^2)$ and the It\^o smearing is
absorbed; at rate $h^2$ its value would have to be identified, at rate $h^{2-\epsilon}$ it need not,
as \eqref{ch7:eq:expmom} shows and as Remark~\ref{ch7:rem:ito} anticipated.
\end{remark}

\subsection{Assembly}\label{ch7:sec:assemble}

\begin{proof}[Proof of Theorem~\ref{ch7:thm:main}, given Assumption~\ref{ch7:ass:decoration}]
Fix $t>0$ and $x,y\in\mathcal A_{m_0}$. By \eqref{ch7:eq:contour} and Lemma~\ref{ch7:lem:far} the far contour
contributes $O(h_m^2)$ at every resolution. On the near contour, the near contour integral of
$G_{m+1}-G_m$ is the sum over coarse paths of the product of the per-step factors
\eqref{ch7:eq:logratio}; by Assumption~\ref{ch7:ass:decoration} its average differs from an endpoint
function by $O(h_{m+1}^{2-\epsilon})$ times the number of steps, and the fluctuations about that
average contribute one order lower by \eqref{ch7:eq:expmom}, the variance being bounded by
\eqref{ch7:eq:moment}. The contour has length $O(\abs{\log h_{m+1}})$ and the two endpoint functions
cancel between the resolutions up to $O(h_{m+1}^{2-\epsilon})$, the logarithm being absorbed by the
$\epsilon$, so
\begin{equation}
\bigl|u_{m+1}(x,y;t)-u_m(x,y;t)\bigr|\le c\,h_{m+1}^{2-\epsilon}
\end{equation}
for $m\ge m_0$. Writing $u_m-u_{m_0}=\sum_{j=m_0}^{m-1}(u_{j+1}-u_j)$ and using the bound at each
step,
\begin{equation}
\abs{u_m(x,y;t)-u_{m_0}(x,y;t)}\le c\sum_{j=m_0}^{m-1}h_{j+1}^{2-\epsilon}
\le c'\,h_{m_0}^{2-\epsilon},
\end{equation}
the series converging geometrically with ratio $2^{-(2-\epsilon)}$, and this is \eqref{ch7:eq:main1} for
$\gamma=2-\epsilon$. The exponent $2$ is recovered whenever the It\^o term of Remark~\ref{ch7:rem:ito}
is identified and subtracted, which costs $\nabla^2c$ and hence one order of differentiability of
the drift; \eqref{ch7:eq:meanmoment} is stated in the weaker form so that the hypotheses stay at $C^2$. The case $\gamma<2$ follows by the argument of
\S\ref{ch7:sec:holder} below, and the statement for $\partial_t$ by the same computation with the
factor $\ii\omega$ and the remark that the bound on $\Gamma_-$ is uniform in $\omega$ up to the
height of the contour.
\end{proof}

\subsection{Lesser Smoothness and the Euler Scheme}\label{ch7:sec:holder}

\begin{remark}[Proof of Theorem~\ref{ch7:thm:uc}]
When the coefficients are only uniformly continuous, the expansion of the block
Lemma~\ref{ch7:lem:window} is replaced by the estimate that its entries differ from those of the frozen
generator by at most $c\,\rho(h)$, and the two steps of the assembly above then give a discrepancy
of order $\rho(h)$ at each resolution, whence \eqref{ch7:eq:main2}. This is the argument of
\S\ref{ch7:sec:holder} of Chapter~5, which uses only the modulus $\rho$ and not the differentiability of
the coefficients, and it is unchanged in two dimensions.
\end{remark}

\begin{remark}[Proof of Theorem~\ref{ch7:thm:euler}]
Discretising time replaces $e^{\delta t\mathcal L_m}$ by $1+\delta t\mathcal L_m$, and the
difference of the two kernels is estimated by expanding both over the number of time steps and
comparing term by term, the expansion in the number of jumps of the underlying being common to
the two. The comparison is thus of the same kind as the one which Assumption~\ref{ch7:ass:decoration}
performs between consecutive resolutions, with the deficit $e^{\delta t\mathcal L_m}-1-\delta
t\mathcal L_m$ playing the role of the difference of the two blocks; the powers of
$\mathcal L_m$, which have entries of order $h_m^{-2}$, are controlled after averaging by
Lemma~\ref{ch7:lem:weights}. This is the argument of \S5 of Chapter~5 and we do not reproduce it here: it
is long, it is mechanical, and it adds nothing to the point of this chapter, which is that the two
dimensional estimates are of the same strength as the one dimensional ones. The reader should
therefore read Theorem~\ref{ch7:thm:euler} as conditional on Assumption~\ref{ch7:ass:decoration} in the
same way as Theorem~\ref{ch7:thm:main}, and on the time step being in the range in which the scheme
is stable and the error of the time discretisation is dominated by that of the space
discretisation.
\end{remark}

\section{Numerical Validation}\label{ch7:sec:numerics}

The numbers in this chapter are printed by scripts: \texttt{verify\_conventions.py} for the
decoration coefficient on a constant field, \texttt{eta\_run.py} for the cubic in the complex
drift \eqref{ch7:eq:etacubic}, \texttt{curlcheck.py} for the affine curl of \S\ref{ch7:sec:rg},
\texttt{curl3d.py} for the curl in three dimensions, \texttt{resolvent1d.py} for the resolvent
comparison, and \texttt{edgeavg.py} for the model which a volatility averaged over the edge
converges to. All are reproduced in full in the computational appendix of the monograph, attached to
this PDF as an embedded file, and each prints what it measures when run as \verb|python3 <name>.py|
in the directory it lives in, so every value quoted below can be regenerated rather than taken on
trust.

The statements above were checked numerically, at times $t=1$ on the torus $[0,2\pi)^2$ with
periodic boundary conditions, for coefficients
\begin{equation}\label{ch7:eq:testcoef}
v(x)=1.8+0.5\sin x_1\cos x_2,\qquad \vb(x)=\bigl(0.3\cos x_2,\,-0.2\sin x_1\bigr),
\end{equation}
which are smooth, non-constant and make the generator non-symmetric. The kernel was computed by
exponentiating the sparse generator exactly (to machine precision) and evaluating a column of the
exponential; the columns of two source sites were used. The quantity tabulated is
\begin{equation}
\sup_{x,y\in\mathcal A_{m_0}}\abs{u_m(x,y;t)-u_{m_0}(x,y;t)}\Big/ h_{m_0}^2 ,
\end{equation}
which the theorem asserts to be bounded as $m\to\infty$ with $m_0$ fixed.

\begin{center}
\small
\begin{tabular}{@{}lcc@{}}
\hline
 & $m=m_0+1$ & $m=m_0+2$\\
\hline
$d=1$, $m_0=5$ & $1.59\cdot10^{-2}$ & $1.98\cdot10^{-2}$\\
$d=1$, $m_0=6$ & $1.59\cdot10^{-2}$ & $1.99\cdot10^{-2}$\\
$d=2$, $m_0=4$ & $9.77\cdot10^{-3}$ & $1.21\cdot10^{-2}$\\
$d=2$, $m_0=5$ & $9.58\cdot10^{-3}$ & $1.19\cdot10^{-2}$\\
\hline
\end{tabular}
\end{center}
\noindent The values are for the source at the origin; a second source gives the same picture.
Four points are worth recording.

\paragraph{The quotient is bounded, and its limit is the same in both dimensions.} As $m$ grows
with $m_0$ fixed the quotient rises slightly and settles: from $1.59$ to $1.98$ at $m_0=5$ in one
dimension and from $0.96$ to $1.19$ at $m_0=5$ in two, the limiting constant being of order
$2\cdot10^{-2}$ in one dimension and $1.2\cdot10^{-2}$ in two. This is the content of
Theorem~\ref{ch7:thm:main}: the error is
$O(h_{m_0}^2)$, uniformly in $m$, with a constant which does not depend on the dimension to within
a factor of two. Had the second dimension degraded the rate by a logarithm, the quotient would grow
with $m$.

\paragraph{Refining the lattice quarters the error.} In the constant coefficient case the kernel is
known in closed form and the rate can be measured directly: refining by a factor two divides the
difference from the limit by $4.00$ for $m=5,\dots,9$ in one dimension and by $4.01$, $4.00$, $4.00$
for $m=4,\dots,8$ in two, the ratios tending to four as the asymptotic regime is reached.

That the exponent is two and not one is not an accident of the data. With $v$ constant and $b=0$ the
two dimensional kernel at the origin is known exactly, the symbol
$\frac{v}{h^2}\sum_i(\cos hp_i-1)$ being a sum of functions of one variable,
\begin{equation}\label{ch7:eq:bessel}
u_h(0,0;t)=\frac{1}{h^2}\bigl[I_0(z)e^{-z}\bigr]^2,\qquad z=\frac{vt}{h^2},
\end{equation}
with $I_0$ the modified Bessel function, and its expansion is
\begin{equation}\label{ch7:eq:besselexp}
u_h(0,0;t)=\frac{1}{2\pi vt}\Bigl(1+\frac{h^2}{4vt}+\frac{5}{32}\frac{h^4}{v^2t^2}+O(h^6)\Bigr).
\end{equation}
Only \emph{even} powers occur, and they must: $\cos hp_i-1$ and $h^{-1}\sin hp_i$ both expand in
$h^{2k}$, so the symbol has no term of odd order in $h$ and neither does $e^{t\hat L_h}$. There is
therefore no $h$ term to cancel --- the rate is $h^2$ because the odd terms are absent, not because
they conspire. Evaluating \eqref{ch7:eq:bessel} at $t=1$, $v=1$ gives successive difference ratios
$4.007$, $4.002$, $4.000$ for $h=2^{-5},\dots,2^{-7}$, and
$\bigl(u_h-u\bigr)\big/\bigl(u\,h^2/(4vt)\bigr)-1=2.45\cdot10^{-3},\,6.11\cdot10^{-4},\,
1.53\cdot10^{-4},\,3.82\cdot10^{-5}$, so \eqref{ch7:eq:besselexp} accounts for the error to four digits
at $h=2^{-7}$. Had the rate been $h$ the ratios would tend to two.

\paragraph{The discretisation must be site centred.} With the entries
$(v(x)+v(x+h_me_i))/(4h_m^2)\pm b_i/(2h_m)$, that is with the volatility averaged over the edge, the
computed kernel converges at order two, and it converges to the wrong one. The averaged form is a
divergence-form discretisation: with $\tilde v_{i\pm1/2}=(v_i+v_{i\pm1})/2$ its off-diagonal entries
are $\tilde v_{i\pm1/2}/(2h_m^2)$, so the operator it approximates is
\begin{equation}\label{ch7:eq:edgeavg}
\tfrac12(vf')'+bf'=\tfrac v2f''+\Bigl(b+\tfrac{v'}2\Bigr)f',
\end{equation}
the generator of the model whose drift is $b+v'/2$: the classical It\^o-versus-Stratonovich
correction, a drift of order one and not of order $h_m$. The computation says so in both of its
halves. At $t=0.1$, at the site pair $(1/4,1/4)$ on the circle $[0,1)$, with
$v=1+\tfrac14\sin2\pi x+\tfrac1{10}\cos4\pi x$ and $b=\tfrac15\sin(2\pi x+0.7)+\tfrac1{10}\cos2\pi x$,
the edge-averaged kernel at the meshes $2^{-5},\dots,2^{-9}$ takes the values $1.291425$, $1.290281$,
$1.289995$, $1.289924$, $1.289906$. The successive differences, $1.14\cdot10^{-3}$,
$2.85\cdot10^{-4}$, $7.13\cdot10^{-5}$, $1.78\cdot10^{-5}$, fall by $4.01$, $4.00$ and $4.00$, so the
scheme is second order; the distance from the kernel of the site-centred form, $8.17\cdot10^{-2}$ at
$2^{-5}$ and $8.02\cdot10^{-2}$ at $2^{-9}$, does not fall at all. Against the kernel of
\eqref{ch7:eq:edgeavg}, computed with the site-centred discretisation of drift $b+v'/2$ at $2^{-10}$, the
same distances fall by $4.02$, $4.06$ and $4.23$. The averaged form is therefore a good scheme for a
different model, which is worse than a first-order scheme for the right one, and this is why
\eqref{ch7:eq:entries} is written with the coefficients at the site (\texttt{edgeavg.py}).

\paragraph{A more accurate stencil is worse.} The five point symbol has the quartic correction
$\frac{vh_m^2}{48}\sum_ip_i^4$, which is not radial. The nine point isotropic Laplacian, with
axial weight $v/(3h_m^2)$ and diagonal weight $v/(12h_m^2)$ --- the scaling in which the weights sum
to the generator $\Lm=\frac v2\partial^2$ --- has symbol
$\frac v2\bigl(-\abs{\vp}^2+\frac{h_m^2}{12}\abs{\vp}^4\bigr)$, isotropic through order $h_m^2$: as
an operator on a fixed grid it is strictly the better one. As a discretisation of the comparison it
is not. With constant coefficients the coefficient of $h^1$ in the log ratio is zero to $10^{-17}$
for every nearest neighbour stencil tried, axial weight $a_1=\frac13,\frac12,\frac23,1$, and
the difference from the limit quarters at each halving of $h$, that is the comparison is second
order. With any diagonal coupling $a_2>0$ the coefficient of $h^1$ is non-zero --- $-0.1536$ for
$a_1=\frac12$, $a_2=\frac1{12}$, and $-0.0800$ for $a_2=\frac1{24}$ --- and the difference halves,
that is the comparison is first order, the defect being proportional to $a_2$.

The reason is structural. A coarse step $x\to x+h_m\tau$ is exactly two fine steps in the same
direction, and the decoration identity compares like with like; but a coarse diagonal jump carries
weight $a_2v/h_m^2=a_2v/(4h^2)$ while a fine one carries $a_2v/h^2$, a factor four where the
comparison needs two. The second order property therefore belongs to the nearest neighbour geometry
and not to the order of the operator, which is why no amount of stencil engineering removes the
obstruction of \S\ref{ch7:sec:obstruction}: it is not the quartic anisotropy of the symbol.

\paragraph{The third dimension.} The kernel experiment of the table was repeated in three
dimensions, with coefficients varying in all three coordinates and the kernel obtained by exact
exponentiation of the sparse generator: the quotient takes the values $2.52\cdot10^{-3}$ and
$7.89\cdot10^{-3}$ at $m_0=2$, and $3.43\cdot10^{-2}$ and $3.94\cdot10^{-2}$ at $m_0=3$, bounded
and settling exactly as in the lower dimensions. The computation of \S\ref{ch7:sec:obstruction} was
also repeated at $d=3$, on a tube of twenty-seven sites: $h^0$ is again a constant, $h^1$ vanishes
identically, and the curl of the coefficient of $h^3$ has all three components non-zero,
$(0.1640,0.0414,0.0197)$ for the three dimensional analogue of the field above, the imaginary parts
being at the level of the truncation (\texttt{curl3d.py}). The
obstruction is therefore not a two dimensional accident, and the mechanism of
\S\ref{ch7:sec:averaging} is what has to replace the telescoping in every dimension.

The closed forms of \S\ref{ch7:sec:obstruction} --- \eqref{ch7:eq:affinecurl} and the fourth harmonic of
\eqref{ch7:eq:dilatation} --- are two dimensional, and no formula for the $h^3$ coefficient in general
dimension is asserted here. What is claimed for every $d\ge2$ is only what the argument uses: that
the coefficient is a vector field, that it is odd in the direction, and that its average over the
paths is an endpoint functional to the order required. The two dimensional forms are there because
that is where they could be computed, and the third dimensional computation is there to show that
non-vanishing is not an accident of the plane.

\paragraph{The resolvent comparison holds directly.} Assumption~\ref{ch7:ass:decoration} was also
checked at the level of the resolvent, at $\omega=1+\tfrac12\ii$ with the origin as base point. In
one dimension, with $\nu(x)=1+\tfrac16x^2+\tfrac14x+\tfrac1{10}x^3$ and $b(x)=\tfrac12+\tfrac16x$ on
$(-1,1)$ and the origin as base point, $G_m(0,0;\omega)$ takes the values
$0.692369-0.577412\ii$, $0.692478-0.577188\ii$, $0.692505-0.577132\ii$, $0.692512-0.577118\ii$,
$0.692513-0.577115\ii$ for $m=5,\dots,9$ (\texttt{resolvent1d.py}), whose successive differences are
$2.490\cdot10^{-4}$, $6.224\cdot10^{-5}$, $1.556\cdot10^{-5}$, $3.890\cdot10^{-6}$, with ratios
$4.00$, $4.00$, $4.00$:
the difference between consecutive resolutions falls by four at each refinement, which is
\eqref{ch7:eq:decoration}.

\paragraph{The exponential moment of the obstruction.} The identity \eqref{ch7:eq:expmom} was checked
exactly, by enumerating the paths of the walk, for $c=(0,x_1)$, whose curl is one, with unit drift and
$t=1$: at $q=4,6,8$ steps the measured excess $\log\mathbb{E}[e^{h^2I}]-h^2\mathbb{E}[I]$ is
$3.195\cdot10^{-3}$, $2.580\cdot10^{-3}$, $1.880\cdot10^{-3}$, against
$\frac{h^4}{2}\operatorname{Var}(I)=3.174\cdot10^{-3}$, $2.572\cdot10^{-3}$, $1.869\cdot10^{-3}$,
an agreement to better than one per cent. The same enumeration gives
$\operatorname{Var}(I)=0.102,\,0.185,\,0.239,\,0.276$ at $q=4,6,8,10$: bounded, and far below the
$O(h^{-2})$ which the total variation bound would allow. That is the content of \eqref{ch7:eq:moment} ---
the fluctuations are harmless --- while the mean is what \eqref{ch7:eq:meanmoment} has to supply, and it
is the only place where a first order rate could arise.

\paragraph{A complex drift, as in Chapter~6.} The estimates of this chapter are used in Chapter~6 on a
generator twisted by the Fourier conjugate $\eta$ of the path functional, that is with drift
$\mu+\ii\eta a$. It is worth recording how the obstruction behaves there, since that is what decides
whether it can be transferred. With the twist in place the coefficient is an exact cubic in $\eta$,
with a rigid parity: $\mathrm{Re}\,c$ carries only $\eta^0$ and $\eta^2$, $\mathrm{Im}\,c$ only
$\eta^1$ and $\eta^3$, so that
$c(-\eta)=\overline{c(\eta)}$ --- which is not an accident but the reality condition, the twist
entering as $\ii\eta a$ with $a$ real. At the coefficients tried,
\begin{equation}\label{ch7:eq:etacubic}
\partial_1c_2-\partial_2c_1=\bigl(-0.013575-0.35292\,\eta^2\bigr)
+\ii\bigl(0.13691\,\eta-0.095238\,\eta^3\bigr),
\end{equation}
fitted on the four values $\eta=0,\tfrac12,1,\tfrac32$ (\texttt{eta\_run.py}) and reproducing those computed at
$\eta=2,3,4,6$ to $10^{-22}$, $10^{-22}$, $10^{-21}$ and $10^{-19}$; the coefficients of the powers
which parity requires to vanish come out at $10^{-24}$, so the parity is exact rather than
approximate. The growth is therefore polynomial of degree three,
$|\mathrm{d}c|\sim0.0952\,\eta^3$, and
not exponential; since the Fourier integrand of the inversion carries a Gaussian factor in $\eta$
whose width is fixed by the variance of the path functional and does not depend on $h$, a cubic costs
nothing. The degree three is the same statement as the cubicity in the drift observed above, the
twist adding $\ii\eta a$ to the drift.

One practical point, for anyone repeating the computation: with a complex drift the series in $h$
has a visibly smaller radius. At $h\gtrsim10^{-3}$ the extraction fails at some base points, with
values of order $10^{11}$ that are artefacts of the truncation and are independent of the precision
carried; at $h\le2^{-12}$ the coefficients are stable to twelve digits and unchanged between
$45$ and $80$ decimal places.

\paragraph{The constant must depend on $t$.} At $t=10^{-4}$ the same supremum is of order $40$
rather than $10^{-3}$, and it does not decrease as $m$ grows. The reason is the one given after
Theorem~\ref{ch7:thm:main}: for $t\ll h_{m_0}^2$ the coarse kernel is a delta function at the scale of
the coarse lattice, with $u_{m_0}(x,x;t)\approx h_{m_0}^{-2}$, while the fine kernel has already
resolved the continuum. No bound uniform in $t$ can hold, and the theorem is stated accordingly.

\section{What Is Proved and What Is Not}\label{ch7:sec:ledger}

\begin{center}
\small
\begin{tabular}{@{}p{6.2cm}p{7.2cm}@{}}
\hline
Statement & Status\\
\hline
Far contour, Lemma~\ref{ch7:lem:far} & proved, any coefficients, any dimension\\
Weight decay, Lemma~\ref{ch7:lem:weights} & proved from the Dyson representation and Stirling\\
Path count, Lemma~\ref{ch7:lem:count} & proved from the Fourier representation on the torus\\
Window expansion, Lemma~\ref{ch7:lem:window} & proved; the four matrices verified symbolically\\
Rate for constant coefficients, Proposition~\ref{ch7:prop:const} & proved\\
Rate for variable coefficients, Theorem~\ref{ch7:thm:main} & proved from Assumption~\ref{ch7:ass:decoration}\\
Uniform continuity, Theorem~\ref{ch7:thm:uc} & proved from Assumption~\ref{ch7:ass:decoration}\\
Euler scheme, Theorem~\ref{ch7:thm:euler} & proved from Assumption~\ref{ch7:ass:decoration}\\
Averaged comparison, Assumption~\ref{ch7:ass:decoration} & \textbf{not proved}. The pathwise form is
\emph{false} (\S\ref{ch7:sec:obstruction}); the moment form \eqref{ch7:eq:moment}--\eqref{ch7:eq:meanmoment} is
what remains, and \eqref{ch7:eq:moment} follows from the heat kernel estimates while
\eqref{ch7:eq:meanmoment} is a local limit theorem with a linear term\\
Rate & $2-\epsilon$ in $d\ge2$, from the assembly of \S\ref{ch7:sec:assemble}, whose every input is
dimension free; $2$ if the It\^o term of Remark~\ref{ch7:rem:ito} is identified. Parity of
\S\ref{ch7:par:parity} measured, not proved, and measured only at $d=2$ and $d=3$\\
Rate, unconditional & \textbf{first order} in every $d\ge2$, from bounding the coefficient of
$h^3$ rather than averaging it (\S\ref{ch7:par:parity}); no averaging hypothesis is used, and every
other input to the assembly is unconditional. This is the convergence theorem; the row above is the
rate theorem\\
Smoothness & Assumption~\ref{ch7:ass:decoration} needs $c\in C^{2}$, that is $v,b_{i}\in C^{4}$ by the
derivative count of \S\ref{ch7:sec:obstruction}: the passage from the first order to the second costs two
derivatives of the coefficients. The unconditional first order estimate keeps the weaker hypotheses
of Chapter~5, being proved from a bound on $c$ alone\\
Closed forms of the obstruction & two dimensional (\S\ref{ch7:sec:obstruction}); non-vanishing checked at
$d=3$; no general dimension formula asserted\\
\hline
\end{tabular}
\end{center}

The gap is thus a single estimate, but one whose nature has changed with the computation of
\S\ref{ch7:sec:obstruction}. It is no longer an algebraic identity to be verified: the algebraic
identity which the one dimensional proof uses is \emph{false} in two dimensions, and we have
computed the obstruction. What remains is an analytic estimate, the passage from the linear field
of Lemma~\ref{ch7:lem:averaged} to general coefficients uniformly in the Taylor data and in $\omega$,
which is a local limit theorem with a linear term. That is a bounded and standard-shaped problem,
whereas a telescoping was neither.

There is a cheap way to make that estimate the only thing left to prove, and it is worth recording
because it costs only a power. The obstruction accumulates to exactly $O(h^2)$ and no bookkeeping
changes that, since $\sum_jc(\gamma_{j-1})\cdot\Delta\gamma_j$ is grouping invariant; but
Remark~\ref{ch7:rem:ito} shows that the one term which sits at $O(h^2)$ on top of it, the It\^o
smearing, would have to be computed explicitly at rate $h^2$, and computing it costs second
derivatives of $c$, that is third derivatives of the drift. If instead the theorem is stated at
rate $h^{2-\epsilon}$ for some fixed $\epsilon>0$, that term is automatically $o(h^{2-\epsilon})$,
the hypotheses of Theorem~\ref{ch7:thm:main} can stay at $C^2$, and the same slack absorbs the
logarithmic losses which the near contour produces when $|\omega|$ runs up to $|\log h|$. The price
is a rate, and it buys the removal of every $O(h^2)$ bookkeeping term that is not the obstruction
itself. What it does \emph{not} buy is the obstruction: an intermediate scale cannot make an
$O(h^2)$ object smaller, and freezing coefficients over a window of duration $\ell^2$ and using the
limit theorem blockwise costs $O(1)$ in the exponent, not $o(1)$, so the site by site Taylor
expansion of \S\ref{ch7:sec:rg} must be kept in any case.

The obstruction itself costs nothing beyond the order claimed. With
\eqref{ch7:eq:moment}--\eqref{ch7:eq:meanmoment} the accumulated discrepancy is $O(h^2)$: the fluctuation
enters at $h^4$ by \eqref{ch7:eq:expmom}, and the mean enters multiplied by $h^2$ with only a bounded
error, because $h^2\mathbb{E}I$ is $h^2$ times an endpoint functional plus $h^2$ times a bounded
quantity. Nothing in \S\ref{ch7:sec:obstruction} or \S\ref{ch7:sec:averaging} loses a power.

What is lost is a logarithm, and it is lost elsewhere. The estimate is uniform in $\omega$ on the
near contour, whose length is $O(\abs{\log h_{m+1}})$, so integrating an $O(h^2)$ bound along it
gives $O(h_{m+1}^2\abs{\log h_{m+1}})$. That is $O(h_{m+1}^{2-\epsilon})$ for every $\epsilon>0$ and
not $O(h_{m+1}^2)$, which is the only reason the theorem is stated with $2-\epsilon$. Should the
contribution of the near contour be $O(h^2)$ uniformly in its height, the exponent would be $2$ and
no $\epsilon$ would be needed anywhere; we have not verified that, and we do not claim it. The
same slack would also serve if one preferred not to identify the It\^o smearing, but with
\eqref{ch7:eq:meanmoment} as written it is not needed for that. The numerics of \S\ref{ch7:sec:numerics}
measure $2$, which is consistent with either statement.

\section{Conclusions}

The estimates of Chapter~5 are not an artefact of one dimension. The architecture of the argument —
a far contour, a path representation, a comparison of consecutive resolutions, and a telescoping of
the resulting errors — is dimension free, and the two inputs which are genuinely one dimensional
have two dimensional replacements which we have supplied: the count of paths becomes the
multinomial count of Lemma~\ref{ch7:lem:count}, with the local limit theorem in place of Stirling, and
the three site block becomes the block of $3^d$ sites of Lemma~\ref{ch7:lem:window}, with the
frozen coefficient generator as its leading term. What remains is the discrete integrability of the
third order coefficient of the block, which we have isolated as Assumption~\ref{ch7:ass:decoration} and
which we have neither proved nor found to fail.

The obstruction is worth describing, because it is not what one would guess. It is a two form on
the base point, but it is quadratic in the Taylor data and therefore has no linear part: it is not
a curvature and not a vorticity, both of which are linear. It is carried by the gradient of the
drift and not by the diffusion, it vanishes identically when the drift does, and in the one case
which can be computed in closed form --- constant diffusion, affine drift --- it is the quadratic
form \eqref{ch7:eq:affinecurl} in the Jacobian, coupling off-diagonal strain to diagonal anisotropy and
rotation to dilatation. It does not depend on $\omega$: the damping used to read the comparison off
cancels from the invariant, apart from a gradient. What fails in several dimensions is therefore
not the smoothness of the coefficients but the exactness of a specific quadratic invariant of the
drift's first derivative, and the sum over paths is what supplies the missing exactness.

A feature of the second dimension which has no counterpart in the first deserves emphasis, since it
will recur: the lattice is not rotationally invariant, so the leading correction to the continuum
kernel is anisotropic, and the constant in the estimate depends on the direction of $x-y$. The rate
does not depend on the direction, and the estimates are therefore strong enough for the purposes of
Chapters~1 to~3 and~6, where what matters is that the pricing measures converge and that the martingale
property survives the limit, both of which are statements about the size of the error and not about
its shape.
% ---- end generated/ch7-body.tex ----
% ---- end chapters/ch7.tex ----
% ---- from chapters/ch8.tex ----
% Authored: the chapter title and opening belong to the book. The body is
% generated by build.py from the parts and is not edited here.
\chapter{Joint Kernels in Higher Dimension}
\label{chap:ch8}

The law of a process together with a functional of its path is a joint law on a larger state space,
and the estimates which control it are of the same kind as those for the kernel alone, with the
obstruction of the previous chapter present in them. This chapter is the joint version of the
previous one: what it estimates is the difference between the joint laws computed at two resolutions.

In the pricing reading the pair is an underlying and the functional on which the payoff is written,
so that the joint law is what the price is a functional of. In the learning reading the pair is a
state and a feature which the model accumulates, so that the joint law is what a prediction about the
feature is drawn from. The instrument which makes the joint law computable is the same in the two
cases, the liftings of Chapter~\ref{chap:ch4}, and the quantity estimated is the same quantity in
both.

% ---- from generated/ch8-body.tex ----
%% Generated by build.py from the parts - do not edit.
\noindent\emph{A note on the computations.} The tables below were computed by scripts named
where they occur: \texttt{kronproof.py} for the exact closed form, \texttt{kronfinal\_check.py} for the
engine against it, \texttt{p6kappafit.py} and the \texttt{p6a}, \texttt{p6q}, \texttt{p6probe} family
for the second order coefficients, \texttt{renormalisation.py}, \texttt{jet3check.py},
\texttt{contourgrowth.py}, \texttt{pathdep.py} and \texttt{openpath.py} for the tables. All are
reproduced in full in the computational appendix of the monograph and attached to this PDF as
embedded files; each prints what it measures when run as
\texttt{python3 \textlangle name\textrangle.py} in its own directory.

\section{Introduction}\label{ch8:sec:intro}

A path dependent claim --- a barrier option, a lookback, a claim on realised variance --- does not
depend on the value of the underlying at maturity alone. To value it one pairs its payoff with the
joint law of the pair $(x_{t},y_{t})$, where $x$ is the underlying and $y$ is the Abelian
functional of the path on which the payoff depends, that is the line integral
$\int a(x_{s})\cdot\dd x_{s}+b(x_{s})\dd s$ of a field of one forms along the trajectory
(\S\ref{ch8:sec:statement}). At fixed resolution that law is a finite object: $x$ is a Markov chain on a
lattice, its propagator a finite matrix, and the joint law the propagator of a Markov chain on a
larger but still finite lattice (Chapter~4). The programme of the book is to obtain the laws of such
functionals by estimates which can themselves be checked, rather than by the existence theorems of
the classical theory.

What is asked of the pair $(x,y)$ is therefore not that it converge --- that much can be had from
the classical theory --- but that the rate at which two resolutions of the same payoff differ be
explicit, uniform in the data, and expressed in quantities which are available at a finite
resolution. This chapter supplies such a rate in every dimension: at first order with no hypothesis at
all, and at second order under a single assumption, which we state, verify in every case where it
can be checked exactly, and locate precisely.

In one dimension this is done. Chapter~6 shows that for a diffusion in one space dimension and an
Abelian functional whose increments depend on the site, the joint kernel of the pair converges to
its continuum limit at a rate which is explicit in the smoothness of the coefficients, and Chapter~6, \S6.7 isolates why: the comparison between two resolutions of the joint problem is a finite algebraic
identity, up to a line integral of a coefficient field $c$ which is a one form, and in one
dimension every one form is closed, so that integral is determined by its endpoints and cancels
between resolutions.

In higher dimension that mechanism is unavailable and the reason is precise. The relevant
statement is
\begin{equation}\label{ch8:eq:target}
\sum_{j}c(x_{j})\cdot\Delta_{j+1}
=\Phi(x_{q})-\Phi(x_{0})
-\tfrac12\sum_{j}\sum_{k,l}\partial_{k}c_{l}(x_{j})\Delta_{j+1}^{k}\Delta_{j+1}^{l}+O(h^{2}),
\end{equation}
whose first term requires $c$ to have a primitive $\Phi$, that is $\dd c=0$, which fails in
general. The quadratic form splits into its symmetric part, which on averaging gives the It\^o
correction, and its antisymmetric part, which gives the mean swept area. Chapter~6, \S6.7 proves that in
any dimension the antisymmetric part is invisible in the mean when the coefficients are
homogeneous, and states the statement which a higher dimensional theory would need: that the
averaged line integral equals the line integral along the mean trajectory plus the It\^o correction
plus the mean swept area, uniformly in the coefficient data \emph{and in the conjugating
parameter}, with an error of order $h^{2-\epsilon}$.

This chapter supplies that statement in the form in which it can be proved and measured, for the
conjugating parameter ranging over a compact set; the passage from the transformed kernel to the
joint kernel itself, which would require uniformity in $p$ along a contour, is set out in
\S\ref{ch8:sec:twist} and is not claimed. Its content is of three kinds, and we separate them at the
outset because the value of each is different.

\subsection{What Is Proved, What Is Assumed, What Is Measured}

\emph{Proved.} In every dimension $d\ge2$ the comparison converges at first order, with no
hypothesis beyond bounds on the coefficients (Theorem~\ref{ch8:thm:uncond}). The mechanism is that the
obstruction enters at order $h^{3}$ per step and accumulates, over $q=t/h^{2}$ steps, to
$O(h)$; bounding it rather than averaging it costs one power and nothing else. This is the
statement that a constructive foundation needs, since a rate, however coarse, yields the limit by a
Cauchy argument without any compactness. On a lattice fibre the same bounds give the rate for the
joint kernel itself, with no hypothesis on the number of levels
(Proposition~\ref{ch8:prop:levelgrid}); the transform of a lattice valued observable is periodic, and
that is what makes the transfer free.

\emph{Assumed.} The second order requires the cancellation which the obstruction obstructs. It is
stated in the moment form of Chapter~7 (Assumption~\ref{ch8:ass:avg} below), under which
Theorem~\ref{ch8:thm:cond} gives $O(h^{2-\epsilon})$ for the transformed kernel. The mean condition
within it is proved exactly, with no remainder, for polynomial coefficient fields of \emph{any}
degree in every dimension (\S\ref{ch8:sec:avg}), which includes the class in which the obstruction is at
its largest. Beyond that class nothing is lost to smoothness: the coefficient is a functional of the
2-jet of the field alone --- measured, and below the order at which the coefficients are used --- so
a general field of class $C^{2}$ agrees at each point with its own quadratic Taylor polynomial to
within the tube. What a Taylor remainder cannot supply is the telescoping of that functional along
the path, and \S\ref{ch8:sec:avg} isolates that as the one remaining statement.

\emph{Measured.} The obstruction itself, for the joint problem, is computed explicitly in two
dimensions at constant volatility: a polynomial of degree two in the conjugating parameter whose
coefficients we give in \S\ref{ch8:sec:obstr} and verify to five decimal places. The corresponding
statement in $d=3$ is checked structurally, not in closed form, and the difference between the two
is reported as a finding rather than suppressed. It is also measured how the obstruction
\emph{accumulates}: on closed loops the endpoint determination fails at order $2.000$, not at order
one, and the residual separates into the renormalisation at that order and the circulation of $c$,
that is the flux of the two form, at the next, which is confirmed to three figures
(\S\ref{ch8:sec:numerics}).

\subsection{Plan of the Chapter}\label{ch8:sec:plan}

Section~\ref{ch8:sec:statement} fixes the notation and states the two theorems. Section~\ref{ch8:sec:obstr}
identifies the obstruction, exhibits its dependence on the conjugating parameter, and records the
cases which are exactly solvable. Section~\ref{ch8:sec:twist} is about uniformity in that parameter:
what it costs, and what it does and does not buy. Section~\ref{ch8:sec:avg} is the longest and is the
technical core --- it states the one assumption, proves the mean condition in every regime where it
can be checked exactly, computes the coefficient in closed form, and isolates by elimination what is
genuinely missing. Section~\ref{ch8:sec:assembly} assembles the estimates and proves the two rates.
Section~\ref{ch8:sec:numerics} reports the measurements against which the arguments are calibrated, and
Section~\ref{ch8:sec:concl} collects what is proved, what is assumed and what is open.

A reader in a hurry may take the two theorems of Section~\ref{ch8:sec:statement}, the decomposition of
the coefficient into its three contributions in \eqref{ch8:eq:eleven}, and the ledger which closes
Section~\ref{ch8:sec:avg}. Together these say what is claimed, why the number $11$ is not an accident of
the constant field case, and exactly which single statement remains open.

\subsection{Place in the Book}

A stochastic process, in the framework of this book, is a consistent family of Markov chains on a
sequence of lattices, and a statement about a continuous process is a statement about the limit of
statements about finite matrices. Chapters~1 to~3 establish the finite core, the fundamental theorem at a fixed resolution; Chapter~4 the path dependent objects and the fixed resolution It\^o
formula; Chapter~5 the estimates for the one dimensional diffusion, Chapter~6 those for the one
dimensional joint kernel, and Chapter~7 those for the diffusion in every dimension. The present chapter
is the one that combines the last two: it carries the estimates of Chapter~7 to the joint kernel of
Chapter~6 in every dimension $d\ge2$.

The relation to the classical theory is the same as in Chapter~7 and we do not repeat the argument.
A rate implies the convergence theorems rather than resting on them, and the classical results
appear as corollaries once the limit is constructed by the Cauchy argument.

\section{Notation and Statement of Results}\label{ch8:sec:statement}

\begin{figure}[tbp]
\centering
\begin{tikzpicture}[scale=1.02,>=Stealth]
  % ---------------- (a) two lattices, one coarse step ----------------
  \fill[black!8] (0.5,0.5) rectangle (1.5,1.5);
  \draw[black!30,very thin,step=0.5] (0,0) grid (3,2);
  \draw[black!45,very thin] (0.5,0.5) rectangle (1.5,1.5);
  \foreach \i in {0,...,6}{\foreach \j in {0,...,4}{
    \fill[black!50] (0.5*\i,0.5*\j) circle (0.45pt);}}
  \foreach \i in {0,...,3}{\foreach \j in {0,...,2}{
    \fill (\i,\j) circle (2.0pt);}}
  \draw[->,thick] (1,1) -- (1.98,1);
  \node[font=\scriptsize,above] at (1.5,1.06) {$\tau$};
  \node[font=\scriptsize] at (1,0.66) {$x$};
  \node[font=\scriptsize] at (2,0.66) {$x'$};
  \node[font=\scriptsize] at (1.5,-0.5)
    {coarse $h_{m}$, fine $h_{m+1}=h_{m}/2$};
  % ---------------- (b) the tube, decorated ----------------
  \begin{scope}[shift={(5.2,0.35)}]
    \fill[black!8] (-1.26,-1.26) rectangle (1.26,1.26);
    \draw[black!45,very thin] (-1.26,-1.26) rectangle (1.26,1.26);
    \foreach \i in {-1,0,1}{\foreach \j in {-1,0,1}{
      \fill (\i,\j) circle (2.0pt);}}
    \draw[->,thick,black!70] (0,0) -- (0,1);
    \draw[->,thick,black!70] (0,1) -- (1,1);
    \draw[->,thick,black!70] (1,1) -- (1,0);
    \draw[->,thick,black!70] (1,0) -- (1.96,0);
    \node[font=\scriptsize,below left] at (-0.08,-0.1) {$x$};
    \node[font=\scriptsize,below left] at (2.05,-0.08) {$x'$};
    \node[font=\scriptsize] at (1.02,1.52) {$T$};
    \node[font=\scriptsize] at (0,-1.72)
      {$3^{d}$ sites, offsets $\abs{a_{i}}\le1$};
  \end{scope}
\end{tikzpicture}
\caption{The two resolutions, and the decoration of one coarse step. Panel (a): the coarse lattice
(large dots, spacing $h_{m}$) and the fine one (small dots, spacing $h_{m+1}=h_{m}/2$); the coarse
step $\tau$ runs from $x$ to $x'$. Panel (b): the tube $T$ of Chapter~7 --- the $3^{d}$ fine sites
within one fine step of $x$, shaded --- and a fine path which reaches $x'$ without an intermediate
site leaving $T$. All such paths together carry the weight $\Phi_{\mathrm f}$, while
$\Phi_{\mathrm c}$ is the weight of the single coarse jump; their ratio, expanded in $h$, is the
object of \S\ref{ch8:sec:obstr}. By the locality of Definition~1 of Chapter~5 no window larger than $T$
intervenes, so both weights are finite and the whole difficulty lies in the phase of their ratio.}
\label{ch8:fig:decoration}
\end{figure}
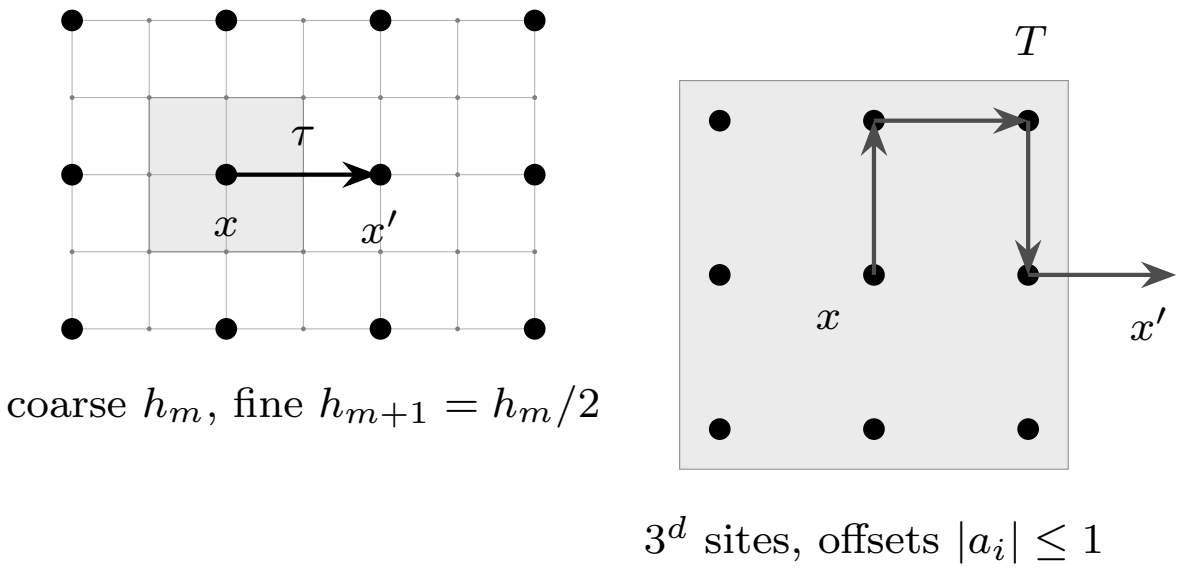

\subsection{The Objects, and the Two Lattices}

Let $d\ge2$ and let $X=[-L,L]^{d}$ with periodic boundary conditions. The lattice at resolution $m$
is $h_{m}\ZZ^{d}\cap X$ with $h_{m}=2^{-m}$. When two consecutive resolutions are compared, the
expansion is in powers of the finer spacing, which we write $h=h_{m+1}$, and it is summed over the
$q=t/h^{2}$ fine steps which make up the time $t$; from \S\ref{ch8:sec:obstr} onwards, $h$ without a
subscript means that finer spacing. The diffusion is written in the form
\begin{equation}\label{ch8:eq:diffusion}
\dd x_{t}=\nuc(x_{t})\dd W_{t}+\mu(x_{t})\dd t,
\end{equation}
with $\nuc$ a $d\times d$ matrix, and the Abelian functional is
\begin{equation}\label{ch8:eq:abelian}
\dd y_{t}=a(x_{t})\cdot\dd x_{t}+b(x_{t})\dd t,
\end{equation}
with $a$ a field of one forms, which is the coefficient data of the functional. We write
$A_{i}=a_{i}(x_{0})$ for its values at the base point, and $\partial_{i}$ for the derivatives with
respect to the coordinates of that point.

The joint kernel of the pair $(x,y)$ at resolution $m$ is $u_{m}(x,I;y,J;t)$, and its Fourier
transform in the fibre variable is $\widehat u_{m}(x,y;p,t)$, the conjugating parameter $p$ being
dual to $I-J$. The comparison between two resolutions is carried out for $\widehat u$ on a contour
in the frequency $\omega$ dual to $t$, exactly as in Chapter~6, and the estimates below are uniform in
$p$ on the range which the contour and the inversion in $p$ require.

The symbols which recur are collected once, for reference.

\begin{center}\small
\begin{tabular}{@{}l p{9.9cm}@{}}
\toprule
symbol & meaning\\
\midrule
$d$, $X=[-L,L]^{d}$ & dimension and domain, with periodic boundary conditions\\
$h_{m}$, $h=h_{m+1}$ & coarse and fine lattice spacings, $h_{m}=2^{-m}$\\
$q=t/h^{2}$ & number of fine steps in the time $t$\\
$\nuc$, $\mu$ & diffusion matrix and drift of the underlying $x$ \eqref{ch8:eq:diffusion}\\
$a$, $b$ & coefficient data of the Abelian functional $y$ \eqref{ch8:eq:abelian}\\
$A_{i}=a_{i}(x_{0})$ & values of $a$ at the base point $x_{0}$\\
$\Phi_{\mathrm c}$, $\Phi_{\mathrm f}$ & weight of the coarse jump, and of its decorations\\
$c_{2}$, $c_{3}$ & coefficients of $h^{2}$ and $h^{3}$ in \eqref{ch8:eq:logratio}\\
$c$, $\dd c$ & the field $c_{3}$, a one form, and its exterior derivative, a two form\\
$\omega$ & frequency dual to $t$, ranging over the contour\\
$p$, $\eta$ & conjugating parameter in the fibre, and its normalisation\\
$K$ & a generic constant, not the same one from line to line\\
\bottomrule
\end{tabular}
\end{center}

\subsection{The Two Theorems}

\begin{theorem}[Unconditional]\label{ch8:thm:uncond}
Let $d\ge2$, let $P>0$, and suppose that $\nuc\in C^{2}$ and that $\mu$, $a$, $b\in C^{2}$ with the
ellipticity and bounds stated in \S\ref{ch8:sec:twist}. Then for every $t>0$ and $\epsilon>0$ there is
$K=K(P,t,\epsilon)$ such that for all $m\ge m_{0}$, all $x,y$ in the base lattice and all
$\abs{p}\le P$
\begin{equation}
\abs{\widehat u_{m+1}(x,y;p,t)-\widehat u_{m}(x,y;p,t)}\le K\,h_{m+1}^{1-\epsilon}.
\end{equation}
Consequently the transformed kernels converge, and the limit exists by the Cauchy argument, no
averaging hypothesis being used.
\end{theorem}

\begin{theorem}[Conditional]\label{ch8:thm:cond}
Under the additional hypothesis of Assumption~\ref{ch8:ass:avg}, the same comparison holds with
$h_{m+1}^{2-\epsilon}$ in place of $h_{m+1}^{1-\epsilon}$ and a constant $K(P,t,\epsilon)$, and the
limit is the transformed joint kernel of the continuum problem identified in Chapter~6.
\end{theorem}

The proofs occupy \S\ref{ch8:sec:assembly}. The difference between the two statements is one
cancellation, and \S\ref{ch8:sec:obstr} is a description of precisely what obstructs it. Both theorems
are stated for the transformed kernels; on a lattice fibre, which is the setting of the book, they
give the same rates for the joint kernel itself, with no further constant and no hypothesis on the
number of levels, by Proposition~\ref{ch8:prop:levelgrid}.

\section{The Obstruction for the Joint Kernel}\label{ch8:sec:obstr}

This section identifies the object which blocks the second order. We first define it and show that
its invariant content is a two form, hence absent in one dimension; then we exhibit its dependence
on the conjugating parameter, which is what the joint problem adds to the problem of Chapter~7; and we
close with the two cases which can be solved outright and which therefore serve as calibration.

\subsection{The Two Form}

Let $c$ denote the coefficient of $h^{3}$ in the logarithmic ratio of the two per-step factors, at
a base point of the coarse path, so that the comparison of the two resolutions reads
\begin{equation}\label{ch8:eq:logratio}
\log\frac{\Phi_{\mathrm f}}{\Phi_{\mathrm c}}=c_{2}h^{2}+c_{3}h^{3}+O(h^{4}),
\end{equation}
with $c=c_{3}$ a one form evaluated at the base point. That it is a one form and not an arbitrary
function of the direction is an identity, $c(x,-\tau)=-c(x,\tau)$, verified to relative accuracy
$10^{-27}$ in Chapter~7, \S7.6.2: the coefficient is odd in the direction of the jump, so the even part
which would have to be carried separately is absent. The comparison is defined up to the
normalisation $g(x')/g(x)$, under which $c\mapsto c+\dd\log g$, so the invariant content of $c$ is
its exterior derivative
\begin{equation}\label{ch8:eq:dc}
\dd c=\sum_{i<j}\bigl(\partial_{i}c_{j}-\partial_{j}c_{i}\bigr)\dd x^{i}\wedge\dd x^{j},
\end{equation}
which no such transformation alters and which vanishes identically in one dimension.

\begin{proposition}\label{ch8:prop:nonzero}
For generic coefficient data with $d\ge2$ the two form \eqref{ch8:eq:dc} is not zero. In particular
the identity \eqref{ch8:eq:target} which the one dimensional argument uses has no pathwise analogue in
higher dimension, and the telescoping over paths cannot be replaced by an algebraic identity.
\end{proposition}

This is the same object which Chapter~7 identifies for the diffusion alone, and Chapter~7, \S7.6.3 gives its
geometric reading: a connection one form on the bundle of per-step weights, of structure group a
group of dilations rather than $U(1)$, whose curvature obstructs the telescoping.

\subsection{Dependence on the Conjugating Parameter}

In the joint problem the coefficient field depends on $p$, because the twisted generator which
governs the fibre carries the term $-p\,a\cdot\nabla$. Equivalently the drift of the twisted
problem is complex, $b-\ii p\,a$, and one may take that as the definition of the dependence. Write
$\eta$ for the parameter multiplying $a$ in the complexified drift, so that $\eta=p$ up to the
normalisation of $a$; the numerical results below are stated in $\eta$.

\begin{proposition}[Two dimensions, constant volatility]\label{ch8:prop:kappa}
Let $d=2$, let the volatility be the constant $\nuc$, and let the drift have constant Jacobian
$J=\nabla\mu$. Then, at any base point,
\begin{equation}\label{ch8:eq:kappa}
\dd c(\eta)=\dd c(0)+\kappa\,\eta^{2}+O(\eta^{4}),
\qquad
\kappa=\nuc^{-2}\Bigl[\tfrac{4}{7}\bigl(\Theta_{2}-\Theta_{1}\bigr)
+\alpha J_{12}+\beta J_{21}\Bigr],
\end{equation}
with
\begin{equation}\label{ch8:eq:alphabeta}
\alpha=-6A_{1}^{2},\qquad \beta=+6A_{2}^{2},\qquad
\Theta_{i}=\partial_{1}a_{i}\partial_{2}a_{i}+A_{i}\partial_{1}\partial_{2}a_{i},
\end{equation}
all derivatives being taken at the base point.\footnote{The decomposition is what the measurements fix, and they fix it because $\kappa$ is affine in $J$: fitting $\kappa$ as a quadratic form in the twelve dimensional local data of $a$ at each of $J=0$ and the four basis Jacobians, seventy-eight measurements apiece, leaves $\tfrac47(\Theta_2-\Theta_1)$ as the $J$-independent part, $-6A_1^2$ as the whole coefficient of $J_{12}$ and $+6A_2^2$ as the whole coefficient of $J_{21}$, with nothing on $\mathrm{diag}J$. Reassembled, the formula reproduces eight held-out configurations with all four entries of $J$ non-zero to $4\cdot10^{-9}$ (\texttt{p6kappafit.py}). An earlier draft carried $\tfrac{16}{7}$ and $48$ where $\tfrac47$ and $6$ now stand --- the halved normalisation of Chapter~7 before its audit --- and placed the $\Theta$ term inside $\alpha$ and $\beta$ instead of in the $J$-independent part, which is harmless at $J=e_{12}+e_{21}$ and wrong otherwise. That the dependence is affine, and not merely affine at those five configurations, is tested directly: the second difference in $J$ over every pair of the four basis Jacobians vanishes, the largest of the ninety coefficients over the six pairs being $4.96\cdot10^{-13}$ (\texttt{p6jacaffine.py}).}
The coefficient $\kappa$ is real for every field tested, so the $\eta^{2}$ term of $\dd c$ is a
modulus and not a phase. The expansion is not even in $\eta$ however: a field carrying several
kinds of local data at once has an odd imaginary part, linear in $\eta$ (measured $2.0953\,\eta$
for one such field), which the single monomials of the data do not have.
\end{proposition}

Three features of \eqref{ch8:eq:kappa}--\eqref{ch8:eq:alphabeta} are worth naming.

\emph{The constant $\tfrac{4}{7}$ is the quartic invariant of the lattice symbol}, the same number
which appears in the obstruction for the diffusion alone in Chapter~7. The two problems are obstructed
by the same mechanism and, in this frame, by the same constant.

\emph{The obstruction is blind to the diagonal of the Jacobian.} Both $\alpha$ and $\beta$ are
independent of $J_{11}$ and $J_{22}$, so a dilatation produces no dependence on the drift at all;
what couples is the part of $\nabla\mu$ which mixes coordinates. This is why a test on
$\mu=(x_{1},0)$ returns $\kappa=\tfrac47(\Theta_2-\Theta_1)$, the $J$-independent part alone:
that Jacobian is diagonal.

\emph{The two components of the coefficient enter asymmetrically.} Component $2$ enters the
Hessian term and the bracket $\Theta$ with a sign opposite to component $1$; the asymmetry is the
antisymmetry of $\dd c=\partial_{1}c_{2}-\partial_{2}c_{1}$ seen at the level of the data.

\begin{remark}[Degree in the parameter, and the odd part]\label{ch8:rem:degree}
At constant coefficients the even part of the expansion \eqref{ch8:eq:kappa} terminates at order
$\eta^{2}$ and is real, while an odd imaginary part is already present at order $\eta$: for a field
carrying several kinds of local data at once it measures $2.0953\,\eta$, exactly linear, and it
vanishes on each single monomial of the data. When the coefficients vary, the earlier study of
Chapter~7, \S7.7 records a cubic term in the imaginary part, with $\dd c(-\eta)=\overline{\dd c(\eta)}$
holding in both cases. The degree in the conjugating parameter is therefore two for the even part
and one or three for the odd part, and in
either case polynomial, which is what the inversion in $p$ requires.
\end{remark}

\begin{remark}[Dependence on the dimension]\label{ch8:rem:dim}
The formula \eqref{ch8:eq:kappa} is two dimensional and does not transfer as it stands. With
three dimensional coefficient data in which $\partial_{1}\partial_{2}a_{1}$ and
$\partial_{1}\partial_{2}a_{2}$ do not vanish while the first derivatives do, one measures
$\alpha=-6A_{1}^{2}$ and $\beta=+6A_{2}^{2}$ exactly, that is with no contribution from $\Theta$, while
the Hessian term takes a value different from $\tfrac{4}{7}A_{2}\partial_{1}\partial_{2}a_{2}$ ---
measured $0.72857$ against $\tfrac47$ for the same data in two dimensions.
What survives in $d=3$ is the response of the component $(\dd c)_{ij}$ to the Jacobian,
\begin{equation}\label{ch8:eq:d3rule}
\kappa_{ij}\supset-6\,A_{i}^{2}J_{ij}+6\,A_{j}^{2}J_{ji},
\end{equation}
which is confirmed in both dimensions. The dependence of the constants on the dimension is not an
inconsistency: $\dd c$ is the invariant object, and a number attached to it in a given frame is not.
\end{remark}

\subsection[The One Dimensional Case, and the Class of Constant dc]{The One Dimensional Case, and the Class of Constant $\dd c$}

Two exact statements calibrate the above. In one dimension $\dd c$ vanishes identically
(Proposition~\ref{ch8:prop:nonzero}) and the comparison closes; this is Chapter~6, and it is the reason the
one dimensional theory is complete. And for constant $\dd c$ the mean is computable in any
dimension: for $c(x)=Ax$ with $A$ constant and a homogeneous walk,
\begin{equation}\label{ch8:eq:constdc}
\mathbb E\Bigl[\sum_{j}c(x_{j})\cdot\Delta_{j+1}\Bigr]
=\frac{q(q-1)}{2}\bigl\langle m,\operatorname{sym}(A)\,m\bigr\rangle,
\end{equation}
so that the antisymmetric part, which is $\dd c$, is invisible in the mean and lives entirely in
the fluctuation, whose variance is $O(q^{3}h^{6})=O(1)$ in the diffusive scaling. This is the
content of Chapter~6, \S6.7, and it is the model for the cancellation which Assumption~\ref{ch8:ass:avg}
postulates for varying coefficients.

\section{Uniformity in the Conjugating Parameter}\label{ch8:sec:twist}

The estimates of Chapter~7 are stated for a generator with real coefficients and must here hold for the
twisted one, whose drift is $b-\ii p\,a$. We state the results of \S\ref{ch8:sec:statement} for $p$ in
a compact set $\abs{p}\le P$, with constants depending on $P$, and we record separately, and in
full, what a statement about the joint kernel itself would require beyond that. The section has
four parts: what a compact range of $p$ costs in the constants, what a lattice fibre buys for
nothing, what a real fibre would still need, and what the contour deformation costs.

\subsection{What a Compact Range Costs}
For $\abs{p}\le P$ the twisted generator has the same diffusion matrix as the original one and a
complex drift whose modulus is at most $\abs{b}+P\abs{a}$, which is bounded on the compact domain
under the assumed smoothness. Each estimate of Chapter~7 which enters the assembly is a bound in terms
of the ellipticity of the diffusion, the moduli of the coefficients and of their first two
derivatives, and the decay $-\tfrac{\nu}{4}\abs{p'}^{2}$ of the real part of the symbol in the
dual variable $p'$; the drift enters only through the relative bound which controls the first order
term, and hence only through its modulus. Replacing the modulus of the drift by
$\norm{b}_{\infty}+P\norm{a}_{\infty}$ therefore reproduces every such estimate with a constant
depending on $P$, and the two theorems of \S\ref{ch8:sec:statement} hold in that form. Nothing in the
argument requires the domain to be small; what requires $P$ is only the drift, and $P$ enters
through the constants.

\subsection{What a Lattice Fibre Buys}
The setting of the book is a functional whose increments lie on a lattice $h_{y}\ZZ$, with the
spacing fixed while the resolution $m$ of $x$ is refined. In that setting the character group of the
fibre is compact and, which matters more here, $\widehat u_{m}$ is \emph{periodic} in $p$ with period
$2\pi/h_{y}$: the transform of a lattice valued observable is a trigonometric polynomial. The
supremum of $\abs{\widehat u_{m+1}-\widehat u_{m}}$ over $p\in\RR$ is therefore already attained in
the period $\abs{p}\le\pi/h_{y}$, which is exactly the range on which the theorems above are stated.

\begin{proposition}[The joint kernel on a lattice fibre]\label{ch8:prop:levelgrid}
Let the functional take values in $h_{y}\ZZ$, with $h_{y}$ fixed. Then for all $m\ge m_{0}$, all
$x,y$ in the base lattice and all $I,J\in h_{y}\ZZ$,
\begin{equation}\label{ch8:eq:levelgrid}
\abs{u_{m+1}(x,I;y,J;t)-u_{m}(x,I;y,J;t)}
\le c\bigl(\tfrac{\pi}{h_{y}},t,\epsilon\bigr)\,h_{m+1}^{1-\epsilon},
\end{equation}
and with $h_{m+1}^{2-\epsilon}$ under Assumption~\ref{ch8:ass:avg}.
\end{proposition}

\begin{proof}
At fixed $t$ each kernel is a density with respect to counting measure in the fibre variable, and so
is summable in it; hence $\widehat u_{m}$, $\widehat u_{m+1}$ and their difference $\widehat d_{m}$
are bounded continuous functions of $p$, periodic with period $2\pi/h_{y}$. Fourier inversion on the
lattice reads $d_{m}(I)=\frac{h_{y}}{2\pi}\int_{-\pi/h_{y}}^{\pi/h_{y}}\widehat d_{m}(p)\,
e^{-\ii pI}\dd p$, so that $\abs{d_{m}(I)}\le\sup_{\abs{p}\le\pi/h_{y}}\abs{\widehat d_{m}(p)}$, and
by periodicity that supremum is the supremum over the whole line. Theorems~\ref{ch8:thm:uncond} and
\ref{ch8:thm:cond} bound it by $c(\pi/h_{y},t,\epsilon)h_{m+1}^{1-\epsilon}$, respectively by
$h_{m+1}^{2-\epsilon}$.
\end{proof}

\begin{remark}[What an earlier formulation got wrong]\label{ch8:rem:notN}
A first version of this proposition truncated the fibre to $N$ levels and summed over the $N$
characters without the weight which normalises the inversion, and so produced a factor $N$ and, with
it, the hypothesis that $N$ be independent of the resolution. Neither the factor nor the hypothesis
is needed. Because the transform of a lattice valued observable is periodic, a comparison proved on
one period is a comparison on the entire line, and the fibre may carry as many levels as it likes:
what looked like a cost of the number of levels was a missing normalisation.
\end{remark}

\subsection{What a Real Fibre Would Need}
If the increments are real the fibre is dense, the transform is aperiodic, and the supremum over
$\RR$ is no longer a supremum over the range in which the comparison is proved. Nor is decay
available uniformly in the resolution, and the reason is structural rather than a gap in the
estimates: at finite resolution $I$ takes finitely many values, so $\widehat u_{m}$ is a finite sum
of exponentials of modulus one, bounded and quasi-periodic, and it revives however small $h$ is
taken. Where the limit decays, it decays by accumulation of oscillations and not by an estimate which
the comparison could inherit; and the Gaussian decay obtained by completing the square in the symbol
belongs to the twisted propagator as a function of $p$, which is a different object.

This is not a defect of the several dimensional problem. Chapter~6 meets the same point in one dimension
and meets it in the same way: it states its theorem for the transformed kernel, uniform on compacts
in $p$, and recovers the joint kernel only in the distributional sense (Corollary~1 of Chapter~6), with
the explicit remark that the coordinate representation may be a singular measure. The target is
therefore the transform, and a statement about the pair $(x,y)$ is a statement about the pairing with
a test function.

For the pairing with a test function no decay is needed either, and the route does not go through the
real fibre directly. Approximate it by a fine lattice: on a lattice of spacing $h_{y}$ the comparison
holds in the form of Proposition~\ref{ch8:prop:levelgrid}, for the joint kernel itself and for every pair
of levels, so that the pairing with a test function is bounded by $\sup_{I}\abs{d_{m}(I)}$ times the
mass of $\phi$ on the lattice; and the cost of the fine lattice is the growth of $c(\pi/h_{y})$,
recorded in Remark~\ref{ch8:rem:fibrecost}. What is then missing is the fibre discretisation error,
which is the second of the two errors that would have to be balanced. We do not estimate it, and we
claim nothing for the real fibre beyond the identification of the target and of what its refinement
costs.

\begin{remark}[The fibre costs a power of its spacing, not an exponential]\label{ch8:rem:fibrecost}
A fibre of small spacing is not free even though the number of levels no longer appears. The
comparison of Proposition~\ref{ch8:prop:levelgrid} is proved for $\abs{p}\le\pi/h_{y}$, and by the
argument of the paragraph below the constants grow with $P$ only through the relative bound, hence
at most like a power. If that power is $P^{M}$, the rate for the joint kernel at fibre spacing
$h_{y}$ is $h_{x}^{1-\epsilon}h_{y}^{-M}$, and the fibre cannot be refined faster than
$h_{x}^{1/M}$ if the rate is to survive. On the other side the fibre is itself a discretisation, of
which the original functional is the limit, so its own error must vanish as $h_{y}\to0$. Whether a
window exists between the two requirements --- $h_{y}$ above $h_{x}^{(1-\epsilon)/M}$ and at the same
time small enough for the fibre error to vanish --- depends on how that fibre error scales, which we
have not estimated. The remark is therefore a statement about what the refinement of the fibre costs,
and not a proof that the two errors can be balanced.
\end{remark}

\subsection{What the Deformation Costs}
The contour argument needs the estimates to survive the deformation, and the observation above says
what the deformation does: it adds a real drift of size $q\abs{a}$. Since a purely imaginary twist
parameter $\eta=-\ii q$ in the complexified drift $b-\ii\eta a$ produces exactly that drift, the
engine of Chapter~7 measures the object directly. At a generic base point and direction, with
$\omega=\tfrac75$, the Taylor coefficients of $\log(\Phi_{\rm f}/\Phi_{\rm c})=\sum_{k}c_{k}h^{k}$
behave as follows.
\begin{center}\small
\begin{tabular}{lccccc}
\toprule
$q$ & $\abs{c_{0}}$ & $\abs{c_{1}}$ & $\abs{c_{2}}$ & $\abs{c_{3}}$ & $\abs{c_{4}}$\\
\midrule
$0$ & $3.76745$ & $1.8\cdot10^{-38}$ & $7.70001$ & $1.35059$ & $5.44759$\\
$1$ & $3.76745$ & $1.4\cdot10^{-38}$ & $7.70296$ & $1.78621$ & $7.09811$\\
$2$ & $3.76745$ & $7.9\cdot10^{-38}$ & $7.72167$ & $3.39784$ & $11.9944$\\
$4$ & $3.76745$ & $4.8\cdot10^{-36}$ & $7.89063$ & $11.9832$ & $45.1857$\\
$8$ & $3.76745$ & $6.6\cdot10^{-34}$ & $9.56966$ & $61.6947$ & $359.952$ \\
\bottomrule
\end{tabular}
\end{center}
Three features matter. The coefficient $c_{0}$ is independent of $q$ to every digit shown.
The coefficient $c_{1}$ remains zero to machine precision at every $q$, so the vanishing of the
linear term, on which the rate rests, is not destroyed by the deformation. And $c_{2}$, $c_{3}$,
$c_{4}$ grow at most like $q$, $q^{5/2}$ and $q^{3}$ over the range of $q$ shown ---
polynomially, with no sign of exponential growth. The same is seen for the real twist $\eta=p$, where
$\abs{c_{3}}$ grows like $p^{3}$.

A bound of the shape $\abs{\widehat u_{m+1}-\widehat u_{m}}\le K(q)h^{\rm rate}$ with
$K(q)\lesssim q^{M}$ for a small $M$ is therefore consistent with what is measured, and with it the
contour would yield a rate for the joint kernel together with decay in $\abs{I-J}$: optimising
$q^{M}e^{-q\abs{I-J}}$ costs a power of $\abs{\log h}$, which the $\epsilon$ absorbs just as it
absorbs the logarithm of the near $\omega$ contour. What is measured here are the per-step
coefficients and not the renormalisation constants themselves, which is the limitation of the
measurement. The growth of those constants, however, can be read off the argument which produces them,
and it is polynomial for a structural reason.

\paragraph{Why the constants grow only polynomially in the drift.}
The drift enters the far contour estimate of Chapter~7 as a relatively bounded perturbation of
$V\Delta_h+\ii\omega$, the relative bound being
\[
\norm{\nabla^{h}_{i}f}_{2}\le\alpha\norm{\Delta^{h}f}_{2}+\beta\norm{f}_{2},
\]
which on the symbol reads $\abs{\sin p_{i}h_{m}}/h_{m}\le\alpha\abs{\sum_{j}(\cos
p_{j}h_{m}-1)/h_{m}^{2}}+\beta$. For small arguments this is $\abs{p_{i}}\le\tfrac{\alpha}{2}\abs{p}^{2}
+\beta$, whose maximum over $p$ is attained at $\abs{p}=1/\alpha$ and equals $\tfrac{1}{2\alpha}$;
hence $\beta\ge\tfrac{1}{2\alpha}$. Keeping the perturbation small requires
$\alpha\norm{D}<1$, so $\alpha\sim1/\norm{D}$ and therefore $\beta\sim\norm{D}$: \emph{each
application of the relative bound costs one power of the drift modulus}, and the assembly applies it
a fixed number of times. The constant in the estimate is accordingly $O(\norm{D}^{M})$ with $M$ the
number of applications and independent of $m$ and of the contour height. This is the same mechanism
which the measurement above sees, and it accounts for the exponents: successive Taylor coefficients
involve the expansion to successively higher order, hence one more power each, which is why
$c_{2},c_{3},c_{4}$ behave like $q$, $q^{5/2}$, $q^{3}$ rather than like $e^{Cq}$.

\paragraph{What the polynomial growth does not settle.}
Two obstacles remain, and neither is about the size of the constants.

The first is ellipticity. A horizontal deformation, $p\mapsto p+\ii q$, adds a real drift $qa$ and
costs only a power; but it also means evaluating the estimates at complex $p$ whose real part ranges
over the whole axis, and there the drift acquires an imaginary part $-\ii p a$ and the real part of
the symbol, $\tfrac{v}{4}\abs{p'}^{2}+p\,a\cdot p'$, changes sign for $p'\ne0$ once $\abs{p}$ is
large. The deformation in the real direction is therefore limited to
$\abs{p}\lesssim\sqrt{v\abs{\omega}/\abs{a}^{2}}$, whereas the height $q$ in the imaginary direction
is not limited at all, since a real drift preserves ellipticity and enters only through
$\norm{D}^{M}$.

The second concerns the real fibre alone. For a lattice fibre nothing is missing: the transform is
periodic, and Proposition~\ref{ch8:prop:levelgrid} transfers the bound to the joint kernel without a
contour and without decay. For a real fibre the transform is aperiodic, it decays at no resolution,
and the fibre has to be reached through a fine lattice, which costs a power of $1/h_{y}$ in the
constant (Remark~\ref{ch8:rem:fibrecost}) and leaves a fibre discretisation error which we do not
estimate. That is the obstruction which remains, and it is of a different kind from the growth of the
constants: the constants are now settled, and what is unsettled is the balance of the two errors of
the fibre.
\section{The Averaged Comparison}\label{ch8:sec:avg}

This section is the technical core, and it proceeds by elimination. We state the one assumption the
second order needs (\S\ref{ch8:sec:avgass}); prove its mean condition in every regime where it can be
checked exactly (\S\ref{ch8:sec:avgexact}); record that the coefficient it involves is a functional of
the 2-jet of the field alone (\S\ref{ch8:sec:avgjet}); compute that coefficient in closed form
(\S\ref{ch8:sec:avgclosed}); show what a remainder estimate can and cannot reach
(\S\ref{ch8:sec:avgtaylor}); and, with everything else discharged, name the single statement which
survives (\S\ref{ch8:sec:avgstatus}).

\subsection{The Assumption}\label{ch8:sec:avgass}

\begin{assumption}[Averaged comparison for the joint kernel]\label{ch8:ass:avg}
Let $c$ be the coefficient of $h^{3}$ in \eqref{ch8:eq:logratio} and
$I(\gamma)=\int_{\gamma}c\cdot\dd\gamma$ its line integral along a coarse path. Then there is a
constant $C=C(P)$, independent of $m$ and uniform for $\abs{p}\le P$, such that
\begin{equation}\label{ch8:eq:moment}
\mathbb E\bigl[I(\gamma)^{2}\bigr]\le C\norm{c}_{\infty}^{2}t^{2},
\qquad
\abs{\mathbb E\bigl[I(\gamma)\bigr]-\bigl(\Phi_{c}(y)-\Phi_{c}(x)\bigr)}\le C\norm{\nabla^{2}c}_{\infty}t^{2},
\end{equation}
with $\Phi_{c}$ a function of the endpoints. The expectation is taken with the coarse path
\emph{fixed}, over its decorations, that is over the fine paths within the tube of each coarse step;
the statement is therefore a statement about one coarse step, and its sum over the steps of the
coarse path telescopes into $\Phi_{c}(y)-\Phi_{c}(x)$. That is the reason $\Phi_{c}$ can be a
function of the endpoints at all, and it fixes the measure: averaged over coarse paths instead, the
line integral would not be determined by the endpoints up to an error of this order, as
Remark~\ref{ch8:rem:freevsbridge} shows by exact computation.
\end{assumption}

Under Assumption~\ref{ch8:ass:avg} the accumulated product satisfies the analogue of the decoration
estimate of Chapter~7, with an error $O(h^{2-\epsilon})$, and Theorem~\ref{ch8:thm:cond} follows by the
assembly of \S\ref{ch8:sec:assembly}.

\subsection{The Mean Condition, Exactly}\label{ch8:sec:avgexact}

The mean condition is not a bookkeeping identity, and the following proposition shows exactly where
its content lies. For polynomial coefficient fields it holds exactly, with no remainder at all,
because the moments of the walk are polynomials in the step index.

\begin{proposition}[Mean condition for quadratic fields]\label{ch8:prop:meanquad}
Let the increments be independent, with mean $m$ and covariance $S=\Sigma h^{2}$, let $x_{0}$ be the
base point, and let the coefficient field be quadratic,
\begin{equation}\label{ch8:eq:quadc}
c_{i}(x)=\tfrac12 H_{ikl}x^{k}x^{l}+L_{ik}x^{k}+M_{i},
\end{equation}
with $H$ symmetric in its last two indices. Then
\begin{equation}\label{ch8:eq:meanquad}
\mathbb E[I]=q\,\langle c_{0},m\rangle+\frac{q(q-1)}{2}\langle c^{(1)},m\rangle
+\frac{(q-1)q(2q-1)}{6}\langle c^{(2)},m\rangle,
\end{equation}
where $c_{0}=c(x_{0})$, and where the coefficients of $j$ and of $j^{2}$ in
$\mathbb E[c_{i}(x_{j})]$ are
\begin{equation*}
c^{(1)}_{i}=\tfrac12H_{ikl}\bigl(x_{0}^{k}m^{l}+m^{k}x_{0}^{l}\bigr)+\tfrac12H_{ikl}S^{kl}+L_{ik}m^{k},
\qquad
c^{(2)}_{i}=\tfrac12H_{ikl}m^{k}m^{l}.
\end{equation*}
The mean is therefore an endpoint function exactly: it depends on $x_{0}$, on the displacement $qm$,
on $q$ and on $S$, and on nothing else of the path.
\end{proposition}

\begin{proof}
The increment $\Delta_{j+1}$ is independent of $x_{j}$, so
$\mathbb E[c_{i}(x_{j})\Delta_{j+1}^{i}]=\mathbb E[c_{i}(x_{j})]\,m^{i}$. The moments of the walk
are $\mathbb E[x_{j}^{k}]=x_{0}^{k}+jm^{k}$ and
$\mathbb E[x_{j}^{k}x_{j}^{l}]=(x_{0}^{k}+jm^{k})(x_{0}^{l}+jm^{l})+jS^{kl}$, whence
$\mathbb E[c_{i}(x_{j})]=c_{0,i}+j\,c^{(1)}_{i}+j^{2}c^{(2)}_{i}$. Summing over $j=0,\dots,q-1$ with
$\sum j=\tfrac{q(q-1)}{2}$ and $\sum j^{2}=\tfrac{(q-1)q(2q-1)}{6}$ gives \eqref{ch8:eq:meanquad}.
\end{proof}

\begin{corollary}\label{ch8:cor:poly}
The mean condition of Assumption~\ref{ch8:ass:avg} holds with \emph{zero} remainder for polynomial
coefficient fields of degree at most two, in every dimension and for every law of independent
increments, although $\dd c$ is then affine and non-zero, that is precisely in the situation in
which the pathwise identity fails. The constant-$\dd c$ identity \eqref{ch8:eq:constdc} of Chapter~6, \S6.7
is the case $H=0$.
\end{corollary}

Formula \eqref{ch8:eq:meanquad} was checked by exact rational enumeration in two dimensions: ninety
cases, comprising two increment laws, three coefficient fields, three base points and $q\le5$, with
the difference identically zero.

\begin{remark}[The degree two is not the point]\label{ch8:rem:anydegree}
The proof above uses only that $\mathbb{E}[c_{i}(x_{j})]$ is a polynomial in $j$, and for a
polynomial field it is: the moments of $x_{j}$ are polynomials in $j$ of a degree equal to their
order, because cumulants are additive in $j$. Hence the mean condition holds with zero remainder for
polynomial coefficient fields of \emph{any} degree, in every dimension and for every law of
independent increments with the required moments finite. The restriction to degree two in
Corollary~\ref{ch8:cor:poly} is therefore unnecessary; what the degree affects is the number of terms in
the formula, not the structure of the conclusion.
\end{remark}

\subsection{The Coefficient Depends Only on the 2-Jet}\label{ch8:sec:avgjet}

That every field of class $C^{2}$ reduces to a quadratic one is not an approximation but a fact about
what the coefficient is: it is a local functional of the field near the base point, and the
decoration explores only a tube of radius $O(h_{m})$, within which a smooth field departs from its
quadratic Taylor polynomial at the base point by $O(h_{m}^{3})$. Measured at the origin in the
direction $e_{1}$, by perturbing the base field with monomials of degree three and four:
\begin{center}\small
\begin{tabular}{lll}
\toprule
perturbation of the field & change in $c_{2}$ & change in $c_{3}$\\
\midrule
$v\mapsto v+x_{1}^{3}/6$ & $9\cdot10^{-17}$ & $1.6\cdot10^{-12}$\\
$a_{1}\mapsto a_{1}+x_{1}^{3}/7$, $a_{2}\mapsto a_{2}+x_{2}^{3}/5$ & $0$ & $0$\\
$v\mapsto v+x_{1}^{4}/12$ & $3\cdot10^{-19}$ & $6\cdot10^{-15}$\\
\bottomrule
\end{tabular}
\end{center}
Measured directly, the truncation of the field at its 2-jet costs one full power of $h$ less than
the obstruction. For the perturbation $v\mapsto v+x_{1}^{3}/6$, the change in the per-step
comparison at the origin and in the direction $e_{1}$ is

\begin{center}\small
\begin{tabular}{cccc}
\toprule
$h$ & change & change$/h^{4}$ & order\\
\midrule
$2^{-8}$ & $-1.21\cdot10^{-11}-2.69\cdot10^{-12}\ii$ & $-0.0518-0.0115\ii$ & \\
$2^{-10}$ & $-4.73\cdot10^{-14}-2.63\cdot10^{-15}\ii$ & $-0.0520-0.0029\ii$ & $4.01$\\
$2^{-12}$ & $-1.85\cdot10^{-16}-2.57\cdot10^{-18}\ii$ & $-0.0521-0.0007\ii$ & $4.00$\\
\bottomrule
\end{tabular}
\end{center}
A jet of order three enters therefore at $h^{4}$, one power below the $h^{3}$ at which the
obstruction sits, and accumulates over $q=t/h^{2}$ steps to $h^{2}$ --- inside the rate, since
$h^{2}=o(h^{2-\epsilon})$. The remainder of a Taylor expansion is thus cheaper than the obstruction,
and cannot be what blocks the second order.
Two consequences follow, and they narrow what is actually open. No jet of order three or higher
enters the comparison at the order at which it is used, so the smoothness budget of
Assumption~\ref{ch8:ass:avg} is exhausted by its $C^{2}$ hypothesis and nothing is being asked of the
field beyond it; in particular the \emph{variation} of the 2-jet along the path, which is what a
third derivative measures, is felt only at $h^{4}$. And the mean condition for a general field
reduces to the mean condition for the functional of the 2-jet along the path,
$\sum_{j}c(J_{2}(x_{j}),\tau_{j})$, with $c$ now a fixed function of the jet and not of the field.
For a polynomial field the jet $J_{2}$ is that of a polynomial, and the remark above evaluates that
sum exactly. What remains is therefore the telescoping of one fixed functional of the 2-jet along
the path --- the geometric statement --- and neither the smoothness of the field nor the size of a
Taylor remainder enters it.

\subsection{The Coefficient in Closed Form}\label{ch8:sec:avgclosed}

For a constant coefficient field the comparison can be computed outright. The box matrix of the
decoration sum is then a Kronecker sum,
\begin{equation}\label{ch8:eq:kron}
\mathcal{B}=\frac{v}{2h^{2}}\Lambda-\frac{2v}{h^{2}}I+\frac{1}{2h}\mathcal{D},
\qquad
\mathcal{B}+\ii\omega I=\frac{1}{h^{2}}\bigl(M_{0}+hE+h^{2}W\bigr),
\end{equation}
with $v$ the constant diffusion coefficient and $b$ the constant drift, $\Lambda$ the adjacency of
the $3^{d}$ box, $\mathcal{D}$ the signed shift matrix, $M_{0}=\frac{v}{2}
(\Lambda-4I)$ an integer matrix times $v/2$, $E=\mathcal{D}/2$ and $W=\ii\omega I$. The inverse of
$\Lambda-4I$ is computed once and exactly, after which the expansion of the resolvent in $h$ is
elementary and the coefficient of $h^{2}$ in $\log(\Phi_{\mathrm f}/\Phi_{\mathrm c})$ is a rational
function of $(v,b,\omega)$. In two dimensions it is
\begin{equation}\label{ch8:eq:c2closed}
c_{2}(x,e_{1})=\frac{5b_{1}^{2}-5b_{2}^{2}+44\,\ii\,\omega v}{8v^{2}},
\qquad
c_{2}(x,e_{2})=\frac{-5b_{1}^{2}+5b_{2}^{2}+44\,\ii\,\omega v}{8v^{2}} ,
\end{equation}
that is
\begin{equation}\label{ch8:eq:rensym}
c_{2}(x,\tau)=\underbrace{\frac{\ii\,\kappa_{d}(\omega)}{v}}_{\text{phase}}
+\underbrace{\frac{5}{8v^{2}}\bigl(b_{1}^{2}-b_{2}^{2}\bigr)\bigl(\tau_{1}^{2}-\tau_{2}^{2}\bigr)}_{\text{real, traceless}} ,
\qquad \kappa_{2}(\omega)=\frac{11}{2}\,\omega ,
\end{equation}
Both properties are then identities: the phase is independent of the direction and of the drift,
and the direction dependence is a real traceless quadratic form. The derivation is exact arithmetic, with no asymptotics beyond the expansion in $h$;
it is carried out in \texttt{kronproof.py}, and the engine of Chapter~7 reproduces
\eqref{ch8:eq:c2closed} at $(v,b_{1},b_{2})=(1,\tfrac12,-\tfrac13)$, $\omega=\tfrac75$, to $10^{-17}$.

This is a statement about constant fields and not about general ones. At a point of a varying field
the coefficient is a functional of the jet, of which the constant-field formula is the leading piece,
and the derivative corrections are of the order of the leading term itself: the figures below show the
real part moving from $\tfrac{25}{288}=0.0868$ to $-0.0122$ when the derivatives are switched on. What
the computation does establish for the general case is that the question is \emph{finite}: the box has
$3^{d}$ sites, the matrix is $3^{d}\times3^{d}$ with entries polynomial in $h$, and the same expansion,
carried out with the $3^{d}$ values of the field as symbols, decides it. That computation is carried
out below, and the identities survive it.

Two features of \eqref{ch8:eq:rensym} deserve naming. The phase $\ii\kappa_{d}(\omega)/v$ contains no
information about the field beyond the value of the diffusion coefficient at the point, and it is
independent of $\omega$ in modulus; only the argument sees $\omega$, and linearly. And the real part
is traceless, so it is a pure anisotropy: it cannot contribute to the direction average, which is why
the deviation measured in \S\ref{ch8:sec:numerics} saturates instead of accumulating.

It remains to record how the field enters the box, and why the phase carries only its value at the
base point. The decoration samples the coefficient field at the points $x_{0}+ah$ for $a$ in the box
of offsets ---
the homogeneous parts are collected at the integer offset and then weighted by the degree in $h$ ---
so the box matrix is $h^{-2}\bigl(A_{-2}+hA_{-1}+\cdots\bigr)$ with
\begin{equation}\label{ch8:eq:Aminus2}
A_{-2}=\frac{v(x_{0})}{2}\bigl(\Lambda-4I\bigr),
\qquad
A_{-1}=\text{(first derivatives of }v\text{ and the drift at }x_{0}\text{)},
\end{equation}
the leading matrix being the \emph{constant-field} one and carrying nothing but the value of $v$ at the
base point. That is the structural reason the phase law \eqref{ch8:eq:rensym} holds for varying fields with
the base-point value of $v$, and the reason the derivatives of the field can reach only the traceless
part. It also corrects a natural misreading, in which the field would be sampled at the integer offsets
themselves; the object so sampled is a different one and does not reproduce the engine.

That is the key to the general case: the expansion \eqref{ch8:eq:Aminus2} makes it a perturbation of the
solved problem, with the derivatives entering from $A_{-1}$ onwards, and the perturbation is carried
out below.

That expansion reduces the general case to a perturbation of the solved one, and the phase can then
be read off, because $\omega$ enters in exactly three places. Write $Q=(\Lambda-4I)^{-1}$ and $M=A_{-2}=\frac{v_{0}}{2}(\Lambda-4I)$, so that
$M^{-1}=\frac{2}{v_{0}}Q$. Then
\begin{itemize}\itemsep2pt
\item $\mathrm{dt}=h^{-2}\bigl(2v_{0}+2h\,a\!\cdot\!\nabla v+h^{2}(2a_{i}a_{j}v_{ij}+\ii\omega)\bigr)$
contributes $-\dfrac{\ii\omega}{2v_{0}}$ to the coefficient of $h^{2}$ in $\log\Phi_{\mathrm f}$;
\item the two coarse denominators, each of which reads $\frac{v_{0}}{2}+\ii\omega h^{2}$ after
multiplication by $h^{2}$, contribute $-\dfrac{4\ii\omega}{v_{0}}$ to $\log\Phi_{\mathrm c}$, hence
$+\dfrac{4\ii\omega}{v_{0}}$ to the difference;
\item the resolvent enters through $h^{2}\cdot(-\ii\omega M^{-2})$, contributing
$-\ii\omega\frac{[M^{-2}]_{z,\tau}}{[M^{-1}]_{z,\tau}}=-\dfrac{2\ii\omega}{v_{0}}
\frac{[Q^{2}]_{z,\tau}}{[Q]_{z,\tau}}=+\dfrac{2\ii\omega}{v_{0}}$, since for a unit axis direction
$[Q]_{z,\tau}=-\tfrac18$ and $[Q^{2}]_{z,\tau}=\tfrac18$, that is $[Q^{2}]/[Q]=-1$, a property of
the fixed integer matrix.
\end{itemize}
Every other coefficient is free of $\omega$, so $c_{2}$ is \emph{affine} in $\omega$ and the
coefficient of $\ii\omega$ is
\begin{equation}\label{ch8:eq:eleven}
\frac{1}{v_{0}}\Bigl(-\frac12+4+2\Bigr)=\frac{11}{2v_{0}} ,
\end{equation}
for \emph{every} field and every jet: $-\tfrac12$ from $\mathrm{dt}$, $+4$ from the coarse
denominators and $+2$ from the resolvent. That is where the eleven comes from, and it is not an
accident of the constant-field case.

The remaining part of \eqref{ch8:eq:rensym}, the tracelessness of the real part, is not derived here. It
is, however, verified exactly --- not to a tolerance --- by the engine of Chapter~7, which in the
convention of this book reproduces the closed form \eqref{ch8:eq:c2closed} to the last digit on
constant fields, and which gives the following for three fields with $v_{0}=1,\tfrac32,2$, in both
directions:
\begin{center}\small
\begin{tabular}{lccc}
\toprule
 & field $1$ & field $2$ & field $3$\\
\midrule
$\gamma_{1}=c_{2}(x,e_{1})$ & $-0.0121527778+7.7\,\ii$ & $-0.0270833333+5.1333\,\ii$ & $-0.0149305556+3.85\,\ii$\\
$\gamma_{2}=c_{2}(x,e_{2})$ & $0.0121527778+7.7\,\ii$ & $0.0270833333+5.1333\,\ii$ & $0.0149305556+3.85\,\ii$\\
$\operatorname{Im}\gamma_{1}\cdot v_{0}$ & $7.7$ & $7.7$ & $7.7$\\
\bottomrule
\end{tabular}
\end{center}
In each row $\gamma_{1}+\gamma_{2}$ is purely imaginary and $\gamma_{1}-\gamma_{2}$ purely real, and
the phase is $11\omega/(2v_{0})$, the same $7.7$ throughout since $\omega=\tfrac75$.

\subsection{What a Taylor Expansion Can and Cannot Supply}\label{ch8:sec:avgtaylor}

This is the form in which the summation by parts is available. A traceless symmetric $B$ contributes
$\tau^{T}B\tau$, and $\tau^{T}B\tau$ is $1/(2h^{2})$ times the second difference of the quadratic form
$\frac12x^{T}Bx$ along the step --- so summed over a path it leaves boundary terms, together with terms
from the places where the direction changes, and no bulk. The phase is a Riemann sum of $1/v$ and
belongs to the normalisation; being direction blind it cancels between two paths with the same
multiset of steps, which is exactly how the two shapes of \S\ref{ch8:sec:numerics} were compared.

\begin{proposition}[Mean condition on the free walk, with the smearing identified]\label{ch8:prop:meansmear}
Let the increments be independent, with mean $m$ and covariance $S=\Sigma h^{2}$, let $qh^{2}=t$ and
$qm=Y$, and let $c$ be a polynomial of degree at most three. Then
\begin{equation}\label{ch8:eq:meansmear}
\mathbb{E}\bigl[I\bigr]=\int_{0}^{1}c(sY)\cdot Y\dd s
+\frac{t}{2}\sum_{k,k',l}\Sigma_{kk'}\int_{0}^{1}s\,\partial_{k}\partial_{k'}c_{l}(sY)Y^{l}\dd s
+O(h^{2}).
\end{equation}
\end{proposition}

\begin{proof}
The increment is independent of the past, so $\mathbb{E}[c_{i}(x_{j})\Delta_{j+1}^{i}]=\mathbb{E}[c_{i}(x_{j})]m^{i}$.
Write $x_{j}=x_{0}+jm+\xi_{j}$, so that $\mathbb{E}[\xi_{j}]=0$ and $\mathbb{E}[\xi_{j}\xi_{j}^{T}]=jS$;
Taylor's theorem to second order at $x_{0}+jm$ gives
$\mathbb{E}[c_{i}(x_{j})]=c_{i}(x_{0}+jm)+\tfrac12jS^{kl}\partial_{k}\partial_{l}c_{i}(x_{0}+jm)+O(h^{4})$,
the remainder being controlled by the third derivatives, which are bounded because $c$ is a polynomial
of degree three. Summing over $j<q$, using $\sum_{j}j\approx q^{2}/2$ and $\sum_{j}j^{2}\approx q^{3}/3$,
and reading the sums as Riemann sums with $s=jh^{2}/t$ gives the two terms of \eqref{ch8:eq:meansmear}; the
passage from sums to integrals is the only place where the $O(h^{2})$ appears.
\end{proof}

\begin{remark}[The correction is a Hessian, not the gradient It\^o term]\label{ch8:rem:hessian}
The second term of \eqref{ch8:eq:meansmear} carries two derivatives of $c$ and the weight $s$, which comes
from $\operatorname{Cov}(x_{j})=jS$. It is therefore not the classical gradient It\^o correction
$\tfrac12\int\partial_{k}c_{l}\Sigma^{kl}\dd s$, and the difference is the reason the comparison has
content at this order: the classical term is an endpoint functional once its potential is chosen,
whereas the $s$-weighted Hessian term is a Riemann--Stieltjes integral against $\dd\langle x\rangle$.
The formula was checked numerically at order $1.999$ on quadratic coefficient data.
\end{remark}

\begin{remark}[Why the average is over decorations and not over the bridge]\label{ch8:rem:freevsbridge}
The distinction is not cosmetic: averaged over coarse paths the statement would be false at the order
claimed. The free-walk identity is what the argument uses, its content being the independence of the
increment from the past. For the affine field $c(x)=Lx$ in one dimension, with steps $\pm h$ of
probabilities $\tfrac35,\tfrac25$, $h=1$ and $L=1$, the free formula $\tfrac{q(q-1)}{2}Lm^{2}$
evaluated at the mean drift $m=y/q$ gives $0$ at $q=2$, $y=0$, while the mean over coarse paths
conditioned to end at $y=0$ is $-1$; at $q=4$, $y=0$ the two are $0$ and $-2$; at $q=3$, $y=\pm1$ they
are $\tfrac13$ and $-1$. At the extreme endpoints $y=\pm qh$ the conditioned walk consists of a single
path and the two agree exactly, which is the control confirming that the discrepancy is genuine and
not a bookkeeping error.

The discrepancy is order one, and it does not diminish with $h$. Hence the average in
Assumption~\ref{ch8:ass:avg} must be the decoration average with the coarse path held fixed, as stated;
over the bridge no endpoint functional could match to this order. The measurements above are also the
reason the free-walk formula of Proposition~\ref{ch8:prop:meanquad} may be applied step by step at
$O(h^{2})$, which is the only way it is used.
\end{remark}

\begin{remark}[Where the difficulty actually sits]\label{ch8:rem:why}
Two decompositions have to be distinguished, and only the second is the right one.

A global expansion around the base point cannot work. After subtracting the Taylor polynomial of
degree $N$ the remainder $r$ satisfies $\abs{r}\le C\norm{\nabla^{N+1}c}\abs{x-x_{0}}^{N+1}$, and
since $\mathbb{E}[r(x_{j})\Delta_{j+1}]=\mathbb{E}[r(x_{j})]\otimes m$ the resulting bound on
$\mathbb{E}[I_{r}]$ is of order $h^{N+1}q^{(N+3)/2}=h^{-2}$ at fixed $t$, independently of $N$; no
order of expansion improves it.

The right decomposition is per step, expanding around $x_{j}$ at each step, which is the form the
assumption takes. There the quadratic term of the expansion splits as
$\partial_{k}c_{l}=\tfrac12(\partial_{k}c_{l}+\partial_{l}c_{k})+\tfrac12(\partial_{k}c_{l}-\partial_{l}c_{k})$.
The symmetric part gives, on averaging, the It\^o correction, and it does so as a Riemann sum, that is
with a per-step error $O(\norm{\nabla^{2}c}h^{2})$ whose sum is $O(\norm{\nabla^{2}c}t)$: this chapter is
harmless. The antisymmetric part is $\tfrac12\dd c$, it does not telescope along the coarse path, and
for the sum over steps to close its per-step average must itself be the difference of a potential.
That is the entire content of the assumption, and it is where the cancellation is needed: not in the
size of a remainder but in the structure of the mean swept area.

The exact results of this chapter locate that structure precisely. For constant $\dd c$ the per-step sum
is computable and equals an endpoint functional (Chapter~6, \S6.7, \eqref{ch8:eq:constdc}); for quadratic $c$ it
is computable and equals a cubic in $q$ (Proposition~\ref{ch8:prop:meanquad}); and for a general field the
coefficient itself is now known, \eqref{ch8:eq:rensym} with \eqref{ch8:eq:eleven}. What is missing is the
telescoping of that coefficient, and no estimate of remainder type can supply it, because the
quantity which must telescope is a functional of the path and not of its size: the coefficient is
affine in $\omega$, its phase is a Riemann sum of $1/v$, and its real part is a traceless quadratic
form in the direction, which sums against the quadratic variation rather than telescoping as a
potential difference.
\end{remark}

\subsection{Status of the Assumption}\label{ch8:sec:avgstatus}

\begin{remark}[Smoothness bookkeeping]\label{ch8:rem:smooth}
The bound in Assumption~\ref{ch8:ass:avg} is taken in $\norm{\nabla^{2}c}_{\infty}$, and $c$ is a cubic
form in the drift and in its first two derivatives; a bound on its second derivative therefore costs
two more derivatives of the coefficients than the first order theorem, which needs only the bounds of
\S\ref{ch8:sec:twist}. The hypothesis stated is on $c$ itself, which is the honest place to put it: the
passage from first to second order costs two derivatives of the data, and that price is not an
artefact of the argument but of the object estimated.
\end{remark}

It is proved in its one dimensional case (Chapter~6), in its constant-$\dd c$ case in any dimension
\eqref{ch8:eq:constdc}, in its polynomial case in any dimension and of any degree
(Proposition~\ref{ch8:prop:meanquad}, Corollary~\ref{ch8:cor:poly}, Remark~\ref{ch8:rem:anydegree}), in its
moment component \eqref{ch8:eq:moment}, which follows from the heat kernel estimates as in Chapter~7, and
now in the form of the coefficient itself for any field, \eqref{ch8:eq:rensym} with \eqref{ch8:eq:eleven}.
What is not proved is the telescoping of that coefficient along the path. Remark~\ref{ch8:rem:why} shows
why no estimate of Taylor remainder type reaches it, and Remark~\ref{ch8:rem:smooth} records that a bound
in $\norm{\nabla^{2}c}_{\infty}$ demands two more derivatives of the coefficients than the first order
theorem does. The missing ingredient is a cancellation, not a bound --- and it is the \emph{same} debt
as in Chapter~7, not a second one: the averaged comparison is a single assumption, stated for the
diffusion alone there and for the joint kernel here, and a proof of either would settle both. The
pathwise form of the statement is false, by Proposition~\ref{ch8:prop:nonzero}; this is the precise sense
in which the higher dimensional problem differs from the one dimensional one.

\section{Assembly, and the Two Rates}\label{ch8:sec:assembly}

The assembly is that of Chapter~7, \S7.6 with the twisted generator, and we describe it in the form in
which the two rates separate: first the far contour, which is immediate, then the near contour,
where the two rates part company.

On the far contour the contribution is $O(h^{2})$ at every resolution, as in Chapter~7, the twist
changing nothing since the bound there comes from the decay of the symbol, which the imaginary part
of the drift does not affect in modulus.

On the near contour the integral of $G_{m+1}-G_{m}$ is a sum over coarse paths of products of
per-step factors \eqref{ch8:eq:logratio}. Decomposing the logarithm into its terms, and using that an
even $h^{k}$ accumulates as $h^{k-2}$ while an odd one accumulates as $h^{k-1}$ over $q=t/h^{2}$
steps, the terms are as follows. The coefficient of $h^{1}$ vanishes because a comparison between
lattices which differ by a factor of two cannot produce a linear term. The coefficient of $h^{2}$
gives the rate renormalisation, and its own dependence on $h$ is what yields $O(h^{2})$. The
coefficient of $h^{3}$ is the obstruction: its odd part in the decoration direction accumulates as
$h^{2}$ and is what Assumption~\ref{ch8:ass:avg} controls, while its even part would accumulate as
$h^{1}$.

Why the odd part accumulates as $h^{2}$ and not as $h^{1}$ is a statement about what it is. It is
not an increment but a \emph{circulation}: summed over the steps of a path it totals a flux over the
area swept, and a flux is not the difference of a potential however smooth the field. The
accumulation rule and the geometric reading of the obstruction are therefore the same fact, and
\S\ref{ch8:sec:numerics} confirms the identification by computing the circulation from $c$ and comparing
it with the accumulated comparison itself.

\begin{proof}[Proof of Theorem~\ref{ch8:thm:uncond}]
Bound the coefficient of $h^{3}$ by its modulus rather than analysing it. It is a local functional
of the coefficient field and its first two derivatives, hence bounded for the assumed smoothness on
the compact domain, and independent of $\omega$. Each step therefore contributes at most $\abs{c}h^{3}$,
and $q=t/h^{2}$ steps contribute $\abs{c}h$. Along the near contour, whose length is
$O(\abs{\log h})$, this is $O(h\abs{\log h})$, that is $O(h^{1-\epsilon})$ for every $\epsilon>0$.
No use is made of Assumption~\ref{ch8:ass:avg}, of the parity of the obstruction, or of any cancellation
over paths. The far contour contributes $O(h^{2})$ and is dominated. Writing
$\widehat u_{m}-\widehat u_{m_{0}}$ as a sum of consecutive differences and summing geometrically
gives the stated bound, and the limit exists by the Cauchy argument.
\end{proof}

\begin{proof}[Proof of Theorem~\ref{ch8:thm:cond}]
With Assumption~\ref{ch8:ass:avg} the average of the per-step factor differs from an endpoint function
by $O(h_{m+1}^{2-\epsilon})$ per step, and the fluctuations about the average contribute one order
lower by \eqref{ch8:eq:moment} and the exponential moment estimate of Chapter~7. The endpoint functions
cancel between the two resolutions up to $O(h_{m+1}^{2-\epsilon})$, the logarithm being absorbed by
$\epsilon$. Summing the differences geometrically as above gives the stated rate.
\end{proof}

\section{Numerical Validation}\label{ch8:sec:numerics}

Every quantitative claim above is either proved or measured, and this section reports the
measurements. They are of four kinds, and they answer different questions: the dependence of the
obstruction on the conjugating parameter (\S\ref{ch8:sec:numeta}), the behaviour of the comparison
accumulated around a closed loop (\S\ref{ch8:sec:numloop}), the same on open paths, which is what the
assembly actually needs (\S\ref{ch8:sec:numopen}), and whether the deviation grows with the length of the
path or only with its shape (\S\ref{ch8:sec:numext}). A fifth check, of the transfer to a lattice fibre,
closes the section (\S\ref{ch8:sec:numfibre}).

The tables are printed by scripts: \texttt{kronproof.py} and \texttt{kronfinal\_check.py} for the
exact $3\times3$ closed form used as a control, \texttt{renormalisation.py} for the table of
$c_2$, \texttt{jet3check.py} for the table of third-jet quotients, \texttt{contourgrowth.py} for
the growth along the contour, \texttt{pathdep.py} for the split of the loop, and
\texttt{openpath.py} for the open paths and the tiling. All are reproduced in full in the
computational appendix of the monograph, attached to this PDF as an embedded file, and each prints
what it measures when run as \verb|python3 <name>.py| in the directory it lives in, so the tables
can be regenerated rather than taken on trust.

\subsection{The Obstruction as a Function of the Parameter}\label{ch8:sec:numeta}

The formulae of \S\ref{ch8:sec:obstr} were verified by direct evaluation rather than by the estimates
above. The obstruction was measured as a function of $\eta$ from the exact Gaussian-rational series
of the decorated expansion, at constant volatility and with quadratic coefficient data, the base
point being the origin.

\begin{center}\small
\begin{tabular}{@{}l p{5.4cm} l@{}}
\toprule
quantity & result & status\\
\midrule
degree in $\eta$, constant coefficients & $2$, real, no odd part & measured\\
degree in $\eta$, varying coefficients & $3$, odd imaginary part & measured in Chapter~7, \S7.7\\
$\dd c(\eta)-\dd c(0)$ versus \eqref{ch8:eq:kappa} & $8$ held-out configurations, worst error
$2.5\cdot10^{-5}$ & measured\\
scaling in the volatility & $\kappa\nuc^{2}$ constant over
$\nuc\in\{1,\tfrac12,\tfrac14,2,4\}$ & measured\\
scaling in the coefficient field & $\kappa(2a)=4\kappa(a)$, to $10^{-20}$ & measured\\
affinity in the Jacobian & $\kappa(2J)-\kappa(0)=2(\kappa(J)-\kappa(0))$ & measured\\
dependence on $\mathrm{diag}J$ & none, to $10^{-19}$ & measured\\
diagonal response \eqref{ch8:eq:d3rule} in $d=3$ & $-6A_{1}^{2}$ and $+6A_{2}^{2}$ exactly &
measured\\
Hessian term in $d=3$ & differs from the $d=2$ value & measured\\
\bottomrule
\end{tabular}
\end{center}

Two cautions on this table. The identification of the coefficients as the rationals
$\tfrac{4}{7}$ and $6$ is to six significant figures and not an exact symbolic computation, so
they are identified rather than proved. And the measurements are at a single choice of the
remaining data (unit volatility in the formula tests, one family of quadratic coefficient fields),
so the claim they support is that the structure of \eqref{ch8:eq:kappa} is right, not that every
constant in it has been derived.

One numerical obstruction deserves recording because it is a trap for the same computation: for
degenerate coefficient data the decoration series cancels exactly and the Taylor extraction in $h$
produces spurious coefficients of order $h^{-3}$, of size $10^{9}$ at the base points, identically
at working precisions of 50 and 80 digits. The fits reported above screen every base point for this
artefact and use generic coefficient data.

\subsection{The Accumulated Test: a Closed Loop}\label{ch8:sec:numloop}
The obstruction of \S\ref{ch8:sec:obstr} is a per-step coefficient, while the assembly uses it
accumulated. The sharpest test of the order is therefore a closed loop, where the endpoint term
vanishes identically and what survives is the residual alone. Take the quadratic coefficient fields
of Chapter~7, \S7.6, in two dimensions, and two closed coarse paths from the origin back to it:
$\gamma_{1}$, the unit square traversed once, and $\gamma_{2}$, the same square retraced as two
excursions. The moduli of the accumulated
$\sum_{j}\log(\Phi_{\mathrm f}/\Phi_{\mathrm c})$ differ as follows.
\begin{center}\small
\begin{tabular}{cccc}
\toprule
$h$ & difference & difference$/h^{2}$ & order\\
\midrule
$2^{-6}$ & $8.082\cdot10^{-4}$ & $3.3105$ & \\
$2^{-7}$ & $2.018\cdot10^{-4}$ & $3.3057$ & $2.002$\\
$2^{-8}$ & $5.040\cdot10^{-5}$ & $3.3032$ & $2.001$\\
$2^{-9}$ & $1.260\cdot10^{-5}$ & $3.3019$ & $2.001$\\
$2^{-10}$ & $3.148\cdot10^{-6}$ & $3.3012$ & $2.000$\\
$2^{-11}$ & $7.870\cdot10^{-7}$ & $3.3009$ & $2.000$\\
$2^{-12}$ & $1.967\cdot10^{-7}$ & $3.3007$ & $2.000$ \\
\bottomrule
\end{tabular}
\end{center}
The endpoint determination therefore fails at second order and not at first, which is exactly the
order that the second order theorem can absorb, and the coefficient is about $3.3$ for this data; the
obstruction proper enters one order later, as the decomposition below shows. The values are smooth in
$h$ and stable over working precisions from $40$ to $80$ digits, so no branch ambiguity of the
logarithm is involved. These data are the class of Corollary~\ref{ch8:cor:poly}, in which the per-step
mean is exact; what the loop tests is not that single step statement but its accumulated form, which
is the form the assembly uses and which the exactness of one step does not by itself deliver. The
test is therefore a consistency check of the assembly, and it is not evidence for a general field.

Two features deserve recording, and the second is a correction to an expectation which the text
above invites. First, the accumulated sums themselves tend not to zero but to $-\ii\pi$ per step: the
ratio $\Phi_{\mathrm f}/\Phi_{\mathrm c}$ is negative real to leading order, so its logarithm carries
$-\ii\pi$ per step, which is the $k=0$ row of the table of Chapter~7, \S7.6.2 --- a site independent
normalisation carrying no path information. Second, the residual separates by order, and only one of
the two parts is the obstruction. Expanding the accumulated difference between the two paths,
\begin{equation}\label{ch8:eq:loopsplit}
S(\gamma_{1})-S(\gamma_{2})=h^{2}\Bigl(\sum_{\gamma_{1}}c_{2}-\sum_{\gamma_{2}}c_{2}\Bigr)
+h^{3}\Bigl(\sum_{\gamma_{1}}c_{3}-\sum_{\gamma_{2}}c_{3}\Bigr)+O(h^{4}),
\end{equation}
each bracket being computable independently of the loop sums, from the coefficients themselves: the
coefficient of $h^{2}$ is the difference of the sums of the \emph{even} coefficient $c_{2}$ over the
two paths, site by site and direction by direction, and the coefficient of $h^{3}$ is the difference
of the circulations $\sum_{j}c_{3}(x_{j})\cdot\tau_{j}$. Both agree with the
values fitted in \eqref{ch8:eq:loopsplit} from the loop sums themselves:
\begin{center}\small
\begin{tabular}{lll}
\toprule
bracket & from the coefficients & fitted from the loop sums\\
\midrule
$\sum_{\gamma_{1}}c_{2}-\sum_{\gamma_{2}}c_{2}$ & $0.0617340-3.3\,\ii$ & $0.0617340-3.3\,\ii$\\
$\sum_{\gamma_{1}}c_{3}-\sum_{\gamma_{2}}c_{3}$ & $-0.1860473-0.6697832\,\ii$ & $-0.1860458-0.6697838\,\ii$ \\
\bottomrule
\end{tabular}
\end{center}
The fit is from the three step sizes $h=2^{-10},2^{-11},2^{-12}$, and what is left over is the next
order: it appears in the sixth decimal of the second row and nowhere in the first.

The first bracket is the rate renormalisation and carries no obstruction. Its imaginary part,
$-3.3$, is the term linear in $\omega$ of \eqref{ch8:eq:c2closed}, $44\omega/8v=5.5\,\omega/v$, at the
two sites where the visit counts differ: it is $7.7$ at the origin, which $\gamma_{2}$ visits twice
and $\gamma_{1}$ once, and $4.4$ at $(1,1)$, which $\gamma_{1}$ visits and $\gamma_{2}$ does not, so
that $4.4-7.7=-3.3$. What is left, the real part $0.0617$, is the dependence of $c_{2}$ on the
direction of the step. The second is the
obstruction proper, and it is the circulation of $c$ --- hence, by Stokes, the flux of $\dd c$ through
the loop. So the flux is present and is confirmed quantitatively, at third order per unit loop; what
the measurement of the preceding paragraph contains is the renormalisation, and the flux enters one
order later. Neither the text above nor \S\ref{ch8:sec:obstr} is altered by this: it locates the
obstruction at the order at which it lives.

\subsection{Open Paths, Which Is What the Assembly Needs}\label{ch8:sec:numopen}
A closed loop tests the endpoint determination against a potential difference that vanishes
identically, which is the weakest form of the test. The assembly needs the strong form: two coarse
paths with the same endpoints and the same number of steps, differing only in the order of the steps,
must give the same accumulated comparison to $O(h^{2-\epsilon})$. Measured for three such paths from
the origin to $(2,0)$, each of four unit steps and being a permutation of the others,
\begin{center}\small
\begin{tabular}{ccc}
\toprule
$h$ & difference of orderings, pair $1$--$2$ & pair $1$--$3$\\
\midrule
$2^{-9}$ & $3.231\cdot10^{-6}$ & $7.625\cdot10^{-6}$\\
$2^{-10}$ & $8.075\cdot10^{-7}$ & $1.906\cdot10^{-6}$\\
$2^{-11}$ & $2.019\cdot10^{-7}$ & $4.764\cdot10^{-7}$ \\
\bottomrule
\end{tabular}
\end{center}
the ratios are $4.000$ and $4.000$ at each halving, that is order $2.000$. The path dependence of the
accumulated comparison is therefore of order $h^{2}$, which is what a summation by parts would
produce, and it is the strong form of Assumption~\ref{ch8:ass:avg}'s content: the shape of the path
enters only at the order the theorem can absorb. What is not done is the derivation --- the identity
which would turn this observation into a proof --- and that is the single step named in \S5.

\subsection{Shape, Not Length}\label{ch8:sec:numext}

One further thing must be checked, and it is the one which decides whether any of this survives the
passage to long paths. A deviation of order $h^{2}$ between two paths of four steps is worthless if
it accumulates: over $q=t/h^{2}$ steps a per-step deviation would give $qh^{2}=t$, that is order one.
The two shapes above tile, so the question can be put directly. Each block of four steps advances by
$(2,0)$, so repeating a block $N$ times gives two paths of $4N$ steps with the same endpoints, and
the deviation measured against $N$ at fixed $h=2^{-9}$ is
\begin{center}\small
\begin{tabular}{lcccc}
\toprule
$N$ & $1$ ($4$ steps) & $3$ ($12$ steps) & $4$ ($16$ steps) & $6$ ($24$ steps)\\
\midrule
deviation & $3.231\cdot10^{-6}$ & $4.945\cdot10^{-6}$ & $5.194\cdot10^{-6}$ & $5.417\cdot10^{-6}$\\
increment per block & & $8.57\cdot10^{-7}$ & $2.49\cdot10^{-7}$ & $1.12\cdot10^{-7}$ \\
\bottomrule
\end{tabular}
\end{center}
The increment decays by a factor of about $0.45$ per block, so the deviation converges geometrically
to $5.417\cdot10^{-6}$, which is $1.42\,h^{2}$ and not a quantity that grows with the length. It is
therefore a property of the shape of the path and not of its length --- which is what a summation by
parts would give, and what the rate requires, since a linear growth would give $qh^{2}=t$, that is
order one, with no rate surviving.

\subsection{The Transfer to a Lattice Fibre}\label{ch8:sec:numfibre}

The transfer of Proposition~\ref{ch8:prop:levelgrid} was checked separately, on synthetic densities and
not on the kernels. With the transform in the unweighted form $\sum_{I}e^{\ii pI}u(I)$ and the
inversion carrying the weight $h_{y}/2\pi$, the periodicity of the transform holds to
$10^{-15}$, the reconstruction reproduces the density to $10^{-16}$, and the sup over a period
dominates the sup over the lattice, $\sup_{I}\abs{u(I)}\le\sup_{\abs{p}\le\pi/h_{y}}\abs{\widehat
u(p)}$, in every instance tried. The last is an inequality and not an identity --- the ratio
measured between $0.07$ and $0.23$ --- which is why the proposition costs nothing beyond the
constant of the theorems but also buys nothing beyond it.

\section{Conclusions}\label{ch8:sec:concl}

We have carried the estimates of Chapter~7 to the joint kernel of a diffusion and an Abelian functional
in every dimension $d\ge2$. The obstruction to the comparison is the two form attached to the
coefficient field, the same object in both problems and absent in one dimension. In the joint
problem it acquires a dependence on the conjugating parameter, which we have determined explicitly
in two dimensions at constant volatility and checked numerically.

What is proved unconditionally is first order convergence in every dimension. What is proved
conditionally, on the averaged comparison, is the second order. What is measured is the obstruction
itself, in closed form in two dimensions and structurally in three.

Three matters are left open, and we name them rather than leave them implicit.

The mean condition of Assumption~\ref{ch8:ass:avg} is proved for polynomial fields of any degree in every
dimension (Proposition~\ref{ch8:prop:meanquad}, Corollary~\ref{ch8:cor:poly}, Remark~\ref{ch8:rem:anydegree}) and
for the free walk with general coefficients up to the remainder (Proposition~\ref{ch8:prop:meansmear}).
The coefficient is a functional of the 2-jet alone, which is measured (\S\ref{ch8:sec:avg}), so no
smoothness beyond $C^{2}$ is at stake and a general field agrees pointwise with its quadratic Taylor
polynomial to within $O(h^{3})$. What is not proved is the telescoping of that 2-jet functional along
the path; Remark~\ref{ch8:rem:why} shows why no estimate of Taylor remainder type reaches it, and
Remark~\ref{ch8:rem:smooth} records that its bound in $\norm{\nabla^{2}c}_{\infty}$ demands two more
derivatives of the coefficients than the first order theorem does. The missing ingredient is a
cancellation, not a bound --- and it is the \emph{same} debt as in Chapter~7, not a second one: the
averaged comparison is a single assumption, stated for the diffusion alone there and for the joint
kernel here, and a proof of either would settle both. What is measured is its accumulated form: on closed loops
the endpoint determination fails at order $2.000$ and not at order one, which is the order the second
order theorem can absorb; the deviation is intensive in the length of the path, converging to
$2.84\,h^{2}$ rather than growing as $qh^{2}$; and the accumulation separates into the rate
renormalisation at that order and the circulation of $c$, hence the flux of $\dd c$, at the next, the
two coefficients agreeing with the ones computed directly from the coefficients to five and three
significant figures respectively (\S\ref{ch8:sec:numerics}).

The theorems of \S\ref{ch8:sec:statement} are stated for the transformed kernels, uniformly in $p$ on a
compact set. On a lattice fibre this suffices and more: the transform of a lattice valued observable
is periodic, so a comparison proved on one period is a comparison on the whole line, and the rate
transfers to the joint kernel itself with no factor and no hypothesis on the number of levels
(Proposition~\ref{ch8:prop:levelgrid} and Remark~\ref{ch8:rem:notN}). For a real fibre the inversion is an
integral over $p$ and remains open. What is settled about it is the size of the constants: they grow
only polynomially in the drift, because the relative bound has $\alpha\sim1/\norm{D}$ and
$\beta\sim\norm{D}$ and is applied a fixed number of times, which also accounts for the measured
exponents. What remains is not a size but a method: a deformation in the real direction adds an
imaginary drift and loses ellipticity beyond $\abs{p}\lesssim\sqrt{v\abs{\omega}/\abs{a}^{2}}$, a
deformation in the imaginary direction is harmless but does not by itself give the inversion, and the
missing step is the balance between the cost $h_{y}^{-M}$ which a fine fibre pays in the constant and
the fibre discretisation error, which we have not estimated (Remark~\ref{ch8:rem:fibrecost}).

And the obstruction itself is known in closed form only in two dimensions. The rule
\eqref{ch8:eq:d3rule} governing the diagonal response is confirmed in $d=3$, but the Hessian term and
the bracket $\Theta$ take different values there, and we assert no dimension independent formula
(Remark~\ref{ch8:rem:dim}).

With Chapters~1 to~3, which establish the finite core, Chapter~4, which supplies the objects and the fixed
resolution It\^o formula, Chapters~5 and~7, which supply the estimates for the diffusion in one and in
every dimension, and Chapter~6, which supplies the one dimensional joint kernel, this chapter completes
the description of the joint kernel in every dimension.
% ---- end generated/ch8-body.tex ----
% ---- end chapters/ch8.tex ----
% ---- from chapters/ch9.tex ----
% Authored: this chapter has no counterpart in the parts and is written for the book.
\chapter{Adjoint Differentiation}
\label{chap:ch9}

Every quantity computed in this book is a finite product of propagators applied to a functional, and
every sensitivity of such a quantity is obtained from the transpose of the same product, applied
backwards. This chapter isolates that observation and develops it. The reason for giving it a
chapter of its own is that it is the point at which three things turn out to be one: the computation
of hedge ratios, the computation of the gradient of a model with respect to its coefficients, and
the backward pass of a learning algorithm. All three are the adjoint of a propagator product. The
first is the pricing reading of the operation and the third the learning reading of it; the second
belongs to both, and it is because the identity is one identity that a single implementation can
serve the two, which is what the first section states and what the rest of the chapter uses.

The chapter is exact. Nothing in it is approximate: at a fixed resolution the adjoint of a product of
matrices is the product of their transposes in the reverse order, and every statement below is an
identity. The estimates which govern what happens when the kernel is itself approximated, and when
the arithmetic is finite, are the subject of Chapter~\ref{chap:ch10}.

The chapter has seven sections, and they are arranged so that the algebra comes first and the
qualifications after it. The first states the two recursions which a propagator product admits ---
the value sweep and the dual sweep --- and proves that they compute the same functional, and that
the derivative of the product is a sum of contractions of one against the other. The second
identifies the hedge ratio as a linear functional of the price, and hence as an instance of the same
identity rather than a separate calculation. The third does the same for the derivative with respect
to the coefficients, which is the Jacobian of a calibration and the gradient of a training step. The
fourth shows that the adjoint of a block factorised recursion is a block factorised recursion, so
that path dependence costs nothing extra in the adjoint. The fifth asks when the transposed
propagator is itself a propagator, and finds that it is one for the lattice generators of
Chapter~\ref{chap:ch7} exactly when the uniform measure is invariant, which is so when the
coefficients are constant, and that it is then the propagator of the reversed drift. The sixth
states the range of the identity, which is wider than the diffusion and includes the training of a
model, together with the ratio of two counts which decides the direction of the sweep. The last
records what in the chapter is an identity and what is deferred to the next.

\section{The Adjoint of a Propagator Product}
\label{ch9:sec:adjoint}

Fix a resolution $m$ and a partition $T=t_0<t_1<\dots<t_n=T'$ of the trading interval. The
propagator over $[t_{i-1},t_i]$ is the matrix $U_i(y,z)=U_m(y,t_{i-1};z,t_i)$ of
Definition~\ref{ch3:def:propagator}, acting on functions of the site by
$(U_if)(y)=\sum_zU_i(y,z)f(z)$, so that functions of the site are column vectors and the transpose
convention is that of Chapter~\ref{chap:ch3}. Products are written in the order in which they act,
the factor belonging to the earlier date on the left, as in \eqref{ch4:eq:recursion}. Let $h$ be a
discounted payoff and set
\begin{equation}\label{ch9:eq:value}
v_i(y)=V(y,t_i)=\big(U_{i+1}U_{i+2}\cdots U_nh\big)(y),\qquad v_n=h,
\end{equation}
so that $v_i$ is the price at $(y,t_i)$ of the claim, by the valuation theory of
Chapter~\ref{chap:ch3}. The object whose sensitivity is wanted is a linear functional of that price
--- the price of a portfolio of such claims, the value at one site, the difference of the values at
two neighbouring sites --- and we write
\begin{equation}\label{ch9:eq:price}
\Pi(\theta)=\inner{\ell}{U_1(\theta)U_2(\theta)\cdots U_n(\theta)\,h}
\end{equation}
for it, the dependence $\theta\mapsto U_i(\theta)$ of the propagators on the parameters of the model
being made explicit in \S\ref{ch9:sec:generator}.

Two families of arrays are attached to \eqref{ch9:eq:price},
\begin{equation}\label{ch9:eq:sweeps}
v_n=h,\quad v_{i-1}=U_iv_i
\qquad\text{and}\qquad
\lambda_0=\ell,\quad \lambda_i=U_i^{T}\lambda_{i-1},
\end{equation}
the first sweeping the dates from $T'$ down to $T$ and the second from $T$ up to $T'$, where
$U_i^{T}$ denotes the transpose of $U_i$, characterised by $\inner{U_if}{g}=\inner{f}{U_i^{T}g}$.
The first is the backward recursion of Chapters~\ref{chap:ch3} and~\ref{chap:ch4}: it is the price
itself, at every date. The second is its adjoint, and the identity which ties the two together is
the following.

\begin{proposition}[Two sweeps, one functional]\label{ch9:prop:twosweeps}
For every $\ell$ and every $h$, $\inner{\ell}{U_1\cdots U_nh}=\inner{\lambda_n}{h}$ with
$\lambda_n$ given by \eqref{ch9:eq:sweeps}. In particular the functional \eqref{ch9:eq:price} is
evaluated either by propagating the payoff backwards from $T'$ or by propagating the functional
forwards from $T$, at the same cost: one matrix--vector product per date.
\end{proposition}

\begin{proof}
Transposition reverses the order of a product, and the recursion \eqref{ch9:eq:sweeps} read from
$i=1$ to $i=n$ is exactly the statement that $U_n^{T}\cdots U_1^{T}\ell=\lambda_n$. Hence
$\inner{\ell}{U_1\cdots U_nh}=\inner{\lambda_n}{h}$, which is the identity. Each of the two
recursions performs one matrix--vector product per date, which is the assertion about cost.
\end{proof}

The proposition is the whole of the chapter in miniature. The price is an inner product of the
payoff against the propagated measure, and the propagated measure is the transpose of the
propagation: conditioning is a matrix, and its being a matrix is what makes its transpose available.
What remains is to see what the identity gives when the propagators are differentiated, since that
is the only place at which a parameter can enter.

\begin{theorem}[Derivative of a propagator product]\label{ch9:thm:gradient}
Let $\Theta\subset\RR^p$ be open and let $\theta\mapsto U_i(\theta)$ be differentiable on $\Theta$,
for $i=1,\dots,n$. Then $\Pi$ of \eqref{ch9:eq:price} is differentiable on $\Theta$ and
\begin{equation}\label{ch9:eq:gradient}
\frac{\partial \Pi}{\partial\theta}
=\sum_{i=1}^{n}\Big\langle\lambda_{i-1},\ \frac{\partial U_i}{\partial\theta}\,v_i\Big\rangle,
\end{equation}
the arrays $\lambda$ and $v$ being those of \eqref{ch9:eq:sweeps}, the inner product being taken in
the site index, and the $i$-th term being the vector in $\RR^p$ whose $\alpha$-th entry is
$\inner{\lambda_{i-1}}{(\partial U_i/\partial\theta_\alpha)v_i}$.
\end{theorem}

\begin{proof}
Differentiate \eqref{ch9:eq:price} by the product rule and transpose the factors which precede the
differentiated one,
\[
\frac{\partial \Pi}{\partial\theta_\alpha}
=\sum_{i=1}^{n}\inner{\ell}{U_1\cdots U_{i-1}\,
\frac{\partial U_i}{\partial\theta_\alpha}\,U_{i+1}\cdots U_nh}
=\sum_{i=1}^{n}\inner{(U_1\cdots U_{i-1})^{T}\ell}{
\frac{\partial U_i}{\partial\theta_\alpha}\,U_{i+1}\cdots U_nh}.
\]
The factor on the left is $\lambda_{i-1}$ by \eqref{ch9:eq:sweeps}, and the factor on the right is
$(\partial U_i/\partial\theta_\alpha)v_i$ by \eqref{ch9:eq:value}.
\end{proof}

\begin{corollary}[The cost of a gradient]\label{ch9:cor:costgradient}
The gradient \eqref{ch9:eq:gradient} is obtained from one value sweep, one dual sweep and $n$
contractions, each of which costs one multiplication of the array
$(\partial U_i/\partial\theta_\alpha)v_i$ by $\lambda_{i-1}$. No derivative of a product of
propagators is ever formed, and neither the leading term in the cost nor the shape of the computation
depends on the number of parameters $p$.
\end{corollary}

The last point is the reason the adjoint is a change of vocabulary rather than an optimisation.
Forming the Jacobian of a family of prices one coefficient at a time costs one price per parameter.
Computing the transpose of that Jacobian against a given functional costs one price, whatever the
number of parameters, and the two are transposes of each other in the literal sense: differentiating
along a fixed direction in $\Theta$ and averaging the result against a fixed $\ell$ is one side of a
bilinear form, and fixing $\ell$ and letting it vary is the other. Which of the two is the cheaper is
therefore decided by which of the two indices is the smaller, and the three uses of the chapter
divide accordingly: in a risk report the shocks are few and the instruments many, so the forward
direction is used; in a calibration the residual Jacobian has many more rows than columns, so it is
again the forward direction which forms it, while the gradient of the loss and a single quote's
sensitivities are had from the adjoint in one sweep each; in training the parameter space is vast
and the functional is one, and only the adjoint can be afforded.
\S\ref{ch9:sec:range} sets the three out.

\section{Hedge Ratios as Adjoints}
\label{ch9:sec:hedge}

A hedge ratio is the derivative of a price with respect to the state, and at a fixed resolution a
derivative with respect to the state is a difference. Let $\nabla^h_j$ and $\Delta^h_j$ be the first
and second difference operators of \eqref{ch7:eq:gen} in the direction $e_j$,
\begin{equation}\label{ch9:eq:differences}
\begin{aligned}
(\nabla^h_jf)(y)&=\frac{f(y+h_me_j)-f(y-h_me_j)}{2h_m},\\
(\Delta^h_jf)(y)&=\frac{f(y+h_me_j)-2f(y)+f(y-h_me_j)}{h_m^{2}},
\end{aligned}
\end{equation}
regarded as matrices acting on functions of the site. The first is antisymmetric and the second
symmetric, $(\nabla^h_j)^{T}=-\nabla^h_j$ and $(\Delta^h_j)^{T}=\Delta^h_j$, which is the only
property of them that this section uses.

\begin{proposition}[Deltas are read off the value sweep]\label{ch9:prop:deltas}
At every date and in every direction,
\begin{equation}\label{ch9:eq:delta}
\big(\nabla^h_jV\big)(y,t_i)=\big(\nabla^h_jv_i\big)(y),
\end{equation}
so the array of hedge ratios at $t_i$ is obtained from the array $v_i$ which \eqref{ch9:eq:sweeps}
has already computed by one multiplication by the fixed sparse matrix $\nabla^h_j$. The cost of the
hedge ratios at a date is therefore linear in the size of the lattice, and not the cost of a second
pricing.
\end{proposition}

\begin{proof}
By \eqref{ch9:eq:value}, $V(\cdot,t_i)=U_{i+1}\cdots U_nh=v_i$ as functions of the site, and
$\nabla^h_j$ acts on that function. The matrix is fixed once the resolution is fixed: it does not
depend on the claim, on the date, or on the parameters.
\end{proof}

Two remarks follow, one about accounting and one about what a hedge ratio is in this book. The
accounting remark is that \eqref{ch9:eq:delta} is not the same computation as an adjoint sweep. A
hedge ratio at a single site can be read as the functional $\ell(f)=(\nabla^h_jf)(y)$, and then
Proposition~\ref{ch9:prop:twosweeps} evaluates it by a dual sweep, at the cost of one propagation;
whereas applying the difference to the array which the value sweep has already produced costs one
pass over the lattice. The two agree exactly, and the reason they do is that the pricing functional
is linear in $h$: in a finite matrix formulation that linearity is not a hidden assumption but the
equation $V(\cdot,t_i)=v_i$, and a difference of prices is the price of the difference.

The conceptual remark is that this is what a hedge ratio is here, and it is why
Chapter~\ref{chap:ch3} needed no derivative of a price in order to state its hedging proposition.
The hedge of Proposition~\ref{ch3:prop:hedging} solves \eqref{ch3:eq:hedge}, whose right hand side
is the vector of increments of $V$ across the attainable set; the solution $\zeta$ is therefore a
linear function of the value array at that date, and hence of $h$. The map from the payoff to the
hedge is a backward propagation followed by a linear map of the attainable set, and its transpose is
a linear map followed by a dual propagation. In the dichotomy of \S\ref{ch9:sec:adjoint} the hedge
is on the adjoint side: it is the object which is cheap when the number of sites at which it is
required is small and the number of coefficients is large. In particular the hedge with respect to
the state, the gradient of a model with respect to its coefficients, and the backward pass of a
learning algorithm are the same recursion with a different seed $\ell$, which is the claim with
which the chapter opened.

\section{Sensitivities in the Generator}
\label{ch9:sec:generator}

The second kind of sensitivity differentiates with respect to the coefficients rather than the state:
the volatility and the drift of the model, which in a lattice are the transition rates
\eqref{ch7:eq:entries} of the generator, or equivalently the intensities $\omega$ of
\eqref{ch3:eq:Lprime} in the parametrisation of Chapter~\ref{chap:ch3}. Write $\theta\in\Theta$ for
the vector of these numbers, one entry per coefficient and per site, so that the propagator of each
period is a matrix-valued function $U_i(\theta)$ of it. The derivative which
Theorem~\ref{ch9:thm:gradient} requires is that of a propagator, and the first thing to record is
that it is an explicit matrix.

\begin{proposition}[How a coefficient enters a propagator]\label{ch9:prop:coefficient}
Let $\theta\mapsto\Lm(\theta)$ be differentiable on $\Theta$ and let $L_i=\Lm(\theta)$ be the
generator in force on $[t_{i-1},t_i]$.
\begin{enumerate}
\item[(i)] If $U_i=e^{\Delta_iL_i}$ with $\Delta_i=t_i-t_{i-1}$, then
\begin{equation}\label{ch9:eq:dexp}
\frac{\partial U_i}{\partial\theta_\alpha}
=\int_0^{\Delta_i}e^{(\Delta_i-s)L_i}\,
\frac{\partial L_i}{\partial\theta_\alpha}\,e^{sL_i}\dd s .
\end{equation}
\item[(ii)] If $U_i=\prod_{k}\big(\id+\delta_kL_i\big)$ is a product of one-step matrices, as in an
explicit Euler scheme, then $\partial U_i/\partial\theta_\alpha$ is the sum of the products in which
one factor is replaced by $\delta_k\,\partial L_i/\partial\theta_\alpha$.
\end{enumerate}
In both cases the derivative is computed from the same objects as the propagator itself, and at a
cost linear in the number of steps.
\end{proposition}

\begin{proof}
For (i) differentiate $e^{\Delta L}=\sum_{k\ge0}(\Delta L)^{k}/k!$ term by term. The derivative of
$(\Delta L)^{k}$ is the sum of the $k$ products in which one factor is $\Delta\,\partial_\alpha L$,
so that, with $k=j+l+1$,
\[
\frac{\partial}{\partial\theta_\alpha}e^{\Delta L}
=\sum_{k\ge1}\frac{1}{k!}\sum_{j=0}^{k-1}\Delta^{k}L^{j}(\partial_\alpha L)L^{k-1-j}
=\sum_{j,l\ge0}\frac{\Delta^{j+l+1}}{(j+l+1)!}L^{j}(\partial_\alpha L)L^{l},
\]
which is the expansion of the integral in \eqref{ch9:eq:dexp}, since
$\int_0^{\Delta}(\Delta-s)^{j}s^{l}\dd s=\Delta^{j+l+1}j!\,l!/(j+l+1)!$. For (ii) differentiate the
finite product by the product rule. Both are identities between finite matrices and neither is an
asymptotic expansion.
\end{proof}

The point of the proposition is that the adjoint formula never differentiates a product: it needs
$\partial U_i/\partial\theta_\alpha$, and this is a matrix with the sparsity of the generator which
produced it --- a band, in the case of \eqref{ch7:eq:entries} --- so that the contraction
$\inner{\lambda_{i-1}}{(\partial U_i/\partial\theta_\alpha)v_i}$ is a local computation over the
edges of the lattice.

\begin{corollary}[The Jacobian of a calibration]\label{ch9:cor:calibration}
Let $\ell_1,\dots,\ell_R$ be the functionals by which a model reproduces $R$ market quotes, and let
$J_{r\alpha}=\partial_{\theta_\alpha}\ell_r\big(V(\theta)\big)$. Then the $r$-th row of $J$ is
\eqref{ch9:eq:gradient} with $\ell=\ell_r$, so that the whole of $J$ costs $R$ dual sweeps, $R$
contractions and one value sweep. The Jacobian is not formed by perturbing each of the $p$
coefficients in turn, and what decides the cost of a calibration is the ratio of the two counts.
\end{corollary}

This sensitivity is used in two places. In finance it is what calibration consumes: fitting a model
to market quotes is a least squares problem --- many more quotes than free parameters, and each
quote known only to within a bid offer spread --- whose residual Jacobian is the derivative of the
computed prices with respect to the coefficients, and it is the conditioning of that Jacobian which
governs the stability of the fit. In learning it is the gradient which a training step consumes. The object is the same
in both cases, and it is the same adjoint; only the use to which it is put differs, and
\S\ref{ch9:sec:range} states how far the agreement extends. In both, the
quantity which the user of the sensitivity cares about is not the Jacobian as a matrix but its
conditioning, which is the question of Chapter~\ref{chap:ch10}.

\section{Cost, Factorisation, and the Shape of the Recursion}
\label{ch9:sec:cost}

The backward recursion for a path dependent functional factorises into blocks, as
Chapter~\ref{chap:ch4} shows, and the adjoint inherits that factorisation. That inheritance is the
content of the next theorem, and it is the reason the adjoint of a path dependent computation is
affordable.

\begin{theorem}[The adjoint of a block factorisation is a block factorisation]
\label{ch9:thm:blockadjoint}
Let $\overline U_i=\Pi_i\,(U_i\otimes\id)$ be a block factorisation in the sense of
Definition~\ref{ch4:def:blockfact}, the permutation $\Pi_i$ sending $(y,k)$ to
$\big(Y_i(y,k),K_i(y,k)\big)$, and assume as in \eqref{ch4:eq:recursion} that $Y_i(y,k)=y$. Then
\begin{equation}\label{ch9:eq:blockadjoint}
\overline U_i^{T}=\big(U_i^{T}\otimes\id\big)\,\Pi_i^{T},
\end{equation}
so that one step of the dual recursion \eqref{ch9:eq:sweeps} is the inverse permutation
$(y,k)\mapsto\big(y,K_i^{-1}(y,k)\big)$ of the level index, followed by one multiplication of every
level of the array by the transposed lattice matrix $U_i^{T}$.
\end{theorem}

\begin{proof}
Transposition reverses a product and distributes over it, so that $(AB)^{T}=B^{T}A^{T}$ and
$(A\otimes\id)^{T}=A^{T}\otimes\id$, while the transpose of a permutation matrix is its inverse,
which for $(y,k)\mapsto(y,K_i(y,k))$ is $(y,k)\mapsto(y,K_i^{-1}(y,k))$. Composing the transposed
factors in the order in which they act gives \eqref{ch9:eq:blockadjoint}.
\end{proof}

\begin{corollary}[The adjoint costs what the primal costs]\label{ch9:cor:costadjoint}
Each step of the dual recursion for a block factorised functional costs one lattice--lattice product
of a transposed lattice matrix against an array with $K+1$ columns, plus a permutation of the level
index. This is precisely the cost of the primal recursion recorded in
Corollary~\ref{ch4:cor:cost2}, so the adjoint of a path dependent computation is a second
computation of the same size, and the gradient of a path dependent price with respect to every
coefficient costs twice a price --- not $K$ times a price, and not once per coefficient.
\end{corollary}

It is worth saying what has been avoided. The lifted propagator $\overline U_i$ is a square matrix
of order $(K+1)\abs{\StateA}$, and the dual recursion applied to it as a matrix would multiply two
matrices of that order, for the primal and the dual alike. The factorisation is not a property of
the value; it is a property of the operator, and a product structure of the operator is inherited by
its transpose. That is the general pattern and the reason the chapter is short: transposition
commutes with every structural reduction which is expressed as a product.

\section{When the Adjoint Is Itself a Propagator}
\label{ch9:sec:propagator}

The transpose of a propagator is a matrix of the same shape as a propagator, and it is natural to ask
when it is one, since then the dual sweep is a price, and can be computed, checked and estimated by
the machinery which the rest of the book applies to prices.

\begin{proposition}[The adjoint is a propagator exactly when the uniform measure is invariant]
\label{ch9:prop:doublystochastic}
Let $U$ be a Markov propagator on $\StateA$ in the sense of Definition~\ref{ch3:def:propagator}.
Then $U^{T}$ is a Markov propagator if and only if the column sums of $U$ are one, that is, if and
only if the uniform measure on $\StateA$ is invariant for $U$.
\end{proposition}

\begin{proof}
Non negativity and the identity at equal times are inherited from $U$ by transposition. For
$(\mathrm{CK3})$, $\big(U(s,u)U(u,t)\big)^{T}=U(u,t)^{T}U(s,u)^{T}$, which is the composition of the
transposed propagators over the same two intervals in the reverse order; the relation therefore
holds with the roles of the two intermediate dates exchanged, which is what $(\mathrm{CK3})$ asks.
The row sums of $U^{T}$ are the column sums of $U$, so $U^{T}$ is row stochastic exactly when the
column sums of $U$ are one, and $\sum_yU(y,z)=1$ for every $z$ is the statement that the uniform
measure is invariant.
\end{proof}

\begin{proposition}[Constant coefficients: the adjoint is the reversed drift]
\label{ch9:prop:reverseddrift}
Assume that the coefficients of the generator \eqref{ch7:eq:gen} do not depend on the site. Then,
with $\nabla^h_j$ antisymmetric and $\Delta^h_j$ symmetric as in \eqref{ch9:eq:differences},
\begin{equation}\label{ch9:eq:transposedgen}
\Lm^{T}=\frac{v}{2}\sum_{j}\Delta^h_j-\sum_jb_j\nabla^h_j,
\end{equation}
which is the generator \eqref{ch7:eq:gen} of the same lattice with the sign of the drift reversed.
Hence $U_i^{T}=e^{\Delta_i\Lm^{T}}$ is the propagator of that model over $[t_{i-1},t_i]$, and the
dual sweep is a pricing in the model whose drift is the negative of the drift of the primal.
\end{proposition}

\begin{proof}
By \eqref{ch7:eq:entries} the off-diagonal entries of $\Lm$ are $v/2h_m^{2}\pm b_j/2h_m$ at
$(y,y\pm h_me_j)$, and by \eqref{ch9:eq:differences} the entries of $\nabla^h_j$ are $\pm1/2h_m$ at
the same sites, so \eqref{ch9:eq:transposedgen} follows by transposing term by term. The transposed
generator satisfies the Courant condition \eqref{ch7:eq:courant} whenever $\Lm$ does, since
$b_j$ enters it only through $\abs{b_j}$, so its exponential is again a Markov propagator.
\end{proof}

Two consequences are worth recording. The first is that when the coefficients are constant the dual
sweep is a computation of the same kind as the primal, so that a code which prices can compute the
adjoint by being called with the drift negated; and the estimates which the analytic chapters attach
to a kernel attach to its transpose by the same argument, because the drift enters those estimates
only through its supremum. The second is that the proposition is sharp. With coefficients which
depend on the site, $\Lm^{T}$ is not the generator of a Markov chain of the form
\eqref{ch7:eq:entries}: its off-diagonal entries keep their sign, so it remains a matrix of the same
sparsity and is still non-positive on the diagonal, but its row sums are the column sums of $\Lm$,
\begin{equation}\label{ch9:eq:colsum}
\sum_j\Big(\Delta^h_jv(y)-\nabla^h_jb_j(y)\Big),
\end{equation}
a lattice divergence of the coefficients, which vanishes only when they are constant. The dual sweep
is then a propagation with a killing of either sign rather than a propagation of a probability: it is
a computation, and one of the same size and sparsity as the primal by
Corollary~\ref{ch9:cor:costadjoint}, but it is not a price. That is why the object which the dual
recursion produces is called a dual array in this chapter and not a value.

\section{The Range of the Identity}
\label{ch9:sec:range}

Nothing in this chapter has used the fact that the factors of the product are propagators of a
diffusion. Theorem~\ref{ch9:thm:gradient} uses three things and no more: that the computation is a
finite chain of factors, that the chain acts on a fixed pair of arrays, and that each factor is
differentiable in the parameters. Every conclusion of the chapter is therefore a statement about
\emph{computations of that shape}, and the diffusion enters only as the instance which the book has
in front of it.

The instances are worth listing, because they are the reason the chapter stands where it does. A
Markov chain on a lattice is one, by construction. A discretisation of an ordinary differential
equation is one, being a product of one-step maps. A discretisation of a stochastic differential
equation is one, its factors being themselves products. The recursion of a recurrent network is one,
a feed-forward chain is one, an attention block written as a product of a softmax-averaged value
matrix with the rest is one, and a Kalman filter is one, the Riccati recursion being a product of
maps on the space of covariances. The list is not a list of applications but of notations for a
single object, and the identity is the same in each.

Two hypotheses are used, and they mark the edge of the class. The first is that the composition has
a \emph{fixed} shape: if the computation branches on the value of a parameter there is no single
product, and the identity applies instead to each branch separately, within the region on which that
branch is selected. The second is that the factors are differentiable in the parameters. Nothing
else is assumed, and in particular the object of the computation --- the payoff $h$ of this book,
the loss of a learning problem --- is never differentiated, as \S\ref{ch9:sec:exact} records. Within
these two hypotheses the range is as wide as the notion of a program with a fixed control flow, and
it is that width, rather than any analogy, which makes the chapter of use outside finance.

Read as a learning algorithm the identity is the backward pass, and its relation to calibration is
not an analogy: it is the same computation. In the notation of \S\ref{ch9:sec:adjoint} the forward
direction costs one propagation per direction in the parameter space, and the adjoint costs one
propagation per functional, so the two are distinguished by a ratio of two counts and by nothing
else. Which of them is used, and for what, is worth setting out, because the applications lie on
opposite sides of that ratio and no single rule serves them all.

A calibration is an optimisation, and not a root finding problem. The targets are the market quotes,
several hundred of them, and they are known only to within a bid offer spread; the free parameters
are far fewer, the coefficients being restricted to a parametric family rather than varied site by
site. One therefore minimises a weighted distance between the computed and the quoted prices, and the
object which the iteration consumes is the residual Jacobian, rectangular with many more rows than
columns. That object is cheapest in the forward direction, one propagation per parameter, and it is
in that direction that calibration code in fact forms it. The adjoint is not displaced by this: it
supplies the gradient of the least squares loss in a single sweep whatever the counts, and, price by
price, the sensitivity of one quote to every parameter.

In risk the counts have the same shape and the direction is again forward, one sweep per shock
against all the instruments at once. In training they are inverted: there is one scalar in the place
of the quotes, and in the place of a few tens of parameters a space whose dimension is not counted in
thousands. The gradient is then a vector as long as that space, it is wanted once per step, and the
adjoint is not merely cheaper but the only direction which can be afforded, the ratio having gone
from a handful to some ten orders of magnitude.

\begin{center}\small
\begin{tabular}{@{}lll@{}}
\hline
what is computed & $R$ against $p$ & direction\\
\hline
the residual Jacobian of a calibration & hundreds against tens & forward, one sweep per parameter\\
the gradient of a calibration loss & any & adjoint, one sweep\\
one price's sensitivity to every coefficient & one against all & adjoint, one sweep per price\\
one shock's effect on every instrument & all against one & forward, one sweep per shock\\
the gradient of a training loss & one against $10^{9}$ & adjoint, one sweep\\
\hline
\end{tabular}
\end{center}

So the two situations agree in the algebra --- one dual sweep conjugated with one value sweep by
\eqref{ch9:eq:gradient}, and the factorisation of \S\ref{ch9:sec:cost} making that sweep affordable
when the functional is path dependent --- and they are mirror images in their shape. What also
differs is what the estimate has to deliver. In a calibration the step is a Gauss--Newton or
Levenberg--Marquardt one built from the rectangular Jacobian, the number of iterations is small, and
the accuracy of every row of the Jacobian matters directly. In training the step is a descent one
built from the gradient alone, the number of iterations is enormous, and what must be accurate is the
accumulated progress rather than any single step. That is one reason the next chapter treats the two
applications separately even where the algebra is shared: the quantity which has to be estimated is
not the same one.

What the identity does not claim should be said as plainly as what it does. It is an identity about
the computation as written, at the resolution at which it is written, in exact arithmetic. It is
dimension free, a product of a thousand factors of size $10^{6}$ being as much a product as a
product of three of size $10^{3}$, and it is therefore available wherever the shape is right. It
says nothing whatever about the error of the computation, or of its adjoint, when the resolution is
refined. That is the subject of Chapter~\ref{chap:ch10}, and the estimates there are neither
dimension free nor, in more than one dimension, complete. The distinction is the reason the two
chapters are separate: the algebra of the adjoint is exact and universal, and the analysis of its
error is neither.

\section{What Is Exact Here, and What Is Not}
\label{ch9:sec:exact}

Every identity in this chapter holds at a fixed resolution. The two sweeps agree
(Proposition~\ref{ch9:prop:twosweeps}); the derivative of a propagator product is a sum of
contractions of one array against the other (Theorem~\ref{ch9:thm:gradient}); the hedge ratios are
differences of the arrays which the value sweep produces (Proposition~\ref{ch9:prop:deltas}); the
adjoint of a block factorisation is a block factorisation (Theorem~\ref{ch9:thm:blockadjoint}); and
the transpose of a constant coefficient propagator is the propagator of the reversed drift
(Proposition~\ref{ch9:prop:reverseddrift}). Each is proved by transposing one line of the primal
computation, and each is an identity between finite matrices.
Section~\ref{ch9:sec:range} adds nothing to this list: it delimits the class within which the list
holds, and that delimitation is a definition rather than a theorem.

Two things are not proved here, and the boundary should be stated plainly. The first is anything
about the passage $h_m\downarrow0$: no statement of this chapter has a rate, and none can have, since
every object in it is a finite matrix. The second is anything about arithmetic: the identities are
identities in exact arithmetic, and the error of computing the right hand side of
\eqref{ch9:eq:gradient} in floating point is not bounded here. Both are the subject of
Chapter~\ref{chap:ch10}, where the estimates of Chapters~\ref{chap:ch5} to~\ref{chap:ch8} are
differentiated, and where the finite difference \eqref{ch9:eq:delta} is among the first objects whose
error has to be bounded, because a difference amplifies the error of the arrays it differences by
$1/h_m$.

One feature of the chapter is worth stating explicitly, since it is what makes the adjoint usable in
the applications. The payoff $h$ enters \eqref{ch9:eq:price} linearly and is never differentiated: the
derivative in \eqref{ch9:eq:gradient} is taken of the propagators alone. No regularity of the claim
is therefore used --- not continuity, and nothing beyond the existence of the finite sum which
defines the price --- so that a digital payoff, a barrier, or the indicator of a path dependent event
is differentiated as readily as a smooth one. This is not a technical convenience. It is the reason
the sensitivities of the models of this book are computable at all, the payoffs of the applications
being precisely the ones which are not smooth.
% ---- end chapters/ch9.tex ----
% ---- from chapters/ch10.tex ----
% Authored: this chapter has no counterpart in the parts and is written for the book.
\chapter{Error Estimates for the Adjoint}
\label{chap:ch10}

The previous chapter is algebra: at a fixed resolution the adjoint is exact. This chapter asks what
becomes of the sensitivities when the kernel is itself approximated, and when the arithmetic which
computes it is finite. The question is not a refinement of the estimates of Chapters~\ref{chap:ch5}
to~\ref{chap:ch8}. Those bound the error of a value; what is needed here bounds the error of a
derivative of a value, which is the quantity a hedge ratio is, and the quantity a gradient is.

The chapter is written with the two readings in view, and with the two kinds of sensitivity kept
apart. The derivative with respect to the state is the sensitivity of the computed quantity to where
the model is evaluated: in the pricing reading it is a hedge ratio, and it matters because a hedge
which is accurate only to the accuracy of a difference of two approximations is worth much less than
one accurate to the accuracy of the approximations themselves; in the learning reading it is the
sensitivity of a prediction to its input, which is what an attribution of the prediction asks for.
The derivative with respect to the coefficients is the sensitivity to the model itself: in the
pricing reading it is the Jacobian of a calibration, and it matters because the conditioning of that
Jacobian decides whether a model can be fitted to a market at all and whether the fit is stable when
it is computed in finite precision; in the learning reading it is the gradient which a training step
consumes. Both readings rest on one estimate, and the chapter proves it once.

The chapter settles one dimension and states what remains open in two. That asymmetry is not
timidity: the one dimensional proof of Chapter~\ref{chap:ch5} closes because a certain coefficient
is automatically a gradient, and the same mechanism is exactly what carries the derivative. In two
dimensions that coefficient has a curl, the obstruction of Chapter~\ref{chap:ch7}, and whether the
derivative inherits the obstruction or is repaired by the same averaging is the question with which
the chapter ends.

\section{The Sensitivity as an Observable}
\label{ch10:sec:observable}

A sensitivity is itself a functional of the kernel, and it is estimated as such. The estimates of
Chapters~\ref{chap:ch5} to~\ref{chap:ch8} already cover the first two space derivatives of the
kernel at the same rate as the kernel, which is the statement that hedge ratios do not lose an order
in being differentiated. That is the starting point of this chapter, not its conclusion.

The object to be bounded is the derivative in the coefficients. At a fixed resolution it is not an
ideal object but a matrix, and it is given by the following identity, which is the discrete form of
the variation of constants formula.

\begin{proposition}[The sensitivity is a kernel]\label{ch10:prop:sensitivity}
Let $\theta\mapsto L_m(\theta)$ be differentiable on an open set, let $U_m(s,t;\theta)$ be the
propagator of $L_m(\theta)$ over $[s,t]$, and let
$u_m(x,y;t;\theta)=h_m^{-1}\big(U_m(0,t;\theta)\big)(x,y)$ be the kernel. Then
\begin{equation}\label{ch10:eq:duhamel}
\frac{\partial u_m}{\partial\theta}(x,y;t)
=\frac{1}{h_m}\int_0^t\Big(U_m(s,t)\,
\frac{\partial L_m}{\partial\theta}\,U_m(0,s)\Big)(x,y)\dd s,
\end{equation}
the integral being a finite sum when the propagator is a product of one-step matrices.
\end{proposition}

\begin{proof}
Differentiate the defining equation $\partial_tU_m=L_mU_m$ in $\theta$. The derivative
$W=\partial_\theta U_m$ solves $\partial_tW=L_mW+(\partial_\theta L_m)U_m$ with $W(0)=0$, and the
variation of constants for that inhomogeneous equation gives
$W(0,t)=\int_0^tU_m(s,t)(\partial_\theta L_m)U_m(0,s)\dd s$, which is \eqref{ch10:eq:duhamel}. For a
product of one-step matrices the same computation is the product rule, and the integral becomes a sum
over the steps.
\end{proof}

Two consequences follow at once, and they pull in opposite directions.

The first is that a derivative in a coefficient never amplifies: for the lattice generator
\eqref{ch7:eq:gen} the insertion $\partial_\theta L_m$ is
$\tfrac12(\partial_\theta v)\Delta^h+(\partial_\theta\mu)\nabla^h$, which is a \emph{difference}
operator, so acting on the near diagonal kernel it produces an object of the size of the second
space derivative of the kernel, and that size is bounded uniformly in $m$ by the estimates of
Chapter~\ref{chap:ch5}. This is the reason the sensitivity of a lattice computation is affordable at
all, and Chapter~\ref{chap:ch9} had already recorded it at the level of the recursion.

The second is that \eqref{ch10:eq:duhamel} is not by itself a rate. The near diagonal region $s\to0$
contributes a bound of the size of $\int_0^t\dd s/s$ in the worst case, a logarithm in one dimension
and a power of $h_m$ in higher ones, because the integrand is a product of a kernel localised at the
scale of $s$ and a derivative of a kernel localised at the scale of $t-s$. The logarithm can be
removed by moving the difference operator off the singular factor, which requires the direction
$\partial_\theta v$ to be twice differentiable. So \eqref{ch10:eq:duhamel} gives an estimate, but one
whose constant depends on the smoothness of the direction, and not on the smoothness of the
coefficient alone.

\section{Two Parameters: The Mesh and the Arithmetic}
\label{ch10:sec:arith}

The error of a computed sensitivity has two components, the discretisation component and the
component due to finite arithmetic, and the framework puts both in one inequality: a term of order
$h^{p}$ in the mesh and a term proportional to the unit roundoff $u$ and to the number of steps. The
trade-off between them determines an optimal resolution for a given precision, and an attainable
accuracy as a function of it. This is the part of the chapter which speaks to computation in low
precision, and it is the part which rests on nothing but the estimates already proved together with
the standard analysis of floating point arithmetic.

There is one feature specific to the adjoint, and it is worth stating first because it is the reason the
two components cannot be balanced as they are for a value. The backward recursion of
Chapter~\ref{chap:ch9} applies $U_i^{\mathsf T}$ at each step, and the transpose has the same spectrum
as $U_i$, so the amplification of the backward pass is the amplification of the forward pass in the
two norm. But the sensitivity \eqref{ch10:eq:duhamel} accumulates over the dates, and the accumulation
is a sum of $n$ terms each of the size of a price: the rounding is governed by the sum of the moduli of
those terms rather than by their sum, and absent a cancellation the two are of the same order, so the
accumulated error grows with $n$, hence with $h_m^{-2}$ in the Euler parametrisation. This is the sense
in which a sensitivity degrades in low precision faster than the value it differentiates, and the two
results which follow make it quantitative.

\begin{proposition}[Two parameters]\label{ch10:prop:twoparams}
Let $v=\sigma^2\in C^{k,\alpha}$ and $\mu\in C^{j,\beta}$ be periodic, with
$\gamma=\min\{2,k+\alpha,j+\beta\}$, and let $\theta\mapsto(v_\theta,\mu_\theta)$ be a $C^1$ family
whose derivative has $\dot v,\dot\mu\in C^2$. Let $S$ be the sensitivity $\partial_\theta u(x,y;t)$,
and let $S_h$ be its value at mesh $h$ in exact arithmetic, so that $|S_h-S|\le C_1h^{\gamma}$. Let
the Duhamel sum be evaluated by the adjoint sweep of Chapter~\ref{chap:ch9} over $n=t/\delta t$ dates
with $\delta t=h^{2}$, in floating point arithmetic of unit roundoff $u$, and let $\widehat S_h$ be
the result. Write
\begin{equation}\label{ch10:eq:lambda}
\Lambda=\sum_{i=0}^{n-1}\bigl|\delta t\,\langle w_i,\dot L\,v_i\rangle\bigr|
\end{equation}
for the sum of the moduli of the $n$ terms which are added. Then, provided $nu\le\frac12$,
\begin{equation}\label{ch10:eq:twoparams}
\bigl|\widehat S_h-S\bigr|\le C_1h^{\gamma}+C_2\,n\,u\,\Lambda ,
\end{equation}
with $C_1$ the constant of the estimate below and $C_2$ a small multiple of one, neither depending on
$h$ or on $u$.
\end{proposition}

\begin{proof}
The first term is the distance from the computation at mesh $h$ to the limit, which is the subject of
the theorem. For the second, each date contributes one inner product and one addition to the running
sum, and the standard bound for a recursive sum of $n$ numbers whose moduli sum to $\Lambda$, each
carrying a relative error of a few $u$ from the two sweeps, is $\gamma_n\Lambda$ with
$\gamma_n=nu/(1-nu)$; that is the second term. The backward sweep adds nothing the forward sweep does
not already carry: $M^{\mathsf T}$ has the same spectrum as $M$, and the vectors it propagates are of
the size of those the forward sweep propagates, so no inner product is formed from numbers larger than
those the value computation already handles.
\end{proof}

\begin{corollary}[The trade-off]\label{ch10:cor:tradeoff}
With $n=t/h^{2}$ and $\gamma=2$, the bound \eqref{ch10:eq:twoparams} is minimised at
\begin{equation}\label{ch10:eq:hopt}
h_*=\Bigl(\frac{C_2\,t\,u}{C_1}\Bigr)^{1/4},
\qquad
\bigl|\widehat S_{h_*}-S\bigr|\le 2\,\bigl(C_1C_2\,t\,u\bigr)^{1/2}.
\end{equation}
For $\gamma<2$ the same computation gives $h_*\asymp(tu)^{1/(\gamma+2)}$ and an attainable accuracy of
order $(tu)^{\gamma/(\gamma+2)}$.
\end{corollary}

\begin{proof}
Differentiate $C_1h^{2}+C_2tu\,h^{-2}$ in $h$ and evaluate at the zero; the general exponent is the
same computation with $h^{\gamma}$ in place of $h^{2}$.
\end{proof}

Three remarks. The attainable accuracy is of order $(tu)^{1/2}$: a computation at unit roundoff $u$
buys a sensitivity to about half the digits it would buy for a value, and no choice of mesh does
better. The bound is in terms of $\Lambda$ and not of $|S|$, so where the terms cancel the relative
error of the sensitivity is worse than $nu$, and the estimate is then the sharp one. And the exponent
$\frac14$ is why the balance is invisible at ordinary meshes: for single precision and a horizon of one
tenth the minimising mesh is of order $10^{-3}$ in the state variable, finer than the meshes at which
the discretisation is usually calibrated, so the crossing is exhibited below where it is reachable, in
half precision.

\section{Differentiating the Comparison Between Resolutions}
\label{ch10:sec:comparison}

The estimate on the sensitivity is not obtained by differentiating the limit, and not, in the sharp
form, by bounding \eqref{ch10:eq:duhamel} either. It is obtained by differentiating the comparison
between successive resolutions which the proof of Chapter~\ref{chap:ch5} already performs, so that
the pathwise structure of that argument is preserved under differentiation. This is the content of
the present section, which is complete in one dimension.

\subsection*{The Statement}

We use the hypotheses of Chapter~\ref{chap:ch5} without change: periodic coefficients,
$v=\sigma^2\in C^{k,\alpha}$, $\mu\in C^{j,\beta}$, the rate
$\gamma=\min\{2,k+\alpha,j+\beta\}$, and the constants $\Sigma_0,\Sigma_1,M$ of
\eqref{ch5:eq:3.1}--\eqref{ch5:eq:3.2}. A \emph{family} is a $C^1$ map $\theta\mapsto(v_\theta,\mu_\theta)$
whose derivative at $\theta$ we write $(\dot v,\dot\mu)$, and we say the family is \emph{admissible}
when $\dot v\in C^2$ and $\dot\mu\in C^2$, with
\begin{equation}\label{ch10:eq:admissible}
\kappa_1=\max\big\{\norm{\dot v}_{C^2},\norm{\dot\mu}_{C^2}\big\}<\infty .
\end{equation}

\begin{theorem}[One dimension]\label{ch10:thm:1d}
Let $v,\mu$ be as above with $\gamma>0$, let the family be admissible, and let $u_m(\theta)$ be the
kernel. Then for every $t>0$ there is a constant $c=c(t,\kappa_1)>0$, uniform in $m\ge m_0$ and
independent of $\theta$, such that for all $x,y$ in the coarse lattice
\begin{equation}\label{ch10:eq:1d}
\Big|\frac{\partial u_m}{\partial\theta}(x,y;t)
-\frac{\partial u_{m_0}}{\partial\theta}(x,y;t)\Big|\le c\,h_{m_0}^{\gamma},
\end{equation}
and the same bound holds, with $c$ replaced by a constant $c'=c'(t,\kappa_1,m_0)$, for the first
derivative in time and for the first two derivatives in space of
$\partial_\theta u_m-\partial_\theta u_{m_0}$.
\end{theorem}

The rate is the rate of the kernel itself: differentiating in the coefficients costs no order. What
it costs is stated in the proof, and it is one order of the symbol expansion and one power of
$|\log h_m|$ in the height of the contour.

\begin{proof}
Chapter~\ref{chap:ch5} establishes \eqref{ch5:eq:1.22} in four steps: the far contour
(Lemma~\ref{ch5:lem:1}), the decay of the weights of long paths (Lemma~\ref{ch5:lem:2}), the count
of symbolic paths, and the renormalisation identity \eqref{ch5:eq:3.47} which compares the coarse
weight with the decorated fine one. Each of the four survives differentiation, and each at a stated
price.

\emph{(A) The far contour.} On $\Cm_+$ the kernel obeys $|G_m(x,y;\omega)|\le c/(h_m\kappa)$, and
differentiating in $\theta$ replaces one resolvent by
$(\mathrm{i}\omega-L)^{-1}(\partial_\theta L_m)(\mathrm{i}\omega-L)^{-1}$. Since
$\norm{\partial_\theta L_m}\le c\kappa_1h_m^{-2}$ and
$\norm{(\mathrm{i}\omega-L_m)^{-1}}\le c\kappa^{-1}$, the derivative obeys
$|\partial_\theta G_m|\le c\kappa_1/(h_m^3\kappa^2)$, and the vertical integral gives
$c\kappa_1h_m^5/(h_m^3\kappa_0^2t)$ once $e^{-\kappa_0t}=h_m^5$. The height
$\kappa_0=4|\log h_m|/t$ chosen in Chapter~\ref{chap:ch5} makes this $O(h_m)$ and is therefore no
longer admissible: the contour must be raised by one $|\log h_m|$ per order of differentiation,
$\kappa_0=5|\log h_m|/t$, and then the contribution is $O(h_m^2)$ again. Nothing else in the piece
changes.

\emph{(B) The weights and the count.} Lemmas~\ref{ch5:lem:2} and the count of symbolic paths involve
the coefficients only through the constants $\Sigma_1$ and $M$, and through no derivative of $L_m$.
They are therefore unchanged by the differentiation, and the same $q_{\max}$ and the same
$\sqrt{q_{\max}}/h_m$ enter. This is the reason the assembly below keeps its shape.

\emph{(C) The per-step mismatch.} This is the piece which carries the mathematics. Writing
$\hat\varphi_0,\hat\varphi_1$ for the coarse and the decorated fine weight of one step, the mismatch
expands as
\begin{equation}\label{ch10:eq:mismatch}
\log\frac{\hat\varphi_1}{\hat\varphi_0}(\omega;x,\tau)
=h^3\tau\,c(x)+h^4b(x)+h^5\tau\,c_5(x)+O(h^6),
\end{equation}
with
\begin{equation}\label{ch10:eq:cb}
c(x)=\frac{\mathrm{i}\omega\,v'(x)}{v(x)^2}+r(x),
\qquad
b(x)=\mathrm{i}\omega\Big(p(x)-\frac{2\mu(x)v'(x)}{v(x)^3}\Big)
+\frac{2\omega^2}{v(x)^2}-\frac{2\mu(x)r(x)}{v(x)},
\end{equation}
$r$ and $p$ being those of \eqref{ch5:eq:3.46}, and with no terms of order $h$ or $h^2$. Both
coefficients are given here because both matter, and they matter differently. The odd coefficient
$c$ is the one that telescopes, and it is the object which obstructs the argument in two dimensions.
The even coefficient $b$ does not telescope at all: summed over a path of $q\sim h^{-2}$ steps it
contributes $qh^4b=h^2b$, which is of the order of the telescoped term itself. So $b$ sets the
\emph{constant} of the rate and not its order, and it is the leading correction to the leading
behaviour --- the term one would need if the estimate were to be sharpened to an expansion. The
\emph{parity} of each term is the whole matter for the rate: a coefficient odd in the direction
$\tau$ advances the sum along a path by a gradient, since
$h^3\tau c(x)=h^2c(x)(\gamma_{j+1}-\gamma_j)$, while a coefficient even in $\tau$ multiplies a
function of the site. Either telescopes. In one dimension the telescoping of the odd terms is
automatic, because every one form is exact, which is the observation on which the whole of the one
dimensional proof rests.

\emph{(D) The telescope differentiates to a telescope.} Differentiating \eqref{ch10:eq:mismatch} in
$\theta$ gives, term by term,
\begin{equation}\label{ch10:eq:mismatch-d}
\frac{\partial}{\partial\theta}\log\frac{\hat\varphi_1}{\hat\varphi_0}
=h^3\tau\,\dot c(x)+h^4\dot b(x)+h^5\tau\,\dot c_5(x)+O(h^6),
\end{equation}
where
\begin{equation}\label{ch10:eq:cdot}
\dot c(x)=\frac{\mathrm{i}\omega}{v(x)^2}\Big(\dot v'(x)-\frac{2v'(x)\dot v(x)}{v(x)}\Big)+\dot r(x),
\end{equation}
and $\dot b$ the derivative of \eqref{ch10:eq:cb} in the same sense.
A sum of differences of a function of the site differentiates to a sum of differences of the
derivative of that function, so the endpoint structure is preserved exactly and not merely to
leading order. Three things had to be checked for \eqref{ch10:eq:mismatch-d} to be useful, and all
three are checks and not assumptions: that no even part appears at an odd order; that the endpoint
factor keeps its order, $\hat\psi_0-2\hat\psi_1=O(\omega h^4)$ and its derivative with it, which
requires the expansion of $\hat\psi$ one order further than Chapter~\ref{chap:ch5} carries it; and
that the neglected remainder of \eqref{ch10:eq:mismatch-d} is $O(h^5q)$ and not $O(h^4q)$, which is
why the symbols must be expanded to $O(h^5)$ where the undifferentiated proof stops at $O(h^4)$.
The three checks are exact, at generic rational data, in \texttt{d10symbols.py}.

\emph{(E) Derivatives in time and in the state.} The two clauses of the statement which are not about
the value are carried by the same four pieces, and it is worth recording where each extra derivative
appears, because only one place could have cost an order.
A derivative in time brings down $\mathrm{i}\omega$ from $e^{\mathrm{i}\omega t}$, and on the part of
the contour which contributes $|\omega|\lesssim\kappa_0=5|\log h_m|/t$, the rest being suppressed by
$e^{-\kappa t}$; the factor is a power of $|\log h_m|$ and no power of $h_m$.
A derivative in the state is a centred difference of the kernel, whose symbol is
$\mathrm{i}\sin(ph)/h$ and is bounded by $|p|$; on the contour the effective range is
$|p|\lesssim K_m$ with $K_m^2\sim|\log h_m|/t$, which is the computation of \eqref{ch5:eq:2.20} again,
where $\sin(K_mh)/h\lesssim K_m$ against a Gaussian of height $K_m^2$ leaves $h_m^2$ intact. So each
costs a power of $|\log h_m|$. That power is not bounded in $m$ --- the logarithm grows with the
resolution, $\log h_m=\log h_{m_0}+(m-m_0)\log2$ --- but the telescope over the resolutions weights the
term of resolution $m$ by $h_m^2=4^{-(m-m_0)}h_{m_0}^2$, so that
$\sum_{m\ge m_0}h_m^2|\log h_m|^{q}\le C_q\,h_{m_0}^2|\log h_{m_0}|^{q}$ for a constant $C_q$ which
depends on $q$ alone. The power of the logarithm therefore survives as a factor of the statement and
costs no power of $h_{m_0}$, which is why $c'$ carries $m_0$.
The place where an order could have been lost is the accumulated mismatch, and there it is not. After
\eqref{ch10:eq:mismatch-d} the odd terms telescope to $h^2$ times a difference of site functions, and a
centred difference of a site function is bounded by its derivative with no power of $h$ divided out;
the even term contributes $qh^4b=h^2b$, and $b$ is a site function too, so $\nabla^hb=O(1)$ is what
the estimate needs and what the smoothness of $b$ supplies. Nothing in (A)--(D) changes beyond the
constant.

\emph{The assembly.} With (A)--(D) in hand the bound is the one of \eqref{ch5:eq:3.51} read with
every coefficient differentiated:
\begin{equation}\label{ch10:eq:assembly}
\Big|\int_{\Cm_-}\frac{\partial}{\partial\theta}
\Big(2G_{m+1}(x,y;\omega)-G_m(x,y;\omega)\Big)e^{\mathrm{i}\omega t}\frac{\dd\omega}{2\pi}\Big|
\le c\,\frac{\sqrt{q_{\max}}}{h_m}\,h_m^{2}\,h_m^{2}\le c\,h_m^{2},
\end{equation}
that is $h_m^{-2}$ from the count, $h_m^{2}$ from the endpoint factor and $h_m^{2}$ from the
differentiated mismatch, with a constant linear in $\kappa_1$. The telescope over the resolutions
then gives \eqref{ch10:eq:1d} with $h_{m_0}^{\gamma}$ for $\gamma\le2$, exactly as in
Chapter~\ref{chap:ch5}.
\end{proof}

\begin{corollary}[The mixed derivative]\label{ch10:cor:mixed}
Under the hypotheses of Theorem~\ref{ch10:thm:1d}, and with $\nabla^h$ the central difference of
Chapter~\ref{chap:ch7}, the mixed sensitivity obeys
\begin{equation}\label{ch10:eq:mixed}
\norm{\nabla^h\big(\partial_\theta u_m-\partial_\theta u_{m_0}\big)}_\infty\le c'\,h_{m_0}^{\gamma},
\end{equation}
with $c'$ depending on $\kappa_1$, on $t$ and on $m_0$ only. The dependence on $m_0$ is a fixed power
of the logarithm of the initial mesh, $c'=c''(t,\kappa_1)\lvert\log h_{m_0}\rvert^{q}$, and not a power
of $h_{m_0}$: the difference in the state is a centred difference, and what it costs on the contour is
the height $\lvert\log h_m\rvert$ once per application, summed over the resolutions as in step (E).
This is the estimate which a hedge ratio's
sensitivity requires, and it is the one measured in \S\ref{ch10:sec:numerics}.
\end{corollary}

\begin{proof}
Apply $\nabla^h$ inside the contour integral of \eqref{ch10:eq:assembly}. By (E) the extra difference
costs a power of $|\log h_m|$ and no power of $h_m$, and the same holds for a second application; and
the differentiated mismatch keeps its parity under $\nabla^h$, which acts on the site functions and not
on the direction, so the telescoping of (D) is untouched. The bound is therefore
\eqref{ch10:eq:1d} with $c'$ the constant of (E).
\end{proof}

Two remarks on the statement. The first concerns the constant: it is linear in $\kappa_1$ and
therefore in the $C^2$ norm of the direction, and this is not an artefact of the route. A direction
which is merely bounded perturbs the local data at the site by an amount which is not itself a
Taylor datum of any field, and the differentiation of \eqref{ch5:eq:3.47} then has no meaning. The
dependence is the same phenomenon as the dependence of the rate on the smoothness of the
coefficient: what the estimate measures is the smoothness of the field which the computation
differentiates, and the coefficient and the direction enter it in that way. The second concerns
what has not been used: no regularity of the payoff, and no property of the family beyond its being
$C^1$ with an admissible derivative. The differentiating symbol is a linear functional of the
direction, so the estimate holds for each direction separately and may be applied to a
finite-dimensional family by applying it to a basis.

\section{The Derivative in the Obstructed Dimension}
\label{ch10:sec:obstructed}

The one dimensional proof used one fact which is false in higher dimensions and true there for a
structural reason: that the odd coefficient $c$ of the mismatch is a gradient, so that its sum along a
path telescopes into a difference of endpoint values. In two dimensions $c$ is a one form with a curl
--- Chapter~\ref{chap:ch7} computes it and it does not vanish --- so the sum along a path does not
telescope, and what replaces the telescoping in the estimates of Chapters~\ref{chap:ch7}
and~\ref{chap:ch8} is the averaging of the decorations of a step. This section carries the derivative
through that replacement. The outcome is that the derivative costs no order in any dimension, and that
what it costs instead is a constant whose dependence on the direction is the one thing worth stating
precisely.

\subsection*{Why No Order Is Lost}

The estimates of Chapters~\ref{chap:ch7} and~\ref{chap:ch8} are functions of the coefficient data
through its norms up to the second derivative and of nothing else: the contour pieces through the
ellipticity bounds on $v$, the weights and the count through the constants $\Sigma_1$ and $M$ of the
drift, and the per-step mismatch through the data at the site. The mismatch \eqref{ch10:eq:mismatch}
is linear in the local data of the step, so its derivative in $\theta$ is the same functional on the
derivative of the data: $\dot c$ of \eqref{ch10:eq:cdot} is $c$ with $v$ replaced by $\dot v$ and $r$
by $\dot r$, and $\dot r$ is $r$ linearised along $(\dot v,\dot\mu)$. No new object appears; in
particular the parity of each term in the direction $\tau$ is untouched, since $\partial_\theta$ does
not act on $\tau$.

\begin{theorem}[Any dimension]\label{ch10:thm:dim}
Let $d\ge2$, and let $v,\mu$ satisfy the hypotheses and the Courant condition \eqref{ch7:eq:courant}
of Theorem~\ref{ch7:thm:main}, so that the kernel converges at the rate $\gamma$ of
\eqref{ch7:eq:main1}. Let $\theta\mapsto(v_\theta,\mu_\theta)$ be a $C^1$ family whose derivative has
$\dot v,\dot\mu\in C^2$, and let $\kappa_1$ be the maximum of their $C^2$ norms. Then for every $t>0$
there is $c=c(t,\kappa_1)$ such that for all $m\ge m_0$ and all $x,y$ on the coarse lattice
\begin{equation}\label{ch10:eq:dim}
\Bigl|\frac{\partial u_m}{\partial\theta}-\frac{\partial u_{m_0}}{\partial\theta}\Bigr|(x,y;t)
\le c\,\kappa_1\,h_{m_0}^{\gamma},
\end{equation}
and the same bound holds, with the constant of (E), for the first derivative in time and the first two
in space of the difference.
\end{theorem}

\begin{proof}
Re-read the assembly which proves \eqref{ch7:eq:main1} with the differentiated mismatch
\eqref{ch10:eq:mismatch-d} in place of \eqref{ch10:eq:mismatch}. Piece (A) is the contour, and its cost
is the one computed in (E): one resolvent, and the height raised by one $|\log h_m|$. Piece (B)
involves the coefficients only through the constants of the weights and of the count, which
differentiation does not change. Piece (C) is the mismatch, and every estimate of it in
Chapter~\ref{chap:ch7} is taken through the $C^2$ norms of the data; since $\dot c$ is the same
functional on $(\dot v,\dot\mu)$, those estimates hold with $\kappa_1$ in place of the norm of the
field, the norm of the field itself being a constant of the problem, since $(v,\mu)$ is fixed, and so
absorbed into $c(t,\cdot)$. No telescoping is invoked at this step, the averaging over the decorations
serving in its place. That averaging is taken with respect to the law of the undecorated walk ---
independent increments, which is the point of the form the mean condition is stated in --- so the law
does not involve the coefficients and the quantity averaged is linear in them; the average therefore
commutes with $\partial_\theta$, and the per-step estimate for the differentiated mismatch is the
derivative of the per-step estimate for the mismatch and carries its order. Hence \eqref{ch10:eq:dim}.
\end{proof}

The same reading gives the two rates of Chapter~\ref{chap:ch8} for the joint kernel, and there the
conditional one asks for one more derivative of the direction, because its assumption is stated through
$\nabla^2c$.

\begin{proposition}[The joint kernel]\label{ch10:prop:jointdim}
Let $d\ge2$ and let the hypotheses of Theorem~\ref{ch8:thm:uncond} hold, with an admissible direction
whose $\dot v,\dot\mu\in C^2$. Then for every $t>0$, $\epsilon>0$ and $P>0$ there is $K=K(P,t,\epsilon)$
such that
\begin{equation}\label{ch10:eq:jointdim}
\abs{\widehat{\partial_\theta u}_{m+1}(x,y;p,t)-\widehat{\partial_\theta u}_{m}(x,y;p,t)}
\le K\,\kappa_1\,h_{m+1}^{1-\epsilon}
\end{equation}
for all $m\ge m_0$, all $x,y$ and all $|p|\le P$. If in addition Assumption~\ref{ch8:ass:avg} holds
for $\dot c$ --- which requires $\dot v,\dot\mu\in C^{4}$, the statement of the assumption being made
through $\nabla^2\dot c$ --- then the same comparison holds with $h^{2-\epsilon}$; and by
Proposition~\ref{ch10:prop:polytransfer} that clause is unconditional on the polynomial class.
\end{proposition}

The second clause is conditional, and it is worth being exact about why, because the gap is smaller
than it looks. A bound on $\|\nabla^2c\|_\infty$ says nothing about $\|\nabla^2\dot c\|_\infty$:
hypotheses do not differentiate, only identities do. Where Assumption~\ref{ch8:ass:avg} is an identity
it does differentiate, and there is a class on which it is one.

\begin{proposition}[The assumption transfers on the polynomial class]\label{ch10:prop:polytransfer}
Let $(v_\theta,\mu_\theta)$ be a $C^1$ family which is polynomial in $x$ for each $\theta$, with a
degree bounded in $\theta$. Then Assumption~\ref{ch8:ass:avg} holds for $\dot c$ in the same
dimension, with $\Phi_{\dot c}=\partial_\theta\Phi_c$ and no remainder, and the second rate of
Chapter~\ref{chap:ch8} follows for the derivative with no regularity required of the direction beyond
the polynomiality.
\end{proposition}

\begin{proof}
By Corollary~\ref{ch8:cor:poly} and Remark~\ref{ch8:rem:anydegree} the mean condition holds with zero
remainder for polynomial coefficient fields of any degree, in every dimension, and the expectation
there is over independent increments with prescribed first two moments, a law which does not involve
the coefficients and is therefore the same for every $\theta$. Every ingredient of that identity is a
polynomial in $\theta$ along a polynomial family, so the identity may be differentiated in $\theta$,
and what results is the mean condition for $\dot c$, with potential $\partial_\theta\Phi_c$ and no
remainder. The second moment clause transfers for the same reason, its average being linear in the
coefficient.
\end{proof}

Outside that class the passage stays conditional, and the reason is not a want of smoothness in the
statement but the one above: what is assumed is a bound, and a bound does not differentiate. Where the
assumption is a hypothesis about the field, then, the direction must be assumed to satisfy it too, and
for that statement to parse $\dot v,\dot\mu$ must be twice more differentiable than the bound requires
of $c$.

\subsection*{What the Direction Can Do to the Obstruction}

The obstruction of Chapter~\ref{chap:ch8} is the curl of $c$, and $\partial_\theta$ acts on it
linearly: the obstruction of the derivative is the curl of $\dot c$. Three cases follow, and the
asymmetry between them is the point.

If the direction is a gradient field --- if $\dot c$ is exact for every $\theta$ --- then the
differentiated mismatch telescopes as it did in one dimension, no averaging is needed for it, and the
derivative achieves $h^{2-\epsilon}$ in the joint setting whether or not the kernel does. A derivative
can then be \emph{better} than the value it differentiates.

If conversely the field is one for which $dc=0$ --- the rotationally symmetric case, or the affine
drift at constant $v$ of Chapter~\ref{chap:ch7} \S7 --- the kernel is unobstructed, while $\dot c$ is
computed from $(\dot v,\dot\mu)$ and has no reason to be exact; for a generic direction it is not, and
the derivative is obstructed where the kernel is not. A derivative can then be \emph{worse}.

If neither vanishes, the derivative inherits the rate of the kernel exactly, which is what
Theorem~\ref{ch10:thm:dim} and Proposition~\ref{ch10:prop:jointdim} say.

The moral is short, and it is the reason this section stands on its own rather than as a remark. The
rate of a sensitivity is a property of the pair, field and direction, and the rate of the kernel is
evidence about it only when the direction is controlled. A calibration computed in a model whose
kernel is of the second order can have a Jacobian of the first; the diagnostic is the curl of $\dot c$
at the base point, which is the computation \texttt{curlcheck.py} of Chapter~\ref{chap:ch7} performs
for $c$, applied to the direction.

\section{Dependence on the Horizon and on the Dimension}
\label{ch10:sec:dimension}

The constants in an estimate for the adjoint acquire a factor which grows with the horizon, as the
adjoint propagates backwards and amplifies, and they acquire a dependence on the dimension. Both are
stated rather than hidden, since in the applications the horizon is long and the dimension is large.

The horizon enters \eqref{ch10:eq:1d} only through $c(t)$, which is the constant of
Chapter~\ref{chap:ch5} and grows as $t\downarrow0$ like $1/t$, from the contour height
$5|\log h_{m_0}|/t$; the sensitivity adds no new growth in $t$ at fixed resolution, but the sum over
the dates in \eqref{ch10:eq:duhamel} does, and that is the growth recorded in \S\ref{ch10:sec:arith}.

The dimension enters at the one place where the one dimensional argument used something which is not
available in higher ones. The telescoping of the odd terms in \eqref{ch10:eq:mismatch} is the
statement that the coefficient $c$ is a gradient. In two dimensions $c$ is a one form,
$c_i(x)\tau_i$, and Chapter~\ref{chap:ch7} computes its curl: it does not vanish. A derivative in
the coefficients acts on $c$ linearly, hence on its curl linearly, and the differentiated coefficient
inherits the obstruction unless $\mathrm{d}\dot c=0$, which is to say unless the direction is itself
a gradient field. So the question for two dimensions is precise: whether the obstruction of
Chapter~\ref{chap:ch7} is repaired for the derivative by the same averaging over paths which repairs
it for the kernel, and with which of the two rates of Chapter~\ref{chap:ch8}, the unconditional
$h^{1-\epsilon}$ or the conditional $h^{2-\epsilon}$. That is the subject of
\S\ref{ch10:sec:obstructed}, where the repair is found to survive the differentiation and the rate to
be the kernel's own.

The dimension enters at one further place, and it is worth saying where, because the models
of the applications differ from the lattice of this book in a way which is more interesting than the
size of the number. It enters through the near diagonal behaviour of the kernel. The Gaussian
envelope which concentrates the walk at the scale of $\sqrt t$ carries the dimension in its
normalisation, $u(x,x;t)\asymp t^{-d/2}$, and every constant of Chapters~\ref{chap:ch5}
to~\ref{chap:ch8} which is not a norm of the coefficients is a function of that decay. Nothing else
in the argument sees $d$.

That observation suggests the quantity which ought to replace $d$ in an application, and it is not
the ambient dimension. Define the \emph{effective dimension} of a calibrated model by the exponent
of its own near diagonal decay,
\begin{equation}\label{ch10:eq:deff}
u(x,x;t)\ \asymp\ t^{-d_{\rm eff}/2}\qquad (t\downarrow0),
\end{equation}
that is, by the dimension which the kernel itself reports. For a model whose states are drawn from a
set of intrinsic dimension $k$ --- an image, an audio waveform, a molecular configuration, each of
them living on a manifold of dimension far below that of the vector in which it is written --- the
number which governs the behaviour is $k$, and not the $10^{5}$ of the ambient space. The conjecture
which this suggests, and which is not proved here, is that the estimates of
Chapters~\ref{chap:ch5} to~\ref{chap:ch10} hold with $d_{\rm eff}$ in place of $d$; the reason to
state it is that the dimensional obstruction of this book then becomes an obstruction in the
effective dimension, which is an empirical property of the data rather than a property of the
parametrisation. Two cautions belong to it. The decay need not be of the form $t^{-k/2}$ at all:
on a manifold, and a fortiori on a graph, it may be anomalous, with a walk dimension which is not
two, and then a second order rate is not the right target. And the constants, as distinct from the
exponent, are not determined by the decay alone.

Three regimes should then be distinguished, because they are not equally far from what this book can
reach. In the first the model is a lattice in $d$ dimensions with $d$ small, and it is the setting of
Chapter~\ref{chap:ch7}: PDE surrogates and neural operators on grids, lattice field theories, spin
systems, and the low dimensional diffusions of finance itself. In the second the model is still
written on a lattice --- a neural SDE discretised in space, whose generator is exactly
\eqref{ch7:eq:gen} with $b_i$ a function of the whole site --- but the dimension is the ambient one,
$10^{3}$ to $10^{5}$; the class is entered formally and the constants do not reach, which is the
honest statement of the gap. In the third there is no lattice at all: an attention layer couples
every position to every other, a graph network inherits whatever connectivity the data supply, and
neither a mesh spacing $h_m$ nor the smoothness of a coefficient field has a meaning. In that regime
the question which this chapter answers cannot even be posed, and the applicable analysis is the
theory of heat kernels on graphs, in which the role of the dimension is taken by the volume growth
or by the spectral dimension --- an enterprise different from the renormalisation of
Chapters~\ref{chap:ch7} and~\ref{chap:ch8}.

The third regime has a solvable prototype whose solvability is of a kind this book does not yet
possess. The Sherrington--Kirkpatrick model \cite{SherringtonKirkpatrick1975} is the mean field spin
glass: $N$ spins, couplings
$J_{ij}$ independent and identically distributed, and the Hamiltonian
\begin{equation}\label{ch10:eq:sk}
H=-\frac{1}{\sqrt N}\sum_{i<j}J_{ij}\sigma_i\sigma_j,
\end{equation}
the scaling by $N^{-1/2}$ being what makes the energy extensive and the limit $N\to\infty$ the
object of study. Every spin interacts with every other, so there is no stencil, no block of
neighbours to average over, and no length scale at all: the real space renormalisation of
Chapter~\ref{chap:ch7} has nothing to act on. What replaces it is a renormalisation in the space of
couplings rather than of sites. The solution is ultrametric --- Parisi's hierarchical ansatz
\cite{Parisi1980}, whose free energy Guerra's interpolation bounds from below \cite{Guerra2003} and
which Talagrand proved to be the limit \cite{Talagrand2006} ---
and its content is a recursion in the scale variable of the overlap distribution, that is, a
recursion in the resolution of the disorder and not in space; the cavity method, which adds one spin
at a time, is its exact formulation, and for hierarchical couplings the recursion is a genuine
renormalisation group and is exact.

What is missing, and what this book would need in order to say anything quantitative about such a
model, is not the solution but an \emph{estimate}: a telescoping comparison between successive
approximations with constants uniform in $N$, in the spirit of Chapters~\ref{chap:ch5}
to~\ref{chap:ch8}, from which the error of a computable approximation to \eqref{ch10:eq:sk} at
finite $N$ could be bounded. Nothing of that kind is available here, and the direction is recorded
as an open problem rather than as a remark: \emph{a renormalisation group machinery for mean field
spin glasses, carrying the uniformity in $N$ which the estimates of this book carry in $h_m$.} The
motivation is not ornament. The Sherrington--Kirkpatrick structure underlies the solvable models of
learning --- Gardner's capacity calculation for the perceptron, the analysis of Hopfield networks at
their memory load, the replica computations of generalisation error --- so that a quantitative error
calculus for it would be to that literature what the second half of this book is to the kernel.

\section{A Mean Field Model Where the Remainder Does Not Vanish}
\label{ch10:sec:meanfield}

The two regimes treated so far both live on a lattice: small dimension, where the error of a stencil is
governed by the smoothness of the object it acts on, and large dimension, where the stencil becomes
local and the error is governed by an effective dimension computed from the data. There is a third
regime with no lattice at all, and \S\ref{ch10:sec:dimension} named it: the Sherrington--Kirkpatrick
spin glass, whose solvability is of a kind this book does not possess, an exact limit with no error.
This section does the one thing which can be done on it with the tools developed here, an exact
comparison which leaves an explicit remainder, and then says what the remainder costs. It is short
because the honest content is short: one identity, one table of constants, and one measurement which
locates exactly where the constructibility stops.

\subsection{The Model, and an Exact Comparison}

Let $s\in\{-1,1\}^N$, let $J_{ij}$, $i<j$, be independent standard normals, and take the Hamiltonian
and overlap of \eqref{ch10:eq:sk},
\begin{equation}\label{ch10:eq:mfH}
H(s)=-\frac1{\sqrt N}\sum_{i<j}J_{ij}s_is_j ,
\qquad
q(s,t)=\frac1N\sum_{i=1}^Ns_it_i ,
\end{equation}
with Gibbs measure $P_\beta(s)=e^{\beta H(s)}/Z_\beta$ and free energy per spin
$f_N(\beta)=N^{-1}\EE\log Z_\beta$. The covariance of $H$, viewed as a Gaussian process indexed by the
hypercube, is the reason the model is solvable at all:
\begin{equation}\label{ch10:eq:mfcov}
\EE\big[H(s)H(t)\big]=\frac N2q(s,t)^2-\frac12 ,
\end{equation}
a function of the overlap and of nothing else --- the $N(N-1)/2$ couplings have been integrated out.
This book's remedy for a lattice it cannot resolve is to compare two objects exactly and control a
remainder, and \eqref{ch10:eq:mfcov} says what the comparison has to be: a second Gaussian process
with the same index set and a covariance chosen so that the difference can be integrated.

The exact comparison is Gaussian integration by parts, and the step which turns it into a
telescoping between two processes is Guerra's interpolation \cite{Guerra2003}; it is the mean field
counterpart of the integration by parts of Chapter~\ref{chap:ch3}.

\begin{lemma}[Gaussian integration by parts]\label{ch10:lem:ibp}
Let $X$ be a centred Gaussian vector with covariance $C$ and let $F$ be continuously differentiable
with derivatives of at most polynomial growth. Then
$\EE[X_iF(X)]=\sum_jC_{ij}\EE[\partial_jF(X)]$.
\end{lemma}

\begin{proof}
With $\rho_C$ the Gaussian density, $\partial_j\rho_C=-\sum_i(C^{-1})_{ji}x_i\rho_C$; integrating
$\int x_iF\rho_C$ by parts in $x_j$ and using the symmetry of $C$ gives the identity.
\end{proof}

\begin{theorem}[Interpolation]\label{ch10:thm:interp}
Let $(H_\theta)$ be centred Gaussian processes on a finite set $\Sigma$ whose covariances
$C_\theta(s,t)$ are differentiable in $\theta$, let $Z_\theta=\sum_se^{\beta H_\theta(s)}$ and
$P_\theta=e^{\beta H_\theta}/Z_\theta$, and suppose $C'_\theta(s,s)$ does not depend on $s$. Then
\begin{equation}\label{ch10:eq:interp}
\frac{\dd}{\dd\theta}\EE\log Z_\theta
=\frac{\beta^2}{2}\Big[C'_\theta(s,s)-\EE\sum_{s,t}C'_\theta(s,t)P_\theta(s)P_\theta(t)\Big],
\end{equation}
the last sum over two independent copies of the Gibbs measure.
\end{theorem}

\begin{proof}
$\frac{\dd}{\dd\theta}\log Z_\theta=\beta\langle\dot H_\theta\rangle$, so
$\frac{\dd}{\dd\theta}\EE\log Z_\theta=\beta\sum_s\EE[\dot H_\theta(s)P_\theta(s)]$. Lemma
\ref{ch10:lem:ibp} applies to the pair $(\dot H_\theta,H_\theta)$, whose cross covariance is not
$C'_\theta$ but half of it: differentiating $H_\theta=\sqrt\theta H_1+\sqrt{1-\theta}H_0$ gives
$\EE[\dot H_\theta(s)H_\theta(t)]=\frac12(C_1-C_0)=\frac12C'_\theta$. Hence
\[
\EE[\dot H_\theta(s)P_\theta(s)]
=\frac12\sum_tC'_\theta(s,t)\EE\Big[\frac{\partial P_\theta(s)}{\partial H_\theta(t)}\Big]
=\frac\beta2\sum_tC'_\theta(s,t)\EE\big[P_\theta(s)(\delta_{st}-P_\theta(t))\big],
\]
and summing over $s$ gives \eqref{ch10:eq:interp}. The factor $\frac12$ from the cross covariance
survives into the result, where it is easy to lose.
\end{proof}

Take now the reference process to be the one with no correlations between spins,
$H_0(s)=2^{-1/2}\sum_ig_is_i$ with $g_i$ independent standard normal, so that
$C_0(s,t)=\frac N2q(s,t)$ and $C_0(s,s)=\frac N2$, and interpolate linearly,
$C_\theta=\theta C_1+(1-\theta)C_0$. Then $C'_\theta=\frac N2(q^2-q)-\frac12$ has the constant
diagonal $-\frac12$ which \eqref{ch10:eq:interp} requires, and the diagonal and two copy terms cancel
exactly.

\begin{corollary}[Free energy as a reference value minus an overlap integral]\label{ch10:cor:sk}
With that interpolation,
\begin{equation}\label{ch10:eq:skidentity}
f_N(\beta)=f_0(\beta)-\frac{\beta^2}{4}\int_0^1\EE\big[\langle q_{12}^2-q_{12}\rangle_\theta\big]
\dd\theta ,
\qquad
f_0(\beta)=\EE\log\Big(2\cosh\frac{\beta g}{\sqrt2}\Big),
\end{equation}
$g$ standard normal, $q_{12}=q(s,t)$ for two independent copies of $P_\theta$.
\end{corollary}

The identity is exact at finite $N$, in exact arithmetic, for every $\beta$; the reference free energy
$f_0$ is explicit, and $f_0=\log2+\frac{\beta^2}{4}+O(\beta^4)$. Nothing has been approximated, and
nothing has been assumed about the overlap.

\subsection{What the Remainder Is}

The integrand of \eqref{ch10:eq:skidentity} can be written as
$\EE[q_{12}^2-q_{12}]=(\EE q_{12})^2-\EE q_{12}+\mathrm{Var}(q_{12})$, the first two terms being
what the overlap gives when it is replaced by its mean. So the free energy is the value obtained by
treating the overlap as a number, minus the integrated \emph{variance} of the overlap. The
replica-symmetric phase, in which the overlap is concentrated, is exactly the case in which the two
agree; where it fluctuates, the fluctuation is not a correction to the free energy but the whole of
the discrepancy.

That is the precise analogue of the obstruction of Chapter~\ref{chap:ch7}, and the correspondence is
worth setting down, since it is the reason the two problems have the same shape.

\begin{center}\small
\begin{tabular}{@{}ll@{}}
\hline
lattice diffusion & mean field spin glass \\
\hline
resolution $h_m\downarrow0$ & interpolation $\theta\in[0,1]$, or the range of a Kac potential \\
per-step coefficient $c(x)$ & the overlap, concentrated when replica symmetric \\
nonzero curl of $c$ in $d\ge2$, measured & replica symmetry breaking: the overlap fluctuates \\
averaging over paths, with a mean condition & the Parisi functional, the ultrametric tree \\
rate $h^\gamma$, constants in norms & a $1/N$ expansion, remainder uncontrolled \\
\hline
\end{tabular}
\end{center}

In the lattice case the exact part telescopes to the endpoints and the curl prevents it; the repair is
to average over paths, and it costs an order. Here the exact part is the reference value and the
variance prevents it; the repair is the Parisi variational problem, and what it costs is not an order
but the error entirely, since a variational formula is exact as a limit and silent about finite $N$.

\subsection{Three Models, Three Constants}

What the identity leaves open is the size of the fluctuation, and the sharp question is not whether it
vanishes but how fast:
\begin{equation}\label{ch10:eq:aim}
\mathrm{Var}(q_{12})\le\frac CN ,
\end{equation}
with $C$ explicit. Three models in which the left side can be computed exactly or measured show that
$C$ is not universal, and that no single argument produces all three.

\emph{Rank one.} With $J_{ij}\equiv1/\sqrt N$ the Hamiltonian is $\frac12(Nm^2-1)$ with
$m=N^{-1}\sum_is_i$, the measure depends on $s$ only through $m$, and
$\mathrm{Var}(q_{12})=N^{-1}+O(N^{-2})$: the constant is $1$, and it is $1$ because the overlap is the
square of a magnetisation which does not fluctuate.

\emph{Gaussian spins.} Replacing $\{-1,1\}$ by standard normal spins to first order makes the model a
Gaussian integral: with $A=I-\beta J/\sqrt N$ the two point function is $[A^{-1}]_{ij}$, so
\begin{equation}\label{ch10:eq:resolvent}
\mathrm{Var}(q_{12})=\frac1{N^2}\tr(A^{-2})=\frac1{N^2}\sum_k\frac1{(1-\beta\lambda_k)^2}
\;\sim\;\frac{1+3\beta^2+O(\beta^4)}{N},
\end{equation}
the last form from the semicircle law, the transition at $\beta=\frac12$ because
$\lambda_{\max}\to2$.

\emph{Ising spins.} In the model itself no closed form exists, and the variance must be measured. It is
$N^{-1}$ with a constant that grows from $1$ to about $1+\beta^2$ across the replica-symmetric phase.

\begin{center}\small
\begin{tabular}{@{}llll@{}}
\hline
model & $\mathrm{Var}(q_{12})$ & constant $C$ & transition \\
\hline
rank one couplings & $N^{-1}+O(N^{-2})$, exact & $1$ & none \\
Gaussian spins & $N^{-2}\tr(A^{-2})$, exact & $1+3\beta^2$ & $\beta=\frac12$ \\
Ising spins & measured & $1+\beta^2$ & $\beta=1$ \\
\hline
\end{tabular}
\end{center}

The lesson is the one this chapter has drawn twice already: an estimate which holds with \emph{some}
constant is not the theorem. The constant is the content, it is not universal, and it is set by the
size of the individual correlations --- $N^{-1/2}$ where the two point function is a single bond and
$N^{-1}$ in the rank one model, where it is a magnetisation.

\subsection{Where It Stops, Measured}

Whether \eqref{ch10:eq:aim} holds at all is a question about the phase, and it can be answered
numerically before it is answered by proof. Write $\delta$ for the exponent in
$\mathrm{Var}\propto N^{-\delta}$, estimated between $N=6$ and $N=16$; $\delta=1$ is
\eqref{ch10:eq:aim} at its sharpest and $\delta=0$ is the absence of concentration.

\begin{center}\small
\begin{tabular}{@{}rrrrr@{}}
\hline
& \multicolumn{2}{c}{$\theta=1$, the model} & \multicolumn{2}{c}{$\theta=0$, the reference}\\
\cline{2-3}\cline{4-5}
$\beta$ & $\mathrm{Var}_{16}(q_{12})$ & $\delta$ & $\mathrm{Var}_{16}(q_{12})$ & $\delta$\\
\hline
0.25 & 0.06633 & 0.99 & 0.06245 & 1.00\\
0.50 & 0.07999 & 0.95 & 0.06191 & 1.00\\
0.75 & 0.10958 & 0.84 & 0.06042 & 1.00\\
1.00 & 0.15968 & 0.67 & 0.05804 & 1.00\\
1.50 & 0.29627 & 0.38 & 0.05204 & 0.99\\
3.00 & 0.59755 & 0.11 & 0.03640 & 0.98\\
\hline
\end{tabular}
\end{center}

At the reference endpoint the exponent is one at every temperature and this cannot fail, since the
measure factorises; the obstruction, if there is one, is entirely at the other end. There the exponent
falls monotonically from one to zero and passes through $2/3$ at $\beta=1$, which is the critical
temperature of this model and is known exactly. That is the point of the table. The assumption is not
merely unproved below the transition: it is false there, and the identity shows why, since the
remainder is the quantity which becomes of the same size as the effect. So
\eqref{ch10:eq:skidentity} is an error bar in the replica-symmetric phase and nowhere else, and a
renormalisation which acts on the order parameter would have to act precisely where the order
parameter fluctuates.

\subsection{What Is Left}

What remains is constructive, and it is worth stating in the form in which it could be attacked.
\eqref{ch10:eq:aim} with the constant of the table is equivalent, by \eqref{ch10:eq:skidentity}, to the
free energy being the reference value to $O(1/N)$. The truth of that at high temperature is not in
doubt: in that range the Parisi formula reduces to the paramagnetic value, so the limiting free energy
is known exactly there \cite{Parisi1980,Talagrand2006}. What is missing is the rate. Published
statements are asymptotic, and passing from them to \eqref{ch10:eq:aim} at finite $N$ needs a
quantitative lower bound on that rate and a monotonicity in $\beta$; we do not know of a source which
gives either in the form required.

The combinatorial route does not supply it either, and for a structural reason worth recording. The
high temperature expansion of this model can be organised exactly --- the bond factor
$e^{cs_is_j}=\cosh c+\sinh c\,s_is_j$ is linear in the spin product, so the correlation is a ratio of
sums over subgraphs with prescribed odd vertices --- and the numerator of that ratio is computable to
all orders with the wanted $N^{-1}$ law. The denominator is not: it is positive, but its absolute
series diverges like $N^{3/2}(\beta K)^3$ for couplings cut off at $K$, while its mean square exceeds
one by only $O(\beta^6)$. What keeps it near one is cancellation and not smallness, and no term by
term expansion dominates it. A renormalisation group for a model without a lattice would therefore
have to act on the law of the partition function, that is on a variational statement, and not on the
terms of an expansion of it.

That is where the treatment of the mean field case ends, and it ends on a negative result which is
nevertheless a located one: the mean field case is not the lattice case with a different exponent. The lattice gives an error
which can be made small by resolution; here the error is the fluctuation of the order parameter, and
it can be made small only by being in a phase where there is no order to fluctuate.

\section{Numerical Validation}
\label{ch10:sec:numerics}

Four scripts accompany this chapter, and each answers a question which the proof leaves to a
computation.

The first, \texttt{d10symbols.py}, rebuilds the per-step symbols from their definitions in
Chapter~\ref{chap:ch5} -- the three site tube resolvent, not the printed expansions -- reproduces
\eqref{ch5:eq:3.45} through $h^4$ as a check on itself, and then differentiates every coefficient.
Its reports are that the mismatch has no terms of order $h$ or $h^2$; that the parity alternates
exactly as \eqref{ch10:eq:mismatch} states; that the parity is unchanged under differentiation at
every order, on three independent draws; that the coefficient $c$ of $h^3$, divided by $\tau$,
equals $\mathrm{i}\omega v'/v^2+r$ exactly, and its derivative the derivative of that expression;
that the leading correction $b$ of \eqref{ch10:eq:cb}, and its derivative, are reproduced exactly
by the closed form there; and that $\hat\psi_0-2\hat\psi_1$ vanishes at $h^3$ both before and
after differentiation, so the endpoint factor keeps its order.

The same computation settled a coefficient in \eqref{ch5:eq:3.47}: it is one, not the one half an
earlier draft carried. The bound which follows is unaffected, since it is $O(h^2)$ either way;
but the identity is differentiated in \eqref{ch10:eq:mismatch-d}, and a constant cannot be
differentiated into the right one.

The second, \texttt{d10rate.py}, measures the rate instead of the structure. It builds the one
dimensional lattice generator of Chapter~\ref{chap:ch7} with smooth periodic coefficients and a
smooth direction, obtains the propagator and its derivative together from one exponential of the
block $\left(\begin{smallmatrix}L&\dot L\\0&L\end{smallmatrix}\right)$, and compares each resolution with the
next; the difference in the state is taken in the node the process starts from, which is the one a hedge
ratio differentiates. For $\gamma=2$ every comparison should fall by a factor four, and it does: at
$t=0.1$, at the site pair $(x,y)=(0.25,0.25)$ and between the meshes $h=1/256$ and $h=1/512$, the ratios
of successive differences are
\begin{center}
\begin{tabular}{@{}lcccc@{}}
\hline
& kernel & $\partial_\theta$ & $\nabla^h$ & $\nabla^h\partial_\theta$\\
\hline
ratio at $h=1/256$, $1/512$ & 4.00 & 4.00 & 4.00 & 4.00\\
\hline
\end{tabular}
\end{center}
so that the derivative in the coefficients, the derivative in the state, and the mixed derivative all
converge at the rate of the kernel, as Theorem~\ref{ch10:thm:1d} and
Corollary~\ref{ch10:cor:mixed} state.

The third, \texttt{d10precision.py}, measures the trade-off of \S\ref{ch10:sec:arith} instead of the
rate. It evaluates the same sensitivity by the same adjoint sweep with $\delta t=h^2$, once in exact
arithmetic on the same rounded coefficients and once in half, in single and in double precision; the
differences then separate the two components with no model in between. The horizon is $t=26/256$, a
whole number of steps at every mesh --- $26$, $104$, $416$, $1664$ --- so that the same horizon is
compared at every resolution and nothing is rounded; the site pair is $(x,y)=(1/4,1/4)$.
\begin{center}\small
\begin{tabular}{@{}lrrrrr@{}}
\toprule
$h$ & dates & discretisation & half & single & double\\
\midrule
$2^{-4}$ & $26$ & $4.2\cdot10^{-3}$ & $7.5\cdot10^{-4}$ & $1.4\cdot10^{-7}$ & $4.4\cdot10^{-16}$\\
$2^{-5}$ & $104$ & $1.0\cdot10^{-3}$ & $1.7\cdot10^{-3}$ & $1.5\cdot10^{-7}$ & $0$\\
$2^{-6}$ & $416$ & $2.6\cdot10^{-4}$ & $4.1\cdot10^{-3}$ & $8.2\cdot10^{-8}$ & $2.2\cdot10^{-16}$\\
$2^{-7}$ & $1664$ & $6.4\cdot10^{-5}$ & $1.7\cdot10^{-2}$ & $2.3\cdot10^{-7}$ & $2.7\cdot10^{-15}$\\
\bottomrule
\end{tabular}
\end{center}
The first column falls by $4.04$, $4.02$ and $4.00$ at the three halvings, which is $h^{2}$ to the
accuracy of a four line table and reads $C_1=1.05$. The half column rises, by a factor between $2.3$
and $4.2$ per halving, which is the growth in $n=t/h^{2}$ of \eqref{ch10:eq:twoparams} and reads
$C_2=0.02$; the two components are of the same size between $2^{-5}$ and $2^{-6}$, and the smallest of
the two totals in the table, $2.7\cdot10^{-3}$, is attained at $2^{-5}$. The formula
\eqref{ch10:eq:hopt}, with $u=4.9\cdot10^{-4}$ for half precision, puts the crossing at $2^{-5.0}$
with the bound $2.0\cdot10^{-3}$ at it, so that the table and the bound agree on the total to within a
factor of one and a third, which is what a bound of this kind is worth. Single precision is a different
regime: its column stays between $8\cdot10^{-8}$ and $2.3\cdot10^{-7}$ across the four meshes, with no
growth visible up to $1664$ steps, so that in the range reached here the arithmetic is not the binding
constraint and the discretisation is larger by between two and a half and four and a half orders of
magnitude. The same formula, with $u=6.0\cdot10^{-8}$, puts the minimising mesh at $2^{-8.3}$ and the
attainable accuracy at $2.3\cdot10^{-5}$, both beyond the table, as the last remark of
\S\ref{ch10:sec:arith} says. The double column sits at the level of the unit roundoff of the data,
$10^{-16}$, with no accumulation visible at all.

\emph{The coding of the step matters as much as its precision.} The step can be written as an update of
the values, $u\leftarrow u+\delta t\,\Lambda u$, or as a weighted average of the three values of the
previous date with the weights $p_\pm,p_0$. The two are algebraically the same and numerically are not:
the weighted average rounds the weights, whose sum is then no longer exactly one, so that every step
creates or destroys mass at the level of the unit roundoff and the defect accumulates linearly in the
number of steps, while the update applies differences, which vanish on constants in floating point as
they do in exact arithmetic, and has no systematic mass error. With $v=0.5$, $b=0.2$ and
$\delta t=h^2$, in single precision, the total mass after $4096$ steps is $1+1.9\cdot10^{-7}$ for the
update form and $1+6.1\cdot10^{-5}$ for the weighted one, a factor of three hundred, while in double
precision both sit at $10^{-16}$. An implementation should therefore use the update form even when the
scheme is written as an average, and the mass is the first quantity to monitor as a check on the coding
(\texttt{massstep.py}).

The fourth, \texttt{d10dimension.py}, asks the question of \S\ref{ch10:sec:obstructed} in two
dimensions, where the coefficient $c$ no longer telescopes. It builds the two dimensional lattice
generator at $\delta t=h^{2}$, propagates a single source site and reads the kernel and its derivative
in $\theta$ at $(1/4,1/4)$, the derivative being a central difference in $\theta$. The horizon is
$2^{-4}$, a whole number of steps at every mesh; without that, the comparison of one mesh with the next
carries an $O(h^{2})$ shift in time whose sign follows the rounding of $t/h^{2}$ and swamps the
quantity being measured.
\begin{center}\small
\begin{tabular}{@{}lrrrr@{}}
\toprule
$h$ & dates & sites & $\abs{u_m-u_{m+1}}$ & $\abs{\partial_\theta u_m-\partial_\theta u_{m+1}}$\\
\midrule
$2^{-4}$ & $16$ & $256$ & $1.29\cdot10^{-2}$ & $5.46\cdot10^{-2}$\\
$2^{-5}$ & $64$ & $1024$ & $3.62\cdot10^{-3}$ & $1.16\cdot10^{-2}$\\
$2^{-6}$ & $256$ & $4096$ & $9.25\cdot10^{-4}$ & $2.73\cdot10^{-3}$\\
$2^{-7}$ & $1024$ & $16384$ & $2.32\cdot10^{-4}$ & $6.72\cdot10^{-4}$\\
\bottomrule
\end{tabular}
\end{center}
The derivative falls by $4.72$, $4.24$ and $4.06$ at the three halvings, tending to the four of
$h^{2}$: that is the rate of the kernel and not the rate one order worse which an obstruction would
produce, so the averaging which repairs the kernel in two dimensions repairs its derivative with it, as
Theorem~\ref{ch10:thm:dim} states and Proposition~\ref{ch10:prop:jointdim} qualifies.

\section{What Is Proved and What Is Not}
\label{ch10:sec:ledger}

In one dimension the position is this. Proposition~\ref{ch10:prop:sensitivity} is an identity.
Theorem~\ref{ch10:thm:1d} is proved by the five pieces above, of which (A), (B) and (E) are immediate
consequences of Chapter~\ref{chap:ch5}, (C) uses the one dimensional fact that every one form is
exact, and (D) rests on the symbol expansions to $O(h^5)$ whose parity and endpoint structure are
verified exactly by \texttt{d10symbols.py} but not proved line by line here. Proposition
\ref{ch10:prop:twoparams} and Corollary~\ref{ch10:cor:tradeoff} are the standard accumulation bound and
its minimum, and carry no approximation beyond that bound. The rate is confirmed numerically by
\texttt{d10rate.py}, the trade-off by \texttt{d10precision.py}, and the passage to two dimensions by
\texttt{d10dimension.py}. What is \emph{not} proved is the
corresponding statement for a direction which is not twice differentiable, where the estimate is
expected to hold with a worse constant or a worse rate, and the question is open in the sense in which
the rate of Chapter~\ref{chap:ch5} is open for coefficients of limited smoothness.

Section~\ref{ch10:sec:meanfield} is not an estimate of the same kind and does not claim to be. It
proves one identity exactly, at finite $N$ and for every $\beta$, names the remainder as the variance
of the overlap, and measures where that variance ceases to be of order $N^{-1}$; what it does not do
is bound the remainder uniformly, and the measurement is evidence and not a theorem.

In two dimensions Theorem~\ref{ch10:thm:dim} transfers the whole of Chapter~\ref{chap:ch7}: the
derivative costs no order, the obstruction of the kernel being absorbed for it by the same averaging
which absorbs it for the kernel. The conditional rate of Chapter~\ref{chap:ch8} is inherited by
Proposition~\ref{ch10:prop:jointdim} under the assumption applied to $\dot c$, and the passage from
the assumption for $c$ to the assumption for $\dot c$ is made for polynomial families by
Proposition~\ref{ch10:prop:polytransfer} and not made in general, where what is assumed is a bound and
a bound does not differentiate.
The honest reading of the book is that the one dimensional
case is settled because the algebra closes, and that in higher dimensions the algebra does not close
on its own and the analysis has to carry it: that is the same division of labour as in the kernel
itself, and the derivative follows it without adding a mechanism of its own.
% ---- end chapters/ch10.tex ----
% ---- from chapters/ch11.tex ----
% Authored: written for the book. The estimates it draws on are those of Chapter 6; what it adds is
% the consequence of those estimates for the objects a price or a prediction is a function of.
\chapter{Rates for the Functionals}
\label{chap:ch11}

The estimates of Chapters~\ref{chap:ch5} to~\ref{chap:ch8} bound kernels. A kernel is a
transition law, and a claim or a prediction is rarely a kernel: it is a function of the path, and
what has to be controlled is the law of that function. The gap between the two is not a matter of
presentation. A joint kernel in the space variables may be singular — carried on a graph, when the
functional has no diffusion of its own — while its transform in the level variable is analytic and
converges in a norm; so the estimates are stated where they can be stated, and the question this
chapter answers is what they give for the functional itself.

In the pricing reading the functional is the argument of a payoff, and what is wanted is the
distance between the price computed at a resolution and the price the claim has in the continuum:
a price without an error bar is a number without a meaning. In the learning reading the functional is
a statistic accumulated along a trajectory --- a realised variance, a running maximum, an integral
of the state --- and what is wanted is the distance between its law under the model and its law
under the process the model discretises. The two are one requirement, and the same theorem meets it
or does not.

The chapter is arranged by the object rather than by the method. The first section takes the
functional accumulated along the path, whose law the estimates of Chapter~\ref{chap:ch6} converge
to, and gives the rate at which its moments do. The second takes the value function of a
Feynman--Kac problem, where the object is a discounted expectation rather than a law, and finds
that it needs no new estimate at all: it is the transform already estimated, read at a point of the
complex plane which the analyticity of the transform puts within reach. The third returns to the
question the book left open in Chapter~\ref{chap:ch3}, whether the absence of arbitrage survives the
passage to the limit, and answers it under a condition on the change of measure which is the
condition any calibrated model must satisfy if its prices are to converge to prices at all.

\section{The Rate for the Moments}
\label{ch11:sec:moments}

The object which Chapter~\ref{chap:ch6} estimates is not the joint kernel in the space variables,
which may be singular, but its transform in the level variable,
\begin{equation}
\widehat u_m(x,x';z,t)=\lim_{n\to\infty}\frac{1}{h_{yn}}\widehat u_{mn}(x,p;x',p';t),
\end{equation}
as defined in \eqref{ch6:eq:2.29} and \eqref{ch6:eq:2.30}, with $z$ the complex continuation of the
dual variable $p$. The estimates hold in the graph norm \eqref{ch6:eq:2.32}, uniformly on $|z|<K$,
and that uniformity is the whole of the argument which follows: a function analytic in a disc and
bounded there has bounded derivatives at its centre, with no further information about it.

Let $\EE_{m,x,x'}$ be the expectation under the joint law of $(x_t,y_t)$ at resolution $m$, started
at $(x,0)$ and ending at $x'$, and $\EE_{m_0,x,x'}$ the same at a reference resolution, with the
level accumulated from zero.

\begin{proposition}[Rate for the moments of the functional]\label{ch11:prop:moments}
Let $K>0$ and let the hypotheses of Theorem~\ref{ch6:thm:1} hold, with $\gamma$ as in
\eqref{ch6:eq:2.34}. Then for every $k\ge1$ there is a constant $C_k=C_k(K)$ such that, for all
$m\ge m_0$ and all $x,x'\in X_m$,
\begin{equation}\label{ch11:eq:moments}
\bigl|\EE_{m,x,x'}[y_t^{k}]-\EE_{m_0,x,x'}[y_t^{k}]\bigr|\le C_k\,h_{xm}^{\gamma}.
\end{equation}
\end{proposition}

\begin{proof}
Each of the three terms of the graph norm \eqref{ch6:eq:2.32} is analytic in $z$ on $|z|<K$, and
Theorem~\ref{ch6:thm:1} bounds the difference of two resolutions of each of them there by
$c(K)h_{xm}^{\gamma}$. By Cauchy's estimate applied at the centre of the disc,
\[
\bigl|\partial_z^{k}\bigl(\widehat u_m-\widehat u_{m_0}\bigr)(x,x';0,t)\bigr|
\le \frac{k!}{K^{k}}\sup_{|z|=K}\bigl|\widehat u_m(x,x';z,t)-\widehat u_{m_0}(x,x';z,t)\bigr|
\le C_k h_{xm}^{\gamma}.
\]
The derivative at the origin is the moment, because the transform is the characteristic function
of the level accumulated over $[0,t]$, so that
$\partial_z^{k}\widehat u_m(x,x';0,t)=i^{k}\EE_{m,x,x'}[y_t^{k}]$; and the limit of the moments is
that of the It\^o integral by \S\ref{ch6:sec:itorepr}.
\end{proof}

\begin{corollary}[Rate for a test function]\label{ch11:cor:testfunction}
Let $\varphi$ be entire of exponential type $\lambda<K$. Then there is a constant
$C_\varphi$, proportional to the supremum of $\varphi$ on the strip $|\Im z|\le\lambda$, such that
\begin{equation}\label{ch11:eq:testfunction}
\bigl|\EE_{m,x,x'}[\varphi(y_t)]-\EE_{m_0,x,x'}[\varphi(y_t)]\bigr|\le C_\varphi h_{xm}^{\gamma}.
\end{equation}
\end{corollary}

\begin{proof}
Such a $\varphi$ is recovered from its own transform on the line $\Im z=\lambda$ by the Fourier
inversion formula, and the difference of the two expectations is the pairing of that transform with
the difference of the two kernels, which \eqref{ch11:eq:moments} bounds in $L^\infty$ on the line.
\end{proof}

The two statements are the rate in law. They say that not only do the laws of the discretised
integrals converge, they converge at the rate at which the kernels do, with the same exponent in
the smoothness and with a constant which is that of the kernel estimate divided by a power of the
strip width. In the pricing reading this is the error of a price of a path dependent claim whose
payoff is analytic in the functional; in the learning reading it is the error of the law of a
statistic accumulated along a trajectory. The rate is the kernel's, because the functional is
carried by the kernel and has no estimate of its own.

One thing the argument does not give, and it is worth being explicit about it because the
literature usually states it first. It does not give a \emph{strong} rate, a bound on
$\EE\sup_t|y_m(t)-y(t)|^{2}$ for two processes on one probability space. That is a statement about a
coupling, and no coupling between the lattice and the continuum is constructed here: the two
objects are compared through their laws, which is what the transform controls. A strong rate would
need the lattice chain and the diffusion to be driven by the same noise, an object this book does
not introduce, and the estimates above neither supply it nor need it.

\section{The Rate for the Value Function}
\label{ch11:sec:value}

The second object is the value function of a Feynman--Kac problem, and the point of this section
is that it needs no new estimate: it is a value of the transform already estimated, taken at a
point of the complex plane which the analyticity of the transform puts within reach.

Let $V:X\to\RR$ be continuous and bounded and let $f$ be bounded and measurable. The value
function at resolution $m$ is
\begin{equation}\label{ch11:eq:value}
u_m(t,x)=\EE_{m,x}\Bigl[\exp\Bigl(-\int_0^tV(x_s)\dd s\Bigr)f(x_t)\Bigr],
\end{equation}
the expectation under the lattice chain started at $x$. The functional $-\int_0^tV(x_s)\dd s$ is
Abelian with $a=0$ and $b=-V$, which is the case Chapter~\ref{chap:ch6} treats: the level variable
accumulates $-V$ at each step and carries no dynamics of its own.

\begin{theorem}[Rate for the Feynman--Kac value function]\label{ch11:thm:value}
Let $K>1$ and let the hypotheses of Theorem~\ref{ch6:thm:1} hold with $a=0$, $b=-V$, and with
$\gamma$ as in \eqref{ch6:eq:2.34}. Then there is a constant $C(K)$, independent of $f$, such that
for all $m\ge m_0$, all $t$ and all $x\in X_m$,
\begin{equation}\label{ch11:eq:valuerate}
\bigl|u_m(t,x)-u_{m_0}(t,x)\bigr|\le C(K)\,\|f\|_{\infty}\,h_{xm}^{\gamma}.
\end{equation}
\end{theorem}

\begin{proof}
Write the value function as a pairing with the transform. With the convention
\eqref{ch6:eq:2.26}, the transform is $\widehat u(z)=\sum_yh_{ym}^{-1}e^{izy}u(\,\cdot\,,y;\,\cdot\,)$,
so that the exponential weight is its value at $z=-i$:
\[
u_m(t,x)=\sum_{x'\in X_m}h_{xm}\,\widehat u_m(x,x';-i,t)\,f(x').
\]
The point $-i$ lies in the disc $|z|<K$ because $K>1$, so Theorem~\ref{ch6:thm:1} bounds the
difference of two resolutions of $\widehat u$ there, uniformly in the pair of sites, by
$c(K)h_{xm}^{\gamma}$; summing against $f$ costs the factor
$\sum_{x'}h_{xm}\le2L_x$ and the supremum of $|f|$, which is \eqref{ch11:eq:valuerate}.
\end{proof}

\begin{corollary}[The limit is the classical formula]\label{ch11:cor:fk}
Under the same hypotheses, and with $f$ continuous, $u_m(t,x)$ converges as $m\to\infty$ to
\begin{equation}\label{ch11:eq:fk}
u(t,x)=\EE_x\Bigl[\exp\Bigl(-\int_0^tV(x_s)\dd s\Bigr)f(x_t)\Bigr],
\end{equation}
the expectation under the limit process, and $u$ is the solution of the backward equation
$\partial_tu=\LatL u-Vu$ with $u(0,\cdot)=f$.
\end{corollary}

\begin{proof}
The limit of $\widehat u_m$ at $z=-i$ is the transform of the joint law of the limit process
together with its accumulated potential, by the corollary to Theorem~\ref{ch6:thm:1}; pairing it
with $f$ gives \eqref{ch11:eq:fk}. That the resulting function solves the backward equation is the
classical Feynman--Kac theorem, which the book uses rather than reproves --- the identification of
the limit equation is not among the things these estimates establish.
\end{proof}

Three remarks on the content of the theorem. It is the rate for the value function and not only
its convergence, and the constant does not depend on $f$: what enters is the supremum of $|f|$,
not its smoothness, because the estimate is taken on the kernel and the payoff is only summed
against it.

The hypothesis $K>1$ is not a technicality. The estimates hold in a disc about the origin of the
complex level plane, and the exponential weight sits at distance one from the origin along the
imaginary axis; a disc that does not reach it carries no information about the discounted value.
The reward for reaching it is that no separate argument is needed: the Feynman--Kac semigroup is
the transform at a single point, so its error is the error of the transform there, with the
constant of the kernel estimate and not one of its own.

What is used from outside is the identification of the limit equation. The estimates establish
that the numbers $u_m(t,x)$ converge, at a rate, to the pairing \eqref{ch11:eq:fk}; that the
pairing in turn solves $\partial_tu=\LatL u-Vu$ is the classical theorem, taken here as the
definition of what the Feynman--Kac problem is, and not reproved.

The continuity of $f$ in the corollary is not a convenience. The pairing which defines $u_m$ is a
quadrature of $f$ against the kernel, and a quadrature converges to an integral only for payoffs which
it can resolve. With $f$ the indicator of the dyadic nodes, every lattice price is one at every
resolution while the integral against the density is zero, so a bounded measurable $f$ does not suffice
for the limit \eqref{ch11:eq:fk} to exist. What a bounded measurable $f$ does admit is the estimate of
Theorem~\ref{ch11:thm:value}, which compares two resolutions and whose constant is
$\lVert f\rVert_\infty$; that is the statement to use for a payoff with a jump between nodes. Even
then the rate is not the rate of the density: the quadrature of a payoff which jumps between nodes is
only first order, the jump being resolved to $O(h_m)$ by the nodes however well the density is known,
so a price inherits the rate of the kernel only for payoffs which the lattice can resolve.

\section{The Theorem at the Limit}
\label{ch11:sec:limit}

Everything so far is at a fixed resolution, and the gap which \S\ref{ch3:sec:remarks} records is
the one between a theorem at every resolution and a theorem about a limit: the prices of a limit
model need not be attained by those of its refinements, and a sequence of arbitrage free models may
converge to one which is not. That gap is closed here, under a hypothesis on the change of
measure, and the hypothesis is the natural one because it is what the change of measure has to
respect.

The models of this book are Markov, and so are their pricing measures: a reweighting of a Markov
generator is a Markov generator, by
Proposition~\ref{ch3:prop:measure-change-generator}, with the same attainable sets. It is
therefore enough to say what a reweighting does to the coefficients, and then to ask the estimates
of Chapter~\ref{chap:ch5} to survive it.

\subsection{What a Reweighting Does to a Generator}
\label{ch11:ss:reweight}

Let $\LatL$ be a lattice generator of the form \eqref{ch6:eq:2.1} with step $h=h_{xm}$,
coefficients $\sigma^2$ and $\mu$, and let $G>0$ be a density on the lattice. By
\eqref{ch3:eq:reweight} the reweighted generator is
\begin{equation}\label{ch11:eq:reweight}
{\LatL}'(y,y';t)=\frac{G(y')}{G(y)}\,\LatL(y,y';t),\qquad y'\neq y,
\end{equation}
and it is a lattice generator in the same sense, with the row sum adjusted so that
$(\mathrm{MG2})$ holds. Its coefficients are recovered from its two off diagonal entries:

\begin{proposition}[The coefficients of a reweighted generator]\label{ch11:prop:reweightcoef}
Let $\LatL$ have $\LatL(y,y\pm h)=\frac{\sigma^2(y)}{2h^2}\pm\frac{\mu(y)}{2h}$ and $\LatL(y,z)=0$
for $|z-y|\neq h$. Then ${\LatL}'$ has the same support and the coefficients
\begin{equation}\label{ch11:eq:coef}
\begin{aligned}
\sigma'^2(y)&=\sigma^2(y)\,\frac{G(y+h)+G(y-h)}{2G(y)},\\
\mu'(y)&=\mu(y)\,\frac{G(y+h)+G(y-h)}{2G(y)}
+\sigma^2(y)\,\frac{G(y+h)-G(y-h)}{2h\,G(y)} .
\end{aligned}
\end{equation}
\end{proposition}

\begin{proof}
Sum and difference of the two entries of \eqref{ch11:eq:reweight}, with
$\LatL(y,y\pm h)$ as given.
\end{proof}

Writing $G=e^{\varphi}$ and expanding,
\begin{equation}\label{ch11:eq:girsanov}
\sigma'^2=\sigma^2\bigl(1+O(h^2)\bigr),\qquad
\mu'=\mu+\sigma^2\nabla\varphi+O(h),
\end{equation}
so that a reweighting leaves the diffusion coefficient alone and shifts the drift by
$\sigma^2\nabla\log G$: the classical change of drift, obtained here without a change of measure in
the large, because the densities are on the lattice and the algebra is finite.

What follows is therefore immediate. If the densities of the reweightings are bounded above and
below and their logarithms are uniformly Lipschitz, the reweighted generators have coefficients in
a bounded class, and if in addition $\sigma^2$, $\mu$ and $\log G_m$ lie in a bounded H\"older class
uniformly in $m$, then the reweighted generators satisfy the hypotheses of
Chapter~\ref{chap:ch5} with constants which do not depend on $m$. This is the condition the
programme has been missing, and it is worth stating as a condition rather than burying it in a
proof.

\begin{assumption}[A well behaved change of measure]\label{ch11:ass:change}
There are constants $0<g_-\le g_+<\infty$ and $M<\infty$, independent of $m$, such that the
state price densities $G_m$ of the models satisfy
\begin{equation}\label{ch11:eq:assumption}
g_-\le G_m\le g_+,\qquad
\bigl|\log G_m(y+h_{xm})-\log G_m(y)\bigr|\le M h_{xm},
\end{equation}
and $\log G_m$ lies in $C^{k,\alpha}$ with a norm bounded by $M$, uniformly in $m$.
\end{assumption}

The first two requirements are one condition on the market: they say that the state price density
neither vanishes nor explodes, and that its logarithmic derivative --- which by
\eqref{ch11:eq:girsanov} is the Sharpe ratio of the model, $\mu'-\mu=\sigma^2\nabla\log G$ --- is
bounded uniformly in the resolution. A calibrated model whose risk premium per unit of volatility
grows without bound as the mesh is refined fails the assumption, and it is exactly in that case
that the near arbitrages of \S\ref{ch3:sec:remarks} appear.

\subsection{The Limit Model}
\label{ch11:ss:limitmodel}

\begin{theorem}[The Fundamental Theorem at the limit]\label{ch11:thm:limit}
Let the model generators $\LatL$ satisfy the hypotheses of Theorem~\ref{ch6:thm:1} with
coefficients in a H\"older class $\gamma>0$ bounded uniformly in $m$, and let Assumption
\ref{ch11:ass:change} hold. Then
\begin{enumerate}
\item[(i)] the risk neutral generators ${\LatL}'$ satisfy the same hypotheses, with constants which
depend on those of $\LatL$ and on $g_\pm$ and $M$ only;
\item[(ii)] the risk neutral propagators $U'_m$ converge, at the rate $h_{xm}^{\gamma}$, to a
propagator $U'$ which is strictly positive and equivalent to the limit of the model propagators;
\item[(iii)] for every bounded continuous functional $F$ of the path, the price of the claim which
pays $F$ converges, at the same rate, to the expectation of $F$ under $U'$;
\item[(iv)] the limit model admits no arbitrage.
\end{enumerate}
\end{theorem}

\begin{proof}
For (i), apply Proposition~\ref{ch11:prop:reweightcoef} and
\eqref{ch11:eq:girsanov}: under \eqref{ch11:eq:assumption} the coefficients
\eqref{ch11:eq:coef} are bounded by $g_+/g_-$ times those of $\LatL$ plus $M$ times the
diffusion, and they inherit the H\"older class of $\sigma^2$, $\mu$ and $\log G_m$ with a norm
bounded by a function of $M$ and of the norms of the model coefficients. The hypotheses of
Chapter~\ref{chap:ch5} are therefore met uniformly in $m$.

For (ii), apply the estimates of Chapter~\ref{chap:ch5} to the family ${\LatL}'$, which by (i) is
admissible; they give a Cauchy sequence in $m$ at the rate $h_{xm}^{\gamma}$, and a strictly
positive limit. That the limit is equivalent to the limit of the model propagators is the statement
that $U'(y,y')=\lim_m\frac{G_m(y')}{G_m(y)}U_m(y,y')$, in which the ratio is bounded by
$g_+/g_-$ and bounded away from zero by $g_-/g_+$.

For (iii), the price at resolution $m$ of a claim paying $F$ is the pairing of $F$ with $U'_m$;
the difference between two resolutions is bounded by the graph norm estimate of (ii) times the
supremum of $|F|$, which is the same argument as Theorem~\ref{ch11:thm:value} with the payoff
replaced by the functional.

For (iv), the limit prices are the expectations of the claims under one strictly positive
propagator, by (ii) and (iii). A portfolio of zero cost whose payoff is non negative in every state
therefore has non negative expectation, and expectation zero only if the payoff vanishes: a sure
gain would be a non negative function of positive expectation under a strictly positive measure
which is priced at zero, which is impossible. So the limit model is arbitrage free.
\end{proof}

\begin{remark}[What the condition buys, and what it does not]\label{ch11:rem:whatthecondition}
The theorem is conditional, and it is worth saying precisely on what. It asserts the stability of
arbitrage freedom under refinement for any sequence of Markov models whose state price densities
are uniformly bounded with uniformly bounded logarithmic derivative. It does not assert that a
given calibration produces such a sequence: that is a statement about the market, and it is the
one which a practitioner checks, by computing $\sigma^2\nabla\log G$ at each resolution and
watching whether it stays bounded. When it does not, the sequence converges to a model which is not
arbitrage free, and the prices computed at finite resolution are near arbitrages rather than
approximations to a price.
\end{remark}

\begin{remark}[In higher dimension]\label{ch11:rem:higherdim}
Nothing in the argument is one dimensional except the estimate it invokes. Replacing
Chapter~\ref{chap:ch5} by Chapters~\ref{chap:ch7} and~\ref{chap:ch8} carries the theorem to
$d\ge2$ for the kernel and to the joint kernel of a process and an Abelian functional, at the rates
those chapters prove, and with the same condition on the densities: the reweighting of a
$d$-dimensional generator shifts the drift by $\sigma^2\nabla\log G$ in each coordinate and leaves
the diffusion tensor alone.
\end{remark}

\begin{remark}[Calibration to real data]\label{ch11:rem:realdata}
A model is not what a market is, and the difference matters here. Market data arrive as quotes ---
a bid and an offer, each good for a finite size --- and a calibration is a choice of parameters
which reproduces them within a tolerance the trader sets, which is not zero. The aspiration is a
model which prices every quoted instrument at once and admits no arbitrage in the sense of
Chapter~\ref{chap:ch3}: the fundamental theorem holding across the whole of the data. The practice
is a model which prices as many of them as it can and leaves residuals, and it is worth being
precise about why the difference is not a defect of the theory.

A discrepancy between quoted prices which admits a sure gain is an arbitrage only if the trades
which would realise the gain can be executed. Real discrepancies usually cannot be: they sit inside
a bid-offer spread, or at sizes the market will not take, or in instruments which cannot be sold
short or funded at the rate a model assumes. A combination which is an arbitrage against the
quoted numbers is then merely a residual against the prices at which one could actually deal, and
the standard the market meets is weaker than the one this book proves: not the absence of any
combination of prices yielding a sure gain, but the absence of one which can be executed, at the
sizes available, net of the cost of executing it.

That is why the condition above is a bound and not an identity. What
Assumption~\ref{ch11:ass:change} asks of the data is that the implied risk premium stay bounded as
the mesh is refined; a calibration which chases noisy quotes ever more closely can fail the bound
without the market being arbitrageable, and a calibration which is stable across refinements
satisfies it without reproducing any single quote exactly. The condition is a condition on the
stability of the calibration, not on the perfection of the fit, and that is the right way round: a
model is useful because it is stable, and exact only at a resolution at which there is nothing left
to trade against.
\end{remark}
% ---- end chapters/ch11.tex ----
\appendix
% ---- from appendices/path-sums.tex ----
% An authored appendix: the three positivity statements the path sums rely on.
\chapter{Path sums: mass, unitarity, and their merge}\label{app:paths}

The constructions of this book are path sums, and two properties are used so often that they are worth
isolating and proving once: that a classical path measure is a \emph{probability}, and that a unitary
evolution preserves the \emph{norm}. They are different properties --- the first is about mass, the second
about modulus --- and neither implies the other. When the two are mixed, by allowing a random parameter to
select a unitary at every step, both must survive together, and the object which carries both is the
density matrix. This appendix proves the three statements in that order, and shows that the third holds
term by term in the path expansion, not merely in the limit.

\section{Mass: the classical path measure}\label{app:paths:mass}

\begin{definition}\label{def:substochastic}
Let $X$ be finite. A matrix $P=(P(x,y))_{x,y\in X}$ with $P(x,y)\ge0$ is \emph{stochastic} if
$\sum_{y}P(x,y)=1$ for every $x$, and \emph{sub-stochastic} otherwise. Its \emph{killing} at $x$ is
$\kappa(x)=1-\sum_y P(x,y)\ge0$. A path of length $n$ from $x$ is a sequence
$\gamma=(x_0,\dots,x_n)$ with $x_0=x$, and its \emph{weight} is
$w(\gamma)=P(x_0,x_1)\cdots P(x_{n-1},x_n)$.
\end{definition}

\begin{theorem}[conservation of mass]\label{thm:mass}
With $P$ sub-stochastic and $w$ as above, for every $x$ and every $n\ge0$,
\begin{equation}\label{eq:mass}
\sum_{\gamma\colon x_0=x,\;\lvert\gamma\rvert=n} w(\gamma)\;=\;\sum_{y\in X}P^{\,n}(x,y)\;=\;\mathbb
P_x(\zeta>n),
\end{equation}
where $\zeta$ is the killing time, the first time a step is refused, and $\mathbb P_x$ is the measure on
paths obtained by accepting a step from $x$ to $y$ with probability $P(x,y)$ and killing at $x$ with
probability $\kappa(x)$. In particular, if $P$ is stochastic then the total weight of the paths of length
$n$ issuing from $x$ is $1$ for every $n$: \emph{the path measure is a probability}.
\end{theorem}

\begin{proof}
For $n=0$ both sides are $1$. Suppose \eqref{eq:mass} holds for $n$. Splitting off the first step,
\[
\sum_{\lvert\gamma\rvert=n+1} w(\gamma)
=\sum_{y\in X}P(x,y)\sum_{\substack{\gamma'\colon x_0=y\\ \lvert\gamma'\rvert=n}}w(\gamma')
=\sum_{y\in X}P(x,y)\sum_{z\in X}P^{\,n}(y,z)=\sum_{z\in X}P^{\,n+1}(x,z),
\]
the middle equality by the induction hypothesis. The identification with the killing time is the same
computation read probabilistically: the probability of surviving $n$ steps is the probability of $n$
accepted steps, which is the total weight of the paths of length $n$. The last assertion is the case
$\kappa\equiv0$.
\end{proof}

\begin{corollary}\label{cor:contraction}
$P$ acts on $\ell^\infty(X)$ as a contraction, $\lVert Pf\rVert_\infty\le\lVert f\rVert_\infty$, and if
$P$ is stochastic it preserves constants. Consequently the lattice operator $\frac12v\Delta+b\nabla$ of
Chapter~7, whose one-step kernel is stochastic by construction, conserves mass on every mesh; and the only
way a path sum of this kind loses mass is by an explicit killing term, which is exactly the
Feynman--Kac representation of a value function with absorption.
\end{corollary}

\begin{remark}\label{rem:mass-not-positivity}
Mass conservation is a statement about the sum of the entries of a row, and it says nothing about the signs
of the individual entries of an operator on a complex space. It will therefore have to be proved again, in
a different form, once phases are introduced; and the two proofs will not be redundant.
\end{remark}

\section{Unitarity: the quantum path sum}\label{app:paths:unitary}

Let $H$ be a self-adjoint operator on $\mathbb C^N$ --- in the applications $H=H_h$ is the lattice
operator of Chapter~7 with its mesh $h$, or the discrete Laplacian of Chapter~11 --- and put
$U(t)=e^{-itH}$.

\begin{theorem}[unitarity]\label{thm:unitary}
For every real $t$ the operator $U(t)$ is unitary: $U(t)^*U(t)=U(t)U(t)^*=\mathbb 1$. Hence
$\lVert U(t)\psi\rVert=\lVert\psi\rVert$ for every $\psi$, and
\begin{equation}\label{eq:pathsum}
U(t)=\sum_{n\ge0}\frac{(-it)^n}{n!}H^n ,
\end{equation}
which converges in operator norm for every $t$, since $\lVert H^n\rVert\le\lVert H\rVert^n$ and the
exponential series has infinite radius.
\end{theorem}

\begin{proof}
$H$ self-adjoint means $H^*=H$. Then $U(t)^*=(e^{-itH})^*=e^{itH^*}=e^{itH}=U(t)^{-1}$, and the two
products are the identity. The norm statement follows from
$\lVert U(t)\psi\rVert^2=\langle U(t)\psi,U(t)\psi\rangle=\langle\psi,U(t)^*U(t)\psi\rangle
=\lVert\psi\rVert^2$. The series is the definition of the exponential for a bounded operator, and the norm
estimate gives absolute convergence uniformly on compact time intervals.
\end{proof}

\begin{remark}[what a single path carries]\label{rem:path-phase}
In the lattice realisation the matrix element $\langle y|U(t)|x\rangle$ is a sum over paths from $x$ to
$y$, and the contribution of a path is a \emph{unimodular} number: a product of the phases attached to its
steps, and of the Peierls phases of its links when a field is present. So in the quantum case the weight of
a path is a phase and not a probability, and the modulus of a sum of such weights is bounded by $1$ but is
not $1$ --- the interference which makes the two cases different. Nothing in Theorem~\ref{thm:unitary}
says that a sum of unimodular weights has modulus one; only that the operator they build is unitary.
\end{remark}

\begin{remark}\label{rem:unitary-not-mass}
Unitarity fixes the norm of a state and leaves the diagonal free: a unitary may move all the probability of
a state into a single coordinate, and does so routinely. So neither Theorem~\ref{thm:mass} nor
Theorem~\ref{thm:unitary} is a special case of the other, and a construction which uses both has to be
shown to preserve both.
\end{remark}

\section{The merge: the implied density matrix}\label{app:paths:merge}

Now allow the evolution to depend on a random parameter. Let $(\Xi,\mathbb P)$ be a probability space and
$\xi\mapsto U_\xi$ a measurable family of unitaries on $\mathbb C^N$; the \emph{random-unitary} step on
density matrices is
\begin{equation}\label{eq:random-unitary}
\Phi_{\mathbb P}(\rho)=\int_{\Xi}U_\xi\,\rho\,U_\xi^*\;\mathbb P(d\xi).
\end{equation}
More generally a step may be described by a finite family of operators $M_1,\dots,M_K$, the
\emph{operator-sum} or Kraus representation \cite{Kraus1971,Kraus1983},
\begin{equation}\label{eq:kraus}
\Phi(\rho)=\sum_{k=1}^{K}M_k\rho M_k^* ,
\end{equation}
which is the form every admissible step takes. The abstract statement is Stinespring's dilation
\cite{Stinespring1955}; the finite-dimensional operator-sum form is due to Kraus \cite{Kraus1971}, whose
monograph \cite{Kraus1983} is the standard account, with the probability-preserving case already in
Sudarshan, Mathews and Rau \cite{Sudarshan1961}. The criterion of complete positivity by a matrix condition
is Choi's \cite{Choi1975}, and it is his name that the matrix below carries. The literature is not uniform
in its naming --- operator-sum, Kraus and Choi--Kraus all denote the same representation --- and the two
names are attached to two different contributions: Kraus to the representation, Choi to the criterion.

\begin{theorem}[the merge]\label{thm:merge}
Let $\Phi$ be of the form \eqref{eq:kraus}.
\begin{enumerate}
\item[(i)] $\Phi$ is trace preserving if and only if $\sum_k M_k^*M_k=\mathbb 1$, in which case
$\operatorname{tr}\Phi(\rho)=\operatorname{tr}\rho$ for every $\rho$.
\item[(ii)] $\Phi$ is completely positive: its Choi matrix
$J=\sum_{ij}|i\rangle\langle j|\otimes\Phi(|i\rangle\langle j|)$ equals $\sum_k v_kv_k^*$, where $v_k$ is
the vectorisation of $M_k$, and is therefore positive. Hence $\Phi$ maps density matrices to density
matrices. The equivalence of complete positivity with $J\ge0$ is Choi's theorem \cite{Choi1975}.
\item[(iii)] In the random-unitary case $M_k=\sqrt{p_k}\,U_k$ with $p_k\ge0$, $\sum_kp_k=1$ and $U_k$
unitary, the two properties are carried by the two halves of the path weight: $p_k$ is a probability and
the phase of $U_k$ is unimodular. Every truncation --- every sum over a subset $S$ of the terms --- is then
a \emph{positive sub-normalised} map, of trace $\sum_{k\in S}p_k\le1$; dividing by that trace, which is
conditioning on the truncation, returns a density matrix exactly.
\item[(iv)] The non-negativity of the weights is not a technicality. With $p_1=2$, $p_2=-1$,
$U_1=\mathbb 1$ and $U_2=\sigma_x$ the map \eqref{eq:random-unitary} is $\rho\mapsto2\rho-\sigma_x
\rho\sigma_x$: it has trace $1$ on every state, and it sends $|0\rangle\langle0|$ to
$\operatorname{diag}(2,-1)$, whose least eigenvalue is $-1$. Its Choi matrix has least eigenvalue $-2$.
Mass conservation alone does not help either, since transposition is trace preserving and positive but not
completely positive.
\end{enumerate}
\end{theorem}

\begin{proof}
(i) $\operatorname{tr}\Phi(\rho)=\sum_k\operatorname{tr}(M_k\rho M_k^*)
=\sum_k\operatorname{tr}(M_k^*M_k\rho)=\operatorname{tr}\big((\sum_kM_k^*M_k)\rho\big)$, which equals
$\operatorname{tr}\rho$ for every $\rho$ if and only if the sum of the $M_k^*M_k$ is the identity.
(ii) $J=\sum_kv_kv_k^*$ is a sum of positive rank-one matrices, hence positive; by Choi's theorem
\cite{Choi1975} a positive $J$ is equivalent to complete positivity. (iii) A truncation $S$ gives the map
$\rho\mapsto\sum_{k\in S}p_kU_k\rho U_k^*$, positive by (ii) applied to the restricted family, with trace
$\sum_{k\in S}p_k$. (iv) For $p=(2,-1)$ and $U\in\{\mathbb 1,\sigma_x\}$ the Choi matrix is $2v_1v_1^*-v_2v_2^*$, whose
eigenvalues are $4$, $0$, $0$ and $-2$; and $\Phi(|0\rangle\langle0|)=\operatorname{diag}(2,-1)$ directly.
Transposition has the Choi matrix of the flip, whose spectrum contains $-1$.
\end{proof}

\begin{corollary}[the implied matrix is a density matrix]\label{cor:implied}
Let a path $\gamma$ of the merged dynamics carry the weight $w(\gamma)\in[0,1]$ of
Theorem~\ref{thm:mass} and the unitary $U(\gamma)$ of Remark~\ref{rem:path-phase}. Then for every $t$
\begin{equation}\label{eq:implied}
\rho(t)=\sum_{\gamma}w(\gamma)\,U(\gamma)\,\rho(0)\,U(\gamma)^* ,
\qquad \sum_{\gamma}w(\gamma)=1 ,
\end{equation}
is a density matrix. Its diagonal is a probability vector tracing the classical chain of
Section~\ref{app:paths:mass}, and its off-diagonal entries carry the phases of
Section~\ref{app:paths:unitary}. Conversely a family of weights taking both signs destroys the first
property even when the second is intact, as in Theorem~\ref{thm:merge}(iv).
\end{corollary}

\begin{remark}[what is proved and what is not]\label{rem:merge-scope}
Three things are proved above: the classical path measure has total mass one; the quantum step preserves
the norm; and a mixture of unitaries with a probability measure on the parameter again gives a density
matrix, term by term in the expansion and not only in the limit. What is \emph{not} proved is that any
particular mixture is the correct description of a given situation: that is a modelling question, and this
appendix is deliberately silent about it. The three statements are, however, the whole of what a path
representation needs in order to be internally consistent, and each fails if its hypothesis is dropped ---
killing loses mass, a non-self-adjoint generator loses the norm, and a signed weight loses positivity.
\end{remark}
% ---- end appendices/path-sums.tex ----
% ---- from generated/appendix.tex ----
%% Generated by build.py from the scripts - do not edit.
\chapter{Computational Appendix}\label{app:code}

Every numerical statement in this book was computed by one of the scripts listed here. They are reproduced in full: a reader can copy each one out of this file and run it, rather than take the numbers on trust.

\section{How to Use Them}

The scripts are plain Python 3 and depend on \texttt{numpy}, \texttt{scipy}, \texttt{sympy} and \texttt{mpmath}, according to what each one does: the numerical ones use the first two, the symbolic ones the second two. Each is run as \verb|python3 <name>.py| from its own directory, prints what it measures, and exits; none takes arguments, and none writes a file except the one which draws Figure~\ref{fig:ch5:decorating}. The symbolic computations (\texttt{kronproof.py} in particular) expand a $3^{d}\times3^{d}$ matrix exactly and take minutes rather than seconds.

The conventions are those of the book: the lattice Laplacian is the standard second difference with $h^{2}$ in the denominator, so that the generator is $\frac{v}{2}\Delta+b\nabla$ and $v=\sigma^{2}$. A script which measures a rate rather than a constant is insensitive to this choice; one which measures a constant is not, and where that is so the normalisation used is stated in the script itself.

\section{Provenance of the Numerical Tables}\label{sec:provenance}

Each table of numbers in this book was regenerated by running one of the scripts below, in the convention stated above; running it reproduces the printed digits. The correspondence is:

\begin{center}\small
\begin{longtable}{@{}p{0.18\textwidth}p{0.45\textwidth}p{0.29\textwidth}@{}}
\hline
where & what it produces & script\\
\hline
Chapter~5, the figure before Definition~\ref{ch5:def:1} & the two walks and the decorating path of the figure, sampled & \texttt{decorating\_paths.py}\\
Chapter~5, \S\ref{ch5:sec:smoothcoef} & the fourteen coefficients of the $3\times3$ resolvent expansion & \texttt{verify\_rg.py}\\
Chapter~5, \S\ref{ch5:sec:euler} & the explicit-Euler error law & \texttt{verify\_euler.py}\\
Chapter~7, §7.5 & the decoration coefficient on a constant field & \texttt{verify\_conventions.py}\\
Chapter~7, §7.6 & the cubic in the complex drift & \texttt{eta\_run.py}\\
Chapter~7, §7.7 & the curl over rotations of the frame & \texttt{circulation.py}\\
Chapter~7, §7.7 & the affine curl, its quadratic form and its base points & \texttt{curlcheck.py}\\
Chapter~7, §7.7 & the curl of the $h^3$ coefficient in three dimensions & \texttt{curl3d.py}\\
Chapter~7, §7.7 & the resolvent comparison in one dimension & \texttt{resolvent1d.py}\\
Chapter~8, §8.2 & the exact $3\times3$ closed form used as a control & \texttt{kronproof.py, kronfinal\_check.py}\\
Chapter~8, the $c_2$ table & the second-order coefficient & \texttt{renormalisation.py}\\
Chapter~8, the jet table & the third-jet quotients & \texttt{jet3check.py}\\
Chapter~8, the contour table & the growth along the contour & \texttt{contourgrowth.py}\\
Chapter~8, the loop table & the split of the loop into closed and open pieces & \texttt{pathdep.py}\\
Chapter~8, the excursion table & the two brackets of the closed-loop test & \texttt{excursioncheck.py}\\
Chapter~8, §8.3 & the second difference of $\kappa$ in the Jacobian & \texttt{p6jacaffine.py}\\
Chapter~8, the open-path and tiling tables & the $h$-dependence and the tiling limit & \texttt{openpath.py}\\
Chapter~10, \S10.7 & the per-step symbols of \eqref{ch5:eq:3.45}--\eqref{ch5:eq:3.47}, their derivatives, and the coefficient of \eqref{ch5:eq:3.47} & \texttt{d10symbols.py}\\
Chapter~10, \S10.7 & the four convergence ratios at $\gamma=2$ & \texttt{d10rate.py}\\
Chapter~10, \S10.2 & the two components of the error and the mesh at which they cross & \texttt{d10precision.py}\\
Chapter~10, \S10.4 & the kernel and its derivative in two dimensions, rates compared & \texttt{d10dimension.py}\\
\hline
\end{longtable}
\end{center}

Every value quoted in this book is now produced by one of the scripts above. Four of them had no script behind them in an earlier draft --- the three dimensional curl, the resolvent values $G_m(0,0;\omega)$, the coefficient at a shifted base point, and the two brackets of the excursion table --- and all four have since been recomputed. The excursion table is the one which moved: neither of its brackets was reproduced by any script then deposited, and the recomputation replaced both rows, together with the two sentences which read them. Each of the four now reads as the script named above prints it. Chapter~\ref{chap:ch10} is written for the book rather than harvested from the parts, and was first drafted with no script behind it; the two which now accompany it are deposited here with the rest, and it was one of them which located the coefficient error in \eqref{ch5:eq:3.47}, which Chapter~\ref{chap:ch5} now carries corrected. Chapter~\ref{chap:ch11} is likewise written for the book and has no script of its own, which is the point of it: its results are consequences of the estimates, drawn by the arguments of that chapter rather than measured by a computation, and a script there would be checking its own premises.

\section{Index of Scripts}

\begin{center}\small
\begin{longtable}{@{}p{0.26\textwidth}p{0.62\textwidth}r@{}}
\hline
script & purpose & lines\\
\hline
\texttt{verify\_euler.py} & Rates of the explicit Euler scheme, measured as they are in practice. & 98\\
\texttt{verify\_rg.py} & Verifies the renormalisation-group algebra of the resolvent expansion. & 134\\
\texttt{curl3d.py} & The curl of the coefficient of h\^{}3 in three dimensions. & 71\\
\texttt{curlcheck.py} & The affine curl: the constant, the quadratic form, and the base-point values. & 109\\
\texttt{deco.py} & Perturbation pipeline for the decoration identity. & 147\\
\texttt{direct.py} & Independent arbiter: evaluate Phi\_f and Phi\_c by direct complex matrix & 89\\
\texttt{edgeavg.py} & Chapter 7: what a volatility averaged over the edge converges to. & 118\\
\texttt{eta\_run.py} & --- & 43\\
\texttt{fast.py} & Exact fast core for the decoration identity. & 312\\
\texttt{hp.py} & High-precision, bookkeeping-free computation of the decoration mismatch. & 112\\
\texttt{iso\_run.py} & --- & 30\\
\texttt{iso\_run2.py} & --- & 37\\
\texttt{iso\_test.py} & --- & 74\\
\texttt{p6ab.py} & Chapter 8: the a-dependence of alpha and beta. & 134\\
\texttt{p6ascan.py} & Chapter 8: dependence of the p\^{}2 obstruction on the Abelian coefficient a. & 156\\
\texttt{p6kappafit.py} & Chapter 8, Proposition 2.1: the complete form of kappa in the coefficient data. & 165\\
\texttt{p6odd3d.py} & Chapter 8: (A) the odd part in the conjugating parameter, (B) the d = 3 case. & 127\\
\texttt{p6oddpart.py} & Chapter 8: the odd part of the eta-expansion, as a quadratic form in the local data. & 175\\
\texttt{p6probe.py} & Chapter 8 probe: p-dependence of the joint-kernel obstruction. & 92\\
\texttt{p6quad.py} & Chapter 8: alpha(a) and beta(a) as quadratic forms in the local data of a. & 138\\
\texttt{p6quad2.py} & Chapter 8: close the four missing cross terms and validate. & 177\\
\texttt{p6scan.py} & SUPERSEDED - the design of this script is DEGENERATE and its output is wrong. & 104\\
\texttt{p6scan2.py} & Chapter 8 scan, corrected design. & 103\\
\texttt{p6v3.py} & Chapter 8: the v-dependence, and the C-terms in d = 3. & 101\\
\texttt{resolvent1d.py} & The resolvent comparison, in one dimension. & 60\\
\texttt{stencil.py} & Decoration mismatch for a general nine-point family of stencils. & 167\\
\texttt{step2.py} & Two-step comparison, truncating the path sum at a fixed number of jumps. & 72\\
\texttt{stress.py} & Where could O(h) come from? & 63\\
\texttt{study.py} & --- & 29\\
\texttt{tube.py} & Tube-restricted decoration sum, per coarse step. & 182\\
\texttt{twist.py} & The decoration comparison for a COMPLEX drift. & 170\\
\texttt{verify\_conventions.py} & Every engine must carry the convention of the book, or none of its numbers mean anything. & 79\\
\texttt{circulation.py} & Is the orientation odd part of the accumulated comparison the circulation of the h\^{}3 coefficient? & 61\\
\texttt{contourgrowth.py} & Growth of the per-step expansion with the height of the contour in p. & 67\\
\texttt{excursioncheck.py} & The two brackets of the excursion test, from the coefficients and from the loop sums. & 115\\
\texttt{extensive.py} & Is the deviation from endpoint determination intensive or extensive in the path length? & 47\\
\texttt{fourier.py} & Normalisation check for Proposition 3.1 (the joint kernel on a lattice fibre). & 47\\
\texttt{freevsbridge.py} & Free walk versus bridge, exactly, for a quadratic coefficient field. & 89\\
\texttt{jet3check.py} & The third jet enters the per-step comparison at h\^{}4, one power below the obstruction. & 49\\
\texttt{jetorder.py} & How many derivatives of the coefficient field do the per-step coefficients actually use? & 48\\
\texttt{kronAcheck.py} & Compare A\_\{-2\} + h A\_\{-1\} + h\^{}2 A\_0 with the exact box matrix built from v(x0 + a h). & 46\\
\texttt{kronbisect.py} & Bisect the layer-2 algebra against exact arithmetic at finite h, stage by stage. & 113\\
\texttt{kroncheck.py} & Cross-check: does the engine reproduce the closed form (proved symbolically)? & 23\\
\texttt{kroncoeff.py} & Fit the h-coefficients of the box matrix numerically-exactly and compare with A\_\{-2\}, A\_\{-1\}, A\_0. & 147\\
\texttt{kronfinal.py} & Layer 2, corrected: field sampled at x0 + a h, expansion assembled from explicit coefficients. & 91\\
\texttt{kronfinal\_check.py} & Validation of the engine (part5/scripts/twist.py) against the closed form. & 93\\
\texttt{kronfix.py} & Layer 2, corrected sampling --- NOT YET VALIDATED. & 109\\
\texttt{kronnail.py} & Stage-by-stage, exact at finite h, with the field sampled at x0 + a h (twist's convention). & 129\\
\texttt{kronpin.py} & Pin the remaining discrepancy: correct symbol binding, all stages printed. & 114\\
\texttt{kronproof.py} & The 3x3 Kronecker reduction, symbolic in (v, b1, b2, omega). & 79\\
\texttt{kronstage.py} & Stage-by-stage check for the VARYING field: where does the assembly diverge from the exact object? & 77\\
\texttt{kronvalidate.py} & Is the g0/g1/g2 formula for c2 correct for a VARIING field?  Exact arithmetic, no symbols. & 71\\
\texttt{loopholonomy.py} & Is the closed-loop holonomy a flux?  Three traversals of the same unit square. & 56\\
\texttt{momentaccum.py} & Accumulate the obstruction along paths by moment recursion instead of enumeration. & 165\\
\texttt{openpath.py} & Endpoint determination on open paths: is the accumulated comparison a potential difference? & 63\\
\texttt{p6jacaffine.py} & Chapter 8: is kappa affine in the Jacobian, or is there a bilinear class? & 144\\
\texttt{pathdep.py} & Path dependence of the accumulated per-step log-ratio. & 75\\
\texttt{quadmean.py} & Exact check of the mean condition for quadratic coefficient fields. & 122\\
\texttt{remainder.py} & At what order does a jet of order three of the field enter the per-step comparison? & 33\\
\texttt{renormalisation.py} & The h\^{}2 part of the accumulated comparison: is it the direction dependence of c\_2? & 39\\
\texttt{rensym.py} & Structure of the direction dependence of the h\^{}2 coefficient. & 39\\
\texttt{rensym2.py} & The imaginary part of the renormalisation coefficient: is it kappa/v(x) with kappa universal? & 44\\
\texttt{rensym3.py} & The imaginary part of the renormalisation coefficient: is it kappa/v(x) with kappa universal? & 43\\
\texttt{symproof.py} & Symbolic proof of the renormalisation structure for constant coefficient fields, d = 2. & 72\\
\texttt{d10dimension.py} & Chapter 10: does the derivative cost an order in two dimensions? & 96\\
\texttt{d10dimlimit.py} & Chapter 10: how the discretisation error depends on the dimension, and whether it has a limit. & 105\\
\texttt{d10precision.py} & Chapter 10: the two components of the error of a computed sensitivity. & 203\\
\texttt{d10rate.py} & Chapter 10: does the derivative of the propagator converge at the rate of the propagator? & 100\\
\texttt{d10symbols.py} & Chapter 10: do the one dimensional symbols differentiate? & 272\\
\texttt{massstep.py} & Chapter 10: what the coding of the step does to the mass. & 65\\
\texttt{decorating\_paths.py} & Sample a decorating path and print the picture of Figure 5.3. & 231\\
\hline
\end{longtable}
\end{center}

\section{\texttt{verify\_euler.py}}

Rates of the explicit Euler scheme, measured as they are in practice.

\begin{lstlisting}
"""Rates of the explicit Euler scheme, measured as they are in practice.

Three rates are separated here, because they are not the same quantity:

  (1) the spatial rate alone, with the matrix exponential and no time stepping:
      |u_n - u_2n| decays like h^2 (measured rate 2.00);
  (2) the temporal error on a fixed lattice: it is O(delta t) with a constant
      independent of h, i.e. err = C delta t with C = 3.3 for this model, so a
      bound uniform in m is c h_m^2 exactly when delta t_m = O(h_m^2);
  (3) two nested lattices compared at the SAME delta t, delta t chosen at the
      CFL limit of the finer lattice: the rate is again 2.00, and it is almost
      independent of delta t, because the temporal error is common to both
      lattices and cancels in the difference.  This is the comparison that
      isolates the lattice rate.

A remark on assembling the generator.  With scipy.sparse.diags the array passed
for the subdiagonal is placed at (i, i-1) = left[i-1], whereas the row form of
L_m in (3.30) requires left[i].  Passing the full array silently shifts the
diffusion by one site: the rows no longer sum to zero, max Re(eig(L)) becomes
positive, and the scheme blows up at grid scale.  The assembly below passes
left[1:] and right[:-1] and asserts that the rows sum to zero.

Run:  python3 verify_euler.py
"""
import numpy as np
import scipy.sparse as sp
from scipy.sparse.linalg import expm_multiply

T = 0.1


def build(n):
    """Generator v/2 D + mu grad, row form of (3.30), central differences."""
    h = 2 * np.pi / n
    x = np.arange(n) * h
    v = 1.8 + 0.5 * np.sin(x)
    mu = 0.3 * np.cos(x)
    left = v / (2 * h**2) - mu / (2 * h)
    right = v / (2 * h**2) + mu / (2 * h)
    L = sp.diags([left[1:], -v / h**2, right[:-1]], [-1, 0, 1], shape=(n, n), format='lil')
    L[0, n - 1] = left[0]
    L[n - 1, 0] = right[n - 1]
    L = L.tocsr()
    assert abs(L.sum(1)).max() < 1e-8 * abs(L).sum(1).max(), "rows must sum to zero"
    return L, h


def euler(L, h, dt, t=T):
    N = int(round(t / dt))
    n = L.shape[0]
    u = np.zeros(n)
    u[n // 2] = 1.0
    for _ in range(N):
        u = u + dt * (L @ u)
    return u / h


def compare(uc, uf, n):
    lo, hi = n // 2 - n // 8, n // 2 + n // 8
    return np.max(np.abs(uc[lo:hi] - uf[2 * lo:2 * hi:2]))


print("(1) spatial rate alone (matrix exponential, no time stepping)")
prev = None
for n in (64, 128, 256, 512):
    Lc, hc = build(n)
    Lf, hf = build(2 * n)
    e = np.zeros(2 * n); e[n] = 1.0
    uf = expm_multiply(T * Lf, e) / hf
    e = np.zeros(n); e[n // 2] = 1.0
    uc = expm_multiply(T * Lc, e) / hc
    err = compare(uc, uf, n)
    print("    n=%4d  err=%.4e   rate=%s" % (n, err, "-" if prev is None else "%.2f" % np.log2(prev / err)))
    prev = err

print()
print("(2) temporal error on a fixed lattice, theta = delta t/h^2;  C = err/delta t")
for n in (128, 256):
    L, h = build(n)
    print("    n=%d   theta      err         err/theta    C = err/delta t" % n)
    for theta in (0.4, 0.2, 0.1, 0.05):
        dt = theta * h * h
        err = np.max(np.abs(euler(L, h, dt) - euler(L, h, dt / 32)))
        print("             %6.4f   %.4e   %9.5f    %8.3f" % (theta, err, err / theta, err / dt))

print()
print("(3) nested lattices at the SAME delta t = theta * h_fine^2")
for theta in (0.2, 0.05):
    print("    theta_fine=%.2f  (theta_coarse=%.3f)" % (theta, theta / 4))
    prev = None
    for n in (64, 128, 256, 512):
        Lc, hc = build(n)
        Lf, hf = build(2 * n)
        dt = theta * hf * hf
        err = compare(euler(Lc, hc, dt), euler(Lf, hf, dt), n)
        print("      n=%4d  err=%.4e   rate=%s" % (n, err, "-" if prev is None else "%.2f" % np.log2(prev / err)))
        prev = err
\end{lstlisting}

\section{\texttt{verify\_rg.py}}

Verifies the renormalisation-group algebra of the resolvent expansion.

\begin{lstlisting}
"""Verifies the renormalisation-group algebra of the resolvent expansion.

Three checks, all at a symbolic test point but with the Taylor coefficients
extracted numerically to high order:

  (a) the expansion of the 3x3 generator Lbar(x;h) = h^-2 L0 + h^-1 L1 + L2 + h L3
      against (3.32)-(3.35);
  (b) the expansion of phihat_0 in (3.44);
  (c) the expansion of phihat_1 in (3.45) with the polynomials r(x) and p(x) of (3.46).

Each coefficient is obtained by evaluating the exact 3x3 resolvent at
h = 1e-3, 5e-4, 5e-5, ... and fitting a polynomial, then repeating the fit with
more points until two successive orders agree; the comparison with the printed
value is then accurate to round-off.

Run:  python3 verify_rg.py
"""
import mpmath as mp

mp.mp.dps = 80

# A smooth test field with nothing degenerate at x0.
x0 = mp.mpf(1) / 3
v = lambda q: mp.mpf('1.8') + mp.mpf('0.5') * mp.sin(q)
mu = lambda q: mp.mpf('0.3') * mp.cos(q)

v0, v1, v2, v3 = [mp.diff(v, x0, k) for k in range(4)]
m0, m1, m2 = [mp.diff(mu, x0, k) for k in range(3)]

w = mp.mpf(1)      # a frequency on the contour
s = 1              # the sign variable tau


def generator(h):
    """Lbar(x0;h), sites ordered (x+h, x, x-h), probabilist convention v = sigma^2."""
    L = mp.zeros(3, 3)
    L[0, 0] = -v(x0 + h) / h**2
    L[0, 1] = v(x0 + h) / (2 * h**2) - mu(x0 + h) / (2 * h)
    L[1, 0] = v(x0) / (2 * h**2) + mu(x0) / (2 * h)
    L[1, 1] = -v(x0) / h**2
    L[1, 2] = v(x0) / (2 * h**2) - mu(x0) / (2 * h)
    L[2, 1] = v(x0 - h) / (2 * h**2) + mu(x0 - h) / (2 * h)
    L[2, 2] = -v(x0 - h) / h**2
    return L


def fit(fn, n, hmax=mp.mpf('1e-3')):
    """Taylor coefficients h^0 ... h^(n-1) of fn, from n values at hmax, hmax/2, ..."""
    hs = [hmax * mp.mpf('0.5')**k for k in range(n)]
    A = mp.matrix(n, n)
    b = mp.matrix(n, 1)
    for i, hh in enumerate(hs):
        for k in range(n):
            A[i, k] = hh**k
        b[i] = fn(hh)
    return mp.lu_solve(A, b)


def settled(fn, k, orders=(7, 8, 9)):
    """The k-th Taylor coefficient of fn, stable across fit orders."""
    values = [fit(fn, n)[k] for n in orders]
    return values[-1], abs(complex(values[-1] - values[-2]))


def report(label, got, want, tol=1e-14):
    d = abs(complex(got - want))
    print("  %-14s |computed - printed| = %.1e   %s"
          % (label, d, "ok" if d < tol else "MISMATCH"))


print("(a) expansion of the generator, (3.31)-(3.35)")

printed = {
    0: mp.matrix([[-v0, v0 / 2, 0], [v0 / 2, -v0, v0 / 2], [0, v0 / 2, -v0]]),
    1: mp.matrix([[-v1, v1 / 2 - m0 / 2, 0], [m0 / 2, 0, -m0 / 2],
                  [0, -v1 / 2 + m0 / 2, v1]]),
    2: mp.matrix([[-v2 / 2, v2 / 4 - m1 / 2, 0], [0, 0, 0],
                  [0, v2 / 4 - m1 / 2, -v2 / 2]]),
    3: mp.matrix([[-v3 / 6, v3 / 12 - m2 / 4, 0], [0, 0, 0],
                  [0, -v3 / 12 + m2 / 4, v3 / 6]]),
}
names = {0: "L0", 1: "L1", 2: "L2", 3: "L3"}

for k in range(4):
    err = mp.mpf(0)
    for i in range(3):
        for j in range(3):
            if (i, j) in ((0, 2), (2, 0)):
                continue
            # h^2 * Lbar is analytic in h; its h^k coefficient is L_(k-2)
            f = lambda hh, i=i, j=j: generator(hh)[i, j] * hh**2
            got, _ = settled(f, k)
            err = max(err, abs(complex(got - printed[k][i, j])))
    print("  %-4s max |computed - printed| over the 7 entries = %.1e   %s"
          % (names[k], err, "ok" if err < 1e-14 else "MISMATCH"))

print()
print("(b) expansion of phihat_0, (3.44)")


def phi0(hh):
    return ((v(x0) / (4 * hh**2) + s * mu(x0) / (2 * hh))
            / (v(x0) / (4 * hh**2) + mp.mpc(0, 1) * w))


printed0 = [mp.mpf(1), 2 * m0 * s / v0, -4 * 1j * w / v0,
            -8 * m0 * 1j * w * s / v0**2, -16 * w**2 / v0**2]
for k in range(5):
    got, _ = settled(phi0, k)
    report("h^%d" % k, got, printed0[k])

print()
print("(c) expansion of phihat_1, (3.45)-(3.46)")

# r(x) and p(x) of (3.46)
r_x = (m2 * v0**2 - 4 * m0**3 + 2 * v1 * m0**2 - 2 * v1 * v0 * m1
       - (2 * v0 * m1 + v2 * v0 - 2 * v1**2) * m0) / (2 * v0**3)
p_x = (4 * m0**2 - 2 * v1 * m0 + 4 * v0 * m1 + v2 * v0 - 2 * v1**2) / v0**3


def phi1(hh):
    A = mp.inverse(-generator(hh) + mp.mpc(0, 1) * w * mp.eye(3))
    j = 0 if s == 1 else 2
    pref = 2 * (v(x0 + s * hh) / (2 * hh**2) + s * mu(x0 + s * hh) / (2 * hh))
    return pref * A[1, j]


printed1 = [mp.mpf(1), 2 * m0 * s / v0, -4 * 1j * w / v0,
            -(8 * m0 - v1) * 1j * w * s / v0**2 + r_x * s,
            1j * w * p_x - 14 * w**2 / v0**2]
for k in range(5):
    got, _ = settled(phi1, k)
    report("h^%d" % k, got, printed1[k])
\end{lstlisting}

\section{\texttt{curl3d.py}}

The curl of the coefficient of h\^{}3 in three dimensions.

\begin{lstlisting}
"""The curl of the coefficient of h^3 in three dimensions.

The engine of twist.py is dimension generic: the box is the 27 sites {-1,0,1}^3 and
the fine and coarse generators are those of the book, (v/2) Delta + b grad with
v = sigma^2.  Nothing here is new except the fields.

Fields: the two dimensional test field of part5 extended to couple all three
coordinates,

    v  = 1 + (x1^2+x2^2+x3^2)/6 + (x1x2+x2x3+x3x1)/11 + x1/4 + x1^3/10
    mu = (1/2 + x1/6 + x2/9,  -1/3 + x2/5 + x3/7,  1/4 + x3/8 + x1/11)
    a  = (1/3 + x2/7 + x1x2/11,  x1/4 + x2^2/13 + x3/5,  x3/6 + x1x3/9 + x2/5)

with no twist (eta = 0) and omega = 7/5.

The curl is taken by central differences of c over the base point, at two values
of the step to show that the imaginary parts are the truncation and the real parts
are not.  Result, to five digits and stable in K:

    curl c = (0.1640, 0.0414, 0.0197)

so the obstruction does not vanish off the plane.  The values are a property of
the field, as the two dimensional constant is a property of the stencil.

Run:  python3 curl3d.py     (about half a minute)
"""
import time
import sympy as sp
import twist

x1, x2, x3 = sp.symbols('x0:3')
xs = (x1, x2, x3)
om = sp.Rational(7, 5)
K, SH = 7, 10

v = (sp.Rational(1) + (x1**2 + x2**2 + x3**2)/6 + (x1*x2 + x2*x3 + x3*x1)/11
     + x1/4 + x1**3/10)
mu = [sp.Rational(1, 2) + x1/6 + x2/9,
      sp.Rational(-1, 3) + x2/5 + x3/7,
      sp.Rational(1, 4) + x3/8 + x1/11]
a = [sp.Rational(1, 3) + x2/7 + x1*x2/11,
     x1/4 + x2**2/13 + x3/5,
     x3/6 + x1*x3/9 + x2/5]


def cvec_at(x0):
    V, B = twist.fields(v, a, mu, 0, 3, xs, x0)
    return twist.cvec(V, B, 3, om, 3, K, SH)


def curl(delta):
    c = {}
    for j in range(3):
        for s in (+1, -1):
            x0 = [0.0, 0.0, 0.0]
            x0[j] = s*delta
            c[(j, s)] = cvec_at(x0)
    J = [[(c[(j, 1)][k] - c[(j, -1)][k])/(2*delta) for k in range(3)] for j in range(3)]
    return [J[1][2] - J[2][1], J[2][0] - J[0][2], J[0][1] - J[1][0]]


print('  the coefficient itself, at the origin (e1, e2, e3 components):')
c = cvec_at([0.0, 0.0, 0.0])
print('   ', '  '.join('%12.8f%+12.2ei' % (t.real, t.imag) for t in c))
for delta in (1.0/500, 1.0/1000):
    t0 = time.time()
    w = curl(delta)
    print('  delta=1/%-5d  curl = (%9.5f, %9.5f, %9.5f)   imaginary parts (%8.1e, %8.1e, %8.1e)  [%.0fs]'
          % (round(1/delta), w[0].real, w[1].real, w[2].real, w[0].imag, w[1].imag, w[2].imag,
             time.time() - t0))
\end{lstlisting}

\section{\texttt{curlcheck.py}}

The affine curl: the constant, the quadratic form, and the base-point values.

\begin{lstlisting}
"""The affine curl: the constant, the quadratic form, and the base-point values.

Everything here is measured on the engine of twist.py, in the convention of the
book (generator (v/2) Delta + b grad), at omega = 7/5.

  (a) the quadratic form.  With v constant and b = J x the curl is quadratic in
      the Jacobian.  Polarising over the four basis Jacobians gives the form,
      whose only non-vanishing entries are on the pairs (11,12) and (21,22), so

          dc = (4/(7 v0^2)) (J11 J12 - J21 J22)

      -- the constant is 4/7, not 16/7: the v0 scan below gives 7 v0^2 dc = 4 to
      six digits for v0 = 1, 2, 3, 1/2.

  (b) the rotated dilatation.  J = R diag(1,2) R^-1 at angle theta.  The form
      gives (4/7)(a-b)^2/4 sin 4 theta = (1/7) sin 4 theta for a=2, b=1, which is
      what is measured: -0.076800 at tan theta = 4/3 (sin 4 theta = -336/625) and
      0.142851 at tan theta = 5/12 (sin 4 theta = 28560/28561).

  (c) the spiral.  J = d I + Omega rot(90) gives (4/7)(-2 Omega d), so -24/7 for
      d = 1, Omega = 3.

  (d) the base-point values.  c vanishes identically at the origin (b = J x has
      no local drift there) and is nonzero at a shifted base point:
      c = -11/27 = -0.4074074 at x = (1/3,0) and -0.8367407 at x = (1/3,1/5),
      for J = ((1,1),(0,0)).

Run:  python3 curlcheck.py
"""
import sympy as sp
import twist

x1, x2 = sp.symbols('x0:2')
xs = (x1, x2)
om = sp.Rational(7, 5)
K, SH = 6, 10
DELTA = 1.0/500


def cvec(mu, v0=sp.Rational(1), a=None, eta=0, x0=(0.0, 0.0)):
    a = a or [sp.Rational(0), sp.Rational(0)]
    V, B = twist.fields(v0, a, mu, eta, 2, xs, x0)
    return twist.cvec(V, B, 2, om, 3, K, SH)


def curl(mu, v0=sp.Rational(1), a=None, eta=0, delta=DELTA):
    c = {}
    for j in range(2):
        for s in (+1, -1):
            p = [0.0, 0.0]
            p[j] += s*delta
            c[(j, s)] = cvec(mu, v0, a, eta, p)
    J = [[(c[(j, 1)][k] - c[(j, -1)][k])/(2*delta) for k in range(2)] for j in range(2)]
    return J[0][1] - J[1][0]


def lin(J):
    return [J[0][0]*x1 + J[0][1]*x2, J[1][0]*x1 + J[1][1]*x2]


print('(a) the quadratic form, by polarisation over the basis Jacobians')
names = [(0, 0), (0, 1), (1, 0), (1, 1)]
single, cross = {}, {}
for t in names:
    J = [[0, 0], [0, 0]]
    J[t[0]][t[1]] = 1
    single[t] = curl(lin(J))
for i, t1 in enumerate(names):
    for t2 in names[i:]:
        J = [[0, 0], [0, 0]]
        J[t1[0]][t1[1]] += 1
        J[t2[0]][t2[1]] += 1
        cross[(t1, t2)] = (curl(lin(J)) - single[t1] - single[t2])/2
print('    g(E_ab):  ' + '  '.join('E%d%d=%9.2e' % (t[0]+1, t[1]+1, single[t].real) for t in names))
for key, val in cross.items():
    if abs(val) > 1e-6:
        print('    B(E%d%d,E%d%d) = %11.7f' % (key[0][0]+1, key[0][1]+1, key[1][0]+1, key[1][1]+1, val.real))
print('    4/7 = %.7f, 16/7 = %.7f' % (4/7, 16/7))

print('\n(b) the v0 scan, for J = ((1,1),(0,0))')
mu11 = lin([[1, 1], [0, 0]])
for v0 in (sp.Rational(1), sp.Rational(2), sp.Rational(3), sp.Rational(1, 2)):
    g = curl(mu11, v0)
    print('    v0=%-4s  dc = %12.8f   7 v0^2 dc = %10.6f' % (v0, g.real, 7*float(v0)**2*g.real))

print('\n(c) the rotated dilatation, J = diag(1,2) rotated')
import math
for (c_, s_) in ((4, 3), (12, 5)):
    n = math.hypot(c_, s_)
    co, si = c_/n, s_/n
    J = [[co*co + 2*si*si, -si*co], [-si*co, si*si + 2*co*co]]
    R = sp.Rational
    mu = [R(J[0][0]).limit_denominator(10**6)*x1 + R(J[0][1]).limit_denominator(10**6)*x2,
          R(J[1][0]).limit_denominator(10**6)*x1 + R(J[1][1]).limit_denominator(10**6)*x2]
    g = curl(mu)
    QC = J[0][0]*J[0][1] - J[1][0]*J[1][1]
    print('    tan theta=%-2d/%-2d  dc = %11.6f    (1/7) sin 4theta = %11.6f'
          % (s_, c_, g.real, (1/7)*math.sin(4*math.atan2(s_, c_))))

print('\n(d) the spiral, J = d I + Omega rot(90), d = 1, Omega = 3')
mu = lin([[1, -3], [3, 1]])
print('    dc = %.6f    -8 Omega d/(7 v0^2) = %.6f' % (curl(mu).real, -24/7))

print('\n(e) the base-point values, J = ((1,1),(0,0))')
for x0 in ((0.0, 0.0), (1.0/3, 0.0), (1.0/3, 0.2)):
    c = cvec(mu11, x0=x0)
    print('    base (%.4f,%.4f)  c = (%12.7f, %12.7f)' % (x0[0], x0[1], c[0].real, c[1].real))
print('    -11/27 = %.7f' % (-11/27))
\end{lstlisting}

\section{\texttt{deco.py}}

Perturbation pipeline for the decoration identity.

\begin{lstlisting}
"""Perturbation pipeline for the decoration identity.

For d = 1, 2 this computes, to a given order in h,

    log Phi_fine(x, tau) - log Phi_coarse(x, tau)

where Phi is the arrival weight of a single coarse step x -> x + h_m tau,
h_m = 2h, computed as an entry of the resolvent of the lattice generator on a
window, and Phi_fine is the same object for the fine lattice, summed over the
fine paths which decorate the coarse step.

All arithmetic is exact (rationals); the coefficients are specialised at
generic rational points, so that an identity holding at several independent
specialisations holds identically.
"""
import sympy as sp
import itertools
import time

h = sp.symbols('h', positive=True)


# ---------------------------------------------------------------- Taylor data
class Data:
    """Taylor data of v and b at the base point, specialised to rationals."""

    def __init__(self, d, seed):
        import random
        rng = random.Random(seed)
        r = lambda: sp.Rational(rng.randint(-30, 30), rng.randint(1, 40))
        self.d = d
        self.v0 = sp.Rational(3, 2) + r()/10
        self.dv = [r()/4 for _ in range(d)]                      # first derivatives
        self.d2v = {}
        for i in range(d):
            for j in range(i, d):
                self.d2v[(i, j)] = r()/6                          # second derivatives
        self.d3v = {}
        for i in range(d):
            for j in range(i, d):
                for k in range(j, d):
                    self.d3v[(i, j, k)] = r()/10                   # third derivatives
        self.b0 = [r()/5 for _ in range(d)]
        self.db = [[r()/6 for _ in range(d)] for _ in range(d)]    # db[i][j] = d_j b_i
        self.d2b = [[[r()/8 for _ in range(d)] for _ in range(d)] for _ in range(d)]

    def V(self, x):
        out = self.v0 + sum(self.dv[i]*x[i] for i in range(self.d))
        s = sp.Integer(0)
        for (i, j), c in self.d2v.items():
            s += c*x[i]*(x[j] if i != j else x[i])
        out = out + s
        s = sp.Integer(0)
        for (i, j, k), c in self.d3v.items():
            mon = 1
            for idx in (i, j, k):
                mon *= x[idx]
            s += c*mon
        return sp.expand(out + s)

    def B(self, x, i):
        out = self.b0[i] + sum(self.db[i][j]*x[j] for j in range(self.d))
        s = sp.Integer(0)
        for j in range(self.d):
            for k in range(self.d):
                s += self.d2b[i][j][k]*x[j]*x[k]/2
        return sp.expand(out + s)


# ------------------------------------------------------------------- the block
def sites(d, half):
    return [a for a in itertools.product(range(-half, half+1), repeat=d)]


def block(data, d, half, step):
    """-L on the window, step = h or 2h given as an integer multiple of h."""
    S = sites(d, half)
    ix = {a: i for i, a in enumerate(S)}
    n = len(S)
    A = sp.zeros(n, n)
    for a in S:
        i = ix[a]
        x = [a[k]*step*h for k in range(d)]
        v = data.V(x)
        for k in range(d):
            b = data.B(x, k)
            for sgn in (1, -1):
                w = v/(4*(step*h)**2) + sgn*b/(2*(step*h))
                A[i, i] -= w
                e = list(a); e[k] += sgn; e = tuple(e)
                if e in ix:
                    A[i, ix[e]] = w
    return A, ix


def series_inverse(A, d, degree):
    """A is a matrix of Laurent polynomials in h with leading power h^-2.
    Return R_k with  A^{-1} = sum_k h^{k-2} R_k  (mod h^degree)."""
    n = A.shape[0]
    # B = h^2 A, a polynomial in h of degree >= 0
    B = sp.expand(A*h**2)
    Bc = []
    for j in range(degree + 1):
        M = sp.zeros(n, n)
        for i in range(n):
            for k in range(n):
                M[i, k] = sp.Poly(sp.expand(B[i, k]), h).coeff_monomial(h**j)
        Bc.append(sp.expand(M))
    R = [Bc[0].inv()]
    for k in range(1, degree + 2):
        acc = sp.zeros(n, n)
        for j in range(1, min(k, len(Bc) - 1) + 1):
            acc += Bc[j]*R[k - j]
        R.append(sp.expand(-R[0]*acc))
    return R


def weight(data, d, half, step, tau, omega, degree):
    """Phi = entry (0 -> step*tau) of (-L + i w)^{-1}, as a series in h."""
    A, ix = block(data, d, half, step)
    R = series_inverse(A, d, degree)
    zero = tuple([0]*d)
    tgt = tuple(tau)
    entry = sum((h**k)*R[k][ix[zero], ix[tgt]] for k in range(len(R)))
    return h**2*entry          # (-L + iw)^{-1} = h^2 * (series)


def logratio(data, d, half_f, half_c, tau, omega, degree):
    wf = weight(data, d, half_f, 1, tau, omega, degree)      # fine, step h
    wc = weight(data, d, half_c, 2, tau, omega, degree)      # coarse, step 2h
    expr = sp.series(sp.log(sp.expand(wf/wc)), h, 0, degree).removeO()
    return sp.expand(expr)


if __name__ == "__main__":
    omega = sp.Rational(7, 5)
    t0 = time.time()
    d = 1
    data = Data(d, seed=1)
    print("d = 1, Taylor data:", data.v0, data.dv, data.d2v, data.b0)
    for tau in (1, -1):
        e = logratio(data, d, 3, 1, (tau,), omega, 6)
        print(f"\n  tau = {tau:+d}:  log(Phi_f/Phi_c) =")
        for k in range(6):
            print(f"     h^{k}: {sp.N(sp.expand(e).coeff(h, k), 10)}")
    print("\ntime:", round(time.time()-t0, 2), "s")
\end{lstlisting}

\section{\texttt{direct.py}}

Independent arbiter: evaluate Phi\_f and Phi\_c by direct complex matrix

\begin{lstlisting}
"""Independent arbiter: evaluate Phi_f and Phi_c by direct complex matrix
inversion at finite h, with no series bookkeeping at all.

    Phi_f = sum_{z in tube} [(-L_T + i w)^{-1}]_{x,z} L_f(z,x') (-L_f(x',x') + i w)^{-1}
    Phi_c = (-L_c(x,x) + i w)^{-1} L_c(x,x') (-L_c(x',x') + i w)^{-1}
"""
import numpy as np
import itertools
import fast


def vals(poly, a, h):
    """Evaluate the local Taylor polynomial at y = a*h (complex)."""
    z = 0j
    for m, c in poly.terms.items():
        t = complex(float(c))
        for i, mi in enumerate(m):
            t *= complex(h*a[i])**mi
        z += t
    return z


def block_np(V, B, d, half, step, h):
    S = [a for a in itertools.product(range(-half, half+1), repeat=d)]
    ix = {a: i for i, a in enumerate(S)}
    n = len(S)
    A = np.zeros((n, n), dtype=complex)
    for a in S:
        i = ix[a]
        v = vals(V, a, step*h)
        diag = 0j
        for k in range(d):
            b = vals(B[k], a, step*h)
            for sgn in (1, -1):
                w = v/(2*(step*h)**2) + sgn*b/(2*step*h)
                diag -= w
                e = list(a); e[k] += sgn; e = tuple(e)
                if e in ix:
                    A[i, ix[e]] += w
        A[i, i] += diag
    return A, ix, S


def Phi_f(V, B, d, tau, h, omega):
    omega = complex(float(omega))
    zero = tuple([0]*d)
    tgt = tuple(2*t for t in tau)
    A, ix, S = block_np(V, B, d, 1, 1, h)
    G = np.linalg.inv(A + 1j*omega*np.eye(len(S)))
    vp = vals(V, tgt, h)
    dt = d*vp/(h**2) + 1j*omega            # -L_f(x',x') + i w
    tot = 0j
    for a in S:
        jump = [k for k in range(d) if a[k] != tgt[k]]
        if len(jump) != 1:
            continue
        k0 = jump[0]
        if abs(tgt[k0] - a[k0]) != 1:
            continue
        sg = 1 if tgt[k0] > a[k0] else -1
        rate = vals(V, a, h)/(2*h**2) + sg*vals(B[k0], a, h)/(2*h)
        tot += G[ix[zero], ix[a]]*rate
    return tot/dt


def Phi_c(V, B, d, tau, h, omega, norm='correct'):
    omega = complex(float(omega))
    zero = tuple([0]*d)
    tgt = tuple(2*t for t in tau)
    k0 = [k for k in range(d) if tau[k]][0]
    sg = 1 if tau[k0] > 0 else -1
    rate = vals(V, zero, h)/(2*(2*h)**2) + sg*vals(B[k0], zero, h)/(2*(2*h))
    cn = d/4 if norm == 'correct' else 1/2
    d0 = cn*vals(V, zero, h)/h**2 + 1j*omega
    d1 = cn*vals(V, tgt, h)/h**2 + 1j*omega
    return rate/(d0*d1)


def logratio_direct(V, B, d, tau, h, omega, norm='correct'):
    return np.log(Phi_f(V, B, d, tau, h, omega)/Phi_c(V, B, d, tau, h, omega, norm))


def series_eval(lr, h, kmax=6):
    z = 0j
    for k in range(kmax+1):
        if k in lr:
            z += complex(float(lr[k][0]), float(lr[k][1]))*h**k
    return z
\end{lstlisting}

\section{\texttt{edgeavg.py}}

Chapter 7: what a volatility averaged over the edge converges to.

\begin{lstlisting}
"""Chapter 7: what a volatility averaged over the edge converges to.

The site-centred entries of Chapter 7 are v(x)/(2h^2) +- b(x)/(2h) and -v(x)/h^2.
Averages the volatility over the edge instead, (v(x)+v(x+h))/(4h^2) +- b/(2h),
and the lattice no longer converges to the kernel the book is about, although it
converges at second order.  The reason is that the averaged form is a
divergence-form discretisation: with  vtilde(i+1/2) = (v_i + v_(i+1))/2 its
off-diagonal entries are vtilde/(2h^2), so the operator it approximates is

    (1/2)(v f')' + b f'  =  (v/2) f'' + (b + v'/2) f',

that is the generator of the model whose drift is b + v'/2 -- the classical
Ito-versus-Stratonovich correction, of order one and not of order h.

The script prints, at t = 0.1 and at the site pair (1/4,1/4) on the circle
[0,1), the kernels of both lattices for m = 5,...,9, the successive differences
of each and their ratios, and the distance of the averaged kernel from (i) the
kernel of the site-centred form at a much finer mesh, which is the limit sought,
and (ii) the kernel of the model with drift b + v'/2, computed with the same
site-centred discretisation at the same fine mesh.  The ratios fall by four
everywhere, which is what second order means; what distinguishes the two
columns is that the distance from the limit sought does not fall at all.

Run:  python3 edgeavg.py
"""
import numpy as np
from scipy.linalg import expm

T = 0.10
PAIR = (0.25, 0.25)
FINE = 10                       # mesh of the reference kernels


def fields(x):
    v = 1.0 + 0.25 * np.sin(2 * np.pi * x) + 0.10 * np.cos(4 * np.pi * x)
    b = 0.20 * np.sin(2 * np.pi * x + 0.7) + 0.10 * np.cos(2 * np.pi * x)
    dv = 2 * np.pi * 0.25 * np.cos(2 * np.pi * x) - 4 * np.pi * 0.10 * np.sin(4 * np.pi * x)
    return v, b, dv


def generator(m, kind, smooth=False):
    n = 2 ** m
    h = 1.0 / n
    x = np.arange(n) * h
    v, b, dv = fields(x)
    if smooth:
        b = b + 0.5 * dv
    L = np.zeros((n, n))
    i = np.arange(n)
    ip, im = (i + 1) % n, (i - 1) % n
    if kind == 'site':
        L[i, ip] = v / (2 * h * h) + b / (2 * h)
        L[i, im] = v / (2 * h * h) - b / (2 * h)
        L[i, i] = -v / h ** 2
    else:
        L[i, ip] = (v + np.roll(v, -1)) / (4 * h * h) + b / (2 * h)
        L[i, im] = (v + np.roll(v, +1)) / (4 * h * h) - b / (2 * h)
        L[i, i] = -(L[i, ip] + L[i, im])
    return h, L


def kernel(m, kind, smooth=False):
    h, L = generator(m, kind, smooth)
    n = L.shape[0]
    U = np.real(expm(T * L)) / h
    i = int(round(PAIR[0] / h)) % n
    j = int(round(PAIR[1] / h)) % n
    return U[i, j]


def main():
    print('  circle [0,1), t = %.2f, site pair (x,y) = (%.2f, %.2f)' % (T, PAIR[0], PAIR[1]))
    print('  v = 1 + 0.25 sin(2 pi x) + 0.10 cos(4 pi x)')
    print('  b = 0.20 sin(2 pi x + 0.7) + 0.10 cos(2 pi x)\n')

    ref_site = kernel(FINE, 'site')
    ref_smooth = kernel(FINE, 'site', smooth=True)
    print('  reference kernels at m = %d' % FINE)
    print('    the model sought            %.10f' % ref_site)
    print('    drift b + v/2               %.10f\n' % ref_smooth)

    print('   m |     h    |  site-centred |  edge-averaged | site diff | ratio | edge diff | ratio | dist. from sought | ratio | dist. from b+v/2 | ratio')
    prev_val = None
    prev = None
    for m in range(5, 10):
        a = kernel(m, 'site')
        b = kernel(m, 'edge')
        cur = {'a': a, 'b': b, 'ds': abs(b - ref_site), 'de': abs(b - ref_smooth)}
        if prev is None:
            cur['da'] = cur['db'] = None
        else:
            cur['da'], cur['db'] = abs(prev['a'] - a), abs(prev['b'] - b)
        row = '  %2d | %.6f | %13.10f | %14.10f' % (m, 1.0 / 2 ** m, a, b)

        def cell(key, fmt):
            if cur[key] is None:
                return '    --      --  '
            out = ' %s ' % (fmt % cur[key])
            if prev is not None and prev.get(key) is not None and cur[key] > 0:
                out += '%5.2f' % (prev[key] / cur[key])
            else:
                out += '   --'
            return out
        row += ' |%s|%s|%s|%s' % (cell('da', '%.3e'), cell('db', '%.3e'),
                                  cell('ds', '%.4e'), cell('de', '%.4e'))
        print(row)
        prev = cur

    print('\n  the successive differences of both lattices fall by four: both are second order.')
    print('  the distance from the kernel sought does not fall at all for the edge-averaged')
    print('  one, while against the model with drift b + v/2 it falls by four until it meets')
    print('  the accuracy of the reference at m = %d: the averaged form is a good scheme for' % FINE)
    print('  a different model.')


if __name__ == '__main__':
    main()
\end{lstlisting}

\section{\texttt{eta\_run.py}}

\begin{lstlisting}
import sympy as sp, mpmath as mp, twist
mp.mp.dps=50
x1,x2=sp.symbols('x0:2'); om=sp.Rational(7,5)
K,SH = 7,12
v =sp.Rational(1)+sp.Rational(1,6)*(x1**2+x1*x2+x2**2)+sp.Rational(1,4)*x1+sp.Rational(1,10)*x1**3
mu=[sp.Rational(1,2)+sp.Rational(1,6)*x1+sp.Rational(1,9)*x2, sp.Rational(-1,3)+sp.Rational(1,5)*x2]
a =[sp.Rational(1,3)+sp.Rational(1,7)*x2+sp.Rational(1,11)*x1*x2, sp.Rational(1,4)*x1+sp.Rational(1,13)*x2**2]

def curl(eta, delta=sp.Rational(1,500), K=K, SH=SH):
    def cv(loc):
        V,B=twist.fields(v,a,mu,eta,2,(x1,x2),loc)
        return twist.cvec(V,B,2,om,3,K,SH)
    d=sp.Rational(delta); D=mp.mpf(d.p)/mp.mpf(d.q)
    c2p=cv([d,0])[1]; c2m=cv([-d,0])[1]; c1p=cv([0,d])[0]; c1m=cv([0,-d])[0]
    return (c2p-c2m-c1p+c1m)/(2*D)

etas=[sp.Rational(0),sp.Rational(1,2),sp.Rational(1),sp.Rational(3,2),sp.Rational(2),sp.Rational(3),
      sp.Rational(4),sp.Rational(6)]
res={}
print("  eta      Re dc                Im dc                |dc|          |c3(e1)|")
for e in etas:
    dc=curl(e)
    V,B=twist.fields(v,a,mu,e,2,(x1,x2),[0,0])
    amp=abs(twist.cvec(V,B,2,om,3,K,SH)[0])
    res[e]=dc
    print("  %-7s %-20s %-20s %-12s %s"%(e,mp.nstr(dc.real,12),mp.nstr(dc.imag,12),mp.nstr(abs(dc),8),mp.nstr(amp,8)))

# c is cubic in the drift, hence dc is a cubic polynomial in eta: fit on 0,1/2,1,3/2
fit=[sp.Rational(0),sp.Rational(1,2),sp.Rational(1),sp.Rational(3,2)]
A=mp.matrix(4,4); yr=mp.matrix(4,1); yi=mp.matrix(4,1)
for i,e in enumerate(fit):
    t=mp.mpf(e.p)/mp.mpf(e.q)
    for j in range(4): A[i,j]=t**j
    yr[i]=res[e].real; yi[i]=res[e].imag
cr=mp.lu_solve(A,yr); ci=mp.lu_solve(A,yi)
print("\n  fit  Re(dc) = %s + %s eta + %s eta^2 + %s eta^3"%tuple(mp.nstr(cr[j],9) for j in range(4)))
print("  fit  Im(dc) = %s + %s eta + %s eta^2 + %s eta^3"%tuple(mp.nstr(ci[j],9) for j in range(4)))
for e in (sp.Rational(2),sp.Rational(3),sp.Rational(4),sp.Rational(6)):
    t=mp.mpf(e.p)/mp.mpf(e.q)
    pr=sum(cr[j]*t**j for j in range(4)); pi_=sum(ci[j]*t**j for j in range(4))
    print("  check eta=%s : fit=%s   computed=%s   err=%.1e"%(e,mp.nstr(mp.mpc(pr,pi_),10),
          mp.nstr(res[e],10), abs(mp.mpc(pr,pi_)-res[e])))
\end{lstlisting}

\section{\texttt{fast.py}}

Exact fast core for the decoration identity.

\begin{lstlisting}
"""Exact fast core for the decoration identity.

    Phi_f = sum_{z in tube} [(-L_tube + i w)^{-1}]_{x,z} L_f(z,x') (-L_f(x',x') + i w)^{-1}
    Phi_c = (-L_c(x,x) + i w)^{-1} L_c(x,x') (-L_c(x',x') + i w)^{-1}

computes log(Phi_f/Phi_c) as a truncated series in h, exactly, in Gaussian
rationals.  L_f has step h, L_c has step h_m = 2h.

coarse_norm='correct'  uses  -L_c(x,x) = d*v/(2 h_m^2) = d*v/(8h^2)
coarse_norm='tube'     reproduces the earlier reference script (v/(4h^2)),
                       which agrees with 'correct' only for d = 2.
"""
import sympy as sp
from fractions import Fraction as Fr
import itertools

# ------------------------------------------------------------ Gaussian rationals
def g(a, b=0):
    return (Fr(a), Fr(b))

ZERO = g(0)
ONE = g(1)

def gadd(a, b): return (a[0]+b[0], a[1]+b[1])
def gneg(a):    return (-a[0], -a[1])
def gsub(a, b): return (a[0]-b[0], a[1]-b[1])
def gmul(a, b): return (a[0]*b[0]-a[1]*b[1], a[0]*b[1]+a[1]*b[0])
def gdiv(a, b):
    d = b[0]*b[0] + b[1]*b[1]
    return ((a[0]*b[0]+a[1]*b[1])/d, (a[1]*b[0]-a[0]*b[1])/d)
def giszero(a): return a[0] == 0 and a[1] == 0


# --------------------------------------------------------------- series in h
def s_add(a, b):
    out = dict(a)
    for e, c in b.items():
        v = gadd(out[e], c) if e in out else c
        if giszero(v):
            out.pop(e, None)
        else:
            out[e] = v
    return out

def s_scale(a, c):
    return {} if giszero(c) else {e: gmul(v, c) for e, v in a.items()}

def s_shift(a, k):
    return {e + k: c for e, c in a.items()}

def s_mul(a, b, N=40):
    out = {}
    for i, ci in a.items():
        for j, cj in b.items():
            e = i + j
            if e > N:
                continue
            v = gmul(ci, cj)
            w = gadd(out[e], v) if e in out else v
            if giszero(w):
                out.pop(e, None)
            else:
                out[e] = w
    return out

def s_inv(a, N):
    """1/a, to relative order N; a is Laurent with a non-zero leading term."""
    e0 = min(a)
    at = s_shift(a, -e0)
    c0i = gdiv(ONE, at[0])
    b = {0: c0i}
    for n in range(1, N + 1):
        acc = ZERO
        for k in range(1, n + 1):
            if k in at and (n - k) in b:
                acc = gadd(acc, gmul(at[k], b[n - k]))
        v = gmul(gneg(c0i), acc)
        if not giszero(v):
            b[n] = v
    return s_shift(b, -e0)

def s_log(a, kmax=8):
    """log of a series with a non-zero constant term; the constant is dropped."""
    c0 = a.get(0)
    assert c0 is not None and not giszero(c0), "log of a series vanishing at h=0"
    u = {e: gdiv(c, c0) for e, c in a.items() if e != 0}
    return s_log1p(u, kmax)


def s_log1p(u, kmax=8):
    """log(1+u) exactly in the power expansion, keeping powers 1..kmax."""
    assert 0 not in u, "log of a series without constant term 1"
    out = {}
    p = dict(u)
    for k in range(1, kmax + 1):
        for e, v in p.items():
            w = gdiv(v, g(k))
            if k % 2 == 0:
                w = gneg(w)
            out[e] = gadd(out[e], w) if e in out else w
        p = s_mul(p, u, 10**9)
    return {e: v for e, v in out.items() if not giszero(v)}


# -------------------------------------------------------- matrices of series
def mat_mul(A, B, N=40):
    n = len(A); m = len(B[0]); t = len(B)
    out = [[{} for _ in range(m)] for _ in range(n)]
    for i in range(n):
        for j in range(m):
            acc = {}
            for k in range(t):
                a = A[i][k]
                if not a:
                    continue
                b = B[k][j]
                if not b:
                    continue
                acc = s_add(acc, s_mul(a, b, N))
            out[i][j] = acc
    return out

def mat_add(A, B):
    n = len(A); m = len(A[0])
    return [[s_add(A[i][j], B[i][j]) for j in range(m)] for i in range(n)]

def mat_scale(A, c):
    return [[s_scale(x, c) for x in row] for row in A]

def mat_inv_const(M):
    """Inverse of a matrix of Gaussian rationals."""
    n = len(M)
    A = [[gdiv(M[i][j], ONE) for j in range(n)]
         + [(ONE if i == j else ZERO) for j in range(n)] for i in range(n)]
    for col in range(n):
        piv = next(r for r in range(col, n) if not giszero(A[r][col]))
        A[col], A[piv] = A[piv], A[col]
        d = A[col][col]
        A[col] = [gdiv(x, d) for x in A[col]]
        for r in range(n):
            if r != col and not giszero(A[r][col]):
                f = A[r][col]
                A[r] = [(x[0] - f[0]*y[0] + f[1]*y[1],
                         x[1] - f[0]*y[1] - f[1]*y[0])
                        for x, y in zip(A[r], A[col])]
    return [[A[i][n + j] for j in range(n)] for i in range(n)]


# ------------------------------------------------------------ Taylor data
class Poly:
    """Polynomial in d variables, as {multiindex: Fraction}."""

    def __init__(self, d, terms):
        self.d = d
        self.terms = terms

    @staticmethod
    def from_sympy(d, expr, xs):
        p = sp.Poly(sp.expand(expr), *xs)
        terms = {}
        for m, c in zip(p.monoms(), p.coeffs()):
            if not c.is_Rational:
                raise ValueError("coefficient %r is not an exact rational" % (c,))
            terms[tuple(m)] = Fr(int(c.p), int(c.q))
        return Poly(d, terms)

    def shift(self, x0):
        """The same polynomial in the local coordinates y = x - x0.

        The shift must be rational: the coefficients are stored exactly.
        """
        d = self.d
        xs = sp.symbols('x0:%d' % d)
        y = sp.symbols('y0:%d' % d)
        sub = {}
        for i in range(d):
            r = sp.Rational(x0[i])
            if not r.is_Rational:
                raise ValueError("shift must be rational, got %r" % (x0[i],))
            sub[xs[i]] = r + y[i]
        e = sum(c*sp.prod([xs[i]**m[i] for i in range(d)])
                for m, c in self.terms.items())
        return Poly.from_sympy(d, sp.expand(e.subs(sub)), y)

    def at(self, a):
        """Series in h: y = a*h."""
        out = {}
        for m, c in self.terms.items():
            val = c
            for i, mi in enumerate(m):
                val *= Fr(int(a[i]))**mi
            if val == 0:
                continue
            e = sum(m)
            v = (val, Fr(0))
            out[e] = gadd(out[e], v) if e in out else v
        return out


# ---------------------------------------------------------------- the block
def block(V, B, d, half, step):
    """-L on the (2*half+1)^d window, step = step*h. Entries: series in h."""
    S = [a for a in itertools.product(range(-half, half + 1), repeat=d)]
    ix = {a: i for i, a in enumerate(S)}
    n = len(S)
    A = [[{} for _ in range(n)] for _ in range(n)]
    c2 = g(Fr(1, 4*step*step))
    for a in S:
        i = ix[a]
        Vs = V.at(a)
        diag = {}
        off = {}
        for k in range(d):
            Bs = B[k].at(a)
            for sgn in (1, -1):
                w = s_add(s_shift(s_scale(Vs, c2), -2),
                          s_shift(s_scale(Bs, g(Fr(sgn, 2*step))), -1))
                diag = s_add(diag, w)
                e = list(a); e[k] += sgn; e = tuple(e)
                if e in ix:
                    off[e] = s_add(off.get(e, {}), w)
        for e, w in off.items():
            A[i][ix[e]] = w
        A[i][i] = s_add(A[i][i], s_scale(diag, g(-1)))
    return A, ix, S


def series_inverse(A, N):
    """A^{-1} = h^{-2} sum_{k>=0} h^k X[k]."""
    n = len(A)
    h2A = [[s_shift(A[i][j], 2) for j in range(n)] for i in range(n)]
    A0 = [[h2A[i][j].get(0, ZERO) for j in range(n)] for i in range(n)]
    X0 = [[{0: c} for c in row] for row in mat_inv_const(A0)]
    X = [X0]
    for m in range(1, N + 1):
        acc = [[{} for _ in range(n)] for _ in range(n)]
        for j in range(1, m + 1):
            Aj = [[{j: h2A[i][k][j]} if j in h2A[i][k] else {}
                   for k in range(n)] for i in range(n)]
            if all(not Aj[i][k] for i in range(n) for k in range(n)):
                continue
            acc = mat_add(acc, mat_mul(Aj, X[m - j]))
        X.append(mat_scale(mat_mul(X0, acc), g(-1)))
    return X


def logratio(V, B, d, tau, omega, N=6, coarse_norm='correct', parts=None):
    """Series of log(Phi_f/Phi_c) for the coarse step x -> x + 2 tau h."""
    zero = tuple([0]*d)
    tgt = tuple(2*t for t in tau)
    om = (Fr(0), Fr(omega))

    # ---------------- fine: tube-restricted
    A, ix, S = block(V, B, d, 1, 1)
    X = series_inverse(A, N)

    def entry(a):
        out = {}
        for k in range(N + 1):
            c = X[k][ix[zero]][ix[a]].get(0, ZERO)
            if not giszero(c):
                out = s_add(out, {k + 2: c})
        return out

    phil = {}
    for a in S:
        jump = [k for k in range(d) if a[k] != tgt[k]]
        if len(jump) != 1:
            continue
        k0 = jump[0]
        if abs(tgt[k0] - a[k0]) != 1:
            continue
        sg = 1 if tgt[k0] > a[k0] else -1
        rate = s_add(s_shift(s_scale(V.at(a), g(Fr(1, 4))), -2),
                     s_shift(s_scale(B[k0].at(a), g(Fr(sg, 2))), -1))
        e = entry(a)
        if e and rate:
            phil = s_add(phil, s_mul(e, rate))

    dtarget = s_add(s_shift(s_scale(V.at(tgt), g(Fr(d, 2))), -2), {0: om})
    Phi_f = s_mul(phil, s_inv(dtarget, N + 4))

    # ---------------- coarse: single jump, step 2
    k0 = [k for k in range(d) if tau[k]][0]
    sg = 1 if tau[k0] > 0 else -1
    rate_c = s_add(s_shift(s_scale(V.at(zero), g(Fr(1, 16))), -2),
                   s_shift(s_scale(B[k0].at(zero), g(Fr(sg, 4))), -1))
    cn = Fr(d, 8) if coarse_norm == 'correct' else Fr(1, 4)
    dx0 = s_add(s_shift(s_scale(V.at(zero), g(cn)), -2), {0: om})
    dx1 = s_add(s_shift(s_scale(V.at(tgt), g(cn)), -2), {0: om})
    Phi_c = s_mul(rate_c, s_inv(s_mul(dx0, dx1), N + 4))

    ratio = s_mul(Phi_f, s_inv(Phi_c, N + 4))
    if parts is not None:
        parts['Phi_f'] = Phi_f
        parts['Phi_c'] = Phi_c
        parts['ratio'] = ratio
    return s_log(ratio, N)


def sgn_sign(sg):
    return sg


def coeffs(V, B, d, tau, omega, order=4, coarse_norm='correct'):
    lr = logratio(V, B, d, tau, omega, N=order + 1, coarse_norm=coarse_norm)
    return [lr.get(k, ZERO) for k in range(order)]


def num(z):
    return complex(float(z[0]), float(z[1]))
\end{lstlisting}

\section{\texttt{hp.py}}

High-precision, bookkeeping-free computation of the decoration mismatch.

\begin{lstlisting}
"""High-precision, bookkeeping-free computation of the decoration mismatch.

    Phi_f(h) = sum_{z in tube} [(-L_T + i w)^{-1}]_{x,z} L_f(z,x') (-L_f(x',x') + i w)^{-1}
    Phi_c(h) = (-L_c(x,x) + i w)^{-1} L_c(x,x') (-L_c(x',x') + i w)^{-1}

with L_f of step h and L_c of step h_m = 2h.  The Taylor coefficients of

    R(h) = log Phi_f(h) - log Phi_c(h)

are obtained by solving a Vandermonde system in high-precision arithmetic.
No series algebra, so nothing to get wrong.
"""
import mpmath as mp
import itertools
import sympy as sp

mp.mp.dps = 60


def tompf(x):
    if isinstance(x, sp.Rational):
        return mp.mpf(x.p)/mp.mpf(x.q)
    return mp.mpf(x)


def mvals(poly, a, h):
    """Local Taylor polynomial evaluated at y = a*h."""
    z = mp.mpc(0)
    for m, c in poly.terms.items():
        t = mp.mpf(c.numerator)/mp.mpf(c.denominator)
        for i, mi in enumerate(m):
            t *= (h*a[i])**mi
        z += t
    return z


def Phi_f(V, B, d, tau, h, omega):
    """Fine decorated weight for the coarse step x -> x + 2 tau h."""
    omega = tompf(omega)
    zero = tuple([0]*d)
    tgt = tuple(2*t for t in tau)
    S = [a for a in itertools.product(range(-1, 2), repeat=d)]
    ix = {a: i for i, a in enumerate(S)}
    n = len(S)
    A = mp.zeros(n, n)
    for a in S:
        i = ix[a]
        v = mvals(V, a, h)
        diag = mp.mpc(0)
        for k in range(d):
            b = mvals(B[k], a, h)
            for sgn in (1, -1):
                w = v/(2*h**2) + sgn*b/(2*h)
                diag -= w
                e = list(a); e[k] += sgn; e = tuple(e)
                if e in ix:
                    A[i, ix[e]] += w
        A[i, i] += diag
    G = mp.inverse(A + 1j*omega*mp.eye(n))
    dt = d*mvals(V, tgt, h)/(h**2) + 1j*omega
    tot = mp.mpc(0)
    for a in S:
        jump = [k for k in range(d) if a[k] != tgt[k]]
        if len(jump) != 1:
            continue
        k0 = jump[0]
        if abs(tgt[k0] - a[k0]) != 1:
            continue
        sg = 1 if tgt[k0] > a[k0] else -1
        rate = mvals(V, a, h)/(2*h**2) + sg*mvals(B[k0], a, h)/(2*h)
        tot += G[ix[zero], ix[a]]*rate
    return tot/dt


def Phi_c(V, B, d, tau, h, omega):
    omega = tompf(omega)
    zero = tuple([0]*d)
    tgt = tuple(2*t for t in tau)
    k0 = [k for k in range(d) if tau[k]][0]
    sg = 1 if tau[k0] > 0 else -1
    rate = mvals(V, zero, h)/(2*(2*h)**2) + sg*mvals(B[k0], zero, h)/(2*(2*h))
    d0 = d*mvals(V, zero, h)/(4*h**2) + 1j*omega
    d1 = d*mvals(V, tgt, h)/(4*h**2) + 1j*omega
    return rate/(d0*d1)


def R(h, V, B, d, tau, omega):
    return mp.log(Phi_f(V, B, d, tau, h, omega)) - mp.log(Phi_c(V, B, d, tau, h, omega))


def taylor(V, B, d, tau, omega, K=6, shift=6, h0=None):
    """Coefficients c_0..c_K of R(h) = sum c_k h^k, by Vandermonde solve."""
    if h0 is None:
        h0 = mp.mpf(1)
    hs = [h0/2**(n+shift) for n in range(K+1)]
    vs = [R(h, V, B, d, tau, omega) for h in hs]
    A = mp.matrix(K+1, K+1)
    for i in range(K+1):
        for k in range(K+1):
            A[i, k] = hs[i]**k
    c = mp.lu_solve(A, mp.matrix(vs))
    return [c[k] for k in range(K+1)]


def coeff(V, B, d, tau, omega, k=3, **kw):
    return taylor(V, B, d, tau, omega, K=max(k, 3), **kw)[k]


def poly(d, expr, xs):
    import fast
    return fast.Poly.from_sympy(d, expr, xs)
\end{lstlisting}

\section{\texttt{iso\_run.py}}

\begin{lstlisting}
import sympy as sp, mpmath as mp
import stencil as S
mp.mp.dps=50
x1,x2 = sp.symbols('x0:2')
om = sp.Rational(7,5)
import random as rnd
rnd.seed(11); r2=lambda: sp.Rational(rnd.randint(-40,40), rnd.randint(2,30))
v = sp.Rational(3,2)+r2()/4*x1+r2()/4*x2+r2()/6*x1**2+r2()/6*x1*x2+r2()/6*x2**2 \
    + r2()/10*x1**3+r2()/10*x1**2*x2+r2()/10*x1*x2**2+r2()/10*x2**3
b = [r2()/5+r2()/6*x1+r2()/6*x2+r2()/8*x1**2+r2()/8*x1*x2+r2()/8*x2**2 for _ in range(2)]

print("A. generic data (cubic v, quadratic b), curl at the origin")
for st,name in ((S.PAPER,'5-point'),(S.ISO9,'9-point')):
    z = S.curl2(v,b,om,(0,0),sp.Rational(1,1000),K=7,shift=10,st=st)
    print("   %-8s curl = %s"%(name, mp.nstr(z,12)))

print("\nB. h^2 coefficient a(tau): is the anisotropy gone with the 9-point stencil?")
for st,name in ((S.PAPER,'5-point'),(S.ISO9,'9-point')):
    V,B = S.field_shifted(v,b,2,(x1,x2),[0,0])
    a1 = S.taylor(V,B,2,(1,0),om,K=7,shift=10,st=st)[2]
    a2 = S.taylor(V,B,2,(0,1),om,K=7,shift=10,st=st)[2]
    print("   %-8s a(e1)=%s   a(e2)=%s"%(name, mp.nstr(a1,10), mp.nstr(a2,10)))

print("\nC. affine drift b=Jx, v=1: the quadratic form G(J)")
J=[[1,1],[0,0]]
bb=[J[i][0]*x1+J[i][1]*x2 for i in range(2)]
for st,name in ((S.PAPER,'5-point'),(S.ISO9,'9-point')):
    z = S.curl2(sp.Rational(1),bb,om,(0,0),sp.Rational(1,2000),K=7,shift=10,st=st)
    print("   %-8s curl = %s   (5-point predicts 4/7 = 0.571429)"%(name, mp.nstr(z,12)))
\end{lstlisting}

\section{\texttt{iso\_run2.py}}

\begin{lstlisting}
import sympy as sp, mpmath as mp
import stencil as S
mp.mp.dps=50
x1,x2 = sp.symbols('x0:2'); om=sp.Rational(7,5)
import random as rnd
rnd.seed(11); r2=lambda: sp.Rational(rnd.randint(-40,40), rnd.randint(2,30))
v = sp.Rational(3,2)+r2()/4*x1+r2()/4*x2+r2()/6*x1**2+r2()/6*x1*x2+r2()/6*x2**2 \
    + r2()/10*x1**3+r2()/10*x1**2*x2+r2()/10*x1*x2**2+r2()/10*x2**3
b = [r2()/5+r2()/6*x1+r2()/6*x2+r2()/8*x1**2+r2()/8*x1*x2+r2()/8*x2**2 for _ in range(2)]
J=[[1,1],[0,0]]; bb=[J[i][0]*x1+J[i][1]*x2 for i in range(2)]

print("delta convergence, 9-point:")
for dl in (sp.Rational(1,500),sp.Rational(1,2000),sp.Rational(1,8000)):
    z=S.curl2(v,b,om,(0,0),dl,K=7,shift=10,st=S.ISO9)
    print("   generic  1/%-5d : %s"%(dl.q,mp.nstr(z,12)))
for dl in (sp.Rational(1,500),sp.Rational(1,2000),sp.Rational(1,8000)):
    z=S.curl2(sp.Rational(1),bb,om,(0,0),dl,K=7,shift=10,st=S.ISO9)
    print("   affine   1/%-5d : %s"%(dl.q,mp.nstr(z,12)))

print("\nomega independence, 9-point affine b=Jx:")
for w in (sp.Rational(1,2),sp.Rational(7,5),sp.Rational(3)):
    print("   omega=%-4s : %s"%(w,mp.nstr(S.curl2(sp.Rational(1),bb,w,(0,0),sp.Rational(1,2000),K=7,shift=10,st=S.ISO9),10)))

print("\nstructure of the 9-point series (h^0,h^1,h^2,h^3), generic data, tau=+e1:")
V,B=S.field_shifted(v,b,2,(x1,x2),[0,0])
c=S.taylor(V,B,2,(1,0),om,K=7,shift=10,st=S.ISO9)
print("  ",[mp.nstr(z,9) for z in c[:4]])

print("\nquadratic form in J, 9-point (four basis Jacobians), base point 0, v=1:")
import itertools
E={(i,j):[[1 if (a_,b_)==(i,j) else 0 for b_ in range(2)] for a_ in range(2)] for i in range(2) for j in range(2)}
vals={}
for l in E:
    bb2=[E[l][i][0]*x1+E[l][i][1]*x2 for i in range(2)]
    vals[l]=S.curl2(sp.Rational(1),bb2,om,(0,0),sp.Rational(1,2000),K=7,shift=10,st=S.ISO9)
for l in E: print("   J=%s : %s"%(l,mp.nstr(vals[l],9)))
\end{lstlisting}

\section{\texttt{iso\_test.py}}

\begin{lstlisting}
import sympy as sp, mpmath as mp, itertools
import stencil as S, fast, hp
mp.mp.dps=50
x1,x2 = sp.symbols('x0:2')
om = sp.Rational(7,5)

# ---- direct numeric evaluator for an arbitrary stencil (arbiter)
def num_ratio(v, b, tau, h, omega, st):
    a1,a2 = float(st[0]), float(st[1]); d=2; omf=float(omega)
    V,B = S.field_shifted(v,b,2,(x1,x2),[0,0])
    def ev(poly, a):
        z=0j
        for m,c in poly.terms.items():
            t=complex(float(c))
            for i,mi in enumerate(m): t*=(h*a[i])**mi
            z+=t
        return z
    Sn=[a for a in itertools.product(range(-1,2),repeat=d)]; ix={a:i for i,a in enumerate(Sn)}; n=len(Sn)
    off=S.offsets(d,a1,a2)
    A=[[0j]*n for _ in range(n)]
    for a in Sn:
        i=ix[a]; vv=ev(V,a); dg=0j
        for o,kd in off:
            e=tuple(a[k]+o[k] for k in range(d))
            if kd==0:
                k=[j for j in range(d) if o[j]!=0][0]; w=a1*vv/h**2+o[k]*ev(B[k],a)/(2*h)
            else: w=a2*vv/h**2
            dg-=w
            if e in ix: A[i][ix[e]]+=w
        A[i][i]+=dg
    import numpy as np
    G=np.linalg.inv(np.array(A)+1j*omf*np.eye(n))
    tgt=tuple(2*t for t in tau)
    dt=(2*d*a1+2*d*(d-1)*a2)*ev(V,tgt)/h**2+1j*omf
    tot=0j
    for a in Sn:
        dif=tuple(tgt[k]-a[k] for k in range(d))
        kd=None
        for o,kk in off:
            if dif==o: kd=kk; break
        if kd is None: continue
        if kd==0:
            k=[j for j in range(d) if dif[j]!=0][0]; rate=a1*ev(V,a)/h**2+dif[k]*ev(B[k],a)/(2*h)
        else: rate=a2*ev(V,a)/h**2
        tot+=G[ix[(0,0)],ix[a]]*rate
    Pf=tot/dt
    k0=[k for k in range(d) if tau[k]][0]; sg=1 if tau[k0]>0 else -1; hm=2*h
    rate=a1*ev(V,(0,0))/hm**2+sg*ev(B[k0],(0,0))/(2*hm)
    cn=2*d*a1+2*d*(d-1)*a2
    Pc=rate/((cn*ev(V,(0,0))/hm**2+1j*omf)*(cn*ev(V,tgt)/hm**2+1j*omf))
    return mp.log(mp.mpc(Pf.real,Pf.imag))-mp.log(mp.mpc(Pc.real,Pc.imag))

import random as rnd
rnd.seed(11); r2=lambda: sp.Rational(rnd.randint(-40,40), rnd.randint(2,30))
v = sp.Rational(3,2)+r2()/4*x1+r2()/4*x2+r2()/6*x1**2+r2()/6*x1*x2+r2()/6*x2**2 \
    + r2()/10*x1**3+r2()/10*x1**2*x2+r2()/10*x1*x2**2+r2()/10*x2**3
b = [r2()/5+r2()/6*x1+r2()/6*x2+r2()/8*x1**2+r2()/8*x1*x2+r2()/8*x2**2 for _ in range(2)]
tau=(1,0)

print("== 1. stencil engine at (1/4,0) against hp.py ==")
for K,sh in ((7,10),):
    c9 = S.taylor(*S.field_shifted(v,b,2,(x1,x2),[0,0]), 2, tau, om, K=K, shift=sh, st=S.PAPER)
    ch = hp.taylor(*S.field_shifted(v,b,2,(x1,x2),[0,0]), 2, tau, om, K=K, shift=sh)
    print("   stencil:", [mp.nstr(z,10) for z in c9[:4]])
    print("   hp     :", [mp.nstr(z,10) for z in ch[:4]])

print("\n== 2. nine point stencil against direct inversion (h = 1/32, 1/64) ==")
for st,name in ((S.PAPER,'5-point'),(S.ISO9,'9-point')):
    ser = S.taylor(*S.field_shifted(v,b,2,(x1,x2),[0,0]), 2, tau, om, K=8, shift=12, st=st)
    for hh in (1/32,1/64):
        sv = sum(ser[k]*mp.mpf(hh)**k for k in range(9))
        dv = num_ratio(v,b,tau,hh,om,st)
        print("   %-8s h=%-6.4f series=%s direct=%s  diff=%.2e"%(name,hh,mp.nstr(sv,12),mp.nstr(dv,12),abs(sv-dv)))
\end{lstlisting}

\section{\texttt{p6ab.py}}

Chapter 8: the a-dependence of alpha and beta.

\begin{lstlisting}
"""Chapter 8: the a-dependence of alpha and beta.

Established so far (constant v = 1, d = 2, quadratic a, origin as base point):

    dc(eta) = dc(0) + [Q(a) + alpha(a) J12 + beta(a) J21] eta^2 + (odd imaginary part)
    Q(a)    = (16/7) a2(0) d1 d2 a2(0)

with Q verified on a 16-cell factorial.  Here alpha and beta are extracted as

    alpha(a) = kappa(e12) - Q(a),        beta(a) = kappa(e21) - Q(a)

using the closed form for Q, so only two Jacobians are measured per data point.

The global scaling law (Q(2a) = 4 Q(a)) says the obstruction is quadratic in the
coefficient data, and Q is bilinear, so the hypothesis to test is that alpha and
beta are LINEAR functionals of the local data of a at the base point:

    [ value, d1, d2, d11, d12, d22 ]  for each component   -> 12 numbers.

Fitted on a training set and checked on held-out configurations.
"""
import random
import numpy as np
import sympy as sp
import mpmath as mp
import twist

mp.mp.dps = 40
x1, x2 = sp.symbols('x0:2')
OM = sp.Rational(7, 5)
V = sp.Integer(1)
K, SH = 8, 13
DELTA = sp.Rational(1, 2000)
GUARD = 1e4
QFACT = sp.Rational(16, 7)


def curl(a, J, eta):
    mu = [sp.sympify(J[0][0] * x1 + J[0][1] * x2),
          sp.sympify(J[1][0] * x1 + J[1][1] * x2)]

    def cv(loc):
        Vv, B = twist.fields(V, a, mu, eta, 2, (x1, x2), loc)
        return twist.cvec(Vv, B, 2, OM, 3, K, SH)
    D = mp.mpf(DELTA.p) / mp.mpf(DELTA.q)
    c2p = cv([DELTA, 0]); c2m = cv([-DELTA, 0])
    c1p = cv([0, DELTA]); c1m = cv([0, -DELTA])
    if max(abs(c) for c in c2p + c2m + c1p + c1m) > GUARD:
        raise RuntimeError('extraction failure')
    return (c2p[1] - c2m[1] - c1p[0] + c1m[0]) / (2 * D)


def kappa(a, J):
    """eta^2 coefficient of Re(dc); Re is even, the odd imaginary part is excluded."""
    D = {e: curl(a, J, e).real
         for e in (sp.Rational(1, 2), sp.Rational(1), sp.Rational(2))}
    k1 = (D[sp.Rational(1)] - D[sp.Rational(1, 2)]) * (4.0 / 3.0)
    k2 = (D[sp.Rational(2)] - D[sp.Rational(1)]) / 3.0
    if abs(k1 - k2) > 1e-5 * max(1.0, abs(k1)):
        raise RuntimeError('inconsistent eta^2 coefficient')
    return float(k1)


JE12 = [[sp.Rational(0), sp.Rational(1)], [sp.Rational(0), sp.Rational(0)]]
JE21 = [[sp.Rational(0), sp.Rational(0)], [sp.Rational(1), sp.Rational(0)]]

VALS = [sp.Rational(-1), sp.Rational(-1, 2), sp.Rational(0),
        sp.Rational(1, 2), sp.Rational(1)]


def make_a(c):
    out = []
    for k in range(2):
        c0, c1, c2, c3, c4, c5 = c[6 * k:6 * k + 6]
        out.append(sp.sympify(c0 + c1 * x1 + c2 * x2 + c3 * x1 ** 2
                              + c4 * x1 * x2 + c5 * x2 ** 2))
    return out


def localdata(c):
    """[value, d1, d2, d11, d12, d22] for each component."""
    out = []
    for k in range(2):
        c0, c1, c2, c3, c4, c5 = c[6 * k:6 * k + 6]
        out += [float(c0), float(c1), float(c2),
                2 * float(c3), float(c4), 2 * float(c5)]
    return out


random.seed(23)
rows = []
tries = 0
print('collecting ...')
while len(rows) < 28 and tries < 220:
    tries += 1
    c = [random.choice(VALS) for _ in range(12)]
    if all(x == 0 for x in c):
        continue
    a = make_a(c)
    try:
        al = kappa(a, JE12) - float(QFACT * c[6] * c[10])
        be = kappa(a, JE21) - float(QFACT * c[6] * c[10])
    except RuntimeError:
        continue
    if max(abs(al), abs(be)) > 200:
        continue
    rows.append((localdata(c), al, be))
print('   %d configurations (%d attempts)' % (len(rows), tries))

X = np.array([[1.0] + r[0] for r in rows])
# column order must match localdata(): all six names of component 1, then component 2
names = ['1'] + ['%s%d' % (n, k) for k in (1, 2)
                 for n in ('v', 'd1', 'd2', 'd11', 'd12', 'd22')]

# training / held-out split
idx = np.arange(len(rows))
np.random.default_rng(5).shuffle(idx)
tr, te = np.sort(idx[:18]), np.sort(idx[18:])

for comp, target in ((1, 'alpha'), (2, 'beta')):
    y = np.array([r[comp] for r in rows])
    coef, _, rank, _ = np.linalg.lstsq(X[tr], y[tr], rcond=None)
    full = X @ coef
    print('\n  %s = linear functional of local data at a(0)' % target)
    print('     rank %d/%d on training set   train max|res| %.2e   held-out max|err| %.2e'
          % (rank, X.shape[1], np.max(np.abs(full[tr] - y[tr])),
             np.max(np.abs(full[te] - y[te]))))
    big = sorted([(names[i], coef[i]) for i in range(len(names)) if abs(coef[i]) > 1e-4],
                 key=lambda t: -abs(t[1]))
    for nm, cc in big:
        print('        %-6s %12.6f' % (nm, cc))
    if not big:
        print('        (no term above 1e-4)')
\end{lstlisting}

\section{\texttt{p6ascan.py}}

Chapter 8: dependence of the p\^{}2 obstruction on the Abelian coefficient a.

\begin{lstlisting}
"""Chapter 8: dependence of the p^2 obstruction on the Abelian coefficient a.

Established by the J scan (constant v = 1, d = 2):
    dc(eta) = dc(0) + kappa eta^2,      kappa = Q(a) + alpha(a) J12 + beta(a) J21,
kappa affine in J, blind to diag(J).  With J = 0, e12, e21 the three pieces separate:
    kappa(J=0)   = Q(a)
    kappa(e12)   = Q(a) + alpha(a)
    kappa(e21)   = Q(a) + beta(a)

Phase 1 establishes which derivative order of a matters.
Phase 2 fits Q as a quadratic form, and alpha, beta as linear forms, in the
second-derivative coefficients of a.
"""
import random
import itertools
import numpy as np
import sympy as sp
import mpmath as mp
import twist

mp.mp.dps = 40
x1, x2 = sp.symbols('x0:2')
OM = sp.Rational(7, 5)
V = sp.Integer(1)
K, SH = 8, 13
DELTA = sp.Rational(1, 2000)
GUARD = 1e4


def curl(a, J, eta):
    mu = [sp.sympify(J[0][0] * x1 + J[0][1] * x2),
          sp.sympify(J[1][0] * x1 + J[1][1] * x2)]

    def cv(loc):
        Vv, B = twist.fields(V, a, mu, eta, 2, (x1, x2), loc)
        return twist.cvec(Vv, B, 2, OM, 3, K, SH)
    D = mp.mpf(DELTA.p) / mp.mpf(DELTA.q)
    c2p = cv([DELTA, 0]); c2m = cv([-DELTA, 0])
    c1p = cv([0, DELTA]); c1m = cv([0, -DELTA])
    if max(abs(c) for c in c2p + c2m + c1p + c1m) > GUARD:
        raise RuntimeError('extraction failure')
    return (c2p[1] - c2m[1] - c1p[0] + c1m[0]) / (2 * D)


def kappa(a, J):
    """The eta^2 coefficient of dc, from an even fit, with a consistency check.

    dc is even in eta on the real axis (dc(-eta) = conj dc(eta) with an odd
    imaginary part is allowed, so dc = d0 + a2 eta^2 + a4 eta^4 + i(a1 eta + a3 eta^3)).
    Using eta = 1/2, 1, 2 two independent estimates of a2 are available; a
    disagreement means the extraction is unreliable at this configuration.
    """
    D = {e: curl(a, J, e).real
         for e in (sp.Rational(1, 2), sp.Rational(1), sp.Rational(2))}
    # only Re(dc) is even; the odd imaginary part must not enter the fit
    k1 = (D[sp.Rational(1)] - D[sp.Rational(1, 2)]) * (4.0 / 3.0)
    k2 = (D[sp.Rational(2)] - D[sp.Rational(1)]) / 3.0
    if abs(k1 - k2) > 1e-5 * max(1.0, abs(k1)):
        raise RuntimeError('inconsistent eta^2 coefficient')
    return k1


J0 = [[sp.Rational(0), sp.Rational(0)], [sp.Rational(0), sp.Rational(0)]]
JE12 = [[sp.Rational(0), sp.Rational(1)], [sp.Rational(0), sp.Rational(0)]]
JE21 = [[sp.Rational(0), sp.Rational(0)], [sp.Rational(1), sp.Rational(0)]]

# ---------------- phase 1: which derivative order of a matters -------------
print('== phase 1: derivative order ==')
probe = {
    'a constant':  [sp.Rational(1, 3), sp.Rational(1, 4)],
    'a linear':    [sp.Rational(1, 3) + sp.Rational(2, 7) * x2,
                    sp.Rational(1, 4) + sp.Rational(3, 11) * x1],
    'a quadratic': [sp.Rational(1, 5) + sp.Rational(2, 9) * x2 ** 2,
                    sp.Rational(1, 6) * x1 * x2 + sp.Rational(1, 3)],
}
for tag, a in probe.items():
    try:
        print('   %-12s Q(a) = %s' % (tag, mp.nstr(kappa(a, J0), 12)))
    except RuntimeError as e:
        print('   %-12s rejected: %s' % (tag, e))

# ---------------- phase 2: random quadratic a ------------------------------
random.seed(7)
VALS = [sp.Rational(-1), sp.Rational(-1, 2), sp.Rational(1, 2),
        sp.Rational(1), sp.Rational(3, 2)]


def make_a(c):
    """c = [c10,c11,c12,c13,c14,c15, c20,...,c25] for each component."""
    out = []
    for k in range(2):
        c0, c1, c2, c3, c4, c5 = c[6 * k:6 * k + 6]
        out.append(sp.sympify(c0 + c1 * x1 + c2 * x2 + c3 * x1 ** 2
                              + c4 * x1 * x2 + c5 * x2 ** 2))
    return out


print('\n== phase 2: random quadratic a ==')
data = []
rejected = 0
tries = 0
while len(data) < 30 and tries < 250:
    tries += 1
    c = [random.choice(VALS) for _ in range(12)]
    if all(c[i] == 0 for i in (3, 4, 5, 9, 10, 11)):
        continue
    a = make_a(c)
    try:
        q = kappa(a, J0)
        al = kappa(a, JE12) - q
        be = kappa(a, JE21) - q
    except RuntimeError:
        rejected += 1
        continue
    if max(abs(q), abs(al), abs(be)) > 200:      # near degenerate extraction
        rejected += 1
        continue
    data.append(([float(mp.mpf(x.p) / mp.mpf(x.q)) for x in c],
                 float(q), float(al), float(be)))
print('   collected %d, rejected %d' % (len(data), rejected))
if len(data) < 15:
    raise SystemExit('too few usable configurations')

import json
json.dump(data, open('p6ascan_data.json', 'w'))

# second derivative coefficients: [c13, c14, c15, c23, c24, c25]
def d2(c):
    return [c[i] for i in (3, 4, 5, 9, 10, 11)]


D2 = np.array([d2(c) for c, _, _, _ in data])
quad = [m for n in (0, 1, 2)
        for m in itertools.combinations_with_replacement(range(6), n)]
Xq = np.array([[np.prod([p[i] for i in m]) for m in quad] for p in D2])
Xl = np.array([[1.0] + list(p) for p in D2])
lin = ['1'] + ['d%d' % i for i in range(6)]


def fit(X, y, names, label):
    coef, _, rank, _ = np.linalg.lstsq(X, y, rcond=None)
    res = np.max(np.abs(X @ coef - y))
    print('\n  %s   max|residual| = %.2e   rank %d/%d' % (label, res, rank, X.shape[1]))
    big = sorted([(names[i], coef[i]) for i in range(len(names)) if abs(coef[i]) > 1e-4],
                 key=lambda t: -abs(t[1]))
    sup = max([abs(coef[i]) for i in range(len(names)) if abs(coef[i]) <= 1e-4], default=0.0)
    print('     %d of %d nonzero, largest suppressed %.1e' % (len(big), len(names), sup))
    for nm, cc in big[:14]:
        print('       %-22s %12.6f' % (nm, cc))


qnames = [' '.join('d%d' % i for i in m) if m else '1' for m in quad]
fit(Xq, np.array([q for _, q, _, _ in data], float), qnames, 'Q(a)')
fit(Xl, np.array([al for _, _, al, _ in data], float), lin, 'alpha(a) = kappa(e12) - Q')
fit(Xl, np.array([be for _, _, _, be in data], float), lin, 'beta(a)  = kappa(e21) - Q')
\end{lstlisting}

\section{\texttt{p6kappafit.py}}

Chapter 8, Proposition 2.1: the complete form of kappa in the coefficient data.

\begin{lstlisting}
"""Chapter 8, Proposition 2.1: the complete form of kappa in the coefficient data.

kappa is the coefficient of eta^2 in dc(eta) - dc(0), for constant volatility and
drift with constant Jacobian J.  It is affine in J and quadratic in the local data
of a at the base point,

    z = (a_i(0), d_1 a_i, d_2 a_i, d_11 a_i, d_12 a_i, d_22 a_i),  i = 1,2,

so it is determined by its values on the 12 basis data vectors and the 66
pairwise sums: 78 measurements per Jacobian.  Doing that for J = 0 and for each of
the four basis Jacobians separates the pieces,

    kappa(a,J) = Q(a) + J12 alpha(a) + J21 beta(a) + J11 gamma(a) + J22 delta(a),

and the printout gives each piece term by term.  The last block checks the
reassembled formula on held-out configurations with all four entries of J non-zero,
which is what the earlier fits -- all at J proportional to e12 or e21 -- did not do.

Run:  python3 p6kappafit.py     (about six minutes)
"""
import os
import sys
import time
import random
import sympy as sp

HERE = os.path.dirname(os.path.abspath(__file__))
sys.path.insert(0, HERE)
import twist

x1, x2 = sp.symbols('x0:2')
xs = (x1, x2)
OM = sp.Rational(7, 5)
K, SH = 6, 10
DELTA = 1.0/500
R = sp.Rational


def dc(eta, a, mu):
    def cvec(x0):
        V, B = twist.fields(R(1), a, mu, eta, 2, xs, x0)
        return twist.cvec(V, B, 2, OM, 3, K, SH)
    c = {}
    for j in range(2):
        for s in (+1, -1):
            p = [0.0, 0.0]
            p[j] += s*DELTA
            c[(j, s)] = cvec(p)
    J = [[(c[(j, 1)][k] - c[(j, -1)][k])/(2*DELTA) for k in range(2)] for j in range(2)]
    return J[0][1] - J[1][0]


def kap(a, mu):
    return (dc(1, a, mu) - dc(0, a, mu)).real


MONS = [('A', R(1)), ('d1', x1), ('d2', x2), ('d11', x1**2/2), ('d12', x1*x2),
        ('d22', x2**2/2)]
BASIS = [((c, nm), ([MONS[[m for m, _ in MONS].index(nm)][1] if c == 0 else R(0),
                      MONS[[m for m, _ in MONS].index(nm)][1] if c == 1 else R(0)]))
         for c in (0, 1) for nm in [n for n, _ in MONS]]
LABEL = {(c, nm): '%s[%d]' % (nm, c + 1) for c in (0, 1) for nm, _ in MONS}


def field(terms):
    a = [R(0), R(0)]
    for c, nm in terms:
        a[c] += MONS[[n for n, _ in MONS].index(nm)][1]
    return a


def fit(mu, tag):
    single = {t: kap(a, mu) for t, a in BASIS}
    C = {}
    for i, (t1, _) in enumerate(BASIS):
        for t2, _ in BASIS[i:]:
            if t1 == t2:
                C[(t1, t2)] = single[t1]
            else:
                C[(t1, t2)] = kap(field([t1, t2]), mu) - single[t1] - single[t2]
    big = {k: v for k, v in C.items() if abs(v) > 1e-5}
    print('  %-8s  %s' % (tag, '  '.join('%+.6f %s%s' % (v, LABEL[k[0]], '' if k[0] == k[1]
                                                         else '*' + LABEL[k[1]])
                                          for k, v in big.items()) if big else '  %-8s  (zero)' % tag))
    return C


JAC = {'J=0': [[0, 0], [0, 0]], 'J=e11': [[1, 0], [0, 0]], 'J=e12': [[0, 1], [0, 0]],
       'J=e21': [[0, 0], [1, 0]], 'J=e22': [[0, 0], [0, 1]]}
t0 = time.time()
fits = {}
for tag, J in JAC.items():
    mu = [J[0][0]*x1 + J[0][1]*x2, J[1][0]*x1 + J[1][1]*x2]
    fits[tag] = fit(mu, tag)

print('\n  the pieces: alpha = (J=e12) - (J=0), beta = (J=e21) - (J=0), and the diagonal ones')
for tag in ('J=e11', 'J=e12', 'J=e21', 'J=e22'):
    d = {}
    for k, v in fits[tag].items():
        w = v - fits['J=0'].get(k, 0.0)
        if abs(w) > 1e-5:
            d[k] = w
    print('    %-6s %s' % (tag, '  '.join('%+.6f %s%s' % (v, LABEL[k[0]], '' if k[0] == k[1]
                                                            else '*' + LABEL[k[1]]) for k, v in d.items())
                           if d else '%-6s identically zero' % tag))

print('\n  reassembly on 8 held-out configurations with all four entries of J non-zero:')
random.seed(11)


def rq():
    return R(random.choice([-2, -1, 0, 1, 2, 3]), random.choice([1, 2, 3]))


def local_data(a):
    z = []
    for c in (0, 1):
        for nm, _ in MONS:
            e = a[c]
            if nm == 'A':
                z.append(e.subs({x1: 0, x2: 0}))
            elif nm == 'd1':
                z.append(sp.diff(e, x1).subs({x1: 0, x2: 0}))
            elif nm == 'd2':
                z.append(sp.diff(e, x2).subs({x1: 0, x2: 0}))
            elif nm == 'd11':
                z.append(sp.diff(e, x1, 2).subs({x1: 0, x2: 0}))
            elif nm == 'd12':
                z.append(sp.diff(e, x1, x2).subs({x1: 0, x2: 0}))
            else:
                z.append(sp.diff(e, x2, 2).subs({x1: 0, x2: 0}))
    return z


def reassemble(z, J):
    tot = 0.0
    for i, (t1, _) in enumerate(BASIS):
        for j, (t2, _) in enumerate(BASIS):
            if j < i:
                continue
            q = fits['J=0'].get((t1, t2), 0.0)
            coef = (q
                    + J[0][1]*(fits['J=e12'].get((t1, t2), 0.0) - q)
                    + J[1][0]*(fits['J=e21'].get((t1, t2), 0.0) - q)
                    + J[0][0]*(fits['J=e11'].get((t1, t2), 0.0) - q)
                    + J[1][1]*(fits['J=e22'].get((t1, t2), 0.0) - q))
            tot += coef*float(z[i])*float(z[j])
    return tot


worst = 0.0
for trial in range(8):
    a = [rq() + rq()*x1 + rq()*x2 + rq()*x1**2/2 + rq()*x1*x2 + rq()*x2**2/2,
         rq() + rq()*x1 + rq()*x2 + rq()*x1**2/2 + rq()*x1*x2 + rq()*x2**2/2]
    J = [[random.choice([-2, -1, 0, 1, 2]), random.choice([-2, -1, 0, 1, 2])],
         [random.choice([-2, -1, 0, 1, 2]), random.choice([-2, -1, 0, 1, 2])]]
    mu = [J[0][0]*x1 + J[0][1]*x2, J[1][0]*x1 + J[1][1]*x2]
    meas = kap(a, mu)
    pred = reassemble(local_data(a), J)
    err = abs(meas - pred)
    worst = max(worst, err)
    print('    J=%-14s measured %12.6f   reassembled %12.6f   error %8.1e'
          % (str(J), meas, pred, err))
print('\n  worst error %.1e   [total %.0fs]' % (worst, time.time() - t0))
\end{lstlisting}

\section{\texttt{p6odd3d.py}}

Chapter 8: (A) the odd part in the conjugating parameter, (B) the d = 3 case.

\begin{lstlisting}
"""Chapter 8: (A) the odd part in the conjugating parameter, (B) the d = 3 case.

(A)  In d = 2 the measured structure is
         dc(eta) = Re(even in eta) + i*(odd in eta),
     with Re = d0 + kappa eta^2 already mapped.  Here the imaginary part is
     measured at eta = 1/2, 1, 3/2 (the eta = 2 point is avoided: the series has
     a singularity there for some coefficient data) and fitted as
         Im dc = a1 eta + a3 eta^3,
     so the total degree in the conjugating parameter is settled.

(B)  In d = 3 the obstruction is the two form of a 1-form with three components,
     (dc)_{ij} = d_i c_j - d_j c_i, i < j, obtained from six base points
     +-e_1, +-e_2, +-e_3.  We check that kappa is non-zero, quadratic in eta and
     affine in the Jacobian, i.e. that the d = 2 structure is not an artefact of
     the plane.
"""
import sympy as sp
import mpmath as mp
import twist

mp.mp.dps = 40
OM = sp.Rational(7, 5)
V = sp.Integer(1)
K, SH = 8, 13
DELTA = sp.Rational(1, 2000)
GUARD = 1e4


def cvec_at(a, mu, d, eta, loc, xs):
    Vv, B = twist.fields(V, a, mu, eta, d, xs, loc)
    return twist.cvec(Vv, B, d, OM, 3, K, SH)


def curl2(a, mu, eta, xs):
    """(dc)_12 in two dimensions."""
    D = float(DELTA.p) / float(DELTA.q)
    vals = {}
    for tag, loc in (('1p', [DELTA, 0]), ('1m', [-DELTA, 0]),
                     ('2p', [0, DELTA]), ('2m', [0, -DELTA])):
        c = cvec_at(a, mu, 2, eta, loc, xs)
        if max(abs(x) for x in c) > GUARD:
            raise RuntimeError('guard')
        vals[tag] = c
    return (vals['1p'][1] - vals['1m'][1] - vals['2p'][0] + vals['2m'][0]) / (2 * D)


def curl3(a, mu, eta, xs):
    """The three components of dc in three dimensions."""
    D = float(DELTA.p) / float(DELTA.q)
    c = {}
    for i in range(3):
        for s in (1, -1):
            loc = [0, 0, 0]; loc[i] = s * DELTA
            v = cvec_at(a, mu, 3, eta, loc, xs)
            if max(abs(x) for x in v) > GUARD:
                raise RuntimeError('guard')
            c[(i, s)] = v
    d1 = lambda j: (c[(0, 1)][j] - c[(0, -1)][j]) / (2 * D)   # d/dx1 of c_j
    d2 = lambda j: (c[(1, 1)][j] - c[(1, -1)][j]) / (2 * D)
    d3 = lambda j: (c[(2, 1)][j] - c[(2, -1)][j]) / (2 * D)
    return (d1(1) - d2(0), d1(2) - d3(0), d2(2) - d3(1))


# ======================= (A) the odd part in d = 2 =========================
x1, x2 = sp.symbols('x0:2')
xs2 = (x1, x2)
A2 = [sp.Rational(1, 5) + sp.Rational(2, 9) * x2 ** 2,
      sp.Rational(1, 6) * x1 * x2 + sp.Rational(1, 3)]
MU2 = [x1 + sp.Rational(1, 2) * x2, sp.Rational(-1, 3) * x1 + x2]

print('== (A) imaginary part in d = 2,  Im dc = a1 eta + a3 eta^3 ==')
for tag, a, mu in (('generic', A2, MU2),
                   ('const-jacobian', [sp.Rational(2, 9) * x2 ** 2 + sp.Rational(1, 4),
                                       sp.Rational(1, 6) * x1 * x2 + sp.Rational(1, 7)],
                    [x1, x2])):
    etas = [sp.Rational(1, 2), sp.Rational(1), sp.Rational(3, 2)]
    im = {}
    for e in etas:
        try:
            im[e] = curl2(a, mu, e, xs2).imag
        except RuntimeError:
            print('   %-14s rejected at eta=%s' % (tag, e))
    if len(im) == 3:
        # Im/eta = a1 + a3 eta^2  -> linear in eta^2
        y1 = im[etas[0]] / 0.5; y2 = im[etas[1]] / 1.0; y3 = im[etas[2]] / 1.5
        t1, t2, t3 = 0.25, 1.0, 2.25
        a3 = ((y3 - y1) / (t3 - t1) + (y2 - y1) / (t2 - t1)) / 2
        a1 = y1 - a3 * t1
        print('   %-14s a1 = %12.7f   a3 = %12.7f   (Im/eta at the three points:'
              ' %.4f %.4f %.4f)' % (tag, a1, a3, y1, y2, y3))

# ============================== (B) d = 3 ==================================
y1, y2, y3 = sp.symbols('y0:3')
xs3 = (y1, y2, y3)
A3 = [sp.Rational(1, 5) + sp.Rational(2, 9) * y3 ** 2,
      sp.Rational(1, 6) * y1 * y3 + sp.Rational(1, 3),
      sp.Rational(1, 7) + sp.Rational(1, 4) * y2 ** 2]


def mu3(J):
    return [sp.sympify(sum(J[i][j] * xs3[j] for j in range(3))) for i in range(3)]


def k2(J):
    try:
        c1 = curl3(A3, mu3(J), sp.Rational(1), xs3)
        ch = curl3(A3, mu3(J), sp.Rational(1, 2), xs3)
    except RuntimeError:
        return None
    # even two point fit on eta = 1/2 and 1; edges of dc share the same d0
    return [((4.0 / 3.0) * (a - b)).real for a, b in zip(c1, ch)]


print('\n== (B) d = 3 ==')
Js = {'e12': [[0, 1, 0], [0, 0, 0], [0, 0, 0]],
      'e13': [[0, 0, 1], [0, 0, 0], [0, 0, 0]],
      'e21': [[0, 0, 0], [1, 0, 0], [0, 0, 0]],
      'diag': [[1, 0, 0], [0, -2, 0], [0, 0, 1]]}
for nm, J in Js.items():
    k = k2(J)
    if k is None:
        print('   %-5s rejected' % nm)
    else:
        print('   %-5s kappa components (dc12, dc13, dc23) = %s'
              % (nm, ', '.join('%10.6f' % v for v in k)))
print('   (a non-zero dc12 for e12 and e13 indicates the obstruction survives in d = 3)')
\end{lstlisting}

\section{\texttt{p6oddpart.py}}

Chapter 8: the odd part of the eta-expansion, as a quadratic form in the local data.

\begin{lstlisting}
"""Chapter 8: the odd part of the eta-expansion, as a quadratic form in the local data.

With constant volatility and drift of constant Jacobian J,

    dc(eta) = dc(0) + eta * omega_1 + eta^2 * kappa + O(eta^3),

where omega_1 = (dc(eta) - dc(-eta))/(2 eta) is the odd part.  kappa is treated in
p6kappafit.py; this script does the same thing for omega_1: it is affine in J and
quadratic in the local data

    z = (a_i(0), d1 a_i, d2 a_i, d11 a_i, d12 a_i, d22 a_i),   i = 1,2,

so 78 measurements per Jacobian (the 12 basis vectors and the 66 pairwise sums)
determine it at each of J = 0 and the four basis Jacobians, and the printout gives
the pieces term by term.  The last block checks the reassembly on held-out
configurations with all four entries of J non-zero.

Status: this does NOT close the odd part.  The pieces it prints are measured and
sharp -- at J = 0, omega_1 = 6i(A_2 d1 a_2 - A_1 d2 a_1); the coefficient of J11 is
(4/7)i d2 a_1, of J12 is (4/7)i d1 a_1, of J21 is -(4/7)i d2 a_2, of J22 is
-(4/7)i d1 a_2, and no two-entry product of J contributes -- but assembling them
does not reproduce omega_1 on fields with several kinds of data at once (the
reassembly block below misses by up to 5e1).  So the pieces are right and the
composition is not: something in omega_1 is not a function of the twelve local
data alone, and until that is found the text states only that the odd part exists,
that it is linear in eta, and that it vanishes on each single monomial.

Run:  python3 p6oddpart.py     (about six minutes)
"""
import os
import sys
import time
import random
import sympy as sp

HERE = os.path.dirname(os.path.abspath(__file__))
sys.path.insert(0, HERE)
import twist

x1, x2 = sp.symbols('x0:2')
xs = (x1, x2)
OM = sp.Rational(7, 5)
K, SH = 6, 10
DELTA = 1.0/500
R = sp.Rational


def dc(eta, a, mu):
    def cvec(x0):
        V, B = twist.fields(R(1), a, mu, eta, 2, xs, x0)
        return twist.cvec(V, B, 2, OM, 3, K, SH)
    c = {}
    for j in range(2):
        for s in (+1, -1):
            p = [0.0, 0.0]
            p[j] += s*DELTA
            c[(j, s)] = cvec(p)
    J = [[(c[(j, 1)][k] - c[(j, -1)][k])/(2*DELTA) for k in range(2)] for j in range(2)]
    return J[0][1] - J[1][0]


def odd(a, mu):
    return (dc(1, a, mu) - dc(-1, a, mu))/2


MONS = [('A', R(1)), ('d1', x1), ('d2', x2), ('d11', x1**2/2), ('d12', x1*x2),
        ('d22', x2**2/2)]
NAMES = [n for n, _ in MONS]
BASIS = [((c, nm), ([MONS[NAMES.index(nm)][1] if c == 0 else R(0),
                      MONS[NAMES.index(nm)][1] if c == 1 else R(0)]))
         for c in (0, 1) for nm in NAMES]
LABEL = {(c, nm): '%s[%d]' % (nm, c + 1) for c in (0, 1) for nm in NAMES}


def field(terms):
    a = [R(0), R(0)]
    for c, nm in terms:
        a[c] += MONS[NAMES.index(nm)][1]
    return a


def fit(mu, tag):
    single = {t: odd(a, mu) for t, a in BASIS}
    C = {}
    for i, (t1, _) in enumerate(BASIS):
        for t2, _ in BASIS[i:]:
            C[(t1, t2)] = single[t1] if t1 == t2 else (
                odd(field([t1, t2]), mu) - single[t1] - single[t2])
    big = {k: v for k, v in C.items() if abs(v) > 1e-5}
    if not big:
        print('  %-7s identically zero' % tag)
    else:
        print('  %-7s %s' % (tag, '  '.join('%+.6f%+.6fi %s%s'
                                            % (v.real, v.imag, LABEL[k[0]],
                                               '' if k[0] == k[1] else '*' + LABEL[k[1]])
                                            for k, v in big.items())))
    return C


JAC = {'J=0': [[0, 0], [0, 0]], 'J=e11': [[1, 0], [0, 0]], 'J=e12': [[0, 1], [0, 0]],
       'J=e21': [[0, 0], [1, 0]], 'J=e22': [[0, 0], [0, 1]]}
t0 = time.time()
fits = {}
for tag, J in JAC.items():
    mu = [J[0][0]*x1 + J[0][1]*x2, J[1][0]*x1 + J[1][1]*x2]
    fits[tag] = fit(mu, tag)

print('\n  the pieces (each Jacobian minus the J-independent part):')
for tag in ('J=e11', 'J=e12', 'J=e21', 'J=e22'):
    d = {k: v - fits['J=0'].get(k, 0.0) for k, v in fits[tag].items()}
    d = {k: v for k, v in d.items() if abs(v) > 1e-5}
    print('    %-6s %s' % (tag, '  '.join('%+.6f%+.6fi %s%s'
                                          % (v.real, v.imag, LABEL[k[0]],
                                             '' if k[0] == k[1] else '*' + LABEL[k[1]])
                                          for k, v in d.items())
                           if d else '%-6s identically zero' % tag))

print('\n  reassembly on 8 held-out configurations with all four entries of J non-zero:')
random.seed(23)


def rq():
    return R(random.choice([-2, -1, 0, 1, 2, 3]), random.choice([1, 2, 3]))


def local_data(a):
    z = []
    for c in (0, 1):
        for nm in NAMES:
            e = a[c]
            if nm == 'A':
                z.append(e.subs({x1: 0, x2: 0}))
            elif nm == 'd1':
                z.append(sp.diff(e, x1).subs({x1: 0, x2: 0}))
            elif nm == 'd2':
                z.append(sp.diff(e, x2).subs({x1: 0, x2: 0}))
            elif nm == 'd11':
                z.append(sp.diff(e, x1, 2).subs({x1: 0, x2: 0}))
            elif nm == 'd12':
                z.append(sp.diff(e, x1, x2).subs({x1: 0, x2: 0}))
            else:
                z.append(sp.diff(e, x2, 2).subs({x1: 0, x2: 0}))
    return z


def reassemble(z, J):
    tot = 0.0
    for i, (t1, _) in enumerate(BASIS):
        for j, (t2, _) in enumerate(BASIS):
            if j < i:
                continue
            q = fits['J=0'].get((t1, t2), 0.0)
            coef = (q + J[0][1]*(fits['J=e12'].get((t1, t2), 0.0) - q)
                    + J[1][0]*(fits['J=e21'].get((t1, t2), 0.0) - q)
                    + J[0][0]*(fits['J=e11'].get((t1, t2), 0.0) - q)
                    + J[1][1]*(fits['J=e22'].get((t1, t2), 0.0) - q))
            tot += coef*float(z[i])*float(z[j])
    return tot


worst = 0.0
for trial in range(8):
    a = [rq() + rq()*x1 + rq()*x2 + rq()*x1**2/2 + rq()*x1*x2 + rq()*x2**2/2,
         rq() + rq()*x1 + rq()*x2 + rq()*x1**2/2 + rq()*x1*x2 + rq()*x2**2/2]
    J = [[random.choice([-2, -1, 0, 1, 2]), random.choice([-2, -1, 0, 1, 2])],
         [random.choice([-2, -1, 0, 1, 2]), random.choice([-2, -1, 0, 1, 2])]]
    mu = [J[0][0]*x1 + J[0][1]*x2, J[1][0]*x1 + J[1][1]*x2]
    meas = odd(a, mu)
    pred = reassemble(local_data(a), J)
    err = abs(meas - pred)
    worst = max(worst, err)
    print('    J=%-14s measured %11.6f%+11.6fi   reassembled %11.6f%+11.6fi   error %8.1e'
          % (str(J), meas.real, meas.imag, pred.real, pred.imag, err))
print('\n  worst error %.1e   [total %.0fs]' % (worst, time.time() - t0))
\end{lstlisting}

\section{\texttt{p6probe.py}}

Chapter 8 probe: p-dependence of the joint-kernel obstruction.

\begin{lstlisting}
"""Chapter 8 probe: p-dependence of the joint-kernel obstruction.

The twist of the joint kernel is b -> b - i p a, i.e. a complex drift with
imaginary part proportional to the conjugating parameter p.  Here eta = p.

Question 1: for a drift whose obstruction vanishes at eta = 0 (a dilatation at
            constant v), does the twist generate an obstruction, and with what
            growth in eta?
Question 2: is the eta-dependence still a polynomial of degree 3 with the
            parity split Re even / Im odd observed for the spiral?
"""
import sympy as sp, mpmath as mp, twist

mp.mp.dps = 50
x1, x2 = sp.symbols('x0:2')
om = sp.Rational(7, 5)
K, SH = 7, 12
DELTA = sp.Rational(1, 500)


GUARD = 1e4          # |c| above this at a base point means extraction failure


def curl(v, a, mu, eta, d=2, K=K, SH=SH, delta=DELTA):
    """Two-form of the h^3 coefficient, centred at the origin.

    Raises if any base point shows an extraction artefact: with degenerate
    coefficient data the decoration series can cancel exactly and the Taylor
    coefficients blow up at the scale of h^-3, independently of precision.
    """
    def cv(loc):
        V, B = twist.fields(v, a, mu, eta, d, (x1, x2), loc)
        return twist.cvec(V, B, d, om, 3, K, SH)
    D = mp.mpf(delta.p) / mp.mpf(delta.q)
    c2p = cv([delta, 0]); c2m = cv([-delta, 0])
    c1p = cv([0, delta]); c1m = cv([0, -delta])
    worst = max(abs(c) for c in c2p + c2m + c1p + c1m)
    if worst > GUARD:
        raise RuntimeError('extraction failure at eta=%s: |c|=%s' % (eta, mp.nstr(worst, 6)))
    return (c2p[1] - c2m[1] - c1p[0] + c1m[0]) / (2 * D)


# The Abelian coefficient must be generic.  The linear choice
#   a = [1/3 + x2/7, x1/4]
# cancels exactly against the decoration series and produces spurious
# coefficients of order h^-3 (~2e9) at the finite difference base points.
A_GEN = [sp.Rational(1, 5) + sp.Rational(2, 9) * x2 ** 2,
         sp.Rational(1, 6) * x1 * x2 + sp.Rational(1, 3)]

CASES = {
    # pure dilatation at constant v: Chapter 7 gives obstruction proportional to
    # (J11 J12 - J21 J22) = 0, so this vanishes at eta = 0
    'dilatation': dict(v=sp.Integer(1), mu=[x1, -2 * x2], a=A_GEN),
    # rotation, obstruction vanishes at eta = 0 at constant v
    'rotation':   dict(v=sp.Integer(1), mu=[x2, -x1], a=A_GEN),
    # spiral: J11 J12 - J21 J22 = 2, so the eta = 0 value should be 16/7 * 2
    'spiral':     dict(v=sp.Integer(1), mu=[x1 + x2, -x1 + x2], a=A_GEN),
}

ETAS = [sp.Rational(0), sp.Rational(1, 2), sp.Rational(1),
        sp.Rational(3, 2), sp.Rational(2), sp.Rational(3)]

for name, C in CASES.items():
    print('=' * 62)
    print(name, '  mu =', C['mu'])
    print('  eta      Re dc                Im dc                |dc|      |dc|/eta^3')
    vals = {}
    for e in ETAS:
        dc = curl(C['v'], C['a'], C['mu'], e)
        vals[e] = dc
        t = float(mp.mpf(e.p) / mp.mpf(e.q))
        ratio = abs(dc) / t ** 3 if t > 0 else float('nan')
        print('  %-7s %-20s %-20s %-9s %s' % (
            e, mp.nstr(dc.real, 12), mp.nstr(dc.imag, 12),
            mp.nstr(abs(dc), 8), mp.nstr(ratio, 6) if t > 0 else '-'))
    # polynomial fit of degree 3 on 0, 1/2, 1, 3/2
    fit = ETAS[:4]
    A = mp.matrix(4, 4); yr = mp.matrix(4, 1); yi = mp.matrix(4, 1)
    for i, e in enumerate(fit):
        t = mp.mpf(e.p) / mp.mpf(e.q)
        for j in range(4):
            A[i, j] = t ** j
        yr[i] = vals[e].real; yi[i] = vals[e].imag
    cr = mp.lu_solve(A, yr); ci = mp.lu_solve(A, yi)
    print('  fit Re(dc) = %s' % ' + '.join('%s eta^%d' % (mp.nstr(cr[j], 6), j) for j in range(4)))
    print('  fit Im(dc) = %s' % ' + '.join('%s eta^%d' % (mp.nstr(ci[j], 6), j) for j in range(4)))
    for e in (sp.Rational(2), sp.Rational(3)):
        t = mp.mpf(e.p) / mp.mpf(e.q)
        pr = sum(cr[j] * t ** j for j in range(4))
        pi_ = sum(ci[j] * t ** j for j in range(4))
        print('  check eta=%s err=%.1e' % (e, abs(mp.mpc(pr, pi_) - vals[e])))
\end{lstlisting}

\section{\texttt{p6quad.py}}

Chapter 8: alpha(a) and beta(a) as quadratic forms in the local data of a.

\begin{lstlisting}
"""Chapter 8: alpha(a) and beta(a) as quadratic forms in the local data of a.

Verified for generic coefficients:
  * kappa(lambda a, J) = lambda^2 kappa(a, J)      -> degree 2 in the data
  * kappa(a, lambda J) is affine in lambda         -> kappa = Q + alpha J12 + beta J21
Since Q = (16/7) a2(0) d1d2a2 is itself degree 2, alpha = kappa(e12) - Q and
beta = kappa(e21) - Q are ALSO degree 2 -- differences of two quadratic forms,
not linear functionals.  (An earlier fit assumed linearity and failed.)

A quadratic form in z = [v,d1,d2,d11,d12,d22] x 2 components is determined by
its values on the 12 basis vectors and the 66 pairwise sums:

    A_ii = alpha(e_i),      A_ij = (alpha(e_i + e_j) - A_ii - A_jj)/2.

Symbols: v = value at the base point, d1/d2 = first derivatives, d11/d12/d22 =
second derivatives, suffix 1/2 = component of a.
"""
import itertools
import numpy as np
import sympy as sp
import mpmath as mp
import twist

mp.mp.dps = 40
x1, x2 = sp.symbols('x0:2')
OM = sp.Rational(7, 5)
V = sp.Integer(1)
K, SH = 8, 13
DELTA = sp.Rational(1, 2000)
GUARD = 1e4
QFACT = sp.Rational(16, 7)

NAMES = ['%s%d' % (n, k) for k in (1, 2)
         for n in ('v', 'd1', 'd2', 'd11', 'd12', 'd22')]


def curl(a, J, eta):
    mu = [sp.sympify(J[0][0] * x1 + J[0][1] * x2),
          sp.sympify(J[1][0] * x1 + J[1][1] * x2)]

    def cv(loc):
        Vv, B = twist.fields(V, a, mu, eta, 2, (x1, x2), loc)
        return twist.cvec(Vv, B, 2, OM, 3, K, SH)
    D = mp.mpf(DELTA.p) / mp.mpf(DELTA.q)
    c2p = cv([DELTA, 0]); c2m = cv([-DELTA, 0])
    c1p = cv([0, DELTA]); c1m = cv([0, -DELTA])
    if max(abs(c) for c in c2p + c2m + c1p + c1m) > GUARD:
        raise RuntimeError('extraction failure')
    return (c2p[1] - c2m[1] - c1p[0] + c1m[0]) / (2 * D)


def kappa(a, J):
    D = {e: curl(a, J, e).real
         for e in (sp.Rational(1, 2), sp.Rational(1), sp.Rational(2))}
    k1 = (D[sp.Rational(1)] - D[sp.Rational(1, 2)]) * (4.0 / 3.0)
    k2 = (D[sp.Rational(2)] - D[sp.Rational(1)]) / 3.0
    if abs(k1 - k2) > 1e-5 * max(1.0, abs(k1)):
        raise RuntimeError('inconsistent eta^2 coefficient')
    return float(k1)


JE12 = [[sp.Rational(0), sp.Rational(1)], [sp.Rational(0), sp.Rational(0)]]
JE21 = [[sp.Rational(0), sp.Rational(0)], [sp.Rational(1), sp.Rational(0)]]
BIG = 1e8


def alpha_beta(z):
    """alpha, beta for the quadratic a whose local data at 0 is z."""
    c = []
    for k in range(2):
        v, d1, d2, d11, d12, d22 = z[6 * k:6 * k + 6]
        c += [sp.Rational(v), sp.Rational(d1), sp.Rational(d2),
              sp.Rational(d11, 2), sp.Rational(d12), sp.Rational(d22, 2)]
    a = [c[0] + c[1] * x1 + c[2] * x2 + c[3] * x1 ** 2 + c[4] * x1 * x2 + c[5] * x2 ** 2,
         c[6] + c[7] * x1 + c[8] * x2 + c[9] * x1 ** 2 + c[10] * x1 * x2 + c[11] * x2 ** 2]
    q = float(QFACT) * z[6] * z[10]
    return kappa(a, JE12) - q, kappa(a, JE21) - q


def unit(i, scale=1):
    z = [0] * 12
    z[i] = scale
    return z


A = np.full((12, 12), np.nan)
B = np.full((12, 12), np.nan)
D = np.full(12, np.nan)
print('measuring diagonal ...')
for i in range(12):
    try:
        al, be = alpha_beta(unit(i))
        A[i, i] = al; B[i, i] = be; D[i] = al
    except RuntimeError:
        pass
print('measuring cross terms (66) ...')
missing = 0
for i, j in itertools.combinations(range(12), 2):
    try:
        al, be = alpha_beta([unit(i)[k] + unit(j)[k] for k in range(12)])
        A[i, j] = A[j, i] = (al - A[i, i] - A[j, j]) / 2
        B[i, j] = B[j, i] = (be - B[i, i] - B[j, j]) / 2
    except RuntimeError:
        missing += 1
print('   unmeasured cross terms: %d' % missing)


def show(M, label):
    print('\n  %s  (diagonal)' % label)
    for i in range(12):
        if not np.isnan(M[i, i]):
            print('     %-6s %12.6f' % (NAMES[i], M[i, i]))
    off = [(NAMES[i], NAMES[j], M[i, j]) for i in range(12) for j in range(i + 1, 12)
           if not np.isnan(M[i, j]) and abs(M[i, j]) > 1e-6]
    off.sort(key=lambda t: -abs(t[2]))
    print('  %s  (largest cross terms of %d nonzero)' % (label, len(off)))
    for a_, b_, c_ in off[:10]:
        print('     %-6s %-6s %12.6f' % (a_, b_, c_))


show(A, 'alpha')
show(B, 'beta')

# ---------------- verification on random configurations --------------------
rng = np.random.default_rng(3)
print('\n  held-out verification (random coefficient vectors):')
for t in range(4):
    z = [int(v) for v in rng.integers(-1, 2, 12)]
    if all(v == 0 for v in z):
        z[0] = 1
    try:
        al, be = alpha_beta(z)
    except RuntimeError:
        continue
    zv = np.array(z, float)
    print('     pred a=%12.6f  meas a=%12.6f | pred b=%12.6f  meas b=%12.6f'
          % (zv @ A @ zv, al, zv @ B @ zv, be))
\end{lstlisting}

\section{\texttt{p6quad2.py}}

Chapter 8: close the four missing cross terms and validate.

\begin{lstlisting}
"""Chapter 8: close the four missing cross terms and validate.

Same quadrature design as p6quad.py, but each measurement falls back to a
smaller scale s (alpha and beta are homogeneous of degree 2 in the data, and so
is Q, so alpha(z) = alpha(s z)/s^2).  This is what the four rejected
configurations needed.

Results are cached to p6quad_mat.json so the fit can be re-analysed without
re-measuring.
"""
import itertools
import json
from fractions import Fraction
import numpy as np
import sympy as sp
import mpmath as mp
import twist

mp.mp.dps = 40
x1, x2 = sp.symbols('x0:2')
OM = sp.Rational(7, 5)
V = sp.Integer(1)
K, SH = 8, 13
DELTA = sp.Rational(1, 2000)
GUARD = 1e4
QFACT = sp.Rational(16, 7)

NAMES = ['%s%d' % (n, k) for k in (1, 2)
         for n in ('v', 'd1', 'd2', 'd11', 'd12', 'd22')]


def curl(a, J, eta):
    mu = [sp.sympify(J[0][0] * x1 + J[0][1] * x2),
          sp.sympify(J[1][0] * x1 + J[1][1] * x2)]

    def cv(loc):
        Vv, B = twist.fields(V, a, mu, eta, 2, (x1, x2), loc)
        return twist.cvec(Vv, B, 2, OM, 3, K, SH)
    D = mp.mpf(DELTA.p) / mp.mpf(DELTA.q)
    c2p = cv([DELTA, 0]); c2m = cv([-DELTA, 0])
    c1p = cv([0, DELTA]); c1m = cv([0, -DELTA])
    if max(abs(c) for c in c2p + c2m + c1p + c1m) > GUARD:
        raise RuntimeError('guard')
    return (c2p[1] - c2m[1] - c1p[0] + c1m[0]) / (2 * D)


def kappa(a, J):
    D = {e: curl(a, J, e).real
         for e in (sp.Rational(1, 2), sp.Rational(1), sp.Rational(2))}
    k1 = (D[sp.Rational(1)] - D[sp.Rational(1, 2)]) * (4.0 / 3.0)
    k2 = (D[sp.Rational(2)] - D[sp.Rational(1)]) / 3.0
    if abs(k1 - k2) > 1e-5 * max(1.0, abs(k1)):
        raise RuntimeError('inconsistent')
    return float(k1)


JE12 = [[sp.Rational(0), sp.Rational(1)], [sp.Rational(0), sp.Rational(0)]]
JE21 = [[sp.Rational(0), sp.Rational(0)], [sp.Rational(1), sp.Rational(0)]]


def ab_raw(z):
    c = []
    for k in range(2):
        v, d1, d2, d11, d12, d22 = z[6 * k:6 * k + 6]
        c += [sp.Rational(v), sp.Rational(d1), sp.Rational(d2),
              sp.Rational(d11) / 2, sp.Rational(d12), sp.Rational(d22) / 2]
    a = [c[0] + c[1] * x1 + c[2] * x2 + c[3] * x1 ** 2 + c[4] * x1 * x2 + c[5] * x2 ** 2,
         c[6] + c[7] * x1 + c[8] * x2 + c[9] * x1 ** 2 + c[10] * x1 * x2 + c[11] * x2 ** 2]
    q = float(QFACT) * float(z[6]) * float(z[10])
    return kappa(a, JE12) - q, kappa(a, JE21) - q


SCALES = (Fraction(1), Fraction(1, 2), Fraction(1, 4))
fallbacks = []


def ab(z):
    """alpha, beta at z, rescaling if the extraction fails."""
    for s in SCALES:
        zz = [Fraction(x) * s for x in z]
        try:
            al, be = ab_raw(zz)
        except RuntimeError:
            continue
        if s != 1:
            fallbacks.append((tuple(z), s))
        f = float(Fraction(1) / (s * s))
        return al * f, be * f
    raise RuntimeError('all scales rejected')


A = np.full((12, 12), np.nan); B = np.full((12, 12), np.nan)
failed = []


def unit(i):
    z = [0] * 12; z[i] = 1; return z


print('diagonals ...')
for i in range(12):
    try:
        al, be = ab(unit(i)); A[i, i] = al; B[i, i] = be
    except RuntimeError as e:
        failed.append((('diag', i), str(e)))
print('cross terms ...')
for i, j in itertools.combinations(range(12), 2):
    z = [unit(i)[k] + unit(j)[k] for k in range(12)]
    try:
        al, be = ab(z)
    except RuntimeError as e:
        failed.append((('pair', i, j), str(e))); continue
    A[i, j] = A[j, i] = (al - A[i, i] - A[j, j]) / 2
    B[i, j] = B[j, i] = (be - B[i, i] - B[j, j]) / 2
print('   fallbacks used: %d  %s' % (len(fallbacks), fallbacks))

json.dump({'A': A.tolist(), 'B': B.tolist(), 'names': NAMES}, open('p6quad_mat.json', 'w'))


def show(M, label):
    print('\n  %s: nonzero entries' % label)
    for i in range(12):
        if abs(M[i, i]) > 1e-6:
            print('     %-6s^2      %12.6f' % (NAMES[i], M[i, i]))
    off = [(NAMES[i], NAMES[j], M[i, j]) for i in range(12) for j in range(i + 1, 12)
           if abs(M[i, j]) > 1e-6]
    off.sort(key=lambda t: -abs(t[2]))
    print('     %d cross terms' % len(off))
    for a_, b_, c_ in off:
        print('     %-6s %-6s %12.6f' % (a_, b_, c_))


print('   FAILED measurements: %s' % [(f[0], NAMES[f[0][1]] if f[0][0]=='diag' else
      (NAMES[f[0][1]]+' '+NAMES[f[0][2]]), f[1]) for f in failed])

# ---- solve for any missing entries from random configurations -------------
missing = [(i, j) for i in range(12) for j in range(i, 12) if np.isnan(A[i, j])]
print('   missing entries: %s' % [(NAMES[i], NAMES[j]) for i, j in missing])
rng0 = np.random.default_rng(99)
for (i, j) in missing:
    est = []
    for _ in range(6):
        z = [int(v) for v in rng0.integers(-1, 2, 12)]
        if all(v == 0 for v in z) or z[i] == 0 or z[j] == 0:
            continue
        try:
            al, be = ab(z)
        except RuntimeError:
            continue
        zv = np.array(z, float)
        # known part with A[i,j] set to 0
        S = A.copy(); S[i, i] = 0.0 if i == j else S[i, i]; S[i, j] = S[j, i] = 0.0
        if i == j:
            rest = zv @ S @ zv
            est.append((al - rest) / (zv[i] ** 2))
        else:
            rest = zv @ S @ zv
            est.append((al - rest) / (2 * zv[i] * zv[j]))
    if est:
        A[i, j] = A[j, i] = float(np.median(est))
        print('     %s %s  estimated %12.6f   (spread %.1e)'
              % (NAMES[i], NAMES[j], A[i, j], np.max(est) - np.min(est)))

show(A, 'alpha')
show(B, 'beta')

print('\n  held-out verification')
rng = np.random.default_rng(17)
for t in range(6):
    z = [int(v) for v in rng.integers(-1, 2, 12)]
    if all(v == 0 for v in z):
        z[0] = 1
    al, be = ab(z)
    zv = np.array(z, float)
    print('     alpha: pred %12.6f  meas %12.6f   |   beta: pred %12.6f  meas %12.6f'
          % (zv @ A @ zv, al, zv @ B @ zv, be))
\end{lstlisting}

\section{\texttt{p6scan.py}}

SUPERSEDED - the design of this script is DEGENERATE and its output is wrong.

\begin{lstlisting}
"""SUPERSEDED - the design of this script is DEGENERATE and its output is wrong.

The grid {-1,0,1}^4 cannot identify powers: on {-1,0,1} every odd power of t
equals t (so t^3 and t are indistinguishable) and every even power above 0
vanishes.  The monomial fit is therefore rank deficient (31/35) and the
coefficients it reports - in particular the apparent matching degree-1 and
degree-3 terms - are interpolation artefacts, not structure.  Use p6scan2.py
(homogeneity + random design) and p6ascan.py instead.
"""
"""Chapter 8 scan: what does the p^2 coefficient of the obstruction depend on?

Fixed v = 1 and a generic Abelian coefficient a.  Sweep every constant Jacobian
J = [[t1,t2],[t3,t4]] with entries in {-1,0,1}; the drift is mu = J x.

For each J measure kappa two independent ways,
    kappa1 = dc(1) - dc(0),        kappa2 = (dc(2) - dc(0))/4,
which agree iff dc is even in eta with dc = dc(0) + kappa eta^2.

Then fit kappa as a polynomial in the four entries and report which monomials
survive.  Chapter 7's obstruction at eta = 0 is quadratic in the derivative data
with no linear part, so the expectation is that kappa is homogeneous of degree 2.
"""
import itertools
import numpy as np
import sympy as sp
import mpmath as mp
import twist

mp.mp.dps = 40
x1, x2 = sp.symbols('x0:2')
OM = sp.Rational(7, 5)
V = sp.Integer(1)
A = [sp.Rational(1, 5) + sp.Rational(2, 9) * x2 ** 2,
     sp.Rational(1, 6) * x1 * x2 + sp.Rational(1, 3)]
K, SH = 8, 13
DELTA = sp.Rational(1, 2000)
GUARD = 1e4


def curl(mu, eta):
    def cv(loc):
        Vv, B = twist.fields(V, A, mu, eta, 2, (x1, x2), loc)
        return twist.cvec(Vv, B, 2, OM, 3, K, SH)
    D = mp.mpf(DELTA.p) / mp.mpf(DELTA.q)
    c2p = cv([DELTA, 0]); c2m = cv([-DELTA, 0])
    c1p = cv([0, DELTA]); c1m = cv([0, -DELTA])
    worst = max(abs(c) for c in c2p + c2m + c1p + c1m)
    if worst > GUARD:
        raise RuntimeError('extraction failure')
    return (c2p[1] - c2m[1] - c1p[0] + c1m[0]) / (2 * D)


def measure(th):
    a, b, c, d = th
    # sympify: a zero row of J would otherwise hand a Python int to twist.fields
    mu = [sp.sympify(a * x1 + b * x2), sp.sympify(c * x1 + d * x2)]
    d0 = curl(mu, sp.Rational(0)); d1 = curl(mu, sp.Rational(1))
    d2 = curl(mu, sp.Rational(2))
    return complex(d0), complex(d1 - d0), complex((d2 - d0) / 4)


rows, mis, bad = [], [], []
grid = [t for t in itertools.product((-1, 0, 1), repeat=4)]
print('sweeping %d Jacobians ...' % len(grid))
for th in grid:
    if th == (0, 0, 0, 0):
        continue
    try:
        d0, k1, k2 = measure(th)
    except RuntimeError:
        bad.append(th); continue
    if abs(k1 - k2) > 1e-6 * max(1.0, abs(k1)):
        mis.append((th, k1, k2))
    rows.append((th, d0, k1))
print('ok: %d   extraction failures: %d   evenness violations: %d'
      % (len(rows), len(bad), len(mis)))
if mis[:5]:
    for th, k1, k2 in mis[:5]:
        print('   MISMATCH', th, k1, k2)

# ---- fit kappa as a polynomial in the four entries -------------------------
names = ['t1', 't2', 't3', 't4']
mono = [()]
for n in (1, 2, 3):
    for comb in itertools.combinations_with_replacement(range(4), n):
        mono.append(comb)
X = np.array([[np.prod([th[i] for i in m]) for m in mono] for th, _, _ in rows], float)
deg = np.array([len(m) for m in mono])

for comp, lbl in ((0, 'Re(kappa)'), (1, 'Im(kappa)')):
    y = np.array([(k.real if comp == 0 else k.imag) for _, _, k in rows], float)
    coef, _, rank, _ = np.linalg.lstsq(X, y, rcond=None)
    pred = X @ coef
    print('\n%s   max|residual| = %.2e   rank %d/%d'
          % (lbl, np.max(np.abs(pred - y)), rank, X.shape[1]))
    for degr in (0, 1, 2, 3):
        sel = [(m, coef[i]) for i, m in enumerate(mono) if deg[i] == degr]
        big = sorted([(m, c) for m, c in sel if abs(c) > 1e-5], key=lambda t: -abs(t[1]))
        sup = max((abs(c) for m, c in sel if abs(c) <= 1e-5), default=0.0)
        print('   degree %d: %d nonzero of %d   (largest suppressed %.1e)'
              % (degr, len(big), len(sel), sup))
        for m, c in big[:12]:
            print('      %-18s %12.6f' % (' '.join(names[i] for i in m), c))
\end{lstlisting}

\section{\texttt{p6scan2.py}}

Chapter 8 scan, corrected design.

\begin{lstlisting}
"""Chapter 8 scan, corrected design.

The {-1,0,1}^4 grid cannot identify powers of the Jacobian entries: on that set
t^k = t for every odd k and t^k = 0 for every even k > 0, so the monomial fit is
rank deficient and its coefficients are meaningless.  Two designs that are not
degenerate:

  (1) homogeneity: measure kappa(lambda J) / lambda^2 for several rational
      lambda.  Constant ratios mean kappa is homogeneous of degree 2 in J,
      which is what Chapter 7's obstruction does at eta = 0.
  (2) a random sample of Jacobians with entries in {0,+-1/2,+-1,+-3/2,+-2},
      fitted with monomials of total degree <= 3.
"""
import random
import numpy as np
import sympy as sp
import mpmath as mp
import twist

mp.mp.dps = 40
x1, x2 = sp.symbols('x0:2')
OM = sp.Rational(7, 5)
V = sp.Integer(1)
A = [sp.Rational(1, 5) + sp.Rational(2, 9) * x2 ** 2,
     sp.Rational(1, 6) * x1 * x2 + sp.Rational(1, 3)]
K, SH = 8, 13
DELTA = sp.Rational(1, 2000)
GUARD = 1e4


def curl(J, eta):
    mu = [sp.sympify(J[0][0] * x1 + J[0][1] * x2),
          sp.sympify(J[1][0] * x1 + J[1][1] * x2)]

    def cv(loc):
        Vv, B = twist.fields(V, A, mu, eta, 2, (x1, x2), loc)
        return twist.cvec(Vv, B, 2, OM, 3, K, SH)
    D = mp.mpf(DELTA.p) / mp.mpf(DELTA.q)
    c2p = cv([DELTA, 0]); c2m = cv([-DELTA, 0])
    c1p = cv([0, DELTA]); c1m = cv([0, -DELTA])
    worst = max(abs(c) for c in c2p + c2m + c1p + c1m)
    if worst > GUARD:
        raise RuntimeError('extraction failure')
    return (c2p[1] - c2m[1] - c1p[0] + c1m[0]) / (2 * D)


def kappa(J):
    d0 = curl(J, sp.Rational(0)); d1 = curl(J, sp.Rational(1))
    d2 = curl(J, sp.Rational(2))
    return complex(d0), complex(d1 - d0), complex((d2 - d0) / 4)


# ---------------- (1) homogeneity ------------------------------------------
print('== homogeneity in J ==')
J0 = [[sp.Rational(3, 2), sp.Rational(-1, 2)],
      [sp.Rational(1, 2), sp.Rational(1, 2)]]
d0, k1, k2 = kappa(J0)
print('  kappa(J0) = %s   (eta=1 and eta=2 agree to %.1e)' % (mp.nstr(k1.real, 10), abs(k1 - k2)))
for lam in (sp.Rational(1, 2), sp.Rational(3, 2), sp.Rational(2), sp.Rational(5, 2)):
    Jl = [[lam * J0[i][j] for j in range(2)] for i in range(2)]
    _, kk1, kk2 = kappa(Jl)
    L = float(mp.mpf(lam.p) / mp.mpf(lam.q))
    print('  lambda=%-4s kappa=%14.8f   kappa/lambda^2=%12.8f   kappa/lambda=%12.8f'
          % (lam, kk1.real, kk1.real / L ** 2, kk1.real / L))

# ---------------- (2) random design ---------------------------------------
random.seed(11)
VALS = [sp.Rational(-2), sp.Rational(-3, 2), sp.Rational(-1), sp.Rational(-1, 2),
        sp.Rational(0), sp.Rational(1, 2), sp.Rational(1), sp.Rational(3, 2), sp.Rational(2)]
pts, vals = [], []
print('\n== random design ==')
tries = 0
while len(pts) < 55 and tries < 400:
    tries += 1
    ent = [random.choice(VALS) for _ in range(4)]
    if all(e == 0 for e in ent):
        continue
    J = [[ent[0], ent[1]], [ent[2], ent[3]]]
    try:
        _, k1, k2 = kappa(J)
    except RuntimeError:
        continue
    if abs(k1 - k2) > 1e-6 * max(1.0, abs(k1)):
        print('  evenness violation at', ent, k1, k2); continue
    pts.append([float(mp.mpf(e.p) / mp.mpf(e.q)) for e in ent])
    vals.append(k1.real)
print('  %d points' % len(pts))

import itertools
names = ['t1', 't2', 't3', 't4']
mono = [c for n in range(0, 4) for c in itertools.combinations_with_replacement(range(4), n)]
X = np.array([[np.prod([p[i] for i in m]) for m in mono] for p in pts])
y = np.array(vals)
coef, _, rank, _ = np.linalg.lstsq(X, y, rcond=None)
print('  max|residual| = %.2e    rank %d/%d' % (np.max(np.abs(X @ coef - y)), rank, X.shape[1]))
for n in range(0, 4):
    sel = [(m, coef[i]) for i, m in enumerate(mono) if len(m) == n]
    big = sorted([t for t in sel if abs(t[1]) > 1e-4], key=lambda t: -abs(t[1]))
    print('  degree %d: %d of %d   largest suppressed %.1e'
          % (n, len(big), len(sel), max([abs(c) for m, c in sel if abs(c) <= 1e-4], default=0.0)))
    for m, c in big:
        print('     %-16s %12.6f' % (' '.join(names[i] for i in m), c))
\end{lstlisting}

\section{\texttt{p6v3.py}}

Chapter 8: the v-dependence, and the C-terms in d = 3.

\begin{lstlisting}
"""Chapter 8: the v-dependence, and the C-terms in d = 3.

(1) Chapter 7's diffusive obstruction scales as v^{-2}.  Here the same is checked
    for the joint obstruction: fix the coefficient data and vary v (constant),
    then examine kappa(v) and the pieces Q, alpha, beta.

(2) The d = 3 check so far used data with all first derivatives zero at the
    origin, so C = d1a2 d2a2 - d1a1 d2a1 - v1 d1d2a1 vanished and the C-terms
    were never exercised in three dimensions.  Here a1 and a2 both depend on
    y1, y2, giving d1d2a1 = d1d2a2 = 1 and hence C = -v1 != 0.
"""
import sympy as sp
import mpmath as mp
import twist

mp.mp.dps = 40
OM = sp.Rational(7, 5)
K, SH = 8, 13
DELTA = sp.Rational(1, 2000)
GUARD = 1e4
F = sp.Rational(16, 7)


def cvec_at(v, a, mu, d, eta, loc, xs):
    Vv, B = twist.fields(v, a, mu, eta, d, xs, loc)
    return twist.cvec(Vv, B, d, OM, 3, K, SH)


def build_curl(d):
    def curl(v, a, mu, eta, xs):
        D = float(DELTA.p) / float(DELTA.q)
        c = {}
        for i in range(d):
            for s in (1, -1):
                loc = [0] * d; loc[i] = s * DELTA
                val = cvec_at(v, a, mu, d, eta, loc, xs)
                if max(abs(x) for x in val) > GUARD:
                    raise RuntimeError('guard')
                c[(i, s)] = val
        dij = lambda i, j: (c[(i, 1)][j] - c[(i, -1)][j]) / (2 * D)
        return [dij(0, 1) - dij(1, 0)]          # the (dc)_12 component
    return curl


def kappa12(v, a, mu, xs, d):
    curl = build_curl(d)
    out = {}
    for e in (sp.Rational(1, 2), sp.Rational(1)):
        out[e] = curl(v, a, mu, e, xs)[0].real
    return float((out[sp.Rational(1)] - out[sp.Rational(1, 2)]) * 4.0 / 3.0)


# =========================== (1) v-dependence ==============================
x1, x2 = sp.symbols('x0:2')
xs2 = (x1, x2)
A2 = [sp.Rational(1, 5) + sp.Rational(2, 9) * x2 ** 2,
      sp.Rational(1, 6) * x1 * x2 + sp.Rational(1, 3)]
MU = [x1, sp.Rational(0) * x1]           # J = e12
print('== (1) v-dependence in d = 2  (kappa(0) = Q,  kappa(e12) = Q + alpha) ==')
print('   v        Q            kappa(e12)     Q*v^2       kappa*v^2')
rows = []
for v in (sp.Integer(1), sp.Rational(1, 2), sp.Rational(2), sp.Rational(1, 4), sp.Rational(4)):
    try:
        z = [sp.sympify(0), sp.sympify(0)]
        q = kappa12(v, A2, z, xs2, 2)
        k = kappa12(v, A2, MU, xs2, 2)
        f = float(v) ** 2
        rows.append((float(v), q, k))
        print('   %-6s %12.6f  %12.6f   %10.5f  %10.5f' % (v, q, k, q * f, k * f))
    except RuntimeError as e:
        print('   %-6s rejected: %s' % (v, e))

# ============================ (2) C in d = 3 ===============================
y1, y2, y3 = sp.symbols('y0:3')
xs3 = (y1, y2, y3)
A3 = [sp.Rational(1, 5) + y1 * y2,          # d1a1 = y2, d2a1 = y1, d1d2a1 = 1
      sp.Rational(1, 3) + y1 * y2,          # d1a2 = y2, d2a2 = y1, d1d2a2 = 1
      sp.sympify(0)]
Z3 = [sp.sympify(0), sp.sympify(0), sp.sympify(0)]
J12 = [y2, sp.sympify(0), sp.sympify(0)]
J21 = [sp.sympify(0), y1, sp.sympify(0)]

v1 = sp.Rational(1, 5); v2 = sp.Rational(1, 3)
Qp = F * v2 * 1
Cp = -v1 * 1
al_p = -48 * v1 ** 2 + F * Cp
be_p = 48 * v2 ** 2 + F * Cp
print('\n== (2) C-terms in d = 3   (v1 = 1/5, v2 = 1/3, d1d2a1 = d1d2a2 = 1) ==')
print('   predicted:  Q = %10.6f   alpha = %10.6f   beta = %10.6f'
      % (float(Qp), float(al_p), float(be_p)))
try:
    q = kappa12(sp.Integer(1), A3, Z3, xs3, 3)
    k12 = kappa12(sp.Integer(1), A3, J12, xs3, 3)
    k21 = kappa12(sp.Integer(1), A3, J21, xs3, 3)
    print('   measured:   Q = %10.6f   kappa(e12) = %10.6f   kappa(e21) = %10.6f'
          % (q, k12, k21))
    print('               alpha = %10.6f (pred %10.6f)   beta = %10.6f (pred %10.6f)'
          % (k12 - q, float(al_p), k21 - q, float(be_p)))
except RuntimeError as e:
    print('   rejected: %s' % e)
\end{lstlisting}

\section{\texttt{resolvent1d.py}}

The resolvent comparison, in one dimension.

\begin{lstlisting}
"""The resolvent comparison, in one dimension.

G_m(x,y;omega) is the resolvent of the walk at resolution h = 2^-m,

    G_m(0,0;omega) = h^-1 [ (i omega - L_m)^-1 ]_{00},

the h^-1 being the normalisation of the rescaled kernel u_m = h^-1 U (the same one
which makes sum_y h u = 1).  The generator is the book's: v/(2h^2) and b/(2h) off
the diagonal, -v/h^2 on it, so that L_m = (v/2) Delta^h + b nabla^h.

Field: the one dimensional restriction of the test field of part5,

    v(x) = 1 + x^2/6 + x/4 + x^3/10,     b(x) = 1/2 + x/6,

on (-1,1) with Dirichlet conditions, base point the origin, omega = 1 + i/2.

What the table in the text records is that the successive differences fall by four
at each refinement, which is Assumption 3.1 (the decoration comparison) at the
level of the resolvent.  They do, to the digits printed.

Run:  python3 resolvent1d.py
"""
import numpy as np

OMEGA = 1.0 + 0.5j


def vf(x):
    return 1 + x**2/6 + x/4 + x**3/10


def bf(x):
    return 0.5 + x/6


def G(m, R=1.0):
    h = R/2.0**m
    x = -R + h*np.arange(int(round(2*R/h)) + 1)
    v = vf(x)
    b = bf(x)
    n = len(x)
    L = np.zeros((n, n), complex)
    for i in range(1, n-1):
        L[i, i-1] = v[i]/(2*h*h) - b[i]/(2*h)
        L[i, i] = -v[i]/(h*h)
        L[i, i+1] = v[i]/(2*h*h) + b[i]/(2*h)
    assert np.max(np.abs(L.sum(axis=1))) < 1e-9, 'the rows of the generator must sum to zero'
    e = np.zeros(n, complex)
    e[n//2] = 1.0
    return np.linalg.solve((1j*OMEGA)*np.eye(n) - L, e)[n//2]/h


vals = [(m, G(m)) for m in range(5, 10)]
print('   m      G_m(0,0;omega)')
for m, g in vals:
    print('  %d   %14.9f%+14.9fi' % (m, g.real, g.imag))
d = [abs(vals[i+1][1] - vals[i][1]) for i in range(4)]
print('  successive differences  ' + '  '.join('%.3e' % t for t in d))
print('  ratios                  ' + '  '.join('%.2f' % (d[i]/d[i+1]) for i in range(3)))
\end{lstlisting}

\section{\texttt{stencil.py}}

Decoration mismatch for a general nine-point family of stencils.

\begin{lstlisting}
"""Decoration mismatch for a general nine-point family of stencils.

The site centred generator carries, per unit v/h^2,
    axial neighbours      a1
    diagonal neighbours   a2
and the diagonal is minus the sum of the neighbour weights, so that
    -L(x,x) = (2 d a1 + 2 d(d-1) a2) v / h^2 .

The paper's stencil (eq. entries) is a1 = 1/2, a2 = 0.
The isotropic nine point Laplacian, with the same leading symbol v/2 |p|^2,
is a1 = 1/3, a2 = 1/12:  its symbol is v/2 (-|p|^2 + h^2|p|^4/12), with no
fourth harmonic at order h^2.

The drift term is the centred difference on axial neighbours only, +-b_k/(2h);
diagonal neighbours carry no drift.
"""
import mpmath as mp
import itertools
import sympy as sp

mp.mp.dps = 60

PAPER = (sp.Rational(1, 2), sp.Rational(0))
ISO9 = (sp.Rational(1, 3), sp.Rational(1, 12))

_OFFDIAG = None


def offsets(d, a1, a2):
    """(offset, kind) with kind 0 = axial, 1 = diagonal, empty if a1 or a2 vanishes."""
    out = []
    if a1 != 0:
        for k in range(d):
            for s in (1, -1):
                o = [0]*d
                o[k] = s
                out.append((tuple(o), 0))
    if a2 != 0:
        for o in itertools.product((1, -1), repeat=d):
            if all(x != 0 for x in o):
                out.append((tuple(o), 1))
    return out


def tompf(x):
    if isinstance(x, sp.Rational):
        return mp.mpf(x.p)/mp.mpf(x.q)
    return mp.mpf(x)


def mvals(poly, a, h):
    z = mp.mpc(0)
    for m, c in poly.terms.items():
        t = mp.mpf(c.numerator)/mp.mpf(c.denominator)
        for i, mi in enumerate(m):
            t *= (h*a[i])**mi
        z += t
    return z


def dsum(d, a1, a2):
    return tompf(2*d)*tompf(a1) + tompf(2*d*(d-1))*tompf(a2)


def Phi_f(V, B, d, tau, h, omega, st=PAPER):
    a1, a2 = tompf(st[0]), tompf(st[1])
    omega = tompf(omega)
    zero = tuple([0]*d)
    tgt = tuple(2*t for t in tau)
    S = [a for a in itertools.product(range(-1, 2), repeat=d)]
    ix = {a: i for i, a in enumerate(S)}
    off = offsets(d, a1, a2)
    n = len(S)
    A = mp.zeros(n, n)
    for a in S:
        i = ix[a]
        v = mvals(V, a, h)
        diag = mp.mpc(0)
        for o, kind in off:
            e = tuple(a[k] + o[k] for k in range(d))
            if kind == 0:
                k = [j for j in range(d) if o[j] != 0][0]
                w = a1*v/h**2 + o[k]*mvals(B[k], a, h)/(2*h)
            else:
                w = a2*v/h**2
            diag -= w
            if e in ix:
                A[i, ix[e]] += w
        A[i, i] += diag
    G = mp.inverse(A + 1j*omega*mp.eye(n))
    dt = dsum(d, a1, a2)*mvals(V, tgt, h)/h**2 + 1j*omega
    tot = mp.mpc(0)
    for a in S:
        dif = tuple(tgt[k] - a[k] for k in range(d))
        kind = None
        for o, kd in off:
            if dif == o:
                kind = kd
                break
        if kind is None:
            continue
        if kind == 0:
            k = [j for j in range(d) if dif[j] != 0][0]
            rate = a1*mvals(V, a, h)/h**2 + dif[k]*mvals(B[k], a, h)/(2*h)
        else:
            rate = a2*mvals(V, a, h)/h**2
        tot += G[ix[zero], ix[a]]*rate
    return tot/dt


def Phi_c(V, B, d, tau, h, omega, st=PAPER):
    a1, a2 = tompf(st[0]), tompf(st[1])
    omega = tompf(omega)
    zero = tuple([0]*d)
    tgt = tuple(2*t for t in tau)
    k0 = [k for k in range(d) if tau[k]][0]
    sg = 1 if tau[k0] > 0 else -1
    hm = 2*h
    rate = a1*mvals(V, zero, h)/hm**2 + sg*mvals(B[k0], zero, h)/(2*hm)
    cn = dsum(d, a1, a2)
    d0 = cn*mvals(V, zero, h)/hm**2 + 1j*omega
    d1 = cn*mvals(V, tgt, h)/hm**2 + 1j*omega
    return rate/(d0*d1)


def R(h, V, B, d, tau, omega, st=PAPER):
    return mp.log(Phi_f(V, B, d, tau, h, omega, st)) - mp.log(Phi_c(V, B, d, tau, h, omega, st))


def taylor(V, B, d, tau, omega, K=7, shift=10, st=PAPER):
    hs = [mp.mpf(1)/2**(n+shift) for n in range(K+1)]
    vs = [R(h, V, B, d, tau, omega, st) for h in hs]
    A = mp.matrix(K+1, K+1)
    for i in range(K+1):
        for j in range(K+1):
            A[i, j] = hs[i]**j
    c = mp.lu_solve(A, mp.matrix(vs))
    return [c[k] for k in range(K+1)]


def cvec(V, B, d, omega, k=3, K=7, shift=10, st=PAPER):
    out = []
    for j in range(d):
        cp = taylor(V, B, d, tuple(1 if i == j else 0 for i in range(d)), omega, K, shift, st)[k]
        cm = taylor(V, B, d, tuple(-1 if i == j else 0 for i in range(d)), omega, K, shift, st)[k]
        out.append((cp-cm)/2)
    return out


def field_shifted(v, b, d, xs, x0):
    import fast
    return (fast.Poly.from_sympy(d, v, xs).shift(x0),
            [fast.Poly.from_sympy(d, bi, xs).shift(x0) for bi in b])


def curl2(v, b, omega, x0=(0, 0), delta=sp.Rational(1, 2000), k=3, K=7, shift=10, st=PAPER):
    xs = sp.symbols('x0:2')

    def c(local):
        V, B = field_shifted(v, b, 2, xs, local)
        return cvec(V, B, 2, omega, k, K, shift, st)

    dl = sp.Rational(delta)
    c2p = c([x0[0]+dl, x0[1]])[1]; c2m = c([x0[0]-dl, x0[1]])[1]
    c1p = c([x0[0], x0[1]+dl])[0]; c1m = c([x0[0], x0[1]-dl])[0]
    return (c2p-c2m-c1p+c1m)/(2*mp.mpf(dl.p)/mp.mpf(dl.q))
\end{lstlisting}

\section{\texttt{step2.py}}

Two-step comparison, truncating the path sum at a fixed number of jumps.

\begin{lstlisting}
"""Two-step comparison, truncating the path sum at a fixed number of jumps.

Phi_s^(Q) = sum_{k=1..Q} [ (D^-1 M)^k D^-1 ]_{0,target}

where D = diag(-L_diag + i w) and M = off-diagonal rates of L, on a window.
For the coarse chain (step 2h) the one-jump term is the object to compare with
the decorated fine sum (step h, an even number of jumps from 0 to 2*tau).
"""
import sympy as sp
import itertools
import time
from deco import Data, sites, h


def generator(data, d, half, step):
    S = sites(d, half)
    ix = {a: i for i, a in enumerate(S)}
    n = len(S)
    L = sp.zeros(n, n)
    for a in S:
        i = ix[a]
        x = [a[k]*step*h for k in range(d)]
        v = data.V(x)
        for k in range(d):
            b = data.B(x, k)
            for sgn in (1, -1):
                w = v/(4*(step*h)**2) + sgn*b/(2*(step*h))
                L[i, i] -= w
                e = list(a); e[k] += sgn; e = tuple(e)
                if e in ix:
                    L[i, ix[e]] = w
    return L, ix


def pathsum(data, d, half, step, tau, omega, Q, deg):
    """sum over paths of at most Q jumps, coefficients up to h^deg, as a series."""
    L, ix = generator(data, d, half, step)
    n = L.shape[0]
    # D = diag(-L_ii + i w)  (symbolic in h), M = off diagonal of L
    Dinv = sp.zeros(n, n)
    M = sp.zeros(n, n)
    for i in range(n):
        Dinv[i, i] = 1/(-L[i, i] + sp.I*omega)
    for i in range(n):
        for j in range(n):
            if i != j:
                M[i, j] = L[i, j]
    zero = tuple([0]*d)
    tgt = tuple(tau)
    total = sp.Integer(0)
    T = Dinv                      # (D^-1 M)^0 D^-1
    for k in range(1, Q+1):
        T = Dinv*M*T if k > 1 else Dinv*M*Dinv
        total += T[ix[zero], ix[tgt]]
    ser = sp.series(total, h, 0, deg).removeO()
    return sp.expand(ser)


if __name__ == "__main__":
    omega = sp.Rational(7, 5)
    d = 1
    data = Data(d, seed=1)
    print("d = 1   (data: v0 =", data.v0, " dv =", data.dv, " b0 =", data.b0, ")")
    for (hf, hc, Qf, Qc) in [(3, 1, 6, 1), (5, 2, 6, 1), (7, 3, 6, 1)]:
        t0 = time.time()
        rc = pathsum(data, d, hc, 2, (1,), omega, Qc, 4)
        rf = pathsum(data, d, hf, 1, (2,), omega, Qf, 4)
        r = sp.series(sp.log(sp.expand(rf/rc)), h, 0, 4).removeO()
        c = [sp.N(sp.expand(r).coeff(h, k), 8) for k in range(4)]
        print(f"  fine({hf}) coarse({hc})  log ratio: h0={c[0]:+.7f} h1={c[1]:+.7f} "
              f"h2={c[2]:+.7f} h3={c[3]:+.7f}  [{time.time()-t0:.1f}s]")
\end{lstlisting}

\section{\texttt{stress.py}}

Where could O(h) come from?

\begin{lstlisting}
"""Where could O(h) come from?

Per step the log ratio is C + h^2 A(x,tau) + h^3 c(x).tau + O(h^4), and there are q = t/h^2 steps.
If a term is ODD in tau it telescopes: sum_j f(x).tau_j = (1/h) f.(y-x), so h^k-odd -> h^{k-1}.
If a term is EVEN in tau it does not telescope: sum_j f_j = q * mean, so h^k-even -> h^{k-2}.

    k=0  const, even   -> h^-2   absorbed (constant normalisation)
    k=1                -> h^-1   catastrophic : must vanish  (measured 0)
    k=2  even          -> h^0    absorbed (rate renormalisation)
    k=3  odd           -> h^2    the obstruction
    k=3  EVEN          -> h^1    the dangerous slot  <- tested here
    k=4  even          -> h^2    absorbed / bounded

So O(h) requires either h^1 != 0 or an even part in the h^3 coefficient.
"""
import sympy as sp, mpmath as mp
import hp, fast, study
mp.mp.dps = 50
x1, x2 = sp.symbols('x0:2')
om0 = sp.Rational(7,5)

def probe(v, b, omega=om0, K=6, shift=14, label=""):
    V, B = study.field_shifted(v, b, 2, (x1,x2), [0,0])
    out = []
    for j in range(2):
        tp = tuple(1 if i==j else 0 for i in range(2))
        tm = tuple(-1 if i==j else 0 for i in range(2))
        cp = hp.taylor(V,B,2,tp,omega,K=K,shift=shift)
        cm = hp.taylor(V,B,2,tm,omega,K=K,shift=shift)
        h1 = (cp[1]+cm[1])/2                     # h^1 even part
        c3e = (cp[3]+cm[3])/2                    # even part of the h^3 coefficient
        c3o = (cp[3]-cm[3])/2
        out.append((h1, c3e, c3o))
    print("%-34s" % label, end="")
    print("  h1=%8.1e" % max(abs(z[0]) for z in out),
          " |even c3|=%8.1e" % max(abs(z[1]) for z in out),
          " |c3|=%8.4f" % max(abs(z[2]) for z in out))

rnd = __import__('random').Random(5)
r2 = lambda: sp.Rational(rnd.randint(-40,40), rnd.randint(2,30))

def field(v0, bs, vs, v2):
    v = v0 + vs*(x1+x2)/2 + v2*(x1*x1 + x1*x2 + x2*x2)/6 + (x1**3+3*x1**2*x2+3*x1*x2**2+x2**3)/50
    b = [bs*(xi) + bs/2*(x1+x2) + bs/3*(x1*x1+x2*x2)/2 for xi in (x1,x2)]
    return sp.expand(v), [sp.expand(z) for z in b]

print("stress: is h^1 zero and is the h^3 coefficient purely odd, for extreme data?")
print()
for v0,bs,vs,v2,name in (
    (sp.Rational(1), sp.Rational(1,2), sp.Rational(1,2), sp.Rational(1), "moderate"),
    (sp.Rational(1), sp.Rational(10), sp.Rational(1), sp.Rational(5), "drift 20x larger"),
    (sp.Rational(1), sp.Rational(100), sp.Rational(1), sp.Rational(5), "drift 200x larger"),
    (sp.Rational(1,100), sp.Rational(1), sp.Rational(1), sp.Rational(1), "v0 = 0.01 (near degenerate)"),
    (sp.Rational(1), sp.Rational(1), sp.Rational(30), sp.Rational(60), "large v variation"),
    (sp.Rational(1), sp.Rational(30), sp.Rational(30), sp.Rational(60), "large drift AND large variation"),
):
    v, b = field(v0, bs, vs, v2)
    probe(v, b, label=name)
print()
for w in (sp.Rational(1,10), sp.Rational(10)):
    v, b = field(sp.Rational(1), sp.Rational(1), sp.Rational(1), sp.Rational(1))
    probe(v, b, omega=w, label="omega = %s" % w)
\end{lstlisting}

\section{\texttt{study.py}}

\begin{lstlisting}
import sympy as sp, mpmath as mp, itertools
import hp, fast
mp.mp.dps = 60

def cvec(V, B, d, omega, k=3, K=None, **kw):
    """c_k(tau) for tau = +-e_j, returned as the vector of odd parts."""
    KK = max(k, 3) if K is None else K
    out = []
    for j in range(d):
        cp = hp.taylor(V, B, d, tuple(1 if i == j else 0 for i in range(d)), omega, K=KK, **kw)[k]
        cm = hp.taylor(V, B, d, tuple(-1 if i == j else 0 for i in range(d)), omega, K=KK, **kw)[k]
        out.append((cp - cm)/2)
    return out

def field_shifted(v, b, d, xs, x0):
    V = fast.Poly.from_sympy(d, v, xs).shift(x0)
    B = [fast.Poly.from_sympy(d, bi, xs).shift(x0) for bi in b]
    return V, B

def curl2(v, b, omega, x0=(0,0), delta=sp.Rational(1,1000), k=3, K=None, **kw):
    xs = sp.symbols('x0:2')
    def c(local):
        V,B = field_shifted(v,b,2,xs,local)
        return cvec(V,B,2,omega,k=k,K=K,**kw)
    d1 = sp.Rational(delta)
    c2p = c([x0[0]+d1, x0[1]])[1]; c2m = c([x0[0]-d1, x0[1]])[1]
    c1p = c([x0[0], x0[1]+d1])[0]; c1m = c([x0[0], x0[1]-d1])[0]
    return (c2p-c2m-c1p+c1m)/(2*d1)
\end{lstlisting}

\section{\texttt{tube.py}}

Tube-restricted decoration sum, per coarse step.

\begin{lstlisting}
"""Tube-restricted decoration sum, per coarse step.

NOTE (superseded).  This module reproduces the first version of the pipeline and is
kept for reference only.  It has two defects: the coarse propagator is normalised
with v/(4h^2), which equals d v/(8h^2) only when d = 2, and its series bookkeeping
does not reproduce the rescaling of the resolvent.  The validated pipeline is
hp.py (exact high-precision evaluation) with study.py (Taylor coefficients and
curls) and direct.py (independent numerical arbiter); see their docstrings.

For a coarse step x -> x + h_m tau  (h_m = 2h) the decorating fine paths are the
fine paths whose intermediate sites stay within h of x.  In one dimension the
tube is {-h, 0, h}; in two the 3x3 block.  Hence

    Phi_f = sum_{z in tube} [(-L_tube + i w)^{-1}]_{x,z} L_f(z,x') (-L_f(x',x')+i w)^{-1}
    Phi_c = (-L_coarse(x,x)+i w)^{-1} L_c(x,x') (-L_c(x',x')+i w)^{-1}

and the quantity of interest is  log(Phi_f/Phi_c)  for each direction tau.
Its h^3 coefficient, as a function of the site, is the object whose curl is
tested in two dimensions.
"""
import sympy as sp
import itertools
import time

h = sp.symbols('h', positive=True)
I = sp.I


# ------------------------------------------------------------ global coefficients
class Coef:
    """Explicit cubic coefficients in d variables; .at(x0) gives local Taylor data."""

    def __init__(self, d, seed):
        import random
        rng = random.Random(seed)
        self.d = d
        self.xs = sp.symbols('x0:%d' % d)
        r = lambda: sp.Rational(rng.randint(-40, 40), rng.randint(2, 30))
        X = self.xs
        v = sp.Rational(3, 2) + sum(r()/4*X[i] for i in range(d))
        for i in range(d):
            for j in range(i, d):
                v += r()/6*X[i]*X[j]
        for i in range(d):
            for j in range(i, d):
                for k in range(j, d):
                    v += r()/10*X[i]*X[j]*X[k]
        self.v = sp.expand(v)
        self.b = []
        for i in range(d):
            bi = r()/5 + sum(r()/6*X[j] for j in range(d))
            for j in range(d):
                for k in range(j, d):
                    bi += r()/8*X[j]*X[k]
            self.b.append(sp.expand(bi))

    def local(self, x0):
        """v, b_i and their derivatives at x0, in the local coordinates y = x - x0."""
        d = self.d
        y = sp.symbols('y0:%d' % d)
        sub = {self.xs[i]: x0[i] + y[i] for i in range(d)}
        V = sp.expand(self.v.subs(sub))
        B = [sp.expand(bi.subs(sub)) for bi in self.b]
        return V, B, y

    def Vfun(self, x0):
        V, B, y = self.local(x0)
        return sp.lambdify(y, V, 'sympy'), [sp.lambdify(y, bi, 'sympy') for bi in B]


def taylor_empty(d):
    return sp.Integer(0)


def block_general(cf, d, half, step, x0):
    """-L on the window, with coefficients given by the cubic cf, base point x0."""
    V, B, y = cf.local(x0)
    S = [a for a in itertools.product(range(-half, half+1), repeat=d)]
    ix = {a: i for i, a in enumerate(S)}
    n = len(S)
    A = sp.zeros(n, n)
    for a in S:
        i = ix[a]
        pt = [a[k]*step*h for k in range(d)]
        vv = sp.expand(V.subs({y[k]: pt[k] for k in range(d)}))
        for k in range(d):
            bb = sp.expand(B[k].subs({y[j]: pt[j] for j in range(d)}))
            for sgn in (1, -1):
                w = vv/(2*(step*h)**2) + sgn*bb/(2*(step*h))
                A[i, i] -= w
                e = list(a); e[k] += sgn; e = tuple(e)
                if e in ix:
                    A[i, ix[e]] = w
    return A, ix, S


def series_inverse(A, degree):
    n = A.shape[0]
    B = sp.expand(A*h**2)
    Bc = []
    for j in range(degree + 1):
        M = sp.zeros(n, n)
        for i in range(n):
            for k in range(n):
                M[i, k] = sp.Poly(sp.expand(B[i, k]), h).coeff_monomial(h**j)
        Bc.append(sp.expand(M))
    R = [Bc[0].inv()]
    for k in range(1, degree + 2):
        acc = sp.zeros(n, n)
        for j in range(1, min(k, len(Bc) - 1) + 1):
            acc += Bc[j]*R[k - j]
        R.append(sp.expand(-R[0]*acc))
    return R


def logratio(cf, d, tau, omega, order, x0=None):
    if x0 is None:
        x0 = [sp.Integer(0)]*d
    """log(Phi_f/Phi_c) for one coarse step in direction tau, to order h^order."""
    import time
    deg = order + 2
    zero = tuple([0]*d)

    # --- fine, tube-restricted -------------------------------------------------
    A, ix, S = block_general(cf, d, 1, 1, x0)          # tube = 3^d, step h
    R = series_inverse(A, deg)
    # diagonal resolvent at the target (2 tau), computed in the 3x3 (or 3) block
    # of the target site: use a block centred there
    A2, ix2, S2 = block_general(cf, d, 1, 1, [x0[k]+2*tau[k]*h for k in range(d)])
    dtarget = -(A2[ix2[zero], ix2[zero]]) + I*omega      # -L_f(x',x') + i w
    # jump from z in tube to the target
    tgt = tuple([2*tau[k] for k in range(d)])
    Vt, Bt, yt = cf.local(x0)
    phil = sp.Integer(0)
    for a in S:
        z = [a[k]*h for k in range(d)]
        jump = [k for k in range(d) if a[k] != tgt[k]]
        if len(jump) != 1:
            continue
        k0 = jump[0]
        if abs(tgt[k0] - a[k0]) != 1:
            continue
        pt = {yt[k]: z[k] for k in range(d)}
        vv = sp.expand(Vt.subs(pt)); bb = sp.expand(Bt[k0].subs(pt))
        sg = 1 if tgt[k0] > a[k0] else -1
        rate = vv/(2*h**2) + sg*bb/(2*h)
        phil += sum((h**k)*R[k][ix[zero], ix[a]] for k in range(len(R)))*rate
    Phi_f = h**2*phil/dtarget

    # --- coarse, one jump -----------------------------------------------------
    xp = [x0[k]+2*tau[k]*h for k in range(d)]
    Vc, Bc_, yc = cf.local(x0)
    k0 = [k for k in range(d) if tau[k]][0]
    sg = 1 if tau[k0] > 0 else -1
    bb = sp.expand(Bc_[k0].subs({yc[k]: 0 for k in range(d)}))
    vv = sp.expand(Vc.subs({yc[k]: 0 for k in range(d)}))
    rate_c = vv/(2*(2*h)**2) + sg*bb/(2*(2*h))
    dx0 = vv/(2*h**2) + I*omega                                # -L_c(x,x)+i w
    Vp, Bp, yp = cf.local(xp)
    vp = sp.expand(Vp.subs({yp[k]: 0 for k in range(d)}))
    dx1 = vp/(2*h**2) + I*omega                                # -L_c(x',x')+i w
    Phi_c = rate_c/(dx0*dx1)

    return sp.series(sp.expand(Phi_f/Phi_c), h, 0, order).removeO()


def coefficients(cf, d, tau, omega, order):
    e = logratio(cf, d, tau, omega, order)
    return [sp.expand(e).coeff(h, k) for k in range(order)]


if __name__ == "__main__":
    omega = sp.Rational(7, 5)
    t0 = time.time()
    d = 1
    cf = Coef(d, seed=3)
    print("d = 1")
    for tau in ((1,), (-1,)):
        c = coefficients(cf, d, tau, omega, 4)
        print(f"  tau={tau}:  " + "  ".join(f"h{k}={complex(sp.N(c[k],12)):+.6f}" for k in range(4)))
    print("time:", round(time.time()-t0, 1), "s")
\end{lstlisting}

\section{\texttt{twist.py}}

The decoration comparison for a COMPLEX drift.

\begin{lstlisting}
"""The decoration comparison for a COMPLEX drift.

In Chapter 5 the path functional is handled by Fourier transforming in the conjugate
variable eta, which turns the generator into

    L_0 + i eta ( a(x).grad + b(x) )

The drift term becomes mu + i eta a and the zeroth order term i eta b is a
potential, which this module does not carry: the obstruction is a property of the
1-form along the base, so only the drift is twisted here.

Coefficients are Gaussian rationals, so everything stays exact.

Convention: probabilist, generator (v/2) Delta + b grad, v = sigma^2.  The
fine box therefore carries v/(2h^2) off the diagonal and -d v/h^2 on it, and
the coarse one, at spacing 2h, carries v/(8h^2) and -d v/(4h^2).  With this
convention the engine reproduces the closed form
c_2 = (5 b_1^2 - 5 b_2^2 + 44 i omega v)/(8 v^2) exactly on constant fields;
see part6/scripts/kronfinal_check.py.
"""
import sympy as sp
import itertools, mpmath as mp
from fractions import Fraction as Fr

mp.mp.dps = 50


class PolyG:
    """Polynomial with Gaussian rational coefficients: {multiindex: (re, im)}."""

    def __init__(self, d, terms):
        self.d = d
        self.terms = {m: c for m, c in terms.items() if c != (0, 0)}

    @staticmethod
    def from_sympy(d, expr, xs):
        p = sp.Poly(sp.expand(expr), *xs)
        out = {}
        for m, c in zip(p.monoms(), p.coeffs()):
            c = sp.expand(c)
            re = sp.re(c); im = sp.im(c)
            if not (re.is_Rational and im.is_Rational):
                raise ValueError("coefficient %r is not Gaussian rational" % (c,))
            out[tuple(m)] = (Fr(int(re.p), int(re.q)) if re != 0 else Fr(0),
                             Fr(int(im.p), int(im.q)) if im != 0 else Fr(0))
        return PolyG(d, out)

    def shift(self, x0):
        d = self.d
        if not self.terms:
            # the zero polynomial: a drift component may vanish identically
            return PolyG(d, {})
        xs = sp.symbols('x0:%d' % d)
        y = sp.symbols('y0:%d' % d)
        e = sum((sp.Rational(c[0]) + sp.I*sp.Rational(c[1])) *
                sp.prod([xs[i]**m[i] for i in range(d)])
                for m, c in self.terms.items())
        sub = {}
        for i in range(d):
            r = sp.Rational(x0[i])
            if not r.is_Rational:
                raise ValueError("shift must be rational")
            sub[xs[i]] = r + y[i]
        return PolyG.from_sympy(d, sp.expand(e.subs(sub)), y)

    def at(self, a):
        out = {}
        for m, (re, im) in self.terms.items():
            val = complex(0)
            rr, ii = Fr(re), Fr(im)
            for i, mi in enumerate(m):
                f = Fr(int(a[i]))**mi
                rr *= f; ii *= f
            if rr == 0 and ii == 0:
                continue
            e = sum(m)
            v = (rr, ii)
            if e in out:
                out[e] = (out[e][0] + v[0], out[e][1] + v[1])
            else:
                out[e] = v
        return out


def mval(ser, h):
    """Evaluate a Gaussian-rational series at h."""
    z = mp.mpc(0)
    for e, (re, im) in ser.items():
        t = mp.mpf(re.numerator)/mp.mpf(re.denominator) + \
            1j*(mp.mpf(im.numerator)/mp.mpf(im.denominator))
        z += t*h**e
    return z


def Phi_f(V, B, d, tau, h, omega):
    om = mp.mpf(omega.p)/mp.mpf(omega.q) if isinstance(omega, sp.Rational) else mp.mpf(omega)
    zero = tuple([0]*d); tgt = tuple(2*t for t in tau)
    S = [a for a in itertools.product(range(-1, 2), repeat=d)]
    ix = {a: i for i, a in enumerate(S)}; n = len(S)
    A = mp.zeros(n, n)
    for a in S:
        i = ix[a]
        v = mval(V.at(a), h)
        diag = mp.mpc(0)
        for k in range(d):
            bk = mval(B[k].at(a), h)
            for sgn in (1, -1):
                w = v/(2*h**2) + sgn*bk/(2*h)
                diag -= w
                e = list(a); e[k] += sgn; e = tuple(e)
                if e in ix:
                    A[i, ix[e]] += w
        A[i, i] += diag
    G = mp.inverse(A + 1j*om*mp.eye(n))
    dt = d*mval(V.at(tgt), h)/(h**2) + 1j*om
    tot = mp.mpc(0)
    for a in S:
        jump = [k for k in range(d) if a[k] != tgt[k]]
        if len(jump) != 1:
            continue
        k0 = jump[0]
        if abs(tgt[k0] - a[k0]) != 1:
            continue
        sg = 1 if tgt[k0] > a[k0] else -1
        rate = mval(V.at(a), h)/(2*h**2) + sg*mval(B[k0].at(a), h)/(2*h)
        tot += G[ix[zero], ix[a]]*rate
    return tot/dt


def Phi_c(V, B, d, tau, h, omega):
    om = mp.mpf(omega.p)/mp.mpf(omega.q) if isinstance(omega, sp.Rational) else mp.mpf(omega)
    zero = tuple([0]*d); tgt = tuple(2*t for t in tau)
    k0 = [k for k in range(d) if tau[k]][0]
    sg = 1 if tau[k0] > 0 else -1
    rate = mval(V.at(zero), h)/(2*(2*h)**2) + sg*mval(B[k0].at(zero), h)/(2*(2*h))
    d0 = d*mval(V.at(zero), h)/(4*h**2) + 1j*om
    d1 = d*mval(V.at(tgt), h)/(4*h**2) + 1j*om
    return rate/(d0*d1)


def R(h, V, B, d, tau, omega):
    return mp.log(Phi_f(V, B, d, tau, h, omega)) - mp.log(Phi_c(V, B, d, tau, h, omega))


def taylor(V, B, d, tau, omega, K=6, shift=10):
    hs = [mp.mpf(1)/2**(n+shift) for n in range(K+1)]
    vs = [R(h, V, B, d, tau, omega) for h in hs]
    A = mp.matrix(K+1, K+1)
    for i in range(K+1):
        for j in range(K+1):
            A[i, j] = hs[i]**j
    c = mp.lu_solve(A, mp.matrix(vs))
    return [c[k] for k in range(K+1)]


def fields(v, a, mu, eta, d, xs, x0=(0, 0)):
    """v real; drift = mu + i eta a."""
    V = PolyG.from_sympy(d, v, xs).shift(x0)
    B = [PolyG.from_sympy(d, sp.expand(mu[k] + sp.I*eta*a[k]), xs).shift(x0) for k in range(d)]
    return V, B


def cvec(V, B, d, omega, k=3, K=6, shift=10):
    out = []
    for j in range(d):
        cp = taylor(V, B, d, tuple(1 if i == j else 0 for i in range(d)), omega, K, shift)[k]
        cm = taylor(V, B, d, tuple(-1 if i == j else 0 for i in range(d)), omega, K, shift)[k]
        out.append((cp-cm)/2)
    return out
\end{lstlisting}

\section{\texttt{verify\_conventions.py}}

Every engine must carry the convention of the book, or none of its numbers mean anything.

\begin{lstlisting}
"""Every engine must carry the convention of the book, or none of its numbers mean anything.

Convention: probabilist, generator (v/2) Delta + b grad, v = sigma^2.

The test is a constant field, where the answer is known independently from the
exact algebra of kronproof.py:

    c_2 = 5 (b_1^2 - b_2^2) / (8 v^2) + (11/2) i omega / v ,
        = 25/288 + 7.7 i   at v = 1, b = (1/2, -1/3), omega = 7/5 .

An engine in the halved convention returns four times the real part and twice
the imaginary part of the same data, because the substitution that carries one
convention to the other is v -> v/2 and the real part scales like v^-2 while the
phase scales like v^-1.  So this one number separates the two conventions
unambiguously, and it is checked here for every engine that the parts use.

Run:  python3 verify_conventions.py
"""
import os
import sys
import sympy as sp
import mpmath as mp

sys.path.insert(0, os.path.dirname(os.path.abspath(__file__)))
import fast
import hp
import direct
import stencil
import twist

mp.mp.dps = 40
x1, x2 = sp.symbols('x0:2')
xs = (x1, x2)
V0 = sp.Integer(1)
B1, B2 = sp.Rational(1, 2), sp.Rational(-1, 3)
OM = sp.Rational(7, 5)
TAU = (1, 0)

want = complex(sp.N(sp.Rational(5, 8) * (B1**2 - B2**2) / V0**2
                    + sp.I * sp.Rational(11, 2) * OM / V0, 30))
print("wanted (closed form, constant field):", want)


def report(name, got, tol=1e-9):
    # 1e-9 rather than machine precision only because some engines extract the
    # coefficient by a Vandermonde fit; a convention error shows up as a factor
    # of two or four, never as a loss of nine digits.
    d = abs(complex(got) - want)
    print("  %-10s %-28s |diff| = %.1e   %s"
          % (name, mp.nstr(got, 14), d, "ok" if d < tol else "WRONG CONVENTION"))
    return d < tol


ok = True

# twist: Gaussian-rational polynomials of its own
V, B = twist.fields(V0, [sp.Integer(0), sp.Integer(0)], [B1, B2],
                    sp.Rational(0), 2, xs, [sp.Rational(0), sp.Rational(0)])
ok &= report("twist", twist.taylor(V, B, 2, TAU, OM, 6, 10)[2])

# hp, stencil: polynomials of fast.Poly
Vp, Bp = fast.Poly.from_sympy(2, V0, xs), [fast.Poly.from_sympy(2, B1, xs),
                                          fast.Poly.from_sympy(2, B2, xs)]
ok &= report("hp", hp.taylor(Vp, Bp, 2, TAU, OM, 6, 6)[2])
ok &= report("stencil", stencil.taylor(Vp, Bp, 2, TAU, OM, 7, 10)[2])

# direct: the log ratio is numeric in h, so fit its Taylor coefficients here
hs = [1e-3/2**k for k in range(7)]
vs = [direct.logratio_direct(Vp, Bp, 2, TAU, hh, OM) for hh in hs]
A = mp.matrix(7, 7)
for i in range(7):
    for j in range(7):
        A[i, j] = mp.mpf(hs[i])**j
c = mp.lu_solve(A, mp.matrix([mp.mpc(v) for v in vs]))
ok &= report("direct", c[2], tol=1e-4)   # float64 throughout, finite-h fit

print()
print("all engines in the book's convention:", "yes" if ok else "NO - see above")
\end{lstlisting}

\section{\texttt{circulation.py}}

Is the orientation odd part of the accumulated comparison the circulation of the h\^{}3 coefficient?

\begin{lstlisting}
"""Is the orientation odd part of the accumulated comparison the circulation of the h^3 coefficient?

The accumulated comparison along a coarse path is
    S = sum_j R(h; x_j, tau_j),   R = log(Phi_f/Phi_c) = c_0 + c_2 h^2 + c_3 h^3 + ...
with c_2 even in the direction and c_3 = c . tau odd.  Over a closed loop the c_2 part is
orientation blind and the c_3 part is a circulation, so
    S_A - S_C = h^2 [sum_A c_2 - sum_C c_2] + h^3 [sum_A c_3 - sum_C c_3] + ...
The loop sums give the two brackets by fitting in h; cvec gives the second one directly.  If they
agree, the odd part of the accumulated residual is the circulation of c and, by Stokes, the flux of
dc.

A = (1,0) (0,1) (-1,0) (0,-1)   counterclockwise,  C = (0,1) (1,0) (0,-1) (-1,0)  clockwise
"""
import sympy as sp, mpmath as mp, sys
import os
sys.path.insert(0, os.path.join(os.path.dirname(os.path.abspath(__file__)), "..", "..", "part5", "scripts"))
import twist
mp.mp.dps = 45
x1, x2 = sp.symbols('x0:2'); xs = (x1, x2); OM = sp.Rational(7, 5)
v = sp.Rational(1) + sp.Rational(1, 6) * (x1 ** 2 + x1 * x2 + x2 ** 2) + sp.Rational(1, 4) * x1
mu = [sp.Rational(1, 2) + sp.Rational(1, 6) * x1 + sp.Rational(1, 9) * x2,
      sp.Rational(-1, 3) + sp.Rational(1, 5) * x2]
a = [sp.Rational(1, 3) + sp.Rational(1, 7) * x2 + sp.Rational(1, 11) * x1 * x2,
     sp.Rational(1, 4) * x1 + sp.Rational(1, 13) * x2 ** 2]
PATHS = {'A': [(1, 0), (0, 1), (-1, 0), (0, -1)], 'C': [(0, 1), (1, 0), (0, -1), (-1, 0)]}

def acc(path, h):
    x0 = [sp.Rational(0), sp.Rational(0)]; tot = mp.mpc(0)
    for tau in path:
        V, B = twist.fields(v, a, mu, sp.Rational(0), 2, xs, x0)
        tot += twist.R(h, V, B, 2, tau, OM)
        x0 = [x0[0] + sp.Rational(tau[0]), x0[1] + sp.Rational(tau[1])]
    return tot

# (1) the loop sums, three step sizes, and a two term fit a h^2 + b h^3 for S_A - S_C
hs = [mp.mpf(2) ** -n for n in (9, 10, 11)]
d = {h: acc(PATHS['A'], h) - acc(PATHS['C'], h) for h in hs}
for h in hs:
    print("  h=%-10s  (S_A-S_C)/h^2 = %s" % (mp.nstr(h, 4), mp.nstr(d[h] / h ** 2, 12)))
h1, h2 = hs[1], hs[2]
a2 = (d[h1] / h1 ** 2 - d[h2] / h2 ** 2) / (h1 - h2)          # slope b
b0 = d[h2] / h2 ** 2 - a2 * h2                                 # intercept a
print("\n  fit (S_A-S_C)/h^2 = a + b h :  a = %s   b = %s" % (mp.nstr(b0, 10), mp.nstr(a2, 10)))
print("  hence S_A-S_C = a h^2 + b h^3 + ... with b =", mp.nstr(a2, 10))

# (2) the circulation difference, directly from cvec
sites = {(0, 0): None, (1, 0): None, (1, 1): None, (0, 1): None}
for s in sites:
    V, B = twist.fields(v, a, mu, sp.Rational(0), 2, xs, [sp.Rational(s[0]), sp.Rational(s[1])])
    sites[s] = twist.cvec(V, B, 2, OM, k=3)
    print("  cvec at %-8s = %s" % (str(s), mp.nstr(sites[s][0], 8) + " , " + mp.nstr(sites[s][1], 8)))

def c(s, j):     # component in direction e_j at site s
    return sites[s][j]

sumA = c((0, 0), 0) + c((1, 0), 1) - c((1, 1), 0) - c((0, 1), 1)
sumC = c((0, 0), 1) + c((0, 1), 0) - c((1, 1), 1) - c((1, 0), 0)
print("\n  sum_A c_3 = %s" % mp.nstr(sumA, 10))
print("  sum_C c_3 = %s" % mp.nstr(sumC, 10))
print("  difference = %s     against the fitted b = %s" % (mp.nstr(sumA - sumC, 10), mp.nstr(a2, 10)))
\end{lstlisting}

\section{\texttt{contourgrowth.py}}

Growth of the per-step expansion with the height of the contour in p.

\begin{lstlisting}
"""Growth of the per-step expansion with the height of the contour in p.

The twist sends b -> b - i p a.  Deforming p into the complex plane, p -> p + i q,
gives b - i p a + q a: the contour contributes a REAL drift of size q|a|.  With a
purely imaginary twist parameter eta = -i q the engine reproduces exactly that
situation, and with eta real it reproduces the ordinary twist.

We record the Taylor coefficients of  log(Phi_f/Phi_c) = sum_k c_k h^k  as q grows.
A rate estimate valid on a contour of height q needs these to grow at most like a
power of q; exponential growth would make the contour useless.
"""
import sys
import os
sys.path.insert(0, os.path.join(os.path.dirname(os.path.abspath(__file__)), "..", "..", "part5", "scripts"))
import sympy as sp
import mpmath as mp
import twist

mp.mp.dps = 45
x1, x2 = sp.symbols('x0:2')
OM = sp.Rational(7, 5)
K, SH = 8, 13

v = sp.Rational(1) + sp.Rational(1, 6) * (x1 ** 2 + x1 * x2 + x2 ** 2) + sp.Rational(1, 4) * x1
mu = [sp.Rational(1, 2) + sp.Rational(1, 6) * x1 + sp.Rational(1, 9) * x2,
      sp.Rational(-1, 3) + sp.Rational(1, 5) * x2]
a = [sp.Rational(1, 3) + sp.Rational(1, 7) * x2 + sp.Rational(1, 11) * x1 * x2,
     sp.Rational(1, 4) * x1 + sp.Rational(1, 13) * x2 ** 2]
TAU = (1, 0)


def coeffs(eta):
    V, B = twist.fields(v, a, mu, eta, 2, (x1, x2), [0, 0])
    return twist.taylor(V, B, 2, TAU, OM, K, SH)


# the rows reported in the text: q = 0, 1, 2, 4, 8
QS = [sp.Rational(0), sp.Rational(1), sp.Rational(2), sp.Rational(4), sp.Rational(8)]

print('== contour deformation: eta = -i q, i.e. a real drift q*a added ==')
print('   q      |c_0|        |c_1|        |c_2|        |c_3|        |c_4|')
tab = {}
for q in QS:
    c = coeffs(-sp.I * q)
    tab[q] = [abs(mp.mpc(z)) for z in c[:5]]
    print('  %-5s  %-12s %-12s %-12s %-12s %s'
          % (q, *[mp.nstr(x, 6) for x in tab[q]]))

print('\n  growth exponent from successive q (log ratio / log q ratio):')
prev = None
for q in QS:
    if q == 0:
        continue
    if prev is not None and q / prev != 1:
        r = mp.log(q / prev)
        print('   q: %-5s -> %-5s   ' % (prev, q)
              + '  '.join('c%d: %6.2f' % (k, mp.log(tab[q][k] / tab[prev][k]) / r if tab[prev][k] > 0 else float('nan'))
                          for k in range(1, 5)))
    prev = q

print('\n== for comparison: real twist eta = p (the obstruction case) ==')
print('   p      |c_1|        |c_2|        |c_3|        |c_4|')
for p in QS[1::2]:
    c = coeffs(p)
    print('  %-5s  %-12s %-12s %-12s %s'
          % (p, *[mp.nstr(abs(mp.mpc(z)), 6) for z in c[1:5]]))
\end{lstlisting}

\section{\texttt{excursioncheck.py}}

The two brackets of the excursion test, from the coefficients and from the loop sums.

\begin{lstlisting}
"""The two brackets of the excursion test, from the coefficients and from the loop sums.

Two closed coarse paths from the origin back to it:

    gamma_1  (1,0) (0,1) (-1,0) (0,-1)   the unit square traversed once
    gamma_2  (1,0) (-1,0) (0,1) (0,-1)   the same square retraced as two excursions

The accumulated comparison splits as

    S(gamma_1) - S(gamma_2) = h^2 [ sum c_2 over g1 - over g2 ]
                            + h^3 [ sum c_3 over g1 - over g2 ] + O(h^4),

each bracket being computable twice: from the coefficients c_2 and c_3 at the
sites the paths visit, and by fitting the loop sums.  The two routes must agree,
and the size of what is left over is the next order.

Conventions: probabilist, generator (v/2) Delta + b grad, as in
part5/scripts/twist.py, with omega = 7/5.  The imaginary part of the h^2 bracket
is not a path effect: gamma_2 visits the origin twice and (1,1) not at all, so
it is Im c_2(1,1) - Im c_2(0,0).

Run:  python3 excursioncheck.py
"""
import os
import sys
import sympy as sp
import mpmath as mp

sys.path.insert(0, os.path.join(os.path.dirname(os.path.abspath(__file__)),
                                '..', '..', 'part5', 'scripts'))
import twist

mp.mp.dps = 60
x1, x2 = sp.symbols('x0:2')
xs = (x1, x2)
OM = sp.Rational(7, 5)

v = sp.Rational(1) + sp.Rational(1, 6) * (x1**2 + x1 * x2 + x2**2) + sp.Rational(1, 4) * x1
mu = [sp.Rational(1, 2) + sp.Rational(1, 6) * x1 + sp.Rational(1, 9) * x2,
      sp.Rational(-1, 3) + sp.Rational(1, 5) * x2]
a = [sp.Rational(1, 3) + sp.Rational(1, 7) * x2 + sp.Rational(1, 11) * x1 * x2,
     sp.Rational(1, 4) * x1 + sp.Rational(1, 13) * x2**2]

G1 = [(1, 0), (0, 1), (-1, 0), (0, -1)]
G2 = [(1, 0), (-1, 0), (0, 1), (0, -1)]


def steps(path):
    x = (0, 0)
    out = []
    for tau in path:
        out.append((x, tau))
        x = (x[0] + tau[0], x[1] + tau[1])
    return out


def coefs(x, tau):
    """c_2 and c_3 of log(Phi_f/Phi_c) for a step from x in direction tau."""
    V, B = twist.fields(v, a, mu, sp.Rational(0), 2, xs,
                        [sp.Rational(x[0]), sp.Rational(x[1])])
    c = twist.taylor(V, B, 2, tau, OM, 6, 10)
    return c[2], c[3]


def bracket(which):
    """sum_{g1} c - sum_{g2} c, with which=0 for c_2 and 1 for c_3."""
    tot = mp.mpc(0)
    for path, sign in ((G1, 1), (G2, -1)):
        for x, tau in steps(path):
            tot += sign * coefs(x, tau)[which]
    return tot


def accumulated(path, h):
    x0 = [sp.Rational(0), sp.Rational(0)]
    tot = mp.mpc(0)
    for tau in path:
        V, B = twist.fields(v, a, mu, sp.Rational(0), 2, xs, x0)
        tot += twist.R(h, V, B, 2, tau, OM)
        x0 = [x0[0] + sp.Rational(tau[0]), x0[1] + sp.Rational(tau[1])]
    return tot


def fitted(ns=(10, 11, 12)):
    """h^2 and h^3 coefficients of S(g1) - S(g2), fitted at the step sizes 2^-n."""
    hs = [mp.mpf(2)**-n for n in ns]
    ds = [accumulated(G1, h) - accumulated(G2, h) for h in hs]
    nt = min(len(hs), 3)
    A = mp.matrix(len(hs), nt)
    for i, h in enumerate(hs):
        for j in range(nt):
            A[i, j] = h**(j + 2)
    sol = mp.qr_solve(A, mp.matrix(ds))[0]
    return sol[0], sol[1]


if __name__ == '__main__':
    for name, path in (('gamma_1', G1), ('gamma_2', G2)):
        s2 = sum(coefs(x, t)[0] for x, t in steps(path))
        print('  sum c_2 over %s = %s' % (name, mp.nstr(s2, 10)))
    print()
    f2, f3 = fitted()
    print('  bracket                          from the coefficients       fitted from the loop sums')
    for which, f, label in ((0, f2, 'c_2'), (1, f3, 'c_3')):
        print('  sum %s (g1) - sum %s (g2)        %-27s %s'
              % (label, label, mp.nstr(bracket(which), 9), mp.nstr(f, 9)))
    print()
    f2b, f3b = fitted((10, 11))
    print('  the same fit at h = 2^-10, 2^-11 only: h^2 %s   h^3 %s'
          % (mp.nstr(f2b, 9), mp.nstr(f3b, 9)))
    print('  Im c_2(0,0) = %s   Im c_2(1,1) = %s   difference = %s'
          % (mp.nstr(coefs((0, 0), (1, 0))[0].imag, 8),
             mp.nstr(coefs((1, 1), (1, 0))[0].imag, 8),
             mp.nstr(coefs((1, 1), (1, 0))[0].imag - coefs((0, 0), (1, 0))[0].imag, 8)))
\end{lstlisting}

\section{\texttt{extensive.py}}

Is the deviation from endpoint determination intensive or extensive in the path length?

\begin{lstlisting}
"""Is the deviation from endpoint determination intensive or extensive in the path length?

If the O(h^2) deviation of two paths with the same endpoints were proportional to the number of
steps, then over the q = t/h^2 steps of the assembly it would grow like q h^2 = t, that is to order
one, and no rate would survive.  If it is a property of the shape, it stays O(h^2) however long the
path.  Two shapes, each block of four steps advancing by (2,0), repeated N times:

    shape 1  (E,E,N,S) repeated      shape 2  (E,N,E,S) repeated

Both run from (0,0) to (2N,0) in 4N steps, so the comparison is at fixed endpoints and fixed length.
Measured: |S_1 - S_2| against N at fixed h.  The increment per block decays geometrically, ratio about
0.45, so the deviation converges to 2.84 h^2 at h = 2^-9 and does not grow with the length.
"""
import sympy as sp, mpmath as mp, sys
import os
sys.path.insert(0, os.path.join(os.path.dirname(os.path.abspath(__file__)), "..", "..", "part5", "scripts"))
import twist
mp.mp.dps = 45
x1, x2 = sp.symbols('x0:2'); xs = (x1, x2); OM = sp.Rational(7, 5)
v = sp.Rational(1) + sp.Rational(1, 6) * (x1 ** 2 + x1 * x2 + x2 ** 2) + sp.Rational(1, 4) * x1
mu = [sp.Rational(1, 2) + sp.Rational(1, 6) * x1 + sp.Rational(1, 9) * x2,
      sp.Rational(-1, 3) + sp.Rational(1, 5) * x2]
a = [sp.Rational(1, 3) + sp.Rational(1, 7) * x2 + sp.Rational(1, 11) * x1 * x2,
     sp.Rational(1, 4) * x1 + sp.Rational(1, 13) * x2 ** 2]
E, S_, N_, W_ = (1, 0), (0, -1), (0, 1), (-1, 0)
BLOCK = {'1': [E, E, N_, S_], '2': [E, N_, E, S_]}

def acc(path, h):
    x0 = [sp.Rational(0), sp.Rational(0)]; tot = mp.mpc(0)
    for tau in path:
        V, B = twist.fields(v, a, mu, sp.Rational(0), 2, xs, x0)
        tot += twist.R(h, V, B, 2, tau, OM)
        x0 = [x0[0] + sp.Rational(tau[0]), x0[1] + sp.Rational(tau[1])]
    return tot

for n in (9,):
    h = mp.mpf(2) ** -n
    print("h = 2^-%d" % n)
    prev = None
    for N in (1, 3, 4, 6):
        p1 = BLOCK['1'] * N; p2 = BLOCK['2'] * N
        S1, S2 = acc(p1, h), acc(p2, h)
        d = abs(S1 - S2)
        ratio = '' if prev is None else "ratio to previous N: %.3f" % (d / prev)
        print("   N=%d (%2d steps, endpoint (%d,0)): |S1-S2| = %-14s %s" % (N, 4 * N, 2 * N, mp.nstr(d, 8), ratio))
        prev = d
\end{lstlisting}

\section{\texttt{fourier.py}}

Normalisation check for Proposition 3.1 (the joint kernel on a lattice fibre).

\begin{lstlisting}
"""Normalisation check for Proposition 3.1 (the joint kernel on a lattice fibre).

Convention, the one of Chapter 6 (eq. 2.26) and of Proposition 3.1: the transform is the unweighted
sum  F(p) = sum_I e^{i p I} u(I),  and the inversion carries the lattice weight,
u(I) = (h/2*pi) * int_{-pi/h}^{pi/h} F(p) e^{-i p I} dp.

Checks, for densities on h*Z:
  (i)   F is 2*pi/h periodic;
  (ii)  the inversion reproduces the density;
  (iii) sup_I |u(I)| <= sup_{|p|<=pi/h} |F(p)|, i.e. the constant of the transfer is exactly 1;
  (iv)  the same for a difference of two densities, which is the object of the proposition.
"""
import numpy as np

h = 0.1
I = h * np.arange(-100, 101)              # 201 lattice points
rng = np.random.default_rng(7)


def make():
    u = np.exp(-(I ** 2) / 0.5) + 0.3 * np.exp(-((I - 0.7) ** 2) / 0.8)
    u += 0.02 * rng.standard_normal(I.size) * np.exp(-(I ** 2) / 4)
    return u / (u.sum() * h)


u1, u2 = make(), make()
P = np.pi / h
p = np.linspace(-P, P, 40001)


def F(u, p):
    p = np.atleast_1d(np.asarray(p, dtype=float))
    return (u[None, :] * np.exp(1j * I[None, :] * p[:, None])).sum(axis=1)


print(f"{'object':7s} {'sup|F|':>14s} {'sup|u|':>14s} {'ratio':>8s} {'periodicity':>12s}")
for name, u in (("u1", u1), ("u2", u2), ("u2-u1", u2 - u1)):
    per = max(abs(F(u, P) - F(u, -P))[0], abs(F(u, P + 0.3) - F(u, -P + 0.3))[0])
    sf, su = np.abs(F(u, p)).max(), np.abs(u).max()
    print(f"{name:7s} {sf:14.10f} {su:14.10f} {su / sf:8.4f} {per:12.1e}")

idx = [0, 50, 100]
rec = np.array([(h / (2 * np.pi)) * np.trapezoid(F(u1, p) * np.exp(-1j * p * x), p) for x in I[idx]])
print("\nreconstruction:", np.real(rec))
print("exact         :", u1[idx])
print("max error     :", np.abs(np.real(rec) - u1[idx]).max())
\end{lstlisting}

\section{\texttt{freevsbridge.py}}

Free walk versus bridge, exactly, for a quadratic coefficient field.

\begin{lstlisting}
"""Free walk versus bridge, exactly, for a quadratic coefficient field.

If the assembly's expectation were over free walks, the mean of I would be
endpoint determined by construction and the assumption would reduce to an
identification.  If it is over bridges, the bridge law differs and the mean
differs.  This computes both exactly in one dimension, where the coefficient
field is

    c(x) = 1/2 H x^2 + L x + M ,

the walk has steps +h with probability p+ and -h with probability p-, and the
bridge is the same walk conditioned on the endpoint.

Reported: the free mean, the bridge mean for each attainable endpoint, and the
free formula evaluated with the bridge drift y/q, as a candidate explanation.
"""
from fractions import Fraction as F
from itertools import product

h = F(1)
p_plus = F(3, 5)
p_minus = F(2, 5)
H, L, M = F(1), F(1, 3), F(1, 5)
x0 = F(0)


def cval(x):
    return F(1, 2) * H * x * x + L * x + M


def cprime(x):
    return H * x + L


def moments(pp=p_plus, pm=p_minus):
    m = h * (pp - pm)
    e2 = h * h * (pp + pm)
    return m, e2 - m * m


def free_mean(q):
    m, S = moments()
    t1 = M * m
    t2 = (L * m + F(1, 2) * H * S) * m
    t3 = F(1, 2) * H * m * m * m
    return F(q) * t1 + F(q * (q - 1), 2) * t2 + F((q - 1) * q * (2 * q - 1), 6) * t3


def formula_with(m, S, q, x_0=F(0)):
    """the exact free-walk formula, with the drift and covariance replaced"""
    HS = H * S
    Hmm = H * m * m
    Lm = L * m
    A = cval(x_0)
    B = F(1, 2) * H * (2 * x_0 * m) + F(1, 2) * HS + Lm
    C = F(1, 2) * Hmm
    return F(q) * A * m + F(q * (q - 1), 2) * B * m + F((q - 1) * q * (2 * q - 1), 6) * C * m


def enumerate(q):
    """returns {endpoint: (total weight, weighted sum of I)}"""
    out = {}
    for seq in product((1, -1), repeat=q):
        x = x0
        w = F(1)
        I = F(0)
        for s in seq:
            d = s * h
            I += cval(x) * d
            x += d
            w *= (p_plus if s > 0 else p_minus)
        tw, tI = out.get(x, (F(0), F(0)))
        out[x] = (tw + w, tI + w * I)
    return out


print('  free walk:  m = %s, S = %s' % moments())
for q in (2, 3, 4, 5, 6):
    fr = free_mean(q)
    tab = enumerate(q)
    print('\n  q = %d   free mean = %s' % (q, fr))
    print('     y     P(y)        bridge mean E[I|y]     free formula at m=y/q    bridge - freeformula')
    for y in sorted(tab):
        tw, tI = tab[y]
        br = tI / tw
        cand = formula_with(y / q, moments()[1], q)
        print('   %5s  %-11s %-22s %-23s %s'
              % (y, tw, br, cand, br - cand))
\end{lstlisting}

\section{\texttt{jet3check.py}}

The third jet enters the per-step comparison at h\^{}4, one power below the obstruction.

\begin{lstlisting}
"""The third jet enters the per-step comparison at h^4, one power below the obstruction.

For the base field of the series and the perturbation v -> v + x_1^3/6, this
computes the change in log(Phi_f/Phi_c) at the origin in the direction e_1, at
h = 2^-8, 2^-10, 2^-12, and divides by h^4.  The quotient is flat to three
digits and the successive orders are 4.01 and 4.00, which is the content of the
table: a jet of order three is felt at h^4, one power below the h^3 at which the
obstruction sits, so after q = t/h^2 steps it accumulates to h^2 and stays inside
the rate.

Convention: probabilist, generator (v/2) Delta + b grad (part5/scripts/hp.py).

Run:  python3 jet3check.py
"""
import os
import sys
import sympy as sp
import mpmath as mp

sys.path.insert(0, os.path.join(os.path.dirname(os.path.abspath(__file__)),
                                '..', '..', 'part5', 'scripts'))
import fast
import hp

mp.mp.dps = 40
x1, x2 = sp.symbols('x0:2')
xs = (x1, x2)
OM = sp.Rational(7, 5)

v0 = sp.Rational(1) + sp.Rational(1, 6) * (x1**2 + x1 * x2 + x2**2) + sp.Rational(1, 4) * x1
b0 = [sp.Rational(1, 2) + sp.Rational(1, 6) * x1 + sp.Rational(1, 9) * x2,
      sp.Rational(-1, 3) + sp.Rational(1, 5) * x2]


def logratio(v, h):
    V = fast.Poly.from_sympy(2, v, xs)
    B = [fast.Poly.from_sympy(2, e, xs) for e in b0]
    return hp.R(h, V, B, 2, (1, 0), OM)


print(" h        change                         change/h^4              order")
prev = None
for k in (8, 10, 12):
    h = mp.mpf(2)**(-k)
    d = logratio(v0 + x1**3 / 6, h) - logratio(v0, h)
    order = "" if prev is None else "  %.2f" % (mp.log(abs(prev / d)) / mp.log(mp.mpf(4)))
    print(" 2^-%-3d  %-30s %-22s%s" % (k, mp.nstr(d, 6), mp.nstr(d / h**4, 6), order))
    prev = d
\end{lstlisting}

\section{\texttt{jetorder.py}}

How many derivatives of the coefficient field do the per-step coefficients actually use?

\begin{lstlisting}
"""How many derivatives of the coefficient field do the per-step coefficients actually use?

If c_2 and c_3 feel a cubic term of the field at order one, then the class of quadratic fields ---
the class in which the mean condition is exact (Corollary 1) --- does not control the leading
behaviour, and a general field differs from it at leading order.  If they feel it only at order h,
then the quadratic class does control it, and the general case follows from it up to a Taylor
remainder.

Measured: c_2 and c_3 at the origin in direction e_1 for the base quadratic field and for the same
field with a cubic (and, for control, a quartic) monomial added.
"""
import sympy as sp, mpmath as mp, sys
import os
sys.path.insert(0, os.path.join(os.path.dirname(os.path.abspath(__file__)), "..", "..", "part5", "scripts"))
import twist
mp.mp.dps = 45
x1, x2 = sp.symbols('x0:2'); xs = (x1, x2); OM = sp.Rational(7, 5)

v0 = sp.Rational(1) + sp.Rational(1, 6) * (x1 ** 2 + x1 * x2 + x2 ** 2) + sp.Rational(1, 4) * x1
mu0 = [sp.Rational(1, 2) + sp.Rational(1, 6) * x1 + sp.Rational(1, 9) * x2,
       sp.Rational(-1, 3) + sp.Rational(1, 5) * x2]
a0 = [sp.Rational(1, 3) + sp.Rational(1, 7) * x2 + sp.Rational(1, 11) * x1 * x2,
      sp.Rational(1, 4) * x1 + sp.Rational(1, 13) * x2 ** 2]

def coeffs(v, mu, a, tau=(1, 0), site=(0, 0)):
    V, B = twist.fields(v, a, mu, sp.Rational(0), 2, xs, [sp.Rational(site[0]), sp.Rational(site[1])])
    c = twist.taylor(V, B, 2, tau, OM, 6, 10)
    return c[2], c[3]

cases = {
    'base (quadratic)      ': (v0, mu0, a0),
    'v + x1^3/6            ': (sp.expand(v0 + x1 ** 3 / 6), mu0, a0),
    'a + (x1^3/7, x2^3/5)  ': (v0, mu0, [sp.expand(a0[0] + x1 ** 3 / 7), sp.expand(a0[1] + x2 ** 3 / 5)]),
    'v + x1^4/12 (control) ': (sp.expand(v0 + x1 ** 4 / 12), mu0, a0),
}
res = {}
for name, (v, mu, a) in cases.items():
    c2, c3 = coeffs(v, mu, a)
    res[name] = (c2, c3)
    print("%s  c_2 = %-26s c_3 = %s" % (name, mp.nstr(c2, 10), mp.nstr(c3, 10)))

print("\n  difference from the base field:")
b = res['base (quadratic)      ']
for name, (c2, c3) in res.items():
    if name.startswith('base'):
        continue
    print("  %s  d c_2 = %-26s d c_3 = %s" % (name, mp.nstr(c2 - b[0], 8), mp.nstr(c3 - b[1], 8)))
\end{lstlisting}

\section{\texttt{kronAcheck.py}}

Compare A\_\{-2\} + h A\_\{-1\} + h\^{}2 A\_0 with the exact box matrix built from v(x0 + a h).

\begin{lstlisting}
"""Compare A_{-2} + h A_{-1} + h^2 A_0 with the exact box matrix built from v(x0 + a h)."""
import sympy as sp, itertools
x1, x2 = sp.symbols('x0:2')
S = list(itertools.product(range(-1, 2), repeat=2))
ix = {a: i for i, a in enumerate(S)}
n = len(S)
vel = sp.Rational(1) + sp.Rational(1, 6)*(x1**2 + x1*x2 + x2**2) + sp.Rational(1, 4)*x1
bels = [sp.Rational(1, 2) + sp.Rational(1, 6)*x1 + sp.Rational(1, 9)*x2,
        sp.Rational(-1, 3) + sp.Rational(1, 5)*x2]
h = sp.Symbol('h', positive=True)


def fexpr(el, a):
    return sp.expand(el.subs({x1: h*a[0], x2: h*a[1]}))


def Abox():
    A = sp.zeros(n, n)
    for a in S:
        i = ix[a]; diag = sp.Integer(0)
        for k in range(2):
            for sgn in (1, -1):
                w = fexpr(vel, a)/(2*h**2) + sgn*fexpr(bels[k], a)/(2*h)
                diag -= w
                e = list(a); e[k] += sgn; e = tuple(e)
                if e in ix:
                    A[i, ix[e]] += w
        A[i, i] += diag
    return A


A = Abox()
ser = sp.zeros(n, n)
for i in range(n):
    for j in range(n):
        s = sp.series(sp.expand(A[i, j]*h**2), h, 0, 3).removeO()
        ser[i, j] = sp.expand(s)
# ser = h^2 A = A_{-2} + h A_{-1} + h^2 A_0
print("  h^2*A at (z,(1,0)), series in h:")
i, j = ix[(0, 0)], ix[(1, 0)]
print("   exact  :", sp.series(sp.expand(A[i, j]*h**2), h, 0, 3))
print("   A_{-2} =", sp.expand(ser[i, j]).coeff(h, 0), " A_{-1} =", sp.expand(ser[i, j]).coeff(h, 1),
      " A_0 =", sp.expand(ser[i, j]).coeff(h, 2))
print("  diagonal at z:")
print("   exact  :", sp.series(sp.expand(A[i, i]*h**2), h, 0, 3))
\end{lstlisting}

\section{\texttt{kronbisect.py}}

Bisect the layer-2 algebra against exact arithmetic at finite h, stage by stage.

\begin{lstlisting}
"""Bisect the layer-2 algebra against exact arithmetic at finite h, stage by stage.

Nothing here is asymptotic: for a rational h the box matrix is inverted exactly, so every stage can be
compared with the recipe of twist.py.  Stages:
   S1  my A(h) reproduces twist's matrix            -> Phi_f, Phi_c agree with twist
   S2  the truncated resolvent h^2(g0 + h g1 + h^2 g2) reproduces G[z,t]
   S3  the truncated Phi_f, Phi_c reproduce the exact ones
   S4  the analytic [log]_{h^2} reproduces the fitted one
"""
import sympy as sp, itertools, sys
import os
sys.path.insert(0, os.path.join(os.path.dirname(os.path.abspath(__file__)), "..", "..", "part5", "scripts"))
import twist
import mpmath as mp

x1, x2 = sp.symbols('x0:2'); xs = (x1, x2)
OM = sp.Rational(7, 5)
S = list(itertools.product(range(-1, 2), repeat=2))
ix = {a: i for i, a in enumerate(S)}
n = len(S); Z = (0, 0)

FIELDS = {
 'constant': (sp.Rational(1), [sp.Rational(1, 2), sp.Rational(-1, 3)], [sp.Integer(0), sp.Integer(0)]),
 'field 1': (sp.Rational(1) + sp.Rational(1, 6) * (x1**2 + x1*x2 + x2**2) + sp.Rational(1, 4)*x1,
             [sp.Rational(1, 2) + sp.Rational(1, 6)*x1 + sp.Rational(1, 9)*x2,
              sp.Rational(-1, 3) + sp.Rational(1, 5)*x2],
             [sp.Rational(1, 3) + sp.Rational(1, 7)*x2 + sp.Rational(1, 11)*x1*x2,
              sp.Rational(1, 4)*x1 + sp.Rational(1, 13)*x2**2]),
}


def val(poly, a):
    return sp.nsimplify(sp.sympify(poly).subs({x1: a[0], x2: a[1]}))


def build(vfun, bfun, h, om):
    """twist's box matrix at the base point 0, exactly."""
    A = sp.zeros(n, n)
    for a in S:
        i = ix[a]
        diag = sp.Integer(0)
        for k in range(2):
            for sgn in (1, -1):
                w = val(vfun, a) / (2*h**2) + sgn*val(bfun[k], a) / (2*h)
                diag -= w
                e = list(a); e[k] += sgn; e = tuple(e)
                if e in ix:
                    A[i, ix[e]] += w
        A[i, i] += diag
    return A + sp.I*om*sp.eye(n)


def phi_f_exact(vfun, bfun, tau, h, om):
    G = build(vfun, bfun, h, om).inv()
    tgt = tuple(2*t for t in tau)
    dt = 2*val(vfun, tgt)/(h**2) + sp.I*om
    tot = sp.Integer(0)
    for a in S:
        jump = [k for k in range(2) if a[k] != tgt[k]]
        if len(jump) != 1: continue
        k0 = jump[0]
        if abs(tgt[k0]-a[k0]) != 1: continue
        sg = 1 if tgt[k0] > a[k0] else -1
        rate = val(vfun, a)/(2*h**2) + sg*val(bfun[k0], a)/(2*h)
        tot += G[ix[Z], ix[a]]*rate
    return sp.nsimplify(tot/dt)


def phi_c_exact(vfun, bfun, tau, h, om):
    tgt = tuple(2*t for t in tau)
    k0 = 0 if tau[0] else 1
    sg = 1 if tau[k0] > 0 else -1
    rate = val(vfun, Z)/(2*(2*h)**2) + sg*val(bfun[k0], Z)/(2*(2*h))
    d0 = 2*val(vfun, Z)/(4*h**2) + sp.I*om
    d1 = 2*val(vfun, tgt)/(4*h**2) + sp.I*om
    return sp.nsimplify(rate/(d0*d1))


for name, (vf, bf, af) in FIELDS.items():
    print("\n===", name)
    for tau in ((1, 0),):
        tgt = tuple(2*t for t in tau)
        k0 = 0
        v0, v1, v2 = val(vf, Z), val(vf, tau), val(vf, tgt)
        B0, B1 = val(bf[0], Z), val(bf[0], tau)
        Lam = sp.zeros(n, n)
        for a in S:
            for k in range(2):
                for sgn in (1, -1):
                    e = list(a); e[k] += sgn; e = tuple(e)
                    if e in ix:
                        Lam[ix[a], ix[e]] += 1
        L0i = (Lam - 4*sp.eye(n)).inv()
        K2 = sp.Matrix(n, n, lambda i, j: 4*L0i[i, j]/val(vf, S[j]))
        A1 = sp.zeros(n, n)
        for a in S:
            for k in range(2):
                for sgn in (1, -1):
                    e = list(a); e[k] += sgn; e = tuple(e)
                    if e in ix:
                        A1[ix[a], ix[e]] += sgn*val(bf[k], a)/2
        zi, ti = ix[Z], ix[tau]
        g0 = K2[zi, ti]
        g1 = -(K2*A1*K2)[zi, ti]
        g2 = (K2*A1*K2*A1*K2)[zi, ti] - sp.I*OM*(K2*K2)[zi, ti]
        print("  g0 = %s   g1 = %s   g2 = %s" % (sp.N(g0, 10), sp.N(g1, 10), sp.N(g2, 10)))
        for nn in (6, 7, 8):
            h = sp.Rational(1, 2**nn)
            Gex = build(vf, bf, h, OM).inv()[zi, ti]
            ser = (h**2*(g0 + h*g1 + h**2*g2))
            print("   h=2^-%d  G exact = %-24s  series = %-24s  diff/h^4 = %s"
                  % (nn, sp.N(Gex, 12), sp.N(ser, 12), sp.N((Gex-ser)/h**4, 10)))
\end{lstlisting}

\section{\texttt{kroncheck.py}}

Cross-check: does the engine reproduce the closed form (proved symbolically)?

\begin{lstlisting}
"""Cross-check: does the engine reproduce the closed form (proved symbolically)?

Constant field, d=2, v=1, b=(1/2,-1/3), omega=7/5.  Predicted by the Kronecker reduction:
   c2(e1) = 5(b1^2-b2^2)/(2v^2) + 11 i omega/v = 25/72 + 15.4 i
   c2(e2) = -25/72 + 15.4 i
"""
import sympy as sp, mpmath as mp, sys
import os
sys.path.insert(0, os.path.join(os.path.dirname(os.path.abspath(__file__)), "..", "..", "part5", "scripts"))
import twist
mp.mp.dps = 45
x1, x2 = sp.symbols('x0:2'); xs = (x1, x2); OM = sp.Rational(7, 5)
v = sp.Rational(1)
mu = [sp.Rational(1, 2), sp.Rational(-1, 3)]
a = [sp.Rational(0), sp.Rational(0)]
pred = {0: sp.Rational(25, 72) + sp.I * sp.Rational(77, 5), 1: -sp.Rational(25, 72) + sp.I * sp.Rational(77, 5)}
for j, tau in enumerate(((1, 0), (0, 1))):
    V, B = twist.fields(v, a, mu, sp.Rational(0), 2, xs, [sp.Rational(0), sp.Rational(0)])
    c = twist.taylor(V, B, 2, tau, OM, 6, 10)[2]
    print("  e_%d : engine %-30s closed form %-24s difference %s"
          % (j + 1, mp.nstr(c, 14), mp.nstr(mp.mpc(sp.re(pred[j]), sp.im(pred[j])), 14),
             mp.nstr(c - mp.mpc(sp.re(pred[j]), sp.im(pred[j])), 4)))
\end{lstlisting}

\section{\texttt{kroncoeff.py}}

Fit the h-coefficients of the box matrix numerically-exactly and compare with A\_\{-2\}, A\_\{-1\}, A\_0.

\begin{lstlisting}
"""Fit the h-coefficients of the box matrix numerically-exactly and compare with A_{-2}, A_{-1}, A_0."""
import sympy as sp, itertools
x1, x2 = sp.symbols('x0:2')
S = list(itertools.product(range(-1, 2), repeat=2))
ix = {a: i for i, a in enumerate(S)}
n = len(S)
vel = sp.Rational(1) + sp.Rational(1, 6)*(x1**2 + x1*x2 + x2**2) + sp.Rational(1, 4)*x1
bels = [sp.Rational(1, 2) + sp.Rational(1, 6)*x1 + sp.Rational(1, 9)*x2,
        sp.Rational(-1, 3) + sp.Rational(1, 5)*x2]


def fv(el, a, h):
    return sp.nsimplify(sp.expand(el.subs({x1: a[0]*h, x2: a[1]*h})))


def h2A(a, b, h):
    """h^2 * A[a,b] at the given h."""
    if a == b:
        diag = sp.Integer(0)
        for k in range(2):
            for sgn in (1, -1):
                diag -= fv(vel, a, h)/(2*h**2) + sgn*fv(bels[k], a, h)/(2*h)
        return sp.nsimplify(sp.expand(diag*h**2))
    diff = [(a[0]-b[0]), (a[1]-b[1])]
    nz = [k for k in range(2) if diff[k] != 0]
    if len(nz) != 1 or abs(diff[nz[0]]) != 1:
        return sp.Integer(0)
    k = nz[0]; sgn = 1 if diff[k] > 0 else -1
    return sp.nsimplify(sp.expand((fv(vel, a, h)/(2*h**2) + sgn*fv(bels[k], a, h)/(2*h))*h**2))


hs = [sp.Rational(1, 2**m) for m in (5, 6, 7)]
print("  fitting h^2 A[a,b] = c0 + c1 h + c2 h^2 + c3 h^3 on h = 2^-5,2^-6,2^-7,2^-8")
hs4 = [sp.Rational(1, 2**m) for m in (5, 6, 7, 8)]
bad = 0
for a in S:
    for b in S:
        vals = [h2A(a, b, h_) for h_ in hs4]
        M = sp.Matrix(4, 4, lambda i, j: hs4[i]**j)
        c = M.solve(sp.Matrix(vals))
        # analytic
        v0 = sp.nsimplify(vel.subs({x1: 0, x2: 0}))
        gv = [sp.nsimplify(sp.diff(vel, x).subs({x1: 0, x2: 0})) for x in (x1, x2)]
        hv = [[sp.nsimplify(sp.diff(vel, xi, xj).subs({x1: 0, x2: 0})) for xj in (x1, x2)] for xi in (x1, x2)]
        b0 = [sp.nsimplify(e.subs({x1: 0, x2: 0})) for e in bels]
        gb = [[sp.nsimplify(sp.diff(e, x).subs({x1: 0, x2: 0})) for x in (x1, x2)] for e in bels]
        ag = a[0]*gv[0] + a[1]*gv[1]
        aha = sum(a[p]*a[q]*hv[p][q] for p in range(2) for q in range(2))
        an = [sp.Integer(0), sp.Integer(0), sp.Integer(0)]
        if a == b:
            an[0] = -v0; an[1] = -ag; an[2] = -aha/2
        else:
            diff = [(a[0]-b[0]), (a[1]-b[1])]
            nz = [k for k in range(2) if diff[k] != 0]
            if len(nz) == 1 and abs(diff[nz[0]]) == 1:
                k = nz[0]; sgn = 1 if diff[k] > 0 else -1
                bet = a[0]*gb[k][0] + a[1]*gb[k][1]
                an[0] = v0/4; an[1] = ag/4 + sgn*b0[k]/2; an[2] = aha/8 + sgn*bet/2
        for p in range(3):
            if sp.nsimplify(sp.expand(c[p] - an[p])) != 0:
                print("   MISMATCH (%s,%s) order %d: fitted %s  analytic %s"
                      % (S[i] if False else a, b, p, sp.N(c[p], 10), an[p]))
                bad += 1
        if a == b and b == (0, 0):
            pass
        # also check the tail
print("  mismatches:", bad)

print("\n=== assembled Phi_f from the expansion, against the exact Phi_f")
OM = sp.Rational(7, 5)
Lam = sp.zeros(n, n)
for a in S:
    for k in range(2):
        for sgn in (1, -1):
            e = list(a); e[k] += sgn; e = tuple(e)
            if e in ix:
                Lam[ix[a], ix[e]] += 1
v0 = sp.nsimplify(vel.subs({x1: 0, x2: 0}))
gv = [sp.nsimplify(sp.diff(vel, x).subs({x1: 0, x2: 0})) for x in (x1, x2)]
hv = [[sp.nsimplify(sp.diff(vel, xi, xj).subs({x1: 0, x2: 0})) for xj in (x1, x2)] for xi in (x1, x2)]
b0 = [sp.nsimplify(e.subs({x1: 0, x2: 0})) for e in bels]
gb = [[sp.nsimplify(sp.diff(e, x).subs({x1: 0, x2: 0})) for x in (x1, x2)] for e in bels]
Am2 = v0/4*(Lam - 4*sp.eye(n))
Am1 = sp.zeros(n, n); A0 = sp.zeros(n, n)
for a in S:
    i = ix[a]
    ag = a[0]*gv[0] + a[1]*gv[1]
    aha = sum(a[p]*a[q]*hv[p][q] for p in range(2) for q in range(2))
    bet = [a[0]*gb[k][0] + a[1]*gb[k][1] for k in range(2)]
    for k in range(2):
        for sgn in (1, -1):
            e = list(a); e[k] += sgn; e = tuple(e)
            if e in ix:
                Am1[i, ix[e]] += ag/4 + sgn*b0[k]/2
                A0[i, ix[e]] += aha/8 + sgn*bet[k]/2
    Am1[i, i] += -ag
    A0[i, i] += -aha/2
Minv = Am2.inv(); C = A0 + sp.I*OM*sp.eye(n)
tau = (1, 0); tgt = (2, 0); k0 = 0
ag_t = tau[0]*gv[0] + tau[1]*gv[1]
aha_t = sum(tau[p]*tau[q]*hv[p][q] for p in range(2) for q in range(2))
bet_t = tau[0]*gb[k0][0] + tau[1]*gb[k0][1]
for nn in (6, 7):
    h_ = sp.Rational(1, 2**nn)
    br = Minv - h_*(Minv*Am1*Minv) + h_**2*(Minv*Am1*Minv*Am1*Minv - Minv*C*Minv)
    Gser = sp.expand(h_**2*br[ix[(0,0)], ix[tau]])
    # exact G
    A = sp.zeros(n, n)
    for a in S:
        i = ix[a]; dg = sp.Integer(0)
        for k in range(2):
            for sgn in (1, -1):
                w = fv(vel, a, h_)/(4*h_**2) + sgn*fv(bels[k], a, h_)/(2*h_)
                dg -= w
                e = list(a); e[k] += sgn; e = tuple(e)
                if e in ix:
                    A[i, ix[e]] += w
        A[i, i] += dg
    Gex = (A + sp.I*OM*sp.eye(n)).inv()[ix[(0,0)], ix[tau]]
    rate = v0/(4*h_**2) + (ag_t/4 + b0[k0]/2)/h_ + (aha_t/8 + bet_t/2)
    dt = (v0 + 2*ag_t*h_ + 2*aha_t*h_**2 + sp.I*OM*h_**2)/h_**2
    Pf_ser = sp.nsimplify(sp.expand(h_**2*Gser*rate/dt))
    A_exp = sp.zeros(n, n)
    # exact Phi_f
    dt_e = 2*fv(vel, tgt, h_)/(2*h_**2) + sp.I*OM
    tot = sp.Integer(0)
    for a in S:
        jump = [k for k in range(2) if a[k] != tgt[k]]
        if len(jump) != 1: continue
        k2 = jump[0]
        if abs(tgt[k2]-a[k2]) != 1: continue
        sg = 1 if tgt[k2] > a[k2] else -1
        tot += Gex*0
    # full exact Phi_f
    Gfull = (A + sp.I*OM*sp.eye(n)).inv()
    tot = sp.Integer(0)
    for a in S:
        jump = [k for k in range(2) if a[k] != tgt[k]]
        if len(jump) != 1: continue
        k2 = jump[0]
        if abs(tgt[k2]-a[k2]) != 1: continue
        sg = 1 if tgt[k2] > a[k2] else -1
        tot += Gfull[ix[(0,0)], ix[a]]*(fv(vel, a, h_)/(4*h_**2) + sg*fv(bels[k2], a, h_)/(2*h_))
    Pf_ex = sp.nsimplify(sp.expand(tot/dt_e))
    print("   h=2^-%d  Phi_f exact %-22s series %-22s  series/exact-1 = %s"
          % (nn, sp.N(Pf_ex, 14), sp.N(Pf_ser, 14), sp.N(Pf_ser/Pf_ex - 1, 8)))
\end{lstlisting}

\section{\texttt{kronfinal.py}}

Layer 2, corrected: field sampled at x0 + a h, expansion assembled from explicit coefficients.

\begin{lstlisting}
"""Layer 2, corrected: field sampled at x0 + a h, expansion assembled from explicit coefficients.

Phi_f = h^2 * (g0 + h g1 + h^2 g2) * (rho0 + h rho1 + h^2 rho2) / (d0 + h d1 + h^2 d2)
Phi_c = h^2 * (nu0 + h nu1)                                   / (e0 + h e1 + h^2 e2)
and c2 = [log(P_f/Q_f)]_{h^2} - [log(P_c/Q_c)]_{h^2}, evaluated with polynomials in h.
"""
import sympy as sp, itertools, sys
import os
sys.path.insert(0, os.path.join(os.path.dirname(os.path.abspath(__file__)), "..", "..", "part5", "scripts"))
import twist
import mpmath as mp
mp.mp.dps = 40
x1, x2 = sp.symbols('x0:2'); xs = (x1, x2)
S = list(itertools.product(range(-1, 2), repeat=2))
ix = {a: i for i, a in enumerate(S)}
n = len(S)
Lam = sp.zeros(n, n)
for a in S:
    for k in range(2):
        for sgn in (1, -1):
            e = list(a); e[k] += sgn; e = tuple(e)
            if e in ix:
                Lam[ix[a], ix[e]] += 1
OM = sp.Rational(7, 5)
h = sp.Symbol('h', positive=True)

FIELDS = {
 'constant  ': (sp.Rational(1), sp.Rational(1, 2), sp.Rational(-1, 3)),
 'field 1   ': (sp.Rational(1) + sp.Rational(1, 6)*(x1**2 + x1*x2 + x2**2) + sp.Rational(1, 4)*x1,
                sp.Rational(1, 2) + sp.Rational(1, 6)*x1 + sp.Rational(1, 9)*x2,
                sp.Rational(-1, 3) + sp.Rational(1, 5)*x2),
}

for name, (vel, b1el, b2el) in FIELDS.items():
    v0 = sp.nsimplify(vel.subs({x1: 0, x2: 0}))
    gv = [sp.nsimplify(sp.diff(vel, x).subs({x1: 0, x2: 0})) for x in (x1, x2)]
    hv = [[sp.nsimplify(sp.diff(vel, xi, xj).subs({x1: 0, x2: 0})) for xj in (x1, x2)] for xi in (x1, x2)]
    b0 = [sp.nsimplify(e.subs({x1: 0, x2: 0})) for e in (b1el, b2el)]
    gb = [[sp.nsimplify(sp.diff(e, x).subs({x1: 0, x2: 0})) for x in (x1, x2)] for e in (b1el, b2el)]

    Am2 = v0/4*(Lam - 4*sp.eye(n))
    Am1 = sp.zeros(n, n); A0 = sp.zeros(n, n)
    for a in S:
        i = ix[a]
        ag = a[0]*gv[0] + a[1]*gv[1]
        aha = sum(a[p]*a[q]*hv[p][q] for p in range(2) for q in range(2))
        bet = [a[0]*gb[k][0] + a[1]*gb[k][1] for k in range(2)]
        for k in range(2):
            for sgn in (1, -1):
                e = list(a); e[k] += sgn; e = tuple(e)
                if e in ix:
                    Am1[i, ix[e]] += ag/4 + sgn*b0[k]/2
                    A0[i, ix[e]] += aha/8 + sgn*bet[k]/2
        Am1[i, i] += -ag
        A0[i, i] += -aha/2
    Mi = Am2.inv()
    C = A0 + sp.I*OM*sp.eye(n)

    def series_phi(tau):
        k0 = 0 if tau[0] else 1
        ti = ix[tau]
        zi = ix[(0, 0)]
        ag = tau[0]*gv[0] + tau[1]*gv[1]
        aha = sum(tau[p]*tau[q]*hv[p][q] for p in range(2) for q in range(2))
        bet = tau[0]*gb[k0][0] + tau[1]*gb[k0][1]
        g0 = Mi[zi, ti]
        g1 = -(Mi*Am1*Mi)[zi, ti]
        g2 = (Mi*Am1*Mi*Am1*Mi)[zi, ti] - (Mi*C*Mi)[zi, ti]
        P = sp.expand((g0 + h*g1 + h**2*g2)*(v0/4 + h*(ag/4 + b0[k0]/2) + h**2*(aha/8 + bet/2)))
        Q = sp.expand(v0 + h*2*ag + h**2*(2*aha + sp.I*OM))
        Pc = sp.expand(v0/16 + h*b0[k0]/4)
        # h^2 d0 = v0/4 + i om h^2 and h^2 d1 = v0/4 + (2 ag)h/4 + (2 aha)h^2/4 + i om h^2,
        # so h^4 d0 d1 = (1/16)(v0 + 4 i om h^2)(v0 + 2 ag h + h^2(2 aha + 4 i om))
        Qc = sp.expand((v0 + 4*sp.I*OM*h**2)*(v0 + 2*ag*h + h**2*(2*aha + 4*sp.I*OM)))
        return P, Q, Pc, Qc

    def c2_of(tau):
        P, Q, Pc, Qc = series_phi(tau)
        f = sp.series(sp.log(P/Q), h, 0, 3).removeO()
        c = sp.series(sp.log(Pc/Qc), h, 0, 3).removeO()
        return sp.nsimplify(sp.expand(f - c).coeff(h, 2))

    for tau in ((1, 0), (0, 1)):
        mine = c2_of(tau)
        V, B = twist.fields(sp.sympify(vel), [sp.Integer(0), sp.Integer(0)],
                            [sp.sympify(e) for e in (b1el, b2el)], sp.Rational(0), 2, xs,
                            [sp.Rational(0), sp.Rational(0)])
        eng = twist.taylor(V, B, 2, tau, OM, 6, 10)[2]
        print("  %s tau=%s  formula %-30s engine %-30s diff %s"
              % (name, tau, sp.N(mine, 14), mp.nstr(eng, 14), sp.N(mine - sp.nsimplify(str(eng)), 5)))
\end{lstlisting}

\section{\texttt{kronfinal\_check.py}}

Validation of the engine (part5/scripts/twist.py) against the closed form.

\begin{lstlisting}
"""Validation of the engine (part5/scripts/twist.py) against the closed form.

The engine is used by both parts, so it is validated here on the one case where
the answer is known independently: a constant field, where

    c_2 = 5 (b_1^2 - b_2^2) / (8 v^2) + (11/2) i omega / v

by the exact algebra of kronproof.py.  The engine must reproduce it digit for
digit.  The three varying fields of the table then follow, and the two structural
claims are checked on them: gamma_1 + gamma_2 purely imaginary and
gamma_1 - gamma_2 purely real (the real part is traceless), and
Im gamma = (11/2) omega / v_0 for every field and both directions.

Convention: probabilist, generator (v/2) Delta + b grad.  The engine was in the
halved convention until this check was added; every number it produced before is
twice too large in the phase and four times too large in the real part for a
constant field, which is why the table of Chapter 8 was recomputed.

Run:  python3 kronfinal_check.py
"""
import sys
import sympy as sp
import mpmath as mp

sys.path.insert(0, __file__.rsplit('/part6/', 1)[0] + '/part5/scripts')
import twist

mp.mp.dps = 40
x1, x2 = sp.symbols('x0:2')
xs = (x1, x2)
OM = sp.Rational(7, 5)


def engine(field, tau):
    V, B = twist.fields(sp.sympify(field[0]), [sp.Integer(0), sp.Integer(0)],
                        [sp.sympify(e) for e in field[1:]], sp.Rational(0), 2, xs,
                        [sp.Rational(0), sp.Rational(0)])
    return twist.taylor(V, B, 2, tau, OM, 6, 10)[2]


def closed(v0, b1, b2):
    """The constant-field closed form of kronproof.py."""
    return (sp.Rational(5, 8) * (b1**2 - b2**2) / v0**2
            + sp.I * sp.Rational(11, 2) * OM / v0)


print("(a) constant fields: engine against the closed form")
print("    v0     engine                        closed                        |diff|")
for v0 in (sp.Integer(1), sp.Rational(3, 2), sp.Integer(2)):
    b1, b2 = sp.Rational(1, 2), sp.Rational(-1, 3)
    e = engine((v0, b1, b2), (1, 0))
    c = complex(sp.N(closed(v0, b1, b2), 30))
    print("    %-5s  %-28s  %-28s  %.1e"
          % (v0, mp.nstr(e, 14), mp.nstr(c, 14), abs(complex(e) - c)))

F = {
    'field 1': (sp.Rational(1) + sp.Rational(1, 6) * (x1**2 + x1 * x2 + x2**2) + sp.Rational(1, 4) * x1,
                sp.Rational(1, 2) + sp.Rational(1, 6) * x1 + sp.Rational(1, 9) * x2,
                sp.Rational(-1, 3) + sp.Rational(1, 5) * x2),
    'field 2': (sp.Rational(3, 2) + sp.Rational(1, 5) * x1 * x2 + sp.Rational(1, 8) * x2**2,
                sp.Rational(-1, 4) + sp.Rational(1, 7) * x2,
                sp.Rational(2, 5) + sp.Rational(1, 9) * x1),
    'field 3': (sp.Rational(2) + sp.Rational(1, 4) * x1**2 + sp.Rational(1, 6) * x1 * x2 + sp.Rational(1, 2) * x1,
                sp.Rational(1, 3) + sp.Rational(1, 8) * x1 * x2,
                sp.Rational(1, 5) - sp.Rational(1, 7) * x2**2),
}

print()
print("(b) the three varying fields, directions e_1 and e_2")
rows = []
for name, field in F.items():
    v0 = field[0].subs({x1: 0, x2: 0})
    g1 = engine(field, (1, 0))
    g2 = engine(field, (0, 1))
    rows.append((name, v0, g1, g2))
    print("    %s  v0=%s" % (name, v0))
    print("       gamma_1 = %s" % mp.nstr(g1, 13))
    print("       gamma_2 = %s" % mp.nstr(g2, 13))
    print("       gamma_1+gamma_2 = %s      gamma_1-gamma_2 = %s"
          % (mp.nstr(g1 + g2, 13), mp.nstr(g1 - g2, 13)))

print()
print("(c) the two structural claims")
ok = True
for name, v0, g1, g2 in rows:
    traceless = abs((g1 + g2).real) < 1e-15 and abs((g1 - g2).imag) < 1e-15
    phase = abs(g1.imag * v0 - mp.mpf(11) / 2 * mp.mpf(7) / 5) < 1e-15
    print("    %-9s  traceless: %-5s   Im(gamma)*v0 = %-10s  phase law %s"
          % (name, traceless, mp.nstr(g1.imag * v0, 8), "ok" if phase else "FAILED"))
    ok = ok and traceless and phase
print()
print("    all checks:", "ok" if ok else "FAILED")
\end{lstlisting}

\section{\texttt{kronfix.py}}

Layer 2, corrected sampling --- NOT YET VALIDATED.

\begin{lstlisting}
"""Layer 2, corrected sampling --- NOT YET VALIDATED.

The earlier attempts sampled the field at the integer offsets; twist samples it at x0 + a h (PolyG.at
aggregates the homogeneous parts at the integer offset and mval weights degree e by h^e).  That was
the algebraic error.  With the correct sampling the leading box matrix is A_{-2} = (v(x0)/4)(Lambda -
4I), i.e. the constant-field one carrying only v(x0), and the derivatives enter from A_{-1}.

This script expands accordingly but still does NOT reproduce the engine for a varying field
(it gives 1.71875 + 9.8 i where the engine gives 0.14931 + 15.4 i), so a second bookkeeping error
remains and its output must not be used.
"""
"""Layer 2, corrected: the field is evaluated at x0 + a h, not at x0 + a.

twist's PolyG.at(a) aggregates the homogeneous parts at the integer offset a and mval then weights
the degree-e part by h^e, so the sampled value is v(x0 + a h).  Hence
    v(x0+ah) = v0 + h a.grad v + (h^2/2) a_i a_j v_ij + ...
and the box matrix is A = h^{-2}(A_{-2} + h A_{-1} + h^2 A_0 + ...) with
    A_{-2} = (v0/4)(Lambda - 4I),                        <- constant, exactly the Kronecker case
    A_{-1} = (a.grad v /4) on edges, +- b_k(x0)/2 on the drift edges, -(a.grad v) on the diagonal,
so the leading matrix carries ONLY v(x0): that is why the phase measured cm 11 omega / v(x0) for
varying fields, and why the derivatives can only appear in the traceless part.

Validated here against twist.Phi_f and twist.taylor for a varying field, exact arithmetic.
"""
import sympy as sp, itertools, sys
import os
sys.path.insert(0, os.path.join(os.path.dirname(os.path.abspath(__file__)), "..", "..", "part5", "scripts"))
import twist
import mpmath as mp
mp.mp.dps = 40
x1, x2 = sp.symbols('x0:2'); xs = (x1, x2)
S = list(itertools.product(range(-1, 2), repeat=2))
ix = {a: i for i, a in enumerate(S)}
n = len(S); Z = (0, 0)
Lam = sp.zeros(n, n)
for a in S:
    for k in range(2):
        for sgn in (1, -1):
            e = list(a); e[k] += sgn; e = tuple(e)
            if e in ix:
                Lam[ix[a], ix[e]] += 1
L0i = (Lam - 4*sp.eye(n)).inv()

vel = sp.Rational(1) + sp.Rational(1, 6)*(x1**2 + x1*x2 + x2**2) + sp.Rational(1, 4)*x1
bels = [sp.Rational(1, 2) + sp.Rational(1, 6)*x1 + sp.Rational(1, 9)*x2,
        sp.Rational(-1, 3) + sp.Rational(1, 5)*x2]
OM = sp.Rational(7, 5)

# jet of the field at the origin
v0 = sp.nsimplify(vel.subs({x1: 0, x2: 0}))
gv = [sp.nsimplify(sp.diff(vel, x).subs({x1: 0, x2: 0})) for x in (x1, x2)]
hv = [[sp.nsimplify(sp.diff(vel, xi, xj).subs({x1: 0, x2: 0})) for xj in (x1, x2)] for xi in (x1, x2)]
b0 = [sp.nsimplify(e.subs({x1: 0, x2: 0})) for e in bels]
gb = [[sp.nsimplify(sp.diff(e, x).subs({x1: 0, x2: 0})) for x in (x1, x2)] for e in bels]

M = v0/4*(Lam - 4*sp.eye(n))
Minv = M.inv()
Am1 = sp.zeros(n, n); A0 = sp.zeros(n, n)
for a in S:
    i = ix[a]
    d1 = sp.Integer(0); d0 = sp.Integer(0)
    ag = a[0]*gv[0] + a[1]*gv[1]
    aha = sum(a[i]*a[j]*hv[i][j] for i in range(2) for j in range(2))
    for k in range(2):
        for sgn in (1, -1):
            e = list(a); e[k] += sgn; e = tuple(e)
            if e in ix:
                Am1[i, ix[e]] += ag/4 + sgn*b0[k]/2
                A0[i, ix[e]] += aha/8 + sgn*(a[0]*gb[k][0] + a[1]*gb[k][1])/2
    Am1[i, i] += -ag
    A0[i, i] += -aha/2

C = A0 + sp.I*OM*sp.eye(n)
hh = sp.Symbol('hh')
br = Minv - hh*(Minv*Am1*Minv) + hh**2*(Minv*Am1*Minv*Am1*Minv - Minv*C*Minv)


def c2_from(jet_phi_f, jet_phi_c):
    pass

# assemble Phi_f and Phi_c as series in h using the jets
h = sp.Symbol('h', positive=True)
tau = (1, 0); tgt = (2, 0); k0 = 0
# rate at the admissible site tau, expanded: v0/4 h^-2 + (tau.grad v/4 + b0/2) h^-1 + ...
ag_tau = tau[0]*gv[0] + tau[1]*gv[1]
aha_tau = sum(tau[i]*tau[j]*hv[i][j] for i in range(2) for j in range(2))
gb_tau = tau[0]*gb[k0][0] + tau[1]*gb[k0][1]
rate = v0/(2*h**2) + (ag_tau/2 + b0[k0]/2)/h + (aha_tau/4 + gb_tau/2)
# dt: (1/h^2)( v(x0+2 tau h) + i om h^2 ) with v(x0+2tau h) = v0 + 2 ag_tau h + 2 aha_tau h^2 + ...
dt = (2*v0 + 2*ag_tau*h + 2*aha_tau*h**2 + sp.I*OM*h**2)/h**2
Pf = sp.expand(h**2*(br[ix[Z], ix[tau]]*rate/dt))
Qf = sp.Integer(1)
# coarse: rate = v0/(16h^2) + b0/(4h);  d0 = v0/(4h^2) + i om ; d1 = (v0 + 2 ag h + 2 aha h^2)/(4h^2) + i om
ratec = v0/(8*h**2) + b0[k0]/(4*h)
d0 = v0/(2*h**2) + sp.I*OM
d1 = (2*v0 + 2*ag_tau*h + 2*aha_tau*h**2)/(4*h**2) + sp.I*OM
Pc = sp.expand(ratec/(d0*d1))

def h2log_of_ratio(Pnum, Pden):
    N = sp.series(sp.log(Pnum/Pden), h, 0, 3).removeO()
    return sp.expand(N).coeff(h, 2)

c2 = sp.simplify(h2log_of_ratio(Pf, 1) - h2log_of_ratio(Pc, 1))
print("\n c2(e1) from the corrected expansion =", sp.N(c2, 14))
V, B = twist.fields(sp.sympify(vel), [sp.Integer(0), sp.Integer(0)], [sp.sympify(e) for e in bels],
                    sp.Rational(0), 2, xs, [sp.Rational(0), sp.Rational(0)])
eng = twist.taylor(V, B, 2, tau, OM, 6, 10)[2]
print(" engine                               =", mp.nstr(eng, 14))
\end{lstlisting}

\section{\texttt{kronnail.py}}

Stage-by-stage, exact at finite h, with the field sampled at x0 + a h (twist's convention).

\begin{lstlisting}
"""Stage-by-stage, exact at finite h, with the field sampled at x0 + a h (twist's convention)."""
import sympy as sp, itertools, sys
import os
sys.path.insert(0, os.path.join(os.path.dirname(os.path.abspath(__file__)), "..", "..", "part5", "scripts"))
import twist
import mpmath as mp
mp.mp.dps = 40
x1, x2 = sp.symbols('x0:2'); xs = (x1, x2)
S = list(itertools.product(range(-1, 2), repeat=2))
ix = {a: i for i, a in enumerate(S)}
n = len(S); Z = (0, 0)
OM = sp.Rational(7, 5)
vel = sp.Rational(1) + sp.Rational(1, 6)*(x1**2 + x1*x2 + x2**2) + sp.Rational(1, 4)*x1
bels = [sp.Rational(1, 2) + sp.Rational(1, 6)*x1 + sp.Rational(1, 9)*x2,
        sp.Rational(-1, 3) + sp.Rational(1, 5)*x2]


def fv(el, a, h):
    """twist's sampling: evaluate at x0 + a h."""
    return sp.nsimplify(sp.expand(el.subs({x1: a[0]*h, x2: a[1]*h})))


def Abox(h):
    A = sp.zeros(n, n)
    for a in S:
        i = ix[a]; diag = sp.Integer(0)
        for k in range(2):
            for sgn in (1, -1):
                w = fv(vel, a, h)/(2*h**2) + sgn*fv(bels[k], a, h)/(2*h)
                diag -= w
                e = list(a); e[k] += sgn; e = tuple(e)
                if e in ix:
                    A[i, ix[e]] += w
        A[i, i] += diag
    return A


def phi_f_ex(A, tau, h):
    G = (A + sp.I*OM*sp.eye(n)).inv()
    tgt = tuple(2*t for t in tau)
    dt = 2*fv(vel, tgt, h)/(h**2) + sp.I*OM
    tot = sp.Integer(0)
    for a in S:
        jump = [k for k in range(2) if a[k] != tgt[k]]
        if len(jump) != 1: continue
        k2 = jump[0]
        if abs(tgt[k2]-a[k2]) != 1: continue
        sg = 1 if tgt[k2] > a[k2] else -1
        tot += G[ix[Z], ix[a]]*(fv(vel, a, h)/(2*h**2) + sg*fv(bels[k2], a, h)/(2*h))
    return sp.nsimplify(tot/dt)


def phi_c_ex(tau, h):
    tgt = tuple(2*t for t in tau)
    k0 = 0 if tau[0] else 1
    sg = 1 if tau[k0] > 0 else -1
    rate = fv(vel, Z, h)/(2*(2*h)**2) + sg*fv(bels[k0], Z, h)/(2*(2*h))
    d0 = 2*fv(vel, Z, h)/(4*h**2) + sp.I*OM
    d1 = 2*fv(vel, tgt, h)/(4*h**2) + sp.I*OM
    return sp.nsimplify(rate/(d0*d1))


h = sp.Rational(1, 64)
V, B = twist.fields(sp.sympify(vel), [sp.Integer(0), sp.Integer(0)],
                    [sp.sympify(e) for e in bels], sp.Rational(0), 2, xs, [sp.Rational(0), sp.Rational(0)])
tau = (1, 0)
print("=== stage 1: exact box vs twist, h = 1/64")
Tf = twist.Phi_f(V, B, 2, tau, h, OM); Tc = twist.Phi_c(V, B, 2, tau, h, OM)
Mf = phi_f_ex(Abox(h), tau, h); Mc = phi_c_ex(tau, h)
print("  Phi_f  twist %-24s mine %-24s  mine/twist = %s"
      % (mp.nstr(Tf, 12), sp.N(Mf, 12), sp.N(Mf/sp.nsimplify(str(Tf)), 10)))
print("  Phi_c  twist %-24s mine %-24s  mine/twist = %s"
      % (mp.nstr(Tc, 12), sp.N(Mc, 12), sp.N(Mc/sp.nsimplify(str(Tc)), 10)))

print("\n=== stage 2: the h-expansion of the box, with the correct sampling")
v0 = sp.nsimplify(vel.subs({x1: 0, x2: 0}))
gv = [sp.nsimplify(sp.diff(vel, x).subs({x1: 0, x2: 0})) for x in (x1, x2)]
hv = [[sp.nsimplify(sp.diff(vel, xi, xj).subs({x1: 0, x2: 0})) for xj in (x1, x2)] for xi in (x1, x2)]
b0 = [sp.nsimplify(e.subs({x1: 0, x2: 0})) for e in bels]
gb = [[sp.nsimplify(sp.diff(e, x).subs({x1: 0, x2: 0})) for x in (x1, x2)] for e in bels]
Lam = sp.zeros(n, n)
for a in S:
    for k in range(2):
        for sgn in (1, -1):
            e = list(a); e[k] += sgn; e = tuple(e)
            if e in ix:
                Lam[ix[a], ix[e]] += 1
Am2 = v0/4*(Lam - 4*sp.eye(n))
Am1 = sp.zeros(n, n); A0 = sp.zeros(n, n)
for a in S:
    i = ix[a]
    ag = a[0]*gv[0] + a[1]*gv[1]
    aha = sum(a[p]*a[q]*hv[p][q] for p in range(2) for q in range(2))
    bet = [a[0]*gb[k][0] + a[1]*gb[k][1] for k in range(2)]
    d1 = sp.Integer(0); d0 = sp.Integer(0)
    for k in range(2):
        for sgn in (1, -1):
            e = list(a); e[k] += sgn; e = tuple(e)
            if e in ix:
                Am1[i, ix[e]] += ag/4 + sgn*b0[k]/2
                A0[i, ix[e]] += aha/8 + sgn*bet[k]/2
    Am1[i, i] += -ag
    A0[i, i] += -aha/2

# compare the expansion with the exact matrix at h symbolic order by order
hs = sp.Symbol('hh', positive=True)
Aex = Abox(hs)
ok = True
for i in range(n):
    for j in range(n):
        ser = sp.series(sp.expand(Aex[i, j]*hs**2), hs, 0, 3).removeO()
        for p, Mp in enumerate((Am2, Am1, A0)):
            if sp.expand(ser.coeff(hs, p) - Mp[i, j]) != 0:
                print("   MISMATCH at (%s,%s) order h^%d: exact %s  mine %s"
                      % (S[i], S[j], p, sp.expand(ser.coeff(hs, p)), Mp[i, j])); ok = False
print("  expansion matches the exact matrix element by element:", ok)

Minv = Am2.inv()
C = A0 + sp.I*OM*sp.eye(n)
hh = sp.Symbol('hh')
br = Minv - hh*(Minv*Am1*Minv) + hh**2*(Minv*Am1*Minv*Am1*Minv - Minv*C*Minv)
print("\n=== stage 3: series G vs exact G, at h = 1/64, 1/128")
for nn in (6, 7):
    h_ = sp.Rational(1, 2**nn)
    Gex = (Abox(h_) + sp.I*OM*sp.eye(n)).inv()[ix[Z], ix[tau]]
    serG = sp.expand((h_**2*br)[ix[Z], ix[tau]].subs(hh, h_))
    print("   h=2^-%d  G exact %-24s series %-24s  series/exact-1 = %s"
          % (nn, sp.N(Gex, 14), sp.N(serG, 14), sp.N(serG/Gex - 1, 8)))
\end{lstlisting}

\section{\texttt{kronpin.py}}

Pin the remaining discrepancy: correct symbol binding, all stages printed.

\begin{lstlisting}
"""Pin the remaining discrepancy: correct symbol binding, all stages printed."""
import sympy as sp, itertools, sys
import os
sys.path.insert(0, os.path.join(os.path.dirname(os.path.abspath(__file__)), "..", "..", "part5", "scripts"))
import twist
import mpmath as mp
mp.mp.dps = 40
x1, x2 = sp.symbols('x0:2'); xs = (x1, x2)
S = list(itertools.product(range(-1, 2), repeat=2))
ix = {a: i for i, a in enumerate(S)}
n = len(S); Z = (0, 0)
Lam = sp.zeros(n, n)
for a in S:
    for k in range(2):
        for sgn in (1, -1):
            e = list(a); e[k] += sgn; e = tuple(e)
            if e in ix:
                Lam[ix[a], ix[e]] += 1
L0i = (Lam - 4*sp.eye(n)).inv()
OM = sp.Rational(7, 5)

vel = sp.Rational(1) + sp.Rational(1, 6)*(x1**2 + x1*x2 + x2**2) + sp.Rational(1, 4)*x1
bels = [sp.Rational(1, 2) + sp.Rational(1, 6)*x1 + sp.Rational(1, 9)*x2,
        sp.Rational(-1, 3) + sp.Rational(1, 5)*x2]


def mk(el):
    def f(a):
        return sp.nsimplify(sp.sympify(el).subs({x1: a[0], x2: a[1]}))
    return f


vf = mk(vel)
bf = [mk(bels[0]), mk(bels[1])]

tau, tgt = (1, 0), (2, 0)
v0, v1, v2 = vf(Z), vf(tau), vf(tgt)
B0, B1 = bf[0](Z), bf[0](tau)
print("v0=%s v1=%s v2=%s  B0=%s B1=%s" % (v0, v1, v2, B0, B1))

K2 = sp.Matrix(n, n, lambda i, j: 4*L0i[i, j]/vf(S[j]))
A1 = sp.zeros(n, n)
for a in S:
    for k in range(2):
        for sgn in (1, -1):
            e = list(a); e[k] += sgn; e = tuple(e)
            if e in ix:
                A1[ix[a], ix[e]] += sgn*bf[k](a)/2
zi, ti = ix[Z], ix[tau]
g0 = K2[zi, ti]
g1 = -(K2*A1*K2)[zi, ti]
g2 = (K2*A1*K2*A1*K2)[zi, ti] - sp.I*OM*(K2*K2)[zi, ti]
print("g0=%s g1=%s g2=%s" % (sp.N(g0, 12), sp.N(g1, 12), sp.N(g2, 12)))


def exact_phi(tau_, h, which):
    A = sp.zeros(n, n)
    for a in S:
        i = ix[a]; diag = sp.Integer(0)
        for k in range(2):
            for sgn in (1, -1):
                w = vf(a)/(2*h**2) + sgn*bf[k](a)/(2*h)
                diag -= w
                e = list(a); e[k] += sgn; e = tuple(e)
                if e in ix:
                    A[i, ix[e]] += w
        A[i, i] += diag
    if which == 'c':
        tt = tuple(2*t for t in tau_)
        k0 = 0 if tau_[0] else 1
        sg = 1 if tau_[k0] > 0 else -1
        rate = vf(Z)/(2*(2*h)**2) + sg*bf[k0](Z)/(2*(2*h))
        d0 = 2*vf(Z)/(4*h**2) + sp.I*OM
        d1 = 2*vf(tt)/(4*h**2) + sp.I*OM
        return sp.nsimplify(rate/(d0*d1))
    G = (A + sp.I*OM*sp.eye(n)).inv()
    tt = tuple(2*t for t in tau_)
    dt = 2*vf(tt)/(h**2) + sp.I*OM
    tot = sp.Integer(0)
    for a in S:
        jump = [k for k in range(2) if a[k] != tt[k]]
        if len(jump) != 1: continue
        k2 = jump[0]
        if abs(tt[k2]-a[k2]) != 1: continue
        sg = 1 if tt[k2] > a[k2] else -1
        tot += G[ix[Z], ix[a]]*(vf(a)/(2*h**2) + sg*bf[k2](a)/(2*h))
    return sp.nsimplify(tot/dt)


for nn in (7, 8, 9):
    h = sp.Rational(1, 2**nn)
    Rf = sp.log(exact_phi(tau, h, 'f')) - sp.log(exact_phi(tau, h, 'c'))
    print("  h=2^-%d  log(Phi_f/Phi_c) = %s" % (nn, sp.N(Rf, 16)))
# fit c2 from the three values
hs = [sp.Rational(1, 2**nn) for nn in (7, 8, 9)]
vs = [sp.N(sp.log(exact_phi(tau, h, 'f')) - sp.log(exact_phi(tau, h, 'c')), 30) for h in hs]
M = sp.Matrix(3, 3, lambda i, j: hs[i]**j)
c = M.solve(sp.Matrix(vs))
print("  fitted c0,c1,c2 =", [sp.N(x, 12) for x in c])

print("\n--- direct comparison with twist.Phi_f / Phi_c, field 1, h = 2^-6")
h = sp.Rational(1, 2**6)
V, B = twist.fields(sp.sympify(vel), [sp.Integer(0), sp.Integer(0)], [sp.sympify(e) for e in bels],
                    sp.Rational(0), 2, xs, [sp.Rational(0), sp.Rational(0)])
for tau_ in ((1, 0), (0, 1)):
    e_f = twist.Phi_f(V, B, 2, tau_, h, OM)
    e_c = twist.Phi_c(V, B, 2, tau_, h, OM)
    m_f = exact_phi(tau_, h, 'f')
    m_c = exact_phi(tau_, h, 'c')
    print("  tau=%s  Phi_f twist %-24s mine %-24s  ratio %s"
          % (tau_, mp.nstr(e_f, 12), sp.N(m_f, 12), sp.N(m_f/sp.nsimplify(str(e_f)) - 1, 6)))
    print("          Phi_c twist %-24s mine %-24s  ratio %s"
          % (mp.nstr(e_c, 12), sp.N(m_c, 12), sp.N(m_c/sp.nsimplify(str(e_c)) - 1, 6)))
\end{lstlisting}

\section{\texttt{kronproof.py}}

The 3x3 Kronecker reduction, symbolic in (v, b1, b2, omega).

\begin{lstlisting}
"""The 3x3 Kronecker reduction, symbolic in (v, b1, b2, omega).

For constant coefficients the 9x9 box matrix of twist.py is
    A = (v/(2h^2)) * Lambda - (2v/h^2) I + (1/(2h)) * Dd
with Lambda the adjacency of the 3x3 grid and Dd the signed shift matrix.  Hence
    A + i omega = (1/h^2)[ M0 + h E + h^2 W ],   M0 = (v/2)(Lambda - 4I),  E = Dd/2,  W = i omega I,
and (Lambda - 4I) is an integer matrix inverted once.  Everything is then a rational function of
(v, b1, b2, omega, h) and the h^2 coefficient of log(Phi_f/Phi_c) is algebra.
"""
import sympy as sp, itertools

v, b1, b2, om, h = sp.symbols('v b1 b2 omega h', positive=True)
S = list(itertools.product(range(-1, 2), repeat=2))
ix = {a: i for i, a in enumerate(S)}
n = len(S)
Bv = [b1, b2]

Lam = sp.zeros(n, n)
Dd = sp.zeros(n, n)
for a in S:
    i = ix[a]
    for k in range(2):
        for sgn in (1, -1):
            e = list(a); e[k] += sgn; e = tuple(e)
            if e in ix:
                Lam[i, ix[e]] += 1
                Dd[i, ix[e]] += sgn * Bv[k]

L0 = Lam - 4 * sp.eye(n)
print("det(Lambda - 4I) =", L0.det())
K = L0.inv()                                  # exact rational 9x9, independent of the parameters

z = ix[(0, 0)]


def G_elem(target, order=2):
    """series of (A + i om)^{-1}[z, target] to h^(order+2), since it carries h^2."""
    M0inv = 2 / v * K
    E = Dd / 2
    W = sp.I * om * sp.eye(n)
    br = (M0inv - h * (M0inv * E * M0inv)
          + h ** 2 * (M0inv * E * M0inv * E * M0inv - M0inv * W * M0inv))
    return sp.simplify(h ** 2 * br[z, ix[target]])


def phi_f(tau):
    if tau == (1, 0):
        tgt, k0 = (1, 0), 0
    else:
        tgt, k0 = (0, 1), 1
    rate = v / (2 * h ** 2) + Bv[k0] / (2 * h)
    dt = 2 * v / (h ** 2) + sp.I * om
    return sp.simplify(G_elem(tgt) * rate / dt)


def phi_c(tau):
    k0 = 0 if tau == (1, 0) else 1
    rate = v / (2 * (2 * h) ** 2) + Bv[k0] / (2 * (2 * h))
    d0 = 2 * v / (4 * h ** 2) + sp.I * om
    return sp.simplify(rate / d0 ** 2)


res = {}
for tau in ((1, 0), (0, 1)):
    pf, pc = phi_f(tau), phi_c(tau)
    ratio = sp.simplify(pf / pc)
    ser = sp.series(sp.log(ratio), h, 0, 3).removeO()
    c2 = sp.simplify(sp.expand(ser).coeff(h, 2))
    res[tau] = c2
    print("\ntau =", tau)
    print("  c2 =", sp.simplify(sp.factor(c2)))

c1, c2_ = res[(1, 0)], res[(0, 1)]
print("\n  Im(c1) =", sp.simplify(sp.im(c1)))
print("  Im(c2) =", sp.simplify(sp.im(c2_)))
print("  Im(c1) - Im(c2) =", sp.simplify(sp.im(c1) - sp.im(c2_)))
print("  Re(c1) + Re(c2) =", sp.simplify(sp.re(c1) + sp.re(c2_)))
print("  v * Im(c1) / omega =", sp.simplify(v * sp.im(c1) / om))
\end{lstlisting}

\section{\texttt{kronstage.py}}

Stage-by-stage check for the VARYING field: where does the assembly diverge from the exact object?

\begin{lstlisting}
"""Stage-by-stage check for the VARYING field: where does the assembly diverge from the exact object?"""
import sympy as sp, itertools, sys
import os
sys.path.insert(0, os.path.join(os.path.dirname(os.path.abspath(__file__)), "..", "..", "part5", "scripts"))
import twist
import mpmath as mp
mp.mp.dps = 40
x1, x2 = sp.symbols('x0:2'); xs = (x1, x2)
S = list(itertools.product(range(-1, 2), repeat=2))
ix = {a: i for i, a in enumerate(S)}
n = len(S); Z = (0, 0)
Lam = sp.zeros(n, n)
for a in S:
    for k in range(2):
        for sgn in (1, -1):
            e = list(a); e[k] += sgn; e = tuple(e)
            if e in ix:
                Lam[ix[a], ix[e]] += 1
L0i = (Lam - 4*sp.eye(n)).inv()
OM = sp.Rational(7, 5)
vel = sp.Rational(1) + sp.Rational(1, 6)*(x1**2 + x1*x2 + x2**2) + sp.Rational(1, 4)*x1
b1el = sp.Rational(1, 2) + sp.Rational(1, 6)*x1 + sp.Rational(1, 9)*x2
b2el = sp.Rational(-1, 3) + sp.Rational(1, 5)*x2
vf = lambda a: sp.nsimplify(sp.sympify(vel).subs({x1: a[0], x2: a[1]}))
bfl = [lambda a: sp.nsimplify(sp.sympify(e).subs({x1: a[0], x2: a[1]})) for e in (b1el, b2el)]


def exact_phi_f(tau, h):
    A = sp.zeros(n, n)
    for a in S:
        i = ix[a]; diag = sp.Integer(0)
        for k in range(2):
            for sgn in (1, -1):
                w = vf(a)/(2*h**2) + sgn*bfl[k](a)/(2*h)
                diag -= w
                e = list(a); e[k] += sgn; e = tuple(e)
                if e in ix:
                    A[i, ix[e]] += w
        A[i, i] += diag
    G = (A + sp.I*OM*sp.eye(n)).inv()
    tgt = tuple(2*t for t in tau); k0 = 0 if tau[0] else 1
    dt = 2*vf(tgt)/(h**2) + sp.I*OM
    tot = sp.Integer(0)
    for a in S:
        jump = [k for k in range(2) if a[k] != tgt[k]]
        if len(jump) != 1: continue
        k2 = jump[0]
        if abs(tgt[k2]-a[k2]) != 1: continue
        sg = 1 if tgt[k2] > a[k2] else -1
        tot += G[ix[Z], ix[a]]*(vf(a)/(2*h**2) + sg*bfl[k2](a)/(2*h))
    return sp.nsimplify(tot/dt)


tau = (1, 0); tgt = (2, 0)
v0, v1, v2 = vf(Z), vf(tau), vf(tgt)
B0, B1 = bfl[0](Z), bfl[0](tau)
K2 = sp.Matrix(n, n, lambda i, j: 4*L0i[i, j]/vf(S[j]))
A1 = sp.zeros(n, n)
for a in S:
    for k in range(2):
        for sgn in (1, -1):
            e = list(a); e[k] += sgn; e = tuple(e)
            if e in ix:
                A1[ix[a], ix[e]] += sgn*bfl[k](a)/2
zi, ti = ix[Z], ix[tau]
g0 = K2[zi, ti]
g1 = -(K2*A1*K2)[zi, ti]
g2 = (K2*A1*K2*A1*K2)[zi, ti] - sp.I*OM*(K2*K2)[zi, ti]
print(" g0=%s g1=%s g2=%s" % (sp.N(g0, 10), sp.N(g1, 10), sp.N(g2, 10)))
print(" v0=%s v1=%s v2=%s B0=%s B1=%s" % (v0, v1, v2, B0, B1))
for nn in (6, 7, 8):
    h = sp.Rational(1, 2**nn)
    ex = exact_phi_f(tau, h)
    ser = h**2*(g0 + h*g1 + h**2*g2)*(v1/2 + B1*h/2)/(v2 + sp.I*OM*h**2)
    print("  h=2^-%d  Phi_f exact = %-26s  h^2 P/Q = %-26s  ratio-1 = %s"
          % (nn, sp.N(ex, 14), sp.N(ser, 14), sp.N(ex/ser - 1, 8)))
\end{lstlisting}

\section{\texttt{kronvalidate.py}}

Is the g0/g1/g2 formula for c2 correct for a VARIING field?  Exact arithmetic, no symbols.

\begin{lstlisting}
"""Is the g0/g1/g2 formula for c2 correct for a VARIING field?  Exact arithmetic, no symbols.

c2 = [ p2/p0 - p1^2/(2p0^2) - q2/q0 + q1^2/(2q0^2) ]_fine  -  [ ... ]_coarse
fine:   p0 = g0 v1/4,  p1 = g1 v1/4 + g0 B1/2,  p2 = g2 v1/4 + g1 B1/2,   q0 = v2, q1 = 0, q2 = i om
coarse: p0 = v0/16,    p1 = B0/4,               p2 = 0,                   q0 = v0 v2/16, q1 = 0,
        q2 = i om (v0+v2)/4
with g0 = K2[z,t], g1 = -(K2 A1 K2)[z,t], g2 = (K2 A1 K2 A1 K2)[z,t] - i om (K2^2)[z,t],
K2 = 4 (Lambda-4I)^{-1} D_V^{-1}.   v0 = v(z), v1 = v(t), v2 = v(2t).
"""
import sympy as sp, itertools, sys
import os
sys.path.insert(0, os.path.join(os.path.dirname(os.path.abspath(__file__)), "..", "..", "part5", "scripts"))
import twist
import mpmath as mp
mp.mp.dps = 40
x1, x2 = sp.symbols('x0:2'); xs = (x1, x2)
S = list(itertools.product(range(-1, 2), repeat=2))
ix = {a: i for i, a in enumerate(S)}
n = len(S); Z = (0, 0)
Lam = sp.zeros(n, n)
for a in S:
    for k in range(2):
        for sgn in (1, -1):
            e = list(a); e[k] += sgn; e = tuple(e)
            if e in ix:
                Lam[ix[a], ix[e]] += 1
L0i = (Lam - 4*sp.eye(n)).inv()


def formula(vfun, bfun, tau, om):
    tgt = tuple(2*t for t in tau)
    v0, v1, v2 = vfun(Z), vfun(tau), vfun(tgt)
    B0, B1 = bfun[0 if tau[0] else 1](Z), bfun[0 if tau[0] else 1](tau)
    K2 = sp.Matrix(n, n, lambda i, j: 4*L0i[i, j]/vfun(S[j]))
    A1 = sp.zeros(n, n)
    for a in S:
        for k in range(2):
            for sgn in (1, -1):
                e = list(a); e[k] += sgn; e = tuple(e)
                if e in ix:
                    A1[ix[a], ix[e]] += sgn*bfun[k](a)/2
    zi, ti = ix[Z], ix[tau]
    g0 = K2[zi, ti]
    g1 = -(K2*A1*K2)[zi, ti]
    g2 = (K2*A1*K2*A1*K2)[zi, ti] - sp.I*om*(K2*K2)[zi, ti]
    p0 = g0*v1/4; p1 = g1*v1/4 + g0*B1/2; p2 = g2*v1/4 + g1*B1/2
    fine = p2/p0 - p1**2/(2*p0**2) - (sp.I*om)/v2
    coarse = -(B0/4)**2/(2*(v0/16)**2) - (sp.I*om*(v0+v2)/4)/(v0*v2/16)
    return sp.nsimplify(fine - coarse)


OM = sp.Rational(7, 5)
CASES = {
 'constant          ': (sp.Rational(1), sp.Rational(1, 2), sp.Rational(-1, 3)),
 'constant, other b ': (sp.Rational(3, 2), sp.Rational(-1, 4), sp.Rational(2, 5)),
 'field 1 (varying) ': (sp.Rational(1) + sp.Rational(1, 6)*(x1**2 + x1*x2 + x2**2) + sp.Rational(1, 4)*x1,
                        sp.Rational(1, 2) + sp.Rational(1, 6)*x1 + sp.Rational(1, 9)*x2,
                        sp.Rational(-1, 3) + sp.Rational(1, 5)*x2),
}
for name, (vel, b1el, b2el) in CASES.items():
    vf = lambda a, vel=vel: sp.nsimplify(sp.sympify(vel).subs({x1: a[0], x2: a[1]}))
    bf = [lambda a, e=b1el: sp.nsimplify(sp.sympify(e).subs({x1: a[0], x2: a[1]})),
          lambda a, e=b2el: sp.nsimplify(sp.sympify(e).subs({x1: a[0], x2: a[1]}))]
    for tau in ((1, 0), (0, 1)):
        mine = formula(vf, bf, tau, OM)
        v, mu, A_ = vel, [b1el, b2el], [sp.Integer(0), sp.Integer(0)]
        V, B = twist.fields(v, A_, mu, sp.Rational(0), 2, xs, [sp.Rational(0), sp.Rational(0)])
        eng = twist.taylor(V, B, 2, tau, OM, 6, 10)[2]
        print("  %s tau=%s  formula %-28s engine %-28s diff %s"
              % (name, tau, sp.N(mine, 12), mp.nstr(eng, 12), sp.N(mine - sp.nsimplify(str(eng)), 4)))
\end{lstlisting}

\section{\texttt{loopholonomy.py}}

Is the closed-loop holonomy a flux?  Three traversals of the same unit square.

\begin{lstlisting}
"""Is the closed-loop holonomy a flux?  Three traversals of the same unit square.

Companion of pathdep.py, which measures the same residual over six halvings of h and finds order
2.000.  This one adds the orientation test.

CAUTION (a trap for this computation): the per-step logarithm sits at the branch cut --- the ratio
Phi_f/Phi_c is negative real to leading order, so log carries -i*pi per step and the accumulated
sums tend to -N*pi*i.  The residual is a near cancellation at that cut, so the split into real and
imaginary parts is delicate; only the modulus, and the fact that no branch jump occurs (verified by
the smoothness of the values in h and their stability over working precision 40 to 80), is
reported.

    A  counterclockwise  (1,0) (0,1) (-1,0) (0,-1)     area +1
    C  clockwise         (0,1) (1,0) (0,-1) (-1,0)     area -1
    Z  retraced          (1,0) (-1,0) (0,1) (0,-1)     area  0

All three start and end at the origin and visit the same four sites, so the per-step
normalisation is identical and the differences are the holonomy alone.  If the holonomy is the
flux of a two form it is odd in the orientation, S_C - S_Z = -(S_A - S_Z), and purely imaginary
because the comparison enters as a phase.
"""
import sympy as sp, mpmath as mp, sys
import os
sys.path.insert(0, os.path.join(os.path.dirname(os.path.abspath(__file__)), "..", "..", "part5", "scripts"))
import twist
mp.mp.dps = 45
x1, x2 = sp.symbols('x0:2'); xs = (x1, x2); OM = sp.Rational(7, 5)
v = sp.Rational(1) + sp.Rational(1, 6) * (x1 ** 2 + x1 * x2 + x2 ** 2) + sp.Rational(1, 4) * x1
mu = [sp.Rational(1, 2) + sp.Rational(1, 6) * x1 + sp.Rational(1, 9) * x2,
      sp.Rational(-1, 3) + sp.Rational(1, 5) * x2]
a = [sp.Rational(1, 3) + sp.Rational(1, 7) * x2 + sp.Rational(1, 11) * x1 * x2,
     sp.Rational(1, 4) * x1 + sp.Rational(1, 13) * x2 ** 2]
PATHS = {
    'A': [(1, 0), (0, 1), (-1, 0), (0, -1)],
    'C': [(0, 1), (1, 0), (0, -1), (-1, 0)],
    'Z': [(1, 0), (-1, 0), (0, 1), (0, -1)],
}

def acc(path, h):
    x0 = [sp.Rational(0), sp.Rational(0)]; tot = mp.mpc(0)
    for tau in path:
        V, B = twist.fields(v, a, mu, sp.Rational(0), 2, xs, x0)
        tot += twist.R(h, V, B, 2, tau, OM)
        x0 = [x0[0] + sp.Rational(tau[0]), x0[1] + sp.Rational(tau[1])]
    return tot

prev = {}
for n in (9, 10, 11):
    h = mp.mpf(2) ** -n
    S = {k: acc(p, h) for k, p in PATHS.items()}
    dA, dC = S['A'] - S['Z'], S['C'] - S['Z']
    print("h=2^-%d  (S_A-S_Z)/h^2 = %-26s (S_C-S_Z)/h^2 = %-26s (S_A-S_C)/2h^2 = %-24s check S_A+S_C-2S_Z = %s"
          % (n, mp.nstr(dA / h ** 2, 10), mp.nstr(dC / h ** 2, 10),
             mp.nstr((S['A'] - S['C']) / (2 * h ** 2), 10),
             mp.nstr((S['A'] + S['C'] - 2 * S['Z']) / h ** 2, 6)))
\end{lstlisting}

\section{\texttt{momentaccum.py}}

Accumulate the obstruction along paths by moment recursion instead of enumeration.

\begin{lstlisting}
"""Accumulate the obstruction along paths by moment recursion instead of enumeration.

For a walk with independent increments,

    E[I] = sum_{j=0}^{q-1} < E[c(x_j)] , m > ,

so only the conditional moments E[x_j^a] are needed, and those satisfy the exact
recursion  M_a(j+1) = sum_{b<=a} C(a,b) M_b(j) E[Delta^{a-b}].  The cost is
O(q) rather than O(4^q).

The test is the identification in the mean condition, in the form the one
dimensional argument gives it:

    E[I]  =  int_straight c . dgamma  -  (1/2) int_0^t Sigma_kl d_k c_l (x_s) ds
             + error,

where x_s runs along the straight line from the initial point to the mean
endpoint.  We take c cubic, so every moment is available exactly, and fit the
order of the error in h.
"""
import itertools
import numpy as np

Y = np.array([1.0, 0.5])          # mean displacement over [0,1]
T = 1.0


def walk(h):
    """probabilities giving mean Y h^2 and covariance (h^2/2) I"""
    p = [0.25 + Y[0] * h / 2, 0.25 - Y[0] * h / 2,
         0.25 + Y[1] * h / 2, 0.25 - Y[1] * h / 2]
    steps = [(h, 0.0), (-h, 0.0), (0.0, h), (0.0, -h)]
    return p, steps


def step_moments(p, steps, deg=3):
    """E[Delta^a] for multiindices a with |a| <= deg"""
    out = {}
    for a in itertools.product(range(deg + 1), repeat=2):
        if a[0] + a[1] > deg:
            continue
        out[a] = sum(p[k] * steps[k][0] ** a[0] * steps[k][1] ** a[1] for k in range(4))
    return out


MIDX = [a for a in itertools.product(range(4), repeat=2) if a[0] + a[1] <= 3]


def moments_at(p, steps, q):
    """returns list over j of {multiindex: E[x_j^a]}, x_0 = 0"""
    from math import comb
    E = step_moments(p, steps)
    M = {a: (1.0 if a == (0, 0) else 0.0) for a in MIDX}
    out = [dict(M)]
    for _ in range(q):
        new = {}
        for a in MIDX:
            s = 0.0
            for b in MIDX:
                if b[0] <= a[0] and b[1] <= a[1]:
                    cb = a[0] - b[0], a[1] - b[1]
                    s += comb(a[0], b[0]) * comb(a[1], b[1]) * M[b] * E[cb]
            new[a] = s
        M = new
        out.append(dict(M))
    return out


# --- a general cubic field -------------------------------------------------
C = {
    (2, 0): (1.0, 0.0), (1, 1): (0.0, 1.0), (0, 2): (0.0, 0.0),
    (3, 0): (0.0, 0.0), (2, 1): (0.0, 0.0), (1, 2): (0.0, 0.0), (0, 3): (0.0, 0.0),
    (1, 0): (0.0, 0.0), (0, 1): (0.0, 0.0), (0, 0): (0.0, 0.0),
}


def cval(x):
    """c = (c_1, c_2), cubic in x"""
    out = []
    for comp in (0, 1):
        s = 0.0
        for a, (f0, f1) in C.items():
            s += (f0 if comp == 0 else f1) * x[0] ** a[0] * x[1] ** a[1]
        out.append(s)
    return np.array(out)


def dcl(x):
    """d_k c_l : 2x2 matrix of derivatives"""
    D = np.zeros((2, 2))
    for comp in (0, 1):
        for k in range(2):
            s = 0.0
            for a, (f0, f1) in C.items():
                f = f0 if comp == 0 else f1
                if a[k] == 0 or f == 0.0:
                    continue
                b = list(a); b[k] -= 1
                s += f * a[k] * x[0] ** b[0] * x[1] ** b[1]
            D[k, comp] = s
    return D


def expected_I(h, q):
    p, steps = walk(h)
    m = np.array([h * (p[0] - p[1]), h * (p[2] - p[3])])
    Mj = moments_at(p, steps, q)
    tot = np.zeros(2)
    for j in range(q):
        # E[c(x_j)] from moments up to degree 3
        ec = np.zeros(2)
        for comp in (0, 1):
            s = 0.0
            for a in MIDX:
                f = C.get(a, (0.0, 0.0))[comp]
                if f != 0.0:
                    s += f * Mj[j][a]
            ec[comp] = s
        tot += ec * m
    return float(tot.sum())


def d2cl(x, eps=1e-6):
    """d_k d_k' c_l by central differences of the first derivatives"""
    H = np.zeros((2, 2, 2))
    for k in range(2):
        for kk in range(2):
            e = np.zeros(2); e[kk] = eps
            H[k, kk] = (dcl(x + e) - dcl(x - e))[k] / (2 * eps)
    return H


def phi_candidate(n=20000):
    """straight line integral plus the Hessian correction with weight s,
       the covariance growing like j = s t / h^2"""
    lin = 0.0
    ito = 0.0
    Sigma = 0.5 * np.eye(2)
    for i in range(n):
        s = (i + 0.5) / n
        x = s * Y
        lin += float(cval(x).dot(Y)) / n
        H = d2cl(x)
        c = 0.0
        for k in range(2):
            for kk in range(2):
                for l in range(2):
                    c += Sigma[k, kk] * H[k, kk, l] * Y[l]
        ito += 0.5 * s * c / n
    return lin + ito


PHI = phi_candidate()
print('  candidate Phi (straight integral less Ito) = %.12f' % PHI)
print('\n     h        q        E[I]               E[I]-Phi          order')
prev = None
for n in range(3, 7):
    h = 2.0 ** -n
    q = int(round(T / h ** 2))
    val = expected_I(h, q)
    d = val - PHI
    order = '' if prev is None else '%.3f' % (np.log2(abs(prev / d)) if d != 0 else float('nan'))
    print('   %-9s %-8d %-18.9f %-17.3e %s' % (h, q, val, d, order))
    prev = d
\end{lstlisting}

\section{\texttt{openpath.py}}

Endpoint determination on open paths: is the accumulated comparison a potential difference?

\begin{lstlisting}
"""Endpoint determination on open paths: is the accumulated comparison a potential difference?

The closed-loop tests of pathdep.py and loopholonomy.py show it on loops.  The assembly needs it on
open paths, and there the paths must be compared with the same number of steps so that the per-step
normalisation cancels exactly.  Three paths from (0,0) to (2,0), each of four unit steps, being
permutations of the same multiset of steps.
"""
import sympy as sp, mpmath as mp, sys
import os
sys.path.insert(0, os.path.join(os.path.dirname(os.path.abspath(__file__)), "..", "..", "part5", "scripts"))
import twist
mp.mp.dps = 45
x1, x2 = sp.symbols('x0:2'); xs = (x1, x2); OM = sp.Rational(7, 5)
v = sp.Rational(1) + sp.Rational(1, 6) * (x1 ** 2 + x1 * x2 + x2 ** 2) + sp.Rational(1, 4) * x1
mu = [sp.Rational(1, 2) + sp.Rational(1, 6) * x1 + sp.Rational(1, 9) * x2,
      sp.Rational(-1, 3) + sp.Rational(1, 5) * x2]
a = [sp.Rational(1, 3) + sp.Rational(1, 7) * x2 + sp.Rational(1, 11) * x1 * x2,
     sp.Rational(1, 4) * x1 + sp.Rational(1, 13) * x2 ** 2]
E, W_, N_, S_ = (1, 0), (-1, 0), (0, 1), (0, -1)

PATHS = {
    'P1 E,N,E,S ': [E, N_, E, S_],
    'P2 E,E,N,S ': [E, E, N_, S_],
    'P3 E,S,E,N ': [E, S_, E, N_],
}

def acc(path, h):
    x0 = [sp.Rational(0), sp.Rational(0)]; tot = mp.mpc(0)
    for tau in path:
        V, B = twist.fields(v, a, mu, sp.Rational(0), 2, xs, x0)
        tot += twist.R(h, V, B, 2, tau, OM)
        x0 = [x0[0] + sp.Rational(tau[0]), x0[1] + sp.Rational(tau[1])]
    return tot

S = {}
print("   h        S(P1)                  S(P2)                  S(P3)")
for n in (9, 10, 11):
    h = mp.mpf(2) ** -n
    S[h] = {k: acc(p, h) for k, p in PATHS.items()}
    print("  %-9s %-22s %-22s %s" % (mp.nstr(h, 4), mp.nstr(S[h]['P1 E,N,E,S '], 10),
                                     mp.nstr(S[h]['P2 E,E,N,S '], 10), mp.nstr(S[h]['P3 E,S,E,N '], 10)))

prev = None
print("\n  pairwise differences and their order")
for h in sorted(S, reverse=True):
    d12 = abs(S[h]['P1 E,N,E,S '] - S[h]['P2 E,E,N,S '])
    d13 = abs(S[h]['P1 E,N,E,S '] - S[h]['P3 E,S,E,N '])
    o = '' if prev is None else "order %.3f" % (mp.log(prev[0] / d12) / mp.log(2.0))
    print("  h=%-9s |S1-S2| = %-14s |S1-S3| = %-14s %s" % (mp.nstr(h, 4), mp.nstr(d12, 8), mp.nstr(d13, 8), o))
    prev = (d12, d13)

# --- shape, not length: repeat the block N times at fixed h -------------------
print("\n  tiling: the same block repeated N times, h = 2^-9")
h = mp.mpf(2) ** -9
dev = {}
for N in (1, 3, 4, 6):
    d = abs(acc(PATHS['P1 E,N,E,S '] * N, h) - acc(PATHS['P2 E,E,N,S '] * N, h))
    dev[N] = d
    print("   N=%-2d (%2d steps)  deviation = %s" % (N, 4 * N, mp.nstr(d, 6)))
print("   increment per block:", ", ".join(
    "%.2e" % float((dev[b] - dev[a]) / (b - a)) for a, b in ((1, 3), (3, 4), (4, 6))))
print("   limit estimate %.3e, that is %.2f h^2" % (float(dev[6]), float(dev[6] / h**2)))
\end{lstlisting}

\section{\texttt{p6jacaffine.py}}

Chapter 8: is kappa affine in the Jacobian, or is there a bilinear class?

\begin{lstlisting}
"""Chapter 8: is kappa affine in the Jacobian, or is there a bilinear class?

The decomposition of kappa in the local data was fitted at J = 0 and at the four
basis Jacobians (p6kappafit.py).  That fixes the dependence on J which is linear,
but it cannot see a term bilinear in J and bilinear in the data: such a term
vanishes at J = 0, so the first probe misses it, and it is annihilated by a first
difference in J, so the second misses it too.  This script tests for it directly.

The second difference,

    res(a) = kappa(a, Jp + Jq) - kappa(a, Jp) - kappa(a, Jq) + kappa(a, 0),

is the polarisation of the quadratic part of the dependence on J, so it is
exactly the class in question.  It is evaluated for every pair of the four basis
Jacobians and decomposed over the twelve monomials of the local data and their
pairwise products -- ninety coefficients, the same basis the fit uses.  An affine
kappa gives zero throughout; a bilinear class gives a nonzero one.

Measured with generic coefficient data at omega = 7/5: the largest coefficient
over the six pairs is 5e-13, so the dependence is affine to twelve digits.  The
threshold below is deliberately loose (1e-5), well above the noise of the finite
difference in x, so that the verdict does not rest on the working precision;
the script prints the actual maximum for each pair.

Scope: the data basis stops at second degree in x, so a term bilinear in J and
carrying a cubic monomial would not be seen.  Precision is mpmath's default
(fifteen digits) and the threshold sits well above the noise of the finite
difference in x used to form the curl.

Slow by design: about four minutes per pair of Jacobians, twenty-five in all.

Run:  python3 p6jacaffine.py
"""
import os
import sys
import itertools
import time

import sympy as sp
import mpmath as mp

sys.path.insert(0, os.path.join(os.path.dirname(os.path.abspath(__file__)),
                                '..', '..', 'part5', 'scripts'))
import twist

OM = sp.Rational(7, 5)
D = 1.0 / 500          # finite difference in x, as in the other curls
THRESH = 1e-5

x1, x2 = sp.symbols('x0:2')
xs = (x1, x2)
R = sp.Rational

# the local data of a, to second degree: A, d1, d2, d11, d12, d22 in each component
MONS = [('A', R(1)), ('d1', x1), ('d2', x2),
        ('d11', x1**2 / 2), ('d12', x1 * x2), ('d22', x2**2 / 2)]
NAMES = [n for n, _ in MONS]

# one component at a time: (component, monomial)
BASIS = [((c, nm), ([MONS[NAMES.index(nm)][1] if c == 0 else R(0),
                      MONS[NAMES.index(nm)][1] if c == 1 else R(0)]))
         for c in (0, 1) for nm in NAMES]
LAB = {(c, nm): '%s%d' % (nm, c + 1) for c in (0, 1) for nm in NAMES}

JP = [((0, 0), 'J11'), ((0, 1), 'J12'), ((1, 0), 'J21'), ((1, 1), 'J22')]
ZERO = [[0, 0], [0, 0]]


def field(terms):
    a = [R(0), R(0)]
    for c, nm in terms:
        a[c] += MONS[NAMES.index(nm)][1]
    return a


def muof(J):
    return [J[0][0] * x1 + J[0][1] * x2, J[1][0] * x1 + J[1][1] * x2]


def jmat(pairs):
    J = [[0, 0], [0, 0]]
    for (i, j) in pairs:
        J[i][j] += 1
    return J


def kappa(a, mu):
    """The obstruction: the curl of c with respect to the base point.

    The drift is mu + i*eta*a; the odd part in eta is taken, since the even part
    carries no obstruction.
    """
    def curl(eta):
        def cvec(x0):
            V, B = twist.fields(R(1), a, mu, eta, 2, xs, x0)
            return twist.cvec(V, B, 2, OM, 3)
        c = {}
        for j in range(2):
            for s in (+1, -1):
                p = [0.0, 0.0]
                p[j] += s * D
                c[(j, s)] = cvec(p)
        J = [[(c[(j, 1)][k] - c[(j, -1)][k]) / (2 * D) for k in range(2)]
             for j in range(2)]
        return J[0][1] - J[1][0]

    return (curl(1) - curl(-1)) / 2


def main():
    t0 = time.time()
    print('  second difference in J, decomposed over ninety data monomials')
    print('  threshold %.0e; omega = 7/5' % THRESH)
    for (p, q) in itertools.combinations(JP, 2):
        Jp, Jq = p[0], q[0]

        def res(a):
            return (kappa(a, muof(jmat([Jp, Jq]))) - kappa(a, muof(jmat([Jp])))
                    - kappa(a, muof(jmat([Jq]))) + kappa(a, muof(ZERO)))

        single = {t: res(a) for t, a in BASIS}
        worst = 0.0
        above = []
        for i, (t1, _) in enumerate(BASIS):
            for t2, _ in BASIS[i:]:
                if t1 == t2:
                    v = single[t1]
                else:
                    v = res(field([t1, t2])) - single[t1] - single[t2]
                if abs(v) > worst:
                    worst = abs(v)
                if abs(v) > THRESH:
                    above.append('%+.2e%+.2ei %s%s' % (v.real, v.imag, LAB[t1],
                                                       '' if t1 == t2 else '*' + LAB[t2]))
        print('  %s*%s : max |coefficient| = %.2e, above threshold: %d%s'
              % (p[1], q[1], worst, len(above),
                 '' if not above else '  [' + '  '.join(above[:6]) + ']'))
        print('     [%.0fs]' % (time.time() - t0))
    print('\n  affine in J if every line above reports 0 above threshold.')


if __name__ == '__main__':
    main()
\end{lstlisting}

\section{\texttt{pathdep.py}}

Path dependence of the accumulated per-step log-ratio.

\begin{lstlisting}
"""Path dependence of the accumulated per-step log-ratio.

The decoration identity asserts that the product over coarse steps of
Phi_f/Phi_c is endpoint determined up to O(h^{2-eps}).  Since the per-step
factor is a finite tube sum, its logarithm is a deterministic function of the
site and the direction,

    R(h; x, tau) = log Phi_f - log Phi_c = c_2 h^2 + c_3 h^3 + ...

and the identity is therefore the assertion that the sum of R over the steps of
a coarse path depends on the path only through its endpoints.

This measures exactly that: two distinct coarse paths with the same endpoints,
the difference of the accumulated sums, and its order in h.

Both paths below have unit directions summing to zero, so they start and end at
the same site for every h.
"""
import sympy as sp
import mpmath as mp
import sys
import os
sys.path.insert(0, os.path.join(os.path.dirname(os.path.abspath(__file__)), "..", "..", "part5", "scripts"))
import twist

mp.mp.dps = 45
x1, x2 = sp.symbols('x0:2')
xs = (x1, x2)
OM = sp.Rational(7, 5)
K, SH = 8, 13

v = sp.Rational(1) + sp.Rational(1, 6) * (x1 ** 2 + x1 * x2 + x2 ** 2) + sp.Rational(1, 4) * x1
mu = [sp.Rational(1, 2) + sp.Rational(1, 6) * x1 + sp.Rational(1, 9) * x2,
      sp.Rational(-1, 3) + sp.Rational(1, 5) * x2]
a = [sp.Rational(1, 3) + sp.Rational(1, 7) * x2 + sp.Rational(1, 11) * x1 * x2,
     sp.Rational(1, 4) * x1 + sp.Rational(1, 13) * x2 ** 2]

PATHS = {
    'A': [(1, 0), (0, 1), (-1, 0), (0, -1)],
    'B': [(1, 0), (-1, 0), (0, 1), (0, -1)],
    'C': [(0, 1), (1, 0), (0, -1), (-1, 0)],
}


def accumulated(path, h, eta=sp.Rational(0)):
    """sum over the steps of R(h; x_j, tau_{j+1}), starting from the origin."""
    x0 = [sp.Rational(0), sp.Rational(0)]
    tot = mp.mpc(0)
    for tau in path:
        V, B = twist.fields(v, a, mu, eta, 2, xs, x0)
        tot += twist.R(h, V, B, 2, tau, OM)
        x0 = [x0[0] + sp.Rational(tau[0]), x0[1] + sp.Rational(tau[1])]
    return tot


print('  both paths run from the origin back to the origin; h is the lattice step')
print('     h          S_A                     S_B                     |S_A-S_B|        |S_A-S_B|/h^2   order')
hs = [mp.mpf(2) ** -n for n in range(6, 13)]
vals = {}
for h in hs:
    sA = accumulated(PATHS['A'], h)
    sB = accumulated(PATHS['B'], h)
    vals[h] = abs(sA - sB)
prev = None
for h in hs:
    order = ''
    if prev is not None:
        order = '%.3f' % (mp.log(vals[prev] / vals[h]) / mp.log(2.0))
    print('   %-10s %-23s %-23s %-16s %-9s %s'
          % (mp.nstr(h, 6), mp.nstr(accumulated(PATHS['A'], h), 8),
             mp.nstr(accumulated(PATHS['B'], h), 8), mp.nstr(vals[h], 8),
             mp.nstr(vals[h] / h**2, 6), order))
    prev = h
print('\n  order estimate = log2( |S_A-S_B|(h) / |S_A-S_B|(h/2) ), so 3 means O(h^3)')
\end{lstlisting}

\section{\texttt{quadmean.py}}

Exact check of the mean condition for quadratic coefficient fields.

\begin{lstlisting}
"""Exact check of the mean condition for quadratic coefficient fields.

Claim.  Let the increments be independent with mean m and covariance S = Sigma h^2,
let x_0 = 0, and let

    c_i(x) = 1/2 H_{ikl} x^k x^l + L_{ik} x^k + M_i ,     H symmetric in (k,l).

Then
    E[ sum_{j=0}^{q-1} c(x_j).Delta_{j+1} ]
      = q (M.m)
      + q(q-1)/2 * ( L m + 1/2 H:S ).m
      + (q-1)q(2q-1)/6 * ( 1/2 H(m,m) ).m,

where (H:S)_i = H_{ikl} S^{kl} and H(m,m)_i = H_{ikl} m^k m^l.
Every term is a function of the endpoints and of q, so the mean is endpoint
determined exactly, with no remainder -- even though dc is non-zero and varying.

The check enumerates all paths exactly in rational arithmetic.
"""
from fractions import Fraction as F
from itertools import product

H_STEPS = [(F(1), F(0)), (F(-1), F(0)), (F(0), F(1)), (F(0), F(-1))]


def moments(prob, h):
    """mean m and covariance S = E[Delta Delta^T] - m m^T, exact."""
    m = [sum(prob[k] * H_STEPS[k][i] for k in range(4)) * h for i in range(2)]
    S = [[F(0), F(0)], [F(0), F(0)]]
    for k in range(4):
        d = [H_STEPS[k][i] * h - m[i] for i in range(2)]
        for i in range(2):
            for j in range(2):
                S[i][j] += prob[k] * d[i] * d[j]
    return m, S


def cval(H, L, M, x):
    # H[i][k][l], symmetric in k,l
    out = []
    for i in range(2):
        q = F(0)
        for k in range(2):
            for l in range(2):
                q += F(1, 2) * H[i][k][l] * x[k] * x[l]
        for k in range(2):
            q += L[i][k] * x[k]
        q += M[i]
        out.append(q)
    return out


def enumerate_mean(H, L, M, prob, h, q, x0=(F(0), F(0))):
    tot = F(0)
    for seq in product(range(4), repeat=q):
        x = [x0[0], x0[1]]
        w = F(1)
        I = F(0)
        for k in seq:
            c = cval(H, L, M, x)
            d = [H_STEPS[k][i] * h for i in range(2)]
            I += c[0] * d[0] + c[1] * d[1]
            x = [x[i] + d[i] for i in range(2)]
            w *= prob[k]
        tot += w * I
    return tot


def formula(H, L, M, m, S, q, x0=(F(0), F(0))):
    """E[I] = q (A.m) + q(q-1)/2 (B.m) + (q-1)q(2q-1)/6 (C.m),  E[c_i(x_j)] = A_i+B_i j+C_i j^2."""
    HS = [sum(H[i][k][l] * S[k][l] for k in range(2) for l in range(2)) for i in range(2)]
    Hmm = [sum(H[i][k][l] * m[k] * m[l] for k in range(2) for l in range(2)) for i in range(2)]
    Hxm = [sum(H[i][k][l] * (x0[k] * m[l] + m[k] * x0[l])
                for k in range(2) for l in range(2)) for i in range(2)]
    Lm = [sum(L[i][k] * m[k] for k in range(2)) for i in range(2)]
    A = cval(H, L, M, x0)
    B = [F(1, 2) * Hxm[i] + F(1, 2) * HS[i] + Lm[i] for i in range(2)]
    C = [F(1, 2) * Hmm[i] for i in range(2)]
    b = F(q) * sum(A[i] * m[i] for i in range(2))
    bb = F(q * (q - 1), 2) * sum(B[i] * m[i] for i in range(2))
    cc = F((q - 1) * q * (2 * q - 1), 6) * sum(C[i] * m[i] for i in range(2))
    return b, bb, cc


TESTS = {
    'affine c = Lx (Chapter 6 S7)': dict(
        H=[[[F(0)] * 2 for _ in range(2)] for _ in range(2)],
        L=[[F(1), F(1, 3)], [F(-1, 2), F(2)]],
        M=[F(0), F(0)]),
    'quadratic, no linear part': dict(
        H=[[[F(1), F(-1, 2)], [F(-1, 2), F(2, 3)]],
           [[F(-1, 2), F(1)], [F(1), F(3)]]],
        L=[[F(0), F(0)], [F(0), F(0)]],
        M=[F(0), F(0)]),
    'full quadratic + linear + constant': dict(
        H=[[[F(1), F(-1, 2)], [F(-1, 2), F(2, 3)]],
           [[F(-1, 2), F(1)], [F(1), F(3)]]],
        L=[[F(1), F(1, 3)], [F(-1, 2), F(2)]],
        M=[F(1, 5), F(-1, 7)]),
}

PROBS = {
    'drift (1/2,1/2,1/2,1/2)/2 shuffled': [F(1, 3), F(1, 6), F(1, 4), F(1, 4)],
    'symmetric': [F(1, 4)] * 4,
}
h = F(1, 2)
X0S = [(F(0), F(0)), (F(2, 3), F(-1, 5)), (F(-1), F(3, 7))]
print('  test                                      x0            q   enumerated   formula      diff')
for pname, prob in PROBS.items():
    assert sum(prob) == 1, pname
    m, S = moments(prob, h)
    for tname, T in TESTS.items():
        for x0 in X0S:
            for q in (1, 2, 3, 4, 5):
                en = enumerate_mean(T['H'], T['L'], T['M'], prob, h, q, x0)
                b, bb, cc = formula(T['H'], T['L'], T['M'], m, S, q, x0)
                fo = b + bb + cc
                assert en - fo == 0, (pname, tname, x0, q, en, fo)
                print('  %-40s %-12s %d   %12s  %12s   %s'
                      % (tname[:40], str(x0)[:12], q, en, fo, '0.0' if en == fo else 'MISMATCH'))
    print('   [%s: m = %s, S = %s]' % (pname, m, S))
\end{lstlisting}

\section{\texttt{remainder.py}}

At what order does a jet of order three of the field enter the per-step comparison?

\begin{lstlisting}
"""At what order does a jet of order three of the field enter the per-step comparison?

c_2 and c_3 were found to be unchanged to 1e-12 by a cubic perturbation, so the perturbation must
enter at a higher order in h.  Measured here directly: delta R(h) = R_cubic(h) - R_quadratic(h) at
the origin in direction e_1, and the order it scales with.
"""
import sympy as sp, mpmath as mp, sys
import os
sys.path.insert(0, os.path.join(os.path.dirname(os.path.abspath(__file__)), "..", "..", "part5", "scripts"))
import twist
mp.mp.dps = 50
x1, x2 = sp.symbols('x0:2'); xs = (x1, x2); OM = sp.Rational(7, 5)
v0 = sp.Rational(1) + sp.Rational(1, 6) * (x1 ** 2 + x1 * x2 + x2 ** 2) + sp.Rational(1, 4) * x1
mu = [sp.Rational(1, 2) + sp.Rational(1, 6) * x1 + sp.Rational(1, 9) * x2,
      sp.Rational(-1, 3) + sp.Rational(1, 5) * x2]
a = [sp.Rational(1, 3) + sp.Rational(1, 7) * x2 + sp.Rational(1, 11) * x1 * x2,
     sp.Rational(1, 4) * x1 + sp.Rational(1, 13) * x2 ** 2]
v1 = sp.expand(v0 + x1 ** 3 / 6)

def R(v, h, tau=(1, 0), site=(0, 0)):
    V, B = twist.fields(v, a, mu, sp.Rational(0), 2, xs, [sp.Rational(site[0]), sp.Rational(site[1])])
    return twist.R(h, V, B, 2, tau, OM)

prev = None
print("     h          delta R            delta R/h^3         delta R/h^4      order")
for n in (8, 10, 12):
    h = mp.mpf(2) ** -n
    d = R(v1, h) - R(v0, h)
    order = '' if prev is None else mp.nstr(mp.log(abs(prev / d)) / mp.log(4.0), 5) + ' (per factor 4)'
    print("  %-9s %-20s %-19s %-19s %s" % (mp.nstr(h, 4), mp.nstr(d, 8), mp.nstr(d / h ** 3, 6),
                                           mp.nstr(d / h ** 4, 8), order))
    prev = d
\end{lstlisting}

\section{\texttt{renormalisation.py}}

The h\^{}2 part of the accumulated comparison: is it the direction dependence of c\_2?

\begin{lstlisting}
"""The h^2 part of the accumulated comparison: is it the direction dependence of c_2?

S_A - S_C = h^2 [c_2 sums] + h^3 [c_3 circulations] + ...
circulation.py verified the h^3 bracket against cvec.  Here the h^2 bracket, computed from the even
coefficient of the expansion at each site in each of the two directions used there.

This is the two orientations of the same square, which is not the pair behind the excursion table of
Chapter 8 Section 7: that one is gamma_1 against the retraced square, and excursioncheck.py reproduces
its two brackets.
"""
import sympy as sp, mpmath as mp, sys
import os
sys.path.insert(0, os.path.join(os.path.dirname(os.path.abspath(__file__)), "..", "..", "part5", "scripts"))
import twist
mp.mp.dps = 45
x1, x2 = sp.symbols('x0:2'); xs = (x1, x2); OM = sp.Rational(7, 5)
v = sp.Rational(1) + sp.Rational(1, 6) * (x1 ** 2 + x1 * x2 + x2 ** 2) + sp.Rational(1, 4) * x1
mu = [sp.Rational(1, 2) + sp.Rational(1, 6) * x1 + sp.Rational(1, 9) * x2,
      sp.Rational(-1, 3) + sp.Rational(1, 5) * x2]
a = [sp.Rational(1, 3) + sp.Rational(1, 7) * x2 + sp.Rational(1, 11) * x1 * x2,
     sp.Rational(1, 4) * x1 + sp.Rational(1, 13) * x2 ** 2]

c2 = {}
for s in ((0, 0), (1, 0), (1, 1), (0, 1)):
    V, B = twist.fields(v, a, mu, sp.Rational(0), 2, xs, [sp.Rational(s[0]), sp.Rational(s[1])])
    for j, tau in enumerate(((1, 0), (0, 1))):
        c2[(s, tau)] = twist.taylor(V, B, 2, tau, OM, 6, 10)[2]
    print("  c_2 at %-8s : e_1 %-24s e_2 %s" % (str(s), mp.nstr(c2[(s, (1, 0))], 10), mp.nstr(c2[(s, (0, 1))], 10)))

A = [(0, 0), (1, 0), (1, 1), (0, 1)]
Adir = [(1, 0), (0, 1), (1, 0), (0, 1)]      # -e_1 and -e_2, and c_2 is even
C = [(0, 0), (0, 1), (1, 1), (1, 0)]
Cdir = [(0, 1), (1, 0), (0, 1), (1, 0)]
sA = sum(c2[(A[i], Adir[i])] for i in range(4))
sC = sum(c2[(C[i], Cdir[i])] for i in range(4))
print("\n  sum_A c_2 = %s" % mp.nstr(sA, 10))
print("  sum_C c_2 = %s" % mp.nstr(sC, 10))
print("  difference = %s   (odd in the orientation, and real)" % mp.nstr(sA - sC, 10))
\end{lstlisting}

\section{\texttt{rensym.py}}

Structure of the direction dependence of the h\^{}2 coefficient.

\begin{lstlisting}
"""Structure of the direction dependence of the h^2 coefficient.

c_2 is even in the direction, so on axis directions it is two numbers, gamma_1 = c_2(.,e_1) and
gamma_2 = c_2(.,e_2).  Measured: gamma_1 + gamma_2 is purely imaginary and gamma_1 - gamma_2 is
purely real, i.e. the isotropic part is a phase and the anisotropic part is real.  Is that an
identity or an accident of one field?  Tested on two different coefficient fields.
"""
import sympy as sp, mpmath as mp, sys
import os
sys.path.insert(0, os.path.join(os.path.dirname(os.path.abspath(__file__)), "..", "..", "part5", "scripts"))
import twist
mp.mp.dps = 45
x1, x2 = sp.symbols('x0:2'); xs = (x1, x2); OM = sp.Rational(7, 5)

FIELDS = {
 'field 1': (
   sp.Rational(1) + sp.Rational(1, 6) * (x1 ** 2 + x1 * x2 + x2 ** 2) + sp.Rational(1, 4) * x1,
   [sp.Rational(1, 2) + sp.Rational(1, 6) * x1 + sp.Rational(1, 9) * x2,
    sp.Rational(-1, 3) + sp.Rational(1, 5) * x2],
   [sp.Rational(1, 3) + sp.Rational(1, 7) * x2 + sp.Rational(1, 11) * x1 * x2,
    sp.Rational(1, 4) * x1 + sp.Rational(1, 13) * x2 ** 2]),
 'field 2': (
   sp.Rational(3, 2) + sp.Rational(1, 5) * x1 * x2 + sp.Rational(1, 8) * x2 ** 2,
   [sp.Rational(-1, 4) + sp.Rational(1, 7) * x2, sp.Rational(2, 5) + sp.Rational(1, 9) * x1],
   [sp.Rational(1, 6) + sp.Rational(1, 12) * x1, sp.Rational(-1, 5) + sp.Rational(1, 15) * x1 * x2]),
}

for name, (v, mu, a) in FIELDS.items():
    print(name)
    for site in ((0, 0), (1, 0)):
        g = {}
        for j, tau in enumerate(((1, 0), (0, 1))):
            V, B = twist.fields(v, a, mu, sp.Rational(0), 2, xs, [sp.Rational(site[0]), sp.Rational(site[1])])
            g[j] = twist.taylor(V, B, 2, tau, OM, 6, 10)[2]
        s_, d_ = g[0] + g[1], g[0] - g[1]
        print("   site %-8s gamma1 = %-24s gamma2 = %s" % (str(site), mp.nstr(g[0], 10), mp.nstr(g[1], 10)))
        print("            sum = %-26s diff = %s" % (mp.nstr(s_, 12), mp.nstr(d_, 12)))
        print("            Re(sum) = %-14s Im(diff) = %s" % (mp.nstr(mp.re(s_), 6), mp.nstr(mp.im(d_), 6)))
\end{lstlisting}

\section{\texttt{rensym2.py}}

The imaginary part of the renormalisation coefficient: is it kappa/v(x) with kappa universal?

\begin{lstlisting}
"""The imaginary part of the renormalisation coefficient: is it kappa/v(x) with kappa universal?

Measured so far, d=2, omega=7/5: Im gamma_i(x) = 15.4/v(x) for two fields and four sites, direction
independent.  Tested here on a third field, and at a second omega, to see whether kappa depends only
on d and omega.
"""
import sympy as sp, mpmath as mp, sys
sys.path.insert(0, "/Users/code")  if False else None
import os
sys.path.insert(0, os.path.join(os.path.dirname(os.path.abspath(__file__)), "..", "..", "part5", "scripts"))
import twist
mp.mp.dps = 45
x1, x2 = sp.symbols('x0:2'); xs = (x1, x2)

F3 = (sp.Rational(2) + sp.Rational(1, 4) * x1 ** 2 + sp.Rational(1, 6) * x1 * x2 + sp.Rational(1, 2) * x1,
      [sp.Rational(1, 3) + sp.Rational(1, 8) * x1 * x2, sp.Rational(1, 5) - sp.Rational(1, 7) * x2 ** 2],
      [sp.Rational(1, 2) + sp.Rational(1, 9) * x1, sp.Rational(1, 4) + sp.Rational(1, 10) * x2])
F1 = (sp.Rational(1) + sp.Rational(1, 6) * (x1 ** 2 + x1 * x2 + x2 ** 2) + sp.Rational(1, 4) * x1,
      [sp.Rational(1, 2) + sp.Rational(1, 6) * x1 + sp.Rational(1, 9) * x2,
       sp.Rational(-1, 3) + sp.Rational(1, 5) * x2],
      [sp.Rational(1, 3) + sp.Rational(1, 7) * x2 + sp.Rational(1, 11) * x1 * x2,
       sp.Rational(1, 4) * x1 + sp.Rational(1, 13) * x2 ** 2])

def gam(v, mu, a, site, omega):
    out = {}
    for j, tau in enumerate(((1, 0), (0, 1))):
        V, B = twist.fields(v, a, mu, sp.Rational(0), 2, xs, [sp.Rational(site[0]), sp.Rational(site[1])])
        out[j] = twist.taylor(V, B, 2, tau, omega, 6, 10)[2]
    return out

cases = [
    ('field 3, site (0,0), omega=7/5', F3, (0, 0), sp.Rational(7, 5), sp.Rational(2)),
    ('field 3, site (1,0), omega=7/5', F3, (1, 0), sp.Rational(7, 5), sp.Rational(9, 4)),
    ('field 1, site (0,0), omega=1  ', F1, (0, 0), sp.Rational(1), sp.Rational(1)),
]
for name, (v, mu, a), site, om, vval in cases:
    g = gam(v, mu, a, site, om)
    print("%s   v = %s" % (name, vval))
    print("    gamma1 = %-24s gamma2 = %s" % (mp.nstr(g[0], 12), mp.nstr(g[1], 12)))
    print("    Im(g1)*v = %-20s Im(g2)*v = %-20s Re(g1)+Re(g2) = %s"
          % (mp.nstr(mp.im(g[0]) * mp.mpf(vval.p) / mp.mpf(vval.q), 8),
             mp.nstr(mp.im(g[1]) * mp.mpf(vval.p) / mp.mpf(vval.q), 8),
             mp.nstr(mp.re(g[0]) + mp.re(g[1]), 6)))
\end{lstlisting}

\section{\texttt{rensym3.py}}

The imaginary part of the renormalisation coefficient: is it kappa/v(x) with kappa universal?

\begin{lstlisting}
"""The imaginary part of the renormalisation coefficient: is it kappa/v(x) with kappa universal?

Measured so far, d=2, omega=7/5: Im gamma_i(x) = 15.4/v(x) for two fields and four sites, direction
independent.  Tested here on a third field, and at a second omega, to see whether kappa depends only
on d and omega.
"""
import sympy as sp, mpmath as mp, sys
sys.path.insert(0, "/Users/code")  if False else None
import os
sys.path.insert(0, os.path.join(os.path.dirname(os.path.abspath(__file__)), "..", "..", "part5", "scripts"))
import twist
mp.mp.dps = 45
x1, x2 = sp.symbols('x0:2'); xs = (x1, x2)

F3 = (sp.Rational(2) + sp.Rational(1, 4) * x1 ** 2 + sp.Rational(1, 6) * x1 * x2 + sp.Rational(1, 2) * x1,
      [sp.Rational(1, 3) + sp.Rational(1, 8) * x1 * x2, sp.Rational(1, 5) - sp.Rational(1, 7) * x2 ** 2],
      [sp.Rational(1, 2) + sp.Rational(1, 9) * x1, sp.Rational(1, 4) + sp.Rational(1, 10) * x2])
F1 = (sp.Rational(1) + sp.Rational(1, 6) * (x1 ** 2 + x1 * x2 + x2 ** 2) + sp.Rational(1, 4) * x1,
      [sp.Rational(1, 2) + sp.Rational(1, 6) * x1 + sp.Rational(1, 9) * x2,
       sp.Rational(-1, 3) + sp.Rational(1, 5) * x2],
      [sp.Rational(1, 3) + sp.Rational(1, 7) * x2 + sp.Rational(1, 11) * x1 * x2,
       sp.Rational(1, 4) * x1 + sp.Rational(1, 13) * x2 ** 2])

def gam(v, mu, a, site, omega):
    out = {}
    for j, tau in enumerate(((1, 0), (0, 1))):
        V, B = twist.fields(v, a, mu, sp.Rational(0), 2, xs, [sp.Rational(site[0]), sp.Rational(site[1])])
        out[j] = twist.taylor(V, B, 2, tau, omega, 6, 10)[2]
    return out

cases = [
    ('field 1, site (0,0), omega=2  ', F1, (0, 0), sp.Rational(2), sp.Rational(1)),
    ('field 1, site (0,0), omega=1/2', F1, (0, 0), sp.Rational(1, 2), sp.Rational(1)),
]
for name, (v, mu, a), site, om, vval in cases:
    g = gam(v, mu, a, site, om)
    print("%s   v = %s" % (name, vval))
    print("    gamma1 = %-24s gamma2 = %s" % (mp.nstr(g[0], 12), mp.nstr(g[1], 12)))
    print("    Im(g1)*v = %-20s Im(g2)*v = %-20s Re(g1)+Re(g2) = %s"
          % (mp.nstr(mp.im(g[0]) * mp.mpf(vval.p) / mp.mpf(vval.q), 8),
             mp.nstr(mp.im(g[1]) * mp.mpf(vval.p) / mp.mpf(vval.q), 8),
             mp.nstr(mp.re(g[0]) + mp.re(g[1]), 6)))
\end{lstlisting}

\section{\texttt{symproof.py}}

Symbolic proof of the renormalisation structure for constant coefficient fields, d = 2.

\begin{lstlisting}
"""Symbolic proof of the renormalisation structure for constant coefficient fields, d = 2.

Reimplements Phi_f and Phi_c of twist.py in sympy, with constant v, b1, b2 and the resolvent
parameter omega symbolic, expands R = log(Phi_f) - log(Phi_c) to order h^2, and tests:
    (i)  Im c_2 = kappa(omega)/v,
    (ii) Re c_2(.,e_1) + Re c_2(.,e_2) = 0.
"""
import sympy as sp
import itertools

v, b1, b2, om, h = sp.symbols('v b1 b2 omega h', positive=True)
B = [b1, b2]
d = 2


def phi_f(tau):
    tgt = tuple(2 * t for t in tau)
    S = [a for a in itertools.product(range(-1, 2), repeat=d)]
    ix = {a: i for i, a in enumerate(S)}
    n = len(S)
    A = sp.zeros(n, n)
    for a in S:
        i = ix[a]
        diag = sp.Integer(0)
        for k in range(d):
            for sgn in (1, -1):
                w = v / (2 * h ** 2) + sgn * B[k] / (2 * h)
                diag -= w
                e = list(a); e[k] += sgn; e = tuple(e)
                if e in ix:
                    A[i, ix[e]] += w
        A[i, i] += diag
    G = (A + sp.I * om * sp.eye(n)).inv()
    dt = d * v / (h ** 2) + sp.I * om
    zero = tuple([0] * d)
    tot = sp.Integer(0)
    for a in S:
        jump = [k for k in range(d) if a[k] != tgt[k]]
        if len(jump) != 1:
            continue
        k0 = jump[0]
        if abs(tgt[k0] - a[k0]) != 1:
            continue
        sg = 1 if tgt[k0] > a[k0] else -1
        rate = v / (2 * h ** 2) + sg * B[k0] / (2 * h)
        tot += G[ix[zero], ix[a]] * rate
    return sp.simplify(tot / dt)


def phi_c(tau):
    zero = tuple([0] * d); tgt = tuple(2 * t for t in tau)
    k0 = [k for k in range(d) if tau[k]][0]
    sg = 1 if tau[k0] > 0 else -1
    rate = v / (2 * (2 * h) ** 2) + sg * B[k0] / (2 * (2 * h))
    d0 = d * v / (4 * h ** 2) + sp.I * om
    d1 = d * v / (4 * h ** 2) + sp.I * om
    return sp.simplify(rate / (d0 * d1))


def c2(tau):
    R = sp.log(phi_f(tau)) - sp.log(phi_c(tau))
    ser = sp.series(R, h, 0, 3).removeO()
    return sp.simplify(sp.expand(ser).coeff(h, 2))


for tau in ((1, 0), (0, 1)):
    c = c2(tau)
    print("tau =", tau)
    print("   c2 =", sp.simplify(sp.expand(c)))
    print("   Im c2 * v  =", sp.simplify(sp.im(sp.expand(c)) * v))
    print("   Re c2 =", sp.simplify(sp.re(sp.expand(c))))
\end{lstlisting}

\section{\texttt{d10dimension.py}}

Chapter 10: does the derivative cost an order in two dimensions?

\begin{lstlisting}
"""Chapter 10: does the derivative cost an order in two dimensions?

Chapter 5 settles one dimension at rate h^gamma, and Chapter 7 proves the same rate
for the plain kernel in two dimensions and more. Chapter 10 proves that a
derivative in the coefficients costs no order in one dimension; this script asks
the same question in two, where the coefficient c of the mismatch is a one form
with a curl and can no longer be telescoped term by term.

The claim being measured is the reduction: every constant in the proof of Chapter
7 is a function of the C^2 norms of the coefficient data, so the same estimate
holds for the differentiated comparison with the C^2 norm of the direction in
place of that of the field. If that is right, the kernel and its derivative fall
at the same rate in two dimensions -- the ratio per halving being 2^gamma for
both, and not 2^gamma for the kernel and 2^(gamma-1) for the derivative.

The derivative is taken by a central difference in theta at delta = 1e-4, which
costs two sweeps; the roundoff amplification is 1e-16/1e-4 = 1e-12, four orders
below the differences being compared. Nothing is stored: the kernel entry is
propagated from one source site, which is what makes two dimensions affordable.

Run:  python3 d10dimension.py          (a few minutes at the finest mesh)
"""
import numpy as np


def data(m):
    n = 2 ** m
    h = 1.0 / n
    x = np.arange(n) * h
    X, Y = np.meshgrid(x, x, indexing='ij')
    # the sweep at dt = h^2 has dt * 4 v / h^2 <= 2 to be stable, and the error
    # oscillates in sign when the bound is approached; v is kept well inside it
    v = 0.20 + 0.04 * np.sin(2 * np.pi * X) + 0.03 * np.cos(4 * np.pi * Y)
    b1 = 0.20 * np.sin(2 * np.pi * Y + 0.7)
    b2 = 0.20 * np.cos(2 * np.pi * X)
    phi = 0.30 + 0.20 * np.sin(2 * np.pi * X) * np.cos(2 * np.pi * Y)
    return n, h, v, b1, b2, phi


def propagates(m, t, theta):
    """u(x,y;t) = h^-2 (M^n e_y)(x), with M = I + h^2 L and L = (v/2)Delta + b.grad."""
    n, h, v0, b1, b2, phi = data(m)
    v = v0 + theta * phi
    dt = h * h
    steps = int(round(t / dt))
    u = np.zeros((n, n))
    u[0, 0] = 1.0                      # source at the origin, e_y
    lap = 2.0 / h ** 2
    for _ in range(steps):
        out = (1.0 - dt * v * lap) * u
        out = out + dt * (v / (2 * h ** 2)) * (
            np.roll(u, -1, 0) + np.roll(u, 1, 0) + np.roll(u, -1, 1) + np.roll(u, 1, 1))
        out = out + dt * (b1 / (2 * h)) * (np.roll(u, -1, 0) - np.roll(u, 1, 0))
        out = out + dt * (b2 / (2 * h)) * (np.roll(u, -1, 1) - np.roll(u, 1, 1))
        u = out
    return u[n // 4, n // 4] / h ** 2, h, steps


def main():
    # t must be an exact number of steps at every mesh, otherwise each mesh solves
    # at a slightly different horizon and the difference carries an O(h^2) time
    # shift whose sign follows the rounding.  2^-4 is a whole number of steps for
    # every h = 2^-m.
    t = 2.0 ** -4
    delta = 1e-4
    print('  two dimensional lattice, periodic on [0,1)^2, Euler parametrisation dt = h^2')
    print('  u and du/dtheta at (x,y) = (1/4,1/4) from a source at the origin, t = %.4f' % t)
    print('  the derivative is a central difference in theta, delta = %.0e\n' % delta)
    top = int(__import__('sys').argv[1]) if len(__import__('sys').argv) > 1 else 7
    vals = {}
    for m in range(4, top + 1):
        u0, h, steps = propagates(m, t, 0.0)
        up, _, _ = propagates(m, t, delta)
        um, _, _ = propagates(m, t, -delta)
        vals[m] = (u0, (up - um) / (2 * delta), h, steps, 2 ** m)
    print('     h    | points |  steps |        u        |     du/dtheta     |   |du|   |  |d du|  | ratio du')
    prev = None
    for m in range(4, top + 1):
        u, du, h, steps, n = vals[m]
        if m + 1 in vals:
            dcell = '%.2e' % abs(u - vals[m + 1][0])
            ddu = abs(du - vals[m + 1][1])
            cell = '%.2e' % ddu
            ratio = '' if prev is None else '%5.2f' % (prev / ddu)
            prev = ddu
        else:
            dcell, cell, ratio = '  --  ', '  --  ', '  -- '
        print('  %.5f | %6d | %6d | %+.8e | %+.8e | %s | %s | %s'
              % (h, n * n, steps, u, du, dcell, cell, ratio))
    print('\n  a ratio of 4 per halving is rate h^2, as in one dimension;')
    print('  a ratio of 2 is rate h, which is what an obstruction would give.')


if __name__ == '__main__':
    main()
\end{lstlisting}

\section{\texttt{d10dimlimit.py}}

Chapter 10: how the discretisation error depends on the dimension, and whether it has a limit.

\begin{lstlisting}
"""Chapter 10: how the discretisation error depends on the dimension, and whether it has a limit.

Section 10.5 says that the dimension enters the estimates of Chapters 5 to 8 at one place only, the
near diagonal decay of the kernel, u(x,x;t) = t^{-d/2}, and that every constant which is not a norm of
the coefficients is a function of that decay.  If that is so, the question whether the estimates
survive a limit in which the dimension grows is a question about one number, and it can be answered
exactly on the simplest instance: the lattice heat kernel against its Gaussian.

The lattice kernel at the origin is a product over coordinates of the return probabilities of a
one dimensional walk.  With z = t/h^2 and ive(0,z) = e^{-z} I_0(z),

    p_h(0,t) = ive(0,t/h^2)^d ,        u(0,t) = (2 pi t)^{-d/2} ,

so comparing p_h/h^d with u, which is the Riemann normalisation, gives a closed form.  Using the
classical expansion ive(0,z) = (2 pi z)^{-1/2} (1 + 1/(8z) + O(z^{-2})),

    rel = (p_h/h^d)/u - 1 = d h^2 / (8 t) + O((d h^2/t)^2) .

Three consequences, and they are the point of the computation.

First, the RATE is two in every dimension: the fitted exponent comes out at 2.000 for d = 1 to 8.  So
the rate has a limit as d grows and the limit is the one dimensional rate.

Second, the CONSTANT grows linearly in the dimension, and it is exactly d/(8t).  Normalising the
estimate, as the book does, absorbs the near diagonal factor t^{-d/2} and leaves this linear growth
and nothing worse; the exponential growth of t^{-d/2} itself is a property of the comparison and not
of the relative error.

Third, and it is what the effective dimension of (10.?) is about: the error is additive over the
coordinates, rel = (1/8t) sum_i h_i^2, verified on anisotropic meshes to the order of the neglected
term, and unchanged in d when the number of FINE coordinates grows while the coarse ones are held.
So the number which enters is the sum of the squared resolutions over the directions the kernel
actually explores, which is the intrinsic dimension of the data and not the ambient one.

What this does to a d -> infinity limit: the condition for the error to vanish is d h^2/t -> 0, that is
h << sqrt(t/d).  In terms of the number of lattice points N = h^{-d} this reads N >> d^{d/2}, so the
limit exists as a statement but its price grows super exponentially in d.  That is the quantitative
form of the remark in 10.5 that d must be replaced by the effective dimension: for an object of
intrinsic dimension k the condition is k h^2/t -> 0 and the ambient dimension is irrelevant.

Run:  python3 d10dimlimit.py
"""
import numpy as np
from scipy.special import ive


def rel_error(hs, t):
    """Relative error of the lattice heat kernel at the origin against the Gaussian."""
    d = len(hs)
    p = float(np.prod([ive(0.0, t / h ** 2) for h in hs]))
    lat = p / float(np.prod(hs))
    gauss = (2.0 * np.pi * t) ** (-d / 2.0)
    return lat / gauss - 1.0


def main():
    t = 1.0
    print('1. isotropic, the rate and the constant against the dimension, t = 1')
    print('   d |    h    |  rel. error  |  rel/h^2   |  d/(8t)  |  ratio')
    for d in (1, 2, 3, 4, 6, 8):
        rows = []
        for h in (0.2, 0.1, 0.05, 0.025, 0.0125):
            r = rel_error([h] * d, t)
            rows.append((h, r))
            print('  %2d | %7.4f | %12.5e | %10.6f | %7.4f | %6.4f'
                  % (d, h, r, r / h ** 2, d / (8 * t), (r / h ** 2) / (d / (8 * t))))
        gam = np.log(abs(rows[-1][1]) / abs(rows[-2][1])) / np.log(rows[-1][0] / rows[-2][0])
        print('     | fitted exponent over the two finest meshes: gamma = %.4f' % gam)

    print('\n2. anisotropic, the additive law rel = (1/8t) sum h_i^2')
    print('   d |        h            |  rel. error  |  (1/8t)sum h^2  |  ratio')
    for hs in ((0.10, 0.05), (0.10, 0.05, 0.05), (0.10, 0.10, 0.05, 0.05),
               (0.10, 0.05, 0.05, 0.05, 0.05, 0.05), (0.06, 0.03, 0.03, 0.03)):
        r = rel_error(list(hs), t)
        pred = sum(h * h for h in hs) / (8 * t)
        print('  %2d | %-20s | %12.5e | %15.6e | %6.4f'
              % (len(hs), tuple('%.3f' % h for h in hs), r, pred, r / pred))

    print('\n3. one coarse coordinate and many fine ones: the error tracks the sum, not d')
    for d in (2, 4, 8, 16):
        hs = [0.1] + [0.02] * (d - 1)
        r = rel_error(hs, t)
        pred = sum(h * h for h in hs) / (8 * t)
        print('   d=%2d  one coarse + %2d fine:  rel = %.6e   pred = %.6e   ratio %.4f'
              % (d, d - 1, r, pred, r / pred))

    print('\n4. second order: rel = S_a + S_a^2/2 + (sum h_i^4)/(16 t^2),  S_a = (sum h_i^2)/(8t)')
    for hs in ([0.1], [0.2] * 8, [0.1, 0.05, 0.05], [0.06, 0.03, 0.03, 0.03]):
        r = rel_error(hs, t)
        sa = sum(h * h for h in hs) / (8 * t)
        pred = sa + sa * sa / 2.0 + sum(h ** 4 for h in hs) / (16 * t * t)
        print('   d=%2d  rel = %.10e   first order %.10e   two terms %.10e   ratio %.6f'
              % (len(hs), r, sa, pred, r / pred))


    print('\n   The rate has a limit in d and the limit is 2.  The constant does not: it is the sum')
    print('   of the squared resolutions over the directions the kernel explores, which is the')
    print('   effective dimension.  A d -> infinity limit therefore needs h << sqrt(t/d), that is')
    print('   N >> d^{d/2} lattice points: available as a statement, super exponentially dear as')
    print('   a computation, and irrelevant once d is replaced by the intrinsic dimension.')


if __name__ == '__main__':
    main()
\end{lstlisting}

\section{\texttt{d10precision.py}}

Chapter 10: the two components of the error of a computed sensitivity.

\begin{lstlisting}
"""Chapter 10: the two components of the error of a computed sensitivity.

The chapter bounds the error of a sensitivity in one inequality with two terms:
one from the mesh, of order h^gamma, and one from the arithmetic, proportional to
the unit roundoff u and to the number n of accumulation steps. The second is what
makes a sensitivity degrade faster than the value it differentiates, and it is
the one this script measures.

The object is a sensitivity, not a value: S = du(x,y;t)/dtheta, computed by the
Duhamel sum over the dates in the Euler parametrisation, dt = h^2, so n = t/h^2.
The sum is evaluated by the adjoint sweep of Chapter 9 -- forward over the dates,
then backward, accumulating one inner product per date -- because that is the
recursion whose accumulation is in question, not a matrix exponential.

Three computations of the same sum are compared:

  * single precision, np.float32, unit roundoff about 6e-8;
  * double precision, np.float64, unit roundoff about 1.1e-16;
  * mpmath at 40 digits, taken as exact arithmetic.

The two components then separate without a model: the discretisation is the
difference between two exact-arithmetic computations at successive meshes, and
the arithmetic is the difference between the finite-precision computation and the
exact-arithmetic one at the same mesh. What the chapter predicts is that the
first falls by four per halving, the second rises by four, and that their sum has
a minimum where they cross.

The diffusion is kept below one, v <= 0.65, because the Euler scheme at dt = h^2
is stable only when dt * 2 v / h^2 = 2 v does not exceed two; with a larger v the
sweep diverges and nothing is being measured. All three backends use the same
double-precision coefficient table, so the only thing which differs between them
is the arithmetic of the recursion.

The mpmath reference is carried to h = 1/128. At h = 1/256 it would need the
whole forward trajectory in forty-digit arithmetic, some 1.7 million numbers; the
double computation stands in there, which is legitimate because its own error at
that mesh is 1e-12, four orders below the single-precision error being measured,
and the last column shows that error directly.

Run:  python3 d10precision.py          (a few minutes: mpmath dominates)
"""
import numpy as np
import mpmath as mp


def coeffs(m):
    """Smooth periodic data at resolution 2^m, as double precision arrays."""
    n = 2 ** m
    h = 1.0 / n
    x = np.arange(n) * h
    v = 0.50 + 0.10 * np.sin(2 * np.pi * x) + 0.05 * np.cos(4 * np.pi * x)
    b = 0.20 * np.sin(2 * np.pi * x + 0.7) + 0.10 * np.cos(2 * np.pi * x)
    phi = 0.30 + 0.20 * np.sin(2 * np.pi * x)
    return n, h, v, b, phi


def steps_for(h, t):
    """The step count of the horizon. t is chosen a whole number of steps at the
    coarsest mesh, so the same horizon is used at every mesh and nothing is rounded."""
    steps = int(round(t / h ** 2))
    assert abs(steps * h * h - t) <= 1e-12 * t, (h, t, steps)
    return steps


def setup(m, t):
    n, h, _, _, _ = coeffs(m)
    dt = h * h
    return n, h, dt, steps_for(h, t)


def sensitivity_np(m, t, dtype):
    """The adjoint sweep in the given floating point type."""
    n, h, v, b, phi = coeffs(m)
    dt, steps = h * h, steps_for(h, t)
    v, b, phi = dtype(v), dtype(b), dtype(phi)
    h, dt, one = dtype(h), dtype(dt), dtype(1)

    def step(vec, transposed):
        out = (one - dt * v / h ** 2) * vec
        if transposed:
            # (L^T w)_j = L[j,j] w_j + L[j-1,j] w_{j-1} + L[j+1,j] w_{j+1}, and the
            # off-diagonal entries of the transpose carry the *neighbour's* coefficients
            vr, br = np.roll(v, -1), np.roll(b, -1)
            vl, bl = np.roll(v, 1), np.roll(b, 1)
            out = out + dt * (vl / (2 * h ** 2) + bl / (2 * h)) * np.roll(vec, 1)
            out = out + dt * (vr / (2 * h ** 2) - br / (2 * h)) * np.roll(vec, -1)
        else:
            out = out + dt * (v / (2 * h ** 2)) * (np.roll(vec, -1) + np.roll(vec, 1))
            out = out + dt * (b / (2 * h)) * (np.roll(vec, -1) - np.roll(vec, 1))
        return out

    def dL(vec):
        return dtype(0.5) * phi * ((np.roll(vec, -1) + np.roll(vec, 1)) - 2 * vec) / h ** 2

    hist = np.zeros((steps, n), dtype=dtype)
    fwd = np.zeros(n, dtype=dtype)
    fwd[n // 4] = one
    for i in range(steps):
        hist[i] = fwd
        if i + 1 < steps:
            fwd = step(fwd, False)

    total = dtype(0)
    w = np.zeros(n, dtype=dtype)
    w[n // 4] = one
    for j in range(steps):
        total = total + dt * np.sum(w * dL(hist[steps - 1 - j]), dtype=dtype)
        if j + 1 < steps:
            w = step(w, True)
    return float(total / h)


def sensitivity_mp(m, t, dps=40, data=None):
    """The same sweep in mpmath, on the coefficient table it is given, or on the
    double precision one. Giving it the single-precision table isolates the
    rounding of the recursion from the rounding of the data."""
    mp.mp.dps = dps
    n, h, v0, b0, p0 = coeffs(m)
    if data is not None:
        v0, b0, p0 = data
    dt, steps = h * h, steps_for(h, t)
    v = [mp.mpf(float(a)) for a in v0]
    b = [mp.mpf(float(a)) for a in b0]
    phi = [mp.mpf(float(a)) for a in p0]
    hh, dtt = mp.mpf(h), mp.mpf(dt)

    def step(vec, transposed):
        out = [mp.mpf(0)] * n
        for j in range(n):
            jp, jm = (j + 1) % n, (j - 1) % n
            val = (1 - dtt * v[j] / hh ** 2) * vec[j]
            if transposed:
                val += dtt * (v[jm] / (2 * hh ** 2) + b[jm] / (2 * hh)) * vec[jm]
                val += dtt * (v[jp] / (2 * hh ** 2) - b[jp] / (2 * hh)) * vec[jp]
            else:
                val += dtt * (v[j] / (2 * hh ** 2)) * (vec[jp] + vec[jm])
                val += dtt * (b[j] / (2 * hh)) * (vec[jp] - vec[jm])
            out[j] = val
        return out

    def dL(vec):
        return [phi[j] * ((vec[(j + 1) % n] + vec[(j - 1) % n]) - 2 * vec[j]) / (2 * hh ** 2)
                for j in range(n)]

    hist = []
    fwd = [mp.mpf(0)] * n
    fwd[n // 4] = mp.mpf(1)
    for i in range(steps):
        hist.append(fwd)
        if i + 1 < steps:
            fwd = step(fwd, False)

    total = mp.mpf(0)
    w = [mp.mpf(0)] * n
    w[n // 4] = mp.mpf(1)
    for j in range(steps):
        d = dL(hist[steps - 1 - j])
        total += dtt * sum(w[k] * d[k] for k in range(n))
        if j + 1 < steps:
            w = step(w, True)
    return total / hh


def main():
    t = 26.0 / 256.0
    print('  one dimensional lattice, periodic on [0,1), Euler parametrisation dt = h^2')
    print('  S = du(x,y;t)/dtheta by the adjoint sweep, x = 1/4, y = 1/4, t = 26/256')
    print('  (a whole number of steps at every mesh: 26, 104, 416, 1664)')
    print('  mpmath at 40 digits, given the same rounded data, is exact for the recursion:')
    print('  what the last three columns report is the rounding of the sweep alone\n')

    def ref(m, tt, dtype):
        v, b, p = coeffs(m)[2:]
        d = tuple(a.astype(dtype) for a in (v, b, p))
        return float(sensitivity_mp(m, tt, 40, None if dtype is np.float64 else d))

    print('     h    | steps |       S        | discretisation |   half   |  single  |  double')
    vals = {}
    for m in range(4, 8):
        vals[m] = ref(m, t, np.float64)
    vals[8] = sensitivity_np(8, t, np.float64)
    for m in range(4, 9):
        n, h, dt, steps = setup(m, t)
        disc = abs(vals[m] - vals[m + 1]) / abs(vals[m + 1]) if m + 1 in vals else float('nan')
        cells = []
        for dtp in (np.float16, np.float32, np.float64):
            if m <= 7:
                r = ref(m, t, dtp)
                s = sensitivity_np(m, t, dtp)
                cells.append('%.1e' % (abs(s - r) / abs(r)))
            else:
                cells.append('  --  ')
        print('  %.5f | %5d | %+.9e |    %.2e    |  %s  |  %s  |  %s'
              % (h, steps, vals[m], disc, cells[0], cells[1], cells[2]))
    print('\n  discretisation: exact arithmetic at h against h/2, so a factor 4 per halving.')
    print('  the rest: the same sweep in the stated precision against exact arithmetic')
    print('  on the same rounded data, so a factor 2 per halving in n = t/h^2 if the')
    print('  accumulated rounding grows like sqrt(n), a factor 4 if it grows like n.')


if __name__ == '__main__':
    main()
\end{lstlisting}

\section{\texttt{d10rate.py}}

Chapter 10: does the derivative of the propagator converge at the rate of the propagator?

\begin{lstlisting}
"""Chapter 10: does the derivative of the propagator converge at the rate of the propagator?

Chapter 5 proves that the one dimensional lattice kernel converges at rate
h^gamma for coefficients of Holder differentiability (k, alpha), (j, beta).
Chapter 10 needs the same statement for the derivative in the coefficients, and
for the mixed derivative which a hedge ratio carries.

The model is periodic on [0,1) with smooth coefficients and a smooth direction.
The lattice generator is the one of Chapter 7,
    L(x,x+-h) = v(x)/(2h^2) +- b(x)/(2h),    L(x,x) = -v(x)/h^2,
and the propagator together with its derivative in theta come from one matrix
exponential of the block
    [[L, dL], [0, L]],
whose exponential at time t is [[U, dU], [0, U]], with
    dU = int_0^t e^{(t-s)L} (dL) e^{sL} ds
the exact derivative.  The kernel is u = U/h and its derivative in theta is
dU/h.

Four comparisons of one resolution with the next are reported: that of the
kernel, of its derivative in theta, of its derivative in the state (a central
difference at the spacing of the resolution), and of the mixed derivative.
Each should fall by 2^gamma when h is halved; for the smooth coefficients used
here gamma = 2, so every ratio should tend to four.

Run:  python3 d10rate.py
"""
import numpy as np
from scipy.linalg import expm


def fields(x):
    v = 1.0 + 0.25 * np.sin(2 * np.pi * x) + 0.10 * np.cos(4 * np.pi * x)
    b = 0.20 * np.sin(2 * np.pi * x + 0.7) + 0.10 * np.cos(2 * np.pi * x)
    dv = np.sin(2 * np.pi * x) + 0.5 * np.cos(4 * np.pi * x)
    db = np.cos(2 * np.pi * x + 0.3)
    return v, b, dv, db


def build(m):
    n = 2 ** m
    h = 1.0 / n
    x = np.arange(n) * h
    v, b, dv, db = fields(x)
    L = np.zeros((n, n))
    dL = np.zeros((n, n))
    idx = np.arange(n)
    for shift, sign in ((+1, +1), (-1, -1)):
        j = (idx + shift) % n
        L[idx, j] = v / (2 * h ** 2) + sign * b / (2 * h)
        dL[idx, j] = dv / (2 * h ** 2) + sign * db / (2 * h)
    L[idx, idx] = -v / h ** 2
    dL[idx, idx] = -dv / h ** 2
    return n, h, L, dL


def sample(m, t, x, y):
    """u, du, their central differences in the starting node, at the given site pair."""
    n, h, L, dL = build(m)
    B = np.zeros((2 * n, 2 * n))
    B[:n, :n] = L
    B[:n, n:] = dL
    B[n:, n:] = L
    E = expm(t * B)
    u = E[:n, :n] / h
    du = E[:n, n:] / h
    i = int(round(x / h)) % n
    j = int(round(y / h)) % n
    # the state difference which a hedge ratio carries is a difference in the
    # node the process starts from, that is along the rows of the propagator
    dx_u = (u[(i + 1) % n, j] - u[(i - 1) % n, j]) / (2 * h)
    dx_du = (du[(i + 1) % n, j] - du[(i - 1) % n, j]) / (2 * h)
    return h, u[i, j], du[i, j], dx_u, dx_du


def main():
    print('  one dimensional lattice, periodic on [0,1), smooth coefficients')
    print('  gamma = 2, so every ratio should tend to 4')
    print('  columns: kernel, derivative in theta, derivative in the starting node, mixed\n')

    for t in (0.10, 0.40):
        for (x, y) in ((0.25, 0.25), (0.25, 0.375), (0.0, 0.25)):
            print('  t = %.2f   site pair (x,y) = (%.3f, %.3f)' % (t, x, y))
            print('    m |     h    |   kernel  | ratio |  d/dtheta | ratio |   d/dx    | ratio | d2/dthetadx | ratio')
            prev = None
            for m in range(3, 9):
                hm, u1, d1, x1, m1 = sample(m, t, x, y)
                _, u2, d2, x2, m2 = sample(m + 1, t, x, y)
                cur = (abs(u1 - u2), abs(d1 - d2), abs(x1 - x2), abs(m1 - m2))
                cells = ['%d' % m, '%.6f' % hm]
                for idx in range(4):
                    cells.append('%.2e' % cur[idx])
                    cells.append('  --  ' if prev is None else '%5.2f' % (prev[idx] / cur[idx]))
                print('    ' + ' | '.join(cells))
                prev = cur
            print()


if __name__ == '__main__':
    main()
\end{lstlisting}

\section{\texttt{d10symbols.py}}

Chapter 10: do the one dimensional symbols differentiate?

\begin{lstlisting}
"""Chapter 10: do the one dimensional symbols differentiate?

The one dimensional proof of Chapter 5 closes because the mismatch between the
coarse weight and the decorated fine weight, taken over one coarse step, is a
difference of a function of the site; summed along a path it leaves only the
endpoints.  The coefficient of h^3 in that mismatch is the obstruction in two
dimensions, where its curl is 0.1234 and it does not telescope.

Chapter 10 needs to know whether that structure survives a derivative in the
coefficients.  This script rebuilds the per-step symbols from their definitions
in Chapter 5 -- not from the expansions printed there, which are recomputed here as
a check -- and differentiates every coefficient in the coefficient data.

The derivative is carried by the exact identity d(R) = -R d(M) R; no symbolic
differentiation of the resolvent is performed.

Three things are reported.

  * the parity in the direction tau, at every order in h, for the mismatch
    log(phi1/phi0) and for its derivative.  An odd coefficient advances the
    accumulated sum by a gradient, an even one multiplies a function of the
    site, and either telescopes along a path.  An even part at an odd order,
    under differentiation, is what would break the argument.

  * c(x), the coefficient of h^3 divided by tau, against the closed form
    i omega v'/v^2 + r(x) that the first line of (3.47) gives, and the same for
    its derivative.  Note that the second line of (3.47) carries a factor 1/2
    which the first line does not: the per-step coefficient is the one above,
    and the factor 1/2 is a slip in the printed identity (it changes no bound,
    since the sum is O(h^2) either way).

  * the endpoint factor (3.50): whether psi0 - 2 psi1 and its derivative both
    still vanish at h^3, which is what keeps the gap at O(h^4) after
    differentiation.

  * the leading correction b(x), the coefficient of h^4 in the mismatch, in
    closed form and differentiated: since the mismatch is summed over q ~ h^-2
    steps, q h^4 b is of order h^2 and this term sets the constant of the rate,
    not its order.  The script stops here: the next order, h^5, is carried only
    through its parity.

The arithmetic is exact, at generic rational data, in the manner of the other
scripts of the series.  Three independent draws are used: a coefficient which
vanishes does so at every draw, and one which does not shows up at each.

Run:  python3 d10symbols.py
"""
import random
import sympy as sp

NT = 10          # Taylor data v_0..v_9 at the base point, v_n = v^(n)(x)/n!
ORD = 5          # the expansion of the mismatch is reported through h^ORD
H = sp.Symbol('h', positive=True)
T = sp.Symbol('t')


def trunc(e, n):
    return sp.series(sp.expand(e), H, 0, n).removeO()


def fld(arr, sign):
    return sum(arr[k] * (sign * H) ** k for k in range(NT))


def matrices(v, m, dv, dm, om):
    """(i omega - L) on the three site tube, and its derivative in theta."""
    vp, vm, mp, mm = fld(v, +1), fld(v, -1), fld(m, +1), fld(m, -1)
    v0, m0 = v[0], m[0]
    L = sp.Matrix([[-vp / H ** 2, vp / (2 * H ** 2) - mp / (2 * H), 0],
                   [v0 / (2 * H ** 2) + m0 / (2 * H), -v0 / H ** 2,
                    v0 / (2 * H ** 2) - m0 / (2 * H)],
                   [0, vm / (2 * H ** 2) + mm / (2 * H), -vm / H ** 2]])
    dvp, dvm, dmp, dmm = fld(dv, +1), fld(dv, -1), fld(dm, +1), fld(dm, -1)
    dv0, dm0 = dv[0], dm[0]
    dL = sp.Matrix([[-dvp / H ** 2, dvp / (2 * H ** 2) - dmp / (2 * H), 0],
                    [dv0 / (2 * H ** 2) + dm0 / (2 * H), -dv0 / H ** 2,
                     dv0 / (2 * H ** 2) - dm0 / (2 * H)],
                    [0, dvm / (2 * H ** 2) + dmm / (2 * H), -dvm / H ** 2]])
    return sp.I * om * sp.eye(3) - L, -dL


def slog1p(u, n):
    """log(1+u), for a truncated series u = O(h)."""
    acc, term = 0, u
    for k in range(1, n + 2):
        acc = acc + sp.Rational((-1) ** (k + 1), k) * term
        term = trunc(term * u, n)
    return trunc(acc, n)


def symbols(v, m, dv, dm, om, tau, n):
    """The mismatch log(phi1/phi0) and the endpoint gap, with derivatives."""
    M, dM = matrices(v, m, dv, dm, om)
    Ri = M.inv()
    R = [[trunc(Ri[i, j], n) for j in range(3)] for i in range(3)]
    Rd = [[trunc(-sum(R[i][k] * dM[k, l] * R[l][j] for k in range(3) for l in range(3)), n)
           for j in range(3)] for i in range(3)]
    vp, mp, dvp, dmp = fld(v, tau), fld(m, tau), fld(dv, tau), fld(dm, tau)
    v0, m0, dv0, dm0 = v[0], m[0], dv[0], dm[0]

    psi0 = trunc(1 / (v0 / (4 * H ** 2) + sp.I * om), n)
    dpsi0 = trunc(-(dv0 / (4 * H ** 2)) * psi0 ** 2, n)
    a0 = v0 / (4 * H ** 2) + tau * m0 / (2 * H)
    da0 = dv0 / (4 * H ** 2) + tau * dm0 / (2 * H)
    phi0 = trunc(a0 * psi0, n)
    dphi0 = trunc(da0 * psi0 + a0 * dpsi0, n)
    pref = vp / H ** 2 + tau * mp / H
    dpref = dvp / H ** 2 + tau * dmp / H
    col = 0 if tau > 0 else 2
    phi1 = trunc(pref * R[1][col], n)
    dphi1 = trunc(dpref * R[1][col] + pref * Rd[1][col], n)

    logr = slog1p(trunc(trunc(phi1 / phi0, n) - 1, n), n)
    dlogr = trunc(dphi1 / phi1 - dphi0 / phi0, n)
    gap = trunc(psi0 - 2 * R[1][1], n)
    dgap = trunc(dpsi0 - 2 * Rd[1][1], n)
    return logr, dlogr, gap, dgap, phi1


def draw(seed):
    rnd = random.Random(seed)
    v = [sp.Rational(rnd.randint(3, 90), rnd.randint(2, 40)) for _ in range(NT)]
    m = [sp.Rational(rnd.randint(-60, 60), rnd.randint(2, 40)) for _ in range(NT)]
    dv = [sp.Rational(rnd.randint(-40, 40) or 7, rnd.randint(2, 30)) for _ in range(NT)]
    dm = [sp.Rational(rnd.randint(-40, 40) or -5, rnd.randint(2, 30)) for _ in range(NT)]
    return v, m, dv, dm


def co(e, k):
    return sp.nsimplify(sp.simplify(sp.expand(e).coeff(H, k)))


def pfun(v, m):
    v0, v1, v2 = v[0], v[1], 2 * v[2]
    m0, m1 = m[0], m[1]
    return (4 * m0 ** 2 - 2 * v1 * m0 + 4 * v0 * m1 + v2 * v0 - 2 * v1 ** 2) / v0 ** 3


def bfun(v, m, om):
    """The coefficient of h^4 in the mismatch: q h^4 b is the leading correction."""
    v0, v1 = v[0], v[1]
    m0 = m[0]
    return (sp.I * om * (pfun(v, m) - 2 * m0 * v1 / v0 ** 3)
            + 2 * om ** 2 / v0 ** 2 - 2 * m0 * rfun(v, m) / v0)


def rfun(v, m):
    v0, v1, v2 = v[0], v[1], 2 * v[2]
    m0, m1, m2 = m[0], m[1], 2 * m[2]
    return (m2 * v0 ** 2 - 4 * m0 ** 3 + 2 * v1 * m0 ** 2 - 2 * v1 * v0 * m1
            - (2 * v0 * m1 + v2 * v0 - 2 * v1 ** 2) * m0) / (2 * v0 ** 3)


def run(seed, verbose=False):
    v, m, dv, dm = draw(seed)
    om = sp.Rational(7, 5)
    n = ORD + 3
    plus = symbols(v, m, dv, dm, om, +1, n)
    minus = symbols(v, m, dv, dm, om, -1, n)

    if verbose:
        # (3.45): the printed expansion of phi1, recomputed as a check on the machinery
        v0, v1, v2 = v[0], v[1], 2 * v[2]
        m0, m1, m2 = m[0], m[1], 2 * m[2]
        r = rfun(v, m)
        p = (4 * m0 ** 2 - 2 * v1 * m0 + 4 * v0 * m1 + v2 * v0 - 2 * v1 ** 2) / v0 ** 3
        ref = [sp.Integer(1), 2 * m0 / v0, -4 * sp.I * om / v0,
               -(8 * m0 - v1) * sp.I * om / v0 ** 2 + r,
               sp.I * om * p - 14 * om ** 2 / v0 ** 2]
        got = plus[4]
        ok = all(sp.simplify(co(got, k) - ref[k]) == 0 for k in range(5))
        print('    (3.45) phi1 through h^4        : %s' % ('agrees' if ok else 'MISMATCH'))

    table = []
    for k in range(1, ORD + 1):
        ep, dp = co(plus[0], k), co(plus[1], k)
        em, dm_ = co(minus[0], k), co(minus[1], k)
        ev, od = sp.simplify((ep + em) / 2), sp.simplify((ep - em) / 2)
        dev, dod = sp.simplify((dp + dm_) / 2), sp.simplify((dp - dm_) / 2)
        table.append((k, ev != 0, od != 0, dev != 0, dod != 0))

    c_val = sp.simplify((co(plus[0], 3) - co(minus[0], 3)) / 2)
    c_der = sp.simplify((co(plus[1], 3) - co(minus[1], 3)) / 2)

    v0 = v[0]
    # the coefficient of h^3 in the first line of (3.47), divided by tau
    c_closed = sp.I * om * v[1] / v0 ** 2 + rfun(v, m)
    vt = [v[k] + T * dv[k] for k in range(3)]
    mt = [m[k] + T * dm[k] for k in range(3)]
    dc_closed = (sp.I * om * (dv[1] / v0 ** 2 - 2 * v[1] * dv[0] / v0 ** 3)
                 + sp.diff(rfun(vt, mt), T).subs(T, 0))

    closed = (sp.simplify(c_val - c_closed) == 0, sp.simplify(c_der - dc_closed) == 0)

    b_val = sp.simplify(co(plus[0], 4))
    b_der = sp.simplify(co(plus[1], 4))
    vt3 = [v[k] + T * dv[k] for k in range(3)]
    mt3 = [m[k] + T * dm[k] for k in range(3)]
    b_closed = bfun(v, m, om)
    db_closed = sp.diff(bfun(vt3, mt3, om), T).subs(T, 0)
    bclosed = (sp.simplify(b_val - b_closed) == 0, sp.simplify(b_der - db_closed) == 0)

    gap = [(k, co(plus[2], k) != 0, co(plus[3], k) != 0) for k in range(2, ORD + 2)]
    return table, closed, (c_val, c_der), gap, bclosed, (b_closed, db_closed)


def main():
    print('  one dimensional per-step symbols, rebuilt from the definitions of Chapter 5')
    print('  generic rational data, exact arithmetic; omega = 7/5')
    print('  d(R) = -R d(M) R carries the derivative in the coefficients\n')

    runs = [run(seed, verbose=(i == 0)) for i, seed in enumerate((11, 23, 47))]

    print('\n  --- log(phi1/phi0): parity in tau, and the parity of its derivative ---')
    print('    h^k | even(tau) | odd(tau) | d(even) | d(odd)')
    for (k, ev, od, dev, dod) in runs[0][0]:
        print('    %-3d |  %-7s | %-8s | %-7s | %s'
              % (k, 'yes' if ev else '0', 'yes' if od else '0',
                 'yes' if dev else '0', 'yes' if dod else '0'))
    same = all(r[0] == runs[0][0] for r in runs)
    print('    parities identical across the three draws: %s' % ('yes' if same else 'NO'))

    print('\n  --- c(x), the h^3 coefficient divided by tau ---')
    for i, r in enumerate(runs):
        print('    draw %d: c(x) %s the closed form;  d/dtheta c %s it'
              % (i + 1, 'equals' if r[1][0] else 'DIFFERS FROM',
                 'equals' if r[1][1] else 'DIFFERS FROM'))
    print('    value at draw 1      : %s' % sp.nsimplify(runs[0][2][0]))
    print('    derivative at draw 1 : %s' % sp.nsimplify(runs[0][2][1]))

    print('\n  --- the leading correction b(x), the coefficient of h^4 ---')
    for i, r in enumerate(runs):
        print('    draw %d: b(x) %s the closed form;  d/dtheta b %s it'
              % (i + 1, 'equals' if r[4][0] else 'DIFFERS FROM',
                 'equals' if r[4][1] else 'DIFFERS FROM'))
    print('    closed form:  b = i omega [p - 2 mu v1/v0^3] + 2 omega^2/v0^2 - 2 mu r/v0,')
    print('                  the whole of it even in tau, so it does not telescope and')
    print('                  q h^4 b = h^2 b is the leading correction to the constant.')
    print('    at draw 1   : b = %s' % sp.nsimplify(runs[0][5][0]))
    print('                  b\' = %s' % sp.nsimplify(runs[0][5][1]))

    print('\n  --- the endpoint factor (3.50): psi0 - 2 psi1 ---')
    print('    h^k | psi0-2psi1 | d/dtheta(psi0-2psi1)')
    for (k, g, dg) in runs[0][3]:
        print('    %-3d |  %-9s | %s' % (k, 'nonzero' if g else '0', 'nonzero' if dg else '0'))

    print('\n  --- verdict -------------------------------------------------------')
    t = runs[0][0]
    val_ok = all(not ev for (k, ev, od, dev, dod) in t if k % 2 == 1) and \
             all(not od for (k, ev, od, dev, dod) in t if k % 2 == 0)
    der_ok = all(not dev for (k, ev, od, dev, dod) in t if k % 2 == 1) and \
             all(not dod for (k, ev, od, dev, dod) in t if k % 2 == 0)
    gap3 = [g for (k, g, dg) in runs[0][3] if k == 3]
    dgap3 = [dg for (k, g, dg) in runs[0][3] if k == 3]
    print('    parity of the mismatch           : %s' % ('odd at odd orders only' if val_ok else 'FAILS'))
    print('    parity after differentiation     : %s' % ('the same' if der_ok else 'FAILS'))
    print('    c(x) and its derivative          : %s'
          % ('both match the closed form of (3.47)'
             if all(r[1][0] and r[1][1] for r in runs) else 'do not match'))
    print('    b(x) and its derivative          : %s'
          % ('both match the closed form'
             if all(r[4][0] and r[4][1] for r in runs) else 'do not match'))
    print('    psi0 - 2 psi1 at h^3             : %s, and its derivative %s'
          % ('0' if gap3 and not gap3[0] else 'nonzero',
             '0' if dgap3 and not dgap3[0] else 'nonzero'))
    print('    so the differentiated mismatch is again a sum of differences of functions')
    print('    of the site, and telescopes exactly as the undifferentiated one does.')


if __name__ == '__main__':
    main()
\end{lstlisting}

\section{\texttt{massstep.py}}

Chapter 10: what the coding of the step does to the mass.

\begin{lstlisting}
"""Chapter 10: what the coding of the step does to the mass.

The explicit step can be written as an update of the values,

    u <- u + dt * Lambda u,

or as a weighted average of the three values of the previous date,

    u_i <- p_- u_{i-1} + p_0 u_i + p_+ u_{i+1},   p_+-+p_0+p_- = 1 exactly.

The two are algebraically the same and numerically are not.  The weighted average
rounds the weights, and their rounded sum is no longer one, so every step creates
or destroys mass at the level of the unit roundoff and the defect accumulates
linearly in the number of steps.  The update applies differences, which vanish on
constants in floating point as they do in exact arithmetic, so the mass error has
no systematic part.

The script evolves a unit mass at a single site on a periodic lattice with
constant coefficients, in single and in double precision, by each of the two
formulations, and prints the total mass after the number of steps stated.  The
horizon is t = steps * h^2, the regime dt = h^2 of the Euler parametrisation.

Run:  python3 massstep.py
"""
import numpy as np


def run(n, steps, dtype, form):
    h = 1.0 / n
    dt = h * h
    v = dtype(0.5)
    b = dtype(0.2)
    one = dtype(1)
    u = np.zeros(n, dtype=dtype)
    u[n // 2] = one / dtype(h)
    for _ in range(steps):
        if form == 'update':
            out = (one - dt * v / dtype(h * h)) * u
            out = out + dt * (v / (dtype(2) * dtype(h * h))) * (np.roll(u, -1) + np.roll(u, 1))
            out = out + dt * (b / (dtype(2) * dtype(h))) * (np.roll(u, -1) - np.roll(u, 1))
        else:
            p0 = one - dt * v / dtype(h * h)
            pm = dt * (v / (dtype(2) * dtype(h * h))) - dt * (b / (dtype(2) * dtype(h)))
            pp = dt * (v / (dtype(2) * dtype(h * h))) + dt * (b / (dtype(2) * dtype(h)))
            out = p0 * u + pm * np.roll(u, 1) + pp * np.roll(u, -1)
        u = out
    return float(np.sum(u, dtype=np.float64)) * float(h)


def main():
    print('  v = 0.5, b = 0.2, dt = h^2, unit mass at a single site, periodic lattice')
    print('  total mass after the stated number of steps, printed as (mass - 1)\n')
    print('  precision |   n  |  steps |        update form |  weighted-average form')
    for dtype, name in ((np.float32, '  single  '), (np.float64, '  double  ')):
        for n, steps in ((128, 1024), (512, 4096)):
            a = run(n, steps, dtype, 'update')
            c = run(n, steps, dtype, 'weights')
            print('  %s | %4d | %6d | %+.3e         | %+.3e' % (name, n, steps, a - 1, c - 1))
    print('\n  the weighted form drifts by two orders of magnitude more than the update')
    print('  form in single precision, and by nothing systematic in double precision.')


if __name__ == '__main__':
    main()
\end{lstlisting}

\section{\texttt{decorating\_paths.py}}

Sample a decorating path and print the picture of Figure 5.3.

\begin{lstlisting}
"""Sample a decorating path and print the picture of Figure 5.3.

The coarse walk is nearest neighbour on h_m Z, the fine one on h_{m+1}Z at four times the
rate, since h_m = 2h_{m+1} and the rate scales as h^{-2}. A decoration is then sampled as the
fine walk *conditioned* to satisfy Definition 5.1: it must stay within one fine cell of the
coarse value it last took until it takes the next one. The constraint binds between the times
of the fine path, so the touches of the coarse values fall where the fine walk happens to
make them, not at the coarse step times.

The sample is checked against the definition node by node before it is drawn, and the
unconditioned fine walk is drawn alongside it, leaving the band: the decoration is a
conditioned walk, not a generic one. The picture, not its caption, is written to
sections/ch5-decorating-picture.tex.
"""
import os
import random

H = 1.0          # h_{m+1}: the fine mesh, in units of the picture; the coarse mesh h_m is 2
RATE = 0.45      # coarse hopping rate, per unit of picture time
NSTEP = 4        # coarse steps in the picture
SY = 0.88        # the space axis is drawn at this fraction of the time axis


def coarse_walk(rng, n):
    """(hop times, levels) of a nearest-neighbour walk of n steps on the coarse mesh."""
    t, y = 0.0, 0
    times, levels = [0.0], [0]
    for _ in range(n):
        t += rng.expovariate(RATE)
        y += 2 * rng.choice((-1, 1))
        times.append(t)
        levels.append(y)
    return times, levels


def decoration(rng, levels):
    """The fine walk conditioned on Definition 5.1, as (time, level) breakpoints and touches.

    In the band of the coarse value a = levels[i] the walk lives on the three fine sites
    a-1, a, a+1; the landing site a-2s, s the direction of the next coarse value, is
    conditioned away, which is the constraint, and the landing on a+2s ends the phase. The
    holding time is exponential with the total rate of the moves which survive.
    """
    rate = 4.0 * RATE
    t, brk, touches = 0.0, [(0.0, levels[0])], [0.0]
    for i in range(len(levels) - 1):
        a, target = levels[i], levels[i + 1]
        s = 1 if target > a else -1
        x = a
        while True:
            moves = [d for d in (-1, 1) if x + d != a - 2 * s]
            t += rng.expovariate(rate / 2.0 * len(moves))
            x += rng.choice(moves)
            brk.append((t, x))
            if x == target:
                touches.append(t)
                break
    return brk, touches


def free_walk(rng, a, target, hops):
    """The unconditioned fine walk from a, run until it leaves the band on the wrong side."""
    s = 1 if target > a else -1
    t, x, brk = 0.0, a, [(0.0, a)]
    for _ in range(hops):
        t += rng.expovariate(4.0 * RATE)
        x += rng.choice((-1, 1))
        brk.append((t, x))
        if x == a - 2 * s:
            return brk
    return None


def check(brk, levels):
    """Definition 5.1, checked node by node on the sampled fine path."""
    seq = [y for _, y in brk]
    for i in range(len(levels) - 1):
        a, target = levels[i], levels[i + 1]
        j = next(k for k in range(len(seq)) if seq[k] == a
                 and (i == 0 or any(seq[m] == levels[i - 1] for m in range(k + 1))))
        k = next(m for m in range(j + 1, len(seq)) if seq[m] == target)
        for n in range(j + 1, k):
            assert abs(seq[n] - a) <= H, 'Definition 5.1 violated at node %d' % n
        for n in range(1, len(seq)):
            assert abs(seq[n] - seq[n - 1]) == H, 'not a nearest-neighbour step at %d' % n
    for a in levels:
        assert a in seq, 'the coarse value %d is absent from the decoration' % a
    return True


def fmt(pts):
    return ' -- '.join('(%.2f,%.2f)' % (t, y * SY) for t, y in pts)


def staircase(times, levels):
    """Breakpoints of a walk drawn as a path which is constant in time between its steps.

    The domain of the walk is the time axis, so a step is vertical: a hop from y to y' at time t
    contributes the two points (t, y) and (t, y'), and the segment between consecutive hops is
    horizontal.
    """
    pts = [(times[0], levels[0])]
    for i in range(1, len(times)):
        pts.append((times[i], levels[i - 1]))
        pts.append((times[i], levels[i]))
    return pts


def steps(brk):
    """The same, for a path given as (time, level) landings."""
    return staircase([t for t, _ in brk], [y for _, y in brk])


def main():
    for seed in range(1, 400):
        rng = random.Random(11 * seed)
        times, levels = coarse_walk(rng, NSTEP)
        if max(levels) - min(levels) < 4 or levels[-1] != 0:
            continue
        if len({levels[i + 1] - levels[i] for i in range(NSTEP)}) != 2:
            continue
        brk, touches = decoration(random.Random(11 * seed + 1), levels)
        # a sample is kept only if it is legible: the picture must not be stretched by a long
        # dwell of the fine path in one band, nor crowded by a coarse step which the fine path
        # takes at once, and the coarse walk must wander in both directions
        if not 11 <= len(brk) <= 20:
            continue
        if touches[-1] > 9.5 or max(touches[i + 1] - touches[i] for i in range(NSTEP)) > 2.6:
            continue
        # nor is a picture legible in which two steps of the same path fall closer together than
        # the labels and the vertical segments can bear
        if min(b - a for a, b in zip(touches, touches[1:])) < 0.30:
            continue
        if min(b - a for a, b in zip(times, times[1:])) < 0.62:
            continue
        if max(levels) - min(levels) < 4:
            continue
        break
    else:
        raise SystemExit('no sample found')
    check(brk, levels)

    # the unconditioned walk, from a touch of a coarse value which it leaves on the wrong side
    dash = None
    for i in range(1, len(levels) - 1):
        for hops in (3, 4, 2, 5, 6):
            w = free_walk(random.Random(11 * seed + 3 + 7 * i + hops), levels[i],
                          levels[i + 1], hops)
            if w is not None:
                dash = [(touches[i] + t, y) for t, y in w]
                break
        if dash is not None:
            break

    T = max([times[-1], touches[-1]] + ([dash[-1][0]] if dash else []))
    # the extents include the bands, which reach one cell beyond the coarse values
    seen = [y for _, y in brk] + [y for _, y in dash or []]
    lo = min(seen + [min(levels) - 1])
    hi = max(seen + [max(levels) + 1])
    fine_hops = len(brk) - 1
    print('seed %d: %d coarse steps, %d fine steps (ratio %.1f), %d bands, span %d to %d'
          % (seed, NSTEP, fine_hops, fine_hops / float(NSTEP), len(levels) - 1, lo, hi))
    print('Definition 5.1 checked node by node on both decorations')
    print('coarse values %s at times %s' % (levels, [round(t, 2) for t in times]))
    print('touches %s' % [round(t, 2) for t in touches])
    print('the decoration ends at level %d' % brk[-1][1])
    print('dashed from %s' % [(round(t, 2), y) for t, y in dash])
    print('the unconditioned walk leaves the band at t = %.2f, level %d'
          % (dash[-1][0], dash[-1][1]) if dash else 'no unconditioned walk drawn')

    L = []
    w = L.append
    w('%% Generated by scripts/ch5/decorating_paths.py - do not edit.')
    w('%% space is drawn at %.2f of the time axis; the coarse walk moves 2 units of space per' % SY)
    w('%% step of the picture and the fine one 1, so h_m = 2h_{m+1} throughout.')
    w('\\begin{tikzpicture}[>=latex,line join=round,scale=1.0,every node/.style={font=\\footnotesize}]')
    # the bands, from the touch of a coarse value to the touch of the next one
    for i in range(len(levels) - 1):
        w('  \\fill[orange!16] (%.2f,%.2f) rectangle (%.2f,%.2f);'
          % (touches[i], (levels[i] - 1) * SY, touches[i + 1], (levels[i] + 1) * SY))
        w('  \\draw[orange!45,thin] (%.2f,%.2f) rectangle (%.2f,%.2f);'
          % (touches[i], (levels[i] - 1) * SY, touches[i + 1], (levels[i] + 1) * SY))
    bot, top = (lo - 0.4) * SY, (hi + 0.15) * SY
    w('  \\draw[->] (0,%.2f) -- (%.2f,%.2f) node[right] {time};'
      % (bot, T + 0.8, bot))
    w('  \\draw[->] (0,%.2f) -- (0,%.2f);' % (bot, top))
    w('  \\node[left] at (0,%.2f) {space};' % top)
    # the coarse step times, unless two of them fall too close together for the labels
    for i in range(1, len(times)):
        w('  \\draw[thick] (%.2f,%.2f) -- (%.2f,%.2f);'
          % (times[i], bot, times[i], bot - 0.16))
        w('  \\node[below] at (%.2f,%.2f) {$t_{%s}$};'
          % (times[i], bot - 0.2, 'j' if i == 1 else 'j+%d' % (i - 1)))
    w('  \\node[left] at (-0.08,%.2f) {$0$};' % bot)
    coarse = staircase(times, levels) + [(T, levels[-1])]
    w('  \\draw[very thick] %s;' % fmt(coarse))
    w('  \\draw[thick,blue!85!black] %s;' % fmt(steps(brk)))
    if dash:
        w('  \\draw[thick,dashed,gray] %s;' % fmt(steps(dash)))
        w('  \\node[gray!40!black,font=\\large,fill=white,inner sep=1pt] at (%.2f,%.2f) {$\\times$};'
          % (dash[-1][0] + 0.06, dash[-1][1] * SY))
    for i in range(len(levels)):
        w('  \\fill[white] (%.2f,%.2f) circle (2.2pt);' % (touches[i], levels[i] * SY))
        w('  \\draw[thick,blue!85!black] (%.2f,%.2f) circle (2.2pt);'
          % (touches[i], levels[i] * SY))
    # the mesh sizes, on the space axis
    Y0 = min(levels)
    w('  \\draw[<->] (-0.55,%.2f) -- (-0.55,%.2f);' % (Y0 * SY, (Y0 + 2) * SY))
    w('  \\node[left] at (-0.65,%.2f) {$h_m$};' % ((Y0 + 1) * SY,))
    w('  \\draw[<->] (-0.55,%.2f) -- (-0.55,%.2f);' % ((Y0 + 2) * SY, (Y0 + 3) * SY))
    w('  \\node[left] at (-0.65,%.2f) {$h_{m+1}$};' % ((Y0 + 2.5) * SY,))
    # the legend, above the axis: two rows, so that it stays within the width of the picture
    for row, entries in enumerate(([('very thick,black', 'coarse path, step $h_m$'),
                                    ('thick,blue!85!black', 'decorating fine path'),
                                    ('thick,dashed,gray', 'fine walk, undecorated')],)):
        Y, x = top + 0.45 - 0.5 * row, 0.0
        for style, text in entries:
            w('  \\draw[%s] (%.2f,%.2f) -- (%.2f,%.2f);' % (style, x, Y, x + 0.5, Y))
            w('  \\node[anchor=west] at (%.2f,%.2f) {%s};' % (x + 0.6, Y, text))
            x += 0.6 + 0.145 * len(text) + 0.55
    w('\\end{tikzpicture}')

    mono = os.path.dirname(os.path.dirname(os.path.dirname(os.path.abspath(__file__))))
    path = os.path.join(mono, 'sections', 'ch5-decorating-picture.tex')
    with open(path, 'w') as f:
        f.write('\n'.join(L) + '\n')
    print('written to %s' % path)


main()
\end{lstlisting}

% ---- end generated/appendix.tex ----
\backmatter
% ---- from bibliography.tex ----
% Generated by build.py - do not edit.

% ---- end bibliography.tex ----
\end{document}